\documentclass[AnnalsOfMath]{PUP-math}

\usepackage{enumitem}
\setlist{noitemsep}
\setlist[enumerate]{align=left, left=0pt, labelsep=0pt}
\setlist[enumerate,1]{label*=(\roman*)}

\theoremstyle{plain} 
\newtheorem{prop}[theorem]{Proposition}
\newtheorem{cor}[theorem]{Corollary}

\theoremstyle{definition}

\newtheorem{rmk}[theorem]{Remark}
\newtheorem{conj}[theorem]{Conjecture}
\newtheorem{fakedefinition}[theorem]{`Definition'}

\newcommand{\fz}{\mathfrak{z}}

\newcommand{\cA}{\mathcal{A}}
\newcommand{\cF}{\mathcal{F}}

\newcommand{\calD}{\mathcal{D}}
\newcommand{\cO}{\mathcal{O}}

\newcommand{\cE}{\mathcal{E}}

\newcommand{\cW}{\mathcal{W}}

\newcommand{\U}{\mathcal{U}}
\newcommand{\calL}{\mathcal{L}}
\newcommand{\cB}{\mathcal{B}}
\newcommand{\B}{\mathcal{B}}
\newcommand{\cN}{\mathcal{N}}

\newcommand{\g}{\mathfrak{g}}
\newcommand{\C}{\mathbb{C}}
\newcommand{\F}{\mathbb{F}}
\newcommand{\A}{\mathcal{A}}
\newcommand{\Orb}{\mathbb{O}}

\newcommand{\unip}{\mathrm{Unip}}

\newcommand{\OHC}{\mathrm{OHC}}
\renewcommand{\AA}{\mathbb{A}}

\newcommand{\OO}{\mathbb{O}}

\newcommand{\CC}{\mathbb{C}}

\newcommand{\PP}{\mathbb{P}}
\newcommand{\QQ}{\mathbb{Q}}
\newcommand{\RR}{\mathbb{R}}

\newcommand{\ZZ}{\mathbb{Z}}

\newcommand{\fs}{\mathfrak{s}}
\newcommand{\fS}{\mathfrak{S}}
\newcommand{\fh}{\mathfrak{h}}
\newcommand{\fu}{\mathfrak{u}}
\newcommand{\fq}{\mathfrak{q}}
\newcommand{\fg}{\mathfrak{g}}

\newcommand{\fL}{\mathfrak{L}}
\newcommand{\fl}{\mathfrak{l}}
\newcommand{\fm}{\mathfrak{m}}

\newcommand{\fn}{\mathfrak{n}}

\newcommand{\fZ}{\mathfrak{Z}}
\newcommand{\fX}{\mathfrak{X}}

\newcommand{\fP}{\mathfrak{P}}
\newcommand{\fo}{\mathfrak{o}}
\newcommand{\fp}{\mathfrak{p}}

\newcommand{\vardbtilde}[1]{\tilde{\raisebox{0pt}[0.85\height]{$\tilde{#1}$}}}

\DeclareMathOperator{\Vect}{Vect}
\DeclareMathOperator{\Ann}{Ann}

\DeclareMathOperator{\Pic}{Pic}

\DeclareMathOperator{\Hom}{Hom}

\DeclareMathOperator{\modd}{-mod}

\DeclareMathOperator{\gr}{gr}

\DeclareMathOperator{\Coh}{Coh}

\DeclareMathOperator{\Spec}{Spec}

\DeclareMathOperator{\Cl}{Cl}

\DeclareMathOperator{\Aut}{Aut}

\DeclareMathOperator{\codim}{codim}
\DeclareMathOperator{\Ind}{Ind}

\DeclareMathOperator{\Ad}{Ad}
\DeclareMathOperator{\ad}{ad}

\DeclareMathOperator{\Irr}{Irr}

\DeclareMathOperator{\Prim}{Prim}

\DeclareMathOperator{\Ker}{Ker}

\DeclareMathOperator{\Stab}{Stab}
\DeclareMathOperator{\HC}{HC}

\DeclareMathOperator{\End}{End}

\DeclarePairedDelimiter\floor{\lfloor}{\rfloor}

\copyrightyear{2024}

\usepackage{natbib}
\makeindex
\begin{document}

\frontmatter

\AnnalsNumber=218

\makehalftitle



\title{Unipotent Ideals and Harish-Chandra Bimodules}


\author{Ivan Losev\\
Lucas Mason-Brown\\
Dmytro Matvieievskyi}

\maketitlepage


\copyrightyear{2025}
\locNumber{}
\isbnNumber{}
\isbnPbkNumber{}
\Editorial{}
\ProductionEditorial{}
\TextDesign{}
\JacketCoverDesign{}
\JacketCoverCredit{}
\Production{}
\Publicity{}
\bookcomposed{}
\typesettingPrintingInformation{}

\copyrightpage




\tableofcontents


\begin{preface}

\section{Unitary representations}\label{SS_Preface_unitary}
Let $G_{\RR}$ be a Lie group. A \emph{unitary representation} of $G_{\RR}$ is a pair $(\pi, V)$ consisting of a complex Hilbert space $V$ and a group homomorphism
$$\pi: G_{\RR} \to U(V),$$
where $U(V)$ denotes the group of invertible unitary operators $V \to V$. We require that the map 
$$G_{\RR} \times V \to V, \qquad (g,v) \mapsto \pi(g)v$$
is continuous. An \emph{invariant subspace} of $(\pi,V)$ is a  linear subspace of $V$ which is preserverd by the operators $\{\pi(g) \mid g \in G_{\RR}\}$. We say that $(\pi,V)$ is \emph{irreducible} if it has no closed invariant subspaces, apart from $\{0\}$ and $V$. Let $\mathrm{Irr}_u(G_{\RR})$ denote the set of all irreducible unitary $G_{\RR}$-representations, up to unitary equivalence. A fundamental unsolved problem in representation theory, with its origins in Gelfand's program of abstract harmonic analysis, is to parameterize this set. This is the problem which motivates the present work. 


Suppose for the moment that $G_{\RR}$ is a \emph{simply connected nilpotent Lie group}. An example is the group of upper triangular real $n$-by-$n$ matrices with ones along the diagonal. In this case, Kirillov (\cite{Kirillov1962}) has found a near-perfect solution to the problem described above. Kirillov's result---now known as the `orbit method'---is as follows. Let $\fg_{\RR}$ denote the Lie algebra of $G_{\RR}$. Then $G_{\RR}$ acts on $\fg_{\RR}$ (by the adjoint action $\Ad$) and hence on its linear dual $\fg^*_{\RR}$. Orbits for the latter action are called `co-adjoint orbits'. For each $\xi \in \fg^*_{\RR}$, there is a maximal Lie subalgebra $\fh_{\xi} \subset \fg_{\RR}$ such that $\xi|_{[\fh_{\xi},\fh_{\xi}]}=0$ (this subalgebra is not unique, but the construction we are describing will not depend on its choice). This Lie subalgebra integrates to a Lie subgroup of $G_{\RR}$, which we will denote by $H_\xi$. Because of our assumptions on $G_{\RR}$, $H_{\xi}$ is simply connected. Hence, there is a unique Lie group homomorphism $\chi: H_{\xi} \to S^1$ such that
$$d\chi(X) = 2\pi i \xi(X), \qquad \forall X \in \fh_{\xi}.$$
There is a notion, due to Mackey (\cite{Mackey1}), of \emph{unitary induction}, which associates to a closed subgroup $H_{\RR} \subset G_{\RR}$ and unitary representation $V$ of $H_{\RR}$, a unitary representation of $G_{\RR}$, denoted $\Ind^{G_{\RR}}_{H_{\RR}} V$. Applying this construction to the closed subgroup $H_{\xi} \subset G_{\RR}$ and its one-dimensional unitary representation $\chi$, we obtain a unitary $G_{\RR}$-representation $\Ind^{G_{\RR}}_{H_{\xi}} \chi$. This representation is irreducible, and is determined, up to unitary isomorphism, by the $G_{\RR}$-orbit of $\xi$. Thus, we obtain a map
\begin{equation}\label{eq:Kir}\mathrm{Kir}: \fg^*_{\RR}/G_{\RR} \to \mathrm{Irr}_u(G_{\RR}), \qquad \mathrm{Kir}(G_{\RR} \cdot \xi) = \Ind^{G_{\RR}}_{H_{\xi}} \chi\end{equation}
from the set of $G_{\RR}$-orbits on $\fg_{\RR}^*$ to $\Irr_u(G_{\RR})$. Kirillov proved in \cite{Kirillov1962} that this map is a bijection. 

This kind of result --- establishing a complete and constructive classification of $\Irr_u(G_{\RR})$ in terms of simple geometric data --- can be pushed a bit beyond the case of simply connected nilpotent Lie groups, but not very far. Kostant and Auslander proved in \cite{AuslanderKostant} an analogous result for simply connected solvable Lie groups. For non-solvable Lie groups, there is no known classification of $\Irr_u(G_{\RR})$ in the vein of Kirillov-Auslander-Kostant. However, for \emph{algebraic Lie groups} (i.e. groups of real points of algebraic groups over $\RR$), Duflo showed in \cite{DufloMackey} that the classification of $\Irr_u(G_{\RR})$ reduces to the classification $\Irr_u(H_{\RR})$ for various \emph{reductive} subgroups $H_{\RR} \subset G_{\RR}$. 

Duflo's result suggests that we should try to extend the orbit method to the case of reductive groups. This idea was explored by Vogan in \cite{Vogan1987,Vogan1992,Vogan_dmkrev}. In contrast to the case of solvable groups, the co-adjoint orbits of a reductive group $G_{\RR}$ exhibit an extraordinary amount of structure. A co-adjoint $G_{\RR}$-orbit is called \emph{nilpotent} if it is stable under the scaling action of $\RR_{>0}$ on $\fg^*_{\RR}$. The set of nilpotent orbits in $\fg^*_{\RR}$ is finite and admits a complete classification in terms of $\mathfrak{sl}_2$-triples. By the \emph{Jordan decomposition}, every co-adjoint orbit is constructed from a nilpotent co-adjoint orbit for a Levi subgroup of $G_{\RR}$ (i.e. the centralizer in $G_{\RR}$ of a semisimple element of $\fg_{\RR}$). Vogan has argued that there should be a parallel statement for unitary representations: there should be a finite set of irreducible unitary representations --- to be called `unipotent representations' --- related to nilpotent co-adjoint orbits, such that every irreducible unitary representation is obtained, by several inductive procedures, from a unipotent representation of a suitable Levi subgroup. Vogan has explained in essence how these inductive procedures should work; conjecturally, the crucial remaining obstacle to classifying $\Irr_u(G_{\RR})$ is defining a suitable set of unipotent representations. There are many previous results on unipotent representations, including definitions and constructions of unipotent representations in various special cases due to Adams, Barbasch, Brylinski, McGovern, Vogan, and others. These works will be reviewed in Chapter \ref{sec:intro}. In this monograph, we will give a general definition of unipotent representations in the case of \emph{complex} reductive groups. 

In the 1950s (\cite{HarishChandraI}), Harish-Chandra developed a powerful array of tools for reducing the study of unitary representations of reductive groups to algebra. Let $K_{\CC}$ be the complexification of a maximal compact subgroup $K_{\RR} \subset G_{\RR}$ and let $\fg_{\CC}$ be the complexification of $\fg_{\RR}$. A {\it $(\fg_{\CC},K_{\CC})$-module} is a $\fg_{\CC}$-module $X$ equipped with a compatible locally-finite $K_{\CC}$-action. Here, `compatible' means that the action of $\mathfrak{k}_\CC = \mathrm{Lie}(K_{\CC})$ on $X$ obtained by differentiating the $K_{\CC}$-action coincides, via restriction, with the action of $\fg_{\CC}$. We say that $X$ is \emph{unitary} if it admits a positive-definite Hermitian form such that $\fg_{\RR}$ acts on $X$ by skew-Hermitian operators. One of Harish-Chandra's remarkable results is that there is a natural bijection
\begin{equation}\label{eq:HCpreface}\Irr_u(G_{\RR}) \longleftrightarrow \{\text{irreducible unitary $(\fg_{\CC},K_{\CC})$-modules}\}/\simeq.\end{equation}
The forward map is `passage to locally-finite smooth  vectors'. The inverse map is completion with respect to the metric defined by the form. In view of the bijection (\ref{eq:HCpreface}), the problem of classifying $\Irr_u(G_{\RR})$ is equivalent to the problem of classifying irreducible unitary $(\fg_{\CC},K_{\CC})$-modules. Similarly, the problem of defining a set of unipotent $G_{\RR}$-representations is equivalent to the problem of defining a set of unipotent $(\fg_{\CC},K_{\CC})$-modules.

Now let us assume that $G_{\RR}$ is the Lie group underlying a complex connected reductive algebraic group $G$. Then a $(\fg_{\CC},K_{\CC})$-module is the same thing as a \emph{Harish--Chandra bimodule}, i.e. a bimodule for the universal enveloping algebra $U(\fg)$ of $\fg = \mathrm{Lie}(G)$ such that the adjoint action of $\fg$ integrates to $G$. In this monograph, we will define the notion of a \emph{unipotent Harish-Chandra bimodule}. Our definition is geometric in nature, relying on the theory of quantizations of conical symplectic singularities. An overview will be provided in Chapter \ref{subsec:quantizationspreface} below. It is not at all obvious that our unipotent bimodules are unitary, but it is now known that they are in complete generality (\cite{DavisMasonBrown}). Thus, our unipotent bimodules correspond (under the bijection (\ref{eq:HCpreface})) to a finite set of unitary representations of $G_{\RR}$. We expect that these representations form the building blocks of $\Irr_u(G_{\RR})$. 

\section{Quantizations}\label{subsec:quantizationspreface}

As mentioned above, our approach to defining unipotent Harish-Chandra bimodules is based on the idea of quantization. In this section, we will elaborate on this idea and its connections to the theory of unitary representations.

`Quantization' is a general term which refers to any of several procedures for passing from the realm of classical mechanics to the realm of quantum mechanics. {In both classical and quantum mechanics, a central role is played by the `phase space', i.e. the space of all possible states of the system.} The phase space of a classical mechanical system is a Poisson manifold, i.e. a smooth manifold $M$ equipped with a Poisson bracket $\{\cdot,\cdot\}:C^\infty(M)\times C^\infty(M)\rightarrow C^\infty(M)$. On the other hand, the phase space of a quantum mechanical system is a complex Hilbert space. In both cases, there is a `Hamiltonian' which governs the dynamics of the system. But we will ignore this additional structure for the purposes of this discussion.

To `quantize' a Poisson manifold $M$ is to construct a Hilbert space $V$ which represents the same underlying system. It is not always possible to do so, and when it is, this Hilbert space may not be unique. Nonetheless, we will denote this (ill-defined, in general) correspondence with a squiggly arrow as follows:
\begin{equation}\label{eq:quantpreface1}\{\text{Poisson manifolds}\} \rightsquigarrow \{\text{Hilbert spaces}\}.\end{equation}
Let us now consider what happens when symmetry is present. If a classical mechanical system is symmetric with respect to a Lie group $G_\RR$, then the phase space $M$ admits a \emph{Hamiltonian} $G_{\RR}$-action. This is, in fact, a pair of additional structures: a $G_{\RR}$-action on $M$ preserving the Poisson bracket and a smooth $G_{\RR}$-equivariant map $\mu: M \to \fg_{\RR}^*$ (called a `moment map') such that for each $\xi \in \fg_{\RR}$ the derivation $\{\mu^*(\xi),\cdot\}$ coincides with the vector field $\xi_M$ obtained by differenting the $G_{\RR}$-action on $M$. By a `Hamiltonian $G_{\RR}$-space' we mean a Poisson manifold equipped with a Hamiltonian $G_{\RR}$-action. If a quantum mechanical system possesses $G_{\RR}$-symmetry, then $G_{\RR}$ acts on the phase space by unitary operators. In other words, the phase space has the structure of a unitary $G_{\RR}$-representation. So the `equivariant' version of (\ref{eq:quantpreface1}) is a correspondence of the form
\begin{equation}\label{eq:quantpreface2}\{\text{Hamiltonian $G_{\RR}$-spaces}\} \rightsquigarrow \{\text{unitary $G_{\RR}$-representations}\}.\end{equation}
Suppose we wish to consider physical systems which are `maximally symmetric'. For Hamiltonian $G_{\RR}$-spaces, one natural interpretation of this requirement is that the $G_{\RR}$-action is \emph{transitive}. We call the corresponding spaces `homogeneous Hamiltonian $G_{\RR}$-spaces'. An analogous requirement on the quantum side is that the unitary $G_{\RR}$-representation is \emph{irreducible}. Thus, it is natural to expect that (\ref{eq:quantpreface2}) restricts to a correspondence of the form
\begin{equation}\label{eq:quantpreface3}\{\text{homogeneous Hamiltonian $G_{\RR}$-spaces}\} \rightsquigarrow \mathrm{Irr}_u(G_{\RR}).\end{equation}
It turns out that homogeneous Hamiltonian $G_{\RR}$-spaces are closely related to co-adjoint $G_{\RR}$-orbits, discussed in Section \ref{SS_Preface_unitary}. First, every co-adjoint $G_{\RR}$-orbit $\OO$ is naturally an example of a homogeneous Hamiltonian $G_{\RR}$-space: $\fg_{\RR}^*$ is a Poisson manifold and $\OO$ is a Poisson submanifold. The co-adjoint action of $G_{\RR}$ on $\OO$ is manifestly transitive, and the inclusion $\OO \hookrightarrow \fg_{\RR}^*$ is a moment map. Thus, $\OO$ can be regarded as a homogeneous Hamiltonian $G_{\RR}$-space. We can generalize this example slightly by considering \emph{equivariant covers}. Namely, choose $\xi \in \OO$ and let $H$ denote the stabilizer of $\xi$ in $G_{\RR}$. Let $H'$ be any subgroup intermediate between $H$ and its identity component, and let $\mu'$ denote the $G_{\RR}$-equivariant map $G_{\RR}/H' \to \OO$ given by $gH' \mapsto \Ad^*(g)\xi$. This map is \'{e}tale by the conditions placed on $H'$. So we can lift the Poisson structure from $\OO$ to $G_{\RR}/H'$, and the composition 
$$G_{\RR}/H' \overset{\mu'}{\to} \OO \hookrightarrow \fg_{\RR}^*$$
is a moment map for the $G_{\RR}$-action on $G_{\RR}/H'$. It is an easy but remarkable fact, first observed by Kostant (\cite[Theorem 5.4.1]{KostantQuant}), that \emph{every} homogeneous Hamiltonian $G_{\RR}$-space arises in this fashion. 

If $G_{\RR}$ is nilpotent, then every co-adjoint $G_{\RR}$-orbit is simply connected. So the set $\fg_{\RR}^*/G_{\RR}$ of co-adjoint orbits is precisely the set of homogeneous Hamiltonian $G_{\RR}$-spaces, and Kirillov's bijection (\ref{eq:Kir}) realizes the quantization correspondence (\ref{eq:quantpreface3}).

The picture for reductive groups is considerably more subtle. The structure of the quantum side is much more complicated. However, there are also complications on the classical side. For example, the co-adjoint orbits are no longer simply connected (this is already the case for $\operatorname{SL}_2(\RR)$ and $\operatorname{SL}_2(\CC)$). This observation is extremely important for the present work. 

The quantization problem, as formulated above, is notoriously difficult (as evidenced, for example, by the failure of Kirillov's Orbit Method for reductive groups). One reason is that the classical and quantum sides are very different in nature. However, there is an \emph{algebraic} version of the quantization problem called `deformation quantization' introduced in \cite{BFFLS}, in which the quantum side is a \emph{deformation} of the classical side. This formulation is closer to the perspective we adopt in this monograph, so we will now recall the main ideas.

Let $M$ be a Poisson manifold and $\hbar$ an indeterminate. By a \emph{star-product} on $M$ we mean an associative $\CC[[\hbar]]$-bilinear map
$$\star: C^\infty(M)[[\hbar]]\times C^\infty(M)[[\hbar]] \rightarrow C^\infty(M)[[\hbar]]$$
subject to the following conditions:
\begin{itemize}
\item[(i)] 
$1\in C^\infty(M)$ is a unit,
\item[(ii)] $f\star g-fg\in \hbar C^\infty(M)$ 
\item[(iii)] and $f\star g-g\star f-\hbar\{f,g\}\in \hbar^2 C^\infty(M)$ for all $f,g\in C^\infty(M)$.
\end{itemize}
Note that (i) and (ii) imply that $\star$ defines on $C^{\infty}(M)[[\hbar]]$ the structure of an associative unital algebra, deforming the usual algebra structure on $C^{\infty}(M)$. For this reason, the pair $(C^{\infty}(M)[[\hbar]],\star)$ is called a \emph{deformation quantization} of $C^{\infty}(M)$.

Since we are taking \emph{formal} power series, we cannot specialize $\hbar$ to any nonzero number. However, under suitable conditions (for example, the existence of a multiplicative group action on $M$ appropriately compatible with $\star$), one can specialize $\hbar$ to $1$ for nice (in examples, polynomial) functions on $M$. This gives rise to the structure of a \emph{filtered} associative algebra on the space of polynomial functions. This structure is called a \emph{filtered quantization}. For example, the universal enveloping algebra $U(\fg)$ arises as a filtered quantization of the symmetric algebra $S(\fg)$, for any Lie algebra $\fg$.

The Poisson algebras we will quantize in this monograph are algebras of polynomial functions on certain coadjoint orbits and their equivariant covers. More precisely, let $G$ be a complex reductive algebraic group (e.g. $\operatorname{SL}_n(\CC), \operatorname{SO}_n(\CC)$ or $\operatorname{Sp}_{2n}(\CC)$) and let $\fg$ be its Lie algebra. A co-adjoint $G$-orbit $\OO$ is said to be \emph{nilpotent} if it is stable under the scaling action of $\CC^{\times}$ on $\fg^*$ (equivalently, if it is identified, by a $G$-equivariant isomorphism $\fg^* \simeq \fg$, with a conjugacy class in $\fg$ consisting of nilpotent operators). For reasons discussed above, we will also consider equivariant covers of such orbits (to be called `nilpotent covers'). As mentioned in the previous paragraph, the existence of a $\CC^{\times}$-action allows us to pass from formal to filtered quantizations.

Thus, we arrive at the problem of quantizing nilpotent orbits and covers. There is an extensive literature on the classification of deformation quantizations of Poisson manifolds, culminating in the celebrated work of Kontsevich (\cite{Kontsevich}). In the algebraic setting, there is an analog of Kontsevich's result due to Yekutieli (\cite{Yekutieli}). However, Yekutieli's result requires very strong cohomology vanishing assumptions on the Poisson variety being quantized, which are not satisfied in our setting. For smooth symplectic varieties, the situation is somewhat improved: the deformation quantizations can be classified under relatively mild cohomology vanishing assumptions (\cite{BK}). Although nilpotent covers are symplectic, they do not satisfy these cohomology vanishing assumptions except in special cases (\cite{Losev5}). 

These difficulties can be resolved in the following manner. Let $\widetilde{\Orb}$ be a $G$-equivariant cover of a nilpotent co-adjoint $G$-orbit. Consider its algebra $\CC[\widetilde{\Orb}]$ of polynomial functions. This algebra is finitely generated. So we can consider the affine variety $\widetilde{X} = \Spec(\CC[\widetilde{\OO}])$. It is Poisson, but in almost all cases singular, so the approaches of \cite{Yekutieli} and \cite{BK} do not apply. However, there is an auxilary Poisson variety $\widetilde{Y}$ (a `$\QQ$-factorial terminalization' of $\widetilde{X}$) together with a Poisson morphism $\rho:\widetilde{Y}\rightarrow \widetilde{X}$ 
which is a partial resolution of singularities (in the sense that $\widetilde{Y}$ is normal and $\rho$ is proper and birational) such that the smooth locus $\widetilde{Y}^{reg}$ satisfies the cohomology vanishing conditions of Bezrukavnikov-Kaledin.

This observation leads to a classification of filtered quantizations of $\CC[\widetilde{\Orb}]$, see \cite{Losev4}. The main result is that the filtered quantizations of $\CC[\widetilde{\OO}]$ are classified, up to isomorphism, by points in a finite-dimensional complex vector space, modulo a linear action of a certain finite group. In particular, and this is of paramount importance for the present work, there is a `canonical quantization' of $\CC[\widetilde{\OO}]$ with quantization parameter $0$. Using such quantizations we will define and study unipotent Harish-Chandra bimodules.

\end{preface}

\begin{acknowledgment}
The authors would like to thank Jeffrey Adams, Dan Barbasch, Alexander Premet and David Vogan   for many helpful conversations. The work of I.L. and D.M. was partially supported by the NSF
under grant DMS-2001139. The work of L.M.B was partially supported by the NSF under grant DMS-2501977.
\end{acknowledgment}

\mainmatter

\chapter{Introduction}\label{sec:intro}


The concept of a unipotent representation has its origins in the representation theory of finite Chevalley groups. Let $G(\F_q)$ be the group of $\F_q$-rational points of a connected reductive algebraic group $G$. Consider the set $\operatorname{Irr}(G(\F_q))$ of complex irreducible representations of $G(\F_q)$. The classification of $\operatorname{Irr}(G(\F_q))$ is a very old problem in representation theory, of fundamental importance. In \cite{DeligneLusztig}, Deligne and Lusztig defined a finite set of representations (called \emph{unipotent representations})
$$\mathrm{Unip}(G(\F_q)) \subset \mathrm{Irr}(G(\F_q)).$$
These representations play a central role in the classification of $\operatorname{Irr}(G(\F_q))$, which was later completed by Lusztig in \cite{Lusztig1984}. The importance of these representations is two-fold

\begin{itemize}
    \item[(i)] They are classified by certain parameters related to complex nilpotent co-adjoint orbits (in particular, their classification is independent of $q$).
    \item[(ii)] The classification of 
    $\mathrm{Irr}(G(\F_q))$ reduces to the classification of $\mathrm{Unip}(G(\F_q))$.

\end{itemize}
For details, we refer the reader to \cite{Lusztig1984}. 

Now, replace $\F_q$ with a local field $k$ and let $G(k)$ be the group of $k$-rational points. Consider the set $\mathrm{Irr}_u(G(k))$ of irreducible \emph{unitary} representations of $G(k)$. The classification $\mathrm{Irr}_u(G(k))$ is a classical problem in representation theory. Unlike its $\F_q$-analog, it remains unsolved in general (for an introduction to the theory of unitary representations, we refer the reader to \cite{Knapp}). An intriguing idea, with its origins in \cite{BarbaschVogan1985}, is to find a finite set of representations (called \emph{unipotent}, by analogy)
$$\mathrm{Unip}(G(k)) \subset \mathrm{Irr}_u(G(k)).$$
which plays the same role with respect to $\mathrm{Irr}_u(G(k))$ that $\unip(G(\F_q))$ plays with respect to $\mathrm{Irr}(G(\F_q))$.
In particular, these representations should satisfy suitable analogs of properties (i) and (ii) above. 

The goal of {this monograph} is to give a definition of $\mathrm{Unip}(G(k))$ in the case when $k=\CC$. The definition we propose is geometric in nature, and well-aligned with the orbit method, which is a guiding philosophy in the subject (see the preface above for a brief account or \cite{Vogan1992} for a more detailed overview). Our approach will shed some light on the more complicated case of  $k=\RR$, see Section \ref{subsec:realgroups}.

In {this monograph}, we will take an algebraic point of view: our main objects of study are Harish-Chandra bimodules (rather than group representations, which are analytic in nature). It is well known that these two formalisms are equivalent. We recall the basic theory of Harish-Chandra bimodules in Section \ref{subsec:HCbimodsclassical}. 

\section{Existing constructions}\label{subsec:existingconstructions}
In \cite{BarbaschVogan1985}, Barbasch and Vogan, following ideas of 
Arthur \cite{Arthur1983}, defined a set of representations called \emph{special unipotent}\index{representation!special unipotent}. We will briefly review their construction and discuss its limitations. 

Let $G$ be a complex reductive algebraic group and let $\fg$ be its Lie algebra. Let $\fg^{\vee}$ be the Langlands dual of $\fg$. If we fix a Cartan subalgebra $\fh \subset \fg$, there is a Cartan subalgebra $\fh^{\vee} \subset \fg^{\vee}$, canonically identified with $\fh^*$. To each nilpotent $\Ad(\fg^{\vee})$-orbit $\mathbb{O}^{\vee} \subset (\fg^{\vee})^*$, we can associate a maximal ideal in the enveloping algebra $U(\fg)$ as follows. First, using an $\Ad(\fg^{\vee})$-invariant isomorphism $(\fg^{\vee})^* \simeq \fg^{\vee}$, we can identify $\OO^{\vee}$ with a nilpotent $\Ad(\fg^{\vee})$-orbit in $\fg^{\vee}$ (still denoted by $\mathbb{O}^{\vee}$). Choose an element $e^{\vee} \in \mathbb{O}^{\vee}$ and an $\mathfrak{sl}(2)$-triple $(e^{\vee},f^{\vee},h^{\vee})$. Acting by $\Ad(\fg^{\vee})$ if necessary, we can arrange so that $h^{\vee} \in \fh^{\vee}$. This element is well-defined modulo the usual action of $W$. Every $W$-orbit in $\fh^*$ determines an infinitesimal character for $U(\fg)$ by means of the Harish-Chandra isomorphism $\mathfrak{Z}(\fg) \simeq \CC[\fh^*]^W$. Consider the infinitesimal character corresponding to the element $\frac{1}{2}h^{\vee} \in \fh^{\vee} \simeq \fh^*$, and let $I_{\mathrm{max}}(\frac{1}{2}h^{\vee}) \subset U(\fg)$ be the (unique) maximal ideal with this given infinitesimal character.

\begin{definition}[Def 1.17, \cite{BarbaschVogan1985}]\label{def:spec_unipotent}
The \emph{special unipotent ideal}\index{ideal!special unipotent} attached to $\mathbb{O}^{\vee} \subset \fg^{\vee}$ is the maximal ideal $I_{\mathrm{max}}(\frac{1}{2}h^{\vee}) \subset U(\fg)$. A \emph{special unipotent bimodule} attached to $\mathbb{O}^{\vee}$ is an irreducible Harish-Chandra bimodule which is annihilated (on both sides) by $I_{\mathrm{max}}(\frac{1}{2}h^{\vee})$. Denote the set of such bimodules by $\mathrm{Unip}^s_{\mathbb{O}^{\vee}}(G)$. 
\end{definition}

There is a distinguished class of nilpotent orbits, first defined in \cite{Lusztig1979}, called \emph{special} nilpotent orbits. Below we will recall one of several equivalent definitions. If $\mathbb{O}^{\vee}$ is special, there is a well-known classification of $\mathrm{Unip}^s_{\mathbb{O}^{\vee}}(G)$. Namely, there is a finite group $\overline{A}(\mathbb{O}^{\vee})$ attached to $\mathbb{O}^{\vee}$ called the \emph{Lusztig canonical quotient} (see \cite[Chp 13]{Lusztig1984}) and a natural bijection
\begin{equation}\label{eq:classificationspecial}\mathrm{Unip}^s_{\mathbb{O}^{\vee}}(G) \simeq \{\text{irreducible representations of } \overline{A}(\mathbb{O}^{\vee})\}. \end{equation}
In the case when $\frac{1}{2}h^{\vee}$ is integral, this statement is contained in \cite[Thm III]{BarbaschVogan1985}. The non-integral case can be handled by similar methods, see \cite{Wong2018}. If $\mathbb{O}^{\vee}$ is \emph{not} special, there is (to our knowledge) no known classification of $\mathrm{Unip}^s_{\mathbb{O}^{\vee}}(G)$. As an easy application of the main results in {this monograph}, we obtain one (see Remark \ref{rmk:classificationspecial}). Our classification is of the same form as (\ref{eq:classificationspecial}), except that $\overline{A}(\mathbb{O}^{\vee})$ is replaced with the 
Galois group of a certain finite cover of a nilpotent co-adjoint $G$-orbit. 

Definition \ref{def:spec_unipotent} has a number of limitations, which were known to Barbasch and Vogan at the time when \cite{BarbaschVogan1985} was written. Recall that to every irreducible Harish-Chandra bimodule, there is an associated nilpotent orbit $\mathbb{O} \subset \fg^*$ (see Section \ref{subsec:assvar} for details). There is also an order-reversing map
\begin{equation}\label{eq:Dintro}\mathsf{D}: \{\text{nilpotent orbits in } (\fg^{\vee})^*\} \to \{\text{nilpotent orbits in } \fg^*\}\end{equation}
called Barbasch-Vogan-Lusztig-Spaltenstein (BVLS) duality (see Section \ref{subsec:BVduality}). In \cite{BarbaschVogan1985} it is shown that the nilpotent orbit associated to a special unipotent bimodule $\B \in \mathrm{Unip}^s_{\mathbb{O}^{\vee}}(G)$ is $
\mathsf{D}(\mathbb{O}^{\vee}) \subset \fg^*$. 

The main limitation of Definition \ref{def:spec_unipotent} is related to the fact that $\mathsf{D}$ is not (usually) surjective. Its image is the set of special nilpotent orbits (this explains the word `special' in `special unipotent'). If $\g=\mathfrak{sl}(n)$, then $\mathsf{D}$ is a bijection and hence every orbit is special, but in all other types, there are non-special orbits. There are many bimodules attached to such orbits which deserve to be called unipotent, for example the `metaplectic' bimodules for $\mathrm{Sp}(2n)$, see Example \ref{ex:unipotent}(iv). Thus, Definition \ref{def:spec_unipotent} is incomplete with respect to the problem posed in the previous section (as mentioned above, Barbasch and Vogan were well aware of this difficulty when they formulated their definition, see e.g. the final paragraphs of \cite[Section 1]{BarbaschVogan1985}. In fact, this issue is the main subject of the second half of \cite{Vogan1987}).

Since the publication of \cite{BarbaschVogan1985}, various attempts have been made to generalize the notion of `special unipotent'. Some examples include: 
\begin{itemize}
    \item Vogan's weakly unipotent ideals (\cite{Vogan1984}).
    \item Barbasch's unipotent representations (\cite{Barbasch1989}). 
    \item McGovern's q-unipotent ideals (\cite{McGovern1994}).
    \item Brylinski's unipotent Dixmier algebras (\cite{Brylinski2003}).
\end{itemize}
None of these approaches is completely satisfactory. Vogan's definition picks out an \emph{infinite} set of ideals, which is inconsistent with the orbit method. The approaches of Barbasch and McGovern are based on ad hoc definitions which make sense only in classical types. Brylinski's approach is appealingly geometric but only applies in special cases. Worse still, the approaches of Vogan, McGovern, and Brylinski allow for non-unitary bimodules (see Section \ref{subsec:examplesideals} for examples).

\section{Vogan's desiderata}\label{subsec:desiderata}

In \cite[Chps 6-12]{Vogan1987}, Vogan lists some expected properties of unipotent bimodules. The most concrete of these expectations are catalogued below:

\begin{enumerate}
    \item $\unip(G)$ is a finite subset of $\mathrm{Irr}_u(G)$, attached in a natural fashion to $G$. See \cite[Chp 6, Problem I.2]{Vogan1987}.
    \item The set $\unip(G)$ is partitioned into `packets' $\unip_{\widetilde{\mathbb{O}}}(G)$ parameterized by equivariant covers of nilpotent co-adjoint $G$-orbits, i.e.
    $$\unip(G) = \bigsqcup_{\widetilde{\mathbb{O}}} \unip_{\widetilde{\mathbb{O}}}(G).$$
    See \cite[`Def' 12.8]{Vogan1987}.
    \item For each such $\widetilde{\mathbb{O}}$, there is maximal, completely prime ideal $I(\widetilde{\mathbb{O}}) \subset U(\fg)$ such that 
    $$\mathrm{LAnn}(\B) = \mathrm{RAnn}(\B) = I(\widetilde{\mathbb{O}}), \qquad \forall \ \B \in \mathrm{Unip}_{\widetilde{\mathbb{O}}}(G).$$
    See \cite[Conj 9.19]{Vogan1987} and \cite[`Def' 9.20]{Vogan1987}.
    \item For every $\B \in \mathrm{Irr}_u(G)$, there is a Levi subgroup $L \subset G$ and a unipotent bimodule $\B' \in \unip(L)$ such that $\B$ is obtained from $\B'$ through parabolic induction and/or a complementary series construction (\cite[Chp 6, Problem I.1]{Vogan1987})
    \item For every $\B \in \unip_{\widetilde{\mathbb{O}}}(G)$ and $e \in \mathbb{O}$, there is a finite-dimensional representation $\chi$ of the component group $G_e/G_e^{\circ}$ such that
    \begin{equation}\label{eq:Viso}\B \simeq_G \mathrm{AlgInd}^G_{G_e}\chi\end{equation}
    where $\mathrm{AlgInd}^G_{G_e}$ denotes induction of algebraic group representations (see \cite[Sec 3.3]{Jantzen}) and $\B$ is regarded as a $G$-representation via the adjoint action of $\fg$.  See \cite[Requirement 11.23]{Vogan1987}. See also \cite[Conj 12.1]{Vogan1991} for a more precise criterion. 
    \item $\unip(G)$ includes (as a proper subset) all special unipotent bimodules (cf. Definition \ref{def:spec_unipotent}). See the remarks following \cite[`Def' 12.8]{Vogan1987} as well as \cite[Chp 8]{Vogan1987}.
    \item For each $\widetilde{\mathbb{O}}$, the set $\mathrm{Unip}_{\widetilde{\mathbb{O}}}(G)$ is in one-to-one correspondence with irreducible representations of a finite group, which is related to the geometry of $\widetilde{\mathbb{O}}$. See \cite[Chp 7]{Vogan1987} as well as \cite[Chp 10]{Vogan1987}. 
\end{enumerate}

\section{Definitions}

The main contribution of {this monograph} is a definition of $\unip(G)$. Our definition is case-free, and geometric in nature. The main technical input is the theory of quantizations of conical symplectic singularities.

Conical symplectic singularities are a remarkable class of affine Poisson varieties with contracting $\CC^{\times}$-actions. Some basic examples include Kleinian singularities (i.e. varieties of the form $\CC^2/\Gamma$, for a finite subgroup $\Gamma \subset \mathrm{Sp}(2)$) and the nilpotent cones of complex reductive Lie algebras. If $X$ is a conical symplectic singularity, it makes sense to talk about its filtered quantizations (for a general discussion of graded Poisson varieties and their filtered quantizations, see Section \ref{subsec:quant}). The classification of filtered quantizations of conical symplectic singularities was completed by the first-named author in \cite{Losev4}. Roughly speaking, he shows that the set of filtered quantizations of $X$ is in bijection with points in a complex vector space $\fP^X$, which is recovered from the geometry of $X$. A consequence of this result, which is of primary importance for our purposes, is that every conical symplectic singularity has a distinguished quantization, i.e. the one corresponding to $0 \in \fP^X$. We call this quantization the \emph{canonical quantization} and denote it by $\cA_0^X$. We conjecture that $\cA_0^X$ is a simple algebra, for all $X$, and gather a great deal of evidence in Section \ref{subsec:conjecturesimplicity}.

Now let $\mathbb{O} \subset \fg^*$ be a nilpotent $G$-orbit and let $\widetilde{\mathbb{O}} \to \mathbb{O}$ be a finite connected $G$-equivariant cover. Note that $\widetilde{\mathbb{O}}$ is a smooth symplectic variety, and its ring of regular functions $\CC[\widetilde{\mathbb{O}}]$ is a finitely-generated algebra. Consider the affine variety $\widetilde{X} := \Spec(\CC[\widetilde{\mathbb{O}}])$. It is not difficult to show that $\widetilde{X}$ is a conical symplectic singularity (see \cite[Lem 2.5]{LosevHC})). 

If $\cA$ is \emph{any} filtered quantization of $\widetilde{X}$, the $G$-action on $\widetilde{X}$ lifts to a $G$-action on $\cA$ (by filtered algebra automorphisms). We show that the latter action is Hamiltonian, i.e. there is a quantum co-moment map $\Phi: U(\fg) \to \cA$ (lifting the classical co-moment map $\CC[\fg^*] \to \CC[\widetilde{\mathbb{O}}]$ induced from the natural map of varieties $\widetilde{\OO} \to \fg^*$). Note that $\Phi$ turns $\cA$ into a Dixmier algebra, in the sense of \cite[Def 2.1]{Vogan1990}. Our main definition is as follows.

\begin{definition}\label{defi:unipot_HC}
The \emph{unipotent ideal} attached to $\widetilde{\mathbb{O}}$ is the two-sided ideal
$$I_0(\widetilde{\mathbb{O}}) := \ker{(\Phi: U(\fg) \to \cA_0^{\widetilde{X}})}.$$
A \emph{unipotent bimodule} attached to $\widetilde{\mathbb{O}}$ is an irreducible Harish-Chandra bimodule with left and right annihilator equal to $I_0(\widetilde{\mathbb{O}})$. Denote the set of such bimodules by $\unip_{\widetilde{\mathbb{O}}}(G)$.
\end{definition}

It is easy to see that $I_0(\widetilde{\mathbb{O}})$ is a completely prime primitive ideal with associated variety $\overline{\mathbb{O}}$. We will see below that it is maximal (we supply a complete proof of this fact for linear classical groups. In the remaining cases, we provide only a sketch, leaving the details to a future paper \cite{MBM}). In fact, we will show that the maximality of $I_0(\widetilde{\OO})$ is equivalent to the simplicity of $\mathcal{A}^{\widetilde{X}}_0$, which, as stated above, is a general expectation. 

Much of {this monograph} is devoted to checking Vogan's desiderata for the definition above. Conditions (2), (5), (6), and (7) are verified for arbitrary $G$ {(Chapters \ref{sec:unipotent} and \ref{sec:duality})}. Condition (1) and (3) are verified for linear classical $G$ {(Chapter \ref{sec:unipbimod} and Appendix \ref{sec:maximality})}. Condition (4) is left more or less untouched (proving this condition is equivalent to the classification of $\mathrm{Irr}_u(G)$ for arbitrary $G$, a problem which we do not solve in {this monograph}). Further motivation for Definition \ref{defi:unipot_HC} will be provided in Section \ref{subsec:motivation}.

\section{Results}
We now summarize the main results of {this monograph}. 

\subsection{Infinitesimal characters and maximality}
In Chapter \ref{sec:centralchars}, we give a recipe for computing the infinitesimal characters of all unipotent ideals. There are two steps:

\begin{enumerate}
    \item Reduce the calculation to \emph{rigid} nilpotent orbits. This requires a fairly detailed analysis of the birational geometry of nilpotent covers. 
    \item Determine the infinitesimal characters attached to rigid nilpotent orbits. For this step we rely heavily on the work of McGovern (\cite{McGovern1994}) and Premet (\cite{Premet2013}).
\end{enumerate}
Following this recipe, one can, in principle, compute the infinitesimal character of an arbitrary unipotent ideal. For linear classical groups, we use this recipe to produce simple combinatorial formulas for unipotent infinitesimal characters, see Proposition \ref{prop:centralcharacterbirigidcover}. In (\cite{MBM}), we use a similar procedure to compute the unipotent infinitesimal characters for spin and exceptional groups. 

If $I \subset U(\fg)$ is a primitive ideal with known infinitesimal character and associated variety, there is a combinatorial condition for deciding whether $I$ is a maximal ideal. We check this condition in many cases, and prove the following theorem.

\begin{theorem}[Propositions \ref{prop:maximalityA}, and \ref{prop:maximalitytypeC}, Corollary \ref{cor:maximalitytypeD} below]\label{thm:maximality}
Suppose that $G$ is a linear classical group, i.e. $G=\operatorname{SL}(n)$, $\operatorname{SO}(n)$,
or $\operatorname{Sp}(2n)$. Let $\widetilde{\Orb}$ be a $G$-equivariant cover of a nilpotent orbit in $\g^*$. Then $I_0(\widetilde{\Orb})$ is maximal. 
\end{theorem}

In \cite{MBM}, we extend this result to arbitrary reductive groups.

\subsection{Classification}\label{subsection: classification}

In Chapter \ref{sec:unipotent}, we give geometric classifications of unipotent ideals and bimodules. First, we explain the classification of unipotent ideals. Let $\widetilde{\mathbb{O}}$ and $\widecheck{\mathbb{O}}$ be $G$-equivariant covers of a common nilpotent orbit $\mathbb{O} \subset \fg^*$. Consider the affine varieties $\widetilde{X} = \Spec(\CC[\widetilde{\mathbb{O}}])$, $\widecheck{X} = \Spec(\CC[\widecheck{\mathbb{O}}])$.
Both $\widetilde{X}$ and $\widecheck{X}$ contain finitely many $G$-orbits and $\widetilde{\OO}$ (resp. $\widecheck{\OO}$) is the unique open orbit in $\widetilde{X}$ (resp. $\widecheck{X}$). Furher, every $G$-equivariant morphism $\widetilde{\mathbb{O}} \to \widecheck{\mathbb{O}}$ extends to a unique finite $G$-equivariant morphism $\widetilde{X} \to \widecheck{X}$. The latter morphism is \'{e}tale over the open orbit $\widecheck{\Orb} \subset \widecheck{X}$.

\begin{definition}\label{defi:almost_etale}
A finite $G$-equivariant morphism $\widetilde{X} \rightarrow \widecheck{X}$ is \emph{almost \'{e}tale} if it is \'{e}tale over the open subset in $\widecheck{X}$ obtained by removing all $G$-orbits of codimension $\geqslant 4$. 
\end{definition}

Write $\widetilde{\mathbb{O}} \succeq \widecheck{\mathbb{O}}$ if there is a $G$-equivariant morphism $\widetilde{\mathbb{O}} \to \widecheck{\mathbb{O}}$ such that the induced morphism $\widetilde{X} \rightarrow \widecheck{X}$ is almost \'{e}tale. This defines a partial order on the set of $G$-equivariant covers of $\mathbb{O}$. Taking the symmetric closure, we get an equivalence relation on the same set. It is not difficult to see that each equivalence class contains a unique maximal element, see Lemma \ref{lem:maximalelement} below. 

\begin{theorem}[Theorem \ref{thm:classificationideals} below]\label{thm:classificationidealsintro}
The ideal $I_0(\widetilde{\OO})$ depends only on the equivalence class of $\widetilde{\OO}$ with respect to the equivalence relation above.
The passage from $\widetilde{\mathbb{O}}$ to $I_0(\widetilde{\mathbb{O}})$ defines a bijection between the following two sets
\begin{itemize}
    \item Equivalence classes of $G$-equivariant covers $\widetilde{\mathbb{O}}$ of $\mathbb{O}$.
    \item Unipotent ideals with associated variety $\overline{\Orb}$.
\end{itemize}
\end{theorem}

Next, we explain the classification of $\unip_{\widetilde{\mathbb{O}}}(G)$. By Theorem \ref{thm:classificationidealsintro}, the set $\unip_{\widetilde{\OO}}(G)$ depends only on the equivalence class $[\widetilde{\OO}]$. Thus, we can assume for this discussion that $\widetilde{\OO}$ is maximal in $[\widetilde{\OO}]$. Let $\Pi$ denote the automorphism group of $\widetilde{\mathbb{O}} \to \mathbb{O}$. This is a finite group, and can be described explicitly as follows. Choose an element
$e\in \Orb$ and $x\in \widetilde{\Orb}$ lying over $e$. Write $G_e$ and $G_x$ for the stabilizers in $G$ of $e$ and $x$, repsectively. Then $G_x$ is a subgroup of finite-index in $G_e$ and $\Pi \simeq N_{G_e}(G_x)/G_x$. 

\begin{theorem}[Theorem \ref{thm:classificationbimods} below]\label{thm:classificationbimodsintro}
There is a natural bijection
$$\{\text{irreducible representations of } \Pi\}  \xrightarrow{\sim}  \operatorname{Unip}_{\widetilde{\Orb}}(G).$$ 
\end{theorem}

\subsection{Unitarity}

The most basic requirement of unipotent bimodules is that they are unitary. We have some partial results in this direction.

\begin{theorem}[Theorem \ref{thm:unipotentunitary} below]\label{thm:unitary}
Suppose that $G$ is a linear classical group, and let $\widetilde{\mathbb{O}}$ be a $G$-equivariant cover of a nilpotent orbit in $\fg^*$. Then every bimodule in $\operatorname{Unip}_{\widetilde{\Orb}}(G)$ is unitary.  
\end{theorem}

The proof is structured as follows. For classical $\fg$, there is a well-know classification of irreducible unitary Harish-Chandra bimodules, due to Barbasch (\cite{Barbasch1989}). The main result of his paper, namely Theorem 12.8, is difficult to apply in general. An intermediate result in the same paper establishes the unitarity of a large collection of bimodules attached to \emph{rigid} nilpotent orbits (see \cite[Prop 10.6]{Barbasch1989}). As an easy consequence of this result, we show that all unipotent bimodules attached to rigid orbits are unitary. To complete the proof of Theorem \ref{thm:unitary}, we show that all of the remaining unipotent bimodules (i.e. those attached to covers of induced nilpotent orbits) are obtained from rigid ones by certain unitarity-preserving procedures, namely parabolic induction and complementary series.

We note that a uniform proof of the unitarity of all unipotent bimodules was recently given in \cite{DavisMasonBrown}.

\subsection{$G$-types}

Assume that $\widetilde{\Orb}$ 
is maximal in its equivalence class. Recall that the bimodules in
$\operatorname{Unip}_{\widetilde{\Orb}}(G)$ are classified by irreducible representation of $\Pi$ (i.e. the automorphism group of the cover $\widetilde{\mathbb{O}} \to \mathbb{O}$). 

\begin{theorem}[Theorem \ref{thm:classificationbimods} below]\label{thm:Vogan_conjecture}
Suppose that $I_0(\widetilde{\Orb})$ is maximal, and let $\B \in\operatorname{Unip}_{\widetilde{\Orb}}(G)$. Let $V$ denote the irreducible representation of $\Pi$ corresponding to $\B$. Then there is a good filtration on 
$\B$ and an isomorphism of graded $G$-equivariant $S(\fg)$-modules
$$\gr (\B)\simeq (\C[\widetilde{\Orb}]\otimes V)^\Pi.$$
\end{theorem}
The isomorphism of $G$-representations (\ref{eq:Viso}) predicted by Vogan is an easy consequence of this result (see Proposition \ref{prop:Gtypes} below).

\subsection{`Special unipotent' implies
`unipotent'} In Chapter \ref{sec:duality} we prove the following result.

\begin{theorem}[Corollary \ref{cor:specialimpliesunipotent} below]\label{thm:spec_unip}
Every special unipotent bimodule is unipotent. 
\end{theorem}

In fact, we prove a more precise result. To each nilpotent orbit $\mathbb{O}^{\vee} \subset (\fg^{\vee})^*$, we associate an equivalence class of covers  $\widetilde{\mathsf{D}}(\mathbb{O}^{\vee})$. This defines an injection
$$\widetilde{\mathsf{D}}: \{\text{nilpotent orbits in } (\fg^{\vee})^*\} \hookrightarrow \{\text{equivalence classes of covers of nilpotent orbits in } \fg^*\},$$
{refining BVLS duality (\ref{eq:Dintro})}. We show that for each $\mathbb{O}^{\vee} \subset\fg^{\vee}$, the unipotent ideal $I_0(\widetilde{\mathsf{D}}(\mathbb{O}^{\vee}))$ coincides with the special unipotent ideal $I_{\mathrm{max}}(\frac{1}{2}h^{\vee})$. Thus, $\mathrm{Unip}_{\mathbb{O}^{\vee}}^s(G) = \mathrm{Unip}_{\widetilde{\mathsf{D}}(\mathbb{O}^{\vee})}(G)$. 

Our construction of $\widetilde{\mathsf{D}}$ should be viewed as a special case of a more general (but still largely conjectural) duality known as \emph{symplectic duality}, see \cite[Sec 10]{BPWII}. In Section \ref{subsec:motivationsymplectic}, we will explain how our construction fits into this picture. 

\begin{rmk}
For $G = \mathrm{Sp}(2n)$, there is a notion, due to Barbasch, Ma, Sun, and Zhu, of a `metaplectic special' unipotent bimodule, see \cite{barbasch2020metaplectic} (these are irreducible bimodules with prescribed annihilators arising via a certain `metaplectic' modification of BVLS duality). It is easy to see that these bimodules are unipotent as well, compare to
Chapter \ref{sec:duality}. Since we will not use this result, we omit the proof.
\end{rmk}

\subsection{Dixmier conjecture}

In Section \ref{subsec:Dixmier}, we prove an old conjecture of Vogan (see \cite[Conj 2.3]{Vogan1990}). To each $G$-equivariant cover $\widetilde{\OO} \to \OO$ of a co-adjoint $G$-orbit $\OO \subset \fg^*$ (not necessarily nilpotent) we attach a canonically defined Dixmier algebra $\mathrm{Dix}(\widetilde{\OO})$ such that
$$\CC[\widetilde{\OO}] \simeq_G \mathrm{Dix}(\widetilde{\OO})$$
as representations of $G$. We prove the following result.

\begin{theorem}[Theorem \ref{thm:Dixmier} below]
The assigment $\widetilde{\OO} \mapsto \mathrm{Dix}(\widetilde{\OO})$ defines an injective correspondence between $G$-equivariant covers of co-adjoint orbits and Dixmier algebras for $G$.
\end{theorem}
Moreover, if $\widetilde{\OO}$ is a \emph{nilpotent} cover, then $\mathrm{Dix}(\widetilde{\OO})$ corresponds to the canonical quantization of $\CC[\widetilde{\OO}]$. This is one of the principal justifications for Definition \ref{defi:unipot_HC}.

\subsection{Discussion of approach}
Here we provide a brief description of our approach in {this monograph}, as it differs substantially from previous work in the subject. 
Whereas existing work in the area of unipotent representations has been largely Lie-theoretic, our approach in {this monograph} is mostly geometric. A recurring theme is that many interesting properties of unipotent bimodules which are very difficult to establish using Lie-theoretic methods, become much more elementary and transparent with a geometric point of view. The most striking examples are Theorems \ref{thm:classificationbimodsintro} and \ref{thm:Vogan_conjecture}.

Consider the problem of classifying unipotent bimodules. As discussed in Section \ref{subsec:existingconstructions}, if $\OO^{\vee}$ is special, the special unipotent bimodules attached to $\OO^{\vee}$ are classified by irreducible representations of the Lusztig quotient $\overline{A}(\OO^\vee)$
(this is essentially \cite[Thm III]{BarbaschVogan1985}).  
The proof appearing in \cite{BarbaschVogan1985} is based on a complicated (and in places, case-by-case) argument involving translation functors, Kazhdan-Lusztig cells, and the Springer correspondence. Even the definition of $\overline{A}(\OO^{\vee})$ is quite involved. Our approach, on the other hand, is geometric in nature. The group which replaces $\overline{A}(\OO^{\vee})$ is easy to define, our result is much more general, and our proof is more elementary. {In the important special case when $\OO^\vee$ is even, the group we define is known to coincide with   $\overline{A}(\OO^{\vee})$; this follows from \cite[Proposition 7.4]{LosevHC}}.

Another (related) example is Vogan's conjecture (Theorem \ref{thm:Vogan_conjecture}). Our proof of this result is easy and conceptual, whereas Lie-theoretic proofs, which exist in several cases, are difficult and ad hoc.

When it comes to results which are more Lie-theoretic in nature, our arguments become more technical. A good example is our (partial) computation of unipotent infinitesimal characters (Chapter \ref{sec:centralchars}). In the Barbasch-Vogan picture, special unipotent ideals are \emph{defined} by their infinitesimal characters, which are given by the elegant formula $\frac{1}{2}h^{\vee}$. In our generalization, there is no such simple formula for unipotent infinitesimal characters, and lots of work is required to compute them. Similarly, in the Barbasch-Vogan picture, special unipotent ideals are \emph{defined} to be maximal, whereas in our generalization, lots of work is required to prove this. We note that the simplicity of the canonical quantization (Conjecture \ref{conj:simplicity}) is equivalent to the maximality of the unipotent ideal (Theorem \ref{prop:simplemaximal}) and seems to be a general phenomenon.  However, proving the simplicity conjecture in general seems to be currently out of reach.   

An interesting question is whether it is possible to come up with a simple \emph{conceptual} description of the infinitesimal characters we obtain, perhaps in terms of the dual group $G^{\vee}$. It is possible that further developments in the area of symplectic duality will shed some light on this question.

\section{Outline of monograph}

Here is an outline of the paper. In Chapters \ref{sec:nilp}, \ref{sec:Lie}, and \ref{sec:symplectic}, we review some preliminary results which are needed in later Chapters. Much of the material in these Chapters already exists in the literature. 

In Chapter \ref{sec:nilp} we collect some elementary facts about nilpotent orbits in complex semisimple Lie algebras and their (finite, connected) covers. 

In Chapter \ref{sec:Lie} we recall some basic facts from Lie theory. This includes a discussion of primitive ideals in enveloping algebras, Harish-Chandra bimodules, Barbasch-Vogan-Lusztig-Spaltenstein (BVLS) duality, and $\cW$-algebras. 

In Chapter \ref{sec:symplectic} we review the key elements of the theory of conical symplectic singularities, their filtered quantizations, and Harish-Chandra bimodules over these filtered quantizations. This Chapter contains two main results. The first is a classification of filtered quantizations of conical symplectic singularities, obtained by the first-named author in \cite{Losev4}. 
To each conical symplectic singularity $X$ we attach a finite-dimensional parameter space $\fP$. Filtered quantizations of $X$ are parameterized by points in this parameter space, modulo a linear action of a certain finite group. The second main result is a classification of Harish-Chandra bimodules with full support over a filtered quantization of a conical symplectic singularity. In Chapter \ref{subsec:daggers} we recall a category equivalence, first obtained in \cite{LosevHC}, between the category of Harish-Chandra bimodules with full support and the category of representations of a certain finite group. In Sections \ref{subsec:coverings}, \ref{subsec:Gammalambda} we provide a (partial) description of this finite group (see, in particular,  Corollary \ref{cor:isotypic}, Theorem \ref{thm:LosevHC}). 

In Chapter \ref{sec:canonical}, we introduce the notion of a canonical quantization: the canonical quantization of a conical symplectic singularity is the filtered quantization corresponding to the zero parameter. If $\widetilde{X} \to X$ is a finite Galois covering of conical symplectic singularities, we get a quantization of $X$ by taking the algebra of invariants in the canonical quantization of $\widetilde{X}$. This quantization of $X$, which is typically not canonical, will play an important role in our description of unipotent ideals attached to nontrivial nilpotent covers. In Chapter \ref{sec:invariantssymplectic}, we compute its quantization parameter. Also in this Chapter we compute the finite group which controls the category of Harish-Chandra bimodules with full support over a canonical quantization, see Proposition \ref{prop:Gamma0}. In fact, we conjecture that canonical quantizations are simple algebras (and hence that \emph{all} Harish-Chandra bimodules are automatically of full support). This conjecture is discussed in Section \ref{subsec:conjecturesimplicity}.

In Chapter \ref{sec:unipotent} we introduce our main objects of study: unipotent ideals and unipotent bimodules. These objects are closely related to canonical quantizations, so the results of Chapter \ref{sec:canonical} are applicable. Our primary goal in Chapters \ref{sec:unipotent}-\ref{sec:unipbimod} is to show that these distinguished ideals and bimodules have all (or at least, many) of the properties one would hope for---in other words, to convince the reader that our definitions are the right ones. As a first step towards this goal, we show in Chapter \ref{sec:unipotent} that our unipotent ideals and bimodules have geometric classifications. Unipotent ideals are classified by equivalence classes of covers with respect to a geometrically defined equivalence relation connected to codimension 2 singularities. Unipotent bimodules are classified by irreducible representations of certain finite groups. The latter classification bears an encouraging resemblance to the Barbasch-Vogan classification of special unipotent bimodules. 

One of the most intriguing questions which arises from the results in this monograph is whether our notion of `unipotent' can be meaningfully extended to all real reductive Lie groups. We discuss this question briefly in Section \ref{subsec:realgroups}. One can naively extend our definition to representations of a (non-complex) Lie group by requiring that the annihilator of the underlying $(\fg,K)$-module is a unipotent ideal. Some very simple examples show that this naive definition is flawed---many of the representations in question are non-unitary. Remarkably, this issue seems to vanish if $\widetilde{\OO}=\OO$ is birationally rigid.

In Chapters \ref{sec:nilpquant}-\ref{sec:unipbimod}, we establish some additional key properties of unipotent ideals and bimodules. Unlike the results in Chapter \ref{sec:unipotent}, these require a detailed understanding of the infinitesimal characters of unipotent ideals. This turns out to be a rather delicate matter. In Chapters \ref{sec:nilpquant} and \ref{sec:centralchars} we present an algorithm for computing the infinitesimal character $\gamma_0(\widetilde{\OO})$ of the unipotent ideal attached to an arbitrary nilpotent cover. The main formula is stated in Chapter \ref{sec:centralchars}: for each nilpotent cover $\widetilde{\OO}$, there is a Levi subgroup $L \subset G$ and a rigid $L$-orbit $\OO_L$ such that 
$$\gamma_0(\widetilde{\OO}) = \gamma_0(\OO_L) + \delta,$$
where $\delta$ is a small `shift' related to the quantization parameter constructed in Chapter \ref{sec:canonical}. Thus, the computation of $\gamma_0(\widetilde{\OO})$ is a two-step process: first, compute the infinitesimal character attached to the rigid nilpotent orbit $\OO_L$, and then compute the element $\delta$. For the first step, we appeal to known results of McGovern (for $\fg$ a classical Lie algebra) and Premet (for $\fg$ exceptional). There are six ``bad'' orbits in exceptional types which are not handled by Premet (we handle these orbit separately in Appendix \ref{SS_m_sing}---they are the most difficult cases). On the other hand, computing the element $\delta$ requires a fairly detailed analysis of the birational geometry of $\widetilde{\OO}$. The required machinery is developed in Chapter \ref{sec:nilpquant}. 

Using our computation of unipotent infinitesimal characters, we prove in Section \ref{subsec:maximality} that all unipotent ideals for linear classical groups are maximal. The argument relies on some rather tedious combinatorics, which is relegated to Appendix \ref{sec:maximality}. 

In Chapter \ref{sec:duality} we introduce our refined BVLS duality. This is a map which takes nilpotent orbits for the Langlands dual group $G^{\vee}$ to nilpotent covers for $G$. One consequence of this duality is that all special unipotent ideals (in the sense of \cite{BarbaschVogan1985}) are unipotent. 

In Chapter \ref{sec:unipbimod}, we turn our attention to the unitarity of unipotent bimodules and the related issue of their construction via parabolic induction. The main result in this chapter, Theorem \ref{thm:unipotentunitary}, is that all unipotent bimodules for linear classical groups are unitary. The proof has two components: first, we show that all unipotent bimodules attached to rigid orbits are unitary. Then, we show that these `rigid' unipotents generate all others (through several unitarity-preserving operations). The first step is an immediate consequence of an old result of Barbasch (\cite{Barbasch1989}) and our infinitesimal character computations. For the second step, we introduce a version of parabolic induction for Harish-Chandra bimodules over filtered quantizations. Under suitable conditions, this construction coincides with the usual induction of Harish-Chandra bimodules (this is proved in Appendix \ref{sec:coincidence}), and thus preserves unitarity.

\chapter{Nilpotent orbits and their covers}\label{sec:nilp}
Let $G$ be a complex connected reductive algebraic group and let $\fg$ be its Lie algebra. For some statements, we will assume that $\fg$ is \emph{classical}, i.e. a simple Lie algebra of type $A$, $B$, $C$, or $D$. Occasionally, we will assume that $G$ is \emph{linear classical}. For us this will mean that $G$ belongs to one of four infinite families: $\mathrm{SL}(n)$, $\mathrm{SO}(2n+1)$, $\mathrm{Sp}(2n)$, and $\mathrm{SO}(2n)$. Note that this condition is strictly stronger than the previous. For example, $\mathrm{Spin}(n)$ is classical, but not linear classical {for $n>6$}.

In this chapter, we will collect some basic facts about nilpotent $G$-orbits and their equivariant covers.

\section{Nilpotent orbits}\label{subsec:orbitsclassical}

If $\fg$ is classical, the nilpotent orbits $\mathbb{O} \subset \fg^*$ are parameterized by (decorated) integer partitions. \index{nilpotent orbit}

\begin{definition}
A partition $p$ is of \emph{type C} (resp type \emph{B/D}) if every odd part (resp. even part) occurs with even multiplicity.  
\end{definition}

The following result is well-known.

\begin{prop}[Section 5.1, \cite{CM}]\label{prop:orbitstopartitions}
Suppose $\fg$ is classical. The set of nilpotent orbits $\mathbb{O} \subset \fg^*$ is parameterized by (decorated) partitions as follows.
\begin{enumerate}[label=(\alph*)]
    \item If $\fg = \mathfrak{sl}(n)$, the set of nilpotent orbits is in one-to-one correspondence with partitions of $n$.
    \item If $\fg = \mathfrak{so}(2n+1)$, the set of nilpotent orbits is in one-to-one correspondence with partitions of $2n+1$ of type B/D.
    \item If $\fg = \mathfrak{sp}(2n)$, the set of nilpotent orbits is in one-to-one correspondence with partitions of $2n$ of type C.
    \item If $\fg = \mathfrak{so}(2n)$, the set of nilpotent orbits is in one-to-one correspondence with partitions of $2n$ of type B/D, except that each \emph{very even} partition (i.e. a partition containing only even parts) corresponds to two nilpotent orbits, labeled $\mathbb{O}^{\mathrm{I}}$ and $\mathbb{O}^{\mathrm{II}}$.
\end{enumerate}
\end{prop}

There is a partial order on nilpotent orbits defined by the relation
$$\mathbb{O}^1 \leq \mathbb{O}^2 \iff \overline{\mathbb{O}}^1 \subseteq \overline{\mathbb{O}}^2.$$
If $\fg$ is classical, this order can be described in terms of partitions. Let $p^1 = (p^1_1,p^1_2,...,p^1_r)$ and $p^2=(p^2_1,p^2_2,...,p^2_s)$ be the partitions corresponding to $\mathbb{O}^1$ and $\mathbb{O}^2$, respectively. If the partitions $p^1$ and $p^2$ coincide, then either $\OO^1=\OO^2$, or $\fg=\mathfrak{so}(2n)$ and $\OO^1$, $\OO^2$ are incomparable very even orbits. If on the other hand $p^1 \neq p^2$, then $\mathbb{O}^1 \leq \mathbb{O}^2$ if and only if
$$\sum_{i=1}^j p^1_i \leq \sum_{i=1}^j p^2_i \qquad \text{for every } j,$$
where we complete partitions with zero parts as necessary. See \cite[Thm 6.2.5]{CM}.

If $\fg$ is exceptional, the classification of nilpotent orbits is more complicated. In {this monograph}, we will use the \emph{Bala-Carter} classification. For details, see Section \ref{subsec:BCinclusion}. Tables can be found in \cite[Sec 8.4]{CM}.

\section{Nilpotent covers}\label{subsec:nilpcovers}

A G\emph{-equivariant nilpotent cover} \index{cover!nilpotent}is a triple consisting of a nilpotent orbit $\mathbb{O} \subset \fg^*$, a homogeneous space $\widetilde{\mathbb{O}}$ for $G$, and a finite $G$-equivariant map $\widetilde{\mathbb{O}} \to \mathbb{O}$. The latter is automatically an \'{e}tale cover of $\mathbb{O}$, since $\widetilde{\mathbb{O}}$ is homogeneous. A \emph{morphism of nilpotent covers} is a $G$-equivariant map $f: \widetilde{\mathbb{O}} \to \widecheck{\mathbb{O}}$ such that the following diagram commutes
\begin{center}
\begin{tikzcd}
\widetilde{\mathbb{O}} \ar[r,"f"] \ar[rd] & \widecheck{\mathbb{O}}\ar[d]\\
& \mathbb{O} 
\end{tikzcd}
\end{center}
Denote the set of morphisms $\widetilde{\mathbb{O}} \to \widecheck{\mathbb{O}}$ by $\Hom_{\mathbb{O}}(\widetilde{\mathbb{O}},\widecheck{\mathbb{O}})$. Since $\widetilde{\mathbb{O}}$ is homogeneous and $\widecheck{\mathbb{O}} \to \mathbb{O}$ is finite, $\Hom_{\mathbb{O}}(\widetilde{\mathbb{O}},\widecheck{\mathbb{O}})$ is a finite set. To describe it, choose $e \in \mathbb{O}$ and a preimage $x \in \widetilde{\mathbb{O}}$. Write $G_e, G_x \subset G$ for the stabilizers. Note that $G_x \subseteq G_e$ is a finite index subgroup (in particular, $G_x^{\circ}=G_e^{\circ}$). The set $\Hom_{\mathbb{O}}(\widetilde{\mathbb{O}},\widecheck{\mathbb{O}})$ can be described as follows
\begin{equation}\label{eq:descriptionofHom}
\Hom_{\mathbb{O}}(\widetilde{\mathbb{O}},\widecheck{\mathbb{O}}) \simeq \{x' \in \widecheck{\mathbb{O}} \mid x' \text{ lies over $e$ and }G_x \subseteq G_{x'}\}.
\end{equation}
As a special case of (\ref{eq:descriptionofHom}), we obtain the following isomorphism
$$\Aut_{\OO}(\widetilde{\OO}) \simeq N_{G_e}(G_x)/G_x.$$
The group above is called the \emph{Galois group} of $\widetilde{\mathbb{O}}$. We say that $\widetilde{\mathbb{O}}$ is \emph{Galois} if $G_x$ is a normal subgroup in $G_e$. In this case, $\OO \simeq \widetilde{\OO}/\Aut_{\OO}(\widetilde{\OO})$.

There is a partial order on the set of isomorphism classes of nilpotent covers, defined by
\begin{equation}\label{eq:partialordercovers}\widetilde{\mathbb{O}} \geq \widecheck{\mathbb{O}} \iff \Hom_{\mathbb{O}}(\widetilde{\mathbb{O}},\widecheck{\mathbb{O}}) \neq \emptyset.\end{equation}
Consider the \emph{universal cover} $\widehat{\mathbb{O}} = G/G_e^{\circ}$ of $\OO$. It follows from (\ref{eq:descriptionofHom}) that $\widehat{\mathbb{O}}$ is Galois and the (unique) maximal cover with respect to the order defined by (\ref{eq:partialordercovers}). 

The \emph{equivariant fundamental group} of $\widetilde{\mathbb{O}}$ is by definition
$$\pi_1^G(\widetilde{\OO}) := \Aut_{\widetilde{\mathbb{O}}}(\widehat{\mathbb{O}}) \simeq G_x/G_e^{\circ}.$$
If $\widetilde{\mathbb{O}} \geq \widecheck{\mathbb{O}}$, then there is an embedding
$$\pi_1^G(\widetilde{\mathbb{O}}) \simeq G_x/G_e^{\circ} \hookrightarrow G_{x'}/G_e^{\circ} \simeq \pi_1^G(\widecheck{\mathbb{O}}).$$
In particular, $\pi_1^G(\widetilde{\mathbb{O}})$ is a subgroup of $\pi_1^G(\mathbb{O})$ (well-defined up to conjugacy). The passage from $\widetilde{\mathbb{O}}$ to $\pi_1^G(\widetilde{\mathbb{O}})$ defines a Galois correspondence
$$\pi_1^G: \{\text{isom. classes of nilpotent covers}\} \xrightarrow{\sim} \{\text{conj. classes of subgroups } H \subseteq \pi_1^G(\mathbb{O})\}.$$
For any subgroup $H \subset \pi_1^G(\mathbb{O})$, the preimage under $G_e \twoheadrightarrow G_e/G_e^{\circ}$ is a finite index subgroup $\widetilde{H} \subset G_e$, which defines a nilpotent cover $\widecheck{\mathbb{O}} = G/\widetilde{H}$ with $\pi_1^G(\widecheck{\mathbb{O}}) \simeq H$. Furthermore, $\widetilde{\mathbb{O}}$ is Galois if and only if $\pi_1^G(\widetilde{\mathbb{O}})$ is a normal subgroup of $\pi_1^G(\mathbb{O})$. If it is, there is a group isomorphism
$$\Aut_{\mathbb{O}}(\widetilde{\mathbb{O}}) \simeq \pi_1^G(\mathbb{O})/\pi_1^G(\widetilde{\mathbb{O}}).$$

If $G$ is classical, the groups $\pi_1^G(\mathbb{O})$ can be described in terms of partitions (see \cite[Cor 6.1.6]{CM}). Suppose $\mathbb{O} \subset \fg^*$ is a nilpotent orbit corresponding to a partition $p = (p_1,p_2,...,p_k)$. Put

\begin{align*}
    a &= \text{number of distinct odd } p_i.\\
    b &= \text{number of distinct even } p_i.\\
    d &= \mathrm{gcd}(p_i).\\
\end{align*}
For example, if $p=(4,3^2,2^4,1)$, then $a=2$, $b=2$, and $d=1$. 

If $G$ is 
linear classical, then $\pi_1^G(\OO)$ is given by the following table:

\vspace{5mm}
\begin{center}
\begin{tabular}{|c|c|} \hline
$G$ & $\pi_1^G(\mathbb{O})$\\ \hline
$\mathrm{SL}(n)$ & $\ZZ_d$ \\ \hline
$\mathrm{SO}(2n+1)$ & $(\mathbb{Z}_2)^{a-1}$\\ \hline
$\mathrm{Sp}(2n)$ & $(\ZZ_2)^b$ \\ \hline
$\mathrm{SO}(2n)$ & $(\mathbb{Z}_2)^{\mathrm{max}(a-1,0)}$ 
\\ \hline
\end{tabular}
\end{center}
\vspace{5mm}

We will also need a description of $\pi_1^G(\mathbb{O})$ for $G=\mathrm{Spin}(m)$. A partition $p$ is \emph{rather odd} if every odd part occurs with multiplicity 1. If $G=\mathrm{Spin}(m)$ and $p$ is \emph{not} rather odd, then $\pi_1^G(\mathbb{O}) \simeq \pi_1^{\mathrm{SO}(m)}(\mathbb{O})$. If $p$ is rather odd, then $\pi_1^G(\mathbb{O})$ is a central extension of $\pi_1^{\mathrm{SO}(m)}(\mathbb{O})$
by $\ZZ_2$.

If $G$ is exceptional and simply connected, the groups $\pi_1^G(\mathbb{O})$ can be found in \cite[Sec 8.4]{CM}.

\section{Lusztig-Spaltenstein induction}\label{subsec:LSinduction}

Let $M \subset G$ be a Levi subgroup of $G$, and let $\mathbb{O}_M$ be a nilpotent $M$-orbit. Fix a parabolic subgroup $Q \subset G$ with a Levi decomposition $Q = MU$. The annihilator of $\fq$ in $\fg^*$ is a $Q$-stable subspace, denoted by $\fq^{\perp} \subset \fg^*$. Choosing a nondegenerate invariant symmetric form on $\fg$, we get a $Q$-invariant identification $\fq^\perp \simeq \fu$. Form the $G$-equivariant fiber bundle $G \times^Q (\overline{\mathbb{O}}_M \times \mathfrak{q}^{\perp})$ over the partial flag variety $G/Q$. There is a proper $G$-equivariant map
$$\mu: G \times^Q (\overline{\mathbb{O}}_M \times \mathfrak{q}^{\perp}) \to \mathfrak{g}^*, \qquad \mu(g,\xi) = \Ad^*(g)\xi.$$
The image of $\mu$ is a closed irreducible $G$-invariant subset of the nilpotent cone $\cN\subset \fg^*$. Hence it is the closure in $\cN$ of a uniquely defined nilpotent $G$-orbit $\mathrm{Ind}^G_M\mathbb{O}_M \subset \fg^*$. The correspondence
$$\mathrm{Ind}^G_M: \{\text{nilpotent } M\text{-orbits}\} \to \{\text{nilpotent } G\text{-orbits}\}$$
is called \emph{Lusztig-Spaltenstein induction}\index{induction!Lusztig-Spaltenstein} (or simply `induction' when there is no risk for  confusion). A nilpotent orbit is said to be \emph{rigid} \index{nilpotent orbit!rigid}if it cannot be induced from a proper Levi subgroup. Of course, every nilpotent orbit can be induced from a rigid nilpotent orbit (typically there are several). Some properties of induction are catalogued below.

\begin{prop}[\cite{LusztigSpaltenstein1979} or \cite{CM}, Sec 7]\label{prop:propertiesofInd}
Induction has the following properties
\begin{itemize}
    \item[(i)] $\mathrm{Ind}^G_M$ depends only on $M$ (and not on $Q$)
    \item[(ii)] If $L \subset M$ is a Levi subgroup of $M$, then
    $$\Ind^G_L = \Ind^G_M \circ \Ind^M_L.$$
    \item[(iii)] If $\mathbb{O}_M \subset \fm^*$ is a nilpotent $M$-orbit and $\mathbb{O} = \Ind^G_M \mathbb{O}_M$, then
    $$\mathrm{codim}(\mathbb{O}_M,\cN_M) = \codim(\mathbb{O},\cN),$$
    where $\cN_M$ stands for the nilpotent cone in $\fm^*$.
\end{itemize}
\end{prop}

For classical $\fg$, one can describe the effect of induction on partitions and the partitions corresponding the rigid nilpotent orbits. We will begin by recalling the classification of Levi subgroups in classical types. Let $\fh$ be a Cartan subalgebra of $\fg$ and let $\Delta \subset \fh^*$ be the set of roots. Choose standard coordinates on $\fh$ (see \cite[Sec 5.2]{CM}) and denote the coordinate functions by $\{e_i\}$.

If $G = \mathrm{SL}(n)$ and $a = (a_1,...,a_t)$ is a partition of $n$, there is a Levi subgroup
$$\mathrm{S}(\mathrm{GL}(a_1) \times ... \times \mathrm{GL}(a_t)) := \mathrm{SL}(n) \cap (\mathrm{GL}(a_1) \times ... \times \mathrm{GL}(a_t)) \subset G$$
corresponding to the roots
$$\{\pm(e_i -e_j)\}_{1 \leq i < j \leq a_1} \cup ... \cup \{\pm (e_i-e_j)\}_{n-a_t+1 \leq i < j \leq n} \subset \Delta$$
Every Levi subgroup in $G$ is conjugate to one of this form, and no two such are conjugate.

If $G=\mathrm{SO}(2n+1)$, $\mathrm{Sp}(2n)$, or $\mathrm{SO}(2n)$, $0 \leq m \leq n$, and $a$ is a partition of $n-m$, there is a Levi subgroup
\begin{equation}\label{eq:Levi1}\mathrm{GL}(a_1) \times ... \times \mathrm{GL}(a_t) \times G(m) \subset G\end{equation}
corresponding to the roots
$$\{\pm (e_i-e_j)\}_{1\leq i < j \leq a_1} \cup ... \cup \{\pm (e_i - e_j)\}_{n-m-a_t+1 \leq i <j \leq n-m} \cup \Delta(m) \subset \Delta$$
where $\Delta(m) \subset \Delta$ has the obvious meaning. If $G=\mathrm{SO}(2n+1)$ or $\mathrm{Sp}(2n)$, then every Levi subgroup in $G$ is conjugate to one of this form, and no two such are conjugate.

If $G=\mathrm{SO}(2n)$, and $a$ is a partition of $n$ with only even parts, there is a Levi subgroup
\begin{equation}\label{eq:Levi2}\mathrm{GL}(a_1) \times ...\times \mathrm{GL}(a_t)' \subset G\end{equation}
corresponding to the roots
$$\{\alpha \in \Delta \mid \alpha(\underbrace{1,...,1}_{a_1},\underbrace{2,...,2}_{a_2},...,\underbrace{t,...,t,-t}_{a_t}) = 0\} \subset \Delta.$$
The prime is included to distinguish this subgroup from the subgroup $\mathrm{GL}(a_1) \times... \times \mathrm{GL}(a_t) \subset G$ defined above (to which it is $O(2n)$-, but not $SO(2n)$-conjugate). Every Levi subgroup in $G$ is conjugate to one of the form (\ref{eq:Levi1}) or (\ref{eq:Levi2}), and no two such are conjugate.

We now proceed to describing the effect of induction on partitions. If $p$ is a partition, write $B(p)$ (resp. $C(p)$) for the unique largest partition of type B/D (resp. C) dominated by $p$, and write $p^t$ for the transpose of $p$. If $q$ is another partition (of any size), define the sum $p+q$ `row by row', e.g. $(5,3,1^2) + (4,3^3) = (9,6,4^2)$. The following proposition is standard. For maximal Levi subgroups, it is proved in \cite[Thm 7.3.3]{CM}. The general case follows from the transitivity of induction (see Proposition \ref{prop:propertiesofInd}(ii)).

\begin{prop}\label{prop:inductionclassical}
The following are true:
\begin{enumerate}
    \item[(i)] Suppose $\fg = \mathfrak{sl}(n)$. Let $\fm = \mathfrak{s}(\mathfrak{gl}(a_1) \times ... \times \mathfrak{gl}(a_t)) \subset \fg$ and let
    $$\mathbb{O}_M = \mathbb{O}_M^1 \times ... \times \mathbb{O}_M^t \subset \cN_{\mathrm{GL}(a_1)} \times ... \times \cN_{\mathrm{GL}(a_t)} = \cN_M.$$
    Write $p^j$ for the partition of $a_j$ corresponding to $\mathbb{O}_M^j$ and define 
    $$p' = \sum_{j=1}^t p^j.$$
    Then $p'$ is the partition corresponding to $\Ind^G_M \mathbb{O}_M$.
    \item[(ii)] Suppose $\fg = \mathfrak{so}(2n+1)$ or $\mathfrak{so}(2n)$. Let $\fm =  \mathfrak{gl}(a_1) \times ... \times \mathfrak{gl}(a_t) \times \fg(m)$ or possibly (if $\fg=\mathfrak{so}(2n)$) $\fm = \mathfrak{gl}(a_1) \times ... \times \mathfrak{gl}(a_t)'$ (with all $a_i$ even) and let
    $$\mathbb{O}_M = \mathbb{O}_M^1 \times ... \times \mathbb{O}_M^t \times \mathbb{O}_M^0  \subset \cN_{\mathrm{GL}(a_1)} \times ... \times \cN_{\mathrm{GL}(a_t)} \times  \cN_{\mathrm{SO}(2m+1)} = \cN_M$$
    Write $p^j$ for the partition corresponding to $\mathbb{O}_M^j$ and define $p'$ by the formula
    $$p' = B(p^0 + 2\sum_{j=1}^t p^j).$$
    Then $p'$ is the partition corresponding to $\Ind^G_M \mathbb{O}_M$. If $\fg=\mathfrak{so}(2n)$ and $p'$ is very even, then its decoration can be deduced from \cite[Cor. 7.3.4]{CM}.
    \item[(iii)] Suppose $\fg=\mathfrak{sp}(2n)$. Let $\fm=  \mathfrak{gl}(a_1) \times ... \times \mathfrak{gl}(a_t) \times \mathfrak{sp}(2m)$ and let
    $$\mathbb{O}_M = \mathbb{O}_M^1 \times ... \times \mathbb{O}_M^t \times \mathbb{O}_M^0  \subset  \cN_{\mathrm{GL}(a_1)} \times ... \times \cN_{\mathrm{GL}(a_t)} \times \cN_{\mathrm{Sp}(2m)}  = \cN_M.$$
    Write $p^j$ for the partition corresponding to $\mathbb{O}_M^j$ and define $p'$ by the formula
    $$p' = C(p^0 + 2\sum_{j=1}^t p^j).$$
    Then $p'$ is the partition corresponding to $\Ind^G_M \mathbb{O}_M$.
\end{enumerate}
\end{prop}

As an easy consequence we obtain

\begin{cor}[\cite{CM}, Cor 7.3.5]\label{cor:rigid}
Suppose $\fg$ is classical, and let $\mathbb{O} \subset \fg^*$ be a nilpotent orbit corresponding to a partition $p$. Then $\OO$ is rigid if and only if one of the following is true:
\begin{itemize}
    \item[(i)] $\fg= \mathfrak{sl}(n)$ and $\mathbb{O} = \{0\}$.
    \item[(ii)] $\fg = \mathfrak{so}(2n+1)$ or $\mathfrak{so}(2n)$, $p$ satisfies
    $$p_i \leq p_{i+1}+1, \qquad \forall i,$$
    and no odd part of $p$ occurs exactly twice.
    \item[(iii)] $\fg =\mathfrak{sp}(2n)$, $p$ satisfies
    $$p_i \leq p_{i+1}+1, \qquad \forall i,$$
    and no even part of $p$ occurs exactly twice.
\end{itemize}
\end{cor}

For exceptional $\fg$, a list of rigid orbits, and the effect of induction on them, can be found in \cite[Sec 4]{deGraafElashvili}.

\section{Birational induction}\label{subsec:birationalinduction}

Choose a Levi subgroup $M \subset G$, a nilpotent $M$-orbit $\OO_M$, and a (finite, connected) $M$-equivariant cover $\widetilde{\OO}_M$ of $\OO_M$. Let $\OO = \Ind^G_M \OO_M$. Consider the affine variety $\widetilde{X}_M := \Spec(\CC[\widetilde{\OO}_M])$. There is an $M$-action on $\widetilde{X}_M$ (induced from the $M$-action on $\widetilde{\mathbb{O}}_M$) and a finite surjective $M$-equivariant map $\mu_M: \widetilde{X}_M \to \overline{\mathbb{O}}_M$. Choose a parabolic subgroup $Q \subset G$ with Levi decomposition $Q=MU$. Let $Q$ act on $\widetilde{X}_M\times \fq^\perp$ as follows: $M$ acts diagonally, and $U$ acts by
$$u \cdot (x,y)=(x,u\mu_M(x)-\mu_M(x)+uy),\quad u\in U, \quad x\in \widetilde{X}_M, \quad y\in \fq^\perp,$$
By our construction of the $Q$-action, the map $\mu_M \times \mathrm{id}: \widetilde{X}_M \times \fq^{\perp} \to \overline{\OO}_M \times \fq^{\perp}$ is $Q$-equivariant. It gives rise to a $G$-equivariant map $G \times^Q (\widetilde{X}_M \times \fq^{\perp}) \to G \times^Q (\overline{\OO}_M \times \fq^{\perp})$ of fiber bundles over $G/Q$. Consider the composition
$$\widetilde{\mu}: G \times^Q (\widetilde{X}_M \times \mathfrak{q}^{\perp}) \to G \times^Q (\overline{\mathbb{O}}_M \times \mathfrak{q}^{\perp}) \overset{\mu}{\to} \overline{\mathbb{O}}.$$
Note that $\widetilde{\mu}^{-1}(\OO) \to \OO$ is a (finite, connected) $G$-equivariant cover. The correspondence
$$\mathrm{Bind}^G_M: \{M\text{-eqvt nilpotent covers}\} \to \{G\text{-eqvt nilpotent covers}\}, \qquad \mathrm{Bind}^G_M(\widetilde{\mathbb{O}}_M) := \widetilde{\mu}^{-1}(\mathbb{O}).$$
is called \emph{birational induction}\index{induction!birational}. By construction, the following diagram commutes
\begin{center}
\begin{tikzcd}
\{M\text{-eqvt nilpotent covers}\} \ar[d] \ar[r,"\mathrm{Bind}^G_M"] & \{G\text{-eqvt nilpotent covers}\} \ar[d]\\
\{\text{nilpotent } M\text{-orbits}\} \ar[r, "\mathrm{Ind}^G_M"] &  \{\text{nilpotent } G\text{-orbits}\}
\end{tikzcd}
\end{center}
A $G$-equivariant nilpotent cover $\widetilde{\mathbb{O}}$ is said to be \emph{birationally rigid}\index{nilpotent cover!birationally rigid}\index{nilpotent orbit!birationally rigid} if it cannot be be birationally induced from a proper Levi subgroup. 

Note that $G$ acts `by conjugation' on pairs of the form $(M,\widetilde{\OO}_M)$. Indeed, if $x \in \widetilde{\OO}_M$, then $\widetilde{\OO}_M \simeq M/M_x$ as $M$-equivariant nilpotent covers. The $G$-action is defined by
$$g \cdot (M, \widetilde{\OO}_M) = (\Ad(g)M, \Ad(g)M/\Ad(g)M_x).$$
Note that $\Ad(g)M/\Ad(g)M_x$ is a cover of $\Ad(g)\OO_M$ and is independent of $x$ (up to isomorphism of $\Ad(g)M$-equivariant covers).

The main properties of birational induction are catalogued below.

\begin{prop}\label{prop:birationalinduction}
Birational induction has the following properties

\begin{enumerate}
    \item[(i)] $\mathrm{Bind}^G_M$ depends only on $M$ (and not on $Q$)
    \item[(ii)] If $L \subset M$ is a Levi subgroup of $M$, then
    $$\mathrm{Bind}^G_L = \mathrm{Bind}^G_M \circ \mathrm{Bind}^M_L.$$
    \item[(iii)] If $\widetilde{\mathbb{O}}$ is a $G$-equivariant nilpotent cover, there is a Levi subgroup $M \subset G$ and a birationally rigid $M$-equivariant nilpotent cover $\widetilde{\OO}_M$ such that 
    $$\widetilde{\mathbb{O}} = \mathrm{Bind}^G_M \widetilde{\mathbb{O}}_M.$$
    The pair $(M,\widetilde{\OO}_M)$ is unique up to conjugation by $G$.
    \item[(iv)] For any nilpotent cover $\widetilde{\OO} \to \OO$, write $\deg(\widetilde{\OO})$ for the degree of the covering map. Then
    $$\deg(\widetilde{\mathbb{O}}_M) \text{ divides } \deg(\mathrm{Bind}^G_M(\widetilde{\mathbb{O}}_M)).$$
\end{enumerate}
\end{prop}

\begin{proof}
Properties (i) and (iii) were established for nilpotent orbits in 
\cite[Lem 4.1]{Losev4} and \cite[Thm 4.4]{Losev4}. The proofs can easily be generalized to nilpotent covers (see \cite[Cor 4.3]{Mitya2020}). Properties (ii) and (iv) are trivial.
\end{proof}

It is possible to classify birationally rigid nilpotent orbits and covers. We will prove several results in this direction in Section \ref{subsec:codim2leaves}. 

We conclude this subsection by describing a large class of orbits which are birationally induced from $\{0\}$. Suppose $\OO \subset \fg^*$ is a nilpotent $G$-orbit. Using an $\Ad(\fg)$-invariant identification $\fg \simeq \fg^*$, we can regard $\OO$ as a nilpotent $G$-orbit in $\fg$. Choose an element $e \in \OO$ and an $\mathfrak{sl}(2)$-triple $(e,f,h)$. The operator $\ad(h)$ defines a $\ZZ$-grading on $\fg$
$$\fg = \bigoplus_{i \in \ZZ} \fg_i, \qquad \fg_i := \{\xi \in \fg \mid \ad(h)(\xi) = i\xi\}.$$
We say that $\OO$ is \emph{even}\index{nilpotent orbit!even} if $\fg_i=0$ for every odd integer $i$. In any case, we can define a parabolic subalgebra
\begin{equation}\label{eq:JMlevi}\mathfrak{p}_{\OO} = \mathfrak{l}_{\OO} \oplus \mathfrak{n}_{\OO}, \qquad \fl_{\OO} := \fg_0, \qquad \fn_{\OO} := \bigoplus_{i \geq 1} \fg_i.\end{equation}
We call $\fp_{\OO}$ (resp. $\fl_{\OO}$) the \emph{Jacobson-Morozov} parabolic (resp. Levi) \index{parabolic!Jacobson-Morozov}associated to $\OO$. Both $\fp_{\OO}$ and $\fl_{\OO}$ are well-defined up to conjugation by $G$. The following result is well-known. The proof is contained in \cite{Kostant1959}, see also \cite[Thm 3.3.1]{CM}.

\begin{prop}\label{prop:even}
Suppose $\OO$ is an even nilpotent $G$-orbit. Then
$$\OO = \mathrm{Bind}^G_{L_{\OO}} \{0\}.$$
\end{prop}

\section{Birational induction and equivariant fundamental groups}\label{subsec:inductionpi1}

Fix the notation of Section \ref{subsec:birationalinduction}. Choose an $M$-equivariant nilpotent cover $\widetilde{\mathbb{O}}_M$ and let $\widetilde{\mathbb{O}} = \mathrm{Bind}^G_M \widetilde{\mathbb{O}}_M$. Consider the fiber bundle $\widetilde{Z}^0 := G \times^Q (\widetilde{\mathbb{O}}_M \times \fq^{\perp})$ over the flag variety $G/Q$. 
Note, there is an inclusion $\widetilde{\mathbb{O}}=\widetilde{\mu}^{-1}(\OO)\subset \widetilde{Z}^0$. Let $i$ denote the inclusion map. It gives rise to a group homomorphism
\begin{equation}\label{eqn:istar}
i_*: \pi_1(\widetilde{\mathbb{O}}) \to \pi_1(\widetilde{Z}^0).\end{equation}
Since $i: \widetilde{\mathbb{O}} \hookrightarrow \widetilde{Z}^0$ is an embedding of smooth complex manifolds, the complement $\widetilde{Z}^0 - i(\widetilde{\mathbb{O}})$ is of \emph{real} codimension $\geq 2$. Hence, the homomorphism (\ref{eqn:istar}) is surjective. We will use it to define a surjective homomorphism $f: \pi^G_1(\widetilde{\mathbb{O}}) \to \pi^M_1(\widetilde{\mathbb{O}}_M)$.

First, note that $\widetilde{Z}^0$ is a vector bundle over the homogeneous space $G \times^Q \widetilde{\mathbb{O}}_M$ with fiber equal to $\fq^{\perp}$. Hence, there is a natural isomorphism
$$\pi_1(\widetilde{Z}^0) \simeq \pi_1(G \times^Q \widetilde{\mathbb{O}}_M).$$
If we fix a base point $x \in \widetilde{\mathbb{O}}$, we get a fibration $G \to \widetilde{\mathbb{O}}$ with fiber $G_x$. This fibration gives rise to an exact sequence of homotopy groups
$$\pi_1(G) \to \pi_1(\widetilde{\mathbb{O}}) \to \pi_0(G_x) \to 1.$$
The final (nontrivial) term is isomorphic to $\pi_1^G(\widetilde{\mathbb{O}})$.

Similarly, if we fix a base point $(1,y) \in G \times^Q \widetilde{\mathbb{O}}_M$, we get a fibration $G \to G \times^Q \widetilde{\mathbb{O}}_M$ with fiber $Q_y$ and hence an exact sequence
$$\pi_1(G) \to \pi_1(G \times^Q \widetilde{\mathbb{O}}_M) \to \pi_0(Q_y) \to 1.$$
Since $U\subset Q_y$, we have $Q_y= M_y \ltimes U$. It follows that $\pi_0(Q_y) \simeq \pi_0(M_y) \simeq \pi_1^M(\widetilde{\mathbb{O}}_M)$. Choose $x \in \widetilde{\OO}$ and $y \in \widetilde{\OO}_M$ such that $x$ maps to $(1,y)$ under $\widetilde{\OO} \subset \widetilde{Z}^0 \twoheadrightarrow G \times^Q \widetilde{\OO}_M$. There is a commutative diagram
\begin{center}
    \begin{tikzcd}
      \pi_1(G) \ar[r] \ar[d,equals] & \pi_1(\widetilde{\mathbb{O}}) \ar[r] \ar[d,twoheadrightarrow,"i_*"] & \pi_1^G(\widetilde{\mathbb{O}}) \ar[r] & 1\\
      \pi_1(G) \ar[r] & \pi_1(\widetilde{Z}^0) \ar[r] & \pi_1^M(\widetilde{\mathbb{O}}_M) \ar[r] & 1
    \end{tikzcd}
\end{center}
where the fundamental groups $\pi_1(G)$, $\pi_1(\widetilde{\OO})$, $\pi_1(\widetilde{Z}^0)$, and $\pi_1^M(\widetilde{\OO}_M)$ are defined with respect to the base points $1 \in G$, $x \in \widetilde{\OO}$, $x \in \widetilde{Z}^0$, and $y \in \widetilde{\OO}_M$, respectively. The following lemma is immediate. 

\begin{lemma}\label{lem:mappi1}
There is a uniquely defined homomorphism $f: \pi_1^G(\widetilde{\mathbb{O}}) \to \pi_1^M(\widetilde{\mathbb{O}}_M)$ which makes the above diagram commute. This homomorphism is surjective.
\end{lemma}

\section{Bala-Carter inclusion}\label{subsec:BCinclusion}

Let $M \subset G$ be a Levi subgroup of $G$. We defined in \cref{subsec:LSinduction} a map $\mathrm{Ind}^G_M$ from nilpotent $M$-orbits to nilpotent $G$-orbits. There is a second map with the same source and target called \emph{Bala-Carter inclusion}\index{saturation}\index{Bala-Carter inclusion}. It is defined by
$$\mathrm{Sat}^G_M: \{\text{nilpotent } M\text{-orbits}\} \to \{\text{nilpotent } G\text{-orbits}\}, \qquad \mathrm{Sat}^G_M(\mathbb{O}_M) = G \cdot \mathbb{O}_M$$
(`Sat' stands for `saturation'). An orbit $\mathbb{O}$ is said to be \emph{distinguished}\index{nilpotent orbit!distinguished} if it is not obtained by Bala-Carter inclusion from a proper Levi subgroup. 

\begin{prop}[\cite{BalaCarter1976}, see also \cite{CM}, Thm 8.1.1]\label{prop:BC1}
If $\mathbb{O}$ is a nilpotent $G$-orbit, there is a Levi subgroup $M \subset G$ and a distinguished nilpotent $M$-orbit $\mathbb{O}_M$ such that
$$\mathbb{O} = \mathrm{Sat}^G_M \mathbb{O}_M.$$
The pair $(M,\mathbb{O}_M)$ is unique up to conjugation by $G$.
\end{prop}

Suppose $Q \subset G$ is a parabolic subgroup with Levi decomposition $Q = MU$. We say that $Q$ is \emph{distinguished}\index{parabolic!distinguished} if $ \dim(M/Z(G)) = \dim(U/[U,U])$ (e.g. a Borel subgroup is distinguished). 

\begin{prop}[\cite{BalaCarter1976}, see also \cite{CM}, Thms 8.2.6, 8.2.8]\label{prop:BC2}
If $\mathbb{O}$ is a distinguished nilpotent $G$-orbit, there is a distinguished parabolic $P =LN \subset G$ such that $\mathbb{O} = \mathrm{Ind}^G_L \{0\}$. Furthermore, $P$ is unique up to conjugation by $G$. 
\end{prop}

\begin{rmk}\label{rmk:distinguishedbirational}
We note that every distinguished nilpotent orbit is even, see the discussion preceding \ref{prop:even}. A convenient reference for this fact is \cite[Thm 8.2.3]{CM}. Thus by Proposition \ref{prop:even}, every distinguished orbit $\OO$ is birationally induced from the $\{0\}$-orbit of its Jacobson-Morozov Levi $L_{\OO} \subset G$. Furthermore, the associated parabolic $P_{\OO}$ is automatically distinguished, see \cite[Thm 8.2.6]{CM}. It follows that the parabolic $P$ appearing in Proposition \ref{prop:BC2} is $G$-conjugate to $P_{\OO}$.
\end{rmk}
Combining Propositions \ref{prop:BC1} and \ref{prop:BC2}, we obtain the following well known result.

\begin{cor}[\cite{BalaCarter1976}, see also \cite{CM}, Thm 8.2.12]\label{cor:BC}
There is a bijection between nilpotent $G$-orbits and $G$-conjugacy classes of pairs $(M,P_M)$ consisting of a Levi subgroup $M \subset G$ and a distinguished parabolic $P_M \subset M$. This bijection is defined by
$$(M,P_M) \mapsto \mathrm{Sat}^G_M \mathrm{Ind}^M_{P_M} \{0\}.$$
\end{cor}
We are now prepared to explain the Bala-Carter notation for nilpotent orbits. A nilpotent orbit $\mathbb{O}$ corresponding to a pair $(M,P_M)$ (as in {Corollary} \ref{cor:BC}) is assigned a \emph{Bala-Carter label} $X(s_i)$, where:
\begin{itemize}
    \item $X$ is the Lie type of $[\fm,\fm]$.
    \item $i$ is the number of simple roots in a Levi subalgebra of $\fp_{\fm}$.
    \item $s$ is a letter (either $a$ or $b$) chosen arbitrarily to distinguish orbits with identical $X$ and $i$.
\end{itemize}
If $i=0$, the parenthetical $(s_0)$ is omitted. In a few cases, the orbit $\mathbb{O}$ is underdetermined by the notation $X(s_i)$ and additional symbols are needed to identify it uniquely. For standard notation in these cases, we refer the reader to \cite[Sec 8.4]{CM}.

For classical $\fg$, one can describe the effect of Bala-Carter inclusion on partitions and the partitions corresponding to distinguished nilpotent orbits. If $p$ and $q$ are partitions (of any size), define the union $p \cup q$ by adding multiplicities, e.g. $(5,3,1^2) \cup (4,3^3) = (5,4,3^4,1^2)$. The following proposition is standard.

\begin{prop}\label{prop:inclusionclassical}
The following are true:
\begin{enumerate}
    \item[(i)] Suppose $\fg = \mathfrak{sl}(n)$. Let  $\mathfrak{m}=\mathfrak{s}(\mathfrak{gl}(a_1) \times ... \times \mathfrak{gl}(a_t))$ and let
    $$\mathbb{O}_M = \mathbb{O}_M^1 \times ... \times \mathbb{O}_M^t \subset \cN_{\mathrm{GL}(a_1)} \times ... \times \cN_{\mathrm{GL}(a_t)} = \cN_L=M.$$
    Write $p^j$ for the partition of $a_j$ corresponding to $\mathbb{O}_M^j$ and define $p'$ by the formula
    $$p' = \bigcup_{j=1}^t p^j.$$
    Then $p'$ is the partition corresponding to $\mathrm{Sat}^G_M\mathbb{O}_M$.
    \item[(ii)] Suppose $\fg = \mathfrak{so}(2n+1)$, $\mathfrak{sp}(2n)$ or $\mathfrak{so}(2n)$. Let $\mathfrak{m}=  \mathfrak{gl}(a_1) \times ... \times \mathfrak{gl}(a_t) \times \mathfrak{g}(m)$ or possibly (if $\fg = \mathfrak{so}(2n)$) $\fm = \mathfrak{gl}(a_1) \times ... \times \mathfrak{gl}(a_t)'$ (with all $a_i$ even) and let
    $$\mathbb{O}_M =  \mathbb{O}_M^1 \times ... \times \mathbb{O}_M^t \times \mathbb{O}_M^0 \subset \cN_{\mathrm{GL}(a_1)} \times ... \times \cN_{\mathrm{GL}(a_t)} \times  \cN_{G(m)} = \cN_M.$$
    Write $p^j$ for the partition corresponding to $\mathbb{O}_M^j$. Define $p'$ by the formula
    $$p' =p^0 \cup \bigcup_{j=1}^t (p^j \cup p^j).$$
    Then $p'$ is the partition corresponding to $\mathrm{Sat}^G_M\mathbb{O}_M$. If $\fg=\mathfrak{so}(2n)$ and $p'$ is very even, then $p^0$ is very even. In this case, $p'$ has the same decoration as $p^0$. 
\end{enumerate}
\end{prop}

\begin{cor}[\cite{CM}, Thm 8.2.14]\label{corollary: distinguished}
If $\fg = \mathfrak{sl}(n)$, then the only distinguished nilpotent $G$-orbit is the principal one. If $\fg = \mathfrak{so}(2n+1)$, $\mathfrak{sp}(2n)$, or $\mathfrak{so}(2n)$, then $\mathbb{O}$ is distinguished if and only if the corresponding partition has no repeated parts.
\end{cor}

\chapter{Lie theory preliminaries}\label{sec:Lie}

In this chapter, we will review some Lie theoretic preliminaries which are needed in later chapters. Let $G$ be a complex connected reductive algebraic group. Choose a Borel subgroup $B \subset G$ and a maximal torus $H \subset B$. Denote the root system by $\Delta$ and the Weyl group by $W$. 

\section{Associated varieties and primitive ideals}\label{subsec:assvar}

Let $I \subset U(\fg)$ be a two-sided ideal. If we equip $U(\fg)$ with its usual filtration, there is a $G$-equivariant graded algebra isomorphism isomorphism $\gr U(\fg) \simeq S(\fg)$ by Poincar\'{e}-Birkhoff-Witt. Under this identification, $\gr(I)$ corresponds to a $G$-invariant ideal in $S(\fg)$, which defines a $G$-invariant subset $V(I) \subset \fg^* = \Spec(S(\fg))$, called the \emph{associated variety}\index{associated variety!of primitive ideal} of $I$.

Recall that $I$ is said to be \emph{primitive}\index{ideal!primitive} if it is the annihilator of a simple left $U(\fg)$-module. Denote the set of such ideals by $\mathrm{Prim}(U(\fg))$. If $I \in \mathrm{Prim}(U(\fg))$, then by a theorem of Joseph (\cite{Joseph1985}), $V(I)$ is the closure of a single nilpotent orbit. 

Let $\mathfrak{Z}(\fg)$ denote the center of $U(\fg)$. If $I$ is primitive, then by Quillen's lemma, the intersection $I \cap \mathfrak{Z}(\fg)$ is the kernel of an algebra homomorphism $\mathfrak{Z}(\fg) \to \CC$, called the \emph{infinitesimal character}\index{infinitesimal character} of $I$. Infinitesimal characters are identified with $W$-orbits in $\fh^*$ via the Harish-Chandra isomorphism $\mathfrak{Z}(\fg) \simeq \CC[\fh^*]^{W}$. For each $\gamma \in\fh^*$, define
$$\mathrm{Prim}_{\gamma}(U(\fg)) := \{I \in \mathrm{Prim}(U(\fg)) \mid W\gamma = \text{infinitesimal character of } I\}.$$
\begin{prop}\label{prop:propertiesofprim}
Let $\gamma \in \fh^*$. Then
\begin{itemize}
    \item[(i)] $\mathrm{Prim}_{\gamma}(U(\fg))$ is a finite, nonempty set. 
    \item[(ii)] $\mathrm{Prim}_{\gamma}(U(\fg))$ contains a unique maximal element $I_{\mathrm{max}}(\gamma)$ with respect to the inclusion ordering on ideals. 
    \item[(iii)] If $I,I'$ are prime (for example, primitive) ideals in $U(\fg)$, then
    $$I \subsetneq I' \implies V(I') \subsetneq V(I).$$
\end{itemize}
\end{prop}

\begin{proof}
(i) and (ii) follow from the results of \cite{Duflo1977}. (iii) follows from \cite[Cor 3.6]{BorhoKraft}.
\end{proof}

\section{BVLS duality}\label{subsec:BVduality}

Let $G^{\vee}$ be the Langlands dual of $G$ and let $\fg^{\vee}$ be the Lie algebra of $G^{\vee}$. By construction, $\fg^{\vee}$ contains a Cartan subalgebra $\fh^{\vee} \subset \fg^{\vee}$ which is naturally identified with $\fh^*$. To every nilpotent $G^{\vee}$-orbit $\mathbb{O}^{\vee} \subset (\fg^{\vee})^*$, one can attach a `dual' nilpotent $G$-orbit $\mathsf{D}(\mathbb{O}^{\vee}) \subset \fg^*$ as follows. First, replace $\mathbb{O}^{\vee}$ with its counterpart in $\fg^{\vee}$ using a $G^{\vee}$-invariant identification $\fg^{\vee} \simeq (\fg^{\vee})^*$. Choose an element $e^{\vee} \in \mathbb{O}^{\vee}$ and an $\mathfrak{sl}(2)$-triple $(e^{\vee},f^{\vee},h^{\vee})$. Replacing this triple with a $G^\vee$-conjugate triple if necessary, we can arrange so that $h^{\vee} \in \fh^{\vee}$. This element is well-defined modulo $W$. Consider the infinitesimal character for $U(\fg)$ corresponding to $\frac{1}{2}h^{\vee} \in \fh^{\vee} \simeq \fh^*$ via the Harish-Chandra isomorphism and let $I_{\mathrm{max}}(\frac{1}{2}h^{\vee}) \subset U(\fg)$ be the unique maximal ideal of infinitesimal character $\frac{1}{2}h^{\vee}$, cf. Proposition \ref{prop:propertiesofprim}(ii). Finally, let
$$\mathsf{D}(\mathbb{O}^{\vee}) := \text{unique open } G\text{-orbit in } V(I_{\mathrm{max}}(\frac{1}{2}h^{\vee})).$$
This defines a map 
$$\mathsf{D}: \{\text{nilpotent } G^{\vee}\text{-orbits}\} \to \{\text{nilpotent } G\text{-orbits}\}$$
called \emph{Barbasch-Vogan-Lusztig-Spaltenstein (BVLS) duality}\index{duality!BVLS}. The definition above is due to Barbasch-Vogan (see \cite{BarbaschVogan1985}). Variants of this map can be found in the earlier work of Lusztig (\cite{Lusztig1984}) and Spaltenstein (\cite{Spaltenstein}). 

A nilpotent $G$-orbit is said to be \emph{special}\index{nilpotent orbit!special} if it lies in the image of $\mathsf{D}$. The map $\mathsf{D}$ restricts to a bijection
$$\mathsf{D}: \{\text{special nilpotent } G^{\vee}\text{-orbits}\} \xrightarrow{\sim} \{\text{special nilpotent } G\text{-orbits}\}.$$
The inverse bijection coincides with the duality map for $G^{\vee}$, see \cite[Prop A2]{BarbaschVogan1985}.

BVLS duality intertwines Bala-Carter inclusion (see Section \ref{subsec:BCinclusion}) and Lusztig-Spaltenstein induction (see Section \ref{subsec:LSinduction}) in the following sense.

\begin{prop}[\cite{BarbaschVogan1985}, Prop A2]\label{prop:inclusioninduction}
Suppose $L \subset G$ is a Levi subgroup of $G$. Then $L^{\vee} \subset G^{\vee}$ is a Levi subgroup of $G^{\vee}$. If $\mathbb{O}^{\vee}_{L}$ is a nilpotent $L^{\vee}$-orbit, then
$$\mathrm{Ind}^G_L \mathsf{D}\left(\mathbb{O}^{\vee}_L\right) = \mathsf{D}\left(\mathrm{Sat}^{G^{\vee}}_{L^{\vee}} \mathbb{O}^{\vee}_L\right).$$
\end{prop}

\begin{rmk}
It is worth noting that Proposition \ref{prop:inclusioninduction} is false if we swap $\Ind$ and $\mathrm{Sat}$. That is, in general $\mathrm{Sat}^G_L \mathsf{D}(\mathbb{O}_L^{\vee}) \neq \mathsf{D}\left(\Ind^{G^{\vee}}_{L^{\vee}}\mathbb{O}^{\vee}_L\right)$.
\end{rmk}

If $\fg$ is classical, one can describe $\mathsf{D}$ in terms of partitions.

\begin{definition}\label{def:lande}
Let $p = (p_1,...,p_k)$ be a partition of $n$. Define partitions $l(p)$ and $e(p)$ of $n-1$ and $n+1$, respectively:
$$l(p) := (p_1,...,p_{k-1},p_k-1), \qquad e(p) := (p_1,...,p_k,1).$$
\end{definition}

\begin{prop}[\cite{McGovern1994}, Thm 5.2]\label{prop:BVduality}
If $\fg$ is classical, $\mathsf{D}$ can be described in terms of partitions as follows
\begin{itemize}
    \item[(i)] If $\fg = \mathfrak{sl}(n)$, then $\fg^{\vee} = \mathfrak{sl}(n)$ and $\mathsf{D}(p) = p^t$.
    \item[(ii)] If $\fg = \mathfrak{so}(2n+1)$, then $\fg^{\vee} =  \mathfrak{sp}(2n)$ and $\mathsf{D}(p) = B(e(p)^t)$.
    \item[(iii)] If $\fg = \mathfrak{sp}(2n)$, then $\fg^{\vee} = \mathfrak{so}(2n+1)$ and $\mathsf{D}(p) = C(l(p^t))$.
    \item[(iv)] If $\fg = \mathfrak{so}(2n)$, then $\fg^{\vee} = \mathfrak{so}(2n)$ and $\mathsf{D}(p) =B(p^t)$. If $n$ is divisible by 4 and $p$ is very even, then $\mathsf{D}(p)$ preserves decoration. If $n$ is not divisible by 4 and $p$ is very even, then $\mathsf{D}(p)$ reverses decoration.
\end{itemize}
\end{prop}

For exceptional $\fg$, a description of $\mathsf{D}$ can be recovered from the information given in \cite[Sec 13.4]{Carter1993}.

\section{Maximal ideals}\label{subsec:assvarmax}

Let $\gamma \in \fh^*$ and let $I \in \mathrm{Prim}_{\gamma}(U(\fg))$. In this section, we will describe a (well-known) criterion for deciding whether $I$ is a maximal ideal.\index{ideal!maximal}

Fix $G^{\vee}$ and $\fh^{\vee}$ as in Section \ref{subsec:BVduality}. Using the identification $\fh^{\vee} \simeq \fh^*$, we can regard $\gamma$ as an element of $\fh^{\vee}$. This element defines two subsystems of $\Delta(\fg^{\vee},\fh^{\vee}) = \Delta^{\vee}$
$$\Delta^{\vee}_{\gamma} := \{\alpha^{\vee} \in \Delta^{\vee}: \langle\gamma, \alpha^{\vee}\rangle \in \ZZ\}, \qquad \Delta^{\vee}_{\gamma,0} := \{\alpha^{\vee} \in \Delta^{\vee}: \langle\gamma, \alpha^{\vee}\rangle =0\} \subset \Delta^{\vee}_{\gamma}$$
called the subsystems of \emph{integral} and \emph{singular} co-roots. These subsystems define reductive subalgebras $\mathfrak{l}^{\vee}_{\gamma}$ and $\mathfrak{l}^{\vee}_{\gamma,0}$ of $\fg^{\vee}$
\begin{equation}\label{eqn:integralsingular}
\mathfrak{l}^{\vee}_{\gamma} = Z_{\fg^{\vee}}(\exp(2\pi i\gamma)), \qquad \mathfrak{l}^{\vee}_{\gamma,0} = Z_{\fg^{\vee}}(\gamma) \subset \mathfrak{l}^{\vee}_{\gamma},
\end{equation}
where $Z_{\fg^\vee}$ stands for the centralizer in $\fg^\vee$.

By definition, $\mathfrak{l}^{\vee}_{\gamma,0}$ is a Levi subalgebra of $\mathfrak{g}^{\vee}$. In general, $\fl^{\vee}$ is not (it is the centralizer in $\fg^{\vee}$ of a semisimple group element). Using the bijection $\Delta^{\vee} \simeq \Delta$, we can produce from $\Delta^{\vee}_{\gamma}$ and $\Delta^{\vee}_{\gamma,0}$ two subsystems of $\Delta$
$$\Delta_{\gamma} := \left(\Delta^{\vee}_{\gamma}\right)^{\vee} \subset \Delta, \qquad \Delta_{\gamma,0} := \left(\Delta^{\vee}_{\gamma,0}\right)^{\vee} \subset \Delta_{\gamma}.$$
These root systems define reductive Lie algebras, denoted by  $\mathfrak{l}_{\gamma}$ and $\mathfrak{l}_{\gamma,0}$. The smaller Lie algebra $\mathfrak{l}_{\gamma,0}$ can be identified with a Levi subalgebra of $\fg$. In general, there is no natural embedding $\fl_{\gamma} \subset \fg$.

Consider the nilpotent orbits
$$\mathbb{O}^{\vee}_{\gamma} := \Ind^{L_{\gamma}^{\vee}}_{L_{\gamma,0}^{\vee}} \{0\} \subset (\fl_{\gamma}^{\vee})^*, \qquad \mathbb{O}_{\gamma} := \mathsf{D}(\mathbb{O}^{\vee}_{\gamma}) \subset \fl_{\gamma}^*.$$
\begin{prop}\label{prop:maximalitycriterion}
Suppose $V(I) = \overline{\mathbb{O}}$. Then $I$ is a maximal ideal if and only if
$$\codim(\mathbb{O},\cN) = \codim(\mathbb{O}_{\gamma},\cN_{L_{\gamma}}).$$
\end{prop}

\begin{proof}
In \cite[Chp 13.3]{Lusztig1984}, Lusztig defines a map
$$J^{G}_{L_{\gamma}}: \{\mathbb{O}_{L_{\gamma}} \subset \cN_{L_{\gamma}}\} \to \{\mathbb{O} \subset \cN\}$$
called \emph{truncated induction}\index{induction!truncated}. In the special case when $\fl_{\gamma}$ is a Levi subalgebra of $\fg$, $J^{G}_{L_{\gamma}}$ coincides with $\Ind^G_{L_{\gamma}}$. In any case, there is an equality
\begin{equation}\label{eqn:Jpreservescodim}
\codim(\mathbb{O}_{L_{\gamma}},\cN_{L_{\gamma}}) = \codim(J^{G}_{L_{\gamma}}\mathbb{O}_{L_{\gamma}},\cN).
\end{equation}
Both statements appear in \cite[Chp 13.3]{Lusztig1984}. 

Consider the maximal ideal $I_{\mathrm{max}}(\gamma) \in \mathrm{Prim}_{\gamma}(U(\fg))$. By \cite[Prop A2]{BarbaschVogan1985}
\begin{equation}\label{eqn:Vmax}
V(I_{\mathrm{max}}(\gamma)) = \overline{J^G_{L_{\gamma}}(\mathbb{O}_{\gamma})}\end{equation}
By Proposition \ref{prop:propertiesofprim}(ii), $I \subseteq I_{\mathrm{max}}(\gamma)$ and therefore
\begin{equation}\label{eqn:inclusion}\overline{J^G_{L_{\gamma}}(\mathbb{O}_{\gamma})} = V(I_{\mathrm{max}}(\gamma)) \subseteq V(I) = \overline{\mathbb{O}}.\end{equation}
Now
\begin{align*}
    I \text{ is maximal} &\iff I=I_{\mathrm{max}}(\gamma) & \text{Proposition \ref{prop:propertiesofprim}(ii)} \\
    &\iff V(I) = V(I_{\mathrm{max}}(\gamma)) & \text{Proposition \ref{prop:propertiesofprim}(iii)}\\
    &\iff \mathbb{O} = J^G_{L_{\gamma}}(\mathbb{O}_{\gamma}) & \text{(\ref{eqn:Vmax})}\\
    &\iff \codim(\mathbb{O},\cN) = \codim(\mathbb{O}_{\gamma},\cN_{L_{\gamma}}). & \text{(\ref{eqn:Jpreservescodim}),(\ref{eqn:inclusion})}
\end{align*}
\end{proof}

\begin{rmk}\label{rmk:maximalitycriterion}
Suppose $\fl_{\gamma}$ is a Levi subalgebra of $\fg$. Then $I \subset U(\fg)$ is maximal if and only if $\mathbb{O} = \Ind^G_{L_{\gamma}} \mathbb{O}_{\gamma}$. This follows easily from the proof of Proposition \ref{prop:maximalitycriterion}. 
\end{rmk}

\section{Harish-Chandra bimodules}\label{subsec:HCbimodsclassical}

Let $\cB$ be a $U(\fg)$-bimodule. A \emph{compatible filtration} on $\cB$ is an ascending filtration by subspaces
$$0 = \B_{-1} \subseteq \B_0 \subseteq \B_1 \subseteq ...,  \qquad \bigcup_{i=0}^{\infty} \B_i = \B$$
such that
$$U_i(\fg)\B_j \subseteq \B_{i+j} \quad \text{and} \quad [U_i(\fg),\cB_j] \subseteq \B_{i+j-1}, \qquad \forall i, j \in \ZZ_{\geq 0}.$$
Here, $\CC=U_0(\fg) \subset U_1(\fg) \subset ...$ denote the usual (PBW) filtration of $U(\fg)$. Under the conditions above, $\gr(\B)$ has the structure of a graded $S(\fg)$-module (because $\gr U(\fg)\simeq S(\fg)$). A compatible filtration is \emph{good} if $\gr(\B)$ is finitely-generated for $S(\fg)$. \index{Harish-Chandra bimodule!for $U(\fg)$}

\begin{definition}\label{def:HCbimodsclassical}
A \emph{Harish-Chandra bimodule} is a $U(\fg)$-bimodule which admits a good filtration. A morphism of Harish-Chandra bimodules is a morphism of $U(\fg)$-bimodules. Denote the category of Harish-Chandra bimodules by $\HC(U(\fg))$.
\end{definition}

Note that a $U(\fg)$-bimodule $\cB$ is Harish-Chandra if and only if it is finitely-generated (as a left or right $U(\fg)$-module) and the adjoint action $\ad$ of $\fg$ on $\cB$ is locally finite. The following lemma is standard, and we omit the proof.

\begin{lemma}\label{lem:HC_fin_length}
For $\cB \in \HC(U(\fg))$, the following conditions are equivalent:
\begin{itemize}
    \item[(i)] Every irreducible $\fg$-representation appears with finite multiplicity in $\cB$ (regarded as a $\fg$-module via $\ad$).
    \item[(ii)] $\cB$ has finite length.
    \item[(iii)] $\cB$ has finite support over $\Spec(\fZ(\fg))$.
\end{itemize}
\end{lemma}
In particular, if $\cB \in \HC(U(\fg))$ has left and right infinitesimal characters, then $\cB$ satisfies properties (i) and (ii) above.

We say that a Harish-Chandra bimodule is $G$-\emph{equivariant} if $\ad$ integrates to an {algebraic action} of $G$. Denote the full subcategory of $G$-equivariant bimodules by  $\HC^G(U(\fg)) \subset \HC(U(\fg))$. Let $G^{\mathrm{sc}}$ denote the universal cover of $G$ and $K$ the kernel of the covering map $G^{\mathrm{sc}} \to G$. Then $\HC(U(\fg)) = \HC^{G^{\mathrm{sc}}}(U(\fg))$ and $\HC^G(U(\fg))$ is the full subcategory consisting of all bimodules $\cB \in \HC^{G^{\mathrm{sc}}}(U(\fg))$ with trivial $K$-action. 

Given a Harish-Chandra bimodule $\B$ and a chosen good filtration, we can define the \emph{associated variety}\index{associated variety!of Harish-Chandra bimodule} of $\B$
$$\mathcal{V}(\B) := \mathrm{Supp}(\gr(\B)) = V(\mathrm{Ann}_{S(\fg)}\gr(\B)) \subset \fg^*,$$
where $V(\bullet)$ denotes the zero locus of an ideal in $S(\fg)=\CC[\fg^*]$. 

If $Z$ is an irreducible component of $\mathcal{V}(\B)$, we can define the \emph{generic multiplicity}\index{generic multiplicity} of $\B$ along $Z$
$$m_Z(\B) := \mathrm{mult}_Z(\gr(\B)) \in \ZZ_{\geq 0}.$$
By a standard argument, both $\mathcal{V}(\B)$ and $m_Z(\B)$ are independent of the filtration. Furthermore, $\mathcal{V}(\B) = V(\mathrm{LAnn}(\B)) = V(\mathrm{RAnn}(\B))$ (where $\mathrm{LAnn}(\cB)$ and $\mathrm{RAnn}(\cB)$ denote the left and right annihilatiors, and$V(\bullet)$ is as defined in Section \ref{subsec:assvar}).

Next, 
we recall the notion of parabolic induction for Harish-Chandra bimodules. Suppose $Q \subset G$ is a parabolic subgroup with Levi decomposition $Q=MU$. If $\alpha_L,\alpha_R \in \fX(\fm)$ and $\alpha_L - \alpha_R$ integrates to a character of $M$, we can define a one-dimensional Harish-Chandra bimodule $\CC(\alpha_L,\alpha_R) \in \HC^M(U(\fm))$ by the formula
\begin{equation}\label{eqn:1dimbimod}
\xi t = \alpha_L(\xi)t \quad t \xi = \alpha_R(\xi)t, \qquad \xi \in \fm, \ t \in \CC(\alpha_L,\alpha_R).\end{equation}
Starting with a bimodule $\B_M \in \HC^M(U(\fm))$, we will construct an induced bimodule $\Ind^G_M \B_M \in \HC^G(U(\fg))$, thus defining a functor
$$\Ind^G_M: \HC^M(U(\fm)) \to \HC^G(U(\fg))$$
called \emph{parabolic induction}\index{induction!of Harish-Chandra bimodules}. First, let $\B_M^{\#} := \B_M \otimes \CC(\rho_{\mathfrak{u}},\rho_{\mathfrak{u}})$, where $\rho_{\mathfrak{u}}$ is the one-dimensional representation of $\mathfrak{m}$ defined by
\begin{equation}\label{eqn:defofrho}
\rho_{\mathfrak{u}}(\xi) = \frac{1}{2}\mathrm{Tr}(\ad(\xi)|_{\mathfrak{u}}), \qquad \xi \in \fm.
\end{equation}
Let $\mathfrak{q}^- = \mathfrak{m} \oplus \mathfrak{u}^- \subset \fg$ be the opposite parabolic, and regard $\B_M^{\#}$ as a $U(\mathfrak{q})-U(\mathfrak{q}^-)$-bimodule via the quotient maps $\mathfrak{q} \to \mathfrak{m}$ and $\mathfrak{q}^- \to \mathfrak{m}$. Next define
$$\widetilde{\Ind}^G_M \B_M = \mathrm{Hom}_{U(\fq) \otimes U(\fq^{-})} \left(U(\fg) \otimes U(\fg), \B_M^{\#}\right)$$
where $\Hom$ is defined using the left action of $\mathfrak{q}$ on the first $U(\fg)$ factor and the right action of $\mathfrak{q}^{-}$ on the second. Regard $\widetilde{\Ind}^G_M \B_M$ as a $U(\mathfrak{g})$-bimodule via
$$u_1fu_2(v_1 \otimes v_2) = f(v_1u_1\otimes u_2v_2), \qquad u_1,u_2,v_1,v_2 \in U(\fg), \ f \in \widetilde{\Ind}^G_M \B_M.$$
Finally, let $\Ind^G_M \B_M$ be the subspace of $\widetilde{\Ind}^G_M \B_M$ consisting of {$\operatorname{ad}(\fg)$-finite elements, i.e. elements of $\widetilde{\Ind}^G_M \B_M$ belonging to a
 finite-dimensional  $\operatorname{ad}(\g)$-stable subspace (this is clearly a subbimodule). }

\begin{prop}[\cite{Vogan1981}, Chapter 6]\label{prop:propsofind}
The construction given above defines an exact covariant functor
$$\Ind^G_M: \HC^M(U(\fm))\to \HC^G(U(\fg)).$$
This functor enjoys the following properties:
\begin{itemize}
    \item[(i)] $\Ind^G_M$ depends only on $M$ (and not on $Q$).
    \item[(ii)] If $L \subset M$ is a Levi subgroup of $M$, then there is a natural isomorphism
    $$\Ind^G_L \simeq \Ind^G_M \circ \Ind^M_L.$$
    \item[(iii)] If $\cB_M \in \HC^M(U(\fm))$, then as representations of $G$
    $$\left(\Ind^G_M \cB_M\right)|_G \simeq \mathrm{AlgInd}^G_M (\cB_M|_G),$$
    where $\mathrm{AlgInd}^G_M$ is the functor of algebraic induction (see \cite[Sec 3.3]{Jantzen}).
    \item[(iv)] If $\fh \subset \fm$ and $\B_M$ has left and right infinitesimal characters $(\gamma_{\ell}, \gamma_r) \in \fh^*/W_M \times \fh^*/W_M$, then $\Ind^G_M\B_M$ has left and right infinitesimal characters $(\gamma_{\ell}, \gamma_r) \in \fh^*/W \times \fh^*/W$.
\end{itemize}
\end{prop}
For what follows, we will need an alternative description of $\Ind^G_M$. For $\B_M \in \HC^M(U(\fm))$, define $U(\fg)$-$U(\fm)$-bimodules
$$\Delta^Q(\B_M) := U(\fg) \otimes_{U(\mathfrak{q})}  \B_M, \qquad \nabla^Q(\B_M) := \Hom^{\mathrm{fin}}_{U(\mathfrak{q}^{-})}(U(\fg), \B_M).$$
Here $\Hom^{\mathrm{fin}}$ denotes the direct sum of Homs from graded components as in the definition
of (parabolic) dual Verma modules. {Here the grading is with respect to any one-parameter subgroup $\mathbb{G}_m\rightarrow G$ such that its centralizer is $M$ and the action on $\mathfrak{q}$ is by characters $t\mapsto t^k$, where $k\in \ZZ_{\geqslant 0}$}. Let $\cA_M = U(\fm)/\mathrm{RAnn}(\B_M^{\#})$ and consider the space of right $\cA$-module homomorphisms
$$\Hom_{\cA_M}(\Delta^Q(\cA_M), \nabla^Q(\B_M^{\#})).$$
This space has the structure of a $U(\fg)$-bimodule via the left actions of $U(\fg)$ on $\Delta^Q(\cA_M)$ and $\nabla^Q(\B_M^{\#})$. 

Consider the {subbimodule of $\ad(\fg)$-finite vectors}
$$\Hom^{\fg-\mathrm{fin}}_{\cA_M}(\Delta^Q(\cA_M), \nabla^Q(\B_M^{\#}))\subset \Hom_{\cA_M}(\Delta^Q(\cA_M), \nabla^Q(\B_M^{\#})).$$

\begin{lemma}\label{lem:bimodinduction}
There is a natural isomorphism of $U(\fg)$-bimodules
$$\Hom_{\cA_M}^{\fg-\mathrm{fin}}(\Delta^Q(\cA_M), \nabla^Q(\B_M^{\#}))\simeq \Ind^G_M \B_M.$$
\end{lemma}

\begin{proof}
By the tensor-hom adjunction
\begin{align*}
    \Hom_{\cA_{M}}^{\fg-\mathrm{fin}}(\Delta^Q(\cA_M), \nabla^Q(\B_M^{\#})) &=
    \Hom_{\cA_M}^{\fg-\mathrm{fin}}(U(\fg) \otimes_{U(\fq)} \cA_M, \Hom_{U(\fq^{-})}^{\mathrm{fin}}(U(\fg),\B_M^{\#}))\\
    &\simeq \Hom_{U(\fq)}^{\fg-\mathrm{fin}}(U(\fg), \Hom_{\cA_M}(\cA_M,\Hom_{U(\fq^{-})}^{\mathrm{fin}}(U(\fg),\B_M^{\#})))\\
    &\simeq \Hom_{U(\fq)}^{\fg-\mathrm{fin}}(U(\fg), \Hom_{U(\fq^{-})}^{\mathrm{fin}}(U(\fg), \B_M^{\#})).
\end{align*}
Under the natural isomorphism
$$\Hom_{\CC}(U(\fg), \Hom_{\CC}(U(\fg),\B_M^{\#})) \simeq \Hom_{\CC}(U(\fg) \otimes_{\CC} U(\fg), \B_M^{\#})$$
the subspace $\Hom_{U(\fq)}(U(\fg), \Hom_{U(\fq^{-})}^{\mathrm{fin}}(U(\fg), \B_M^{\#}))$ maps to $\Hom_{U(\fq) \otimes U(\fq^{-})}(U(\fg) \otimes U(\fg), \B_M^{\#}) = \widetilde{\Ind}^G_M \B_M$. The lemma follows immediately. 
\end{proof}

We conclude this section by recalling the Langlands classification of irreducible Harish-Chandra bimodules and listing some easy properties of the resulting parameterization.\index{Langlands parameter}

\begin{definition}\label{def:Langlands}
A \emph{Langlands parameter} for $G$ is a pair $(\lambda_{\ell}, \lambda_r) \in \mathfrak{h}^* \times \mathfrak{h}^*$ such that $\lambda_{\ell} - \lambda_r$ integrates to a character of $H$. For any such pair, let $I_H^G(\lambda_{\ell},\lambda_r) = \mathrm{Ind}^G_H \CC(\lambda_{\ell}, \lambda_r)$ (cf. (\ref{eqn:1dimbimod})).
\end{definition}

Note: in the literature, the term `Langlands parameter' is often reserved for the pair $(\lambda_{\ell}-\lambda_r, \lambda_{\ell} + \lambda_r)$, where $\lambda_{\ell}-\lambda_r$ and $\lambda_{\ell}+\lambda_r$ are called, respectively, the `discrete' and `continuous' parameters. We prefer the more symmetric formulation above.

\begin{theorem}[\cite{Zelobenko}]\label{thm:Langlands}
The following are true:
\begin{itemize}
    \item[(i)] If $(\lambda_{\ell}, \lambda_r)$ is a Langlands parameter for $G$, then $I^G_H(\lambda_{\ell},\lambda_r)$ has a unique irreducible subquotient $\overline{I}_H^G(\lambda_{\ell}, \lambda_r)$ with the following property: $\overline{I}_H^G(\lambda_{\ell}, \lambda_r)$ contains the irreducible $G$-representation of extremal weight $\lambda_{\ell}-\lambda_r$.
    \item[(ii)] Every irreducible $G$-equivariant Harish-Chandra bimodule is isomorphic to $\overline{I}^G_H(\lambda_{\ell},\lambda_r)$ for some Langlands parameter $(\lambda_{\ell},\lambda_r)$.
    \item[(iii)] If $(\lambda_{\ell},\lambda_r)$ and $(\lambda_{\ell}', \lambda_r')$ are Langlands parameters, then $\overline{I}^G_H(\lambda_{\ell}, \lambda_r) \simeq \overline{I}^G_H(\lambda_{\ell}', \lambda_r')$ if and only if there is an element $w \in W$ such that
$$\lambda_{\ell}' = w\lambda_{\ell}, \qquad \lambda_r' = w\lambda_r.$$
\end{itemize}
\end{theorem}

\begin{lemma}\label{lem:Langlands}
Suppose $M \subset G$ is a Levi subgroup and $(\lambda_{\ell},\lambda_r)$ is a Langlands parameter. The following are true:
\begin{itemize}
    \item[(i)] $\Ind^G_M I^M_H(\lambda_{\ell},\lambda_r) \simeq I^G_H(\lambda_{\ell},\lambda_r)$.
    \item[(ii)] If $\Ind^G_M \overline{I}^M_H(\lambda_{\ell},\lambda_r)$ is irreducible, then $\Ind^G_M \overline{I}^M_H(\lambda_{\ell},\lambda_r) \simeq \overline{I}^G_H(\lambda_{\ell},\lambda_r)$.
\end{itemize}
\end{lemma}

\begin{proof}
(i) follows immediately from the transitivity of induction, see Proposition \ref{prop:propsofind}(ii). We proceed to proving (ii). Since $\overline{I}^M_H(\lambda_{\ell},\lambda_r)$ is a subquotient of $I^M_H(\lambda_{\ell},\lambda_r)$ and $\Ind^G_M$ is exact, $\Ind^G_M \overline{I}^M_H(\lambda_{\ell},\lambda_r)$ is a subquotient of $\Ind^G_M I(\lambda_{\ell},\lambda_r) \simeq I^G_H(\lambda_{\ell},\lambda_r)$. Let $V_G$ (resp. $V_M$) denote the irreducible representation of $M$ (resp. $G$) of extremal weight $\lambda_{\ell} - \lambda_r$. It suffices to show that $V_G \hookrightarrow \Ind^G_M \overline{I}^M_H(\lambda_{\ell},\lambda_r)$. By definition, $V_M \hookrightarrow \overline{I}^M_H(\lambda_{\ell},\lambda_r)$. Thus by Proposition \ref{prop:propsofind}(iii), $\mathrm{AlgInd}^G_M V_M \hookrightarrow \Ind^G_M \overline{I}^M_H(\lambda_{\ell},\lambda_r)$. By Frobenius reciprocity, there is an isomorphism
$$\Hom_G(V_G,\mathrm{AlgInd}^G_M V_M) \simeq \Hom_M(V_G,V_M)$$
Note that an extremal weight for $G$ in $V_G$ is also an extremal weight for $M$. Thus, $\Hom_M(V_G,V_M)\neq 0$, as desired.
\end{proof}

\section{${W}$-algebras}\label{subsec:W}
To each nilpotent orbit $\OO$, one can attach a finitely generated associative algebra, called the $W$-\emph{algebra} for $\OO$. In this section, we will recall some basic facts about these algebras and their connection to Harish-Chandra bimodules. For details and proofs, we refer the reader to \cite{GanGinzburg}, \cite{Premet2002}, and \cite{LosevICM}.\index{$W$-algebra}

View $\OO$ as an orbit in $\fg$ 
using a $G$-invariant identification $\fg \simeq \fg^*$.
Pick $e\in \OO$, and fix an $\mathfrak{sl}(2)$-triple $(e,f,h)$. The \emph{Slodowy slice}\index{Slodowy slice} associated to $e$ is the affine subspace of $\fg$ defined by
$$S=e+\Ker (\ad f) \subset \fg.$$
Using the same identification $\fg \simeq \fg^*$, we can regard $S$ as an affine subspace of $\fg^*$. The ring of regular functions $\CC[S]$ has a number of additional structures. First, the Poisson structure on $S(\fg) \simeq \CC[\fg^*]$ (defined using the Lie bracket) induces a Poisson structure on $\CC[S]$, see \cite[Sec 3]{GanGinzburg}. Second, the reductive group $R:= Z_G(e,f,h)$ stabilizes $S$ and thus acts on $\CC[S]$, see \cite[Sec 2]{LosevICM}. This action is Hamiltonian (cf. Section \ref{subsec:equivariant})---the moment map is the restriction to $S$ of the projection $\fg^* \to \mathfrak{r}^*$. Lastly, there is a positive algebra grading on $\CC[S]$. To define it, consider the co-character $\gamma: \CC^\times\rightarrow G$ corresponding to $h$, and let $\CC^{\times}$ act on $\fg^*$ by
$$t \cdot \xi=t^{-2}\Ad^*(\gamma(t))(\xi), \qquad \xi\in \fg^*.$$

This action (called the \emph{Kazhdan action})\index{Kazhdan action} stabilizes $S$, fixes $e$, and contracts the former onto the latter. Thus, it defines a positive grading (called the \emph{Kazhdan grading}) on $\CC[S]$, see \cite[Sec 4]{GanGinzburg}.

The $W$-algebra for $\OO$ is a certain filtered associative algebra $\cW$ such that $\gr(\cW) \simeq \CC[S]$ (as graded Poisson algebras). {It can be defined as a quantum Hamiltonian reduction of $U(\fg)$}{: there is a subalgebra $\mathfrak{m}$ and a character $\psi:\mathfrak{m}\rightarrow \CC$ such that $\cW=\operatorname{End}(U(\fg)\otimes_{U(\mathfrak{m})}\CC_\psi)^{\mathrm{opp}}$, where $\CC_\psi$ denotes the 1-dimensional $\fm$-representation on which $\mathfrak{m}$ acts via $\psi$}. For details, we refer the reader to \cite{Premet2002}, see also \cite[Sec 2.3]{LosevICM}. For what follows, it is important to note that $\cW$ carries a Hamiltonian action of $R$ lifting that on $\CC[S]$, see \cite[Lemma 2.4]{Premetstabilizer}. 

Similarly to Definition \ref{def:HCbimodsclassical}, one can define the category of $R$-equivariant Harish-Chandra bimodules for $\cW$, denoted $\HC^R(\cW)$ (see \cite[Sec. 2.5]{Losev2011})\index{Harish-Chandra bimodule!for $\mathcal{W}$}. Write $\HC^R_{\mathrm{fin}}(\cW) \subset \HC^R(\cW)$ for the subcategory of finite-dimensional bimodules. In \cite[Sec 3.4]{Losev2011}, the first-named author defines a pair of functors\index{extension functor}\index{restriction functor}
$$\bullet_{\dagger}: \HC^G(U(\fg))\to \HC^R(\cW), \qquad \bullet^{\dagger}: \HC^R_{\mathrm{fin}}(\cW) \to  \HC^G_{\overline{\OO}}(U(\fg)).$$
Here and below we write $\HC^G_{\overline{\OO}}(U(\fg))$  for the full subcategory in $\HC^G(U(\fg))$ consisting of all bimodules $\mathcal{B}$ with $\mathcal{V}(\mathcal{B})\subset \overline{\OO}$.

Some properties of these functors are catalogued below.

\begin{prop}\label{prop:propsofdagger}
The following are true:
\begin{itemize}
    \item[(i)] $\bullet_\dagger$ is an exact monoidal functor.
    \item[(ii)] A good filtration on $\cB$ gives rise to a good filtration on $\cB_\dagger$ with the property that the graded $\CC[S]$-module  $\gr (\cB_\dagger)$ is $R$-equivariantly identified with the pullback of $\gr\cB$ to $S$. 
    \item[(iii)] In particular,  $\OO\not\subset\mathcal{V}(\cB)$ if and only if $\cB_\dagger=0$. And if $\overline{\OO}$ is an irreducible component of $\mathcal{V}(\cB)$, then $\cB_\dagger$ is finite dimensional, and  $\dim \cB_\dagger=m_{\overline{\OO}}(\cB)$.
    \item[(iv)] $\bullet^\dagger$ is right adjoint to the functor $\bullet_{\dagger}: \HC^G_{\overline{\OO}}(U(\fg))\rightarrow \HC^R_{\mathrm{fin}}(\cW)$. 
    \item[(v)] If $\mathcal{V}(\cB) = \overline{\OO}$, both kernel and cokernel of the adjunction morphism $\cB\to (\cB_\dagger)^\dagger$ are supported in $\partial \mathbb{O}$. 
    \item[(vi)] We have a bi-functorial isomorphism $\operatorname{Hom}_{U(\g)}(\cB^1,\cB^2)_\dagger\xrightarrow{\sim} \operatorname{Hom}_{\cW}(\cB^1_\dagger,\cB^2_\dagger)$.
\end{itemize}
\end{prop}
\begin{proof}
Part (ii) follows from \cite[Lemma 3.3.2]{Losev2011}. (i) and (iii)-(v) are parts of
\cite[Proposition 3.4.1]{Losev2011}. (vi) is proved similarly to \cite[Lemma 3.10]{BezLosev}.
\end{proof}

We will make repeated use of the following lemma in our study of unipotent ideals.

\begin{lemma}\label{lem:adjunctionmorphismiso}
Let $\cB\in \HC^G_{\overline{\OO}}(U(\fg))$. Assume $\mathcal{V}(\cB)=\overline{\OO}$ and $I=\mathrm{LAnn}_{U(\fg)}(\cB)$ is a maximal ideal. Then the adjunction unit $\cB\to (\cB_\dagger)^\dagger$ is an isomorphism. 
\end{lemma}

\begin{proof}
Write $K, C\in \HC^G(U(\fg))$ for the kernel and the cokernel of $\cB\to (\cB_\dagger)^\dagger$. Note that $\mathcal{V}(K),\mathcal{V}(C) \subseteq \partial \OO$ by  \cref{prop:propsofdagger}(v). In particular, $\mathcal{V}((\cB_{\dagger})^{\dagger}) = \overline{\OO}$. We will show that $K=C=0$.

Since $K \subset \cB$,
$I =\mathrm{LAnn}_{U(\fg)}(\cB) \subseteq \mathrm{LAnn}_{U(\fg)}(K)$. This inclusion must be strict since $\mathcal{V}(K) \subseteq \partial \OO \subsetneq \overline{\OO} = \mathcal{V}(\cB)$. Hence, $\mathrm{LAnn}_{U(\fg)}(K)=U(\fg)$, i.e. $K=0$. 

Now we show that $C=\{0\}$. As in the previous paragraph, it is enough to show that $IC=\{0\}$. Since $I\mathcal{B}=\{0\}$ and $\mathcal{B}\hookrightarrow 
(\cB_\dagger)^\dagger$, we see that the multiplication map $I\otimes_{U(\fg)}
(\cB_\dagger)^\dagger\rightarrow (\cB_\dagger)^\dagger$ factors through 
$I\otimes_{U(\fg)}
C\rightarrow (\cB_\dagger)^\dagger$. The image is a sub-bimodule supported on $\partial \OO$. From (iii) and (iv) of Proposition \ref{prop:propsofdagger} it follows that any bimodule in the image of $\bullet^\dagger$ has no nonzero subbimodules supported on $\partial\OO$. In particular, $I(\cB_\dagger)^\dagger=\{0\}$ and hence $IC=\{0\}$. 
\end{proof}

Let $\mathrm{Prim}_{\overline{\OO}}(U(\fg))$ denote the set of primitive ideals $I \subset U(\fg)$ such that the associated variety $V(I)$ is equal to $\overline{\OO}$. Let $\mathrm{Id}_{\mathrm{fin}}(\cW)$ (resp. $\mathrm{Prim}_{\mathrm{fin}}(\cW)$) denote the set of ideals $J \subset \cW$ (resp. primitive ideals $J \subset \cW$) of finite codimension. Note that a primitive ideal is of finite codimension if and only if it is the annihilator of a finite-dimensional irreducible module.

If $I \in \mathrm{Prim}_{\overline{\OO}}(U(\fg))$, then $I_\dagger \in \mathrm{Id}_{\mathrm{fin}}(\cW)$ by (i) and (iii) of Proposition \ref{prop:propsofdagger}. Conversely, suppose $J \in \mathrm{Prim}_{\mathrm{fin}}(\cW)$. Since the $R$ action on $\cW$ is Hamiltonian, $J$ is stable under the identity component $R^{\circ}$. Let $\underline{J}$ denote the intersection of all $R$-conjugates of $J$. Then $\cW/\underline{J}\in \HC_{\mathrm{fin}}^R(\cW)$. By (i) of Proposition \ref{prop:propsofdagger}, $U(\fg)_\dagger=\cW$. Applying adjunction to the map $\cW\twoheadrightarrow \cW/\underline{J}$, we get a bimodule homomorphism  $U(\fg)\rightarrow (\cW/\underline{J})^{\dagger}$.  Its kernel is an element in $\mathrm{Prim}_{\overline{\OO}}(U(\fg))$, see \cite[Thm 1.2.2(iv),(v) and (vi)]{Losev3}, which we will henceforth denote by $J^\ddag$. Thus, we get maps
$$\bullet_{\dagger}: \mathrm{Prim}_{\overline{\OO}}(U(\fg)) \to \mathrm{Id}_{\mathrm{fin}}(\cW), \qquad \bullet^{\ddag}: \mathrm{Prim}_{\mathrm{fin}}(\cW) \to \mathrm{Prim}_{\overline{\OO}}(U(\fg)).$$ 
We will need the following result from \cite{Losev2011}.

\begin{theorem}[Thm 1.2.2, Section 4.2, \cite{Losev2011}]\label{thm:daggersurjective}
The following are true:
\begin{enumerate}
\item 
The map $\bullet^{\ddag}: \mathrm{Prim}_{\mathrm{fin}}(\cW) \to \Prim_{\overline{\mathbb{O}}}(U(\fg))$ is surjective. Every fiber is a single orbit for the $R$-action on $\mathrm{Prim}_{\mathrm{fin}}(\cW)$. 
\item The map $\bullet_{\dagger}: \mathrm{Prim}_{\overline{\OO}}(U(\fg)) \to \mathrm{Id}_{\mathrm{fin}}(\cW)$ is injective. Moreover,
for every $I\in \Prim_{\overline{\mathbb{O}}}(U(\fg))$, we have 
$I_{\dagger} = \bigcap J$, where the intersection is taken over all $J\in  \mathrm{Prim}_{\mathrm{fin}}(\cW)$ such that $J^{\ddag} = I$.
\end{enumerate}
\end{theorem}


In view of Theorem \ref{thm:daggersurjective}, we can make the following definition.

\begin{definition}\label{def:Wdimension}
If $I \in \mathrm{Prim}_{\overline{\mathbb{O}}}(U(\fg))$, define
$$\operatorname{\operatorname{\cW-\dim}}(I) := \sqrt{\dim(\cW/J)}$$
for $J \in \mathrm{Prim}_{\mathrm{fin}}(\cW)$ such that $I = J^{\ddag}$. The number $\operatorname{\cW-\dim}(I)$ is called the $\cW$-\emph{dimension} of $I$. It is well-defined by Theorem \ref{thm:daggersurjective} and coincides with the dimension of the irreducible $\cW$-module annihilated by $J$.\index{$\mathcal{W}$-dimension}
\end{definition}

Let $J\subset \cW$ be a two-sided ideal of codimension $1$ and define $\underline{J} \subset \cW$ as in the paragraph preceding Theorem \ref{thm:daggersurjective}. Let $R_1$ denote the stabilizer of $J$ in $R$ and $H_1$ its preimage in $G_e$, a finite index subgroup. The Harish-Chandra bimodule 
$(\cW/\underline{J})^{\dagger}$ has a natural algebra structure, see the discussion preceding \cite[Lemma 5.2]{Losev4}. The homomorphism $U(\fg)\rightarrow (\cW/\underline{J})^\dagger$ is one of algebras. By the construction of \cite[Section 3.3]{Losev2011}, this algebra has a distinguished algebra filtration. The next lemma follows from the proof of \cite[Lemma 3.3.3]{Losev2011}.

\begin{lemma}\label{lem:upper_dag_1dim}
There is an inclusion of $G$-equivariant graded algebras
$\gr[(\cW/\underline{J})^\dagger]\hookrightarrow \mathbb{C}[G/H_1]$, and hence an inclusion of $G$-representations $(\cW/\underline{J})^\dagger\hookrightarrow \mathbb{C}[G/H_1]$.
\end{lemma}

\begin{cor}\label{Cor:completely_prime} Let $J$ be a codimension $1$ ideal in $\cW$.  Then the ideal
$J^\ddag$ is completely prime.
\end{cor}
\begin{proof}
Let $H_1,\underline{J}$ have the same meaning as in Lemma  \ref{lem:upper_dag_1dim}.
Since the variety $G/H_1$ is irreducible, the algebra $\CC[G/H_1]$ is a domain. Thus, by Lemma \ref{lem:upper_dag_1dim}, the same is true of  $\gr[(\cW/\underline{J})^\dagger]$, and hence of $(\cW/\underline{J})^\dagger$. Note that there is an algebra embedding $U(\fg)/J^\ddag\hookrightarrow (\cW/\underline{J})^\dagger$. The corollary follows.
\end{proof}

\chapter{Deformations and quantizations of conical symplectic singularities}\label{sec:symplectic}

In this chapter, we will review some basic facts about Poisson deformations and filtered quantizations of conical symplectic singularities.

\section{Poisson deformations and filtered quantizations}\label{subsec:quant}

Let $d \in \ZZ_{>0}$ and let $A$ be a \emph{graded Poisson algebra} of degree $-d$. By this, we will mean a finitely-generated commutative associative unital algebra equipped with two additional structures: an algebra grading 
$$A = \bigoplus_{i=-\infty}^{\infty} A_i$$
and a Poisson bracket $\{\cdot, \cdot\}$ of degree $-d$
$$\{A_i, A_j\}\subset A_{i+j-d}, \qquad i,j \in \ZZ.$$
To any algebra of this form, one can associate various classes of filtered associative algebras which `deform' the Poisson bracket in a suitable sense. In {this monograph}, we will consider two such classes of algebras: \textit{filtered Poisson deformations} and \textit{filtered quantizations}. The definitions are below.\index{quantization!filtered}\index{Poisson deformation!filtered}

\begin{definition}\label{def:poissondef}
\leavevmode
\begin{itemize}
    \item A \emph{filtered Poisson deformation} of $A$ is a pair $(\cA^0,\theta)$ consisting of
\begin{itemize}
    \item[(i)] a Poisson algebra $\cA^0$, equipped with a complete and separated filtration by subspaces
    $$\cA^0 = \bigcup_{i=-\infty}^{\infty} \cA^0_{\leq i}, \qquad ... \subseteq \cA^0_{-1} \subseteq \cA^0_0 \subseteq \cA^0_1 \subseteq ...$$
    such that
    $$\{\cA^0_{\leq i},\cA^0_{\leq j}\} \subseteq \cA^0_{\leq i+j-d}, \qquad i,j \in \ZZ,$$
    and
    \item[(ii)] an isomorphism of graded Poisson algebras
    $$\theta: \gr(\cA^0) \simeq A,$$
    where the Poisson bracket on $\gr(\cA^0)$ is defined by
    $$\{a+\cA^0_{\leq i-1}, b + \cA^0_{\leq j-1}\} := \{a,b\} + \cA^0_{\leq i + j -d-1}, \qquad a \in \cA_{\leq i}^0, \ b \in \cA_{\leq j}^0.$$
\end{itemize}
\item An isomorphism of filtered Poisson deformations $(\cA^0_1, \theta_1) \to (\cA^0_2, \theta_2)$ is an isomorphism of filtered Poisson algebras $\phi: \cA^0_1 \to \cA^0_2$ such that $\theta_1 = \theta_2 \circ \gr(\phi)$. Denote the set of isomorphism classes of filtered Poisson deformations of $A$ by $\mathrm{PDef}(A)$. 
\end{itemize}
\end{definition}

\begin{definition}\label{def:filteredquant}
\leavevmode
\begin{itemize}
    \item A \emph{filtered quantization} of $A$ is a pair $(\cA,\theta)$ consisting of
\begin{itemize}
    \item[(i)] an associative algebra $\cA$ equipped with a complete and separated filtration by subspaces
    $$\cA = \bigcup_{i=-\infty}^{\infty} \cA_{\leq i}, \qquad ... \subseteq \cA_{\leq -1} \subseteq \cA_{\leq 0} \subseteq \cA_{\leq 1} \subseteq ...$$
    such that
    $$[\cA_{\leq i}, \cA_{\leq j}] \subseteq \cA_{\leq i+j-d} \qquad i,j \in \ZZ,$$
    and
    \item[(ii)] an isomorphism of graded Poisson algebras
    $$\theta: \gr(\cA) \xrightarrow{\sim} A,$$
    where the Poisson bracket on $\gr(\cA)$ is defined by
    $$\{a+\cA_{\leq i-1}, b+\cA_{\leq j-1}\}=[a,b]+\cA_{\leq i+j-d-1}, \qquad a \in \cA_{\leq i}, \ b \in \cA_{\leq j}.$$
\end{itemize}
\item An isomorphism of filtered quantizations $(\cA_1, \theta_1) \xrightarrow{\sim} (\cA_2, \theta_2)$ is an isomorphism of filtered algebras $\phi: \cA_1 \xrightarrow{\sim} \cA_2$ such that $\theta_1 = \theta_2 \circ \gr(\phi)$. Denote the set of isomorphism classes of quantizations of $A$ by $\mathrm{Quant}(A)$.
\end{itemize}
\end{definition}

\begin{rmk}
In many cases we will consider, the grading on $A$ is by nonnegative integers. In such cases, the filtrations considered in Definitions \ref{def:filteredquant} and \ref{def:poissondef} are automatically complete and separated. 
\end{rmk}

Now let $X$ be a \emph{graded Poisson variety}. By this, we will mean a normal quasi-projective variety with an {algebraic} $\CC^{\times}$-action and a Poisson bracket $\{ \cdot, \cdot\}: \cO_X \otimes \cO_X \to \cO_X$ of degree $-d$. We will consider the so-called \emph{conical topology} on $X$\index{conical topology}. The open subsets in this topology are the Zariski-open subsets $U \subset X$ which are stable under $\CC^{\times}$. Every point in $X$ admits a $\CC^{\times}$-stable open affine neighborhood. {Indeed, since $X$ is quasi-projective and normal, it admits a $\CC^\times$-equivariant embedding into the projective space $\mathbb{P}(V)$ for some finite dimensional rational $\CC^\times$-representation $V$. This follows, for example, by combining \cite[Proposition 3.2.6]{Brion} with \cite[Theorem 5.2.1]{Brion}. For any point $x\in X$ we can then find a homogeneous and $\CC^\times$-semi-invariant polynomial $P\in \CC[V]$ such that $P(x)\neq 0$, while $P$ is zero on $\overline{X}\setminus X$, where the closure of $X$ is taken in $\mathbb{P}(V)$. The intersection with $X$ of the non-vanishing locus of $P$ is a $\CC^\times$-stable affine open neighborhood of $x$ in $X$.}

{Note that }when viewed in the conical topology, $\cO_X$ is a sheaf of graded Poisson algebras.\index{quantization!filtered} 

\begin{definition}\label{def:filteredquantX}
\leavevmode
\begin{itemize}
    \item A \emph{filtered quantization} of $X$ is a pair $(\mathcal{D},\theta)$ consisting of
\begin{itemize}
    \item[(i)] a sheaf $\mathcal{D}$ of associative algebras in the conical topology on $X$, equipped with a complete and separated filtration by subsheaves of vector spaces
    $$\mathcal{D} = \bigcup_{i=-\infty}^{\infty} \mathcal{D}_{\leq i}, \qquad ...\subseteq \mathcal{D}_{\leq -1} \subseteq \mathcal{D}_{\leq 0} \subseteq \mathcal{D}_{\leq 1} \subseteq...$$
    such that
    $$[\mathcal{D}_{\leq i},\mathcal{D}_{\leq j}] \subseteq \mathcal{D}_{\leq i+j-d},$$
    and
    \item[(ii)] an isomorphism of sheaves of graded Poisson algebras
    $$\theta: \gr(\mathcal{D}) \xrightarrow{\sim} \cO_X.$$
\end{itemize}
\item An isomorphism of filtered quantizations $(\mathcal{D}_1, \theta_1) \xrightarrow{\sim} (\mathcal{D}_2, \theta_2)$ is an isomorphism of sheaves of filtered algebras $\phi: \mathcal{D}_1 \to \mathcal{D}_2$ such that $\theta_1 = \theta_2 \circ \gr(\phi)$. Denote the set of isomorphism classes of filtered quantizations of $X$ by $\mathrm{Quant}(X)$.
\end{itemize}
\end{definition}

Filtered Poisson deformations of $X$ can be defined in a similar way. Often, the isomorphism $\theta$ is clear from the context, and will be omitted from the notation. However, the reader should keep in mind that a Poisson deformation $(\cA^0,\theta)$ (resp. filtered quantization $(\cA,\theta)$) is \emph{not} determined up to isomorphism by $\cA^0$ (resp. $\cA$) alone.

Suppose $A=\CC[X]$, where $X$ is a graded {\it affine} Poisson variety. Of course, $A$ is a graded Poisson algebra of degree $-d$. If $\mathcal{D}$ is a filtered quantization of $X$, then $\Gamma(X,\mathcal{D})$ is a filtered quantization of $A$. The correspondence $\mathcal{D} \mapsto \Gamma(X,\mathcal{D})$ defines a bijection $\mathrm{Quant}(X) \xrightarrow{\sim} \mathrm{Quant}(A)$, see e.g. \cite[Sec 1]{Ginzburg1986}. The inverse bijection $\mathrm{Quant}(A) \xrightarrow{\sim} \mathrm{Quant}(X)$ can be defined as follows. Given $\cA \in \mathrm{Quant}(A)$, form the corresponding Rees algebra\index{Rees algebra}
$$\cA_{\hbar} := \bigoplus_{i \in \ZZ} \cA_{\leq i}\hbar^i.$$
Note that $\cA_{\hbar}/(\hbar-1) \simeq \cA$ and $\cA_{\hbar}/(\hbar) \simeq \gr(\cA) \simeq A$. For every homogeneous element $f \in A$ and $i \in \ZZ$, choose a lift $\widetilde{f}_i \in \cA_{\hbar, i} := \cA_{\hbar}/(\hbar^i)$ of $f$. It is easy to check that the multiplicative system $\{ \widetilde{f}_i^n\} \subset \cA_{\hbar,i}$ satisfies the Ore conditions (the key observation is that $\ad(\widetilde{f}_i)^i=0$, which follows from the commutativity of $A$). Hence, there is a localization $\cA_{\hbar,i}[\widetilde{f}_i^{-1}]$, a graded $\CC[\hbar]/(\hbar^i)$-algebra.  It is easy to check that this localization is independent of the lift. Consider the inverse limit $\cA_{\hbar}[f^{-1}] := \lim_{i \to \infty}\cA_{\hbar,i}[\widetilde{f}_i^{-1}]$ (with respect to the natural homomorphisms $\cA_{\hbar,i+1}[\widetilde{f}_{i+1}^{-1}] \to \cA_{\hbar,i}[\widetilde{f}_i^{-1}]$) in the category of graded algebras. Then $\cA[f^{-1}] := \cA_{\hbar}[f^{-1}]/(\hbar-1)$ is a filtered quantization of the graded Poisson algebra $A[f^{-1}]$. The algebras $\cA[f^{-1}]$, as $f$ runs over the homogeneous elements in $A$, form a sheaf on $X$, denoted by $\mathcal{D}$. It is easy to see that $\mathcal{D} \in \mathrm{Quant}(X)$ and $\Gamma(X,\mathcal{D}) \simeq \cA$. This sheaf is called the \emph{microlocalization} of $\cA$ over $X$.\index{microlocalization}

\section{Kleinian Singularities and the McKay correspondence}\label{subsec:Mckay}

\emph{Kleinian singularities}\index{singularity!Kleinian} are an important class of 2-dimensional graded Poisson varieties. In this chapter, we will review some basic facts about such varieties. For details and proofs, we refer the reader to \cite{McKay} and \cite{GonzalezVerdier}.

Let $\Gamma$ be a nontrivial finite subgroup of $\mathrm{Sp}(2)$. Write $\{V_1,...,V_n\}$ for the nontrivial irreducible representations of $\Gamma$ and $V$ for the 2-dimensional tautological representation of $\Gamma \subset \mathrm{Sp}(2)$. For $i,j \in \{1,...,n\}$, define a nonnegative integer
$$m_{ij}:= \dim \Hom_{\Gamma}(V_j, V_i \otimes V).$$
The \emph{McKay graph}\index{McKay graph} of $\Gamma$ is the graph with vertices $\{V_1,...,V_n\}$ and $m_{ij}$ (non-oriented) edges between $V_i$ and $V_j$. 

\begin{theorem}[\cite{McKay}, Props 3,4]\label{thm:McKay}
The McKay graph of $\Gamma$ is a simply laced Dynkin diagram. The passage from $\Gamma$ to its McKay graph defines a bijection
$$\{\text{nontrivial finite subgroups } \Gamma \subset \mathrm{Sp}(2)\} \xrightarrow{\sim} \{\text{simple root systems of type } ADE\}.$$
\end{theorem}
With $\Gamma$ as above, consider the quotient variety $\Sigma := \CC^2/\Gamma$. The variety $\Sigma$ is called the \emph{Kleinian singularity} corresponding to $\Gamma$. 
The symplectic form on $\CC^2$ induces a Poisson bracket on $\CC[\Sigma] = \CC[\CC^2]^{\Gamma}$ of degree $-2$, making $\Sigma$ a graded Poisson variety of the type described in Section \ref{subsec:quant}. 

{Recall the following classical fact, see, for example, \cite[Theorem 2.16]{Kollar}:  every surface has a unique minimal resolution of singularities}. 
Consider the mininal resolution $\rho:\mathfrak{S} \to \Sigma$. The exceptional divisor $\rho^{-1}(0)$ is a union of components $\{C_1,...,C_m\}$, each isomorphic to $\PP^1(\CC)$, intersecting transversely, see e.g. \cite[Sec 6.1]{Slodowy}. The \emph{intersection graph} of $\mathfrak{S}$ is the (undirected) graph with vertices $\{C_1,...,C_n\}$ and an edge between $C_i$ and $C_j$ if and only if $C_i \cap C_j \neq\emptyset$. The correspondence of Theorem \ref{thm:McKay} has a geometric interpretation, due to Gonzalez-Springberg and Verdier.

\begin{theorem}[\cite{GonzalezVerdier}, Thm 2.2]\label{thm:McKaygeometric}
The McKay graph of $\Gamma$ is isomorphic to the intersection graph of $\mathfrak{S}$. 
\end{theorem}

Let $\fg$ be the simple complex Lie algebra of type ADE corresponding to $\Sigma$ via Theorem \ref{thm:McKay}. Choose a Cartan subalgebra $\fh \subset \fg$, and let $\Delta, \Lambda \subset \fh^*$ be the root system and weight lattice, respectively. Choose simple roots $\{\alpha_1,...,\alpha_n\} \subset \Delta$. By Theorems \ref{thm:McKay} and \ref{thm:McKaygeometric}, there are bijections
\begin{equation}\label{eq:bijectionsMcKay}\{\alpha_1,...,\alpha_n\} \simeq \{C_1,...,C_n\} \simeq \{V_1,...,V_n\}\end{equation}
well-defined up to a diagram automorphism. Arrange the indices so that $\alpha_i$ corresponds to $C_i$ and $V_i$. The following proposition is an easy consequence of \cite[Thm 2.2]{GonzalezVerdier}.

\begin{prop}\label{prop:PicSigma}
There is a group isomorphism
$$\sigma: \Lambda \xrightarrow{\sim} \Pic(\mathfrak{S})$$
with the following property: if $\sigma(\lambda)=\mathcal{L}$, then
$$\langle \mathcal{L}, C_i\rangle = \langle \lambda, \alpha_i^{\vee}\rangle, \qquad 1 \leq i \leq n,$$
{where in the left hand side we have the degree of the restriction of $\mathcal{L}$ to $C_i$.}
\end{prop}

\begin{rmk}\label{rmk:PforKleinian}
The map
$$c_1: \Pic(\mathfrak{S}) \to H^2(\mathfrak{S},\ZZ)$$
which takes each line bundle $\mathcal{L}$ to its first Chern class $c_1(\mathcal{L})$ is an isomorphism. Indeed $\mathfrak{S}$ comes equipped with a $\CC^\times$-action which contracts it onto the exceptional divisor of the minimal resolution $\mathfrak{S} \to \Sigma$.
This gives an identification of $H^2(\mathfrak{S},\ZZ)$ with the second cohomology group of the exceptional divisor, which has basis indexed by the irreducible components. The claim that $c_1$ is an isomorphism easily follows.
\end{rmk}

\section{Symplectic singularities}\label{subsec:symplectic}

Let $X$ be a normal Poisson variety.\index{singularity!symplectic}

\begin{definition}[\cite{Beauville2000}, Def 1.1]\label{def:symplecticsing}
We say that $X$ has \emph{symplectic singularities} if
\begin{itemize}
    \item[(i)] the regular locus $X^{\mathrm{reg}} \subset X$ is symplectic; denote the symplectic form by $\omega^{\mathrm{reg}}$. 
    \item[(ii)] there is a resolution of singularities $\rho: Y \to X$ such that $\rho^*(\omega^{\mathrm{reg}})$ extends to a regular (not necessarily symplectic) $2$-form on $Y$.
\end{itemize}
\end{definition}

If (ii) holds for one resolution, it holds for all others, \cite[Sec 1.2]{Beauville2000}. In particular, if a normal Poisson variety $X'$ admits a birational projective Poisson morphism to $X$, then $X'$ has symplectic singularities as well. 

\begin{lemma}[\cite{Beauville2000}, Prop 1.3]\label{lem:rationalsingularities}
Suppose $X$ has symplectic singularities. Then $X$ has rational singularities. That is, for any proper birational map $\eta: Z \to X$ the following are true:
\begin{itemize}
    \item[(i)] $\eta_*\cO_Z \xrightarrow{\sim} \cO_X$,
    \item[(ii)] $R^i\eta_*\cO_Z = 0$ for $i>0$.
\end{itemize}
{Also, $X$ is Gorenstein meaning that the dualizing sheaf is a line bundle.}
\end{lemma}

In {this monograph}, we will give special attention to so-called \emph{conical symplectic singularities}.\index{singularity!conical symplectic}

\begin{definition}\label{def:conicalsymplecticsing}
A \emph{conical symplectic singularity} is a normal graded Poisson variety $X$ with symplectic singularities such that the $\CC^{\times}$-action on $X$ contracts it to a point.
\end{definition}

Note that a conical symplectic singularity $X$ is in particular a graded Poisson variety of the type considered in Section \ref{subsec:quant}.  Furthermore, the contracting $\CC^{\times}$-action guarantees that $X$ is affine. 

\begin{example}\label{example:symplecticsingularity}
The following are examples of conical symplectic singularities
		\begin{itemize}
		    \item[(i)] Let $\Gamma\subset \mathrm{Sp}(2)$ be a finite subgroup. Then the Kleinian singularity $\Sigma = \CC^2/\Gamma$ is a conical symplectic singularity, see  \cite[Prop 2.4]{Beauville2000}. For $\rho$ we take the minimal resolution $\mathfrak{S} \to \Sigma$.
		    \item[(ii)] Let $\fg$ be a complex reductive Lie algebra and let $\mathbb{O} \subset \fg^*$ be a nilpotent orbit. Then $\Spec(\CC[\mathbb{O}])$ is a conical symplectic singularity, see  \cite[Sec 2.5]{Beauville2000}. 
		    \item[(iii)] In the setting of $(ii)$, let $\widetilde{\mathbb{O}} \to \mathbb{O}$ be a connected finite \'{e}tale cover. Then $\Spec(\CC[\widetilde{\mathbb{O}}])$ is a conical symplectic singularity, see \cite[Lem 2.5]{LosevHC}.
		\end{itemize}
	 \end{example}	
	 
For an arbitrary variety $X$, define the subvarieties $X_0,X_1,X_2,...$ as follows: $X_0:=X$ and  $X_{k+1}:=X_k- X_k^{\mathrm{reg}}$. If $X$ is Poisson, then all $X_k$ are Poisson subvarieties of $X$.

\begin{definition}\label{def:fin_many_leaves} We say that $X$ has finitely many (symplectic) leaves if  $X_k^{\mathrm{reg}}$ is a symplectic variety for all $k$. By a symplectic leaf of $X$ we mean an irreducible (i.e. connected) component of $X_k^{\mathrm{reg}}$ for some $k$.\index{symplectic leaf}
\end{definition}

The following lemma is elementary.

\begin{lemma}\label{lem:fin_many_leaves_equiv}
Let $X$ be a Poisson variety with finitely many leaves. Then the leaves are in bijection with prime Poisson ideals in $\CC[X]$; to a sympletic leaf $\mathcal{L}\subset X$ one assigns the ideal of all functions which vanish on $\mathcal{L}$. Furthermore, if $I\subset \CC[X]$ is an arbitrary Poisson ideal, the variety of zeroes $V(I) \subset X$ is a union of symplectic leaves. 
\end{lemma}

The following important result is Theorem 2.3 in \cite{Kaledin2006}.

\begin{prop}\label{Prop:fin_many_leaves_Kaledin}
Suppose $X$ has symplectic singularities. Then $X$ has finitely many leaves. 
\end{prop}


The next proposition characterizes Kleinian singularities among all conical symplectic singularities.

\begin{prop}\label{Prop:Kleinian_uniqueness}
Let $X$ be a 2-dimensional conical symplectic singularity such that the bracket has degree $-d$ for $d>0$. Then there is a unique (up to conjugation) subgroup $\Gamma\subset \operatorname{Sp}(2)$ such that there is a Poisson isomorphism $X\cong \CC^2/\Gamma$. Moreover, this isomorphism can be chosen to be $\CC^\times$-equivariant, where the action of $\CC^2/\Gamma$ is as follows:
\begin{itemize}
\item[(i)] the rescaled usual action of $\CC^\times$ on $\CC^2/\Gamma$ in the case when $\Gamma$ is not of type $A$;
\item[(ii)] the action induced by the following action on $\CC[\CC^2/\Gamma]$: $t.(xy)=t^d xy, t.(x^n)=t^e x^n, t.(y^n)=t^{nd-e}y^n$, where $e\in \{1,\ldots,nd-1\}$, and we assume that $\Gamma$ is diagonalizable in the standard basis $x,y$ of $\CC^2$. Here $\Gamma$ is of type $A_{n-1}$.
\end{itemize}
\end{prop}

This proposition is quite standard. However, we were not able to locate a proof in the literature, so we provide one for the reader's convenience.

\begin{proof}
The proof is in several steps. We write $x_0\in X$ for the cone point.

{\it Step 1}. Note that $x_0$ is a canonical singularity in the sense of \cite[Definition 2.11]{KollarMori1998}: (1) there holds because $X$ is Gorenstein, see
Lemma \ref{lem:rationalsingularities} and (2) there, where one can take $m=1$, is a special case of 
\cite[Theorem 5.10]{KollarMori1998}. Therefore, by \cite[Theorem 4.20]{KollarMori1998}, a complex analytic neighborhood of $x_0$ in $X$, denote it by $U$, is isomorphic, as a complex analytic space, to a neighborhood of $0$ in $\CC^2/\Gamma$, denote it by $U'$. We can assume that $U'=D/\Gamma$, for a disc $D$ around $0$ in $\CC^2$.
Note that $\Gamma$ is uniquely determined, as the fundamental group of $U\setminus\{x_0\}$ acting on the tangent space $T_0D$.

{\it Step 2}. Transfer the Poisson structure from $U$ to $U'$. The resulting structure lifts to a $\Gamma$-invariant Poisson structure on $D$. Hence, by the Darboux theorem, it is $\Gamma$-equivariantly isomorphic to the Poisson structure coming from the standard symplectic form (possibly after shrinking $D$). In other words, we can assume that $U,U'$ are isomorphic as Poisson complex analytic spaces, where the Poisson structure on $U$ is restricted from $X$, and the Poisson structure on $U'$ is restricted from $\Sigma:=\CC^2/\Gamma$. 

{\it Step 3}. Let $\mathsf{eu}$ and $\mathsf{eu}'$ denote the Euler vector fields on $U,U'$ induced by the $\CC^\times$-actions on $X$ and $\Sigma$ (both are taken so that the degrees of the Poisson brackets are equal to $-d$; in the case of $\Sigma$ we rescale the usual action appropriately). 
Let $\xi:=\mathsf{eu}-\mathsf{eu}'$, this is a Poisson vector field (where we identify $U$ and $U'$ by their Poisson isomorphism). Then $\xi$ lifts to $D$ and hence is given by $\{f,\cdot\}$ for an analytic function $f$ on $U'$. Starting from now it is more convenient to work with formal neighborhoods: we have identified (complete local) Poisson rings $\CC[X]^{\wedge_{x_0}}$ and $\CC[\Sigma]^{\wedge_0}$ with Euler derivations $\mathsf{eu}$ and $\mathsf{eu}'$, respectively, and $f\in \CC[\Sigma]^{\wedge_0}$ such that $\mathsf{eu}=\mathsf{eu}'+\{f,\cdot\}$. We can assume that $f$ is in the maximal ideal $\mathfrak{m}$.

{\it Step 4}. Assume now that $\Gamma$ is not of type $A$. Then $\{g,\cdot\}: \CC[\Sigma]^{\wedge_0}\rightarrow \CC[\Sigma]^{\wedge_0}$ is topologically nilpotent for any $g\in \mathfrak{m}$, one can see this from analyzing the grading on $\CC[\Sigma]$. Therefore, $\exp(\{g,\cdot\})$ gives a well-defined Poisson automorphism of $\CC[\Sigma]^{\wedge_0}$.
From here and the observation that the eigenvalues of $\mathsf{eu}'$ on $\mathfrak{m}/\mathfrak{m}^k$ for all $k>1$ are positive integers, it is easy to deduce that there is a unique element $g\in \mathfrak{m}$ such that $\mathsf{eu}'+\{f,\cdot\}=\exp(\{g,\cdot\})\mathsf{eu}'$. Twisting the isomorphism $\CC[X]^{\wedge_{x_0}}\cong \CC[\Sigma]^{\wedge_0}$ with $\exp(\{g,\cdot\})$, we can assume that the isomorphism intertwines the Euler derivations as well. Hence it uniquely extends to a graded Poisson isomorphism $\CC[X]\cong \CC[\Sigma]$ finishing the proof.

{\it Step 5}. Now assume $\Gamma$ is of type $A_n$. The cases $n=1$ and $n>1$ are somewhat different, we consider the latter and leave the former to the reader. The algebra $\C[\Sigma]^{\wedge_0}$ is 
$\CC[e,h,f]/(ef-h^{n+1})$ with brackets of generators given by $\{h,e\}=ne, \{h,f\}=-nf, \{f,e\}=ch^n$ for a suitable nonzero scalar $c$. 
We can twist the $\CC^\times$-action on $\Sigma$ with a Hamiltonian torus action (a rescaling of the action given by $t.(x,y)=(tx,t^{-1}y)$). The subspace of $g\in \mathfrak{m}$ such that $\{g,\cdot\}$ is topoplogically nilpotent is $\C e\oplus \C f\oplus \mathfrak{m}^2$. So we can find suitable $g$ such that $\mathsf{eu}'-\exp(\{g,\cdot\})\mathsf{eu}=a\{h,\cdot\}$ for $a\in \C$. 
We observe that $\mathsf{eu}'$ and $\{h,\cdot\}$ commute.
Then using the integrality of eigenvalues it is easy to see that $\exp(\{g,\cdot\})\mathsf{eu}$ is obtained from $\mathsf{eu}$ by adding the derivation coming from a Hamiltonian torus action on $\CC^2/\Gamma$ that commutes with the standard contracting action. This gives the conclusion of (ii).
\end{proof}

\section{$\mathbb{Q}$-factorial terminalizations}\label{subsec:Qfactorial}

Let $X$ be a normal Poisson variety with symplectic singularities. Recall that a normal variety $Y$ is $\QQ$-\emph{factorial} if every Weil divisor has a (nonzero) integer multiple which is Cartier. The following is a consequence of \cite{BCHM} (see  \cite[Prop 2.1]{LosevSRA} for a proof).

\begin{prop}\label{prop:terminalization}
There is a birational projective morphism $\rho: Y\to X$ such that
\begin{itemize}
    \item[(i)] $Y$ is an irreducible, normal, Poisson variety  (in particular, $Y$ has symplectic singularities).
    \item[(ii)] $Y$ is $\QQ$-factorial.
    \item[$(iii)$] $Y$ has terminal singularities.
\end{itemize}
\end{prop}

\begin{rmk}\label{rmk:terminalization}
Modulo (i), (iii) is equivalent to the condition that the singular locus of $Y$ is of codimension $\geq 4$, see \cite[Main Thm]{Namikawa_note}. In practice, the latter condition is often easier to check.
\end{rmk}

The map $\rho: Y\to X$ in the proposition above (or the variety $Y$ itself, if the map is understood) is called a $\QQ$-\emph{factorial terminalization}\index{$\QQ$-factorial terminalization}. If $X$ is conical, then $Y$ admits a $\CC^{\times}$-action such that $\rho$ is $\CC^{\times}$-equivariant, see \cite[A.7]{Namikawa3}. 

\begin{example}\label{Ex:Springer_resolution}
Let $\fg$ be a complex reductive Lie algebra and let $\cN \subset \fg^*$ be its nilpotent cone. By Example \ref{example:symplecticsingularity}(ii) (and the normality of $\cN)$, $X:=\cN$ is a conical symplectic singularity. For $\rho: Y \to X$ we take the Springer resolution $T^*(G/B) \to X$. In fact, this example can be generalized to all varieties of the form $X:=\operatorname{Spec}(\CC[\widetilde{\OO}])$, where $\widetilde{\OO}$ is a nilpotent cover, see Section \ref{subsec:terminalizationcover}.
\end{example}

{\begin{example}\label{Ex:minimal_Kleinian} 
Let $X=\Sigma$, a Kleinian singularity. Then its minimal resolution, $\mathfrak{S}$, is also a $\QQ$-factorial terminalization. Combining this observation with 
Proposition \ref{Prop:Kleinian_uniqueness}, we see that every $\mathbb{Q}$-factorial terminalization of a 2-dimensional conical symplectic singularity is, in fact, the minimal resolution of a Kleinian singularity. 
\end{example}}

Let $\fL_1,...\fL_t \subset X$ be the symplectic leaves of codimension 2 -- there are finitely many by Proposition \ref{Prop:fin_many_leaves_Kaledin}. For any proper birational morphism $\rho: Z \to X$ {from a normal Poisson variety $Z$}, consider the open subset
$$Z^2 := \rho^{-1}(X^{\mathrm{reg}} \cup \fL_1 \cup ... \cup \fL_t) \subset Z$$

\begin{lemma}\label{lem:codim2}
$\codim(Z-Z^2,Z) \geq 2$.
\end{lemma}

\begin{proof}
If $\rho:Z \to X$ is a $\QQ$-factorial terminalization, a proof is contained in \cite[Prop 2.14]{Losev4}. In general, let $\overline{\rho}:\overline{Z} \to Z$ be a $\QQ$-factorial terminalization of $Z$. Then $\rho \circ \overline{\rho}: \overline{Z} \to X$ is a $\QQ$-factorial terminalization of $X$. Note that $\overline{\rho}^{-1}(Z^2) \subseteq (\rho \circ \overline{\rho})^{-1}(X^{\mathrm{reg}} \cup \fL_1 \cup ... \cup \fL_t)$. Thus, 
$$\codim(Z - Z^2,Z) \geq \codim(\overline{Z} - \overline{\rho}^{-1}(Z^2),\overline{Z}) \geq \codim(\overline{Z} - (\rho \circ \overline{\rho})^{-1}(X^{\mathrm{reg}} \cup \fL_1 \cup ... \cup \fL_t),\overline{Z}) \geq 2$$
\end{proof}

\begin{definition}\label{defi:Namikawa_space}
Let $X$ be a conical symplectic singularity and $Y$ a $\QQ$-factorial terminalization of $X$.
The \emph{Namikawa space}\index{Namikawa space} associated to $X$ is the complex vector space 
$$\fP := H^2(Y^{\mathrm{reg}},\CC).$$
\end{definition}
We will see below that $\fP$ depends only on $X$ (and not on $\rho: Y \to X$), justifying the terminology. In Section \ref{subsec:quantsymplectic} we will see that $\fP$ plays a central role in the classification of Poisson deformations and filtered quantizations of $X$, and this is the main reason we consider it.

Note that $\fP$ comes equipped with a natural $\QQ$-form, namely $H^2(Y^{\mathrm{reg}},\QQ)$. So it makes sense to consider the rational and real Namikawa spaces $\fP_{\QQ}:=H^2(Y^{\mathrm{reg}},\QQ)$ and $\fP_{\RR}:=H^2(Y^{\mathrm{reg}},\RR)$. Consider the first Chern class map
$$c_1: \Pic(Y^{\mathrm{reg}}) \to H^2(Y^{\mathrm{reg}},\mathbb{Z}).$$
The following {lemma} provides an alternative characterization of $\fP_{\QQ}$.

\begin{lemma}\label{lem:c1complexification}
$c_1$ induces an isomorphism of vector spaces
\begin{equation}\label{eq:h2complexification}\Pic(Y^{\mathrm{reg}})\otimes_{\mathbb{Z}}\QQ  \xrightarrow{\sim} H^2(Y^{\mathrm{reg}},\QQ).
\end{equation}
\end{lemma}

\begin{proof}
Henceforth, we will use the superscript `an' for complex analytic objects. For example,  $\cO^{\mathrm{an}}_{Y^{\mathrm{reg}}}$ will denote the sheaf of complex analytic functions on $Y^{\mathrm{reg}}$. The idea of the proof is as follows. We will first show that the the lemma holds if we replace the algebraic Picard group $\Pic(Y^{\mathrm{reg}})$ with its analytic counterpart  $\Pic^{\mathrm{an}}(Y^{\mathrm{reg}})$ -- this is easy. Passing from $\Pic^{\mathrm{an}}(Y^{\mathrm{reg}})$ to $\Pic(Y^{\mathrm{reg}})$ requires a technical argument. Namely, we will show that every analytic line bundle on $Y^{\mathrm{reg}}$ admits a multiple with a $\CC^\times$-equivariant structure for the $\CC^\times$-action on $Y^{\mathrm{reg}}$ restricted from the contracting $\CC^\times$-action on $Y$. An equivariant analytic line bundle on $Y^{\mathrm{reg}}$ can be pushed forward to obtain a $\CC^\times$-equivariant analytic coherent sheaf on $Y$. {By appealing} to a version of the GAGA principle{, we will then see that every analytic line bundle on $Y^{reg}$ with a $\CC^\times$-equivariant structure is algebraic in a unique way. This will show the claim of the lemma.}

{\it Step 1}. Recall that $\operatorname{Pic}^{\mathrm{an}}(Y^{\mathrm{reg}})\xrightarrow{\sim} H^1(Y^{\mathrm{reg}}, \cO^{\mathrm{an},\times})$. Consider the exponential exact sequence of sheaves on $Y^{\mathrm{reg}}$
$$0\rightarrow \mathbb{Z}\rightarrow \cO^{\mathrm{an}}\rightarrow \cO^{\mathrm{an},\times}\rightarrow 0$$
and the corresponding long exact sequence in cohomology
$$H^1(Y^{\mathrm{reg}}, \cO^{\mathrm{an}})
\rightarrow \operatorname{Pic}^{\mathrm{an}}(Y^{\mathrm{reg}})\xrightarrow{c_1}H^2(Y^{\mathrm{reg}}, \mathbb{Z})\rightarrow H^2(Y^{\mathrm{reg}}, \cO^{\mathrm{an}}).$$
In Step 2, we will show that 
the first and the fourth terms vanish, and hence that $c_1$ is an isomorphism 
$\operatorname{Pic}^{\mathrm{an}}(Y^{\mathrm{reg}})\xrightarrow{\sim}H^2(Y^{\mathrm{reg}},\mathbb{Z})$.

{\it Step 2}.  Next we prove that 
\begin{equation}\label{eq:cohom_vanish}
H^j(Y^{\mathrm{reg}},\mathcal{O}^{\mathrm{an}})=0 \text{ for }  j=1,2.
\end{equation}
The variety $Y$ is singular symplectic, hence Cohen-Macaulay (singular symplectic implies Gorenstein by {Lemma \ref{lem:rationalsingularities}}; it is a standard fact that Gorenstein implies Cohen-Macaulay). This in particular implies that the stalks of $\mathcal{O}^{\mathrm{an}}_Y$ are Cohen-Macaulay.  {Indeed, these stalks are Noetherian. Now we apply the observation that a local Noetherian ring is Cohen-Macaulay if and only if its completion at the maximal ideal is Cohen-Macaulay, see \cite[Lemma 15.43.3]{stacks-project}.}

It follows that the cohomology with support $H^i_{Y^{\mathrm{reg}}}(Y,\mathcal{O}^{\mathrm{an}})$ is $0$ for $i\leqslant 3$. Thus
\begin{equation}\label{eq:cohom_coincidence}
H^j(Y^{\mathrm{reg}},\mathcal{O}^{\mathrm{an}})\xrightarrow{\sim} H^j(Y,\mathcal{O}^{\mathrm{an}}), \text{ for } j\leqslant 2. 
\end{equation}
Since $Y$ is singular symplectic, $\mathcal{O}^{\mathrm{an}}$ coincides with the dualizing sheaf $K^{\mathrm{an}}$. By the Grauert-Riemenschneider theorem, 
$R^i\rho_* K^{\mathrm{an}}=0$ for $i>0$. Since $X$ is affine and hence Stein, this implies that $H^k(Y, K^{\mathrm{an}})=0$ for all $k>0$. This, together with the isomorphism $\mathcal{O}^{\mathrm{an}}\simeq K^{\mathrm{an}}$ and  
(\ref{eq:cohom_coincidence}), implies
(\ref{eq:cohom_vanish}) and hence the isomorphism
$c_1:\operatorname{Pic}^{\mathrm{an}}(Y^{\mathrm{reg}})\xrightarrow{\sim}H^2(Y^{\mathrm{reg}},\mathbb{Z})$.

{\it Step 3}. Next we show that every analytic line bundle on $Y^{\mathrm{reg}}$ admits a multiple with a $\CC^\times$-equivariant structure. 
Pick a sufficiently large integer $n$ and consider the complex analytic manifold $E^n:=\CC^n- \{0\}$
(to be understood as an approximation of $E(\CC^\times)$). This manifold has a free action of $\CC^\times$ with quotient $\mathbb{P}^{n-1}$. 

Consider the product $\mathbb{Y}:=Y^{\mathrm{reg}}\times E^n$. This too admits a free action of $\CC^\times$, i.e. the diagonal one, and the quotient variety $\mathbb{Y}/\CC^\times$ is a fiber bundle over $\mathbb{P}^{n-1}$ with fiber $Y^{\mathrm{reg}}$.

{\it Step 4}. Let $\pi$ denote the projection $\mathbb{Y}\twoheadrightarrow Y^{\mathrm{reg}}$. We claim that every analytic line bundle $\mathcal{L}$ on $\mathbb{Y}$ is isomorphic to the pullback under $\pi$ of a line bundle on $Y^{\mathrm{reg}}$. Since $H^i(Y^{reg},\mathcal{O}^{\mathrm{an}})=0$ for $i=1,2$ (see Step 2) and $H^i(E^n, \mathcal{O}^{\mathrm{an}})=0$ for $0<i<n-1$, we have 
$H^i(\mathbb{Y},\mathcal{O}^{\mathrm{an}})=0$ for $i=1,2$. Hence, the exponential sequence for $\mathbb{Y}$ induces an isomorphism 
$\operatorname{Pic}^{\mathrm{an}}(\mathbb{Y})
\xrightarrow{\sim} H^2(\mathbb{Y},\mathbb{Z})$. By the K\"{u}nneth formula, $\pi^*$ induces an isomorphism 
$H^2(Y^{\mathrm{reg}},\mathbb{Z})\xrightarrow{\sim} H^2(\mathbb{Y},\mathbb{Z})$. Since the exponential sequence is functorial with respect to pullbacks, we see that $\pi^*$ is an isomorphism 
\begin{equation}\label{eq:Pic_iso}
\Pic^{\mathrm{an}}(Y^{\mathrm{reg}})\simeq \Pic^{\mathrm{an}}(\mathbb{Y}).
\end{equation} 
An inverse for (\ref{eq:Pic_iso}) can be
constructed as follows. Let $\iota$ denote the inclusion $\mathbb{Y}\hookrightarrow Y^{\mathrm{reg}}\times \CC^n$. Thanks to (\ref{eq:Pic_iso}), we see that 
$\iota_*\mathcal{L}$ is a line bundle on $Y^{\mathrm{reg}}\times \C^n$. Restricting $\iota_*\mathcal{L}$ to $Y^{\mathrm{reg}}\times \{0\}$ defines an inverse for $\pi^*$.

{\it Step 5}. We have  $H^i(\mathbb{Y}/\CC^\times,
\mathcal{O}^{\mathrm{an}})=
H^i(\mathbb{Y},
\mathcal{O}^{\mathrm{an}})^{\CC^\times}$. By Step 4, the right hand side vanishes. 
Therefore   $c_1:\operatorname{Pic}^{\mathrm{an}}(\mathbb{Y}/\CC^\times)\rightarrow H^2(\mathbb{Y}/\mathbb{C}^\times,\mathbb{Z})$ is an isomorphism. Now we compute the base change of the target to $\CC$.

By \cite[Lemma 2.15]{Losev4} $H^1(Y^{\mathrm{reg}}, \CC)=0$. Also $H^1(\mathbb{P}^{n-1},\CC)=0$. Using the Serre spectral sequence for the fiber bundle $\mathbb{Y}/\CC^\times\rightarrow \mathbb{P}^{n-1}$, we see that 
$H^2(\mathbb{Y}/\CC^\times,\CC)\simeq 
H^2(Y^{\mathrm{reg}},\CC)\oplus H^2(\mathbb{P}^{n-1},\CC)$. Equivalently, the following sequence is exact
\begin{equation}\label{eq:cohomology_exact}
0\rightarrow \C=H^2(\mathbb{P}^{n-1},\CC)\rightarrow H^2(\mathbb{Y}/\CC^\times,\C)\rightarrow H^2(Y^{\mathrm{reg}},\C)\rightarrow 0.
\end{equation}
Here the first map is the pullback under $\mathbb{Y}/\CC^\times\rightarrow \mathbb{P}^{n-1}$ and the second is the pullback under the fiber inclusion $Y^{\mathrm{reg}}\hookrightarrow \mathbb{Y}/\CC^\times$. 

{\it Step 6}. Consider the equivariant Picard group $\Pic^{\mathrm{an}}_{\CC^\times}(Y^{\mathrm{reg}})$. Forgetting equivariance, we get a group homomorphism
$$\Pic^{\mathrm{an}}_{\CC^\times}(Y^{\mathrm{reg}}) \to \Pic^{\mathrm{an}}(Y^{\mathrm{reg}})$$
We claim that the cokernel is torsion.

Since the first Chern class map commutes with pullbacks from (\ref{eq:cohomology_exact}), the cokernel of the homomorphism
\begin{equation}\label{eq:Pic_restr_map}
\Pic^{\mathrm{an}}(\mathbb{Y}/\CC^\times)\rightarrow \Pic^{\mathrm{an}}(Y^{\mathrm{reg}})
\end{equation}
is torsion. Furthermore $\Pic^{\mathrm{an}}(\mathbb{Y}/\CC^\times)\xrightarrow{\sim} \Pic^{\mathrm{an}}_{\CC^\times}(\mathbb{Y})$. 
By Step 4, $\iota_*\mathcal{L}$ is a line bundle on $Y^{\mathrm{reg}}\times \C^n$ for all line bundles $\mathcal{L}$ on $\mathbb{Y}$. It follows that there is an isomorphism $\Pic^{\mathrm{an}}_{\CC^\times}(\mathbb{Y})\xrightarrow{\sim} \Pic^{\mathrm{an}}_{\CC^\times}(Y^{\mathrm{reg}}\times \CC^n)$ given by push-forward. Under this identification, (\ref{eq:Pic_restr_map}) becomes the pullback map under restriction to $Y^{\mathrm{reg}}\times \{0\}$ and forgetting equivariance. In particular, the cokernel of  
$\Pic^{\mathrm{an}}_{\CC^\times}(Y^{\mathrm{reg}})\rightarrow \Pic^{\mathrm{an}}(Y^{\mathrm{reg}})$ is torsion. Equivalently, every line bundle on $Y^{\mathrm{reg}}$ has a multiple which admits a $\CC^\times$-equivariant structure. 

{\it Step 7}. 
Let $\iota$ denote the inclusion $Y^{\mathrm{reg}}\hookrightarrow Y$. Since $\operatorname{codim}_Y (Y^{\mathrm{sing}})\geqslant 3$, the push-forward $\mathcal{F}:=\iota_* \mathcal{L}$ is an analytic coherent sheaf, see \cite[Theorem 5]{Siu}. It comes with a natural $\CC^\times$-equivariant structure. As usual, let $\rho$ denote the projective morphism $Y\rightarrow X$. In Steps 8-9, we will show that the restriction of $\mathcal{F}$ to $Y- \rho^{-1}(0)$ is the analytification of a {unique $\mathbb{C}^\times$-equivariant} algebraic coherent sheaf. We will now explain how this implies the lemma. 

If $\dim Y=2$ (equivalently, {by Example \ref{Ex:minimal_Kleinian}}, if $\rho:Y \to X$ is the minimal resolution of a Kleinian singularity), the claim of the lemma follows from Remark \ref{rmk:PforKleinian}. Thus, we can assume $\dim Y\geqslant 4$. In this case $\rho^{-1}(0)$ has codimension at least $2$, see Lemma \ref{lem:codim2}. Let $\mathcal{F}_1$ denote a $\CC^\times$-equivariant algebraic coherent sheaf on $Y- \rho^{-1}(0)$ with $\mathcal{F}_1^{\mathrm{an}}\simeq \mathcal{F}$.
Let $\mathcal{L}_1$ denote the restriction of $\mathcal{F}_1$ to $Y^{\mathrm{reg}}- \rho^{-1}(0)$ so that the restriction of $\mathcal{L}$ to $Y^{\mathrm{reg}}- \rho^{-1}(0)$ coincides with the analytification of $\mathcal{L}_1$.
Write $\iota_1$ for the inclusion $Y^{\mathrm{reg}}- \rho^{-1}(0)\hookrightarrow Y^{\mathrm{reg}}$. By Step 2 of the proof of \cite[Proposition 3.2]{BoixedaAlvarezLosev}, $\mathcal{L}$, which coincides with the analytic pushforward of its restriction to $Y^{\mathrm{reg}}- \rho^{-1}(0)$, {is isomorphic to} the analytification of $\iota_{1*}\mathcal{L}_1$. It follows that $\iota_{1*}\mathcal{L}_1$ is a line bundle. So $\mathcal{L}$ is the analytification of a {unique} algebraic $\mathbb{C}^\times$-equivariant line bundle. {It follows that $\operatorname{Pic}_{\CC^\times}(Y^{reg})\xrightarrow{\sim} \operatorname{Pic}^{an}_{\CC^\times}(Y^{reg})$. Since every algebraic line bundle on $Y^{\operatorname{reg}}$ admits a $\CC^\times$-equivariant structure, the analytification map $\operatorname{Pic}(Y^{reg})\rightarrow \operatorname{Pic}^{an}(Y^{reg})$ is injective. And Step 6 combined with the isomorphism $\operatorname{Pic}_{\CC^\times}(Y^{reg})\xrightarrow{\sim} \operatorname{Pic}^{an}_{\CC^\times}(Y^{reg})$ shows that the cokernel of $\operatorname{Pic}(Y^{reg})\rightarrow \operatorname{Pic}^{an}(Y^{reg})$ is torsion implying the claim of the lemma.} 

It remains to show that {every $\CC^\times$-equivariant analytic coherent sheaf} $\mathcal{F}$ on $Y- \rho^{-1}(0)$ is the analytification of a {unique $\mathbb{C}^\times$-equivariant} algebraic coherent sheaf. This will be accomplished using the GAGA principle in Steps 8-9. 

{\it Step 8}. Consider a more general situation: let $\rho:Y\rightarrow X$ a projective morphism onto an affine variety $X$. Suppose that $\rho^*:\CC[X]\rightarrow\CC[Y]$ is an isomorphism. Further, assume that $X,Y$ are equipped with $\CC^{\times}$-actions such that $\CC[X]$ is positively graded and $\rho$ is $\CC^{\times}$-equivariant. In particular, there is a point $0\in X$. We claim that every $\CC^\times$-equivariant analytic coherent sheaf $\mathcal{F}$ on $Y- \rho^{-1}(0)$ is the analytification of
a {unique} algebraic $\CC^\times$-equivariant coherent sheaf. 

Assume first that $\CC[X]$ is generated by degree $1$ elements. 
Because of this, the action of $\CC^\times$ on $X- \{0\}$ is free with quotient $\operatorname{Proj}(\CC[X])$. We can form the GIT quotient of $Y$ by the action of $\CC^\times$ 
(with the trivial line bundle). This quotient is $(Y- \rho^{-1}(0))/\C^\times$. The latter variety is projective over $(X- \{0\})/\CC^\times$, hence is projective. Form the analytic coherent sheaf $\underline{\mathcal{F}}$ on
$(Y- \rho^{-1}(0))/\C^\times$ obtained from $\mathcal{F}$ by equivariant descent. By the GAGA principle, it comes from a unique algebraic coherent sheaf
on $(Y- \rho^{-1}(0))/\CC^\times$, to be denoted by $\underline{\mathcal{F}}_1$. Let $\mathcal{F}_1$ denote the pullback of $\underline{\mathcal{F}}_1$ to 
$Y- \rho^{-1}(0)$. The analytification of $\mathcal{F}_1$ is $\mathcal{F}$. The uniqueness (of $\underline{\mathcal{F}}_1$ and hence of $\mathcal{F}_1$) again follows from the GAGA principle: the analytification functor on a projective variety is a category equivalence.  

{\it Step 9}.
Finally, we will prove the {the existence and uniqueness of the algebraization of $\mathcal{F}$} in the general case, i.e. when $\CC[X]$ need not be generated by elements of degree $1$.
There is $e\geqslant 1$ such that the subalgebra $\CC[X]_{(e)}:=\bigoplus_{i\geqslant 0}\CC[X]_{ie}$ (where $\CC[X]_{ie}$
denotes the graded component of degree $ie$) is generated by elements of degree $e$. Let $\Gamma\subset \CC^\times$ be the subgroup of $e$th roots of unity. Consider the quotients $X/\Gamma,Y/\Gamma$. Let $\rho_{(e)}:Y/\Gamma\rightarrow X/\Gamma$
be the natural morphism and $\varpi:Y\rightarrow Y/\Gamma$ the quotient morphism. Consider the sheaf $\varpi_* \mathcal{F}$. It splits as a direct sum $\bigoplus \mathcal{F}(i)$, where $\mathcal{F}(i)$ is the eigensheaf for the action of $\Gamma$ with eigenvalue $t\mapsto t^i, t\in \Gamma, i=0,\ldots,e-1$. The sheaf $\mathcal{F}(i)$ is $\CC^\times/\Gamma$-equivariant---we twist the $\CC^\times$-action on  
$\mathcal{F}(i)$ by $t\mapsto t^{-i}$. The action of $\CC^\times/\Gamma$ on $X/\Gamma$ satisfies the conditions of Step 8. So we can find a unique $\CC^\times/\Gamma$-equivariant algebraic coherent sheaf $\mathcal{F}_1(i)$ on $Y/\Gamma- \rho_{(e)}^{-1}(0)$
whose analytification is $\mathcal{F}(i)$. Form the sheaf $\mathcal{F}_1=\bigoplus_{i=0}^{e-1}\mathcal{F}_1(i)$ on $Y/\Gamma- \rho_{(e)}^{-1}(0)$. Note that the analytification of $\mathcal{F}_1$ is $\varpi_*\mathcal{F}$, {and $\mathcal{F}_1$ is uniquely characterized by this property}. By the functoriality of the construction, 
$\mathcal{F}_1$ is a sheaf of $\varpi_*\mathcal{O}_{Y-\rho^{-1}(0)}$-modules. So we can regard it as a sheaf on $Y- \rho^{-1}(0)$. The analytification is
still the restriction of $\mathcal{F}$
to $Y- \rho^{-1}(0)$.
\end{proof}

If we choose a different $\QQ$-terminalization $Y_1 \to X$, there is a canonical isomorphism $\Pic(Y_1^{\mathrm{reg}}) \simeq \Pic(Y^{\mathrm{reg}})$, see \cite[Prop 2.18]{BPW}. In particular, in view of Lemma \ref{lem:c1complexification}, the Namikawa space $\fP_{\QQ}$ is independent of the choice of $\QQ$-factorial terminalization $\rho: Y \to X$.

\section{Structure of Namikawa space}\label{subsec:structurenamikawa}

Let $X$ be a conical symplectic singularity, and choose a $\QQ$-factorial terminalization $\rho: Y \to X$ as in Section \ref{subsec:Qfactorial}. Following \cite{Namikawa} and \cite{Losev4}, we will provide a description of $\fP$ in terms of the geometry of $X$.

Let $\fL_1, ..., \fL_t \subset X$ be the symplectic leaves of codimension 2. For each $\fL_k$, the formal slice to $\fL_k \subset X$ is {a 2-dimensional symplectic singularity and so, by an argument analogous to Steps 1 and 2 of the proof of Proposition \ref{Prop:Kleinian_uniqueness}, is} identified {as a Poisson formal scheme} with the formal neighborhood at $0$ in a Kleinian singularity $\Sigma_k = \CC^2/\Gamma_k$. {More precisely, the following holds. Choose a point $x\in \fL_k$ for some $k=1,\ldots,t$. Set $V:=T_x\mathfrak{L}_k$. We write $\CC[X]^\wedge$ for the completion of $\C[X]$ with respect to the maximal ideal of $x$. Similarly, we write $\CC[\Sigma_k]^\wedge$ for the completion of $\CC[\Sigma_k]$ with respect to the maximal ideal of $0$. Then we have a Poisson algebra isomorphism
\begin{equation}\label{eq:decompositionleaf_classical}\CC[X]^{\wedge}\simeq \CC[V]^{\wedge} \ \widehat{\otimes} \  \CC[\Sigma_k]^\wedge,\end{equation}
cf. \cite[Theorem 2.3]{Kaledin2006}. Here, and elsewhere, we write $\widehat{\otimes}$ for the completed tensor product.
 }

We will call $\Sigma_k$ the `singularity' of the leaf $\fL_k \subset X$. Let $\fg_k$ be the corresponding simple complex Lie algebra of type ADE (cf. Theorem \ref{thm:McKay}). Fix a Cartan subalgebra $\fh_k \subset \fg_k$, and let $\Delta_k \subset \fh_k^*$, $\Lambda_k \subset \fh_k^*$, and $W_k$ be the root system, weight lattice, and Weyl group, respectively. Let $\rho_k:\mathfrak{S}_k\rightarrow \Sigma_k$ be the minimal resolution. Recall, see Remark \ref{rmk:PforKleinian}, that $\rho_k^{-1}(0)$ is homotopically equivalent to $\mathfrak{S}_k$ due to the contracting $\mathbb{C}^\times$-action on the latter. In particular, there is a natural identification $H^2(\mathfrak{S}_k,\ZZ) \simeq H^2(\rho_k^{-1}(0),\ZZ)$. Now choose a point $x\in \mathfrak{L}_k$. Since $\rho: Y\rightarrow X$ is a $\mathbb{Q}$-factorial terminalization, it gives a symplectic resolution over codimension $2$ leaves. Hence there is an identification $\rho^{-1}(x)\simeq \rho_k^{-1}(0)$. The fundamental group $\pi_1(\fL_k)$ acts by monodromy on the irreducible components $\{C_i\}$ of $\rho_k^{-1}(0)$, and this action preserves the intersection graph of $\mathfrak{S}_k$. Thus, $\pi_1(\fL_k)$ acts by diagram automorphisms on $\Delta_k$---acting on the simple roots via the bijection (\ref{eq:bijectionsMcKay})---and hence also on $\Delta_k$, $\Lambda_k$, $\fh_{k}$, and $W_k$. The \emph{partial Namikawa space}\index{Namikawa space!partial} for $\fL_k$ is the space of monodromy invariants
$$\fP_{k} := (\fh_{k}^*)^{\pi_1(\fL_k)}.$$
Since $\rho: Y \to X$ is a $\QQ$-factorial terminalization, $\rho^{-1}(\fL_k) \subset Y^{\mathrm{reg}}$. The embedding  $\rho_k^{-1}(0) \hookrightarrow Y^{\mathrm{reg}}$ gives rise to a map on cohomology
$$\fP := H^2(Y^{\mathrm{reg}},\CC) \to  H^2(\mathfrak{S}_k,\CC)^{\pi_1(\fL_k)} \simeq \fP_{k}.$$
This map is defined over $\QQ$. Also, define
$$\fP_{0} := H^2(X^{\mathrm{reg}},\CC).$$
The embedding $X^{\mathrm{reg}} \hookrightarrow Y^{\mathrm{reg}}$ gives rise to a map
$$\fP = H^2(Y^{\mathrm{reg}},\CC) \to H^2(X^{\mathrm{reg}},\CC) = \fP_{0},$$
also defined over $\QQ$.

\begin{prop}[\cite{Losev4}, Lem 2.8]\label{prop:partialdecomp}
The maps $\fP \rightarrow \fP_{k}$ defined above induce a linear isomorphism
\begin{equation}\label{eq:partialdecomp}\fP \simeq\bigoplus_{k=0}^t \fP_{k}, \qquad \lambda \mapsto (\lambda_0,\lambda_1,...,\lambda_t).\end{equation}
\end{prop}


We will now define the \emph{Namikawa Weyl group}\index{Namikawa Weyl group} associated to $X$ . For each codimension 2 leaf $\fL_k \subset X$, consider the subgroup of monodromy invariants $W_k^{\pi_1(\fL_k)} \subset W_k$. Note that there is a natural action of $W_k^{\pi_1(\fL_k)}$ on $\fP_{k} = (\fh_{k}^*)^{\pi_1(\fL_k)}$. The Namikawa Weyl group associated to $X$ is the product 
$$W := \prod_{k=1}^t W_k^{\pi_1(\fL_k)}.$$
$W$ acts on $\fP$ via the isomorphism (\ref{eq:partialdecomp}) (the action on $\fP_{0}$ is trivial).

\section{Finite covers of conical symplectic singularities}\label{subsec:finitecovers}

Let $X$ be a conical symplectic singularity. In this section, we will define the notion of a \emph{finite cover} of $X$. Let $p': \widetilde{X}' \to X^{\mathrm{reg}}$ be a finite \'{e}tale cover of the regular locus $X^{\mathrm{reg}} \subset X$. Rescaling if necessary, we can arrange so that the $\CC^{\times}$-action on $X^{\mathrm{reg}}$ lifts to $\widetilde{X}'$. Consider the composition $\widetilde{X}' \overset{p'}{\to} X^{\mathrm{reg}} \hookrightarrow X$ and its Stein factorization
\begin{center}
\begin{tikzcd}
\widetilde{X}' \ar[r,hookrightarrow] \ar[d,"p'"] & \widetilde{X} \ar[d, "p"]\\
X^{\mathrm{reg}} \ar[r,hookrightarrow] & X
\end{tikzcd}
\end{center}
Note that $\widetilde{X}$ is affine and $\widetilde{X}'$ embeds into $\widetilde{X}$ as an open subvariety. Since $\operatorname{codim}(X^{\mathrm{sing}},X)\geqslant 2$ and $p:\widetilde{X} \to X$ is finite, we have that $\codim(\widetilde{X}-\widetilde{X}',\widetilde{X}) \geq 2$. Thus the algebra $\CC[\widetilde{X}']$ is finitely generated and
$ \widetilde{X}= \Spec(\CC[\widetilde{X}'])$. In particular, the $\CC^{\times}$-action on $\widetilde{X}'$ extends to $\widetilde{X}$. In fact, $\widetilde{X}$ is a conical symplectic singularity, see \cite[Lemma 2.5]{LosevHC}. A map $p:\widetilde{X} \to X$ obtained in this fashion is called a \emph{finite cover} of $X$. We say that $p$ is \emph{Galois} if its restriction to $\widetilde{X}'$ is Galois. 

\begin{lemma}\label{lem:leaftoleaf}
Let $p: \widetilde{X} \to X$ be a finite cover and $\fL' \subset \widetilde{X}$ a symplectic leaf of codimension $k$. Then there is a symplectic leaf $\fL \subset X$ of codimension $k$ such that $p(\overline{\fL}') = \overline{\fL}$.
\end{lemma}

\begin{proof}
The assignment $\fL \mapsto I(\overline{\fL})$ defines a bijection between the set of symplectic leaves in $X$ (resp. $\widetilde{X}$) and the set of prime Poisson ideals in $\CC[X]$ (resp. $\CC[\widetilde{X}]$), 
Lemma \ref{lem:fin_many_leaves_equiv}.  If $I$ is a prime Poisson ideal in $\CC[\widetilde{X}]$, then $I \cap \CC[X]$ is a prime Poisson ideal in $\CC[X]$. Thus, $p(\overline{\fL}')$ is the closure of a symplectic leaf in $X$, say $\fL$. Since $p$ is finite, it preserves codimension. Thus, $\fL \subset X$ is of codimension $k$.
\end{proof}

Next, we relate the dimension 2 singularities in $X$ and $\widetilde{X}$. Let $\fL_k \subset X$ be a codimension 2 leaf and $\Sigma_k = \CC^2/\Gamma_k$ the corresponding Kleinian singularity. Let $\Gamma$ denote the algebraic fundamantal group $\pi^{\mathrm{alg}}_1(X^{\mathrm{reg}})$ of $X^{\mathrm{reg}} \subset X$. By \cite[Cor 5.2]{SGA}, $\Gamma$ is the profinite completion of the usual fundamental group $\pi_1(X^{\mathrm{reg}})$. And, by \cite[Main Thm]{Namikawa_fundgroups}, $\Gamma$ is finite. Thus, $\Gamma$ is the maximal finite quotient of $\pi_1(X^{\mathrm{reg}})$. 

{From (\ref{eq:decompositionleaf_classical}) we get the} morphism $\Sigma_k^{\wedge}\hookrightarrow X^\wedge \to X$.
 Since $\Sigma_k^{\wedge, \mathrm{reg}}$ maps to $X^{\mathrm{reg}}$, there is an induced homomorphism
\begin{equation}\label{eq:defofphik}
\phi_k: \Gamma_k = \pi_1^{\mathrm{alg}}(\Sigma_k^{\wedge, \mathrm{reg}}) \to \pi_1^{\mathrm{alg}}(X^{\mathrm{reg}}) = \Gamma.\end{equation}
Choose a point $x \in \fL_k$ and let $p^{-1}(x) = \{\widetilde{x}_j\}$. Write $X^{\wedge x}$ (resp. $\widetilde{X}^{\wedge \widetilde{x}_j}$) for the completion of $X$ (resp. $\widetilde{X}$) at $x$ (resp. $\widetilde{x}_j$), i.e. $X^{\wedge_x}=\operatorname{Spec}(\mathbb{C}[X]^{\wedge_x})$. There is a Cartesian diagram of schemes
\begin{equation}\label{diag:fiberproduct}
\begin{tikzcd}
 \bigsqcup_j \widetilde{X}^{\wedge \widetilde{x}_j} \ar[d] \ar[r] & \widetilde{X} \ar[d]  \\
 X^{\wedge  x} \ar[r] &  X
\end{tikzcd}
\end{equation}
Now suppose $\widetilde{X} \to X$ is Galois with Galois group $\Pi$ so that $\Pi$ acts faithfully on $\widetilde{X}$ and $\widetilde{X}/\Pi \xrightarrow{\sim} X$. Choose a preimage $\widetilde{x}_j$, and let $\fL'_k\subset \widetilde{X}$ be the connected component of $\widetilde{x}_j$ in $p^{-1}(\fL_k)$. The formal slice to $\fL'_k$ at $\widetilde{x}_j$ is isomorphic to the completion at $0$ of a Kleinian singularity $\Sigma_k' := \CC^2/\Gamma_k'$ (with, possibly, $\Gamma_k'=1$).   Consider the induced homomorphism 
$$\Gamma_k \overset{\phi_k}{\to} \Gamma \twoheadrightarrow \Pi.$$
 By construction, $\Gamma_k'$ is the kernel of this map. In particular, $\Gamma_k'$ is normal in $\Gamma_k$. 

Let $\Gamma'$ denote the algebraic fundamental group of $\widetilde{X}$ (since the cover $\widetilde{X} \to X^{reg}$ is Galois, $\Gamma'$ is the kernel of the epimorphism $\Gamma \twoheadrightarrow \Pi$). Similarly to (\ref{eq:defofphik}), there are group homomorphisms $\phi_k': \Gamma_k' \to \Gamma'$. The diagram of schemes (\ref{diag:fiberproduct}) gives rise to a diagram of groups
\begin{equation}\label{diag:galoiscovering}
\begin{tikzcd}
    \Gamma_k' \ar[d, hookrightarrow] \ar[r,"\phi'_k"] & \Gamma' \ar[d, hookrightarrow]  \\
 \Gamma_k \ar[r,"\phi_k"] &  \Gamma
\end{tikzcd}
\end{equation}
\section{Deformations and quantizations of symplectic singularities}\label{subsec:quantsymplectic}

Let $X$ be a conical symplectic singularity and choose a $\QQ$-factorial terminalization $\rho: Y \to X$ as in Section \ref{subsec:Qfactorial}. Recall the  Namikawa space $\fP = H^2(Y^{\mathrm{reg}},\CC)$. Below in this section we will see that this vector space parameterizes both Poisson deformations and filtered quantizations of $Y$. 

The following proposition was obtained in \cite{Losev4} but is a direct consequence of results of Namikawa (\cite{Namikawa3}).

\begin{prop}[\cite{Losev4}, Proposition 2.6]\label{prop:universaldef}
There is a graded Poisson scheme $Y_{\mathrm{univ}}$ over $\fP$ such that
\begin{itemize}
    \item[(i)] The fiber of $Y_{\mathrm{univ}}$ over $0$ is $Y$.
    \item[(ii)] Suppose $B$ is a conical affine scheme, $Z$ is a graded Poisson scheme over $B$, and (i) holds for the morphism $Z \to B$. Then there is a unique  $\CC^{\times}$-equivariant morphism $B \to \fP$ with the following property: there is an isomorphism of graded Poisson schemes $Z \simeq Y_{\mathrm{univ}}\times_\fP B$ over $B$ which is the identity over $0 \in B$.
\end{itemize}
\end{prop}

The scheme $Y_{\mathrm{univ}}$ in Proposition \ref{prop:universaldef} is called the \emph{universal graded Poisson deformation} of $Y$. 

\begin{rmk}\label{rmk:period_map}
Below, we will need a geometric interpretation of the unique morphism $B\rightarrow \fP$. We provide one in the special case when $B$ is a vector space with $\CC^\times$-action given by $t\mapsto t^{-d}$, as on $\fP$ so that the morphism $B\rightarrow \fP$
is linear. {Let $B^\wedge$ (resp. $\fP^\wedge$) denote the formal neighborhood of $0$ in $B$ (resp. $\fP$). Let $Z^\wedge$ (resp. $Y_{\mathrm{univ}}^\wedge$) denote the formal neighborhood of $Y^{\mathrm{reg}}$ in $Z$ (resp. $Y_{\mathrm{univ}}$).} These  are formal deformations of $Y^{\mathrm{reg}}$ over $B^\wedge,\fP^\wedge$. The de Rham cohomology group $H^2_{dR}(Z^\wedge/B^\wedge)$ decomposes as a tensor product $H^2_{dR}(Y^{\mathrm{reg}})\otimes \CC[B^\wedge]$ thanks to the Gauss-Manin connection. In particular, the class of the relative symplectic form on $Z^\wedge/B^\wedge$ can be viewed as 
a map $B^\wedge\rightarrow \fP^\wedge$.
Because of the $\CC^\times$-actions, this map comes from a linear map $B\rightarrow \fP$. It was shown in \cite[Theorem 3.6]{Kaledin-Verbitsky} that the map $B^\wedge\rightarrow \fP^\wedge$ is the unique such map inducing an isomorphism $Z^\wedge\xrightarrow{\sim} B^\wedge\times_{\fP^\wedge}Y_{\mathrm{univ}}^\wedge$ of formal Poisson deformations of $Y^{\mathrm{reg}}$.  A connection between the Namikawa and Kaledin-Verbitsky approaches is explained in \cite[Remark 2.7]{Losev4}. In particular, the map $B\rightarrow \fP$ described in this remark is the same as the one in Proposition \ref{prop:universaldef}.
\end{rmk}

Next, we state an analog of Proposition \ref{prop:universaldef} for filtered quantizations. For any graded smooth symplectic variety $V$, there is a (non-commutative) \emph{period map}\index{period map}
$$\mathrm{Per}: \mathrm{Quant}(V) \to H^2(V,\CC),$$
see \cite[Section 4]{BK}, \cite[Section 2.3]{Losev_isofquant}.

\begin{prop}[\cite{Losev4}, Proposition 3.1(1)]\label{prop:P=Quant(Y)}
The maps 
$$\mathrm{Quant}(Y) \overset{|_{Y^{\mathrm{reg}}}}{\to} \mathrm{Quant}(Y^{\mathrm{reg}}) \overset{\mathrm{Per}}{\to} H^2(Y^{\mathrm{reg}},\CC) = \fP$$
are bijections.
\end{prop}

For $\lambda \in \fP$, write  $\mathcal{D}_{\lambda}$ for the corresponding  filtered quantization of $Y$. There is also 
a universal version of $\mathcal{D}_\lambda$.

\begin{prop}[\cite{Losev4}, Cor 3.2(2)]\label{prop:universal for Q-term}
There is a sheaf $\calD^{Y,\mathrm{univ}}$ of associative $\CC[\fP]$-algebras on $Y$ such that, for each $\lambda \in \fP$, the filtered quantization $\calD_\lambda$ of $Y$ is obtained  as
$\calD^{Y_{\mathrm{univ}}}\otimes_{\CC[\fP]}\C$, where the homomorphism $\CC[\fP]\rightarrow \CC$ corresponds to $\lambda$.
\end{prop}

Next, we explain the relation between Poisson deformations (resp. quantizations) of $Y$ and Poisson deformations (resp. quantizations) of $X$. For $\lambda \in \fP$, set $ \cA_{\lambda} := \Gamma(Y,\mathcal{D}_{\lambda})$. Let $Y_\lambda$ denote the fiber of $Y_{\mathrm{univ}}\rightarrow \fP$ over $\lambda$ and let $\cA^0_\lambda=\CC[Y_\lambda]$.

\begin{theorem}\label{thm:defsofsymplectic}
The following are true:

\begin{itemize}
\item[(i)] For every $\lambda \in \fP$, the algebra $\cA_\lambda^0$ is a filtered Poisson deformation of $\CC[X]$.
    \item[(ii)] Every filtered Poisson deformation of $\CC[X]$ is isomorphic to $\cA^0_{\lambda}$ for some $\lambda \in \fP$.
    \item[(iii)] For every $\lambda, \lambda' \in \fP$, we have $\cA^0_{\lambda} \simeq \cA^0_{\lambda'}$ if and only if $\lambda' \in W\cdot \lambda$.
\end{itemize}
Hence, the map $\lambda \mapsto \cA^0_{\lambda}$ induces a bijection
$$\fP/W \simeq \mathrm{PDef}(X), \qquad W \cdot \lambda \mapsto \cA^0_{\lambda}.$$
\end{theorem}
\begin{proof}
These results are essentially due to Namikawa, \cite{NamikawaQ2}. (i) is a consequence of 
\cite[Prop 2.9(1)]{Losev4}. (ii) and (iii) follow from \cite[Prop 2.12, Cor 2.13]{Losev4}. 
\end{proof}

\begin{theorem}[\cite{Losev4}, Prop 3.3, Thm 3.4]\label{thm:quantsofsymplectic}
Direct analogs of  (i)-(iii) of Theorem \ref{thm:defsofsymplectic} hold for filtered quantizations. In particular, there is a bijection
$$\fP/W \simeq \mathrm{Quant}(X), \qquad W\cdot \lambda \mapsto \cA_{\lambda}.$$
\end{theorem}

Next we discuss the universal graded deformation of $X$. The algebra $\CC[Y_{\mathrm{univ}}]$ is a graded deformation of $\CC[Y]=\CC[X]$ over $\CC[\fP]$. The algebra $\CC[Y_{\mathrm{univ}}]$ carries a $W$-action which is the identity on $\CC[Y]$ and makes the homomorphism $\CC[\fP]\rightarrow \CC[Y_{\mathrm{univ}}]$ equivariant.
These two claims constitute 
\cite[Prop 2.9]{Losev4}. Similarly to the case of the $\QQ$-factorial terminalizations, all Poisson deformations $\cA^0_{\lambda}$ can be recovered as the fibers of a certain Poisson scheme over $\fP/W$.

\begin{prop}[\cite{Losev4}, Prop 2.12, Cor 2.13]\label{prop:universaldefsymplectic}
There is a graded affine Poisson scheme $X_{\mathrm{univ}}$ over $\fP/W$ such that:
    \begin{itemize}
        \item[(i)] The fiber of $X_{\mathrm{univ}}$ over $0$ is $X$.
        \item[(ii)] Suppose $B$ is a conical affine scheme, $Z$ is a graded Poisson scheme over $B$, and (i) holds for the morphism $Z \to B$. Then there is a unique $\CC^{\times}$-equivariant morphism $B \to \fP/W$ with the following property: there is an isomorphism of graded Poisson schemes $Z \xrightarrow{\sim} X_{\mathrm{univ}}\times_{\fP/W} B$  over $B$ which is the identity over $0 \in B$.
        \item[(iii)] We have $\CC[X_{\mathrm{univ}}]=\CC[Y_{\mathrm{univ}}]^W$.
    \end{itemize} 
\end{prop}

The scheme $X_{\mathrm{univ}}$ in Proposition \ref{prop:universaldefsymplectic} is called the \emph{universal Poisson deformation}\index{Poisson deformation!universal} of $X$. Note that (iii) yields a projective morphism $Y_{\mathrm{univ}}\rightarrow X_{\mathrm{univ}}$.

Let $X_\lambda$ denote the fiber of $W\lambda$ in $X_{\mathrm{univ}}$.  The morphism $Y_{\mathrm{univ}}\rightarrow X_{\mathrm{univ}}$ restricts to $Y_\lambda\rightarrow X_\lambda$. 
It follows from the proof of \cite[Proposition 2.9]{Losev4} that 
$Y_{\mathrm{univ}}\rightarrow \fP\times_{\fP/W}X_{\mathrm{univ}}
\rightarrow X_{\mathrm{univ}}$
is the Stein factorization. In particular, $Y_\lambda\rightarrow X_\lambda$ is a projective morphism with connected fibers between normal varieties of the same dimension. Hence it is birational and $\CC[X_\lambda]\xrightarrow{\sim}\CC[Y_\lambda]$.


In the present paper, filtered quantizations of conical symplectic singularities will play the leading role---Poisson deformations enter only as a tool for understanding quantizations. Below, we 
highlight an important example of Theorem \ref{thm:quantsofsymplectic}.

\begin{example}\label{ex:quantizationnilcone}
Let $G$ be a complex simple algebraic group and let $\cN \subset \fg^*$ be its nilpotent cone. Choose a Borel subgroup $B \subset G$ and a maximal torus $H \subset B$. As noted in Example \ref{Ex:Springer_resolution}, $X:=\cN$ is a conical symplectic singularity. For $\rho: Y \to X$, we take the Springer resolution $\rho: T^*(G/B) \to X$. There is a natural identification
$$\fP := H^2(T^*(G/B), \CC) \simeq \fh^*.$$
The regular locus of $X$ is the principal nilpotent orbit $\mathbb{O} \subset X$ and there is a single codimension 2 leaf, namely the subregular orbit $\mathbb{O}' \subset X$. If $\fg$ is simply laced, the singularity of this leaf is of the same type as $\fg$. Otherwise, this singularity corresponds to the \emph{unfolding} of the Dynkin diagram of $\fg$ (e.g. if $\fg$ is of type $G_2$, then the singularity is of type $D_4$). In the latter case, $\pi_1(\mathbb{O}')$ acts by the group of automorphisms of the simply laced Dynkin diagram which folds into $\fg$. In either case, $W$ is identified with the (Lie-theoretic) Weyl group for $\fg$. A convenient reference for these facts is \cite[Sec 8.3-8.4]{Slodowy}.

In this example, the objects $Y_{\mathrm{univ}}$, $X_{\mathrm{univ}}$, and $\cA_{\lambda}$ can be easily described. The universal graded Poisson deformation $Y_{\mathrm{univ}}$ is the homogeneous vector bundle $G\times^B \mathfrak{n}^\perp$ with its natural projection to $\fP=\fh^*$, where $\mathfrak{n}\subset \mathfrak{b}$ denotes the nilpotent radical. The variety $X_{\mathrm{univ}}$ is $\fg^*$ and the map 
$X_{\mathrm{univ}}\rightarrow \fP/W$ is the quotient morphism.   

Finally, the quantization $\cA_{\lambda}$ of $\CC[X]$ with quantization parameter $\lambda \in \fh^*$ is the corresponding \emph{central reduction} of $U(\fg)$. More precisely, the $W$-orbit of $\lambda$ defines a maximal ideal $\mathfrak{m}_{\lambda} \subset \mathfrak{Z}(\fg)$ by means of the Harish-Chandra isomorphism, and $\cA_{\lambda} \simeq U(\fg)/\mathfrak{m}_{\lambda}U(\fg)$. 
\end{example}

Fix a quantization parameter
$$\lambda = (\lambda_0,\lambda_1,...,\lambda_t) \in \bigoplus_{k=0}^t \fP_k \simeq \fP.$$
We will explain how to (partially) recover $\lambda$ from the quantization $\cA_\lambda$.

For a codimension 2 leaf $\fL_k \subset X$, consider the filtered quantization $\cA^{\Sigma_k}_{\lambda_k}$ of $\CC[\Sigma_k]$ corresponding to the quantization parameter $\lambda_k \in \fP_k$. Then $\cA_{\lambda}$ and $\cA^{\Sigma_k}_{\lambda_k}$ are related as follows. Choose a point $x\in \fL_k$. Let $\cA_{\lambda,\hbar}^{\wedge}$ (resp. $\cA_{\lambda_k, \hbar}^{\wedge}$) denote the formal completion of the Rees algebra of $\cA_{\lambda}$ (resp. $\cA_{\lambda_k}^{\Sigma_k}$) at the ideal defined by $x$ (resp. $0$), and let $\AA_{\hbar}(V)$ denote the homogeneous Weyl algebra associated to the symplectic vector space $V:=T_x\fL_k$. Then there is a $\C[[\hbar]]$-algebra isomorphism
\begin{equation}\label{eq:decompositionleaf}\cA_{\lambda,\hbar}^{\wedge}\simeq \AA_{\hbar}(V)^{\wedge} \ \widehat{\otimes}_{\C[[\hbar]]} \ \cA_{\lambda_k,\hbar}^\wedge,\end{equation} 
{which 
modulo $\hbar$ reduces to 
(\ref{eq:decompositionleaf_classical})
See, e.g., \cite[(3.2)]{Losev4}.}   

If (\ref{eq:decompositionleaf_classical}) is fixed, $W_k\lambda_k$ is uniquely determined by (\ref{eq:decompositionleaf}), see \cite[Proposition 3.6]{Losev4}.

\section{Quantizations of Kleinian singularities}\label{subsec:quantKleinian}

We will now elaborate on the classification of filtered quantizations in the case of Kleinian singularities.
On the one hand, this is an interesting special case, where the quantizations can be constructed algebraically. On the other hand, thanks to (\ref{eq:decompositionleaf}), understanding the case of Kleinian singularities is essential for understanding the general case.

Let $\Sigma = \CC^2/\Gamma$ be a Kleinian singularity, and fix the notation of Section \ref{subsec:Mckay}, e.g. $\fg$, $\fh$, and so on. By Remark \ref{rmk:PforKleinian}, there is a natural identification $\fP \simeq \fh^*$ and $W$ is identified with the Weyl group for $\fg$. Thus by Theorem \ref{thm:quantsofsymplectic}, there is a bijection
$$\fh^*/W \xrightarrow{\sim} \mathrm{Quant}(\Sigma), \qquad W\cdot \lambda \mapsto \cA_{\lambda}.$$
For Kleinian singularities, there is an explicit algebraic construction of filtered quantizations due to Crawley-Boevey and Holland (\cite{CrawleyBoevey1998}). For each element $c \in \CC[\Gamma]^{\Gamma}$ of the form
$$c = 1+ \sum_{\gamma \neq 1} c_{\gamma}\gamma$$
there is an associated filtered quantization $eH_ce$ of $\CC[\Sigma]$, constructed as follows. First, form the algebra
$$H_c := \CC\langle x,y\rangle \# \Gamma /(xy-yx-c).$$
The grading on $\CC\langle x,y \rangle \# \Gamma$ induces a filtration on $H_c$. Let $e \in \CC[\Gamma]$ denote the averaging idempotent and consider the subalgebra $eH_ce \subset H_c$. By \cite[Lem 1.1]{CrawleyBoevey1998}, $eH_ce$ is a filtered quantization of $\CC[\Sigma]$.  

The quantization parameter (in $\fh^*/W$) of $eH_ce$ can be recovered from the element $c$ as follows. Recall that the nontrivial irreducible representations $\{V_1,...,V_n\}$ of $\Gamma$ are in bijection with the simple roots $\{\alpha_1,...,\alpha_n\}$ for $\fg$. Write $\{\omega_1,...,\omega_n\}$ for the corresponding fundamental weights and consider the parameter
\begin{equation}\label{eq:defoflambdac}
\lambda^c := \frac{1}{|\Gamma|}\sum_{i=1}^n \mathrm{tr}_{V_i}(c) \omega_i \in \fh^*.\end{equation}
By \cite[Prop 3.17]{Losev4}, there is an isomorphism of filtered quantizations
\begin{equation}\label{eq:lambdatoc}
eH_ce \simeq \cA_{\lambda^c}.\end{equation}
It is worth highlighting several special cases of (\ref{eq:lambdatoc}). If $c=|\Gamma|e(=\sum_{\gamma\in \Gamma}\gamma)$, then $\mathrm{tr}_{V_i}(c)=0$ for $1\leq i \leq n$. Hence $\lambda^{|\Gamma|e}=0$ and
\begin{equation}\label{eq:lambdaGammae}eH_{|\Gamma|e}e \simeq \cA_0.\end{equation}
If $c=1$, then $\mathrm{tr}_{V_i}(c) = \dim(V_i)$ for $1 \leq i \leq n$. Hence
\begin{equation}\label{eq:barycenter}\lambda^1 = \frac{1}{|\Gamma|}\sum_{i=1}^n \dim(V_i) \omega_i.\end{equation}
The parameter $\lambda^1$ will play a particularly important role in subsequent sections. We call it the `weighted barycenter' parameter for $\Sigma$\index{barycenter parameter}. If $\Gamma=\ZZ_{n+1}$ (i.e $\Sigma$ is of type $A_n$), then $V_i$ are one-dimensional and $\lambda^1$ is the usual (unweighted) barycenter of the fundamental alcove for $\fg$.

\section{Ample cones and $\QQ$-terminalizations}\label{amplecones}
We will now recall some facts from \cite{Namikawa_bir} about the classification of $\QQ$-factorial terminalizations. Let $X$ be a conical symplectic singularity. Recall, see Section 
\ref{subsec:structurenamikawa}, that each partial Namikawa space $\fP_{\RR,k}$ is the real span of a root system, namely $\Delta_k^{\pi_1(\fL_k)}$. Let $\fP_{\RR,k}^{\geq 0} \subset \fP_{\RR,k}$ denote the fundamental chamber, and
$$\fP_{\RR}^{\geq 0} := \fP_{\RR,0} \oplus \bigoplus \fP_{\RR,k}^{\geq 0}.$$
Note that $\fP_{\RR}^{\geq 0}$ is a fundamental domain for the $W$-action on $\fP_{\RR}$. Following Namikawa (\cite[Section 1]{Namikawa_bir}), we will construct a finite set of polyhedral cones (indexed by $\QQ$-factorial terminalizations of $X$) which partition $\fP_{\RR}^{\geq 0}$.

Let $Y \to X$ be a $\QQ$-factorial terminalization. A line bundle $\mathcal{L} \in \Pic(Y)$ is said to be \emph{relatively ample} if it is ample with respect to the projective morphism $Y \to X$. Write $\Pic^a(Y) \subset \Pic(Y)$ for the semigroup of relatively ample line bundles. We can regard $\Pic^a(Y)$ as a subset of $\fP_{\RR}$ via the sequence of maps
\begin{equation}\label{eq:sequence_of_maps_Pic}
\Pic^a(Y) \to \Pic(Y^{\mathrm{reg}}) \to \Pic(Y^{\mathrm{reg}}) \otimes_{\ZZ} \RR \simeq H^2(Y^{\mathrm{reg}},\RR) \simeq \fP_{\RR}.
\end{equation}
The first map is by restriction, and the isomorphism $\Pic(Y^{\mathrm{reg}}) \otimes_{\ZZ} \RR \simeq H^2(Y^{\mathrm{reg}},\RR)$ is due to Lemma \ref{lem:c1complexification}. Note that the kernel of the map $\Pic(Y^{\mathrm{reg}})\rightarrow \fP_{\RR}$ consists of torsion line bundles on $Y^{\mathrm{reg}}$, and no torsion line bundle is relatively ample. Thus, (\ref{eq:sequence_of_maps_Pic}) is injective.

The \emph{ample cone}\index{ample cone} of $Y$ is defined to be the $\RR_{\geq 0}$-invariant subset
\begin{equation}\label{eq:Amp_def}
\mathrm{Amp}(Y) := \overline{\mathbb{R}_{\geq 0}\Pic^a(Y)} \subset \fP_{\RR}
\end{equation}
By \cite[Proposition 2.17]{BPW} we have $\mathrm{Amp}(Y) \subset \fP_{\RR}^{\geq 0}$.

In fact, the cones $\mathrm{Amp}(Y)$ are cut out by a collection of rational hyperplanes in $\fP$ which admit an alternative characterization, described below. Recall that for $\lambda\in \fP$ we let $X_\lambda$ (resp. $Y_\lambda$) denote the fiber of the universal deformation $X_{\mathrm{univ}}$ (resp. $Y_{\mathrm{univ}}$) over $W\lambda$ (resp. $\lambda$). As was discussed after Proposition \ref{prop:universaldefsymplectic}, we have a projective birational morphism $Y_\lambda\rightarrow X_\lambda$. 
Let $\fP^{\mathrm{sing}} \subset \fP$ be the subset consisting of all $\lambda \in \fP$ for which this morphism fails to be an isomorphism,
(equivalently, for which $Y_\lambda$ is not affine).

\begin{theorem}[\cite{Namikawa_bir}, Lem. 1, Main Thm.]\label{thm:Q-term cones}
The following are true:

\begin{itemize}
    \item[(i)] There are finitely many isomorphism classes of $\QQ$-terminalizations of $X$.
    \item[(ii)] Each $\QQ$-terminalization is determined up to isomorphism by its ample cone.
    \item[(iii)] The ample cone of each $\QQ$-terminalization is polyhedral, and of full dimension in $\fP_{\RR}^{\geq 0}$.
    \item[(iv)] The ample cones of all $\QQ$-terminalizations partition $\fP_{\RR}^{\geq 0}$.
    \item[(v)] The subset $\mathfrak{P}^{\mathrm{sing}}\subset \mathfrak{P}$ is the $W$-stable union of rational hyperplanes in $\fP$, including the walls corresponding to the $W$-action. The ample cones $\mathrm{Amp}(Y)$ are cut out by these hyperplanes.
\end{itemize}
\end{theorem}

\section{Automorphisms}\label{subsec:automorphisms}

We will need some information about the automorphism groups of Poisson deformations and filtered quantizations of conical symplectic singularities. Let $X$ be a conical symplectic singularity so that $\CC[X]$ is a graded Poisson algebra. Consider the group $\overline{\mathcal{G}}$ of graded Poisson automorphisms of $\CC[X]$. Since $\CC[X]$ is positively graded, $\overline{\mathcal{G}}$ is an algebraic group.  Choose a maximal reductive subgroup $\mathcal{G} \subset \overline{\mathcal{G}}$. Note that $\overline{\mathcal{G}}$ acts on $\mathrm{PDef}(X)$ by
$$g \cdot (\cA^0,\theta) = (\cA^0, g \circ \theta), \qquad g \in \overline{\mathcal{G}}, \ (\cA^0,\theta) \in \mathrm{PDef}(X).$$
In a similar manner, one can define a $\overline{\mathcal{G}}$-action on $\mathrm{Quant}(X)$. These actions induce two (a priori distinct) $\overline{\mathcal{G}}$-actions on $\fP/W$ via Theorems \ref{thm:defsofsymplectic} and \ref{thm:quantsofsymplectic}, called the Poisson and the quantum actions. For example, the Poisson action is defined by regarding $\fP/W$ as the parameter space for filtered Poisson deformations of $\CC[X]$.


Write $\Aut(\cA^0_{\lambda})$ (resp. $\Aut(\cA_{\lambda})$) for the group of filtered Poisson algebra automorphisms of $\cA^0_{\lambda}$ (resp. filtered algebra automorphisms of $\cA_{\lambda}$). Both are algebraic groups. The following result, obtained in \cite[Sec 3.7]{Losev4}, describes the relationship between $\Aut(\cA_{\lambda}^0)$ and $\Aut(\cA_{\lambda})$.

\begin{prop}\label{prop:automorphisms}
The following are true:
\begin{itemize}
\item[(i)] The Poisson and quantum actions of $\overline{\mathcal{G}}$ coincide and are trivial on $\overline{\mathcal{G}}^{\circ}$.
\item[(ii)] Write  $\mathcal{G}_\lambda$ for the stabilizer of $W\lambda$ in $\mathcal{G}${, this is a reductive subgroup}. There is an epimorphism 
$\Aut(\cA_{\lambda})\twoheadrightarrow \mathcal{G}_\lambda$ with unipotent kernel. It is given by taking the associated graded and then projecting onto the maximal reductive subgroup.
\item[(iii)] Similarly, there is an epimorphism 
$\Aut(\cA^0_{\lambda})\twoheadrightarrow \mathcal{G}_\lambda$ with unipotent kernel.
    %
\end{itemize}
%
\end{prop}

\begin{proof}
(i) is \cite[Prop 3.21]{Losev4}. The claim that the map defined in (ii) is an epimorphism is easy, compare to the proof of \cite[Lem 3.19]{Losev4}. The kernel consists of automorphisms of $\cA_\lambda$ such that the assoiciated graded automorphism is unipotent. Such automorphisms are unipotent. So the kernel is unipotent.  The proof of (iii) is the same. 
%
%

\end{proof}

\begin{rmk}\label{rmk:equiv_automor}
Let $G\subset \mathcal{G}$ be a connected reductive subgroup. Thanks
to (i) {and (ii) (resp., (iii))} of Proposition \ref{prop:automorphisms}, {the $G$-action on $\C[X]$ lifts to an action of $\cA_\lambda$ (resp., $\cA^0_\lambda$) by filtered algebra (resp., filtered Poisson algebra) automorphisms so that the isomorphism $\gr\cA_\lambda\xrightarrow{\sim} \CC[X]$ is $G$-equivariant (and the analogous condition for $\cA^0_\lambda$). Moreover, thanks to (ii) and basic facts from the theory of algebraic groups, the lifted action on $\cA_\lambda$ is unique up to conjugating by an automorphism which is trivial on the associated graded algebra. The analogous claim holds for $\cA_\lambda^0$.} 

So there is a natural bijection between the following two sets:
\begin{itemize}
    \item The set of filtered  quantizations of $\CC[X]$ considered up to $G$-equivariant filtered algebra automorphisms.
    \item The set of filtered Poisson deformations of $\CC[X]$ considered up to $G$-equivariant filtered Poisson algebra isomorphisms.
\end{itemize}
By (ii) and (iii) of Proposition \ref{prop:automorphisms},
the maximal reductive subgroups in $\Aut_G(\cA_{\lambda}), 
\Aut_G(\cA^0_{\lambda})$ coincide with $Z_{\mathcal{G}}(G)\cap \mathcal{G}_\lambda$.  
%
%
%
%
\end{rmk}

\section{Hamiltonian quantizations}\label{subsec:equivariant}

Let $A$ be a graded Poisson algebra of the type described in Section \ref{subsec:quant}, and let $G$ be an algebraic group. Suppose $G$ acts rationally on $A$ by graded Poisson automorphisms. The $G$-action gives rise to a Lie algebra homomorphism
$$\fg \to \mathrm{Der}(A), \qquad \xi \mapsto \xi_A,$$
where $\mathrm{Der}(A)$ is the Lie algebra consisting of derivations of $A$. We say that $A$ is \emph{Hamiltonian} if there is a $G$-equivariant map $\varphi: \mathfrak{g} \to A_d$ (called a \emph{classical co-moment map})\index{co-moment map!classical} such that
$$\{\varphi(\xi), a\} = \xi_A(a), \qquad \xi \in \fg, \quad a \in A.$$ 
Now suppose $(\cA,\theta)$ is a filtered quantization of $A$.
We say that $(\cA,\theta)$ is $G$-\emph{equivariant}\index{quantization!equivariant} if $G$ acts rationally on $\cA$ by filtered algebra automorphisms and the isomorphism $\theta: \gr(\cA) \xrightarrow{\sim} A$ is $G$-equivariant. In this setting (as above) we get a Lie algebra homomorphism
$$\fg \to \mathrm{Der}(\cA), \qquad \xi \mapsto \xi_{\cA}.$$

\begin{definition}\label{def:hamiltonian}
Suppose $A$ is a graded Poisson algebra equipped with a Hamiltonian $G$-action.
\begin{itemize}
    \item A \emph{Hamiltonian} quantization\index{quantization!Hamiltonian} of $A$ is a triple $(\cA,\theta,\Phi)$ consisting of
 \begin{itemize}
     \item[(i)] a $G$-equivariant filtered quantization $(\cA,\theta)$ of $A$, and
     \item[(ii)] a $G$-equivariant map $\Phi: \fg \to \cA_{\leq d}$ (called a \emph{quantum co-moment map})\index{co-moment map!quantum} such that
     $$[\Phi(\xi),a] = \xi_{\cA}(a), \qquad \xi \in \fg, \quad a \in \cA.$$
 \end{itemize}
 \item An isomorphism $(\mathcal{A}_1,\theta_1,\Phi_1) \xrightarrow{\sim} (\mathcal{A}_2,\theta_2,\Phi_2)$ of Hamiltonian quantizations of $A$ is a $G$-equivariant isomorphism of filtered algebras $\phi: \cA_1 \to \cA_2$ such that $\theta_1 = \theta_2 \circ \gr(\phi)$ and $\Phi_2 = \phi \circ \Phi_1$. Denote the set of isomorphism classes of Hamiltonian quantizations of $A$ by $\mathrm{Quant}^G(A)$.
\end{itemize}
\end{definition}

Note that if $\cA$ is a $G$-equivariant quantization and $\Phi: \fg \to \cA_{\leq d}$ is a quantum co-moment map, then $\gr(\Phi): \fg \to A_d$ is a classical co-moment map. In the other direction, we have the following result.

\begin{lemma}\label{lem:co-momentunique}
Let $G$ be a connected reductive algebraic group, and let $X=\Spec(A)$ be a conical symplectic singularity with a Hamiltonian $G$-action {such that the connected component of the center, $Z(G)^\circ$, acts trivially}. Then the following are true:
\begin{itemize}
    \item[(i)] There is a unique classical co-moment map $\varphi: \fg \to A_d$.
    \item[(ii)] For each $G$-equivariant quantization $\cA$ of $A$ and Lie algebra homomorphism $\chi:\fg \to \CC$, there is a unique quantum co-moment map
    $$\Phi_{\chi}: \fg \to \cA_{\leq d}$$
   such that $\Phi_{\chi}|_{\fz(\fg)} = \chi$.
\end{itemize}
\end{lemma}

\begin{proof}
Suppose $\varphi_1,\varphi_2:\fg \to A_d$ are classical co-moment maps. Then $\varphi_1-\varphi_2: \fg \to A_d$ takes values in the Poisson center of $A$. Since $\Spec(A)$ is generically symplectic, the Poisson center is $\CC = A_0$, and therefore $\varphi_1-\varphi_2=0$. This proves (i).

Now suppose $\cA$ is a $G$-equivariant quantization of $A$. If $a \in \cA_k$ is central in $\cA$, then $a+\cA_{k-1} \in A_k$ is Poisson-central in $A$. Thus, the center of $\cA$ is $\CC=\cA_{\leq 0}$. Suppose $\Phi_1,\Phi_2:\fg \to \cA_{\leq d}$ are quantum co-moment maps. Then $\Phi_1-\Phi_2: \fg \to \cA_{\leq d}$ is a Lie algebra homomorphism with values in $Z(\cA) =\CC$. That is, $\Phi_1$ and $\Phi_2$ differ by a character of $\fz(\fg)$. {This establishes the uniqueness part of (ii).} It remains to show that the (unique) classical co-moment map $\varphi: \fg \to A_d$ lifts to a quantum co-moment $\Phi:\fg \to \cA_{\leq d}${: $Z(G)^\circ$ acts trivially on $A$, hence on $\cA$, and hence any quantum comoment map is constant on $\mathfrak{z}(\g)$}. We can assume in the proof that $G$ acts faithfully on $A$ {(in particular, $G$ is semisimple)} and so $\fg\hookrightarrow \operatorname{Der}(\cA)$. 

First, observe that every derivation $\Delta$ of $\cA$ satisfying $\Delta(\cA_{\leq i}) \subseteq \cA_{\leq k+i}$ is of the form $[\widetilde{a},\bullet]$ for some $\widetilde{a} \in \cA_{\leq k+d}$. Indeed, let $\delta = \gr(\Delta)$, a Poisson derivation of degree $k$. By \cite[Proposition 2.14]{Losev4}, all such derivations are inner. Thus, $\delta=\{a_0, \bullet\}$ for some element $a_0 \in A_{k+d}$. Choose a lift $\widetilde{a}_0 \in \cA_{\le k+d}$ of $a_0$ and let $\Delta^1=\Delta-[\widetilde{a}_0, \bullet]$. Note that $\Delta^1$ is a filtered derivation of degree $k-1$. Iterating this procedure, we obtain an element $\widetilde{a} \in \cA_{\le k+d}$ such that $\Delta-[\widetilde{a}, \bullet]$ is a filtered derivation of degree $-d-1$. Since the grading on $A$ is by nonnegative integers and $\gr(\Delta-[\widetilde{a}, \bullet])$ is inner, $\gr(\Delta-[\widetilde{a}, \bullet])=0$. Thus, $\Delta = [\widetilde{a},\bullet]$, as desired.

There is a $G$-equivariant Lie algebra embedding $\cA_{\leqslant d}/\cA_{\leqslant 0}\hookrightarrow \operatorname{Der}(\cA)$ given by $a\mapsto [a,\cdot]$. By the above argument, the image of $\fg$ in $\operatorname{Der}(\cA)$ lies in that of $\cA_{\leqslant d}/\cA_{\leqslant 0}$. So the map $\xi\mapsto \xi_{\cA}$ can be viewed as a $G$-equivariant Lie algebra embedding $\fg\hookrightarrow \cA_{\leqslant d}/\cA_{\leqslant 0}$. Since $G$ is {is assumed to be semisimple}, this embedding lifts to a $G$-equivariant Lie algebra embedding $\fg\hookrightarrow \cA_{\leqslant d}$. This lift is our quantum comoment map. 

\end{proof}

If $X$ is a conical symplectic singularity with a Hamiltonian $G$-action {such that $Z(G)^\circ$ acts trivially}, we will consider the \emph{extended Namikawa space}\index{Namikawa space!extended}
$$\overline{\fP} := \fP \oplus \fz(\fg)^*.$$
This space should be viewed as an equivariant version of the Namikawa space for $X$. There is a natural $W$-action on $\overline{\fP}$ defined via the decomposition above (the $W$-action on the second factor is trivial). For each $\lambda = (\lambda',\chi) \in \fP \oplus \fz(\fg)^* = \overline{\fP}$, there is an associated Hamiltonian quantization $(\cA_{\lambda'},\Phi_{\chi})$, where $\cA_{\lambda'}$ is the filtered quantization of $A$ considered in Theorem \ref{thm:quantsofsymplectic} and $\Phi_{\chi}: \fg \to \cA_{\lambda'}$ is the quantum co-moment map considered in Lemma \ref{lem:co-momentunique}(ii). For notational convenience, we will often write $(\cA_{\lambda},\Phi_{\lambda}) := (\cA_{\lambda'},\Phi_{\chi})$. The following is an easy consequence of Theorem \ref{thm:quantsofsymplectic} and Lemma \ref{lem:co-momentunique}.

\begin{prop}\label{prop:classificationHamiltonian}
In the setting of Lemma \ref{lem:co-momentunique}, there is a bijection
$$\overline{\fP}/W \xrightarrow{\sim} \mathrm{Quant}^G(X), \qquad W \cdot \lambda \mapsto (\cA_{\lambda},\Phi_{\lambda}).$$
\end{prop}

\section{Harish-Chandra bimodules for filtered quantizations}\label{subsec:HCbimods}

Let $A$ be a {nonnegatively} graded Poisson algebra {with bracket of degree $-d$} and let $\cA$ be a filtered quantization of $A$. Let $\cB$ be an $\cA$-bimodule. A \emph{compatible filtration} on $\cB$ is an increasing filtration by subspaces
$$0 = \B_{-1} \subseteq \B_0 \subseteq \B_1 \subseteq ...,  \qquad \bigcup_{i=0}^{\infty} \B_i = \B$$
such that
$$\cA_{\leq i} \B_{\leq j} \subseteq \B_{\leq i+j} \quad \text{and} \quad [\cA_{\leq i},\cB_{\leq j}] \subseteq \B_{\leq i+j-d}, \qquad \forall i, j \in \ZZ_{\geq 0}.$$
Under these conditions, $\gr(\B)$ has the structure of a graded $A$-module. A compatible filtration is \emph{good} if $\gr(\B)$ is finitely-generated for $A$.

\begin{definition}\label{def:HCbimodfiltered}
A \emph{Harish-Chandra bimodule}\index{Harish-Chandra bimodule!for a quantization} for $\cA$ is an $\cA$-bimodule which admits a good filtration. A morphism of Harish-Chandra bimodules is a homomorphism of $\cA$-bimodules. Denote the category of Harish-Chandra bimodules for $\cA$ by $\HC(\cA)$.
\end{definition}

Note that if $\cA = U(\fg)$, for $\fg$ a reductive Lie algebra, we recover Definition \ref{def:HCbimodsclassical}. 

For $\mathcal{B} \in \HC(\cA)$, we can define the associated variety $\mathcal{V}(\mathcal{B}) \subset \Spec(A)$\index{associated variety!of Harish-Chandra bimodule} and generic multiplicities $m_Z(\mathcal{B}) \in \ZZ_{>0}$\index{generic multiplicity} as in the Lie-theoretic context (see the remarks following Definition \ref{def:HCbimodsclassical}). We say that $\mathcal{B}$ has \emph{full support} if $\mathcal{V}(\mathcal{B}) = \Spec(A)$. Note that the regular bimodule $\cA$ is Harish-Chandra and satisfies this condition. Write $\HC_{\partial}(\cA) \subset \HC(\cA)$ for the Serre subcategory of Harish-Chandra bimodules which are \emph{not} of full support. Note that $\HC(\cA)$ is a monoidal category (with respect to $\otimes_{\cA}$) and $\HC_{\partial}(\cA)$ is a monoidal ideal (see \cite[Lem 2.13]{LosevHC}). Consider the quotient category
$$\overline{\HC}(\cA) = \HC(\cA)/\HC_{\partial}(\cA).$$
We call $\overline{\HC}(\cA)$ the category of Harish-Chandra bimodules with full support for $\cA$. It is a monoidal category because $\HC_{\partial}(\cA)$ is a monoidal ideal.

\begin{lemma}\label{lem:simplenoboundary}
If the algebra $\cA$ is  simple, then $\HC_{\partial}(\cA)=0$.
\end{lemma}

\begin{proof}
Let $\cB\in \HC_{\partial}(\cA)$, and consider the algebra $\End_{\cA}(\cB)$ of left $\cA$-module endomorphisms of $\cB$. 
There is a natural filtration on $\End_{\cA}(\cB)$ defined by
$$\End_{\cA}(\cB)_{\leq i} := \{f \in \End_{\cA}(\cB) \mid f(\cB_{\leq k}) \subseteq \cB_{\leq k+i} \quad \forall k \in \ZZ\}$$
It is straightfoward to check that the above filtration is good. In particular, $\End_{\cA}(\cB) \in \HC(\cA)$. There is a natural filtered algebra homomorphism 
$$\varphi: \cA \to \End_{\cA}(\cB)^{\mathrm{op}}, \qquad \varphi(a)(b)=ba.$$
Note that $\mathcal{V}(\mathrm{End}_{\cA}(\cB)^{\mathrm{op}}) \subseteq \mathcal{V}(\cB) \subsetneq \Spec(A) = \mathcal{V}(\cA)$. Thus, $\ker{(\varphi)} = \mathrm{RAnn}(\cB)$ is a nonzero two-sided ideal in $\cA$. Since $\cA$ is simple, this implies that $\mathrm{RAnn}(\cB)=\cA$, and hence that $\cB=0$.
\end{proof}

\begin{rmk}\label{rmk:HC_sheaf}
Similarly, if $\calD$ is a sheaf quantization of a Poisson variety $X$ (cf. Definition \ref{def:filteredquantX}), it makes sense to consider Harish-Chandra bimodules for $\calD$. This category will be denoted by $\HC(\calD)$.
\end{rmk}

\section{Restriction and extension functors for Harish-Chandra bimodules}\label{subsec:daggers}

Let $X$ be a conical symplectic singularity and let $\cA_{\lambda}$ be a quantization of $\CC[X]$. Let $\Gamma$ be the algebraic fundamental group of the regular locus $X^{\mathrm{reg}} \subset X$. Recall from Section \ref{subsec:finitecovers} that $\Gamma$ is the maximal finite quotient of $\pi_1(X^{\mathrm{reg}})$. There is a close relationship between Harish-Chandra bimodules with full support for $\cA_{\lambda}$ and (finite-dimensional) representations of $\Gamma$, which we will now describe.

In \cite[Sec 4]{LosevHC}, the first-named author constructs an exact monoidal functor\index{restriction functor}
$$\bullet_{\dagger}: \HC(\cA_{\lambda}) \to \Gamma\modd.$$
The image of $\cB\in \HC(\cA_\lambda)$ under this functor is characterized by the following property: let $\widehat{X}' \to X^{\mathrm{reg}}$ be the universal algebraic cover, i.e. the unique Galois cover with Galois group $\Gamma$. Then for $\cB \in \HC(\cA_{\lambda})$ (and any good filtration), there is an isomorphism of Poisson coherent sheaves on $X^{\mathrm{reg}}$
\begin{equation}\label{eqn:grdagger}
(\gr \cB)|_{X^{\mathrm{reg}}} \simeq \left(\cB_{\dagger} \otimes \cO_{\widehat{X}'}\right)^{\Gamma}.\end{equation}
We will give a formal construction of $\bullet_{\dagger}$ towards the end of this section.

From (\ref{eqn:grdagger}) it follows that $\ker{(\bullet_{\dagger})} = \HC_{\partial}(\cA_{\lambda})$. Thus, $\bullet_{\dagger}$ descends to an embedding
$$\bullet_{\dagger}: \overline{\HC}(\cA_{\lambda}) \hookrightarrow \Gamma\modd.$$
This embedding is full and monoidal, and its image in $\Gamma\modd$ is closed under taking direct summands, see \cite[Lem 4.6]{LosevHC}. Thus, $\mathrm{Im}(\bullet_{\dagger})=\Gamma/\Gamma(\lambda)\modd$ for a normal subgroup $\Gamma(\lambda) \subseteq \Gamma$, which is uniquely determined by $\lambda$. We will provide a more explicit description of $\Gamma(\lambda)$ in Section \ref{subsec:Gammalambda}.

In \cite[Sec 4.4]{LosevHC}, it is shown that $\bullet_{\dagger}$ admits a right adjoint\index{extension functor}
$$\bullet^{\dagger}: \Gamma\modd \to \HC(\cA_{\lambda}).$$
For any $\cB \in \HC(\cA_{\lambda})$, the kernel and cokernel of the unit homomorphism $\cB \to (\cB_{\dagger})^{\dagger}$ are supported on $X^{\mathrm{sing}}$, see \cite[Lem 4.6]{LosevHC}. Hence, $\bullet^{\dagger}: \Gamma \modd \to \overline{\HC}(\cA_{\lambda})$ is a left inverse for $\bullet_{\dagger}: \overline{\HC}(\cA_{\lambda}) \to \Gamma\modd$, and there are quasi-inverse equivalences
\begin{equation}\label{eq:daggers_inverse}
\bullet_{\dagger}: \overline{\HC}(\cA_{\lambda}) 
\xrightarrow{\sim}\Gamma/\Gamma(\lambda)\modd, \qquad \bullet^{\dagger}: \Gamma/\Gamma(\lambda)\modd \xrightarrow{\sim} \overline{\HC}(\cA_{\lambda}).
\end{equation}
Next, we will give constructions of $\bullet_{\dagger}$ and $\bullet^{\dagger}$, following \cite[Section 4]{LosevHC} with some modifications. Consider the universal cover $\widehat{X}' \to X^{\mathrm{reg}}$, and let $\widehat{X}=\Spec(\CC[\widehat{X}'])$. The Poisson variety $\widehat{X}$ is a finite Galois cover of $X$ in the sense of Section \ref{subsec:finitecovers}. As we mentioned there, it is a conical symplectic singularity. Consider the \'{e}tale lift of $\mathcal{D}^{X^{\mathrm{reg}}}_\lambda$ to $\widehat{X}'$, compare to 
\cite[Section 2.6]{LosevHC}. Denote the resulting sheaf of algebras by $\widehat{\mathcal{D}}'$. It has an action of $\Gamma$ by filtered algebra automorphisms with $(\widehat{\mathcal{D}}')^\Gamma\simeq \mathcal{D}^{X^{\mathrm{reg}}}_\lambda$. Consider the category of $\Gamma$-equivariant Harish-Chandra bimodules over $\widehat{\mathcal{D}}'$, denoted $\HC^\Gamma(\widehat{\mathcal{D}}')$.
Taking $\Gamma$-invariants defines an equivalence of categories
$\HC^\Gamma(\widehat{\mathcal{D}}')
\rightarrow \HC(\mathcal{D}^{X^{\mathrm{reg}}}_\lambda)$, compare to the penultimate paragraph in \cite[Section 4.2]{LosevHC}. There is a full embedding $\Gamma\operatorname{-mod}
\hookrightarrow \HC^\Gamma(\widehat{\mathcal{D}}')$ 
given by $V\mapsto \widehat{\mathcal{D}}'\otimes V$. This embedding is an equivalence: using a $\Gamma$-invariant regluing as in \cite[Section 4.2]{LosevHC} we get an equivalence of categories $\HC^\Gamma(\widehat{\mathcal{D}}') \xrightarrow{\sim} \HC^\Gamma(\widetilde{\cA}^0_{\lambda^1})$, {where $\lambda^1=(\lambda_0,0,\ldots,0)$ under (\ref{eq:partialdecomp}), $\widetilde{\cA}_{\lambda^1}$ denotes the quantization of $\widehat{X}$ with parameter $\lambda^1$, $\widetilde{\cA}_{\lambda^1}$ is its microlocalization to $\widehat{X}'$, and} $\HC^\Gamma(\widetilde{\cA}^0_{\lambda^1})$ is {the category of $\Gamma$-equivariant HC $\widetilde{\cA}^0_{\lambda^1}$-bimodules} as defined in \cite[Section 4.2]{LosevHC}. {Namely, these are HC bimodules $\mathcal{B}$ equipped with $\Gamma$-representations making the structure map $\widetilde{\cA}^0_{\lambda^1}\otimes \mathcal{B}\otimes \widetilde{\cA}^0_{\lambda^1}\rightarrow \mathcal{B}$ equivariant.}

Furthermore, $\HC^\Gamma(\widetilde{\cA}^0_{\lambda^1})$ is equivalent to $\Gamma\operatorname{-mod}$, see \cite[Section 4.3]{LosevHC}. The functor $\bullet^{\dagger}:\Gamma\operatorname{-mod}\rightarrow \HC(\cA_\lambda)$
is then given by
\begin{equation}\label{eq:upper_dagger_construction}
    V\mapsto \Gamma(\widehat{X}', \widehat{\mathcal{D}}'\otimes V)^\Gamma.
\end{equation}
This construction is equivalent to that of \cite[Section 4]{LosevHC}. There the equivalence 
$\Gamma\operatorname{-mod}\xrightarrow{\sim}
\HC(\mathcal{D}|_{X^{\mathrm{reg}}})$ was the composition
\begin{equation}\label{eq:compos_equiv0}\Gamma\operatorname{-mod}\xrightarrow{\sim}
\HC^\Gamma(\widetilde{\cA}^0_{\lambda^1})\xrightarrow{\sim} \HC((\widetilde{\cA}^0_{\lambda^1})^\Gamma)\xrightarrow{\sim}\HC(\mathcal{D}|_{X^{\mathrm{reg}}}),\end{equation}
where the third equivalence is regluing. Here we use the composition  
\begin{equation}\label{eq:compos_equiv}
\Gamma\operatorname{-mod}\xrightarrow{\sim}
\HC^\Gamma(\widetilde{\cA}^0_{\lambda^1})\xrightarrow{\sim}\HC^\Gamma(\widehat{\mathcal{D}}') \xrightarrow{\sim}\HC(\mathcal{D}|_{X^{\mathrm{reg}}}),
\end{equation}
where the second equivalence is regluing. Since we are regluing with the same cocycles, (\ref{eq:compos_equiv0}) and (\ref{eq:compos_equiv}) define isomorphic tensor functors. {We emphasize that we do not need the regluing procedure to define the functor $\Gamma\operatorname{-mod}\rightarrow \HC(\calD|_{X^{reg}})$; we only need it to prove that this functor is an equivalence.}

The functor $\bullet_{\dagger}$ can be constructed as the composition of the microlocalization functor $\HC(\A_\lambda)\rightarrow \HC(\mathcal{D}|_{X^{\mathrm{reg}}})$
and the inverse of the composition (\ref{eq:compos_equiv}). This construction shows that $\bullet_{\dagger}$ is tensor, since it is the composition of tensor functors.

Finally, it is evident from the constructions that $\bullet^\dagger$ and $\bullet_{\dagger}$ respect good filtrations. For example, a good filtration on $\mathcal{B}\in \HC(\A_\lambda)$ gives rise to a filtration on $\mathcal{B}_\dagger$. To see this, we note that the microlocalization and global section functors as well as all equivalences in (\ref{eq:compos_equiv}) preserve good filtrations (and the functors (\ref{eq:compos_equiv}) are equivalences of the categories of filtered objects).

\section{Finite covers of filtered quantizations}\label{subsec:coverings}

Let $X$ be a conical symplectic singularity and let $\cA^X=\cA^X_{\lambda}$ be a filtered quantization of $\CC[X]$ (here we introduce the convention, adopted throughout the paper, of indicating the underlying variety of a filtered quantization as a superscript). Let $p: \widetilde{X} \to X$ be a finite cover, as defined in Section \ref{subsec:finitecovers}.

\begin{definition}\label{defi:quant_cover}  
A \emph{cover of} $\cA^X$ (with respect to the map $p:\widetilde{X} \to X$) is 
a filtered quantization $\cA^{\widetilde{X}}$ of $\CC[\widetilde{X}]$ which admits a filtered algebra embedding $\cA^X\hookrightarrow \cA^{\widetilde{X}}$ such that the associated graded homomorphism is the pullback map $\CC[X]\hookrightarrow \CC[\widetilde{X}]$.
\end{definition}


Recall the normal subgroup $\Gamma(\lambda) \subset \Gamma$ defined in Section \ref{subsec:daggers}, and consider the Galois cover $\breve{X}$ of $X$ with Galois group $\Gamma/\Gamma(\lambda)$. Let $\breve{\mathcal{D}}'$ denote the \'{e}tale lift of $\mathcal{D}|_{X^{\mathrm{reg}}}$ to $\breve{X}'$, and let $\breve{\cA}:=\Gamma(\breve{X}', \breve{\mathcal{D}}')$. Note that $\Gamma/\Gamma(\lambda)$ acts on 
$\breve{\cA}$ and $\breve{\cA}^{\Gamma/\Gamma(\lambda)}$
is identified with $\cA^X$. Also note that we have a natural inclusion $\gr \breve{\cA}\hookrightarrow \CC[\breve{X}]$. This inclusion may fail to be an isomorphism, so $\breve{\cA}$ may fail to be a cover of $\cA$.

The following lemma explains the relationship between $\breve{\cA}$ and $\bullet^{\dagger}$.

\begin{lemma}\label{lem:breve_A}
The following are true:
\begin{itemize}
    \item[(i)] For $V\in \Gamma/\Gamma(\lambda)\operatorname{-mod}$, we have $V^\dagger=(\breve{\cA}\otimes V)^{\Gamma/\Gamma(\lambda)}$.
    \item[(ii)] The natural morphism 
    $\breve{X}\rightarrow \operatorname{Spec}(\gr \breve{\cA})$ is an isomorphism over $X^{\mathrm{reg}}$.
\end{itemize}
\end{lemma}
\begin{proof}
In the notation of (\ref{eq:upper_dagger_construction}), there is an isomorphism of filtered quantizations of $\breve{X}'$, $\breve{\mathcal{D}}\simeq \widehat{\mathcal{D}}^{\Gamma(\lambda)}$. Now (i) follows from  (\ref{eq:upper_dagger_construction}).

Next we prove (ii). {Note that $\CC[X]\hookrightarrow \gr \breve{\cA}\hookrightarrow \CC[\breve{X}]$. Since $\CC[\breve{X}]$ is finite over $\CC[X]$, we see that $\gr\breve{\cA}$ is finite over $\CC[X]$.}

{We claim that $\Spec(\gr\breve{A})\rightarrow X$ is \'{e}tale over $X^{reg}$.
Note that $\CC[\breve{X}]$ is a Poisson $\CC[X]$-algebra and $\gr\breve{\cA}$ is a Poisson subalgebra. Consider a point $x\in X^{reg}$. The completion $\CC[X]^{\wedge_x}$ is the algebra of formal power series with the Darboux bracket. Therefore, every finitely generated Poisson $\CC[X]^{\wedge_x}$-module $M$ is of the form $\CC[X]^{\wedge_x}\otimes M_0$, where $M_0$ is the Poisson centralizer of $\CC[X]^{\wedge_x}$ in $M$.
Consider the Poisson $\CC[X]^{\wedge_x}$-algebras $(\gr \breve{\cA})^{\wedge_x}(:=\CC[X]^{\wedge_x}\otimes_{\CC[X]})\gr\breve{\cA}\subset \CC[\breve{X}]^{\wedge_x}$.
The Poisson centralizers of $\CC[X]^{\wedge_x}$ in these completions  are subalgebras, denote them by $C_1\subset C_2$. Since $\breve{X}\rightarrow X$ is \'{e}tale over $X^{reg}$, we see that $C_2$ is the direct sum of several copies of $\CC$. It follows that $C_1$ is also the direct sum of several copies of $\CC$, hence $\Spec(\gr\breve{\cA})\rightarrow X$ is \'{e}tale over $X^{reg}$.} 

So, it suffices to show that the generic degree of $\operatorname{Spec}(\gr \breve{\cA})$ over $X$ is $|\Gamma/\Gamma(\lambda)|$. By (\ref{eqn:grdagger}), it is enough to show that 
\begin{equation}\label{eq:dim_dagger}
\dim \breve{\cA}_\dagger=|\Gamma/\Gamma(\lambda)|. 
\end{equation}
Thanks to (i), $\breve{\cA}=V^\dagger$ for $V=\C[\Gamma/\Gamma(\lambda)]$. 
Now (\ref{eq:dim_dagger}) follows from (\ref{eq:daggers_inverse}).
\end{proof}

The next proposition relates covers of $\cA$ and their categories of Harish-Chandra bimodules to $\breve{\cA}$.

\begin{prop}\label{Prop:quant_covers}
Let $\cA^{\widetilde{X}}$ be a cover of $\cA^X$. The following are true:
\begin{itemize}
    \item[(i)] There is a subgroup $\Gamma^0\subset \Gamma/\Gamma(\lambda)$ with the following properties:
    \begin{itemize}
        \item[$\bullet$] $\widetilde{X} \simeq \breve{X}/\Gamma^0$.
        \item[$\bullet$] There is an isomorphism of filtered algebras $\cA^{\widetilde{X}}\xrightarrow{\sim}\breve{\cA}^{\Gamma^0}$ with associated graded equal to the natural isomorphism $\CC[\widetilde{X}]
    \xrightarrow{\sim} \CC[\breve{X}]^{\Gamma^0}$.
    \end{itemize}
    \item[(ii)] There is an equivalence of categories $\Gamma^0\operatorname{-mod}\xrightarrow{\sim} \overline{\HC}(\cA^{\widetilde{X}})$ given by $V\mapsto (\breve{\cA}\otimes V)^{\Gamma^0}$. 
\end{itemize}
\end{prop}
\begin{proof}
First we prove (i). Let $\Gamma_0$ be a subgroup in $\Gamma$ (defined up to conjugacy) such that $\widetilde{X}\simeq \widehat{X}/\Gamma_0$ and $\widetilde{\mathcal{D}}'\simeq \widehat{\mathcal{D}}'^{\Gamma_0}$. It follows that $\cA^{\widetilde{X}}_\dagger\simeq\CC[\Gamma/\Gamma_0]$.
The action of $\Gamma$ on the image of every Harish-Chandra $\cA^X$-bimodule under $\bullet_\dagger$ factors through $\Gamma/\Gamma(\lambda)$ by the construction of $\Gamma(\lambda)$ in Section \ref{subsec:daggers}. Thus $\Gamma(\lambda)\subset \Gamma_0$. Set $\Gamma^0=\Gamma_0/\Gamma(\lambda)$. We wish to establish an isomorphism $\cA^{\widetilde{X}}\xrightarrow{\sim} \breve{\cA}^{\Gamma^0}$. We begin by exhibiting a filtered algebra isomorphism 
\begin{equation}\label{eq:dagger_alg_isom}
\cA^{\widetilde{X}}_\dagger\xrightarrow{\sim} (\breve{\cA}^{\Gamma^0})_\dagger.
\end{equation} 
Note that both sides carry natural filtered algebra structures. This is evident from the construction of $\bullet_{\dagger}$ at the end of Section \ref{subsec:daggers}. The left hand side comes with a filtration with the property that there is a graded $\Gamma$-equivariant algebra isomorphism $\gr(\cA^{\widetilde{X}}_{\dagger})\xrightarrow{\sim}\C[\Gamma/\Gamma_0]$, where the target is concentrated in degree $0$. It follows that there is a filtered algebra isomorphism $\cA^{\widetilde{X}}_{\dagger}\xrightarrow{\sim}\C[\Gamma/\Gamma_0]$. Similarly, $\breve{\cA}_\dagger\xrightarrow{\sim} \C[\Gamma/\Gamma(\lambda)]$ and hence $(\breve{\cA}^{\Gamma^0})_\dagger\xrightarrow{\sim} \C[\Gamma/\Gamma_0]$. This proves (\ref{eq:dagger_alg_isom}).

Next we will show that (\ref{eq:dagger_alg_isom}) implies $\cA^{\widetilde{X}} \xrightarrow{\sim} \breve{\cA}^{\Gamma^0}$. Indeed, by the construction of $\bullet_\dagger$, (\ref{eq:dagger_alg_isom}) comes from a filtered algebra isomorphism between the microlocalizations of $\cA^{\widetilde{X}}, \breve{\cA}^{\Gamma^0}$ to $X^{\mathrm{reg}}$. Passing to the global sections of the microlocalizations, we get a filtered algebra isomorphism $\cA^{\widetilde{X}} \xrightarrow{\sim} \breve{\cA}^{\Gamma^0}$. This completes the proof of (i).

Finally we prove (ii). Let $\widetilde{\Gamma}\subset\Gamma^0$ be the normal subgroup defined analogously to $\Gamma(\lambda)$ so that 
$\Gamma_0/\widetilde{\Gamma}\operatorname{-mod}\xrightarrow{\sim}\overline{\HC}(\cA^{\widetilde{X}})$, and let $\bullet^{\widetilde{\dagger}}, \bullet_{\widetilde{\dagger}}$ denote the equivalences. (ii) is equivalent to the equality $\widetilde{\Gamma}=\Gamma(\lambda)$. 

Every Harish-Chandra bimodule for $\cA^{\widetilde{X}}$ is also a Harish-Chandra bimodule for $\cA^X$, and a good filtration for $\cA^{\widetilde{X}}$ is also a good filtration for $\cA^X$. Set $\widetilde{\mathcal{B}}:=\CC[\Gamma_0/\widetilde{\Gamma}]^{\widetilde{\dagger}}$. Note that $\widetilde{\mathcal{B}}$ comes with a natural good filtration and $\gr(\widetilde{\mathcal{B}})$ embeds into $\CC[\widehat{X}]^{\widetilde{\Gamma}}$
with cokernel supported on $X^{\mathrm{sing}}$. It follows that $\widetilde{\mathcal{B}}_\dagger=\CC[\Gamma/\widetilde{\Gamma}]$, and so $\Gamma(\lambda) \subset \widetilde{\Gamma}$. On the other hand, $\breve{\cA}$ is a Harish-Chandra $\cA^{\widetilde{X}}$-bimodule with 
$\breve{\cA}_{\widetilde{\dagger}}=\CC[\Gamma_0/\Gamma(\lambda)]$. So 
$\Gamma(\lambda) \supset \widetilde{\Gamma}$. This finishes the proof of (ii). 
\end{proof}

Combining Lemma \ref{lem:breve_A} and 
Proposition \ref{Prop:quant_covers} we deduce the following result, which we will use repeatedly. 

\begin{cor}\label{cor:isotypic}
Let $\cA^{\widetilde{X}}$ be a cover of $\cA^X$ and assume $\overline{\HC}(\cA^{\widetilde{X}}) \simeq \Vect$. Then
\begin{itemize}
    \item[(i)] $\widetilde{X} \simeq \breve{X}$
and $\cA^{\widetilde{X}} \simeq \breve{\cA}$. In particular, $\widetilde{X}$ is a Galois cover of $X$ with Galois group $\Gamma/\Gamma(\lambda)$, and the action of $\Gamma/\Gamma(\lambda)$ on $\CC[\widetilde{X}]$ lifts to $\cA^{\widetilde{X}}$.
    \item[(ii)] The equivalence
    $$\bullet^{\dagger}: \Gamma/\Gamma(\lambda)\modd \xrightarrow{\sim} \overline{\HC}(\cA^X)$$
    is given by $V\mapsto (\cA^{\widetilde{X}}\otimes V)^{\Gamma/\Gamma(\lambda)}$.
\end{itemize}
\end{cor}

\section{Normal subgroup $\Gamma(\lambda) \subset \Gamma$}\label{subsec:Gammalambda}

Let $X$ be a conical symplectic singularity and let $\cA_{\lambda}$ be a filtered quantization of $\CC[X]$. In this section, we will give a description of the normal subgroup $\Gamma(\lambda) \subset \Gamma$, following \cite{LosevHC}. 

First, suppose  $X = \CC^2/\Gamma$ is a Kleinian singularity. Let $\fg$ be the complex simple Lie algebra corresponding to $\Gamma$ (cf. Section \ref{subsec:Mckay}). Fix a Cartan subalgebra $\fh \subset \fg$, and let $\Lambda_r$, $\Lambda$, $W^a = \Lambda_r \rtimes W$, and $W^{ae} = \Lambda \rtimes W$ denote the root lattice, weight lattice, affine Weyl group, and extended affine Weyl group, respectively. Let $\Gamma_{\cN}$ denote the fundamental group of the principal nilpotent orbit in $\fg^*$. There is a well-known isomorphism $\Gamma/[\Gamma,\Gamma] \simeq \Gamma_{\cN}$, see e.g. \cite[Sec 5.1]{LosevHC}. Write $\fX(H)$ for the characters of a group $H$ and $G^{\mathrm{sc}}$ for the connected simply connected group corresponding to $\fg$. Then
$$\Gamma_{\cN} \simeq Z(G^{\mathrm{sc}}) \simeq \fX(\Lambda/\Lambda_r) \simeq \fX(W^{ae}/W^a),$$
Hence, there is an isomorphism
\begin{equation}\label{eqn:charsofGammai}
\fX(\Gamma) \simeq W^{ae}/W^a.
\end{equation}
If $\lambda \in \fh^*$, write $W^{ae}_{\lambda}$ (resp. $W^a_{\lambda}$) for the stabilizer of $\lambda$ under the natural action of $W^{ae}$ (resp. $W^a$) on $\fh^*$. 

\begin{prop}[\cite{LosevHC},  Prop 5.3]\label{prop:Gamma0Kleinian}
Let $X$ be a Kleinian singularity. Then $\fX(\Gamma/\Gamma(\lambda)) \subset \fX(\Gamma)$ corresponds to $W^{ae}_{\lambda}/W^a_{\lambda} \subset W^{ae}/W^a$ under the isomorphism (\ref{eqn:charsofGammai}).
\end{prop}

If $\Gamma$ is not of type $E_8$, Proposition \ref{prop:Gamma0Kleinian} can be used to compute $\Gamma(\lambda)$, see \cite[Section 5]{LosevHC}. This result is rather technical and will not be used in {this monograph}, so we omit it.

In some sense, the general case can be reduced to the case of Kleinian singularities. Indeed, let $X$ be an arbitrary conical symplectic singularity, and $\fL_1,...,\fL_t \subset X$ the codimension 2 leaves. For each $\fL_k$, there is a Kleinian singularity $\Sigma_k = \CC^2/\Gamma_k$, see Section \ref{subsec:structurenamikawa}, and a group homomorphism $\phi_k: \Gamma_k \to \Gamma$, see Section \ref{subsec:finitecovers}. The quantization parameter $\lambda \in \fP$ determines an element $\lambda_k \in \fP_k \simeq (\fh_k^*)^{\pi_1(\fL_k)}$ via the restriciton map $\fP \to \fP_k$, see Proposition \ref{prop:partialdecomp}, and this element defines a normal subgroup $\Gamma_k(\lambda_k) \subset \Gamma_k$.

\begin{theorem}[\cite{LosevHC}, Thm 6.1]\label{thm:LosevHC}
$\Gamma(\lambda)$ is the smallest normal subgroup of $\Gamma$ containing $\bigcup_{k=1}^t \phi_k(\Gamma_k(\lambda_k))$. 
\end{theorem}

\chapter{Canonical quantizations of conical symplectic singularities}\label{sec:canonical}

Let $X$ be a conical symplectic singularity, and let $\fP$ denote the associated (complex) Namikawa space. For each $\lambda \in \fP$, write $\cA_{\lambda}$ for the corresponding filtered quantization of $\CC[X]$, see
Section \ref{subsec:quantsymplectic}.

\begin{definition}
The \emph{canonical quantization} of $\CC[X]$ is the filtered quantization $\cA_0$. \index{quantization!canonical}
\end{definition}

Canonical quantizations of nilpotent covers are closely related to our main objects of study---unipotent ideals and unipotent bimodules, to be defined in Chapter \ref{sec:unipotent}. In this chapter, we will establish some general facts about canonical quantizations for use in later chapters. For many of these results, we will impose the following condition on $X$
\begin{equation}\label{eq:noE8condition}X\text{ does not contain a}\text{ 2-dimensional formal slice of type } E_8\end{equation}
For our purposes, this is a relatively harmless assumption. For nilpotent covers, (\ref{eq:noE8condition}) is satisfied in all but one case (the principal orbit in type $E_8$). This isolated exception can be handled separately by classical methods.

\section{Canonical quantizations}\label{subsec:canonical}

First, we study the group $\Aut(\cA_0) $ of filtered algebra automorphisms of $\cA_0$. As in Section \ref{subsec:automorphisms}, let $\mathcal{G}$ denote the reductive part of the group of graded Poisson automorphisms of $\CC[X]$.

\begin{prop}\label{prop:canonicalauts}
The $\mathcal{G}$-action on $\CC[X]$ lifts to a $\mathcal{G}$-action on $\cA_0$ by filtered algebra automorphisms. Furthermore, $\mathcal{G} \subset \Aut(\cA_0)$ is a maximal reductive subgroup.
\end{prop}

\begin{proof}
The action of $\mathcal{G}$ on $\fP/W$ (see Proposition \ref{prop:automorphisms}(i))
stabilizes $0$. The assertions follow from Proposition \ref{prop:automorphisms}(ii). 
\end{proof}

Next, we describe the normal subgroup $\Gamma(0) \subset \Gamma$ defined in Section \ref{subsec:daggers}. Let $\fL_1,...,\fL_t \subset X$ be the codimension 2 leaves. Recall that for each $\fL_k$ there is a group homomorphism $\phi_k : \Gamma_k\to \Gamma$, see (\ref{eq:defofphik}).

\begin{prop}\label{prop:Gamma0}
Suppose $X$ has property (\ref{eq:noE8condition}). Then $\Gamma(0) \subset \Gamma$ is the minimal normal subgroup containing $\phi_k(\Gamma_k)$ for all $k$.
\end{prop}

\begin{proof}
By Theorem \ref{thm:LosevHC}, it suffices to show that $\Gamma_k(0) = \Gamma_k$ for $1 \leq k \leq t$. By Proposition \ref{prop:Gamma0Kleinian}, the characters of $\Gamma_k/\Gamma_k(0)$ are in bijection with $W^{ae}_0/W^a_0 = W/W = 1$. Thus, $\Gamma_k/\Gamma_k(0)$ has no nontrivial characters. Because $X$ satisfies (\ref{eq:noE8condition}), we have that $\Gamma_k$, and hence $\Gamma_k/\Gamma_k(0)$, is solvable. A solvable group with no nontrivial characters is trivial.
\end{proof}

\begin{rmk}
In this remark we explain the special status of the quantization parameter $0\in \fP$.
Let $Y$ be a $\QQ$-factorial terminalization of $X$, and consider the set of isomorphism classes of graded formal quantizations of $Y$ (for definitions, see e.g. \cite[Sec 2.2]{Losev_isofquant}). By \cite[Cor 2.3.3]{Losev_isofquant}, this set is in natural bijection with $\fP$. 
Recall that a graded formal quantization $\mathcal{D}_{\hbar}$ comes with a lattice $\mathcal{D}_{\hbar^d}$ over $\CC[[\hbar^d]]$. We say that $\mathcal{D}_{\hbar}$ is \emph{even}\index{quantization!even} if there is a graded anti-automorphism $\mathcal{D}_\hbar \to \mathcal{D}_\hbar$ which restricts to an involution of $\mathcal{D}_{\hbar^d}$ sending $\hbar^d$ to $-\hbar^d$. It follows from the results of \cite[Sec 2.3]{Losev_isofquant} that the graded formal quantization of $Y$ corresponding to the parameter $0 \in \fP$ is the unique even such. 
\end{rmk}

\section{Invariants in canonical quantizations: case of Kleinian singularities}\label{subsec:invariantsKleinian}

Suppose $\widetilde{X} \to X$ is a finite Galois cover of conical symplectic singularities (cf. Section \ref{subsec:coverings}). The Galois group $\Pi$ acts on $\CC[\widetilde{X}]$ by graded Poisson automorphisms and $\CC[\widetilde{X}]^\Pi\xrightarrow{\sim} \CC[X]$. By Proposition \ref{prop:automorphisms}, the $\Pi$-action lifts to the canonical quantization $\cA_0^{\widetilde{X}}$ of $\CC[\widetilde{X}]$. The algebra of invariants $(\cA_0^{\widetilde{X}})^{\Pi}$ has the structure of a filtered quantiation of $X$. It is important to note that this quantization is typically not canonical. In this section and the next, we will compute its quantization parameter. First, we handle the case of Kleinian singularities.

Let $\Sigma = \CC^2/\Gamma$ be a Kleinian singularity, let $\Gamma'$ be a normal subgroup of $\Gamma$, and consider the induced Galois cover 
$$\Sigma' := \CC^2/\Gamma' \twoheadrightarrow \CC^2/\Gamma = \Sigma$$
The Galois group $\Pi := \Gamma/\Gamma'$ acts on $\CC[\Sigma']$ by graded Poisson automorphisms and $\CC[\Sigma']^{\Pi} \xrightarrow{\sim} \CC[\Sigma]$. 
We will compute the quantization parameter of $(\cA_0^{\Sigma'})^{\Pi}$, this is a special case of the situation of the previous paragraph. 

Let $\fg$ be the complex simple Lie algebra corresponding to $\Gamma$, see Section \ref{subsec:Mckay}. Choose a Cartan subalgebra $\fh \subset \fg$. Let $\{V_1,...,V_n\}$ be the nontrivial irreducible representations of $\Gamma$ and let $\{\omega_1,...,\omega_n\} \subset \fh^*$ be the corresponding fundamental weights. Define a dominant weight $\epsilon(\Gamma') \in \fh^*$ by the formula
\begin{equation}\label{eq:defofepsilonGamma}
\epsilon(\Gamma') := \frac{|\Gamma'|}{|\Gamma|} \sum_{i=1}^n \dim (V_i^{\Gamma'}) \ \omega_i.\end{equation}

\begin{prop}\label{prop:parameterofinvariants}
There is an isomorphism of filtered quantizations
$$(\cA_0^{\Sigma'})^{\Pi} \simeq \cA_{\epsilon(\Gamma')}^{\Sigma}.$$
\end{prop}

\begin{proof}
{Let $H = \CC\langle x,y\rangle\#\Gamma/(xy-yx-|\Gamma'|e')$ and $H'=\CC \langle x,y\rangle\#\Gamma'/(xy-yx-|\Gamma'|e')$}, where $e' \in\CC[\Gamma']$ is the averaging idempotent. Note that $H'$ embeds as a subalgebra in $H$ and $\Pi$ acts on $e'H'e'$ by algebra automorphisms. If  $z\in (e'H'e')^{\Pi}$, then $ze\in eHe$, and the map $z \mapsto ze$ defines a filtered algebra homomorphism $(e'H'e')^{\Pi} \to eHe$. Its associated graded homomorphism is the natural identification
$$\CC[\Sigma']^{\Pi}\xrightarrow{\sim} \CC[\Sigma].$$
Thus, $(e'H'e')^{\Pi} \simeq eHe$ as filtered quantizations. Now we have isomorphisms of filtered quantizations
$$(\cA_0^{\Sigma'})^{\Pi} \simeq (e'H'e')^{\Pi} \simeq eHe \simeq \cA_{\lambda^{|\Gamma'|e'}}^{\Sigma},$$
where $\lambda^{|\Gamma'|e'} \in \fh^*$ is the parameter defined in Section \ref{subsec:quantKleinian}. The first isomorphism is by (\ref{eq:lambdaGammae}), and the third is by (\ref{eq:lambdatoc}). It remains to show that $\lambda^{|\Gamma'|e'} = \epsilon(\Gamma')$. If $V'$ is an irreducible representation of $\Gamma'$, then
$$\mathrm{tr}_{V'}(e') = \begin{cases} 
      1 & V' \simeq \CC \\
      0 & \text{otherwise} 
   \end{cases}$$
Hence, 
$$\mathrm{tr}_{V_i}(|\Gamma'|e') = |\Gamma'| \dim (V_i^{\Gamma'}),  \qquad  1 \leq i \leq n,$$
and therefore
$$\langle \lambda^{|\Gamma'|e'}, \alpha_i^{\vee}\rangle = \frac{1}{|\Gamma|}\mathrm{tr}_{V_i}(|\Gamma'|e') = \frac{|\Gamma'|}{|\Gamma|}\dim (V_i^{\Gamma'})  = \langle \epsilon(\Gamma'), \alpha_i^{\vee}\rangle, \qquad 1 \leq i \leq n.$$
This completes the proof.
\end{proof}

\begin{example}\label{ex:barycenterKleinian}
Suppose $\Sigma' = \CC^2$, so that $\Gamma'=1$. Then $\epsilon(1)$ is the weighted barycenter parameter $\lambda^1 \in \fh^*$ defined in Section \ref{subsec:quantKleinian}:
$$\epsilon(1) = \lambda^1= \frac{1}{|\Gamma|}\sum_{i=1}^n \dim(V_i) \  \omega_i \in \fh^*.$$
\end{example}

\section{Invariants in canonical quantizations: general case}\label{sec:invariantssymplectic}

Fix $p:\widetilde{X} \to X$ and $\Pi$ as in the first paragraph of Section \ref{subsec:invariantsKleinian}. In this section, we will compute the quantization parameter of $(\cA_0^{\widetilde{X}})^{\Pi}$.

Let $\fL_k \subset X$ be a codimension 2 leaf and let $x \in \fL_k$. Fix the notation of Section \ref{subsec:finitecovers}, i.e. $\Sigma_k = \CC^2/\Gamma_k$, $\widetilde{x}_j$, $\Sigma_k'=\CC^2/\Gamma_k'$, $\phi_k: \Gamma_k \to \Gamma$ and so on. Let $\fg_k$ be the complex simple Lie algebra corresponding to $\Gamma_k$, and let $\fh_k \subset \fg_k$ be a Cartan subalgebra. Let $\{V_1(k),...,V_{n(k)}(k)\}$ be the nontrivial irreducible representations of $\Gamma_k$ and $\{\omega_1(k),...,\omega_{n(k)}(k)\} \subset \fh_k^*$ the corresponding fundamental weights. Define the element
\begin{equation}\label{eq:defofepsilonk}
\epsilon_k := \epsilon(\Gamma_k') = \frac{|\Gamma_k'|}{|\Gamma_k|} \sum_{i=1}^{n(k)} \dim (V_i(k)^{\Gamma_k'}) \omega_i(k) \in \fh_k^*.\end{equation}
\begin{prop}\label{prop:parameterofinvariantssymplectic}
The element $\epsilon_k \in \fh_k^*$ is a fixed point for the monodromy action of $\pi_1(\fL_k)$ on $\fh_k^*$, and hence an element of $\fP_k^X := (\fh_k^*)^{\pi_1(\fL_k)}$. Define the quantization parameter
$$\epsilon:= (0,\epsilon_1,\epsilon_2,...,\epsilon_t) \in \bigoplus_{k=0}^t \fP_k^X \simeq \fP^X.$$
There is an isomorphism of filtered quantizations
$$(\cA_0^{\widetilde{X}})^{\Pi} \simeq \cA_{\epsilon}^X.$$
\end{prop}

\begin{proof}
By Theorem \ref{thm:quantsofsymplectic}, there is a parameter $\lambda \in \fP^X$ (well-defined modulo $W$) such that
$$(\cA_0^{\widetilde{X}})^{\Pi} \simeq \cA_{\lambda}^X.$$
We will show that $\lambda_0 =0$ and $W_k\lambda_k=W_k\epsilon_k$ for $ 1\leq k \leq t$; this will imply the claims of the proposition. 

{\it Step 1}.
First, we show that $\lambda_0 = 0$. Fix a $\QQ$-terminalization $\rho: \widetilde{Y} \to \widetilde{X}$. The identification $\mathrm{Quant}(\widetilde{Y}) \xrightarrow{\sim} \fP^{\widetilde{X}}$ is the composition
$$\mathrm{Quant}(\widetilde{Y}) \to \mathrm{Quant}(\widetilde{Y}^{\mathrm{reg}}) \overset{\mathrm{Per}}{\to} \fP^{\widetilde{X}}$$
where the first map is restriction, see \cref{prop:P=Quant(Y)}. Let $\mathcal{D}^{\widetilde{Y}}$ denote the quantization of $\widetilde{Y}$
with parameter $0$ so that $\cA_0^{\widetilde{X}} = \Gamma(Y,\mathcal{D}^{\widetilde{Y}})$. {Let  $X'$ denote the preimage of $X^{\mathrm{reg}}$ under $p: \widetilde{X} \to X$. Let $\mathcal{D}^{\widetilde{Y}^{\mathrm{reg}}}$, $\mathcal{D}^{X'}$ denote the restrictions of $\mathcal{D}^{\widetilde{Y}}$ to $\widetilde{Y}^{\mathrm{reg}}$ and $X'$ and let $\mathcal{D}^{X^{\mathrm{reg}}} = (\mathcal{D}^{X'})^{\Pi}$. Note that $\mathcal{D}^{X^{\mathrm{reg}}}$ is a quantization of $X^{\mathrm{reg}}$ and $\mathcal{D}^{X'}$ is its (\'{e}tale) lift to $X'$. Since the period map is functorial, see \cite[Sec 4]{BK}, $\mathrm{Per}(\mathcal{D}^{X'}) \in H^2(X',\CC)$ is the image of $\mathrm{Per}(\mathcal{D}^{\widetilde{Y}^{\mathrm{reg}}}) = 0 \in H^2(\widetilde{Y}^{\mathrm{reg}},\CC)$ under the restriction map $H^2(\widetilde{Y}^{\mathrm{reg}},\CC) \to H^2(X',\CC)$ and of  $\mathrm{Per}(\mathcal{D}^{X^{\mathrm{reg}}}) \in H^2(X^{\mathrm{reg}},\CC)$ under the pullback map $H^2(X^{\mathrm{reg}},\CC) \hookrightarrow H^2(X',\CC)$. Thus, $\mathrm{Per}(\mathcal{D}^{X^{\mathrm{reg}}})=0$.
Since $\mathcal{D}^{X^{\mathrm{reg}}}$ is the restriction of $(\cA_0^{\widetilde{X}})^{\Pi}$ to $X^{\mathrm{reg}}$, this implies $\lambda_0=0$}. In steps 2 and 3, we will treat the case of $k>0$.

{\it Step 2}. Consider the codimension 2 leaf $\fL_k \subset X$ and the fiber product (\ref{diag:fiberproduct}). As noted there, the natural map $\bigsqcup_j \widetilde{X}^{\wedge \widetilde{x}_j} \to \widetilde{X}$ is $\Pi$-equivariant. For each $\widetilde{x}_j$ let $\fL'_j$ denote the symplectic leaf in $\widetilde{X}$ containing $\widetilde{x}_j$. There is an isomorphism
\begin{equation*}\label{eq:decompositionleafcover}
    \cA_{0,\hbar}^{\widetilde{X},\wedge}\simeq \AA_{\hbar}(V_j)^{\wedge} \widehat{\otimes}_{\CC[[\hbar]]} \cA_{0,\hbar}^{{\Sigma}',\wedge},
\end{equation*}
where $V_j=T_{\widetilde{x}_j}\fL'_j$, see (\ref{eq:decompositionleaf}). Here $\cA_{0,\hbar}^{\widetilde{X},\wedge}$ (resp. $\cA_{0,\hbar}^{{\Sigma}',\wedge}$) denotes the completion of the Rees algebra of the canonical quantization $\cA_0^{\widetilde{X}}$ (resp. $\cA_0^{{\Sigma}'}$) at $\widetilde{x}_j \in \widetilde{X}$ (resp. $0 \in {\Sigma}'$). Since taking invariants commutes with completions, there is a natural isomorphism
\begin{equation}\label{eq:iso1}
\cA_{\lambda, \hbar}^{X,\wedge_x}\xrightarrow{\sim} \left(\bigoplus \AA_{\hbar}(V_j)^{\wedge} \widehat{\otimes}_{\CC[[\hbar]]}  \cA_{0,\hbar}^{{\Sigma}',\wedge}\right)^{\Pi}
\end{equation}
%
%
Fix a preimage $\widetilde{x}$ of $x$ in $\widetilde{X}$ and let $\Pi':=\Pi_{\widetilde{x}}, V:=V_j$. Since $\Pi$ acts transitively on $\{\widetilde{x}_j\}$, we have
$$\left(\bigoplus \AA_{\hbar}(V_j)^{\wedge} \widehat{\otimes}_{\CC[[\hbar]]}  \cA_{0,\hbar}^{{\Sigma}',\wedge}\right)^{\Pi}\xrightarrow{\sim} \left(\AA_{\hbar}(V)^{\wedge} \widehat{\otimes}_{\CC[[\hbar]]}  \cA_{0,\hbar}^{{\Sigma}',\wedge}\right)^{\Pi'}.$$ 
and therefore
\begin{equation}\label{eq:iso2}
\cA_{\lambda, \hbar}^{X,\wedge_x}\xrightarrow{\sim} \left(\AA_{\hbar}(V)^{\wedge} \widehat{\otimes}_{\CC[[\hbar]]}  \cA_{0,\hbar}^{{\Sigma}',\wedge}\right)^{\Pi'}.\end{equation}

{\it Step 3}. Since $\mathcal{L}_k$ is a smooth and symplectic and $\mathcal{L}_k=p^{-1}(\mathcal{L}_k)/\Pi$, $\Pi'$ acts trivially on $V$ (otherwise, $V/\Pi'$ would be singular because the action of $\Pi'$ on $V$ is symplectic, while $\mathcal{L}_k$ would have singularity $V/\Pi'$ at $x$). 
In particular, there is an identification  $T_x\mathcal{L}\simeq V$ and also 
\begin{equation}\label{eq:classical_iso3}
    \CC[X]^{\wedge_x}\xrightarrow{\sim} 
    \CC[V]^\wedge\widehat{\otimes}
    \CC[\Sigma'_k]^\wedge.
\end{equation}
Furthermore
$$(\AA_{\hbar}(V)^{\wedge} \widehat{\otimes}_{\CC[[\hbar]]}  \cA_{0,\hbar}^{{\Sigma}',\wedge})^{\Pi'}=\AA_{\hbar}(V)^{\wedge} \widehat{\otimes}_{\CC[[\hbar]]}  (\cA_{0,\hbar}^{{\Sigma}',\wedge})^{\Pi'}.$$
So (\ref{eq:iso2}) becomes 
\begin{equation}\label{eq:iso3}
\cA_{\lambda, \hbar}^{X,\wedge_x}\xrightarrow{\sim} \AA_{\hbar}(V)^{\wedge} \widehat{\otimes}_{\CC[[\hbar]]}  (\cA_{0,\hbar}^{{\Sigma}',\wedge})^{\Pi'}.\end{equation}
This isomorphism lifts (\ref{eq:classical_iso3}).
As noted at the end of Section \ref{subsec:quantsymplectic}, the isomorphism (\ref{eq:iso3}) determines $W_k\lambda_k$ uniquely. By Proposition 
\ref{prop:parameterofinvariants}, we have $W_k\lambda_k=W_k\epsilon_k$.

\end{proof}

\begin{example}\label{ex:barycentersymplectic}
Suppose $\widetilde{X}$ has no codimension 2 leaves. Then for each $\Sigma_k$, the quantization parameter $\epsilon_k \in \fh_k^*$ (cf. (\ref{eq:defofepsilonk})) is the `weighted barycenter' parameter considered in Example \ref{ex:barycenterKleinian}
$$\epsilon_k = \frac{1}{|\Gamma_k|} \sum_{i=1}^{n(k)} \dim V_i(k) \ \omega_i(k) \in \fh_k^*.$$
We call the parameter
$$\epsilon = (0,\epsilon_1,\epsilon_2,...,\epsilon_t) \in \bigoplus_{k=0}^t \fP_k^X \simeq \fP^X$$
the \emph{weighted barycenter}\index{barycenter parameter} of $\fP^X$. By Proposition \ref{prop:parameterofinvariantssymplectic}, $\epsilon$ is the quantization parameter of the $\Pi$-invariants in the canonical quantization $\cA_0^{\widetilde{X}}$. This parameter will play an important role in Chapter \ref{sec:centralchars} in the computation of infinitesimal characters.
\end{example}

\section{Almost \'{e}tale covers of symplectic singularities}\label{subsec:almostetale}

In this section, we will define and study finite covers of a very special type. The ideas in this section will be used in Chapter \ref{sec:unipotent} for the classification of unipotent ideals and bimodules. Let $p: \widetilde{X} \to X$ be a finite cover and fix the notation of Section \ref{subsec:finitecovers} (e.g. $\Gamma$, $\Gamma'$, $\Gamma_k$, $\Gamma_k'$ and so on).

\begin{prop}\label{prop:almostetale}
The following conditions are equivalent:
\begin{itemize}
    \item[(i)] $p$ is \'{e}tale over the open subset
    $$X^{\mathrm{reg}} \cup \bigcup_{k=1}^t \fL_k \subset X.$$
    \item[(ii)] For each codimension 2 leaf $\fL_k \subset X$, there is an equality $\Gamma_k'=\Gamma_k$.
\end{itemize}
If either of these equivalent conditions is satisfied, we say that $p:\widetilde{X} \to X$ is \emph{almost \'{e}tale}.\index{cover!almost-\'{e}tale}
\end{prop}

\begin{proof}
First, we prove (i) $\Rightarrow$ (ii). Choose $x \in \fL_k$ and $\widetilde{x} \in p^{-1}(x)$. By assumption, $p$ is \'{e}tale at $\widetilde{x}$. Thus $p$ induces an \'{e}tale map (and hence an isomorphism) $\widetilde{X}^{\wedge \widetilde{x}}\to X^{\wedge x}$. Consider the symplectic vector space $V:=\CC^2\oplus T_x\mathcal{L}_k$. There are identifications $\widetilde{X}^{\wedge \widetilde{x}}\simeq V^{\wedge_0}/\Gamma_k'$ and 
$X^{\wedge x}\simeq V^{\wedge_0}/\Gamma_k$. In particular, $\Gamma_k,\Gamma_k'$ are recovered as the algebraic fundamental groups of the regular loci of $X^{\wedge x}$ 
and $\widetilde{X}^{\wedge \widetilde{x}}$, respectively. Since 
$\widetilde{X}^{\wedge \widetilde{x}}\xrightarrow{\sim} X^{\wedge x}$, we conclude that $\Gamma_k=\Gamma_k'$.

Next, we prove (ii)$\Rightarrow$(i).
A neighborhood of $x\in X$ in the analytic topology is identified with $D/\Gamma_k$ for a disc $D \subset V$
around $0$. Similarly, a neighborhood
of $\widetilde{x}$ in $\widetilde{X}$ is identified with $\widetilde{D}/\Gamma_k$. Shrinking $D,\widetilde{D}$ if necessary  we can assume that $p$ restricts to  $\widetilde{D}/\Gamma_k\rightarrow D/\Gamma_k$. This restriction is a quasi-finite morphism of complex analytic spaces. We claim it is \'{e}tale. Let $D^0:=D- (\{0\} \times T_x\mathcal{L}_k)$. Note that $D^0$ is the locus in $D$ where $\Gamma_k$ acts freely. Define $\ \widetilde{D}^0$ similarly. Let $\widetilde{D}^1\subset \widetilde{D}^0$ be the preimage of $D^0/\Gamma_k$ in $\widetilde{D}$. 
The complement to $\widetilde{D}^1$ in $\widetilde{D}$ is of codimension $2$ because $p$ is quasi-finite. In particular, $\widetilde{D}^1$ is simply connected. The map $\widetilde{D}^1/\Gamma_k\rightarrow D^0/\Gamma_k$ is a Poisson morphism of symplectic complex analytic manifolds, and is therefore \'{e}tale. Thus it lifts to 
a  Poisson morphism of universal covers, which coincide with $\widetilde{D}^1$ and $D^0$, respectively. By the Hartogs theorem, the composition $\widetilde{D}^1\rightarrow D^0\hookrightarrow D$ extends to $\widetilde{D}$. This extension is \'{e}tale outside of codimension $2$, and hence \'{e}tale. This implies that $\widetilde{D}/\Gamma_k\rightarrow D/\Gamma_k$ is \'{e}tale, as desired.
\end{proof}

If $p: \widetilde{X} \to X$ is Galois, there is a third characterization of almost \'{e}tale. Recall the group homomorphisms $\phi_k: \Gamma_k \to \Gamma$, see (\ref{eq:defofphik}).

\begin{prop}\label{prop:almost etale through Gamma}
Suppose  $p: \widetilde{X} \to X$ is Galois. Then $p$ is almost \'{e}tale if and only if $\Gamma'\supset \phi_k(\Gamma_k)$ for all $k$.
\end{prop}

\begin{proof}
\vspace{3mm}
First, suppose $p$ is almost \'{e}tale. For each codimension 2 leaf $\fL_k \subset X$, there is an inclusion $\phi'_k(\Gamma'_k) \subseteq \phi_k(\Gamma_k) \cap \Gamma'$, see (\ref{diag:galoiscovering}). By Proposition \ref{prop:almostetale}, $\Gamma_k' = \Gamma_k$, and so $\phi_k(\Gamma_k) \subset \Gamma'$. 

Conversely, suppose $\phi_k(\Gamma_k) \subset \Gamma'$. 
Let $\mathfrak{L}_k \subset X$ be a codimension 2 leaf and let $x\in \mathfrak{L}_k$. Choose a small contractible open neighborhood $D$ of $x$ in $X$. We wish to show that $p:p^{-1}(D)\rightarrow D$ is unramified. 
Let $\widetilde{D}$ be a connected component of $p^{-1}(D)$, and choose a base point $y\in D - \mathfrak{L}$ for $\pi_1(D - \mathfrak{L})$. Any loop in $D- \mathfrak{L}$ is conjugate to an element in the preimage of $\phi_k(\Gamma_k)\subset \Gamma=\pi_1^{\mathrm{alg}}(X^{\mathrm{reg}})$ in $\pi_1(X^{\mathrm{reg}})$ by the construction of $\phi_k$. Therefore it lies in the preimage of $\Gamma'$ in $\pi_1(X^{\mathrm{reg}})$. It follows that the cover $\widetilde{D}- p^{-1}(\mathfrak{L})
\rightarrow D- \mathfrak{L}$ is an isomorphism. Hence the restriction of $p$ to $\widetilde{D}$ is an isomorphism. This completes the proof. 
\end{proof}

If $\widetilde{X},\widetilde{X}' \to X$ are finite covers of $X$, a \emph{morphism of covers} is a map $f: \widetilde{X} \to \widetilde{X}'$ such that the following diagram commutes
\begin{center}
    \begin{tikzcd}
    \widetilde{X} \ar[dr] \ar[r,"f"] & \widetilde{X}' \ar[d]\\
    & X
    \end{tikzcd}
\end{center}
This defines a partial order on the set of isomorphism classes of finite covers of $X$, analogously to (\ref{eq:partialordercovers}).

\begin{prop}\label{prop:maximaletale}
There is a unique \emph{maximal} almost \'{e}tale cover  $p:\breve{X} \to X$. This cover is Galois. The corresponding Galois group is $\Gamma/\underline{\Gamma}$, where $\underline{\Gamma}$ is the minimal normal subgroup of $\Gamma$ containing $\phi_k(\Gamma_k)$ for all $k$.
\end{prop}

\begin{proof}
First, we will show that $X$ admits a maximal almost \'{e}tale cover. Since $\Gamma$ is finite (and the identity map $X \to X$ is almost \'{e}tale), it suffices to show that for any pair of almost \'{e}tale covers $\widetilde{X}, \widetilde{X}' \to X$, there is an almost \'{e}tale cover $\breve{X} \to X$ which covers both $\widetilde{X}$ and $\widetilde{X}'$. First, form the fiber product $\widetilde{X}\times_X\widetilde{X}'$. The natural morphisms $\widetilde{X}\times_X \widetilde{X}'\rightarrow X,\widetilde{X},\widetilde{X}'$ are finite and \'{e}tale outside of codimension 4. Note that $\widetilde{X}\times_X\widetilde{X}'$
may fail to be irreducible. Let $\breve{X}$
denote the normalization of an irreducible component. Then $\breve{X}$ is an almost \'{e}tale cover of $X$,$\widetilde{X}$, and $\widetilde{X}'$. This proves both the existence and uniqueness of the maximal almost \'{e}tale cover.

It remains to show that the unique maximal almost \'{e}tale cover $\breve{X} \to X$ is Galois. The claim regarding the Galois group then follows from  Proposition 
\ref{prop:almost etale through Gamma}. {The argument that $\breve{X}$ is Galois is similar to the previous paragraph. By the maximality of $\breve{X}$, any irreducible component of $\breve{X}\times_X \breve{X}$ maps to either copy of $\breve{X}$ (generically) isomorphically. 
Let $\widehat{X}$ be the universal algebraic cover of $X$ and $H\subset \Gamma:=\pi_1^{alg}(X^{reg})$ be such that $\breve{X}=\widehat{X}/H$. 
Note that for any $\gamma\in \Gamma$, the variety $\breve{X}\times_X\breve{X}$ contains a component 
whose normalization is $\widehat{X}/(H\cap \gamma H\gamma^{-1})$. This shows that $H$ is normal and finishes the proof.}
\end{proof}

\begin{rmk}\label{rmk:maximaletale}
Note that if $X$ satisfies (\ref{eq:noE8condition}), then $\Gamma(0) = \underline{\Gamma}$ by 
Proposition \ref{prop:Gamma0}. Thus, the Galois group of the cover $\breve{X} \to X$ is $\Gamma/\Gamma(0)$.
\end{rmk}

The next proposition relates almost \'{e}tale covers to canonical quantizations.

\begin{prop}\label{prop:almostetaleinvariants}
Let $p:\widetilde{X} \to X$ be a finite Galois almost \'{e}tale cover, and let $\Pi = \Gamma/\Gamma'$ be its Galois group. There is an isomorphism of filtered quantizations
$$(\cA_0^{\widetilde{X}})^{\Pi} \simeq \cA_0^X.$$
\end{prop}

\begin{proof}
By Proposition \ref{prop:parameterofinvariantssymplectic}, there is an isomorphism of filtered quantizations
$$(\cA_0^{\widetilde{X}})^{\Pi} \simeq \cA_{\epsilon}^X.$$
where $\epsilon \in \fP^X$ is the quantization parameter defined in the statement of that proposition. By Proposition \ref{prop:almostetale}, we have $\Gamma_k' = \Gamma_k$ for $1 \leq k \leq t$. Hence
$$\epsilon_k = \epsilon(\Gamma_k')= \frac{|\Gamma_k'|}{|\Gamma_k|} \sum_{i=1}^{n(k)} \dim (V_i(k)^{\Gamma_k}) \ \omega_i(k) = 0.$$
\end{proof}

\section{Simplicity conjecture}\label{subsec:conjecturesimplicity}

In this section, we gather evidence for the following general conjecture:

\begin{conj}\label{conj:simplicity}
Suppose $X$ is a conical symplectic singularity. Then the canonical quantization $\cA_0$ of $\CC[X]$ is a simple algebra. 
\end{conj}

Note that if $X$ is $\QQ$-factorial and terminal, this conjecture is a special case of \cite[Conj 3.1]{LosevSRA}.

\begin{example}\label{ex:covers}
Let $G$ be a complex connected reductive algebraic group, let $\widetilde{\mathbb{O}}$ be a $G$-equivariant nilpotent cover, and let $X = \Spec(\CC[\widetilde{\mathbb{O}}])$. In Chapter \ref{sec:unipotent}, we consider the primitive ideal $I_0(\widetilde{\mathbb{O}}) = \ker{\Phi} \subset U(\fg)$ (where $\Phi: U(\fg) \to \cA_0$ is the quantum co-moment map for $\cA_0$). Such ideals are called `unipotent' and play a central role in {this monograph}. By Proposition \ref{prop:simplemaximal} below, the simplicity of $\cA_0$ is equivalent to the maximality of $I_0(\widetilde{\mathbb{O}})$. In Appendix \ref{sec:maximality}, we show that $I_0(\widetilde{\mathbb{O}})$ is always maximal if $G$ is linear classical. This argument is generalized to arbitrary groups in the paper \cite{MBM}. 
\end{example}

The other cases where we know Conjecture \ref{conj:simplicity} to 
be true are not otherwise relevant to {this monograph}, so we will be brief. 

\begin{example}\label{ex:slodowy_cover}
Fix the notation of Example \ref{ex:covers}. Let $\mathbb{O}'\subset \partial\mathbb{O}$ 
be an orbit, and $\bullet_{\dagger'}$ be the restriction functor, see Section \ref{subsec:W}. Let $S'$ be the Slodowy slice to $\mathbb{O}'$ and $X'$ be the preimage of $S'$
in $X$. If $\cA_0$ is the canonical quantization of $X$, one can show that $(\cA_0)_{\dagger'}$ is
the canonical quantization of $X'$. The simplicity of $(\cA_0)_{\dagger'}$ should 
follow from the simplicity of $\cA_0$ (together with a transitivity property 
for restriction functors).  
\end{example}

\begin{example}\label{ex:quiver}
Let $X$ be an affine Nakajima quiver variety of the form $\mathcal{M}(v,w)$
in the notation of \cite[Sec 2.1]{BezLosev} (slightly modified: to each loop we assign the space of traceless matrices). The formal slices to symplectic leaves in 
$X$ are of the same form for ``smaller'' quivers. The quantization $\cA_0$
is a simple algebra if and only if none of the quantum slice algebras (which are also
canonical quantizations) have finite dimensional representations. So
Conjecture \ref{conj:simplicity} in this case is equivalent to the assertion that
the canonical quantization of a Nakajima quiver variety has no 
finite dimensional representations. If the underlying quiver $Q$ is 
of finite or affine type, this follows from the results of \cite{BezLosev}, \cite{Losevtotallyaspherical}, \cite{Losevgieseker}, and \cite{Losevwallcrossing}.
First, assume that $Q$ has no loops. 
The argument of \cite[Section 3]{Losevwallcrossing} can be used to prove
that if \cite[Conjecture 1.1]{BezLosev} holds (in the stronger form which also 
includes the claim that the map $\mathsf{CC}$ from there is injective),
then $\cA_0$ has no finite dimensional representations. If
$Q$ is of finite or affine type, the stronger form of \cite[Conjecture 1.1]{BezLosev}
was proved in \cite{BezLosev},\cite{Losevwallcrossing}. And in the case when $Q$ is a Jordan 
quiver, $\cA_0$ has no finite dimensional representations by 
\cite{Losevgieseker}. Thus Conjecture \ref{conj:simplicity}
holds when $Q$ is of finite or affine type. 
\end{example}

\begin{example}\label{ex:quotient_sing}
Let $X=V/\Gamma$, where $V$ is a symplectic vector space and $\Gamma$ is a finite 
group of linear symplectomorphisms of $V$. Consider the normal subgroup $\Gamma' \subset \Gamma$ generated by all symplectic reflections. Let $X_0=V/\Gamma'$ and let $\cA_0^{X_0}$ be the canonical
quantization of $V/\Gamma'$. Then $\cA_0=(\cA_0^{X_0})^{\Gamma/\Gamma'}$ by Proposition \ref{prop:almostetaleinvariants}. It is easy to see that the simplicity of $\cA_0^{X_0}$ implies that of
$\cA_0$. So for Conjecture \ref{conj:simplicity} we can assume 
that $\Gamma$ is generated by symplectic reflections. An important case
is when $V=(\CC^2)^{\oplus n},\Gamma=S_n\ltimes \Gamma_1^n$, where $\Gamma_1\subset 
\mathrm{SL}(2)$. The variety $V/\Gamma$ in this case is a Nakajima 
quiver variety as in Example \ref{ex:quiver}. So in this case Conjecture 
\ref{conj:simplicity} should follow. Another infinite family of examples are
complex reflection groups $G(\ell,r,n)$ (for $r=1$ we recover $S_n\ltimes \Gamma_1^n$,
where $\Gamma_1\simeq \mathbb{Z}_{\ell}$). In this case, the conjecture should also 
follow from the approach of \cite[Section 3]{Losevtotallyaspherical}. This leaves only exceptional groups generated by symplectic reflections. 
\end{example}

\section{Harish-Chandra bimodules for different parameters}\label{SS_HC_different} 
Let $X$ be a conical symplectic singularity. Let $\Gamma$ denote the algebraic fundamental group of $X^{reg}$.
For $\lambda\in \fP/W$, we write $\A_\lambda$ for the corresponding quantization of $X$.
We say that a quantization parameter $\lambda\in \fP/W$ is {\it unipotent} if it is the parameter of the quantization of 
the form $\tilde{\A}_0^{\underline{\Gamma}}$, where $\underline{\Gamma}$ is a quotient of $\Gamma$ and $\tilde{\A}_0$ is the 
quantization of the Galois cover $\tilde{X}$ of $X$ with Galois group $\underline{\Gamma}$ corresponding to quantization parameter $0$. The terminology ``unipotent'' will be justified later. So we have assigned a parameter, $\lambda$, to each quotient of $\Gamma$.

Consider the algebra $\A_{\fP}:=\Gamma(\calD^{Y,\mathrm{univ}})$, where $\calD^{Y,\mathrm{univ}}$ appeared in Proposition \ref{prop:universal for Q-term} and $\fP$ is the Cartan space for $X$, so that $\A_{\lambda}=\A_{\fP}\otimes_{\C[\fP]}\C_\lambda$. Thanks to Definition \ref{def:HCbimodfiltered}, we can talk about HC $\A_{\fP}$-bimodules. Let $\HC(\A_{\lambda^1},\A_{\lambda^2})$
denote the full subcategory in $\HC(\A_{\fP})$,
where the left action of $\CC[\fP]$ is via $\lambda^1$ and the right action is via $\lambda^2$, its objects are called HC $\A_{\lambda^1}$-$\A_{\lambda^2}$-bimodules. It makes sense to speak about the associated variety of such a bimodule in $X$, in particular, about HC bimodules with full support.

The main result of this section is the following proposition. It will play an important role in classifying unipotent bimodules in Section \ref{subsec:classificationbimods}.  

\begin{prop}\label{Prop:diff_unip_HC}
Suppose $X$ satisfies (\ref{eq:noE8condition}).
Let $\lambda^1,\lambda^2$ be two unipotent parameters such that there is a HC $\A_{\lambda^1}$-$\A_{\lambda^2}$-bimodule with full support.
Then $\lambda^1=\lambda^2$.
\end{prop}
\begin{proof}
The proof is in several steps. 

{\it Step 1}. We reduce the claim to the case when $X$ is a Kleinian singularity. Recall, Section \ref{subsec:quantsymplectic}, that any quantization parameter $\lambda$ is of the form 
$(\lambda_0,\lambda_1,\ldots,\lambda_r)$, where $\lambda_0\in H^2(X^{reg},\mathbb{C})$, while the remaining components $\lambda_i$ are quantization
parameters for the slice Kleinian singularities $\Sigma_i$. Observe that if $\lambda$ is unipotent, then $\lambda_0=0$, while the slice quantization
$\A^i_{\lambda_i}$ of $\Sigma_i$ is obtained as the invariants in the canonical quantization of some cover of $\Sigma_i$, cf. Section \ref{subsec:finitecovers}. Hence $\lambda_i$
is unipotent. 

Now let $\mathcal{B}$ be a HC $\A_{\lambda^1}$-$\A_{\lambda^2}$-bimodule with full support. Let $x\in \fL_i$. We can form its restriction $\cB_{\dagger,x}$ to $\Sigma_i$, see \cite[Section 3.3]{Losevwallcrossing}, and it is an $\cA^{\Sigma_i}_{\lambda^1_i}-\cA^{\Sigma_i}_{\lambda^2_i}$-bimodule with full support, see \cite[Lemma 3.5]{Losevwallcrossing}. This completes the reduction to the Kleinian singularity case.    

{\it Step 2}. Until the further notice, $X=\C^2/\Gamma$ and $\Gamma^i:=\Gamma(\lambda^i)$ is the normal subgroup assigned to the parameter $\lambda^i$, see Section \ref{subsec:Gammalambda}. In \cite[Theorem 3.4.5]{sraco}, the first named author constructed an enhanced restriction functor from the category of Harish-Chandra bimodules over 
a suitable version of the quantization of $\C[\C^2]^\Gamma$. There are slight differences: \cite{sraco} deals with the full symplectic reflection algebra, not its spherical subalgebra, and with the graded situation. The filtered setting follows from the graded one, while the case of spherical subalgebras is completely parallel to the case of full symplectic reflection algebras.

We can consider the category of HC $\C[\fP]$-bimodules and its $\Gamma$-equivariant version, the latter will be denoted by $\HC^\Gamma(\C[\fP])$. 

Now consider the open leaf in $\C^2/\Gamma$. There is the restriction functor $\bullet_\dagger:\HC(\A_{\fP})\rightarrow \HC^\Gamma(\C[\fP])$
with the following properties:
\begin{enumerate}
\item It is monoidal; this is part of \cite[Theorem 3.4.5]{sraco}.
\item It intertwines the internal Hom functors: $\Hom_{\A_{\fP}}(M,N)_\dagger\xrightarrow{\sim} \Hom_{\C[\fP]}(M_\dagger,N_\dagger)$, and the same for Hom's over $\A_\fP^{opp}$. This is completely parallel to (1). 
\item It is $\C[\fP]$-bilinear, this follows immediately from the construction.
\item On the full subcategory $\HC(\A_\lambda)$, the functor $\bullet_\dagger$ becomes the restriction functor $\HC(\A_\lambda)\rightarrow \operatorname{Rep}(\Gamma)$ from \cite{LosevHC}, recalled in Section \ref{subsec:daggers}. This is because, for symplectic quotient singularities, the 
    functor from \cite{LosevHC} is a special case of the functor from \cite{sraco}. 
\item The functor $\operatorname{HC}(\A_{\lambda},\A_{\lambda'})\rightarrow \operatorname{Rep}(\Gamma)$  factors through
a full embedding from $\overline{\operatorname{HC}}(\A_{\lambda},\A_{\lambda'})$ for all quantization parameters $\lambda,\lambda'$.
This is a part of \cite[Theorem 3.4.6]{sraco}.
\end{enumerate}  

{\it Step 3}. Combining Proposition \ref{prop:Gamma0} with Corollary \ref{cor:isotypic} we see that $\A_{\lambda^i}$ is the subalgebra of $\Gamma/\Gamma^i$-invariants
in the canonical quantization of $\C^2/\Gamma^i$. It follows that the image of $\overline{\operatorname{HC}}(\A_{\lambda^i})$ in 
$\operatorname{Rep}(\Gamma)$ is $\operatorname{Rep}(\Gamma/\Gamma^i)$. Let $\mathcal{C}$ be the image of $\overline{\operatorname{HC}}(\A_{\lambda^1},\A_{\lambda^2})$ in $\operatorname{Rep}(\Gamma)$. Assume it is nonzero. 
By (1) in Step 2, it is closed under tensoring with representations that are trivial on $\Gamma^i$ for $i=1,2$.
And by (2), for all $U,V\in \mathcal{C}$, the representation $\operatorname{Hom}_{\CC}(U,V)$ is trivial on $\Gamma^i$
for $i=1,2$. Since $\operatorname{Hom}_{\CC}(V,V)$ contains the trivial summand, we see that for any representation $V^1$ trivial on $\Gamma^1$, the direct summand $V^1\subset V^1\otimes \operatorname{Hom}_{\CC}(V,V)=\operatorname{Hom}_{\CC}(V,V^1\otimes V)$ is trivial on $\Gamma^2$. Similarly, any representation trivial on $\Gamma^2$ is also trivial on $\Gamma^1$. It follows that $\Gamma^1=\Gamma^2$, and since $\lambda^i$ is recovered
from $\Gamma^i$, we have $\lambda^1=\lambda^2$.  
This is a contradiction.
\end{proof}

\chapter{Unipotent ideals and bimodules}\label{sec:unipotent}

Let $G$ be a connected reductive algebraic group and let $\widetilde{\mathbb{O}}$ be a $G$-equivariant nilpotent cover. Recall that the affine variety $\widetilde{X} := \Spec(\CC[\widetilde{\OO}])$ is a conical symplectic singularity, see Example \ref{example:symplecticsingularity}. There is a $G$-action on $\CC[\widetilde{\OO}]$ by graded Poisson automorphisms and a classical co-moment map
$$\varphi: \fg \to \CC[\widetilde{\OO}].$$
induced from the natural map of varieties $\widetilde{\OO} \to \OO \subset \fg^*$.

Consider the extended Namikawa space $\overline{\fP}^{\widetilde{X}} = \fP^{\widetilde{X}} \oplus \fz(\fg)^*$, see Section \ref{subsec:equivariant}. By Proposition \ref{prop:classificationHamiltonian}, there is a bijection
$$\overline{\fP}^{\widetilde{X}}/W^{\widetilde{X}} \simeq \mathrm{Quant}^G(\widetilde{X}), \qquad W \cdot \lambda \mapsto (\cA_{\lambda}^{\widetilde{X}},\Phi_{\lambda}^{\widetilde{X}})$$
where $\Phi_{\lambda}^{\widetilde{X}}: U(\fg) \to \cA_{\lambda}^{\widetilde{X}}$ is the unique quantum co-moment map which restricts to the character of $\fz(\fg)$ determined by $\lambda$. 

For each $\lambda \in \overline{\fP}^{\widetilde{X}}$, consider the two-sided ideal
$$I_{\lambda}(\widetilde{\mathbb{O}}) := \ker{(\Phi_{\lambda}^{\widetilde{X}}: U(\fg) \to \cA^{\widetilde{X}}_{\lambda})} \subset U(\fg).$$
 Ideals of this form enjoy a number of favorable properties, see Proposition \ref{prop:propsofIbeta} below. In particular, every $I_{\lambda}(\widetilde{\OO})$ is primitive, completely prime, and has associated variety $\overline{\OO}$.

\begin{definition}\label{def:unipotentideals}
The \emph{unipotent ideal}\index{ideal!unipotent} attached to $\widetilde{\mathbb{O}}$ is the primitive ideal $I_0(\widetilde{\mathbb{O}}) \subset U(\fg)$. 
\end{definition}
The set of unipotent ideals includes, as a proper subset, all \emph{special} unipotent ideals, cf. Definition \ref{def:spec_unipotent}. We will give a proof of this containment in Section \ref{subsec:refinedBVLS}. More examples (and non-examples) are given in Section \ref{subsec:examplesideals}. In Section \ref{subsec:motivation}, we will provide further motivation for Definition \ref{def:unipotentideals}. In particular, we will demonstrate that our definition follows naturally from Vogan's desiderata, cf. Section \ref{subsec:desiderata}, and the orbit method philosophy. 

In many cases, non-isomorphic covers give rise to the same unipotent ideal. In Section \ref{subsec:classificationideals}, we will introduce an equivalence relation on nilpotent covers. Two covers are equivalent if, roughly speaking, they have the same dimension 2 singularities (see Definition \ref{def:equivalencerelation} for a more precise condition). We will show that the map $\widetilde{\OO} \mapsto I_0(\widetilde{\OO})$ induces a bijection
$$\{\text{equivalence classes of } G\text{-eqvt nilpotent covers}\} \xrightarrow{\sim} \{\text{unipotent ideals}\},$$
proving Theorem \ref{thm:classificationidealsintro} from the introduction. In Section \ref{subsec:classificationbimods}, we will turn our attention to bimodules. In view of Definition \ref{def:unipotentideals}, it is natural to define:

\begin{definition}\label{def:unipotentbimods}
A \emph{unipotent bimodule}\index{Harish-Chandra bimodule!unipotent} attached to $\widetilde{\mathbb{O}}$ is an irreducible Harish-Chandra bimodule $\cB \in \HC^G(U(\fg))$ such that
$$\mathrm{LAnn}(\cB) = \mathrm{RAnn}(\cB) = I_0(\widetilde{\mathbb{O}}).$$
Denote the set of (isomorphism classes of) unipotent bimodules by $\unip_{\widetilde{\mathbb{O}}}(G)$.
\end{definition}

In Section \ref{subsec:classificationbimods}, we will give a geometric classification of unipotent bimodules. The main result is as follows. The equivalence class of $\widetilde{\OO}$ contains a unique maximal element, see Lemma \ref{lem:maximalelement}. Since Definition \ref{def:unipotentbimods} depends only on the equivalence class, we can assume that $\widetilde{\OO}$ is maximal. Let $\Pi = \Aut_{\OO}(\widetilde{\OO})$. For each irreducible $\Pi$-representation $V$, there is a Harish-Chandra bimodule $(\cA_0^{\widetilde{X}} \otimes V)^{\Pi}$. This bimodule has left and right infinitesimal characters since $I_0(\widetilde{\OO})$ is primitive. So, by Lemma \ref{lem:HC_fin_length}, it has finite length. Below we will show that $(\cA_0^{\widetilde{X}} \otimes V)^{\Pi}$ has a unique composition factor with associated variety $\overline{\OO}$.  Denote it by $\cB_V$. We will show that the map $V \mapsto \cB_V$ defines a bijection
$$\{\text{irreducible representations of } \Pi\} \xrightarrow{\sim} \unip_{\widetilde{\OO}}(G),$$
proving Theorem \ref{thm:classificationbimodsintro} from the introduction. In fact, assuming that $I_0(\widetilde{\OO})$ is maximal (a condition which we check in many cases; the remaining cases are checked in \cite{MBM}), the bimodule $(\cA_0^{\widetilde{X}} \otimes V)^{\Pi}$ is irreducible. 

Every Harish-Chandra bimodule carries the structure of a $G$-representation (via the adjoint action of $\fg$). In Section \ref{subsec:Gtypes}, we will show that if $\cB \in \unip_{\widetilde{\OO}}(G)$ and $I_0(\widetilde{\OO})$ is maximal,  then there is a subgroup $H \subset G$ and a finite-dimensional $H$-representation $\chi$ such that 
$$\cB \simeq_G \mathrm{AlgInd}^G_H \chi.$$
This result, which follows easily from the classification above, provides an affirmative answer, in the complex case, to an old conjecture of Vogan (Desideratum (5) from Section \ref{subsec:desiderata}).

\section{Hamiltonian quantizations of nilpotent covers}

Let $\widetilde{\OO}$ be a $G$-equivariant nilpotent cover. Choose $\lambda \in \overline{\fP}^{\widetilde{X}}$, and let $(\cA,\Phi):=(\cA^{\widetilde{X}}_{\lambda}, \Phi^{\widetilde{X}}_{\lambda})$ be the associated Hamiltonian quantization of $\CC[\widetilde{\OO}]$. As above, let $I_{\lambda}(\widetilde{\OO}) = \ker{\Phi} \subset U(\fg)$. Below, we will collect some basic facts about $\cA$, $\Phi$, and $I_{\lambda}(\widetilde{\OO})$. 

Choose $e \in \OO$ and $x \in \widetilde{\OO}$ over $e$. Write $R$ (resp. $R_x$) for the reductive part of the stabilizer $G_e$ (resp. $G_x$). Note that $R_x \subset R$ is a finite-index subgroup. Let $\cW$ denote the $W$-algebra associated to $\OO$. Recall the functors defined in Section \ref{subsec:W}
$$\bullet_{\dagger}: \HC^G(U(\fg))\to \HC^R(\cW), \qquad \bullet^{\dagger}: \HC^R_{\mathrm{fin}}(\cW) \to  \HC^G_{\overline{\OO}}(U(\fg)).$$
as well as the maps
$$\bullet_{\dagger}: \mathrm{Prim}_{\overline{\OO}}(U(\fg)) \to \mathrm{Id}_{\mathrm{fin}}(\cW), \qquad \bullet^{\ddag}: \mathrm{Prim}_{\mathrm{fin}}(\cW) \to \mathrm{Prim}_{\overline{\OO}}(U(\fg))$$
We will need the following result.

\begin{lemma}\label{lem:Adagger}
The following are true:
\begin{itemize}
    \item[(i)] $\cA_{\dagger}$ is a filtered algebra, $R$-equivariantly isomorphic to the $\CC[R/R_x]$ (with the trivial filtration).
    \item[(ii)] There is an $R_x$-stable ideal $J\subset \cW$ of codimension 1 such that 
    $$\cA_{\dagger} \simeq \bigoplus_{r \in R/R_x} \cW/rJ.$$
    \item[(iii)] The adjunction map $\cA \to (\cA_{\dagger})^{\dagger}$ is an isomorphism of filtered algebras.
    \item[(iv)] The ideal $I_\lambda(\widetilde{\OO})$ is primitive. Moreover,  $I_{\lambda}(\widetilde{\OO}) = J^{\ddag}$, where $J \subset \cW$ is the ideal in (ii). 
\end{itemize}
\end{lemma}

\begin{proof}
(i) is \cite[Lem 5.2(1)]{Losev4}. (ii) is an immediate consequence of (i). (iii) is \cite[Lem 5.2(3)]{Losev4}. We proceed to proving (iv). The proof that $I_\lambda(\widetilde{\OO})$ is primitive repeats the beginning of \cite[Section 5.3]{Losev4}. By Theorem \ref{thm:daggersurjective}, $J^\ddag$ is the unique primitive ideal $I\subset U(\fg)$ such that  $I_\dagger$ is the intersection of the $R$-conjugates of $J$. The ideal $I_\lambda(\widetilde{\OO})$ satisfies this property by the definition of $J$, so $I_\lambda(\widetilde{\OO})=J^\ddag$, as desired.
%
%
\end{proof}

\begin{prop}\label{prop:propsofIbeta}
The ideal $I_{\lambda}(\widetilde{\mathbb{O}}) \subset U(\fg)$ enjoys the following properties:
\begin{itemize}
    \item[(i)] $V(I_{\lambda}(\widetilde{\OO})) = \overline{\OO}$.
    \item[(ii)] $\operatorname{\operatorname{\cW-\dim}}(I_{\lambda}(\widetilde{\mathbb{O}}))=1$ (cf. Definition \ref{def:Wdimension}).
    \item[(iii)] Suppose $\widetilde{\mathbb{O}} \to \mathbb{O}$ is Galois and the $\Aut_{\OO}(\widetilde{\OO})$-action on $\CC[\widetilde{\OO}]$ lifts to $\cA_{\lambda}^{\widetilde{X}}$. Then $m_{\overline{\mathbb{O}}}(U(\fg)/I_{\lambda}(\widetilde{\mathbb{O}}))=1$.
    \item[(iv)] $I_{\lambda}(\widetilde{\OO})$ is completely prime.
\end{itemize}
\end{prop}

\begin{proof}
Properties (i) and  (ii) are immediate from Lemma \ref{lem:Adagger}. For (iii), let $\cA'$ denote the algebra of $\Aut_{\OO}(\widetilde{\OO})$-invariants in $\cA$. Note that $\cA'$ is a Hamiltonian quantization of $\CC[\OO]$. By Lemma \ref{lem:co-momentunique} its co-moment map coincides with the restriction of $\Phi$ to $\cA' \subset \cA$. Thus, we get an embedding $U(\fg)/I_{\lambda}(\widetilde{\OO})\hookrightarrow \cA'$. Taking multiplicities, we get an inequality $m_{\overline{\OO}}(U(\fg)/I_{\lambda}({\widetilde{\OO}}))\le m_{\overline{\OO}}(\cA')$. Since $\cA'$ is a quantization of $\CC[\OO]$, we have $m_{\overline{\OO}}(\cA')=1$. Property (iii) follows. (iv) follows from Corollary \ref{Cor:completely_prime} combined with (iv) of Lemma \ref{lem:Adagger}.
\end{proof}


\section{Examples of unipotent ideals and bimodules}\label{subsec:examplesideals}

In this section, we will collect some examples (and non-examples) of unipotent ideals and bimodules. 

\begin{example}\label{ex:unipotent}
\leavevmode
\begin{itemize}
    \item[(i)] Let $\mathbb{O} = \{0\}$. Then $I_0(\mathbb{O})$ is the \emph{augmentation ideal} (i.e. the annihilator of the trivial $U(\fg)$-module).
    
    \item[(ii)] More generally, let $\mathbb{O}^{\vee}$ be a nilpotent orbit for the Langlands dual group $G^{\vee}$ and consider the \emph{special unipotent ideal} $I_{\mathrm{max}}(\frac{1}{2}h^{\vee}) \subset U(\fg)$ (cf. Definition \ref{def:spec_unipotent}). In Section \ref{subsec:refinedBVLS}, we will show that $I_{\mathrm{max}}(\frac{1}{2}h^{\vee})$ is unipotent (and describe the cover it is attached to).
    
    \item[(iii)] Let $\fg$ be a simple Lie  algebra and let $\mathbb{O}_{\mathrm{min}}$ be the minimal nilpotent orbit. If $\fg$ is not of type $A$, there is a unique (maximal) completely prime ideal $I_0 \subset U(\fg)$ such that $V(I_0) = \overline{\mathbb{O}}_{\mathrm{min}}$ called the \emph{Joseph ideal} (see \cite{Joseph1976})\index{ideal!Joseph}. This ideal can be constructed in a number of equivalent ways, see e.g. \cite{GarfinkleThesis} and \cite{LevasseurSmithStafford}. It follows from (iv) of Proposition \ref{prop:propsofIbeta} that $I_0(\OO_{\mathrm{min}}) = I_0$. 
    
    \item[(iv)] As a special case of (iii), consider the Joseph ideal $I_0$ for $\fg=\mathfrak{sp}(2n)$. There are two irreducible Harish-Chandra bimodules annihilated by $I_0$---they are the irreducible constituents of the metaplectic representation. Both are unitary. These bimodules are among the most familiar examples of unipotent bimodules which are not special unipotent.

    \item[(v)] Let $G=\mathrm{SL}(2n)$ and let $\widetilde{\mathbb{O}}$ be the universal cover of the so-called `model' nilpotent orbit $\mathbb{O} = \mathbb{O}_{(2^n)}$. We will see below (see Proposition \ref{prop:centralcharacterbirigidcover}) that $I_0(\widetilde{\mathbb{O}}) = I_{\mathrm{max}}(\frac{\rho}{2})$. This ideal plays a central role in the determination of the unitary dual of the universal cover of $\mathrm{SL}(n,\RR)$ (see \cite{Huang1990}).
\end{itemize}

\end{example}

\begin{example}\label{ex:nonexample}
Let $G=\mathrm{Sp}(2n)$ and let $\mathbb{O} = \mathbb{O}_{(2^n)}$ be the `model' nilpotent orbit. In \cite{McGovern1994}, McGovern attaches to $\OO$ several so-called \emph{q-unipotent}\index{ideal!$q$-unipotent} ideals (all are \emph{weakly unipotent}\index{ideal!weakly unipotent} in the sense of \cite{Vogan1984}). One of these ideals is the maximal ideal $I_{\mathrm{max}}(\frac{\rho}{2})$ (it is completely prime, like all of McGovern's ideals). We note that if $n \geq 2$, this ideal is \emph{not} unipotent in the sense of Definition \ref{def:unipotentideals}. This is evident from the infinitesimal character calculations in Section \ref{subsec:centralcharclassical} (see Proposition \ref{prop:centralcharacterbirigidcover}). Unitarity considerations justify its exclusion. For example, an {\tt atlas} computation shows that the spherical Harish-Chandra bimodule $U(\fg)/I_{\mathrm{max}}(\frac{\rho}{2})$ is non-unitary for $2 \leq n \leq 7$ (for $n > 7$, the unitarity calculation cannot be completed in a reasonable amount of time). 
\end{example}

\section{Motivation}\label{subsec:motivation}

In this section, we will provide some evidence in favor of Definition \ref{def:unipotentideals}.

\subsection{Unitarity}\label{SSS_motivation_unitarity}

The most basic requirement of unipotent bimodules is that they are unitary (we will recall what this means in Section \ref{subsec:unitarydef}). If $\cB$ is an irreducible unitary Harish-Chandra bimodule, then by a classical result of Vogan (see \cite[Prop 7.12]{Vogan1986a}), the left and right annihilators of $\cB$ are completely prime ideals. Thus, a minimal requirement of all unipotent ideals is that they are completely prime.

The set of completely prime primitive ideals in $U(\fg)$ includes all primitive ideals of $\cW$-dimension 1, 
see Corollary \ref{Cor:completely_prime}. Conversely, every completely prime primitive ideal with integral infinitesimal character is of $\cW$-dimension 1 (with the possible exception of primitive ideals attached to a single orbit in type $E_8$), see \cite{Losev3}, \cite{Losev2014}. At non-integral infinitesimal character, there are completely prime primitive ideals with $\operatorname{\cW-\dim}(I) > 1$. For example, if $\fg = \mathfrak{sp}(2n)$ and $n \geq 2$, the maximal ideal $I_{\mathrm{max}}(\rho/2)$ is completely prime, whereas
    $$\operatorname{\cW-\dim}(I_{\mathrm{max}}(\rho/2)) = 2^{\floor{\frac{n-1}{2}}},$$ 
    see \cite{Premet2010}, \cite[Proposition 4.1]{LosevPanin}. However, as noted in Example \ref{ex:nonexample}, the Harish-Chandra bimodules $U(\fg)/I_{\mathrm{max}}(\rho/2)$ are non-unitary in all cases we have checked. These examples suggest the following conjecture.

\begin{conj}\label{conj:unitaryWdim}
The left and right annihilators of 
unitary Harish-Chandra bimodules are of $\cW$-dimension $1$.
\end{conj}
In view of this conjecture, and the examples supporting it, we restrict our attention to primitive ideals of $\cW$-dimension 1. This set includes all ideals of the form $I_{\lambda}(\widetilde{\mathbb{O}})$, see Proposition \ref{prop:propsofIbeta}. In classical types, all primitive ideals of $\cW$-dimension 1 are of this form, see \cite[Thm 1.2]{Topley}. In exceptional types, there are exceptions (a few will be encountered below). Nonetheless, ideals of the form $I_\lambda(\widetilde{\OO})$ constitute a particularly nice class of primitive ideals of $\cW$-dimension $1$.

\subsection{Orbit method}\label{subsec:orbitmethod}

The considerations above suggest that all unipotent ideals are of the form $I_{\lambda}(\widetilde{\OO})$. 
The orbit method philosophy suggests that one should further restrict to $\lambda=0$. 

\begin{definition}[\cite{Vogan1990}, Def 2.1]
A \emph{Dixmier algebra}\index{Dixmier algebra} for $G$ is a pair $(\cA,\varphi)$ consisting of an associative algebra $\cA$ and an algebra homomorphism $\varphi: U(\fg) \to \cA$ such that $\cA$, regarded as a $U(\fg)$-bimodule via $\varphi$, is Harish-Chandra (cf. Definition \ref{def:HCbimodsclassical}). An isomorphism of Dixmier algebras $(\cA_1, \varphi_1) \xrightarrow{\sim} (\cA_2, \varphi_2)$ is an isomorphism of algebras $\phi: \cA_1 \to \cA_2$ which intertwines $\varphi_1$ and $\varphi_2$.
\end{definition}

One version of the orbit method, articulated in \cite{Vogan1990}, is the following. For each  $G$-equivariant cover $\widetilde{\OO}^1 $ of a coadjoint orbit $ \OO^1$ (not necessarily nilpotent), there should be an associated Dixmier algebra, denoted $\mathrm{Dix}(\widetilde{\OO}^1)$. As a representation of $G$, $\mathrm{Dix}(\widetilde{\OO}^1)$ should coincide with the ring of regular functions $\CC[\widetilde{\OO}^1]$. Finally, the passage from $\widetilde{\OO}^1$ to $\mathrm{Dix}(\widetilde{\OO}^1)$ should define an injective correspondence
$$\mathrm{Dix}: \{G\text{-eqvt covers of coadjoint orbits}\} \hookrightarrow \{\text{Dixmier algebras for }G\}.$$
This existence of this injection, with the property described above, is known as the \emph{Dixmier conjecture}, see \cite[Conj 2.3]{Vogan1990}. In the terminology of \cite{Vogan1990}, a \emph{unipotent Dixmier algebra}\index{Dixmier algebra!unipotent} is one which corresponds to a \emph{nilpotent} cover under this (not-yet-defined) correspondence. A \emph{unipotent ideal} is the kernel in $U(\fg)$ of the co-moment map for a unipotent Dixmier algebra. 

In Section \ref{subsec:Dixmier}, we will give a definition of $\mathrm{Dix}$, extending a result of the first-named author, see \cite[Thm 5.3]{Losev4}. Our map has the following property: If $\widetilde{\OO}$ is a nilpotent cover, then $\mathrm{Dix}(\widetilde{\OO})$ is the canonical quantization $\cA_0^{\widetilde{X}}$ (regarded as a Dixmier algebra via the co-moment map $\Phi^{\widetilde{X}}_0: U(\fg) \to \cA_0^{\widetilde{X}}$). Thus, Definition \ref{def:unipotentideals} is consistent with the orbit method philosophy.

\subsection{Maximality}

It is generally believed that unipotent ideals should be maximal (see Desideratum (3) from Section \ref{subsec:desiderata} and the references therein). In Appendix \ref{sec:maximality}, we will show that for linear classical groups, the ideal $I_0(\widetilde{\OO})$ is always maximal. This result is extended to arbitrary groups in the paper \cite{MBM}.

Now we show that the maximality of $I_0(\widetilde{\OO})$ is equivalent to the simplicity of $\cA_0^{\widetilde{X}}$. The latter is a general expectation of canonical quantizations of conical symplectic singularities, see Section \ref{subsec:conjecturesimplicity}.

\begin{prop}\label{prop:simplemaximal}
The following conditions are equivalent:
\begin{enumerate}
\item $\cA_0^{\widetilde{X}}$ is simple.
\item $I_0(\widetilde{\mathbb{O}})$ is a maximal ideal.
\end{enumerate}
\end{prop}

\begin{proof}
To simplify notation, let $I:=I_0(\widetilde{\mathbb{O}})$ and $\cA:=\cA_0^{\widetilde{X}}$. We first prove (1)$\Rightarrow$(2). Suppose $I$ is not maximal. Let $\cW$ be the $W$-algebra associated to $\mathbb{O}$ and let $R$ be the reductive part of the centralizer of $e \in \mathbb{O}$. Consider the functors $\bullet_{\dagger}: \HC^G(U(\fg)) \to \HC^R(\cW)$ and $\bullet^{\dagger}: \HC^R_{\mathrm{fin}}(\cW) \to \HC^G_{\overline{\OO}}(U(\fg))$ constructed in Section \ref{subsec:W}. Note that $\cA_{\dagger}$ is a completely reducible 
$R$-equivariant $\mathcal{W}$-bimodule (in fact, it is the direct sum of $1$-dimensional $\mathcal{W}$-bimodules, Lemma 
\ref{lem:Adagger}). Since $U(\fg)/I\subset \cA$ is an inclusion of Harish-Chandra bimodules, and $\bullet_\dagger$ is exact (Proposition \ref{prop:propsofdagger}), $(U(\fg)/I)_\dagger$ is an $R$-stable subbimodule of $\cA_{\dagger}$.  

Consider the algebra $\cA'=[(U(\fg)/I)_\dagger]^{\dagger}$. Recall that
$\cA=(\cA_\dagger)^\dagger$, see Lemma \ref{lem:Adagger}(iii). Thus $\cA'$ embeds into 
$\cA$ as a subalgebra and as a $U(\fg)$-bimodule direct summand. Moreover, the injective map $U(\fg)/I\to \cA=(\cA_\dagger)^\dagger$ factors through the natural map $U(\fg)/I\to ([U(\fg)/I]_\dagger)^{\dagger}=\cA'$. It follows that the latter map is injective, and thus by \cref{prop:propsofdagger}(iv) the bimodule $\cA'$ is an extension of $U(\fg)/I$ by a bimodule supported on $\partial\mathbb{O}$. Note that for every ideal $\widetilde{I}\subset U(\g)$ 
properly containing $I$ we have $\mathcal{V}(U(\fg)/\widetilde{I})\subsetneq \mathcal{V}(U(\fg)/I)$, since $I$ is a prime ideal, see (iii) of Proposition \ref{prop:propertiesofprim}.
Since $U(\g)/I$ has finite length, see Lemma \ref{lem:HC_fin_length}, there is a minimal ideal $\widetilde{I}\subset U(\g)$ properly containing $I$. 

We claim that $\widetilde{I}/I$ is a 2-sided  ideal in $\cA'$. Let us show that it is a right ideal (the proof that it is a left ideal is analogous). The associated variety $\mathcal{V}(\cA'/[\widetilde{I}/I])$ of the $U(\fg)$-bimodule $\cA'/(\widetilde{I}/I)$ is a proper subvariety of $\mathcal{V}(\cA')=\mathcal{V}(U(\g)/I)$. Hence $\mathcal{V}(\widetilde{I}/I)=\mathcal{V}(\cA')$. 
If $\B:=(\widetilde{I}/I)\otimes_{U(\g)}(\cA'/[\widetilde{I}/I])$ is nonzero, then  $\mathcal{V}(\B)\subsetneq \mathcal{V}(\cA')$.  By the tensor-Hom adjunction, the simple $U(\g)$-bimodule $\widetilde{I}/I$ embeds into $\Hom_{U(\g)}(\cA'/[\widetilde{I}/I],\B)$, in fact, into the locally finite part, and the associated variety of the latter is contained in the associated variety of $\B$. This contradicts the fact that $\mathcal{V}(\widetilde{I}/I)=\mathcal{V}(\cA')$. It follows that $\B=\{0\}$. The right multiplication in  $\cA'$ gives a map $\widetilde{I}/I\otimes_{U(\g)} \cA'\rightarrow \cA'$, since $\mathcal{B}=\{0\}$, we get  
$[\widetilde{I}/I]\A'\subseteq \widetilde{I}/I$. Thus $\widetilde{I}/I$ is a nonzero right (and hence two-sided) ideal in $\cA'$.

Since $\cA'$ is a direct bimodule summand of $\cA$, we have that $\cA/[\cA (\widetilde{I}/I)]\neq \{0\}$. Let ${I}'$ be the left annihilator of $\cA/\cA (\widetilde{I}/I)$. Note that 
$\mathcal{V}(\cA/\cA (\widetilde{I}/I))\subset \mathcal{V}(U(\fg)/\widetilde{I})
\subsetneq \mathcal{V}(\cA)$. So $I'$ is a  two-sided ideal in $\cA$, different from $\cA$. Since  $\cA/I'$ is a 
HC bimodule over $U(\fg)/\widetilde{I}$ we see that $\mathcal{V}(\cA/I') \subset \mathcal{V}(U(\fg)/\widetilde{I})$.
On the other hand, $\mathcal{V}(U(\fg)/\widetilde{I})\subsetneq \mathcal{V}(\cA)$. It follows that $I'\neq \{0\}$.  Since $\cA$ is simple, we get a contradiction. This proves (1)$\Rightarrow$(2). 

Next we prove (2)$\Rightarrow$(1). Suppose $\cA$ is not simple. Then there is a nonzero prime ideal $\widetilde{I} \subset \cA$. Since $\cA$ is a Noetherian domain, the GK dimension of the quotient $\cA/\widetilde{I}$ is less than that of $\cA$, and hence the GK dimension of the image of $U(\fg)/I$ in $\cA/\widetilde{I}$
is less than that of $U(\fg)/I$. Since $I$ is maximal, $U(\fg)/I$
has no such algebra quotient, a contradiction. This completes the proof.
\end{proof}

\subsection{Applicability to real reductive groups}\label{subsec:realgroups}

In this section, we will discuss the applicability of Definition \ref{def:unipotentbimods} to the more general setting of real reductive Lie groups. This issue is quite subtle and deserving of careful consideration. It is further explored in \cite{MBM} and \cite{DavisMasonBrown}.

Let $G_{\RR}$ be a real reductive Lie group. For this informal discussion, the precise meaning of `real reductive Lie group' will not matter enormously. It should be broad enough to include non-linear covering groups, such as $\mathrm{Mp}(2n,\RR)$, but narrow enough so that the theory of Harish-Chandra modules is applicable. In particular, $G_{\RR}$ should have finitely many connected components and the identity component should be a finite cover of an algebraic group. Let $\fg$ denote the complexification of the Lie algebra $\fg_{\RR}$ of $G_{\RR}$ and $K$ the complexification of a maximal compact subgroup $K_{\RR} \subset G_{\RR}$. A $(\fg,K)$-\emph{module} is a $U(\fg)$-module with a compatible algebraic $K$-action. A $(\fg,K)$-module is \emph{admissible} if each irreducible $K$-representation appears with finite multiplicity. Under suitable conditions on $G_{\RR}$, there is a bijection between the set of irreducible admissible $G_{\RR}$-representations (up to infinitesimal equivalence) and the set of irreducible admissible $(\mathfrak{g},K)$-modules (up to isomorphisn), see \cite[Chp 0.3]{Vogan1981} for an overview. We note that if $G_{\RR}$ is complex, then a $(\fg,K)$-module is the same thing as a $G_{\RR}$-equivariant Harish-Chandra bimodule for the (complex) Lie algebra $\fg_{\RR}$. As explained in Chapter \ref{sec:intro}, the set $\mathrm{Irr}_u(G_{\RR})$ of irreducible unitary $G_{\RR}$-representations should contain a finite set of `building blocks' $\mathrm{Unip}(G_{\RR}) \subset \mathrm{Irr}_u(G_{\RR})$ called `unipotent representations'. In the complex case, $\mathrm{Unip}(G_{\RR})$ should be the set of Definition \ref{def:unipotentbimods}. Our problem is to formulate a natural generalization of Definition \ref{def:unipotentbimods} for arbitrary $G_{\RR}$. 

If $G_{\RR}$ is algebraic, there is a notion of a \emph{special unipotent representation}, due to Adams, Barbasch, and Vogan. The definition in \cite[Sec 27]{AdamsBarbaschVogan} is stated in the language of $L$-groups and Arthur parameters, but is equivalent to the following (see \cite[Cor 27.13]{AdamsBarbaschVogan}). 

\begin{definition}[\cite{AdamsBarbaschVogan}]\label{def:realunipotent}
Suppose $G_{\RR}$ is an algebraic group and let $\OO^{\vee}$ be a nilpotent orbit for the Langlands dual Lie algebra $\fg^{\vee}$. A \emph{special unipotent representation}\index{representation!special unipotent} attached to $\OO^{\vee}$ is an irreducible $(\fg,K)$-module $\cB$ such that $\Ann(\cB) = I_{\mathrm{max}}(h^{\vee}/2)$.
\end{definition}

\begin{conj}[\cite{Arthur1989}]\label{conj:Arthur}
Suppose $G_{\RR}$ is an algebraic group. Then all special unipotent representations of $G_{\RR}$ are unitary.
\end{conj}

Conjecture \ref{conj:Arthur} is now known in all cases. More precisely, it is proved
\begin{itemize}
    \item for real exceptional groups in \cite{AdamsMillerVogan};
    \item for quasisplit symplectic and special orthogonal groups in \cite{Arthur2013};
    \item for complex classical groups in \cite{Barbasch1989};
    \item for all classical groups in \cite{BMSZunitarity};
    \item for all complex groups in \cite{DavisMasonBrown}.
\end{itemize}

Notably, if we allow $G_{\RR}$ to be a non-linear group, Conjecture \ref{conj:Arthur} is false.

\begin{example}
Let $G_{\RR} = \mathrm{Mp}(4,\RR)$. Hence, $\fg = \mathfrak{sp}(4)$ and $K$ is the `square-root of determinant' cover of $GL(2)$
$$K = \left\{(g,z) \in GL(2) \times \CC^{\times} \mid \det(g) = z^2\right\}$$
Let $\OO^{\vee}$ be the subregular nilpotent orbit for $\fg^{\vee} = \mathfrak{so}(5)$. Then in standard coordinates on $\fh^*$, we have $\frac{1}{2}h^{\vee} = (1,0)$. There are two degenerate principal series representations of $G_{\RR}$ which are annihilated by the ideal $I_{\mathrm{max}}(\frac{1}{2}h^{\vee})$. They are irreducible by \cite[Thm 1]{KudlaRallis}. In the same (standard) coordinates, the $K$-types are as follows:

$$2\ZZ \times 2\ZZ + (\frac{1}{2},\frac{1}{2}) \qquad \text{or} \qquad 2\ZZ \times 2\ZZ - (\frac{1}{2},\frac{1}{2}). $$
By an easy application of the Dirac inequality, we see that neither representation is unitary. We are grateful to Jing-Song Huang for suggesting this approach.
\end{example}

We will show in Chapter \ref{sec:duality} that the set of unipotent ideals includes (as a proper subset) all special unipotent ideals. This suggests the following generalization of Definition \ref{def:realunipotent}.
\begin{fakedefinition}\label{def:unipotent1}
Suppose $G_{\RR}$ is an algebraic group and let $\widetilde{\OO}$ be a finite cover of a nilpotent co-adjoint $\Ad(\fg)$-orbit. A \emph{unipotent representation}\index{representation!unipotent} of $G_{\RR}$ attached to $\widetilde{\OO}$ is an irreducible $(\fg,K)$-module $\cB$ such that $\mathrm{Ann}(\cB) = I_0(\widetilde{\OO})$. 
\end{fakedefinition}

We use quotation marks to indicate that this is \emph{not} a reasonable definition---if we take the ideal $I_0(\widetilde{\OO})$ to be non-special, there can be non-unitary $(\fg,K)$-modules which are annihilated by $I_0(\widetilde{\OO})$ even if $G_{\RR}$ is algebraic. This problem is already apparent for $\mathrm{SL}(2,\RR)$.

\begin{example}\label{ex:Mp2R}
Let $G_{\RR} = \mathrm{SL}(2,\RR)$. Hence, $\fg=\mathfrak{sl}(2)$ and $K=SO(2)$. The unipotent ideal attached to the (universal) 2-fold cover of the principal nilpotent $\Ad(\fg)$-orbit is the maximal ideal in $U(\fg)$ of infinitesimal character $\frac{\rho}{2}$, see (v) of Example \ref{ex:unipotent}. There are two irreducible $(\fg,K)$-modules annihilated by this ideal: one spherical and one non-spherical. The spherical representation is unitary (it is the midpoint of a unitary complementary series). The non-spherical representation, however, is \emph{not} unitary. In particular, it should not be regarded as a unipotent representation of $G_{\RR}$. 
\end{example}

We believe that `Definition' \ref{def:unipotent1} can be fixed in every case by replacing $G_{\RR}$ with an appropriate non-linear covering group. This covering group should depend in a predictable way on the infinitesimal character of $I_0(\widetilde{\OO})$. In Example \ref{ex:Mp2R}, the infinitesimal character of $I_0(\widetilde{\OO})$ is half-integral. This suggests that we should look for genuine representations of the (non-linear) 2-fold covering group $\mathrm{Mp}(2,\RR)$. There are two genuine irreducible representations of $\mathrm{Mp}(2,\RR)$ annihilated by the ideal $I_0(\widetilde{\OO})$, namely the irreducible constituents of the oscillator representation. Both are unitary, and worthy of the label `unipotent'.\index{representation!oscillator}

\begin{rmk}
Further evidence for this idea can be found in \cite{barbasch2020metaplectic}. There, the authors define a class of maximal ideals for $\fg=\mathfrak{sp}(2n)$ called `metaplectic special unipotent'. One can check that all of their ideals are unipotent in our sense, compare to Chapter \ref{sec:duality}. They prove in a later paper that all genuine irreducible representations of $\mathrm{Mp}(2n,\RR)$ annihilated by a metaplectic special unipotent ideal are unitary, see \cite[Theorem 2.1]{BMSZunitarity}.
\end{rmk}

However, if we assume that $\OO$ is birationally rigid, these subtleties seem to vanish. In \cite{DavisMasonBrown}, the second-named author and Dougal Davis have proved the following result.

\begin{theorem}[{\cite[Corollary 5.23]{DavisMasonBrown}}]\label{thm:realunipotent}
Suppose $\OO$ is a birationally rigid nilpotent $\Ad(\fg)$-orbit. Then every irreducible $(\fg,K)$-module $\cB$ with $\Ann(\cB)=I_0(\OO)$ is unitary.
\end{theorem}

\begin{example}
Let $G_{\RR}$ be the split real form of the (unique) simple group of type $E_8$ and consider the nilpotent orbits
$$\OO_1 = 2A_2+2A_1, \qquad \OO_2=  A_4+A_3.$$
Both orbits are rigid, see \cite{deGraafElashvili}. We will explain in Chapter \ref{sec:centralchars} how the ideals $I_0(\OO)$ can be computed---in these cases we have
$$I_0(\OO_1) = I_{\mathrm{max}}(\rho/3), \qquad I_0(\OO_2) = I_{\mathrm{max}}(\rho/5).$$
An {\tt atlas} computation shows that in each case, there is a unique irreducible $(\fg,K)$-module annihilated by $I_0(\OO_i)$, and it is unitary. 
\end{example}

\begin{example}
Let $G_{\RR}$ be the non-linear group $\mathrm{Mp}(2n,\RR)$ (for $n \geq2$) and let $\OO$ be the minimal orbit. As noted in Example \ref{ex:unipotent}(iv), $I_0(\OO)$ is the Joseph ideal. There are two irreducible $(\fg,K)$-modules annihilated by this ideal---they are the irreducible constituents of the oscillator representation. It is a classical fact that both are unitary. 
\end{example}

Theorem \ref{thm:realunipotent} suggests that if $\OO$ is birationally rigid, `Definition' \ref{def:unipotent1} is correct. We propose the following.

\begin{definition}\label{def:unipotent2}
Suppose $\OO$ is a birationally rigid nilpotent $\Ad(\fg)$-orbit. A \emph{unipotent representation}\index{representation!unipotent} of $G_{\RR}$ attached to $\OO$ is an irreducible $(\fg,K)$-module $\cB$ such that $\mathrm{Ann}(\cB) = I_0(\OO)$. Write $\unip_{\OO}(G_{\RR})$ for the set of (equivalences classes of) such representations.
\end{definition}


As further vindication of Definition \ref{def:unipotent2}, we offer the following remark. It is shown in \cite{Losev_Yu} that the set $\unip_{\OO}(G_{\RR})$ admits a geometric parameterization (analogous to our parameterization of unipotent bimodules, see Theorem \ref{thm:classificationbimods} below). More precisely, they prove that for $G_{\RR}$ algebraic, $\unip_{\OO}(G_{\RR})$ is in bijection with pairs $(\OO_K, \mathcal{L})$, where $\OO_K$ is a $K$-orbit in $\OO \cap (\fg/\mathfrak{k})^*$ and $\mathcal{L}$ is an irreducible twisted local system on $\OO_K$
with half-canonical twist.

\section{Equivariant bimodules for Hamiltonian quantizations of nilpotent covers}\label{subsec:equivariantbimods}

In this section we will establish some preliminary facts about Harish-Chandra bimodules for quantizations of nilpotent covers, specializing (and modifying) the general results and constructions of Sections \ref{subsec:daggers}, \ref{subsec:coverings}.

Let $\widetilde{\OO}$ be a $G$-equivariant nilpotent cover and let $(\cA,\Phi):=(\cA^{\widetilde{X}}_{\lambda},\Phi^{\widetilde{X}}_{\lambda})$ be a Hamiltonian quantization of $\CC[\widetilde{\OO}]$. On any Harish-Chandra bimodule $\cB \in \HC(\cA)$, there is an adjoint action of $\fg$, defined by
$$\ad: \fg \to \End(\cB), \qquad \ad(\xi)(b) = \Phi(\xi)b-b\Phi(\xi), \qquad \xi \in \fg, b\in \cB.$$

\begin{definition}\label{def:eqvtbimods}
A Harish-Chandra bimodule $\cB \in \HC(\cA)$ is $G$-\emph{equivariant} if $\ad$ integrates to a {rational} action $\Ad$ of $G$. Denote the full subcategory of $G$-equivariant Harish-Chandra bimodules by $\HC^G(\cA) \subset \HC(\cA)$ and write $\overline{\HC}^G(\cA)$ for the quotient category $\HC^G(\cA)/\HC^G_{\partial}(\cA)$.
\end{definition}

Let $\Gamma = \pi_1(\widetilde{\OO})$. By (\ref{eq:daggers_inverse}), there is an equivalence of categories
$$\bullet_{\dagger}: \overline{\HC}(\cA) \xrightarrow{\sim} \Gamma/\Gamma(\lambda)\modd$$
In this section, we will provide an analogous description of the category $\overline{\HC}^G(\cA)$. Let $\Omega = \pi_1^G(\widetilde{\OO})$ and let $\Omega(\lambda)$ be the image of $\Gamma(\lambda)$ under the quotient map $\Gamma \twoheadrightarrow\Omega$.

\begin{prop}\label{prop:classificationeqvtbimods}
The functor $\bullet_{\dagger}: \overline{\HC}(\cA) \xrightarrow{\sim} \Gamma/\Gamma(\lambda)\modd$ restricts to an equivalence
$$\bullet_{\dagger}:\overline{\HC}^G(\cA) \xrightarrow{\sim} \Omega/\Omega(\lambda)\modd$$
\end{prop}
\begin{proof}
    First, note that the adjoint action of $\fz(\fg)$, and hence of $Z(G)^{\circ}$, on any bimodule $\cB\in \HC^G(\cA)$ is trivial. Replacing $G$ with $G/Z(G)^\circ$ if necessary, we can assume that $G$ is semisimple. Let $G^{\mathrm{sc}}$ denote the universal cover of $G$, and $K$ the kernel of the quotient map $G^{\mathrm{sc}}\to G$. An object in $\overline{\HC}^G(\cA_{\lambda})$ is a bimodule $\cB\in \overline{\HC}(\cA_\lambda)$ on which the action of $K$ is trivial. Since   $\widetilde{\OO}$ is a $G$-equivariant cover, we see that $K\subset G^{\mathrm{sc}}_x$ for all $x\in \widetilde{\OO}$. So $K$ maps to $\Gamma$ with cokernel equal to $\Omega$. 
    
    We claim that the $K$-action on $\cB\in \HC(\cA_\lambda)$ is trivial if and only if the $K$-action on $\cB_{\dagger}$ is trivial. 
    Clearly, $K$ acts trivially on $\cB$ if and only if it acts trivially on $\gr \cB$. By (\ref{eqn:grdagger}), we have $(\gr\cB)|_{X^{\mathrm{reg}}}\simeq (\cB_\dagger\otimes \cO_{\widehat{X}^{\mathrm{reg}}})^\Gamma$. Since $\Gamma$ acts freely on $\widehat{X}^{\mathrm{reg}}$, $K$ acts trivially on $(\gr\cB)|_{X^{\mathrm{reg}}}$ if and only if it acts trivially on $\cB_{\dagger}$. This completes the proof.
\end{proof}


The next lemma will be useful for computing $\cB_\dagger$ when $\cB\in \HC^G(\cA)$. Since the image of $\fg$ in $\cA$ lives in degree $\leqslant d$, every good filtration on $\cB$ is stable under $\ad(\fg)$, and hence $\Ad(G)$. In particular, $\gr(\cB)$ is a $G$-equivariant $\CC[\widetilde{X}]$-module and the $G$-action is uniquely recovered from the Poisson structure.

\begin{lemma}\label{lem:determineddagger}
For any $\cB\in \HC^G(\cA)$ (and any good filtration), there is an isomorphism of $G$-equivariant coherent sheaves on $\widetilde{\OO}$
\begin{equation}\label{eqn:inductiondagger1}
(\gr \cB)|_{\widetilde{\mathbb{O}}} \simeq (\cB_{\dagger} \otimes \cO_{\widehat{\mathbb{O}}})^{\Omega}.
\end{equation}
This isomorphism uniquely determines the $\Omega$-module $\cB_\dagger$.
\end{lemma} 
\begin{proof}
We can assume in the proof that $G$ is semisimple and simply connected and hence that $\HC^G(\cA)=\HC(\cA)$. Recall that for a smooth symplectic variety, a Poisson coherent sheaf is the same thing as a local system, see Step 3 in the proof of \cite[Lemma 3.9]{B_ineq}.

By (\ref{eqn:grdagger}) (and using the notation therein) there is a  local system isomorphism
\begin{equation}\label{eq:local_system_iso1}
(\gr \cB)|_{\widetilde{X}^{\mathrm{reg}}} \simeq (\cB_{\dagger} \otimes \cO_{\widehat{X}'})^{\Gamma}.
\end{equation}
On both sides of (\ref{eq:local_system_iso1}), the $\fg$-action is recovered from the local system structure. It follows that the isomorphism (\ref{eq:local_system_iso1}) is $\fg$- and hence $G$-equivariant.

Since $\operatorname{codim}(\widetilde{X}^{\mathrm{reg}}- \widetilde{\OO}, \widetilde{X}^{\mathrm{reg}})\geqslant 2$, (\ref{eq:local_system_iso1}) is equivalent to a $G$-equivariant local system isomorphism
\begin{equation}\label{eq:local_system_iso2}
(\gr \cB)|_{\widetilde{\OO}} \simeq (\cB_{\dagger} \otimes \cO_{\widehat{\OO}})^{\Gamma}.
\end{equation}
This implies (\ref{eqn:inductiondagger1}). 
It remains to show that (\ref{eqn:inductiondagger1}) determines $\cB_\dagger$ uniquely. For this, it suffices to show that the forgetful functor from the category of $G$-equivariant local systems on $\widetilde{\OO}$ to the category of $G$-equivariant coherent sheaves on $\widetilde{\OO}$ is a full embedding. Choose $x\in \widetilde{\OO}$. Via restriction to $x$, these categories are equivalent to $G_x/G_x^{\circ}\operatorname{-mod}$ and $G_x\operatorname{-mod}$, respectively. The forgetful functor is the pullback under $G_x\twoheadrightarrow G_x/G_x^{\circ}$, which is indeed a full embedding. This completes the proof. 
%
\end{proof}

We will also need a $G$-equivariant version of Corollary \ref{cor:isotypic}. Let $\widecheck{\OO} \to \widetilde{\OO}$ be a $G$-equivariant Galois cover and let $\Pi := \Aut_{\widetilde{\OO}}(\widecheck{\OO})$. Note that the $\Pi$-action on $\widecheck{\OO}$ is by $G$-equivariant Hamiltonian automorphisms.  Suppose $\cA^{\widecheck{X}}$ is a $\Pi$-equivariant Hamiltonian quantization of $\CC[\widecheck{\OO}]$  such that $\cA \simeq (\cA^{\widecheck{X}})^\Pi$ as Hamiltonian quantizations. In view of Proposition \ref{prop:classificationeqvtbimods}, the following result is completely analogous to Corollary \ref{cor:isotypic}, and is proved in a similar fashion.

\begin{cor}\label{cor:isotypiceqvt}
Suppose $\overline{\HC}^G(\cA^{\widecheck{X}}) \simeq \Vect$. Then $\bullet_{\dagger}$ induces an equivalence of monoidal categories
$$\bullet_{\dagger}: \overline{\HC}^G(\cA) \xrightarrow{\sim} \Pi\modd.$$
The inverse equivalence, $\bullet^{\dagger}$, is the functor of \emph{isotypic components}
$$V^{\dagger} \simeq (\cA^{\widecheck{X}}\otimes V)^{\Pi}, \qquad  V \in \Pi\modd.$$
\end{cor}

Our next task is to relate the restriction and extension functors of Section 
\ref{subsec:W} to those of 
Section \ref{subsec:daggers}.
Let  $\cA:=\cA^{\widetilde{X}}_\lambda, I:=I_\lambda(\widetilde{\OO}) \subset U(\fg)$ and $\underline{J} := I_\dagger\subset \cW$. Consider the full subcategory of $\HC^G(U(\fg))$
\begin{equation}\label{eq:HC_Ug_I}
\HC^G(U(\fg)/I) := \{\cB \in \HC^G(U(\fg)) \mid I \subseteq \mathrm{LAnn}(\cB), \  I \subseteq \mathrm{RAnn}(\cB)\}
\end{equation}
and the Serre subcategory
\begin{equation}\label{eq:Serresubcategory}
\HC^G_{\partial}(U(\fg)/I) :=  \{\cB \in \HC^G(U(\fg)) \mid I \subsetneq \mathrm{LAnn}(\cB), \  I \subsetneq \mathrm{RAnn}(\cB)\}
\subset \HC^G(U(\fg)/I)\}.
\end{equation}
Note that $\HC_{\partial}^G(U(\fg)/I)=\HC^G_{\partial \OO}(U(\fg)) \cap \HC^G(U(\fg)/I)$, this follows from (iii) of Proposition \ref{prop:propertiesofprim}. Form the quotient category
\begin{equation}\label{eq:barHC_Ug_I}
\overline{\HC}^G(U(\fg)/I) := \HC^G(U(\fg)/I)/\HC^G_{\partial}(U(\fg)/I).
\end{equation}
Every bimodule $\cB\in \HC^G(\cA)$ can be regarded as a $G$-equivariant Harish-Chandra $U(\fg)$-bimodule via the  co-moment map $\Phi: U(\fg) \to \cA$. The resulting bimodule is contained in the subcategory $\HC^G(U(\fg)/I) \subset \HC^G(U(\fg))$. This defines a forgetful functor
$$\Phi^*: \HC^G(\cA) \to \HC^G(U(\fg)/I).$$
This functor takes $\HC^G_{\partial}(\cA)$ to $\HC^G_{\partial \OO}(U(\fg))$ and thus descends to a functor (still denoted by $\Phi^*$)
$$\Phi^*: \overline{\HC}^G(\cA) \to \overline{\HC}^G(U(\fg)/I).$$
On the other hand, there is a restriction functor
$$\bullet_{\dagger}: \HC^G(U(\fg)) \to \HC^R(\cW),$$
see Section \ref{subsec:W}. By \cref{prop:propsofdagger}, this functor gives rise to a full embedding $\overline{\HC}^G(U(\fg)/I)\hookrightarrow \HC^R(\cW/\underline{J})$. To avoid notational confusion, we will denote this embedding by $\bullet_{\dagger'}$ for the remainder of this section. The next proposition shows that $\bullet_{\dagger}:\overline{\HC}^G(\cA) \to \Omega\modd$ and $\bullet_{\dagger'}:\overline{\HC}^G(U(\fg)/I) \to \HC^R(\cW/\underline{J})$ are compatible in a certain precise sense.

Choose $x\in\widetilde{\OO}$ as in Lemma \ref{lem:Adagger}. Note that $\Omega= R_x/R^{\circ}$. By Lemma \ref{lem:Adagger}(ii), we can choose an $\Omega$-stable ideal $J\subset \cW$ of codimension 1 such that $\underline{J}$ is the intersection of the $R$-conjugates of $J$, and furthermore 
$$\cA_\dagger\simeq \bigoplus_{r\in R/R_x}\cW/rJ.$$
The filtered algebra $\cA_{\dagger}$ comes equipped with a natural $R$-action as well as an $R$-equivariant algebra homomorphism $\cW/\underline{J}\hookrightarrow \cA_\dagger$. 

A representation $V$ of $\Omega$ is the same thing as an $R_x/R^{\circ}$-equivariant module for $\cW/J$. Consequently, there is a uniquely defined $R$-equivariant $\cA_{\dagger}$-module $\mathsf{B}_V$ with fiber $V$ over the summand $\cW/J$. We can view $\mathsf{B}_V$ as an object in $\HC^R(\cW/\underline{J})$.


\begin{lemma}\label{lem:twodaggers} 
The following diagram of functors commutes (up to isomorphism)
    \begin{center}
    \begin{tikzcd}
    \overline{\HC}^G(\cA) \ar[d, "\Phi^*"] \ar[r,"\bullet_\dagger"] & \Omega\modd \ar[d, "\mathsf{B}_\bullet"]\\
     \overline{\HC}^G(U(\fg)/I) \ar[r, "\bullet_{\dagger'}"]& \HC^R(\cW/\underline{J})
    \end{tikzcd}
\end{center}
\end{lemma}
\begin{proof}
We can assume in the proof that $G$ is semisimple and simply connected: the general case is obtained from this one by passing to full subcategories. In particular, the $G$-equivariance condition for a Harish-Chandra $\cA$-bimodule is vacuous, and $\Omega=\Gamma$. Also all categories involved are semisimple, so it is enough to establish an isomorphism on the level of objects. 

Recall, (\ref{eq:daggers_inverse}), that $\bullet_\dagger$ defines an equivalence $\overline{\HC}(\cA)\xrightarrow{\sim} \Gamma/\Gamma(\lambda)\operatorname{-mod}$.
Define $\breve{X}$ and $\breve{\cA}$ (for $\widetilde{X}$ and $\cA$) as in the discussion preceding Lemma \ref{lem:breve_A}. The functor $V\mapsto (\breve{\cA}\otimes V)^\Gamma$ is isomorphic to  $\bullet^\dagger$ by Lemma \ref{lem:breve_A}(i) and so is a quasi-inverse equivalence for $\bullet_\dagger$, see (\ref{eq:daggers_inverse}). Thus, it suffices to exhibit a functorial isomorphism (we omit $\Phi^*$ from the notation)
\begin{equation}\label{eq:fun_isomorphism}
[(\breve{\cA}\otimes V)^\Gamma]_{\dagger'}\xrightarrow{\sim} \mathsf{B}_V.
\end{equation}
Since the $\Gamma$-action on $\breve{\cA}$ fixes the image of the quantum comoment map,
we have a $\Gamma$-equivariant isomorphism $(\breve{\cA}\otimes V)_{\dagger'}
\xrightarrow{\sim} \breve{\cA}_{\dagger'}\otimes V$ of objects in $\HC^R(\cW/\underline{J})$.  
This gives rise to an isomorphism
\begin{equation}\label{eq:fun_isomorphism1}
[(\breve{\cA}\otimes V)^\Gamma]_{\dagger'}\xrightarrow{\sim} (\breve{\cA}_{\dagger'}\otimes V)^\Gamma.
\end{equation}
While $\breve{\cA}$ is not, in general, a quantization of $\CC[\breve{X}]$, it quantizes an algebra whose normalization is $\CC[\breve{X}]$, see (ii) of Lemma \ref{lem:breve_A}. So (i) and (ii) of Lemma \ref{lem:Adagger} hold for $\breve{\cA}$ nonetheless. Thus
$$\breve{\cA}_{\dagger'}\xrightarrow{\sim} \bigoplus_{r\in R/\breve{R}}\cW/rJ,$$ 
where we write $\breve{R}$ for the preimage in $R$ of $\Gamma(\lambda)\subset R/R^{\circ}$. From this description it is easy to see that the right hand side of (\ref{eq:fun_isomorphism1}) coincides with $\mathsf{B}_V$.

    %
\end{proof}

\section{Classification of unipotent ideals}\label{subsec:classificationideals}

Let $\widetilde{\mathbb{O}}$ and $\widecheck{\mathbb{O}}$ be $G$-equivariant covers of a common nilpotent orbit $\mathbb{O}$. As usual, let $\widetilde{X}= \Spec(\CC[\widetilde{\mathbb{O}}])$ and $\widecheck{X} = \Spec(\CC[\widecheck{\mathbb{O}}])$. If $f: \widetilde{\mathbb{O}} \to \widecheck{\mathbb{O}}$ is a morphism of $G$-equivariant covers, there is an induced morphism of spectra $\widetilde{X} \to \widecheck{X}$. The latter morphism is a finite $G$-equivariant cover, which is \'{e}tale over the subset $\widecheck{X}^{\mathrm{reg}} \subset \widecheck{X}$ by Lemma \ref{lem:leaftoleaf}. Recall, Definition \ref{defi:almost_etale}, that this morphism is \emph{almost \'{e}tale} if it is \'{e}tale over the locus
$$\widecheck{X}^{\mathrm{reg}} \cup \bigcup_k \fL_k \subset \widecheck{X}.$$
Define a binary relation $\succeq$ on the set of $G$-equivariant nilpotent covers as follows
$$\widetilde{\mathbb{O}} \succeq \widecheck{\mathbb{O}} \iff \exists \ f: \widetilde{\mathbb{O}} \to \widecheck{\mathbb{O}} \text{ such that } \widetilde{X} \to \widecheck{X} \text{ is an almost \'{e}tale cover}$$
\begin{definition}\label{def:equivalencerelation}
Let $\sim$ be the equivalence relation on the isomorphism classes of $G$-equivariant nilpotent covers defined by taking the symmetric closure of $\succeq$. Denote the equivalence class of $\widetilde{\mathbb{O}}$ by $[\widetilde{\mathbb{O}}]$.
\end{definition}
Below, we will show that unipotent ideals are parameterized by equivalence classes of covers. Our main result requires a bit of preparation. 

\begin{lemma}\label{lem:noE8}
Suppose $G$ is simple. Then $\widetilde{X}$ satisfies condition (\ref{eq:noE8condition}) except in the following case: $G$ is the (unique) simple group of type $E_8$ and $\widetilde{\OO}$ is the principal nilpotent orbit.
\end{lemma}

\begin{proof}
If $\widehat{\OO} \to \widetilde{\OO}$ is a $G$-equivariant cover, then the subgroups $\Gamma_k' \subset \mathrm{Sp}(2)$ corresponding to the codimension 2 leaves in $\widehat{X}$ are contained in the subgroups $\Gamma_k \subset \mathrm{Sp}(2)$ for $\widetilde{X}$, see (\ref{diag:galoiscovering}). Thus, we can reduce to the case when $\widetilde{\OO}$ is an orbit. The singularities corresponding to the codimension 2 leaves in the closures of nilpotent orbits were described in the classical cases by Kraft and Procesi (\cite[Main Thm]{Kraft-Procesi1981},\cite[Thm 2]{Kraft-Procesi}) and in the exceptional cases by Fu et.al. (\cite[Sec 13]{fuetal2015}). The lemma follows by inspection of the references above.
\end{proof}

We note that it is possible to give a more conceptual proof of Lemma \ref{lem:noE8}, at least in exceptional types. The proof requires some machinery which is developed in Section \ref{subsec:terminalizationcover}. We refer the reader to Remark \ref{rmk:noE8}. 

\begin{lemma}\label{lem:maximalelement}
Let $[\widetilde{\OO}]$ be an equivalence class of $G$-equivariant nilpotent covers. Then\index{cover!maximal}
\begin{itemize}
    \item[(i)] $[\widetilde{\mathbb{O}}]$ contains a unique maximal element $\widehat{\mathbb{O}}$.
    \item[(ii)] $\widehat{\OO}$ is Galois over every cover in $[\widetilde{\mathbb{O}}]$.
    \item[(iii)] There is an equivalence of categories $\overline{\HC}^G(\cA_0^{\widehat{X}}) \simeq \Vect$.
\end{itemize}
\end{lemma}

\begin{proof}
For (i), it suffices to show that $\widetilde{X}$ admits a unique maximal almost \'{e}tale $G$-equivariant cover $\widehat{X} \to \widetilde{X}$. For $G$ simply connected, this is a special case of Proposition \ref{prop:maximaletale}. The general case is proved in the same way. For similar reasons, $\widehat{X} \to \widetilde{X}$ is Galois, proving (ii).

For (iii), note that $G$ decomposes as a product $G = H \times E_8 \times ... \times E_8$, where $H$ is a reductive group containing no simple factors of type $E_8$ and $E_8$ is the (unique) simple group of type $E_8$. Indeed, every connected reductive group admits a central isogeny from a product of simple groups and a torus, and the simply connected group of type $E_8$ is adjoint. Thus, it suffices to prove (iii) separately for $G=H$ and $G=E_8$. 

First suppose $G=H$. Then $\widehat{X}$ satisfies property (\ref{eq:noE8condition}) by Lemma \ref{lem:noE8}. Let $\Omega=\pi_1^G(\widehat{\OO}), 
\Gamma=\pi_1(\widehat{\OO})$. Since $\widehat{\OO}$ is maximal in its equivalence class, $\Gamma=\Gamma(0)$, see Proposition \ref{prop:maximaletale} and Remark \ref{rmk:maximaletale}. Thus $\Omega=\Omega(0)$ and by Proposition \ref{prop:classificationeqvtbimods} there is an equivalence 
$$\overline{\HC}^G(\cA_0^{\widehat{X}}) \simeq \Omega/\Omega(0) \simeq \Vect,$$
proving (iii) in this case. 
Next suppose $G=E_8$. If $\OO$ is the principal nilpotent orbit, then $\pi_1(\OO)=1$. So $\widehat{\OO}=\OO$ and $\Gamma=1$. By the results of Section \ref{subsec:daggers}, there is an embedding $\overline{\HC}(\cA_0^{\widehat{X}})\hookrightarrow \Gamma\modd \simeq \Vect$. (iii) follows at once. If $\widehat{\OO}$ is not the principal orbit, then $\widehat{X}$ satisfies property (\ref{eq:noE8condition}), see again Lemma \ref{lem:noE8}. So we can argue as in the case of $G=H$. 
\end{proof}

\begin{prop}\label{prop:equivalencerelation}
Let $\widetilde{\mathbb{O}}$, $\dot{\mathbb{O}} \to \mathbb{O}$ be $G$-equivariant covers. Then
$$[\dot{\mathbb{O}}] = [\widetilde{\mathbb{O}}] \iff I_0(\dot{\mathbb{O}}) = I_0(\widetilde{\mathbb{O}}).$$
Moreover, if $\widetilde{\OO}$ is maximal in its equivalence class, and $\Pi:=\operatorname{Aut}_{\OO}(\widetilde{\OO})$, then there is a $G$-equivariant algebra isomorphism 
\begin{equation}\label{eq:full_invar_charact}
[(U(\fg)/I_0(\widetilde{\OO}))_\dagger]^\dagger\xrightarrow{\sim}(\cA_0^{\widetilde{X}})^\Pi.
\end{equation}
\end{prop}

\begin{proof}
Note that every $G$-equivariant nilpotent cover is a homogeneous space for the universal cover of the derived subgroup of $G$. Thus, we can reduce to the case when $G$ is semisimple and simply connected. The proof has several steps. 

{\it Step 1}. First we prove the easy implication
$[\dot{\mathbb{O}}] = [\widetilde{\mathbb{O}}] \Rightarrow I_0(\dot{\mathbb{O}}) = I_0(\widetilde{\mathbb{O}}).$
By Lemma \ref{lem:maximalelement}, we can assume $\dot{\mathbb{O}}$ is the maximal cover in $[\widetilde{\mathbb{O}}]$. Hence, there is a Galois cover $\dot{\mathbb{O}} \to \widetilde{\mathbb{O}}$ such that the induced map  $\dot{X} = \Spec(\CC[\dot{\OO}]) \to \Spec(\CC[\widetilde{\OO}]) = \widetilde{X}$ is almost \'{e}tale. By Proposition \ref{prop:almostetaleinvariants}, there is a filtered algebra embedding $i:\cA_0^{\widetilde{X}} \hookrightarrow \cA_0^{\dot{X}}$, and by Lemma \ref{lem:co-momentunique}, $i \circ \Phi_0^{\widetilde{X}} = \Phi_0^{\dot{X}}$. Thus
$$I_0(\widetilde{\mathbb{O}}) = \ker{(\Phi_0^{\widetilde{X}})} = \ker{(i \circ \Phi_0^{\widetilde{X}})} = \ker{(\Phi_0^{\dot{X}})} = I_0(\dot{\mathbb{O}}),$$
as desired.

{\it Step 2}. In Steps 3-7, we will prove the hard implication $I_0(\dot{\mathbb{O}}) = I_0(\widetilde{\mathbb{O}}) \Rightarrow [\dot{\mathbb{O}}] = [\widetilde{\mathbb{O}}]$. By Step 1, it suffices to show that every $G$-equivariant nilpotent cover $\widetilde{\mathbb{O}}$, maximal in its equivalence class, is determined by the ideal $I_0(\widetilde{\mathbb{O}})$. Our proof is somewhat indirect. We will first show that $I_0(\widetilde{\OO})$ determines a pair $(\widecheck{\OO},\widecheck{\cA})$ consisting of a nilpotent cover $\widecheck{\OO}$, covered by $\widetilde{\OO}$, and a filtered quantization $\widecheck{\cA}$ of $\CC[\widecheck{\OO}]$. The quantization $\widecheck{\cA}$ is constructed so that $I_0(\widetilde{\OO})$ coincides with the kernel of the (unique) quantum co-moment map $U(\fg) \to \widecheck{\cA}$ and $\widecheck{\OO}$ is the smallest cover of $\OO$ admitting a quantization with this property. Set $\widecheck{X}:=\operatorname{Spec}(\CC[\widecheck{\OO}])$. The pair $(\widecheck{X},\widecheck{\cA})$ determines a unique Galois cover $\breve{\OO}$ of $\widecheck{\OO}$ (see the discussion preceding Lemma \ref{lem:breve_A}) such that the image of $\overline{\HC}(\widecheck{\cA})$ in $\pi_1(\widecheck{\OO})\modd$ consists precisely of representations with trivial restriction to $\pi_1(\breve{\OO})$. To prove that 
$I_0(\widecheck{\mathbb{O}}) = I_0(\widetilde{\mathbb{O}}) \Rightarrow [\widecheck{\mathbb{O}}] = [\widetilde{\mathbb{O}}]$ we will show that $\breve{\OO}\simeq \widetilde{\OO}$.

{\it Step 3}.
Let $R$ denote the reductive part of the centralizer of $e\in \OO$, and let $\Gamma=\pi_1(\widetilde{\Orb})$.
Recall the functors $\bullet_\dagger: \HC^G(U(\fg)) \to \HC^R(\cW)$ and $\bullet^\dagger: \HC^R_{\mathrm{fin}}(\cW) \to \HC^G(U(\fg))$  and the maps $\bullet_{\dagger}: \mathrm{Prim}_{\overline{\OO}}(U(\fg)) \to \mathrm{Id}_{\mathrm{fin}}(\cW), \bullet^{\ddag}: \mathrm{Prim}_{\mathrm{fin}}(\cW) \to \mathrm{Prim}_{\overline{\OO}}(U(\fg))$ from Section \ref{subsec:W}.
By (iv) of Lemma \ref{lem:Adagger}, applied to $\cA_0^{\widetilde{X}}$, there is a $\Gamma$-invariant ideal $J\subset \cW$ of codimension 1 such that $I_0(\widetilde{\mathbb{O}})=J^\ddag$. Consider the finite index subgroup $\operatorname{Stab}_R(J)\subset R$
and let $\widecheck{\mathbb{O}} \to \mathbb{O}$ be the cover corresponding to $\operatorname{Stab}_R(J)/R^{\circ}\subset \pi_1(\OO)$. Thus, $\pi_1(\widecheck{\mathbb{O}}) \simeq \Stab_R(J)/R^{\circ} \subseteq R/R^{\circ} \simeq \pi_1(\mathbb{O})$. Since $J$ is $\pi_1(\widetilde{\mathbb{O}})$-stable, we get the inclusion $\pi_1(\widetilde{\mathbb{O}}) \subseteq \pi_1(\widecheck{\mathbb{O}})$. Thus $\widetilde{\mathbb{O}}$ is a cover of $\widecheck{\mathbb{O}}$. 

{\it Step 4}.
Recall from Section \ref{subsec:W} that $R^{\circ}$ preserves $J$. Let $\underline{J}$ denote the intersection of the $R$-conjugates of $J$. This is an $R$-stable ideal of finite codimension. Consider the Dixmier algebra $\widecheck{\cA}:=(\cW/\underline{J})^\dagger$. We claim that $\widecheck{\cA}$ is a filtered quantization of $\CC[\widecheck{\mathbb{O}}]$. This will be proved in steps 5-6. 

{\it Step 5}.
By Lemma \ref{lem:upper_dag_1dim}, $\gr \widecheck{\cA}$ embeds into $\mathbb{C}[\widecheck{\mathbb{O}}]$ as a $G$-stable subalgebra. This embedding intertwines the natural homomorphisms $\mathbb{C}[\fg^*] \to \gr \widecheck{\cA}$ and $\CC[\fg^*] \to \CC[\widecheck{\OO}]$, and hence the extension $\gr \widecheck{\cA} \hookrightarrow \mathbb{C}[\widecheck{\mathbb{O}}]$ is finite.
By Proposition \ref{prop:propsofdagger} (v), $\widecheck{\cA}_\dagger \simeq \cW/\underline{J}$. The dimension of the right hand side is $|R/\operatorname{Stab}_R(J)|$. Thus $\dim(\widecheck{\cA}_\dagger)$ coincides with the degree of the covering map $\widecheck{\mathbb{O}}\rightarrow \mathbb{O}$. In particular, the induced map $\Spec(\CC[\widecheck{\OO}]) \to \Spec(\gr(\widecheck{\cA}))$ is birational. A finite birational morphism onto a normal variety is an isomorphism. Thus,  $\gr(\widecheck{\cA}) \xrightarrow{\sim} \CC[\widecheck{\OO}]$ is equivalent to the normality of $\Spec(\gr \widecheck{\cA})$. The latter will be established in Step 6. 

{\it Step 6}.
By Step 5, $\widecheck{\mathbb{O}}$ embeds as an open $G$-orbit in $\operatorname{Spec}(\gr \widecheck{\cA})$ with boundary codimension $\geq 2$. Thus, $\Spec(\gr \widecheck{\cA})$ is normal if $\gr \widecheck{\cA}$ is Cohen-Macaulay. 
Note that the homomorphism $\cW\rightarrow (\cA_0^{\widetilde{X}})_\dagger$ factors through $\cW\rightarrow \cW/\underline{J}$. 
Also $\cW/\underline{J}$, regarded as an object in $\HC_{\mathrm{fin}}^R(\cW)$, embeds as a direct summand in 
$(\cA_0^{\widetilde{X}})_\dagger$.
Therefore, $\widecheck{\cA}=(\cW/\underline{J})^\dagger$, regarded as a filtered $U(\fg)$-bimodule, embeds as a direct summand in $[(\cA_0^{\widetilde{X}})_\dagger]^\dagger \simeq \cA_0^{\widetilde{X}}$ (so that $\cA_0^{\widetilde{X}}$ viewed as a filtered bimodule decomposes into the direct sum of $\widecheck{\cA}$ and another filtered bimodule). The last isomorphism follows from Lemma \ref{lem:Adagger}(iii). Passing to associated graded bimodules, we deduce that $\gr\widecheck{\cA}$, regarded as a $\CC[\fg^*]$-module, embeds as a direct summand into $\gr\cA_0^{\widetilde{X}} \simeq \CC[\widetilde{\OO}]$. The latter, and hence the former, is a Cohen-Macaulay module. So $\gr\widecheck{\cA}$ is a Cohen-Macaulay algebra. As noted in Step 5, this implies that $\widecheck{\cA}$ is a filtered quantization of $\CC[\widecheck{\mathbb{O}}]$.

{\it Step 7}.
Form $\breve{\OO}$ and $\breve{\cA}$ as in the discussion preceding Lemma \ref{lem:breve_A} (for $\operatorname{Spec}(\CC[\widecheck{\OO}])$
and its quantization $\widecheck{\cA}$). 
The construction in Step 6 shows that $\cA_0^{\widetilde{X}}$ is a cover of $\widecheck{\cA}$ in the sense of Definition \ref{defi:quant_cover}. Also, since $\widetilde{\OO}$ is maximal in its equivalence class, we have 
$\HC(\cA_0^{\widetilde{X}})\simeq \operatorname{Vect}$, see (iii) of Lemma \ref{lem:maximalelement}. Thus by Corollary \ref{cor:isotypic}, we have that $\breve{\OO}\simeq \widetilde{\OO}, \breve{\cA} \simeq \cA_0^{\widetilde{X}}$. This proves the implication $I_0(\dot{\mathbb{O}}) = I_0(\widetilde{\mathbb{O}})\Rightarrow 
[\dot{\mathbb{O}}] = [\widetilde{\mathbb{O}}]$.

{\it Step 8}. It remains to prove (\ref{eq:full_invar_charact}). 
By Step 7, $\widetilde{\OO}=\breve{\OO}$. By Proposition \ref{Prop:quant_covers}, $\widetilde{\OO}$ is a Galois cover of $\widecheck{\OO}$. Let $\Pi'\subset \Pi$ be the Galois group of $\widetilde{\OO} \to \widecheck{\OO}$. By Proposition \ref{Prop:quant_covers}, $\widecheck{\cA}=(\cA_0^{\widetilde{X}})^{\Pi'}$. In particular, there is an inclusion 
{$(\cA_0^{\widetilde{X}})^{\Pi}\hookrightarrow \widecheck{\cA}$} intertwining the quantum comoment maps. Consider the chain of inclusions
\begin{equation}\label{eq:inclusions1}U(\fg)/I\hookrightarrow (\cA_0^{\widetilde{X}})^{\Pi}\hookrightarrow \widecheck{\cA},\end{equation}
where $I:=I_0(\widetilde{\OO})$. Recall that the functor $\bullet_\dagger$ is exact by (i) of Proposition \ref{prop:propsofdagger}. Applying $\bullet_{\dagger}$ to (\ref{eq:inclusions1}) we obtain another chain of inclusions
$$(U(\fg)/I)_\dagger\hookrightarrow (\cA_0^{\widetilde{X}})^{\Pi}_\dagger\hookrightarrow \widecheck{\cA}_\dagger.$$
By the construction of $\widecheck{\cA}$ in Step 4, the composition of these inclusions is an isomorphism. On the other hand, by (i) of  Lemma \ref{lem:Adagger}, we have 
$$\dim(\cA_0^{\widetilde{X}})^{\Pi}_\dagger=\operatorname{deg}(\widetilde{\OO}/\Pi), \qquad 
\dim \widecheck{\cA}_\dagger=\operatorname{deg}(\widecheck{\OO}),$$
where `deg' denotes the degree of a finite cover of $\OO$.  {Since $\widecheck{\OO}=\widetilde{\OO}/\Pi'$, we}  conclude that $\Pi'=\Pi$. This proves (\ref{eq:full_invar_charact}).
\end{proof} 

The following is an immediate consequence of Proposition \ref{prop:equivalencerelation}.

\begin{theorem}\label{thm:classificationideals}
Let $\OO$ be a fixed nilpotent $G$-orbit. The passage from $\widetilde{\mathbb{O}}$ to $I_0(\widetilde{\mathbb{O}})$ defines a bijection between the following two sets
\begin{itemize}
    \item Equivalence classes of $G$-equivariant covers $\widetilde{\mathbb{O}} \to \mathbb{O}$.
    \item Unipotent ideals $I \subset U(\fg)$ with $V(I) = \overline{\mathbb{O}}$.
\end{itemize}
\end{theorem}

\begin{cor}\label{Cor:check_A_description} 
Suppose $\widetilde{\OO}$ is maximal in its equivalence class and let $\Pi:=\operatorname{Aut}_{\OO}(\widetilde{\OO})$. If $I_0(\widetilde{\OO})$ is a maximal ideal, then the quantum comoment map $U(\fg)\rightarrow \cA_0^{\widetilde{X}}$ induces an algebra algebra isomorphism
$$U(\fg)/I_0(\widetilde{\OO})\xrightarrow{\sim} (\cA_0^{\widetilde{X}})^{\Pi}.$$
\end{cor}
\begin{proof} 
 Lemma \ref{lem:adjunctionmorphismiso} applied to $\cB=U(\fg)/I_0(\widetilde{\OO})$
 implies $U(\fg)/I_0(\widetilde{\OO}) \xrightarrow{\sim} ([U(\fg)/I_0(\widetilde{\OO})]_\dagger)^\dagger$. On the other hand, $([U(\fg)/I_0(\widetilde{\OO})]_\dagger)^\dagger \xrightarrow{\sim} (\cA_0^{\widetilde{X}})^{\Pi}$ by (\ref{eq:full_invar_charact}). This completes the proof.
\end{proof}

\section{Classification of unipotent bimodules}\label{subsec:classificationbimods} Let $\widetilde{\OO}$ be a $G$-equivariant nilpotent cover, maximal in its equivalence class, and let $I:=I_0(\widetilde{\OO})$. Recall the category $\HC^G(U(\fg)/I)$ and its quotient $\overline{\HC}^G(U(\fg)/I)$ defined in (\ref{eq:HC_Ug_I}) and (\ref{eq:barHC_Ug_I}), respectively. Note that the irreducible objects in $\overline{\HC}^G(U(\fg)/I)$ are in bijection with $\mathrm{Unip}_{\widetilde{\mathbb{O}}}(G)$ --- this is a consequence of Definition \ref{def:unipotentbimods} together with (iii) of Proposition \ref{prop:propertiesofprim}.

Let $\Pi := \Aut_{\mathbb{O}}(\widetilde{\mathbb{O}})$, and $\widecheck{\mathbb{O}} := \widetilde{\mathbb{O}}/{\Pi}$. The main result of this section is the following theorem. 

\begin{theorem}\label{thm:classificationbimods}
Suppose $I$ is a maximal ideal. Then
there is an equivalence of monoidal categories
$$\Pi\modd \xrightarrow{\sim} \HC^G(U(\fg)/I), \qquad V \mapsto \Phi^*(\cA_0^{\widetilde{X}} \otimes V)^{\Pi}.$$
In particular, there is a bijection
$$\{\text{irreducible representations of } \Pi\} \xrightarrow{\sim} \mathrm{Unip}_{\widetilde{\mathbb{O}}}(G).$$
\end{theorem}

The proof of this theorem will require some preparation. 

By Proposition \ref{prop:canonicalauts}, the $\Pi$-action on $\CC[\widetilde{\mathbb{O}}]$ lifts to the canonical quantization $\cA_0^{\widetilde{X}}$. Let $\widecheck{X} = \Spec(\CC[\widecheck{\mathbb{O}}])$ and set 
$$\cA := (\cA_0^{\widetilde{X}})^{\Pi}.$$
An object $\cB \in \HC^G(\cA)$ can be regarded as a $G$-equivariant Harish-Chandra bimodule for $U(\fg)/I$ via the co-moment map $\Phi:=\Phi_{\delta}^{\widecheck{X}}: U(\fg) \to \cA_{\delta}^{\widecheck{X}}$. Recall, Corollary \ref{Cor:check_A_description}, that $\Phi$ induces an isomorphism  $U(\fg)/I\xrightarrow{\sim} \cA$.
So we get the forgetful functor 
\begin{equation}\label{eq:forget_fun_HC}
\HC^G(\cA)\rightarrow 
\HC^G(U(\fg)/I);
\end{equation} it is fully faithful. Combining (iii) of Lemma \ref{lem:maximalelement} with Corollary \ref{cor:isotypiceqvt}, we get an equivalence $\HC^G(\cA)\cong \Pi\operatorname{-mod}$. It remains to 
show that (\ref{eq:forget_fun_HC}) is essentially surjective, in other words, that 
\begin{itemize} \item[(*)] any Harish-Chandra bimodule over $U(\g)/I$ in the classical sense, equivalently, for the filtration induced from the PBW filtration on $U(\g)$, is also a Harish-Chandra bimodule with respect to the filtration on $\cA$ turning it into a filtered quantization of $\CC[\widecheck{\Orb}]$. 
\end{itemize}

We first highlight an easy special case, which is also the most common case. 

\begin{lemma}\label{Lem:quant_bimodule_easy}
Suppose $\widecheck{\Orb}=\Orb$, equivalently, $\widetilde{\Orb}$ is a Galois nilpotent cover. Then (*) holds. 
\end{lemma}
\begin{proof}
Take an irreducible object $\cB\in \HC^G(U(\g)/I)$ with its good filtration (compatible with the PBW filtration). Note that $\A$ is a filtered quantization of $\CC[\Orb]$. It follows that $\gr (U(\g)/I)$ is a closed subscheme of $\g^*$ whose reduced subscheme is $\overline{\Orb}$ and whose multiplicity on $\Orb$ is $1$. 

Consider the bimodule $(\cB_\dagger)^\dagger$ over $(\cA_\dagger)^\dagger$, where we use the restriction and extension functors from Section \ref{subsec:W}. By Lemma \ref{lem:adjunctionmorphismiso},
we have a $G$-equivariant algebra isomorphism $\cA\xrightarrow{\sim}(\cA_\dagger)^\dagger$ and an $\cA$-bimodule isomorphism 
$\cB\xrightarrow{\sim}(\cB_\dagger)^\dagger$. The algebra
$(\cA_\dagger)^\dagger$ is filtered. Thanks to the $G$-equivariant isomorphism $\cA\xrightarrow{\sim}(\cA_\dagger)^\dagger$ and Lemma 
\ref{lem:upper_dag_1dim}, we see that it is a filtration of a quantization of $\C[\Orb]$. 

Now we claim that $(\cB_\dagger)^\dagger$ is an $\cA$-bimodule with a good filtration.  
Recall, \cite[Section 3.4]{Losev3}, that 
$(\cB_\dagger)^\dagger$ is obtained as follows:
\begin{enumerate}
\item we take the Rees bimodule $\cB_\hbar$ of $\cB$ for the doubled PBW filtration, then $\cB$ carries a natural $\fg$-action given by $\xi.b:=\frac{1}{\hbar^2}[\xi,b]$; 
\item we complete $\cB_\hbar$ at the point $\chi\in \fg^*$ corresponding to $e$, denote the completion by $\cB_\hbar^{\wedge_\chi}$;
\item we take the locally finite vectors for $\fg$ in $\cB_\hbar^{\wedge_\chi}$;
\item then we take the locally finite vectors for the Kazhdan action of $\CC^\times$ on the resulting space, the result carries a grading; 
\item then we take the invariants for a suitable action of $R/R^\circ$, denote the resulting bimodule by $(\cB_{\hbar,\dagger})^\dagger$. 
\item and finally, we specialize $\hbar=1$ (which gives rise to the filtration). 
\end{enumerate}
In particular, this gives a filtered $\cA$-bimodule.
(2)-(5) make sense for $\gr\cB$ as well, let $([\gr\cB]_\dagger)^{\dagger}$ denote the resulting $G$-equivariant $\C[\Orb]$-module. Since $\gr\cB|_{\Orb}$ is scheme theoretically supported on $\Orb\subset \Orb^2$, we have $([\gr\cB]_\dagger)^{\dagger}=\Gamma(\gr\cB|_{\Orb})$, this follows from \cite[Section 3.2]{Losev3}.  And the construction yields an exact sequence
$$0\rightarrow (\cB_{\hbar,\dagger})^\dagger\xrightarrow{\hbar\cdot} (\cB_{\hbar,\dagger})^\dagger\rightarrow 
([\gr\cB]_\dagger)^{\dagger}.$$
This exact sequence shows that $\gr (\cB_\dagger)^\dagger\hookrightarrow \Gamma(\gr\cB|_{\Orb})$. Therefore, the $\cA$-bimodule filtration on $(\cB_\dagger)^\dagger$ is good.  

%
%
\end{proof}

We now proceed to the general case. Let $\widetilde{X}:=\operatorname{Spec}(\C[\widetilde{\Orb}])$, and $\Pi$ have the same meaning as before Theorem \ref{thm:classificationbimods}. Let $X:=\widetilde{X}/\Pi$, and $\A_\lambda=\tilde{\A}_0^\Pi$ be the quantization of $X$ (so that $\lambda$ is a unipotent parameter in the sense of Section \ref{SS_HC_different}). Let $I$ denote the kernel of $U(\fg)\twoheadrightarrow \A$.  Let $A$ denote the fundamental group of $\Orb$ and let 
$H_0\subset H\subset A$ denote the subgroups corresponding to the covers $\widetilde{X},X$, respectively, so that $\Pi=H/H_0$. Since $\widetilde{\Orb}$ is maximal in its equivalence class, we have $H=N_A(H_0)$. 

Consider the functor $\bullet_\dagger$ from Section \ref{subsec:W} and restrict it to $\HC^G(U(\fg)/I)$. By \cite[Lemma 5.2]{Losev4}, we have $ \A_\dagger=\C[A/H]$, so  the target category for this restriction is the category $\operatorname{Sh}^A((A/H)^2)$ of $A$-equivariant sheaves of finite dimensional vector spaces on the finite set $(A/H)^2$, monoidal with respect to convolution. 
Its simple objects are labelled by triples $(x,y,V)$, where $x,y\in A/H$ and $V$ is an irreducible representation of
$A_{(x,y)}$, up to the $A$-action. This monoidal category is rigid: on the level of labels of simples the duality sends $(x,y,V)$
to $(y,x,V^*)$. 

\begin{lemma}\label{Lem:image_dagger_restrictions}
We have the following:
\begin{enumerate}
\item The image of $\HC^G(U(\fg)/I)$ is closed under convolutions and taking duals. 
\item The simple $(1,1,V)$ lies in the image if and only if $H_0$ acts trivially on $V$ (here we write $1$ for $1H\in A/H$).
\end{enumerate}
\end{lemma}
\begin{proof}
(1) follows from (i) and (vi) of Proposition \ref{prop:propsofdagger}. To prove (2) we argue as follows. Let $\sigma$ be the simple $Q$-equivariant $\cW$-bimodule corresponding to the form $(1,1,V)$. Consider the bimodule $\sigma^\dagger$, in fact, it is a bimodule over $(\cA_\dagger)^\dagger=\cA$. The associated graded $\C[X]$-bimodule $\gr (\sigma^\dagger)$ embeds into the global sections of the $G$-equivariant vector bundle on the open $G$-orbit $\widecheck{\Orb}$ with fiber $V$, this is proved completely analogously to Lemma \ref{Lem:quant_bimodule_easy}. In particular, $\sigma^\dagger$ is HC as an $\cA$-bimodule. Recall that $\widetilde{\Orb}$ is maximal in its equivalence class. Combining Proposition \ref{prop:Gamma0} with the description of $\Gamma(\lambda)$ in Section \ref{subsec:Gammalambda}, we see that $\HC(\A^{\widetilde{X}}_0)\cong \operatorname{Vect}$. Now (2) follows from Corollary \ref{cor:isotypiceqvt}.
\end{proof}

We are going to define a subgroup $\hat{H}\subset A$ containing $H$ as follows. Let $S$ denote the union of supports of the images under $\bullet_\dagger$ of 
objects from $\operatorname{HC}^G(U(\fg)/I)$. Since this image is closed under convolution, $(x,y),(y,z)\in S\Rightarrow (x,z)\in S$.
Since it is closed under the duality, $(x,y)\in S\Rightarrow (y,x)\in S$. And, of course, $(x,x)\in S$ for all $x$. It follows 
that we have equivalence relation $\sim$ on $A/H$ defined by $x\sim y$ if and only if $(x,y)\in S$. It is $A$-invariant. Define 
$\hat{H}$ as the set of all elements  $a\in A$ such that $1\sim aH$. This is the required subgroup. 
Thanks to 2) of Lemma \ref{Lem:image_dagger_restrictions}, Theorem \ref{thm:classificationbimods} is equivalent to $\hat{H}=H$.

\begin{lemma}\label{Lem:intersection_conjugates}
For all $g\in \hat{H}$, we have $H_0\subset H\cap gHg^{-1}$.
\end{lemma}
\begin{proof}
Consider a simple object $L$ of $\operatorname{Sh}^A((A/H)^2)$ that lies in the image of $\HC^G(U(\fg)/I)$. Then its dual $L^*$ also lies in the image. Suppose $L$ is  labelled by $(1,g,V)$ for some $V$. Then $L^*$ is labelled by $(g,1,V^*)$. Consider the convolution of these objects. Then for any direct summand $U$ in 
$\operatorname{Ind}^H_{H\cap gHg^{-1}}(V\otimes V^*)$ we have that $(1,1,U)$ is in the image of  $\HC^G(U(\fg)/I)$. Hence $H_0$ acts trivially on $U$. This implies 
$H_0\subset H\cap gHg^{-1}$. 
\end{proof}  

Our next goal is to use results of  Section \ref{SS_HC_different} to prove the following proposition. 

\begin{prop}\label{Prop:self_normalizing}
The subgroup $H$ is self-normalizing in $\hat{H}$. 
\end{prop}
\begin{proof}
Assume that $g\in \hat{H}\setminus H$ normalizes $H$. Then $g$ gives rise to a $G$-equivariant automorphism of $X$ to be denoted by the same letter. Let $\lambda$ be the parameter of $\cA$, so that $\cA=\cA_\lambda$.

Take an irreducible object $\sigma$ in the image of $\HC^G(U(\fg)/I)$. Let  $\sigma$ be represented by a triple $(g,1,V)$,  and let 
$\cB$ be the corresponding irreducible HC $U(\fg)/I$-bimodule. Equip it with its filtration from the isomorphism $\cB\xrightarrow{\sim} (\cB_\dagger)^\dagger$. The associated graded is a $\C[X]$-bimodule but the two actions of $\C[X]$ do not coincide if $gH\neq 1$. Rather since the two actions of 
$\C[X]_\dagger$ on $\cB_\dagger$ differ by $g$, the same is true for the actions of $\C[X]$ on $\gr \cB$. It follows that if
we twist the left action of $\A_\lambda$ by $g$ (getting an action of $\A_{g\lambda}$), we get a HC bimodule over quantizations
of $\C[X]$, where $\A_{g\lambda}$ acts on the left, while $\A_\lambda$ acts on the right.  Proposition \ref{Prop:diff_unip_HC}
implies that $g\lambda=\lambda$. So $g$ gives an automorphism of $\A_\lambda$. This automorphism is $G$-equivariant contradicting the surjectivity of $U(\fg)\rightarrow \A_\lambda$. 
\end{proof}

\begin{proof}[Proof of Theorem \ref{thm:classificationbimods}]
Recall the subgroups $H_0\subset H\subset \hat{H}$ of $A$. Also consider $H_1:=\bigcap_{g\in \hat{H}}gHg^{-1}$. By Lemma 
\ref{Lem:intersection_conjugates}, $H_0\subset H_1$. By the construction, $H$ is the normalizer of $H_0$, while 
$H_1$ is normal in $\hat{H}$.  If the claim of the theorem is 
false, then the inclusions $H\subset \hat{H}, H_0\subset H_1$ are proper. Also, by Proposition \ref{Prop:self_normalizing}, $\hat{H}\setminus H$
does not contain elements normalizing $H$, so the inclusions $H_1\subset H\subset \hat{H}$ are proper. This can only hold if $H$ is not normal in $A$, which excludes the cases when $A$ is abelian or when it is a central extension of an abelian group (because $H$ has to contain the center).

The remaining cases are $A=S_3,S_3\times \mathbb{Z}/2\mathbb{Z}, S_4,S_5$. The former two groups do not have chains of subgroups
$\{1\}\subsetneq H_0\subsetneq H_1\subsetneq H\subsetneq \hat{H}$. 

Consider the case $A=S_4$. The order is $24$, the product of  four prime factors (with multiplicities). It follows that $\hat{H}=S_4$, while 
the order of $H$ is either $12$ or $8$. The only possibility is when the nontrivial element in $H_0$ has cycle type $(2,2)$,
$H_1$ is the Klein 4-subgroup, and $H$ is the normalizer of $H_0$, the semidirect product of $H_1$ with a transposition. The only orbit $\OO$ with $A=S_4$ in a simple Lie algebra is the orbit $F_4(a_3)$ in $F_4$. By \cite[Table 2]{Westaway2020}, the cover $\widetilde{\OO}$ of $\OO$ corresponding to $H_0$ (denoted in loc. cit. by $\mathrm{tw}(S_2)$) as well as the universal cover $\widehat{\OO}$ of $\OO$ are birationally rigid. It follows from Corollary \ref{cor:criterionbirigid} that $\Spec(\CC[\widetilde{\OO}])$ and $\Spec(\CC[\widehat{\OO}])$ have no codimension 2 leaves, and hence that the map $\Spec(\CC[\widehat{\OO}]) \to \Spec(\CC[\widetilde{\OO}])$ is almost \'{e}tale. Therefore $\widetilde{\OO}$ is not maximal in its equivalence class. 

Consider finally the case $A=S_5$. The only nontrivial normal subgroup is the alternating group $\mathfrak{A}_5$. So if $\hat{H}=S_5$, then $H_1=\mathfrak{A}_5$, contradiction with $H_1\subsetneq H\subsetneq \hat{H}$.  
It follows, in particular, that the cardinalities of $H_0,H_1,H,\hat{H}$ have $1,2,3,4$ prime factors with multiplicities. 
Consider the following cases:

1) $|H_0|=5$. Here the normalizer is $H_0$, a contradiction.

2) $|H_0|=3$. Here the normalizer is $H=S_3\times S_2$, a maximal subgroup in $S_5$. A contradiction.

3) $|H_0|=2$. If the nontrivial element in $H_0$ is a transposition, then $H=N_A(H_0)\cong S_2\times S_3$. This subgroup is maximal, so 
$\hat{H}=S_5$, a contradiction. Now assume the nontrivial element in $H_0$ has cycle type (2,2,1). The normalizer has order $8$
and the situation is similar to the case of $S_4$, in particular, $\hat{H}=S_4$. 

What remains to prove is that the cover $\widetilde{\OO}$ corresponding to this $H_0$ (whose only nontrivial element is, say, $(12)(34)$)
is equivalent to the universal cover $\widehat{\OO}$. But by \cite[Table 5]{Westaway2020}, both $\widetilde{\OO}$ and the $\widehat{\OO}$ are birationally rigid (in loc. cit., $H_0$ is denoted by $\mathrm{tw}(S_2)$). Arguing as in the case of $A=S_4$, we conclude that $\widehat{\OO}$ is equivalent to $\widetilde{\OO}$. This completes the proof.
\end{proof}

\section{A proof of a conjecture of Vogan}\label{subsec:Gtypes}

Recall from Section \ref{subsec:HCbimodsclassical} that every Harish-Chandra bimodule $\B \in \HC^G(U(\fg))$ can be regarded as a $G$-representation via the adjoint action of $\fg$. If $\B$ is unipotent, then this $G$-representation is of a very special form.

\begin{prop}\label{prop:Gtypes}
Assume $I_0(\widetilde{\mathbb{O}})$ is a maximal ideal. Let $\B \in \unip_{\widetilde{\mathbb{O}}}(G)$ and let $e \in \mathbb{O}$. Then there is a finite-dimensional representation $\chi$ of $G_e/G_e^{\circ}$ such that
$$\B \simeq_G \mathrm{AlgInd}^G_{G_e} \chi.$$
\end{prop}

\begin{proof}
Assume $\widetilde{\mathbb{O}}$ is maximal in its equivalence class, and let $\Pi =\operatorname{Aut}_{\OO}(\widetilde{\OO})$. By Theorem \ref{thm:classificationbimods}, there is a uniquely defined irreducible $V \in \Pi\modd$ such that $\cB \simeq \Phi^*(\cA_0^{\widetilde{X}} \otimes V)^{\Pi}$. The filtration on $\cA_0^{\widetilde{X}}$ induces a good filtration on $\B$, and there is a natural isomorphism of graded $G$-equivariant $\CC[\mathbb{O}]$-modules
\begin{equation}\label{eq:grBdescription}
\gr(\B) \simeq \gr(\cA_0^{\widetilde{X}} \otimes V)^{\Pi} \simeq (\CC[\widetilde{\mathbb{O}}] \otimes V)^{\Pi}.\end{equation}
Fix $\widecheck{\mathbb{O}}$ as in Section \ref{subsec:classificationbimods}. Choose $x \in \widecheck{\mathbb{O}}$ over $e$ and $y \in \widetilde{\mathbb{O}}$ over $x$. Thus, $G_y$ is a normal subgroup of finite-index in $G_x$ and $\Pi \simeq G_x/G_y$. Regard $V$ as an irreducible representation of $G_x$, trivial on $G_y$, and form the $G$-equivariant vector bundle $\mathcal{V} := G \times^{G_x} V \to G/G_x \simeq \widecheck{\mathbb{O}}$. Then $(\CC[\widetilde{\mathbb{O}}] \otimes V)^{\Pi}$ can be identified with $\Gamma(\widecheck{\mathbb{O}},\mathcal{V})$. As $G$-representations
$$\Gamma(\widecheck{\mathbb{O}},\mathcal{V}) \simeq_G \mathrm{AlgInd}^G_{G_x}V \simeq_G \mathrm{AlgInd}^G_{G_e}(\mathrm{AlgInd}^{G_e}_{G_x}V).$$
Now set $\chi = \mathrm{AlgInd}^{G_e}_{G_x}V$ to complete the proof.
\end{proof}

\begin{rmk}
In Chapter \ref{sec:duality}, we will show that all special unipotent ideals $I_{\mathrm{max}}(\frac{1}{2}h^{\vee})$ are unipotent. Thus, the conclusion of Proposition \ref{prop:Gtypes} holds for all special unipotent bimodules. In particular, Proposition \ref{prop:Gtypes} provides an affirmative answer to a conjecture of Vogan (\cite[Conj 12.1]{Vogan1991}) in the case of complex groups. Special cases of this conjecture were previously established by \cite{masonbrown2018}, \cite{Barbasch_dualpairs}, and \cite{WongLVB}.
\end{rmk}

Proposition \ref{prop:Gtypes} is the `complex case' of a more general conjecture about irreducible representations of real reductive Lie groups, see \cite[Sec 12]{Vogan1991}. We conclude this section with a discussion of this more general conjecture. A proof of this more general conjecture was proved by Dougal Davis and the second-named author in \cite{DavisMasonBrown}.

Let $K$ be a symmetric subgroup of $G$, i.e. the fixed point locus of a an algebraic involution. Suppose $\mathcal{B}$ is an irreducible $(\fg,K)$-module such that $V(\mathrm{Ann}(\mathcal{B})) = \overline{\OO}$. Then the associated variety of $\mathcal{B}$ is a union of components of $\overline{\OO} \cap (\fg/\mathfrak{k})^*$, and each such component is the closure of a $K$-orbit on $\OO \cap (\fg/\mathfrak{k})^*$, see \cite[Section 4]{Vogan1991}. Now suppose that the associated variety of $\mathcal{B}$ contains a $K$-orbit $\OO_K \subset \OO \cap (\fg/\mathfrak{k})^*$ such that
\begin{equation}\label{eq:codim_conditionK}\codim(\partial \OO_K,\overline{\OO}_K) \geq 2. 
\end{equation} 
Then the associated variety of $\mathcal{B}$ is \emph{equal} to the closure of $\OO_K$, see \cite[Theorem 4.6]{Vogan1991} (we note that Condition \ref{eq:codim_conditionK} is automatic if $\codim(\partial \OO, \overline{\OO})\geq 4$).

Finally suppose that $\Ann(\mathcal{B}) = I_0(\OO)$. In this situation, one can associate to $\mathcal{B}$ a certain $K$-equivariant vector bundle $\mathcal{L}_K$ on $\OO_K$, see \cite[Sec 2]{Vogan1991}. In fact, this vector bundle has the additional structure of a twisted local system (with `twist' corresponding to a square root of the canonical bundle on $\OO_K$, see \cite[Thm 8.7]{Vogan1991}). In this setting, Vogan's conjecture is as follows.

\begin{conj}\label{conj:realVoganconj}
Let $\mathcal{B}$ be an irreducible $(\fg,K)$-module such that $\Ann(\mathcal{B}) = I_0(\OO)$. Suppose that the associated variety of $\mathcal{B}$ is $\overline{\OO}_K$, where $\OO_K$ is a $K$-orbit on $\OO \cap (\fg/\mathfrak{k})^*$ satisfying Condition \ref{eq:codim_conditionK}. Then there is an isomorphism of $K$-representations
$$\cB \simeq_K \Gamma(\OO_{\cB}, \mathcal{L}_{\cB})$$
\end{conj}
As mentioned above, this conjecture was recently proved in \cite[Theorem 5.19]{DavisMasonBrown}.

\chapter{Parabolic Induction of Hamiltonian quantizations}\label{sec:nilpquant}

Let $G$ be a connected reductive algebraic group. In Chapter \ref{sec:unipotent}, we defined, for each $G$-equivariant nilpotent cover $\widetilde{\OO}$, a completely prime primitive ideal $I_0(\widetilde{\OO}) \subset U(\fg)$, Definition \ref{def:unipotentideals}, and a finite set of irreducible Harish-Chandra bimodules $\unip_{\widetilde{\OO}}(G)$, Definition \ref{def:unipotentbimods}. Recall that in Definition \ref{def:equivalencerelation} we have defined an equivalence relation on the set of nilpotent covers.

We proved that the set $\unip_{\widetilde{\OO}}(G)$ is in natural bijection with irreducible representations of the (finite) automorphism group of the maximal cover in the equivalence class $[\widetilde{\OO}]$, see Theorem \ref{thm:classificationbimods}, and that all bimodules in $\unip_{\widetilde{\OO}}(G)$ are, as $G$-representations, of the form conjectured by Vogan {under the assumption that $I_0(\widetilde{\OO})$ is maximal}, see Proposition \ref{prop:Gtypes}. However, three properties of unipotent ideals and bimodules were left more or less unaddressed, namely:
\begin{itemize}
    \item All unipotent bimodules are unitary.
    \item All unipotent ideals are maximal.
    \item All special unipotent bimodules are unipotent.
\end{itemize}
Compare to Desiderata (1), (3), and (6) from Section \ref{subsec:desiderata}. We will verify these properties {for classical groups} in Chapters \ref{sec:centralchars}-\ref{sec:unipbimod}. To do so, we will need a method for computing the infinitesimal character $\gamma_0(\widetilde{\OO}) \in \fh^*/W$ of a unipotent ideal $I_0(\widetilde{\OO})$. This method will be developed in Chapter \ref{sec:centralchars}. In this chapter, we establish some preliminary facts about filtered quantizations of nilpotent covers which are needed for these computations.

We will summarize the method here in order to motivate what follows. Assume for simplicity that $\widetilde{\OO}$ is birationally rigid, see Section \ref{subsec:birationalinduction},  and maximal in its equivalence class (it is not particularly difficult to reduce the computation to this special case and we will do so in Section \ref{subsec:unipotentcentralchars}). These assumptions guarantee that $\widetilde{\OO}$ is the universal $G$-equivariant cover of $\OO$. Taking $\Aut_{\OO}(\widetilde{\OO})$-invariants in $\cA_0^{\widetilde{X}}$, we get a Hamiltonian quantization $\cA_{\epsilon}^X$ of $\CC[\OO]$ corresponding to a special parameter $\epsilon \in \fP^X$ constructed in Example \ref{ex:barycentersymplectic}. Moreover,  $I_{\epsilon}(\OO)=I_0(\widetilde{\OO})$. Choose a Levi subgroup $L \subset G$ and a birationally rigid $L$-orbit $\OO_L$ such that $\OO=\mathrm{Bind}^G_L \OO_L$. Since $\OO_L$ is birationally rigid, $\CC[\OO_L]$ admits a unique filtered quantization, $\cA_0^{X_L}$. In Section \ref{subsec:inductionquantizations}, we will define the notion of parabolic induction for Hamiltonian quantizations. For every `induction parameter' $\lambda \in \fX(\fl)$, we construct an `induced' quantization of $\CC[\OO]$, denoted $\Ind^G_L \cA_{\lambda}^{X_L}$. Every Hamiltonian quantization of $\CC[\OO]$ arises in this fashion. More precisely, there is a linear isomorphism $\eta:\fX(\fl)\xrightarrow{\sim}\overline{\fP}^X$ such that $\Ind^G_L \cA_{\lambda}^{X_L}$ coincides with the Hamiltonian quantization $\cA^X_{\eta(\lambda)}$. Given an induction parameter $\lambda \in \fX(\fl)$ and the infinitesimal character of $I_0(\OO_L)$, it is easy to determine the infinitesimal character of $I_{\eta(\lambda)}(\OO)$. So we need to determine the element $\eta^{-1}(\epsilon)\in \fX(\fl)$. Computing $\eta$ is a nontrivial problem, to which the bulk of this chapter is devoted. In Sections \ref{subsec:descriptionpartial}-\ref{subsec:identification}, we solve this problem under certain conditions on $\OO$. This reduces the problem of computing unipotent infinitesimal characters to the case of birationally rigid orbits. By a similar argument, one can further reduce to the case of rigid nilpotent orbits. For such orbits, we can compute the infinitesimal characters using the results of McGovern (\cite{McGovern1994}) and Premet (\cite{Premet2013}). 


We now state the results of these computations with references to the relevant parts of the text.

{\it Case 1}. Let $G=\operatorname{SL}(n)$. Only the zero orbit is (birationally) rigid, (i) of Proposition \ref{prop:birigidorbitclassical}. Moreover, the universal cover of an orbit $\Orb$ is birationally rigid if and only if the partition $p$ of $\Orb$ is of the form $p=(d^m)$, where $d,m$ are positive integers satisfying $md=n$, see Proposition \ref{prop:birigidcoverA}. The corresponding orbit $\Orb$ is birationally induced from the zero orbit in the Levi $L=S(\operatorname{GL}(m)^d)$, where, recall $S$ indicates that we take the intersection with $\operatorname{SL}(n)$, see (i) of Proposition \ref{prop:Lforbirigidcover}. 
For this cover $\widetilde{\Orb}$, we have 
$$\gamma_0(\widetilde{\mathbb{O}}) = (\frac{n-1}{2d}, \frac{n-3}{2d}, ..., \frac{1-n}{2d}) = \frac{\rho}{d}.$$ 
See (i) of Proposition \ref{prop:centralcharacterbirigidcover} for this computation.

{\it Case 2}. Let $G=\operatorname{SO}(n)$ or $\operatorname{Sp}(2n)$. Let $\Orb$ be a nilpotent orbit in $\fg$ and let $p=(p_1,p_2,\ldots)$ be the corresponding partition of $n$ or $2n$. Then $\Orb$ is birationally rigid if and only if 
\begin{itemize}
\item $p_{i}\leqslant p_{i+1}+1$ for all $i$, 
\item and, in the case when $\fg=\mathfrak{so}(n)$, $p$ is not of the form $(2^m,1^2)$ for some $m$ (in which case $n$ is even).
\end{itemize}
See (ii) and (iii) of Proposition \ref{prop:birigidorbitclassical}. It is possible to completely describe all orbits with birationally rigid covers, however we do not do this here. Instead we consider a broader class of covers, `2-leafless' ones, for which $\operatorname{Spec}(\CC[\widetilde{\Orb}])$ is smooth in codimension $2$, see 
Definition \ref{defi:2_leafless}. These turn out to have simpler classification. It turns out, Proposition \ref{prop:nocodim2leaves}, that $\Orb$ admits a 2-leafless cover if and only if $p_i\leqslant p_{i+1}+2$ for all $i$ and the inequality is strict whenever $p_i$ is even for $\fg=\mathfrak{so}(n)$ and odd for $\fg=\mathfrak{sp}(n)$. One can completely describe the birational induction data for the orbits $\Orb$ admitting birationally rigid covers, although this is technical, see (ii) and (iii) of Proposition \ref{prop:Lforbirigidcover}. 

We now describe the computation of $\gamma_0(\widetilde{\Orb})$ for a birationally rigid $G$-equivariant cover $\widetilde{\Orb}$. This requires some combinatorial preparation. The discussion of the previous paragraph shows that if $\Orb$ admits a birationally rigid cover, then each part in the transposed partition $p^t$ occurs with multiplicity not exceeding $2$. We write $x$ and $y$ for the subpartitions of $p^t$ consisting of all parts with multiplicities $1$ and $2$, respectively. 
We then can form new partitions $f_?(x)$ and $g(y)$ for $?=B$ or $C$ as explained in Definitions \ref{def:fBfC} and \ref{def:xy}. 
Finally, if $q=(q_1,q_2,\ldots)$ is a partition of $n$, we can form an element $\rho^+(q)\in (\frac{1}{2}\ZZ_{\geqslant 0})^{\lfloor n/2\rfloor}$ by appending the positive entries of the form $\frac{q_i-j}{2}$ one for each possible $i$ and $j$ and a suitable number of $0$'s to get $\lfloor n/2\rfloor$ entries. Then we take $q=f_?(x)\cup g(y)$, where $?=B$ for $\fg=\mathfrak{so}(n)$ and $?=C$ for $\fg=\mathfrak{sp}(n)$, and the union sign indicates the union of the partitions so that the multiplicity of any part $\ell$ in $q$
is the sum of the multiplicities of $\ell$ in $f_?(x)$ and $g(y)$. Our final result on computing the infinitesimal characters of the unipotent ideals associated to $G$-equivariant birationally rigid covers is that 
$\gamma_0(\widetilde{\Orb})=\rho^+(q)$ for $q$ as above, see (ii) and (iii) of Proposition
\ref{prop:Lforbirigidcover}. Here $\rho^+(q)$ is viewed as an element of $\fh^*$, whose coordinates in a natural basis are the entries of $\rho^+(q)$,  since we only care about the $W$-orbit of this element, the order of coordinates turns out to be not important. 

{\it Case 3}. $G=\operatorname{Spin}(n)$. Compared to the case of $G=\operatorname{SO}(n)$, here one has more birationally rigid covers, hence more unipotent ideals and corresponding infinitesimal characters. In particular, we will see (Proposition \ref{prop:nocodim2leavesspin}) that a nilpotent orbit $\OO$ with partition $p$ admits a 2-leafless cover, which is not $\operatorname{SO}(n)$-equivariant, if and only if the following conditions are satisfied:
\begin{itemize}
\item $p$ is rather odd, i.e., every odd part occurs with mutliplicity $1$,
\item $p_i\leqslant p_{i+1}+1$ if $p_i$ is even, and $p_i\leqslant p_{i+1}+4$ is $p_i$ is odd,
\item $p_i\neq p_{i+1}+3$ for all $i$.
\end{itemize}
One can also describe the induction data for $\Orb$, see Proposition \ref{prop:nocodim2leavesspinLM}, and reduce the compututation of $\gamma_0(\widetilde{\Orb})$ to the case of $\operatorname{SO}(n)$-equivariant covers,  Section \ref{subsec:centralcharspin}, both results are somewhat technical so we do not provide details here. 

{\it Case 4}. $G$ is exceptional. Here we only provide partial results, see \cite{MBM} for the complete calcluations. We do not give a complete classification of birationally rigid orbits or orbits admitting birationally rigid (or 2-leafless) covers. Instead, we observe that the computation of $\gamma_0(\widetilde{\Orb})$ is relatively easy when the codimension $2$ singularities of the variety $\operatorname{Spec}(\CC[\OO])$ are all of type $A_1$ (this is a typical situation across both classical and exceptional types, the case of the birationally rigid covers for $G=\operatorname{SL}(n)$ provides a notable exception). We determine all orbits $\Orb$ in exceptional Lie algebras admitting 2-leafless covers, Proposition \ref{prop:almost all A1}.  There are four that have dimension $2$ singularities different from $A_1$, three in $E_6$ and one in $E_8$, and for all of them the universal cover is birationally rigid. We describe the codimension $2$ singularities for each of these orbits.  

We will also perform some computations of infinitesimal characters $\gamma_0(\widetilde{\OO})$ for birationally rigid covers for exceptional groups in Sections \ref{subsec:centralcharexceptional} and
\ref{subsec:dualdistinguished}.

\section{Partial resolutions of nilpotent covers}\label{subsec:picard}

Choose a Levi subgroup $L \subset G$ and an $L$-equivariant cover $\widetilde{\mathbb{O}}_L$ such that $\widetilde{\mathbb{O}} = \mathrm{Bind}^G_L \widetilde{\mathbb{O}}_L$. Fix a parabolic subgroup $P \subset G$ with Levi decomposition $P= LN$, and form the $G$-equivariant fiber bundle
$$\pi: \widetilde{Y} := G \times^P (\widetilde{X}_L \times \fp^{\perp}) \to G/P,$$
see Section \ref{subsec:birationalinduction}. There is a proper map $\widetilde{Y} \to \overline{\mathbb{O}}$, which factors through a partial resolution $\rho: \widetilde{Y} \to \widetilde{X}$. 

In this preparatory section, we will compute the Picard group and class group of $\widetilde{Y}$. We begin with a general observation.

\begin{lemma}\label{lem:normalaffinecone}
Let $V$ be a normal affine cone. Then $\Pic(V) = 0$.
\end{lemma}

\begin{proof}
Choose a line bundle $\mathcal{L} \in \Pic(V)$. Since $V$ is normal, $\mathcal{L}$ can be made $\CC^{\times}$-equivariant, see e.g. \cite[Proposition 2.4]{Knop_linear}. Hence, $\mathcal{L}$ has the structure of a graded $\CC[V]$-module. Let $0 \in V$ be the fixed point for the $\CC^{\times}$-action and consider the maximal ideal $\fm_0 \subset \CC[V]$. Since $\mathcal{L}$ is a line bundle, the fiber $\mathcal{L}_0 = \mathcal{L}/\fm_0\mathcal{L}$ is one-dimensional. Choose a nonzero element $\overline{s} \in \mathcal{L}_0$. This element lifts to a section $s \in \mathcal{L}$ which generates $\mathcal{L}$ as a $\CC[V]$-module
by the graded version of Nakayama's lemma. The action map defines a surjection $\CC[V] \to \mathcal{L}$, which is manifestly injective, since $\mathcal{L}$ is torsion-free.
\end{proof}

Let $\fX(L)$ denote the character lattice of $L$. If $\chi \in \fX(L)$, write $\mathcal{L}_{G/P}(\chi)$ for the $G$-equivariant line bundle on $G/P$ with fiber $\chi$ at $[P] \in G/P$. 

\begin{prop}\label{prop:descriptionofpic}
The following are true:
\begin{itemize}
    \item[(i)] The map $\pi^*: \Pic(G/P) \rightarrow \Pic(\widetilde{Y})$ is an isomorphism. If $G$ is semisimple and simply connected, there is a further isomorphism $\mathcal{L}_{G/P}: \fX(L) \xrightarrow{\sim} \Pic(G/P)$.

    \item[(ii)] There is a short exact sequence
    $$0 \to \Pic(\widetilde{Y}) \overset{\mathrm{div}}{\to} \Cl(\widetilde{Y}) \to \Pic(\widetilde{\mathbb{O}}_L) \to 0,$$
    where $\Cl$ denote the class group, and $\mathrm{div}$ is the natural map from line bundles to divisors.
\end{itemize}
\end{prop}

\begin{proof}
First, we prove (ii). Let $C \subset G/P$ be the open Bruhat cell and let $\widetilde{C} = \pi^{-1}(C)$. {Writing $N^-$ for the unipotent radical of the opposite parabolic $P^-$ of $P$}, there is an $N^-$-invariant isomorphism 
\begin{equation}\label{eq:decompC}\widetilde{C} \simeq N^- \times \widetilde{X}_L \times \fp^{\perp}\end{equation}
Taking class groups, we get an identification $\Cl(\widetilde{C}) \simeq \Cl(\widetilde{X}_L)$. Since $\codim(\widetilde{X}_L - \widetilde{\mathbb{O}}_L,\widetilde{X}_L) \geq 2$, $\Cl(\widetilde{X}_L) \simeq \Cl(\widetilde{\mathbb{O}}_L)$ and since $\widetilde{\mathbb{O}}_L$ is smooth, $\Cl(\widetilde{\mathbb{O}}_L) \simeq \Pic(\widetilde{\mathbb{O}}_L)$. Thus, $\Cl(\widetilde{C}) \simeq \Pic(\widetilde{\mathbb{O}}_L)$. 

Since $\widetilde{C}$ is an open subset of $\widetilde{Y}$ there is a surjective map $\mathrm{res}: \Cl(\widetilde{Y}) \to \Cl(\widetilde{C})$. We will show that the following sequence is exact
\begin{equation}\label{eq:sesZOM}0 \to \Pic(\widetilde{Y}) \overset{\mathrm{div}}{\to} \Cl(\widetilde{Y}) \overset{\mathrm{res}}{\to} \Cl(\widetilde{C}) \to 0.\end{equation}
First, note that $\mathrm{div}: \Pic(\widetilde{Y}) \to \Cl(\widetilde{Y})$ is injective, since $\widetilde{Y}$ is normal. Next, observe that the composition 
$$\Pic(\widetilde{Y}) \overset{\mathrm{div}}{\to} \Cl(\widetilde{Y}) \overset{\mathrm{res}}\to \Cl(\widetilde{C})$$
coincides with the composition 
$$\Pic(\widetilde{Y}) \overset{\mathrm{res}}{\to} \Pic(\widetilde{C}) \overset{\mathrm{div}}{\to} \Cl(\widetilde{C}),$$

By (\ref{eq:decompC}), $\Pic(\widetilde{C}) \simeq \Pic(\widetilde{X}_L)$, and by Lemma \ref{lem:normalaffinecone}, $\Pic(\widetilde{X}_L) \simeq 0$. Thus, $\Pic(\widetilde{C}) \simeq 0$ and  $\mathrm{Im}(\mathrm{div}) \subseteq \ker(\mathrm{res})$. 

Next we prove the reverse inclusion: $\mathrm{Im}(\mathrm{div}) \supseteq \ker(\mathrm{res})$. Let $D=G/P-C$ and $\widetilde{D}=\pi^{-1}(D) = \widetilde{Y} -\widetilde{C}$. Denote the irreducible components of $D$ by $D_1,...,D_r$. Note that each $D_i$ is of codimension 1 and the set $\{D_1,...,D_r\}$ is parameterized by simple roots for $\fg$ which are not contained in $\fl$, see \cite[Thm 21.14]{Borel}. For each $D_i$, let $\widetilde{D}_i = \pi^{-1}(D_i)$. Since $\pi: \widetilde{Y} \to G/P$ is a fiber bundle with irreducible fibers, $\widetilde{D}_1,...,\widetilde{D}_r$ are the irreducible components of $\widetilde{D}$. If $[E] \in \ker(\mathrm{res})$, then $[E] = \sum_{i=1}^r c_i [\widetilde{D}_i]$ for integers $c_1,...,c_r$. Hence
$$[E] = \sum_{i=1}^r c_i[\widetilde{D}_i] = \mathrm{div}(\pi^*\cO(\sum_{i=1}^r c_i[D_i])),$$
where $\cO(\sum_{i=1}^r c_i[D_i])) \in \Pic(G/P)$ is the line bundle corresponding to the class $\sum_{i=1}^r c_i[D_i] \in \Cl(G/P)$. In particular, $\ker(\mathrm{res}) \subseteq \mathrm{Im}(\mathrm{div})$. Thus, we have shown that (\ref{eq:sesZOM}) is exact. Using the isomorphism $\Cl(\widetilde{C}) \simeq \Pic(\widetilde{\mathbb{O}}_L)$ described in the first paragraph of the proof, we get a short exact sequence
$$0 \to \Pic(\widetilde{Y}) \overset{\mathrm{div}}{\to} \Cl(\widetilde{Y}) \to \Pic(\widetilde{\mathbb{O}}_L) \to 0,$$
as desired.

Next, we prove (i). For $G$ semisimple, the isomorphism $\mathcal{L}_{G/P}: \fX(L) \simeq \Pic(G/P)$ is standard, see e.g. \cite[Prop 3.2]{Knop1989}. It remains to show that $\pi^*: \Pic(G/P) \xrightarrow{\sim} \Pic(\widetilde{Y})$. 

Restricting along the zero section $G/P\subset \widetilde{Y}$ defines a left inverse for $\pi^*$. In particular, $\pi^*$ is injective. To prove it is surjective, choose a line bundle $\mathcal{L} \in \Pic(\widetilde{Y})$. By the proof of the inclusion $\mathrm{Im}(\mathrm{div}) \supseteq \ker(\mathrm{res})$ above, 
%
$$\mathrm{div}(\mathcal{L}) = \sum_{i=1}^r c_i[\widetilde{D}_i] = \mathrm{div}(\pi^*\cO(\sum_{i=1}^r c_i[D_i])).$$
Since $\mathrm{div}$ is injective, this implies $\mathcal{L} \simeq \pi^*\cO(\sum_{i=1}^r c_i[D_i])$. Thus, $\pi^*$ is surjective, as desired. This completes the proof.
\end{proof}

We say that a line bundle $\mathcal{L} \in \Pic(\widetilde{Y})$ is \emph{relatively ample} if it is ample with respect to the projective morphism $\rho: \widetilde{Y} \to \widetilde{X}$. Write $\Pic^a(\widetilde{Y}) \subset \Pic(\widetilde{Y})$ for the semigroup of relatively ample line bundles and $\fX(L)^{>0} \subset \fX(L)$ for the semigroup of (strictly) $P$-dominant weights.

\begin{prop}\label{prop:descriptionofpica}
The isomorphism $\pi^*: \Pic(G/P) \simeq \Pic(\widetilde{Y})$ of Proposition \ref{prop:descriptionofpic} restricts to a semigroup isomorphism $\Pic^a(G/P) \simeq \Pic^a(\widetilde{Y})$. If $G$ is semisimple and simply connected, there is a further isomorphism $\mathcal{L}_{G/P}:\fX(L)^{>0} \simeq \Pic^a(G/P)$. 
\end{prop}

\begin{proof}
First, observe that there is a finite morphism 
$$\widetilde{Y} \to G \times^P \widetilde{X} \simeq G/P \times \widetilde{X},$$
and $\rho: \widetilde{Y} \to \widetilde{X}$ coincides with the composition $\widetilde{Y} \to G/P \times \widetilde{X} \to \widetilde{X}$ (where the second map is the projection). If $\mathcal{L} \in \Pic(G/P)$ is ample, then $\pi^*\mathcal{L}$ coincides with the pullback to $\widetilde{Y}$ of the relatively ample line bundle $\mathcal{L} \otimes \cO_{\widetilde{X}} \in \Pic(G/P \times \widetilde{X})$. Since $\widetilde{Y} \to G/P \times \widetilde{X}$ is finite, this pullback is relatively ample. 

Conversely, suppose that $\pi^*\mathcal{L}$ is relatively ample. The zero section $G/P \subset \widetilde{Y}$ is the preimage under $\rho: \widetilde{Y} \to \widetilde{X}$ of $0 \in \widetilde{X}$. Hence, $\mathcal{L} \simeq (\pi^*\mathcal{L})|_{G/P}$ is ample. 

Now suppose $G$ is semisimple and simply connected. It is a standard fact that for a character $\chi \in \fX(L)$, the line bundle $\mathcal{L}_{G/P}(\chi)$ is ample if and only if $\chi$ is strictly dominant for $P$. Thus,  $\mathcal{L}_{G/P}:\fX(L) \xrightarrow{\sim} \Pic(G/P)$ restricts to an isomorphism
$\fX(L)^{>0} \xrightarrow{\sim} \Pic^a(G/P)$, as asserted.
\end{proof}

\section{$\QQ$-terminalizations of nilpotent covers}\label{subsec:terminalizationcover}

Fix the notation of Section \ref{subsec:picard}, i.e. $P=LN$, $\widetilde{\OO}_L$, $\widetilde{X}_L$, $\rho:\widetilde{Y} \to \widetilde{X}$, $\pi:\widetilde{Y} \to G/P$, and so on. For the remainder of this section, we will also assume:
$$\widetilde{\OO}_L \text{ is birationally rigid}.$$
Under this assumption, we will see below that the partial resolution $\rho: \widetilde{Y} \to \widetilde{X}$ is a $\QQ$-factorial terminalization. In the case when $\widetilde{\mathbb{O}}$ is a nilpotent orbit, this was proved by the first-named author (see \cite[Cor 4.6]{Losev4}). For universal covers in classical types, a similar result was obtained by Namikawa (see \cite{NamikawaQ}). The result below is due to the third-named author (see \cite[Cor 4.3]{Mitya2020}).

\begin{theorem}\label{thm:Qfactorialcover}
The partial resolution $\rho: \widetilde{Y} \to \widetilde{X}$ is a $\QQ$-factorial terminalization.
\end{theorem}

Theorem \ref{thm:Qfactorialcover} will allow us to provide purely Lie-theoretic descriptions of the Namikawa space $\fP^{\widetilde{X}}$ and Namikawa Weyl group $W^{\widetilde{X}}$ attached to $\widetilde{X}$ (see Section \ref{subsec:structurenamikawa} for definitions). These descriptions are as follows. Let $\underline{\pi}$ denote the restriction of $\pi: \widetilde{Y} \to G/P$ to $\widetilde{Y}^{\mathrm{reg}}$, and consider the pullback map on cohomology $\underline{\pi}^*: H^2(G/P, \CC) \to H^2(\widetilde{Y}^{\mathrm{reg}},\CC) = \fP^{\widetilde{X}}$. Note that  $H^2(G/P,\CC)$ is identified with $\fX(\fl \cap [\fg,\fg])$. Consider the composition
\begin{equation}\label{eq:defofeta}\eta: \fX(\fl \cap [\fg,\fg]) \simeq  H^2(G/P,\CC) \overset{\underline{\pi}^*}{\to} H^2(\widetilde{Y}^{\mathrm{reg}},\CC) = \fP^{\widetilde{X}},\end{equation}
We will see below that this map is an isomorphism. 

Next, we will describe $W^{\widetilde{X}}$ in terms of the data $(L,\widetilde{\OO}_L)$. Let $\Ad^*$ denote the co-adjoint action of $N_G(L)$ on $\fl^*$ and let $\mu_L: \widetilde{X}_L \to \fl^*$ denote the moment map. Consider the group
$$N_G(L, \widetilde{\OO}_L) := \{ (\eta, \zeta) \in N_G(L) \times \Aut(\widetilde{X}_L) \,|\, \Ad^*(\eta) \circ \mu_L = \mu_L \circ \zeta \}.$$
We can regard $L$ as a normal subgroup of $N_G(L,\widetilde{\OO}_L)$ via the natural embedding
$$L \hookrightarrow N_G(L, \widetilde{\OO}_L), \qquad l \mapsto (l, l).$$
Consider the quotient group
      $$ \widetilde{W}^{\widetilde{X}} := N_G(L, \widetilde{\OO}_L)/L. $$
Since both $N_G(L)/L$ and $\Aut(\widetilde{X}_L)$ are finite, so is $\widetilde{W}^{\widetilde{X}}$. Moreover, $\widetilde{W}^{\widetilde{X}}$ acts on $\fX(\fl \cap [\fg,\fg])$ via the natural projection $\widetilde{W}^{\widetilde{X}} \to N_G(L)/L$. 

There is a natural homomorphism
$$\kappa: \widetilde{W}^{\widetilde{X}} \to \Aut_G(\widetilde{X}).$$
For nilpotent orbits, $\kappa$ was defined in Step (ii) of the proof of \cite[Prop 4.7]{Losev4}. For nilpotent covers, the definition is analogous. We note that $\kappa$ is surjective by \cite[Proposition 13]{NamikawaQ2}.

\begin{prop}
\label{prop:namikawacovers}
The following are true:
\begin{itemize}
\item[(i)] The map $\eta: \fX(\fl \cap [\fg,\fg]) \to \fP^{\widetilde{X}}$ is an isomorphism.
\item[(ii)] $W^{\widetilde{X}}$ is identified with the normal subgroup
$$\ker{\kappa} \subset \widetilde{W}^{\widetilde{X}}$$
with its natural action on $\fX(\fl \cap [\fg,\fg])$.
\item[(iii)] The universal graded deformation $\widetilde{Y}_{\mathrm{univ}}$ (see Proposition \ref{prop:universaldef}) is identified with the Poisson scheme {$G\times^P(\fX(\fl \cap [\fg,\fg])\times \widetilde{X}_L\times \fp^\perp)$ over $ \fX(\fl \cap [\fg,\fg])$} for a uniquely determined isomorphism {$\fX(\fl \cap \fg,\fg])\xrightarrow{\sim} H^2(\widetilde{Y}^{\mathrm{reg}},\CC)$}. This isomorphism coincides with $\eta$.
\end{itemize}
\end{prop}
\begin{proof}
(i) is proved for nilpotent orbits in \cite[Prop 4.7]{Losev4}. The proof there can be easily generalized to arbitrary nilpotent covers. (ii) is \cite[Proposition 10]{NamikawaQ2}.

The isomorphism $G\times^P(\fX(\fl\cap [\fg,\fg])\times \widetilde{X}_L\times \fp^\perp)\xrightarrow{\sim} \widetilde{Y}_{\mathrm{univ}}$ is proved in the same way as \cite[Prop 4.7(2)]{Losev4}, see also \cite[Thm 9]{NamikawaQ2}. It remains to check that the isomorphism {$\fX(\fl \cap [\fg,\fg])\xrightarrow{\sim} \fP^{\widetilde{X}}$} in Proposition \ref{prop:universaldef} coincides with the map $\eta$ defined above. The variety $\widetilde{X}^{\mathrm{reg}}=G\times^P(\widetilde{X}^{\mathrm{reg}}\times \fp^\perp)$ is obtained as the Hamiltonian reduction of the product variety $T^*G\times \widetilde{X}^{\mathrm{reg}}$ under the action of $P$ and the deformation {$G\times^P(\fX(\fl \cap [\fg,\fg])\times \widetilde{X}_L\times \fp^\perp)$} is realized as the universal Hamiltonian reduction for this action, i.e. as $\mu_P^{-1}(\fX(\fp\cap [\fg,\fg]))/P$ (where $\mu_P$ is the moment map). It remains to show that $\eta$ coincides with the period map from Remark \ref{rmk:period_map}. This is a special case of \cite[Proposition 3.2.1]{Losev_isofquant}.
\end{proof}


Next, we explain the classification of $\QQ$-factorial terminalizations of $\widetilde{X}$. Recall the subset $\fP^{\mathrm{sing}} \subset \fP$ defined in Section \ref{amplecones}. Let $\fX(\fl \cap [\fg,\fg])_{\RR}^{\geq 0} \subset \fX(\fl \cap [\fg,\fg])_{\RR}$ denote the cone of $\fp$-dominant weights.

\begin{lemma}\label{lem:parab_sing_locus}
The following are true:
\begin{itemize}
\item[(i)] Under the identification $\eta: \fX(\fl \cap [\fg,\fg]) \simeq \fP$, the subset $\fP^{\mathrm{sing}}$ corresponds to the subset
$$\fX(\fl \cap [\fg,\fg])^{\mathrm{sing}} := \{\lambda \in \fX(\fl \cap [\fg,\fg]) \mid L \subsetneq G_{\lambda}\}.$$
\item[(ii)] Under the identification $\eta: \fX(\fl \cap [\fg,\fg])_{\RR} \simeq \fP_{\RR}$, the ample cone $\mathrm{Amp}(G \times^P (\widetilde{X}_L \times \fp^{\perp}))$ (defined by (\ref{eq:Amp_def})) corresponds to the subset $\fX(\fl \cap [\fg,\fg])_{\RR}^{\geq 0}$.
\item[(iii)] Every $\QQ$-terminalization of $\widetilde{X}$ is of the form $G \times^{P_1} (\widetilde{X}_{L}  \times \fp_1^{\perp})$, where $P_1$ is a parabolic subgroup of $G$ with Levi subgroup $L$. Two parabolic subgroups give rise to the same $\QQ$-factorial terminalization if and only if their dominant cones are $W^{\widetilde{X}}$-conjugate.
\end{itemize}
\end{lemma}

\begin{proof}
By definition $\fP^{\mathrm{sing}}$ consists of all parameters $\lambda \in \fP$ for which $\widetilde{Y}_\lambda$ is not affine. Viewing $\lambda$ as an element of $\fX(\mathfrak{l} \cap [\fg,\fg])$, we have $\widetilde{Y}_\lambda=G\times^P(\{\lambda\}\times \mathfrak{p}^\perp\times \widetilde{X}_L)$. If $G_\lambda=L$, this variety is affine. If $G_\lambda\neq L$, then $\widetilde{Y}_\lambda$ contains the projective variety $G_\lambda/(G_\lambda\cap P)$ and is therefore not affine. This completes the proof of (i). (ii) follows immediately from Proposition \ref{prop:descriptionofpica}.  Cones of the form $\fX(\fl \cap  [\fg,\fg])^{\geq 0}_{\RR}$ partition $\fP_{\RR}$. Now (iii) follows from Theorem \ref{thm:Q-term cones}.
\end{proof}

\begin{rmk}\label{Rem:Q-termin_classif}
In \cite[Thm 7]{NamikawaQ2}, Namikawa provides an alternative proof of Lemma \ref{lem:parab_sing_locus}(iii). His proof does not use Theorem \ref{thm:Q-term cones}; he deduces the first part of (iii) from the fact that ample cones of the form $\mathrm{Amp}(G \times^P (\widetilde{X}_L \times \fp^{\perp}))$ cover the fundamental chamber $\fP_{\RR}^{\geq 0}$. The second part of (iii) is standard.
\end{rmk}

Our definition of $\eta$ depends on the choice of $P$. The next result is essential for our computation of $\eta$.

\begin{prop}\label{prop:independence}
For two choices of $P$, the corresponding maps $\eta: \mathfrak{X}(\mathfrak{l} \cap [\fg,\fg]) \xrightarrow{\sim}\fP^{\widetilde{X}}$ coincide up to the action of $W^{\widetilde{X}}$ on the target. 
\end{prop}

\begin{proof}
Replacing $G$ with its commutator subgroup, we can assume that $G$ is semisimple. The proof has several steps.

{\it Step 1}.
Let $P_1,P_2$ be parabolics with Levi factor $L$ and write $\eta_1,\eta_2$ for the corresponding isomorphisms. We identify $\CC[G\times^P(\widetilde{X}_L\times \fp^\perp)]$ with $\CC[\widetilde{X}]$ via the pullback map under the terminalization morphism 
$G\times^P(\widetilde{X}_L\times \fp^\perp)\twoheadrightarrow \widetilde{X}$.

{\it Step 2}.
Choose a Zariski-generic element $\lambda\in \fX(\fl)$. Thanks to (iii) of Theorem \ref{thm:defsofsymplectic}, it suffices to show that $\CC[\widetilde{X}_{\eta_1(\lambda)}]$
and $\CC[\widetilde{X}_{\eta_2(\lambda)}]$ are isomorphic as filtered Poisson deformations of $\CC[\widetilde{X}]$. 

{\it Step 3}.
Let $\ell=\C\lambda$. Consider the variety $\widetilde{Y}_{\ell,i}:=G\times^{P_i}(\ell\times \widetilde{X}_L\times \fp_i^\perp)$ and the natural map $\rho_i: \widetilde{Y}_{\ell,i} \to \fg^*$. Let $\ell^\times:=\ell- \{0\}$ and $\widetilde{Y}^\times_{\ell,i}:=\ell^\times\times_{\ell}\widetilde{Y}_{\ell,i}$. There is a $\C^\times\times G$-equivariant identification  $\widetilde{Y}^\times_{\ell,i}\simeq G\times^L(\ell^\times\times \widetilde{X}_L)$ intertwining the maps to $\g^*$. 

{\it Step 4}.
The pullback $\widetilde{X}_{\eta_i(\ell)}$ of the universal deformation $\widetilde{X}_{\mathrm{univ}}$ (see Proposition 
\ref{prop:universaldefsymplectic}) to $\eta_i(\ell)$ coincides with $\Spec(\CC[\eta_i(\ell)\times_{\fP}\widetilde{Y}_{\mathrm{univ}}])$ as a Poisson scheme over $\ell$ with $\CC^\times$-action. So there is the Stein factorization $Y^\times_{\ell,i}\rightarrow \widetilde{X}_{\eta_i(\ell)}\rightarrow \fg^*$. This gives rise to a $\CC^\times$-equivariant isomorphism of Poisson $\ell$-schemes 
$\widetilde{X}_{\eta_1(\ell)}\xrightarrow{\sim} \widetilde{X}_{\eta_2(\ell)}$. By the construction in Step 1, this isomorphism gives the identity on $\CC[\widetilde{X}]$, i.e. is an isomorphism of filtered Poisson deformations. It follows that 
$\CC[\widetilde{X}_{\eta_1(\lambda)}]$
and $\CC[\widetilde{X}_{\eta_2(\lambda)}]$
are isomorphic as filtered Poisson deformations, completing the proof.
\end{proof}

By Proposition \ref{prop:partialdecomp}, $\fP^{\widetilde{X}}$ decomposes into partial Namikawa spaces
\begin{equation}\label{eq:partialdecomp2}\fP^{\widetilde{X}} \simeq \fP^{\widetilde{X}}_0 \oplus \bigoplus_{k=1}^t \fP^{\widetilde{X}}_k.\end{equation}
Here, $\fP_0^{\widetilde{X}} = H^2(\widetilde{X}^{\mathrm{reg}},\CC)$ and $\fP_k^{\widetilde{X}}$ is the space of $\pi_1(\fL_k)$-invariants in the dual Cartan subalgebra $\fh_k^*$ associated to the Kleinian singularity $\Sigma_k$.

\begin{rmk}\label{rmk:noE8}
We note that Proposition \ref{prop:namikawacovers} together with (\ref{eq:partialdecomp2}) provides a more conceptual proof of Lemma \ref{lem:noE8}, at least in exceptional types. Let $\OO$ be a nilpotent orbit in a simple exceptional Lie algebra and suppose $X=\Spec(\CC[\OO])$ contains a codimension 2 leaf $\fL_k \subset X$ such that the corresponding singularity is of type $E_8$. Note that $E_8$ has no nontrivial diagram automorphisms. Hence, $\dim \fP^X_k = \dim \fh_k^* = 8$. Choose a Levi subgroup $L \subset G$ and a birationally rigid $L$-orbit $\OO_L$ such that $\OO=\mathrm{Bind}^G_L \OO_L$. Then by Proposition \ref{prop:namikawacovers} and the decomposition (\ref{eq:partialdecomp2}), we have $\dim(\fX(\fl)) = \dim \fP^X \geq \dim \fP^X_k=8$. Thus, $G=E_8$ and $(L,\OO_L)=(T,\{0\})$. That is, $\OO$ is the principal nilpotent orbit in $E_8$.
\end{rmk}

For our computation of $\eta$, we will need Lie-theoretic descriptions of the spaces $\fP_k^{\widetilde{X}}$ analogous to Proposition \ref{prop:namikawacovers}. First, we will describe the space $\fP_0^{\widetilde{X}}$. We will describe the spaces $\fP_k^{\widetilde{X}}$, $k \geq 1$, in Section \ref{subsec:descriptionpartial}, under some additional restrictions. 
For the following lemma, suppose $G$ is simply connected. Let $R_x$ denote the reductive part of the stabilizer of $x \in \widetilde{\OO}$ and let $\mathfrak{r}$ be its Lie algebra. Note that the adjoint action of $R_x$ on $\fX(\mathfrak{r})$ factors through $R_x/R_x^{\circ} \simeq \pi^G_1(\widetilde{\OO})$. 

\begin{lemma}\label{lem:computeH2} Suppose that $G$ is semisimple and simply connected.
The following are true:
\begin{itemize}
\item[(i)] There are natural identifications $\Pic(\widetilde{X}^{\mathrm{reg}})\xrightarrow{\sim}\Pic(\widetilde{\OO})$ and $\Pic(\widetilde{\OO})\xrightarrow{\sim}\fX(R_x)$.
\item[(ii)] There are isomorphisms 
$$\Pic(\widetilde{\OO})\otimes_{\ZZ}\CC\xrightarrow{\sim} H^2(\widetilde{\OO},\CC)
\xrightarrow{\sim} \fX(R_x)\otimes_{\ZZ}\CC=\fX(\mathfrak{r})^{\pi_1(\widetilde{\OO})},$$
where the first map is $c_1$.
\item[(iii)] Restriction along the embedding $\widetilde{\OO}\hookrightarrow 
\widetilde{X}^{\mathrm{reg}}$ induces an isomorphism $H^2(\widetilde{X}^{\mathrm{reg}},\CC)\xrightarrow{\sim} H^2(\widetilde{\OO},\CC)$. 
\end{itemize}
\end{lemma}

\begin{proof}
(i): the first isomorphism follows from the inequality $\operatorname{codim}
(\widetilde{X}^{\mathrm{reg}}- \widetilde{\OO},\widetilde{X}^{\mathrm{reg}})\geqslant 2$ and the fact that $\widetilde{X}^{\mathrm{reg}}$ is smooth. The second isomorphism is standard, see e.g. 
\cite[Proposition 3.2]{Knop1989}.

(ii): The isomorphism $H^2(\widetilde{\OO},\CC)
\xrightarrow{\sim} \fX(\mathfrak{r})^{\pi_1(\widetilde{\OO})}$ is also standard, see e.g.
\cite[Theorem 3.3]{BISWAS}. Under the identifications $H^2(\widetilde{\OO},\CC)
\xrightarrow{\sim} \fX(\mathfrak{r})^{\pi_1(\widetilde{\OO})}$ and $\Pic(\widetilde{\OO})\xrightarrow{\sim}\fX(R_x)$, the Chern class map $c_1:\Pic(\widetilde{\OO})\rightarrow 
H^2(\widetilde{\OO}, \CC)$ coincides with the natural map $\fX(R_x)\rightarrow \fX(R_x)\otimes_{\ZZ}\CC=\fX(\mathfrak{r})^{\pi_1(\widetilde{\OO})}$. (ii) follows at once.

\{(iii): For $i \geq 0$, let $\widetilde{X}^i$ denote the union of $\widetilde{\OO} \subset \widetilde{X}^{\mathrm{reg}}$ and all $G$-orbits in $\widetilde{X}$ of codimension $\leq 2i$. Since $G$ acts on $\widetilde{X}$ with finitely many orbits, $\widetilde{X} = \widetilde{X}^N$ for some $N \geq 1$. Thus, it suffices to show that $H^2(\widetilde{X}^{i+1},\CC) \to H^2(\widetilde{X}^i, \CC)$ is an isomorphism for all $i$. Let $Y_1,...,Y_m  \subset \widetilde{X}^{i+1}$ denote the $G$-orbits of codimension $2(i+1)$, so that $\widetilde{X}^{i+1}= \widetilde{X}^i \cup \bigcup_{k=1}^m Y_k$. For each $k$, choose a tubular neighborhood $U_k \subset \widetilde{X}^{i+1}$ of $Y_k$ and set $U_k^{\times} := U_k \cap \widetilde{X}^i$. There is a Mayer-Vietoris sequence
\begin{align*}
        \ldots \to H^1(\widetilde{X}^i, \CC)\oplus \bigoplus H^1(U_k, \CC) &\to \bigoplus H^1(U_k^\times, \CC)\to H^2(\widetilde{X}^{i+1}, \CC) \\
        \to H^2(\widetilde{X}^i, \CC)\oplus \bigoplus H^2(U_k, \CC) &\to \bigoplus H^2(U_k^\times, \CC)\to \ldots
\end{align*}
Note that $U_k$ and $U_k^\times$ are fibrations over $Y_k$ with fibers $D_{2i}$ and $D_{2i}^\times$ respectively, where $D_{2i}$ is the $2(i+1)$-dimensional disc and $D_{2i}^{\times}$ is its puncture. Since 
$$H^1(D_{2i}, \CC)=H^2(D_{2i}, \CC)=H^1(D_{2i}^\times, \CC)=H^2(D_{2i}^\times, \CC)=0,$$
the maps 
$$H^j(Y_k, \CC)\to H^j(U_k,\CC), \quad H^j(Y_k, \CC)\to H^j(U_k^\times,\CC), \qquad j=1,2$$
are isomorphisms. It follows that the restriction maps $H^j(U_i,\CC)\to H^j(U_i^\times,\CC)$ are isomorphisms, and therefore $H^2(\widetilde{X}^{i+1},\CC)\xrightarrow{\sim} H^2(\widetilde{X}^i, \CC)$, as asserted.
\end{proof}

\section{Parabolic induction of Hamiltonian quantizations}\label{subsec:inductionquantizations}

In this section, we will define the notion of parabolic induction for Hamiltonian quantizations of nilpotent covers. 

Choose a parabolic subgroup $P \subset G$ with Levi decomposition $P=LN$. Let $\widetilde{\mathbb{O}}_L$ be an $L$-equivariant nilpotent cover and let $\widetilde{\OO} = \mathrm{Bind}^G_L\widetilde{\OO}_L$. Consider the partial resolution $\rho:\widetilde{Y} := G \times^P (\widetilde{X}_L \times \fp^{\perp}) \to \widetilde{X}$ defined in Section \ref{subsec:picard}.
    
Starting with a Hamiltonian quantization $\cA_L$ of $\CC[\widetilde{\mathbb{O}}_L]$, we will produce a Hamiltonian quantization $\Ind^G_L \cA_L$ of $\CC[\widetilde{\mathbb{O}}]$, thus defining a map
$$\Ind^G_L: \mathrm{Quant}^L(\CC[\widetilde{\mathbb{O}}_L]) \to \mathrm{Quant}^G(\CC[\widetilde{\mathbb{O}}])$$
called \emph{parabolic induction}.\index{induction!of Hamiltonian quantizations} 

Our construction requires a bit of additional notation. 

\begin{itemize}
    \item Let $\mu_L: \widetilde{X}_L \to \overline{\OO}_L \subset \fl^*$ denote the (classical) moment map for the $L$-action on $\widetilde{X}_L$.
    \item Let $\Phi_L: \fl \to \cA_L$ denote the (quantum) co-moment map for the $L$-action on $\cA_L$. 
    \item Let $\mathcal{D}_L$ denote the microlocalization of $\cA_L$ over $\widetilde{X}_L$. 
    \item Let $\mathfrak{D}_{G/N}$ denote the sheaf of differential operators on $G/N$ (regarded as a sheaf in the conical topology on $T^*(G/N)$). Since $L$ normalizes $N$, there is a right $L$-action  on $G/N$ and hence on $T^*(G/N)$ and $\mathfrak{D}_{G/N}$.
    \item Let
    $$\mu_{G/N}:T^*(G/N) \simeq G \times^N (\mathfrak{g}/\mathfrak{n})^* \to \mathfrak{l}^*$$
    denote the (classical) moment map for the  $L$-action on $T^*(G/N)$ (from the right): the map $\mu^*_{G/N}$ sends any element of $\mathfrak{l}$ to the corresponding vector field viewed as a function on $T^*(G/N)$. {Explicitly, a point in $T^*(G/N)$ is $N.(g,x)$ for $x\in (\fg/\mathfrak{n})^*\subset \fg^*$, the map $\mu_{G/N}$ sends this point to $-x|_{\fl}$.}
    \item Let $\Phi^0_{G/N}: \fl \to \Gamma(T^*(G/N),\mathfrak{D}_{G/N})$  denote the usual (quantum) co-moment map for the $L$-action on $\mathfrak{D}_{G/N}$, i.e. $\Phi^0_{G/N}$ takes $\xi \in \fl$ to the action vector field on $G/N$ determined by $\xi$. 
    \item Let $\Phi_{G/N} := \Phi_{G/N}^0 + \rho_{\mathfrak{n}}$, where $\rho_{\fn}$ denotes one-dimensional representation of $\fl$ defined in (\ref{eqn:defofrho}). 
\end{itemize}

Consider the product variety $T^*(G/N) \times \widetilde{X}_L$ with the diagonal left $L$-action
$$l(d,x) := (dl^{-1},lx), \qquad l \in L, \ d \in T^*(G/N), \ x \in \widetilde{X}_L.$$
This action has a moment map
$$\mu := \mu_{G/N} + \mu_{\widetilde{X}_L}: T^*(G/N) \times \widetilde{X}_L \to \fl^*.$$ 
Note that
$$\mu^{-1}(0) \simeq G \times^N (\widetilde{X}_L \times \fp^{\perp}) \subset G \times^N ( \widetilde{X}_L \times \fn^{\perp}) \simeq T^*(G/N) \times \widetilde{X}_L.$$
and therefore
$$\mu^{-1}(0)/L \simeq G \times^P (\widetilde{X}_L \times \fp^{\perp}) = \widetilde{Y}.$$
Form the completed tensor product
$$\mathcal{D}' := \mathfrak{D}_{G/N} \  \widehat{\otimes} \ \mathcal{D}_L.$$
Here, the completion is defined with respect to the filtrations on the factors. Note that $\mathcal{D}'$ is an $L$-equivariant filtered quantization of $T^*(G/N) \times \widetilde{X}_L$ with a naturally defined quantum co-moment map 
$$\Phi := \Phi_{G/N} \otimes 1 + 1 \otimes \Phi_L: \fl \to \Gamma(T^*(G/N) \times \widetilde{X}_L, \mathcal{D}').$$ 
Let $\mathrm{Ind}^G_L \mathcal{D}_L$ denote the quantum Hamiltonian reduction of $\mathcal{D}'$ at $0\in \fl^*$. More precisely, form the left ideal
$$I := \mathcal{D}' \ \mathrm{Span}\{\Phi(\xi) \mid \xi \in \fl\} \subset \mathcal{D}'.$$
The quotient $\mathcal{D}'/I$ is an $L$-equivariant sheaf of $\mathcal{D}'$-modules on $T^*(G/N) \times \widetilde{X}_L$, set-theoretically supported on $\mu^{-1}(0)$. Consider the quotient morphism $q: \mu^{-1}(0) \to \widetilde{Y}$ and set
$$\Ind^G_L \mathcal{D}_L := \left(q_* [\mathcal{D}'/I] \right)^L.$$
Since the $L$-action on $\mu^{-1}(0)$ is free, $\Ind^G_L\mathcal{D}_L$ is a filtered quantization of $\widetilde{Y}$, see e.g. \cite[Lem 3.3.1]{Losev_isofquant}. It has a left $G$-action (coming from the left $G$-action on $\mathfrak{D}_{G/N}$) and a quantum co-moment map (coming from the quantum co-moment map $\Phi_{G,G/N}\otimes 1:\fg \to \Gamma(\mathfrak{D}_{G/N}\widehat{\otimes}\mathcal{D}_L)$, where $\Phi_{G,G/N}:\fg\rightarrow \Gamma(\mathfrak{D}_{G/N})$ is the quantum comoment map for the $G$-action, taking $\xi\in \fg$ to the corresponding vector field on $G/N$).

Finally, define
$$\Ind^G_L \cA_L := \Gamma(\widetilde{Y}, \Ind^G_L \mathcal{D}_L).$$ 

We now show that $\Ind^G_L \cA_L$ is a quantization of $\CC[\widetilde{\OO}]=\CC[\widetilde{X}]$.  The following lemma is standard.

\begin{lemma}\label{lem:quant_pushforward}
Let $X_1,X_2$ be normal graded Poisson varieties and let $\rho:X_1\rightarrow X_2$ be a $\C^\times$-equivariant Poisson morphism such that $\mathcal{O}_{X_2}\xrightarrow{\sim }\rho_*\mathcal{O}_{X_1}$ and $R^1\rho_*\mathcal{O}_{X_1}=0$. Then for any filtered quantization $\mathcal{D}_1$ of $X_1$, $\rho_*\mathcal{D}_1$ is a filtered quantization of $X_2$.
\end{lemma}

We provide a proof for the reader's convenience. 
\begin{proof}
The proof reduces to the case when $X_2$ is affine. The algebra $\Gamma(X_1,\mathcal{D}_1)$ comes with a complete and separated filtration induced from the filtration on $\mathcal{D}_1$. We wish to show that the embedding $\gr\Gamma(\mathcal{D}_1)\hookrightarrow \CC[X_1](=\CC[X_2])$ is an isomorphism or, equivalently, that every homogeneous element lies in the image. Cover $X_1$ by $\CC^\times$-stable open affine subvarieties $U_1,\ldots,U_k$. For a homogeneous element $f\in \CC[X_2]$ of degree $d$ we can find elements $\widetilde{f}_i\in \Gamma(U_i,\mathcal{D}_{1,\leqslant d})$ whose principal symbols satisfy $\gr \widetilde{f}_i=\rho^*(f)|_{U_i}$. The elements 
$\widetilde{f}_i-\widetilde{f}_j$ lie in $\Gamma(U_i\cap U_j, \mathcal{D}_{1,\leqslant d-1})$. Let $f_{ij}\in \C[U_i\cap U_j]$ denote the degree $d-1$ component of $\widetilde{f}_i-\widetilde{f}_j$. It is a 1-cocycle. Since $H^1(X_1,\mathcal{O}_{X_1})=0$, it is a coboundary: we can find degree $d-1$ elements $g_i\in \C[U_i]$ with $f_{ij}=g_i-g_j$. Lift $g_i$ to an  element $\widetilde{g}_i\in \Gamma(U_i,\mathcal{D}_{1,\leqslant d-1})$. Replacing $\widetilde{f}_i$ with $\widetilde{f}_i-\widetilde{g}_i$ if necessary we can arrange so that the elements $\widetilde{f}_i-\widetilde{f}_j$ live in degree $\leqslant d-2$. We can iterate this procedure, each time reducing the degree by $1$. This process converges because the filtration on $\mathcal{D}_1$ is complete and separated. Thus we can assume that $\widetilde{f}_i-\widetilde{f}_j=0$, i.e. that the elements $\widetilde{f}_i$ glue to $\widetilde{f}\in \Gamma(\mathcal{D}_1)$. The principal symbol of this element is $f$ and we are done.
\end{proof}
By Lemma \ref{lem:rationalsingularities}, we have $\cO_{\widetilde{X}} \xrightarrow{\sim} \rho_*\cO_{\widetilde{Y}}$ and $R^i\rho_*\cO_{\widetilde{Y}}=0$ for $i>0$. Hence by Lemma \ref{lem:quant_pushforward}, the filtered algebra $\Ind^G_L \cA_L$ is a quantization of $\CC[\widetilde{\OO}]$. Since $\Ind^G_L \mathcal{D}_L$ has a Hamiltonian $G$-action, so does $\Ind^G_L\cA_L$. Thus, we have defined the desired map $\Ind^G_L: \mathrm{Quant}^L(\CC[\widetilde{\OO}_L]) \to \mathrm{Quant}^G(\CC[\widetilde{\OO}])$. In fact, below we will see that $\Ind^G_L\cA_L$ does not depend on the choice of $P$ justifying the notation.

Next, we will show that parabolic induction of Hamiltonian quantizations is transitive. Choose a parabolic subgroup $Q = MU \subset G$ such that $P \subset Q$ and $L \subset M$. Let $\widetilde{\mathbb{O}}_M=\mathrm{Bind}^M_L \widetilde{\OO}_L$ and $\widetilde{X}_M = \Spec(\CC[\widetilde{\OO}_M])$.

\begin{lemma}\label{lem:induct_transitivity}
Let $\cA_L$ be a Hamiltonian quantization of $\CC[\widetilde{\OO}_L]$. Then there is an isomorphism of Hamiltonian quantizations
$$\Ind^G_L \cA_L\simeq \Ind^G_M(\Ind^M_{L}\cA_L).$$
\end{lemma}

\begin{proof}
It is clear from the construction that the central part of the co-moment map is preserved under induction. Thus, we can reduce to the case when $\fg$ is semisimple. In this case, the Hamiltonian structure is determined uniquely, see Lemma \ref{lem:co-momentunique}. So it suffices to exhibit a $G$-invariant isomorphism of filtered quantizations.

Consider the partial resolution $\overline{\rho}: \widetilde{Z} := G \times^Q (\widetilde{X}_M \times \fq^{\perp}) \to \widetilde{X}$. Note that $\rho: \widetilde{Y} \to \widetilde{X}$ factors through a proper birational map $\rho': \widetilde{Y} \to \widetilde{Z}$.  By the transitivity of pushforwards, it suffices to show that
\begin{equation}\label{eq:transit_induct_iso}
\rho'_*(\Ind^G_L \mathcal{D}_L)\simeq
\Ind^G_M \mathcal{D}_M.
\end{equation}
as filtered quantizations. Let $\rho_M: \widetilde{Y}_M =M\times^{M\cap P}(\widetilde{X}_L\times \fp^{\perp}) \to \widetilde{X}_M$ denote the partial resolution. Consider the variety $T^*(G/U)\times \widetilde{Y}_M$ and the filtered quantization
\begin{equation}\label{eq:aux_quantization}
\mathfrak{D}_{G/U} \ \widehat{\otimes} \ \operatorname{Ind}_L^M \mathcal{D}_L.
\end{equation}
The quantization $\mathfrak{D}_{G/U}\widehat{\otimes}\mathcal{D}_M$ appearing in the construction of $\Ind^G_M \mathcal{D}_M$ is
the pushforward of (\ref{eq:aux_quantization})
under $\operatorname{id}_{T^*(G/U)}\times 
\rho_M$. On the other hand, (\ref{eq:aux_quantization}) is the reduction of 
\begin{equation}\label{eq:aux_quantization1} 
\mathfrak{D}_{G/U} \ \widehat{\otimes} \ 
\mathfrak{D}_{M/M\cap N} \ \widehat{\otimes} \ \mathcal{D}_L
\end{equation}
under the diagonal action of $L$ on the second and third factors.

Note that the $M$-action on (\ref{eq:aux_quantization1}) is trivial on the third factor. Thus, if we reduce (\ref{eq:aux_quantization1}) under the $M$-action, we recover the filtered quantization
$\mathfrak{D}_{G/N} \ 
\widehat{\otimes} \ \mathcal{D}_L$. 
The quantum co-moment map for the $L$-action on $\mathfrak{D}_{G/N} \ 
\widehat{\otimes} \ \mathcal{D}_L$ is induced from the quantum co-moment map on (\ref{eq:aux_quantization1}), and is therefore shifted by the character
$\rho_{\mathfrak{n}}=\rho_{\mathfrak{u}}|_{\mathfrak{l}}+\rho_{\mathfrak{u} \cap \fm}$. So if we reduce 
(\ref{eq:aux_quantization1}) under the action of $M\times L$ we recover the quantization $\Ind^G_L \mathcal{D}_L$.

Thus, the left hand side of (\ref{eq:transit_induct_iso}) is obtained from (\ref{eq:aux_quantization}) by first pushing forward along the map $\operatorname{id}_{T^*(G/U)}
\times \rho_M$ and then reducing under the $M$-action. The right hand side of (\ref{eq:transit_induct_iso}) is obtained from 
(\ref{eq:aux_quantization}) by first reducing under the $M$-action and then pushing forward along 
$\rho'$. Note that $\rho'$ is the morphism induced by $\widehat{\rho}:=\operatorname{id}_{T^*(G/U)} \times \rho_M$ under classical Hamiltonian reduction. The group $M$ acts freely on both $T^*(G/U)\times \widetilde{X}_M$ and $T^*(G/U)\times \widetilde{Y}_M$. Note that the morphisms $\widehat{\rho}$ and $\rho'$ satisfy the conditions of Lemma \ref{lem:quant_pushforward}. Write $\bullet/\!/\!/ M$ for the quantum Hamiltonian reduction of a Hamiltonian quantization by a free action of $M$. To prove (\ref{eq:transit_induct_iso}), it suffices to show for an arbitrary Hamiltonian quantization $\mathcal{D}$ of $T^*(G/U)\times \widetilde{Y}_M$, there is an isomorphism of quantizations 
$(\widehat{\rho}_* \mathcal{D})/\!/\!/ M\xrightarrow{\sim}\rho'_*(\mathcal{D}/\!/\!/M)$. We note that Hamiltonian reduction sends quantizations to quantizations because the action of $M$ is free, while $\rho'_*,\widehat{\rho}_*$ send quantizations to quantizations thanks to Lemma \ref{lem:quant_pushforward}. Note that there is a natural homomorphism 
$(\widehat{\rho}_* \mathcal{D})/\!/\!/ M\rightarrow\rho'_*(\mathcal{D}/\!/\!/M)$ with associated graded equal to the identity map on $\mathcal{O}_{\widetilde{Z}}$. Since both sides are filtered quantizations, this homomorphism is an isomorphism.

%
%
%
%
%

\end{proof}

Now suppose $\widetilde{\OO}_L$ is birationally rigid. By Theorem \ref{thm:quantsofsymplectic}, quantizations of $\CC[\widetilde{\OO}]$ are parameterized by $W^{\widetilde{X}}$-orbits on $\fP^{\widetilde{X}}$. By Proposition \ref{prop:namikawacovers}, the latter space is identified with $\fX(\fl \cap [\fg,\fg])$. Thus, we get a bijection
\begin{equation}\label{eqn:Quantcenter}
\fX(\fl \cap [\fg,\fg])/W^{\widetilde{X}} \xrightarrow{\sim} \mathrm{Quant}(\CC[\widetilde{\mathbb{O}}]), \qquad W^{\widetilde{X}} \cdot \beta \mapsto \cA_{\eta(\beta)}^{\widetilde{X}}.\end{equation}
Abusing notation slightly, we will often abbreviate $\cA_{\eta(\beta)}^{\widetilde{X}}$ by writing $\cA_{\beta}^{\widetilde{X}}$. The isomorphism $\eta: \fX(\fl \cap [\fg,\fg]) \simeq \fP^{\widetilde{X}}$ of Proposition \ref{prop:namikawacovers} extends to an isomorphism (still denoted by $\eta$)
$$\eta: \fX(\fl) \xrightarrow{\sim} \fX(\fl \cap [\fg,\fg]) \oplus \fz(\fg)^* \xrightarrow{\sim} \fP^{\widetilde{X}} \oplus \fz(\fg)^* =: \overline{\fP}^{\widetilde{X}}$$
Hence by Proposition \ref{prop:classificationHamiltonian}, there is a bijection
\begin{equation}\label{eq:reductive1}
\fX(\fl)/W^{\widetilde{X}} \xrightarrow{\sim} \mathrm{Quant}^G(\CC[\widetilde{\OO}]), \qquad W^{\widetilde{X}} \cdot \beta \mapsto (\cA^{\widetilde{X}}_{\beta},\Phi^{\widetilde{X}}_{\beta}) \end{equation}
Here, $\Phi_{\beta}^{\widetilde{X}}: U(\fg) \to \cA_{\beta}^{\widetilde{X}}$ is the (unique) quantum co-moment map satisfying
$$\Phi_{\beta}^{\widetilde{X}}|_{\fz(\fg)} = \beta|_{\fz(\fg)},$$
see Section \ref{subsec:equivariant}. The next proposition shows that the quantization parameter $\beta$ is preserved under $\Ind^G_M$.

\begin{prop}\label{prop:quantizationparaminduction}
For each $\beta \in \fX(\fl)$, there is an isomorphism of Hamiltonian quantizations 
\begin{equation}\label{eq:induct_param}
\Ind^G_M \cA_{\beta}^{\widetilde{X}_M} \simeq \cA_{\beta}^{\widetilde{X}}.
\end{equation}
\end{prop}

\begin{proof}
As in the proof of Lemma \ref{lem:induct_transitivity}, we can assume that $\fg$ is semisimple. Thanks to Lemma \ref{lem:induct_transitivity}, we can also assume that $M=L$. Consider the $\QQ$-factorial terminalization $\widetilde{Y}:=G \times^P (\widetilde{X}_L \times \fp^{\perp}) \to \widetilde{X}$, and its universal graded Poisson deformation
$$\widetilde{Y}_{\mathrm{univ}} = G\times^P(\fX(\fl)\times \widetilde{X}_L\times \fp^\perp),$$ 
see Proposition \ref{prop:namikawacovers}(iii). Let $\mathcal{D}^{\widetilde{Y}_{\mathrm{univ}}}$ (resp. $\mathcal{D}^{\widetilde{X}_{L}}$) denote the canonical quantization of $\widetilde{Y}_{\mathrm{univ}}$ (resp. $\widetilde{X}_{L}$), see Proposition \ref{prop:universal for Q-term}. Consider the filtered quantization $\overline{\mathcal{D}}^{\widetilde{Y},\mathrm{univ}}$ of $\widetilde{Y}_{\mathrm{univ}}$ defined by
$$\overline{\mathcal{D}}^{\widetilde{Y},\mathrm{univ}} := \left([\mathfrak{D}_{G/N} \ \widehat{\otimes} \ \mathcal{D}^{{\widetilde{X}_{L}}}]/ 
[\mathfrak{D}_{G/N} \ \widehat{\otimes} \  \mathcal{D}^{\widetilde{X}_{L}}]\Phi([\mathfrak{l},\mathfrak{l}])\right)^L,$$
By \cite[Theorem 5.4.1]{Losev_isofquant}, $\overline{\mathcal{D}}^{\widetilde{Y},\mathrm{univ}}$ is an even quantization (more precisely, its restriction to $\widetilde{Y}^{\mathrm{reg}}_{\mathrm{univ}}$ is the specialization of a graded even quantization). Thus, by the results of \cite[Section 2.3]{Losev_isofquant}, there is an isomorphism of filtered quantizations
$$\overline{\mathcal{D}}^{\widetilde{Y},\mathrm{univ}} \simeq \mathcal{D}^{\widetilde{Y},\mathrm{univ}}$$
Specializing both sides to $\beta \in \fX(\fl)$, we get an isomorphism of filtered quantizations
$$\Ind^G_L \mathcal{D}_{\beta}^{\widetilde{X}_L} \simeq \mathcal{D}_{\beta}^{\widetilde{Y}}.$$
Now, the result follows by taking global sections.

\end{proof}

Using the classification of filtered quantizations of $\CC[\widetilde{X}]$, see 
Theorem \ref{thm:quantsofsymplectic}, and 
Proposition \ref{prop:quantizationparaminduction} we obtain the following.

\begin{cor}
Let $\cA_M$ be a Hamiltonian quantization of $\widetilde{X}_M$. Then, up to an isomorphism of Hamiltonian quantizations, $\Ind^G_M \cA_M$ depends only on $M$ (and not on $Q$).
\end{cor}

\section{Namikawa spaces vs parabolic induction}\label{subsec:inductioncompatibility}
    In this section we begin the task of computing the isomorphism $\eta: \fX(\fl \cap [\fg,\fg]) \xrightarrow{\sim} \fP^{\widetilde{X}}$. Here we establish some basic facts about the behavior of Namikawa spaces under parabolic induction.

    Choose a parabolic subgroup $P=LN \subset G$ and a birationally rigid $L$-equivariant nilpotent cover $\widetilde{\OO}_L$ such that $\widetilde{\OO}=\mathrm{Bind}^G_L \widetilde{\OO}_L$. Choose also a parabolic subgroup $Q =MU \subset G$ such that $P\subset Q$ and $L \subset M$. Let $\widetilde{\OO}_M = \mathrm{Bind}^M_L \widetilde{\OO}_L$ and $P_M = P \cap M$, a parabolic in $M$. Consider the $\QQ$-terminalizations
    $$\rho: \widetilde{Y} = G \times^P (\widetilde{X}_L \times \fp^{\perp}) \to \widetilde{X}, \qquad \rho_M: \widetilde{Y}_M = M \times^{P_M} (\widetilde{X}_L \times \fp_M^{\perp}) \to \widetilde{X}_M$$
    Note that $\widetilde{Y}_M$ is an  $M$-equivariant fiber bundle over $M/P_M=Q/P$. Hence $G\times^Q(\widetilde{Y}_M\times \fq^\perp)$ is a $G$-equivariant fiber bundle over $G/P$. Also $G\times^Q(\widetilde{Y}_M\times \fq^\perp)$ has a $G$-equivariant projective birational morphism to $G\times^Q(\widetilde{X}_M\times \fq^\perp)$ and hence to $\widetilde{X}$. Note that both varieties $\widetilde{Y}, G\times^Q(\widetilde{Y}_M\times \fq^\perp)$
come with natural $G\times\CC^\times$-actions
(on $G\times^Q(\widetilde{X}_M\times \fq^\perp)$ the factor $\CC^\times$ acts fiberwise via $t.(y,z)=(ty, t^{d}z)$ for $t\in \CC^\times, y\in \widetilde{Y}_M, z\in \fq^\perp$, where $d$ is the weight of the Poisson bracket on $\widetilde{Y}_M$), and the morphisms to $\widetilde{X}$ are $G\times\C^\times$-equivariant. 
    \begin{lemma}\label{lem:compatibility1}
    There is a $G\times \C^\times$-equivariant isomorphism
    $$\widetilde{Y} \xrightarrow{\sim} G \times^Q (\widetilde{Y}_M \times \fq^{\perp})$$
    intertwining the morphisms to $\widetilde{X}$ and to $G/P$.
    \end{lemma}
    \begin{proof}
    The variety $G\times^Q(\widetilde{Y}_M\times \fq^\perp)$ is $\QQ$-factorial terminal and so is a $\QQ$-factorial terminalization of $\widetilde{X}$. The fiber of $0\in \widetilde{X}$ is $G/P$ so restricting a relatively ample line bundle from $G\times^Q(\widetilde{Y}_M\times \fq^\perp)$ to $G/P$ gives an ample line bundle. Thus, by Proposition \ref{prop:descriptionofpica}, the ample cone of $G\times^Q(\widetilde{Y}_M\times \fq^\perp)$
    is contained in that of $\widetilde{Y}$.
    By Lemma \ref{lem:parab_sing_locus}, there is an isomorphism of $\QQ$-factorial terminalizations of $\widetilde{X}$ \begin{equation}\label{eq:terminalizations_iso}\widetilde{Y}\xrightarrow{\sim} G\times^Q(\widetilde{Y}_M\times \fq^\perp).\end{equation}
    Such an isomorphism is necessarily unique because the terminalization morphisms are birational. In particular, it is $G\times \C^\times$-equivariant. The morphisms $\widetilde{Y}, G\times^Q(\widetilde{Y}_M\times \fq^\perp)
    \rightarrow G/P$ are recovered as the categorical quotient morphisms for the action of $\C^\times$. Hence, they are intertwined by (\ref{eq:terminalizations_iso}) as well. 
    \end{proof}

    For the next result, we will assume that $G$ is semisimple in order to simplify the statements. Let ${\pi}'$ denote the composition of $\pi:\widetilde{Y}^{\mathrm{reg}}\rightarrow G/P$ with the projection $G/P\rightarrow G/Q$.
    Consider the composition
    \begin{equation}\label{eq:defofbareta}
    {\eta}': \fX(\fm) \simeq  H^2(G/Q,\CC) \overset{{\pi}^{'*}}{\to} H^2(\widetilde{Y}^{\mathrm{reg}},\CC) = \fP^{\widetilde{X}}.
    \end{equation}
    Note that ${\eta}'$ coincides with the composition 
    \begin{equation}\label{eq:defofbareta1}\fX(\fm) \hookrightarrow \fX(\fl) \xrightarrow{\eta} \fP^{\widetilde{X}}.\end{equation}
    Thanks to Lemma \ref{lem:compatibility1}, we can view $\widetilde{Y}^{\mathrm{reg}}$ as a fiber bundle over $G/Q$ with fiber $\widetilde{Y}_M^{\mathrm{reg}}\times \fq^\perp$. In particular, there are pullback maps
    \begin{equation}\label{eq:pullback_maps}H^2(G/Q,\CC)\rightarrow H^2(\widetilde{Y}^{\mathrm{reg}},\CC),\qquad H^2(\widetilde{Y}^{\mathrm{reg}},\CC) \rightarrow 
    H^2(\widetilde{Y}_M^{\mathrm{reg}},\CC).
    \end{equation} 
    We wish to describe (\ref{eq:pullback_maps}) under the identifications
    \begin{equation}\label{eq:bunch_of_identifications}
    \fX(\fm)\xrightarrow{\sim} H^2(G/Q,\CC), \quad 
    \eta: \fX(\fl)\xrightarrow{\sim} H^2(\widetilde{Y}^{\mathrm{reg}},\CC), \quad
    \eta_M:\fX(\fl\cap [\fm,\fm])\xrightarrow{\sim} H^2(\widetilde{Y}_M^{\mathrm{reg}},\CC).
    \end{equation}

    
    %
    %
    
    \begin{prop}\label{prop:compatibility2}
    Under the identifications in (\ref{eq:bunch_of_identifications}), the pullback maps in (\ref{eq:pullback_maps}) correspond to the inclusion $\fX(\fm)\hookrightarrow \fX(\fl)$ and the projection $\fX(\fl)\twoheadrightarrow \fX(\fl\cap [\fm,\fm])$, respectively.  
%
    \end{prop}
    \begin{proof}
   Since $\eta'$ coincides with (\ref{eq:defofbareta1}), $H^2(G/Q)\rightarrow H^2(\widetilde{Y}^{\mathrm{reg}},\CC)$ is indeed the inclusion $\fX(\fm)\hookrightarrow \fX(\fl)$. Next we show that $H^2(\widetilde{Y}^{\mathrm{reg}},\CC)\rightarrow 
    H^2(\widetilde{Y}_M^{\mathrm{reg}},\CC)$ coincides with $\fX(\fl)\twoheadrightarrow 
    \fX(\fl\cap [\fm,\fm])$. There is a commutative diagram
    \begin{center}
        \begin{tikzcd}
          \widetilde{Y}_M^{\mathrm{reg}} \ar[d] \ar[r] & \widetilde{Y}^{\mathrm{reg}} \ar[d] \\
          M/P_M \ar[r] & G/P
        \end{tikzcd}
    \end{center}
    Here the vertical maps are the projections, while the horizontal maps are the fiber inclusions (for fiber bundles over $G/Q$). We have $H^2(G/P,\CC)\simeq \fX(\fl), H^2(M/P_M,\CC)\simeq \fX(\fl\cap [\fm,\fm])$, and the pullback map for the bottom arrow is the projection $\fX(\fl)\twoheadrightarrow \fX(\fl\cap [\fm,\fm])$. Since the pullback maps for the vertical arrows are $\eta_M$ and $\eta$, we are done.
   
        
%
\end{proof}

\section{Partial resolutions from codimension 2 leaves}\label{subsec:descriptionpartial} 

In this section, we will provide a Lie-theoretic description of the partial Namikawa spaces $\fP_k^{\widetilde{X}}$, $k \geq 1$, under some additional restrictions. Assume $G$ is semisimple in order to simplify the statements.

Choose a Levi subgroup $L \subset G$ and a birationally rigid $L$-equivariant nilpotent cover $\widetilde{\OO}_L$ such that $\widetilde{\OO} = \mathrm{Bind}^G_L \widetilde{\OO}_L$. Choose a codimension 2 leaf $\fL_k \subset \widetilde{X}$ and let $\Sigma_k$ be the corresponding Kleinian singularity. There is a closed embedding $\Sigma_k \hookrightarrow \widetilde{X}$, constructed as follows. The image of $\overline{\fL}_k$ under the moment map $\mu: \widetilde{X} \to \overline{\OO} \subset \fg^*$ is the closure of a codimension $2$ orbit $\OO' \subset \overline{\OO}$. Consider the Slodowy slice $S'$ to $\OO'$ in $\g^*$. Recall from Section \ref{subsec:W} that $S'$ comes with a natural contracting $\CC^\times$-action. Since $S'$ is transverse to $\OO'$, $\mu^{-1}(S')$ is transverse to every leaf in $\mu^{-1}(\OO')$. Thanks to the contracting $\CC^\times$-action on $S'$, $\mu^{-1}(S')$ splits into a disjoint union of connected components indexed by the points in $\mu^{-1}(S'\cap \OO')$. Choose a point in this set lying in $\fL_k$ and let $\Sigma_k$ be a  connected component of $\mu^{-1}(S')$ containing this point. Note that $\Sigma_k$ is a Kleinian singularity and $\operatorname{Spec}(\C[\Sigma_k]^{\wedge})$ is a formal slice to $\fL_k$. 

Assume for the remainder of this section that $H^2(\widetilde{\OO},\CC)=0$. By Lemma \ref{lem:computeH2}, this is equivalent to the condition $\fP_0^{\widetilde{X}}=0$. Under this assumption, we will define, for each codimension 2 leaf $\fL_k \subset \widetilde{X}$, a Levi subgroup $M_k \subset G$ containing $L$ which is `adapted' to $\fL_k$. We will show that this Levi subgroup is uniquely determined by the following property: for any parabolic $Q \subset G$ with Levi factor $M_k$, there is a partial resolution $\overline{\rho}_k:\widetilde{Z}_k = G\times^Q (\widetilde{X}_{M_k} \times \fq^{\perp}) \to \widetilde{X}$ such that:
\begin{itemize}
\item $\overline{\rho}_k$ is a minimal resolution over $\Sigma_k\subset \widetilde{X}$.
\item $\overline{\rho}_k$ is an isomorphism over $\Sigma_j \subset \widetilde{X}$ $j\neq k$.
\item $H^2(\widetilde{\OO}_{M_k},\CC)=0$.
\end{itemize}
In other words, $\rho_k$ resolves $\Sigma_k$ and preserves all other slices, while preserving the cohomology vanishing condition for the cover. Replacing $G$ with a covering group if necessary, we can (and will) assume that $G$ is simply connected. To simplify the notation we write $\fP$ instead of $\fP^{\widetilde{X}}$ and $\fP_{k}$ instead of $\fP^{\widetilde{X}}_{k}$. 

Recall, see Proposition \ref{prop:namikawacovers}, that a choice of parabolic $P \subset G$ with Levi factor $L$ gives rise to an identification $\eta: \fX(\fl) \xrightarrow{\sim} \fP$. Different parabolics $P$ give rise to different identifications, but these identifications coincide up to the $W^{\widetilde{X}}$-action on the target, see Proposition \ref{prop:independence}.

\begin{lemma}\label{lem:prelimfaces}
There is a parabolic $P$ with Levi factor $L$ such that the intersection $\mathfrak{C}:=\fX(\fl)_{\RR}^{\geqslant 0}\cap \eta^{-1}(\fP_{\RR,k}^{\geqslant 0})$ is a face in $\fX(\fl)_{\RR}^{\geqslant 0}$ of dimension $\dim \fP_{k}$.
\end{lemma}

\begin{proof}
Since $\fP_{\RR,0} \simeq H^2(\widetilde{\OO},\RR)=0$, the product $\prod_{i}\fP_{\RR,i}^{\geq 0}$ is a fundamental chamber for the $W^{\widetilde{X}}$-action on $\fP_{\RR}$, and $\fP_{\RR,k}$ is a face of this chamber. By Theorem \ref{thm:Q-term cones}, the ample cones of all $\QQ$-terminalizations of $\widetilde{X}$ partition the fundamental chamber. Hence, there is a $\QQ$-terminalization $\widetilde{Y}$ of $\widetilde{X}$ such that the intersection $\eta (\mathfrak{C}) := \mathrm{Amp}(\widetilde{Y}) \cap \fP_{\RR,k}$ is of dimension $\dim \fP_{\RR,k}$. Note that $\eta(\mathfrak{C})$ is automatically a face of $\mathrm{Amp}(\widetilde{Y})$. By Lemma \ref{lem:parab_sing_locus}, $\widetilde{Y} = G \times^P (\widetilde{X}_L \times \fp^{\perp})$ for a parabolic $P$ with Levi factor $L$ and $\mathrm{Amp}(\widetilde{Y}) = \eta( \fX(\fl)_{\RR}^{\geq 0})$. Hence, $\mathfrak{C} = \fX(\fl)_{\RR}^{\geq 0} \cap \eta^{-1}(\fP_{\RR,k}^{\geq 0})$ is a face of $\fX(\fl)_{\RR}^{\geq 0}$ of dimension $\dim \fP_{\RR,k}$.
\end{proof}
Choose $\mathfrak{C}$ as in Lemma \ref{lem:prelimfaces}. Let $\fm_k$ be the Levi subalgebra 
\begin{equation}\label{eq:defofmk}\fm_k := \fh \oplus \bigoplus_{\alpha \in \Delta(\mathfrak{C})} \fg_{\alpha}, \qquad \Delta(\mathfrak{C}) := \{\alpha \in \Delta \mid \langle \alpha^{\vee},\mathfrak{C}\rangle =0\},\end{equation}
and let $M_k$ be the corresponding Levi subgroup of $G$. By construction, $\fX(\fm_k)_{\RR} = \RR\mathfrak{C}$ and $\eta(\RR\mathfrak{C}) = \fP_{\RR,k}$. Thus, $\eta$ restricts to an isomorphism
\begin{equation}\label{eq:etakdef}\eta_k: \fX(\fm_k) \xrightarrow{\sim} \fP_{k}\end{equation}
\begin{rmk}\label{rmk:Mkunique}
Note that the Levi subgroup $M_k \subset G$ is uniquely determined by the triple $(L,\widetilde{\OO}_L,\fL_k)$. Indeed, $\fm_k$ is the Levi subalgebra corresponding to the co-roots which vanish on the subspace $\eta^{-1}(\fP_{k})$. And by Proposition \ref{prop:independence}, this subspace is independent of $P$. 
\end{rmk}

Next, we will give a geometric characterization of the Levi subgroup $M_k$. It will be convenient to introduce some additional terminology. 

\begin{definition}
A \emph{resolution datum}\index{resolution datum} for $L$ is a triple $(P,M,Q)$ consisting of a parabolic $P \subset G$ with Levi factor $L$, a Levi subgroup $M \subset G$ containing $L$, and a parabolic $Q \subset G$ containing $P$ with Levi factor $M$. A resolution datum is said to be \emph{adapted to the leaf } $\fL_k \subset \widetilde{X}$ if $M=M_k$. 
\end{definition}

\begin{rmk}\label{rmk:resdatum}
If $M$ and $Q$ are fixed, we can always find a resolution datum $(P,M,Q)$. However, if $M$ and $P$ are fixed, it may not be possible to do so. Nonetheless, we can always find a resolution datum $(P',M,Q)$ such that the identification $\eta': \fX(\fl) \xrightarrow{\sim} \fP$ corresponding to $P'$ coincides with $\eta$. Indeed, $\eta$ is determined by a choice of fundamental chamber for the $W^{\widetilde{X}}$-action on $\fX(\fl)_{\RR}$. 
Take $P'$ such that $\fX(\fl)_{\RR}^{\geq 0}$ intersects $\fX(\fm_k)$. Since every fundamental chamber intersects $\fX(\fm_k)$, such a $P'$ exists.
\end{rmk}
A choice of resolution datum $(P,M,Q)$ gives rise to the following package of data:
\begin{itemize}
    \item A $\QQ$-terminalization $\rho:\widetilde{Y} := G \times^P (\widetilde{X}_L \times \fp^{\perp}) \to \widetilde{X}$.
    \item A projection $\pi: \widetilde{Y} \to G/P$.
    \item An $M$-equivariant nilpotent cover $\widetilde{\OO}_M := \mathrm{Bind}^M_L \widetilde{\OO}_L$.
    \item A partial resolution $\overline{\rho}: \widetilde{Z} := G \times^Q (\widetilde{X}_M \times \fq^{\perp}) \to \widetilde{X}$.
    \item A projection $\overline{\pi}: \widetilde{Z} \to G/Q$.
    \item A $\QQ$-terminalization $\rho': \widetilde{Y} \to \widetilde{Z}$ such that $\rho = \overline{\rho} \circ \rho'$.
\end{itemize}
For an $\fL_k$-adapted resolution datum $(P,M_k,Q)$, we write $\widetilde{Z}_k$ for $\widetilde{Z}$, $\overline{\rho}_k$ for $\overline{\rho}$, $\pi_k$ for $\overline{\pi}$, and $\rho_k'$ for $\rho'$. 

\begin{rmk}\label{rmk:XMXL}
Suppose $(P,M,Q)$ is a resolution datum for $L$. Under the identifications 
$\operatorname{Pic}(\widetilde{Z})\simeq \mathfrak{X}(M)$ and $\operatorname{Pic}(\widetilde{Y})\simeq \mathfrak{X}(L)$ constructed in Section \ref{subsec:picard}, the pullback map
 $\rho^{'*}: \operatorname{Pic}(\widetilde{Z})\rightarrow
 \operatorname{Pic}(\widetilde{Y})$ coincides with the inclusion $\mathfrak{X}(M)\hookrightarrow 
 \mathfrak{X}(L)$ induced by the restriction of characters. 
\end{rmk}

\begin{prop}\label{prop:part_resol_slice}
Suppose $(P,M,Q)$ is a resolution datum for $L$ and recall that we assume that $H^2(\widetilde{\OO},\CC)=\{0\}$. Then $M=M_k$ if and only if the following conditions are satisfied:
\begin{itemize}
    \item[(i)] The map $\overline{\rho}: \overline{\rho}^{-1}(\Sigma_k) \to \Sigma_k$ is a minimal resolution, and ${\rho}': \rho^{-1}(\Sigma_k) \to \overline{\rho}^{-1}(\Sigma_k)$ is an isomorphism. 
    \item[(ii)] If $j\neq k$, the map $\overline{\rho}: \overline{\rho}^{-1}(\Sigma_j) \to \Sigma_j$ is an isomorphism, and ${\rho}': \rho^{-1}(\Sigma_j) \to \overline{\rho}^{-1}(\Sigma_j)$ is a minimal resolution. 
    \item[(iii)] $H^2(\widetilde{\OO}_M,\CC)=\{0\}$.
\end{itemize}
 \end{prop}

 \begin{proof}
 The proof is in several steps. 

 {\it Step 1}. We will show that (i) and (ii) are satisfied for any $\fL_k$-adapted resolution datum $(P,M_k,Q)$. Choose a strictly dominant weight $\tau \in \fX(M_k)^{>0}$ and consider the relatively ample line bundle $\mathcal{L}_{\widetilde{Z}_k}(\tau) = \pi_k^*\mathcal{L}_{G/Q}(\tau)$ on $\widetilde{Z}_k$. Write $\mathcal{L}_{\widetilde{Z}_k}(\tau)_i$ for the restriction of $\mathcal{L}_{\widetilde{Z}_k}(\tau)$ to $\overline{\rho}_k^{-1}(\Sigma_i)$. Similarly, define $\mathcal{L}_{\widetilde{Y}}(\tau)$ and $\mathcal{L}_{\widetilde{Y}}(\tau)_i$. 
 
Consider the line bundle $\mathcal{L}_{\widetilde{Y}}(\tau)_k$. Note that $c_1(\mathcal{L}_{\widetilde{Y}}(\tau)_k)=\eta_k(\tau)$. Hence, $\mathcal{L}_{\widetilde{Y}}(\tau)_k$ is relatively ample. There are isomorphisms
\begin{equation}\label{eq:iso}\mathcal{L}_{\widetilde{Y}}(\tau)_k\simeq [\overline{\rho}_k^*\mathcal{L}_{\widetilde{Z}_k}(\tau)]_k\simeq \overline{\rho}_k^*[\mathcal{L}_{\widetilde{Z}_k}(\tau)_k].\end{equation}
Assume that $\rho'_k: \rho_k^{-1}(\Sigma_k) \to \overline{\rho}_k^{-1}(\Sigma_k)$ is \emph{not} an isomorphism, and let $C$ be a simple curve in its exceptional divisor. Then by (\ref{eq:iso}), we have $\langle \mathcal{L}_{\widetilde{Y}}(\tau)_k, C \rangle=0$, contradicting the relative ampleness of $\mathcal{L}_{\widetilde{Y}}(\tau)_k$. Thus, $\rho'_k: \rho_k^{-1}(\Sigma_k) \to \overline{\rho}_k^{-1}(\Sigma_k)$ is an isomorphism and $\overline{\rho}_k:\overline{\rho}_k^{-1}(\Sigma_k)
 \rightarrow \Sigma_k$ is a minimal resolution. 
 
On the other hand, $\mathcal{L}_{\widetilde{Y}}(\tau)_i$ is the trivial line bundle on $\rho^{-1}(\Sigma_i)$, pulled back 
 from a relatively ample line bundle on $\overline{\rho}_k^{-1}(\Sigma_i)$. It follows that 
 $\overline{\rho}_k:\overline{\rho}_k^{-1}(\Sigma_i)\xrightarrow{\sim}\Sigma_i$ is an isomorphism and hence that 
 $\rho'_k: \rho^{-1}(\Sigma_i)\rightarrow 
 \overline{\rho}_k^{-1}(\Sigma_i)$ is a minimal resolution.

{\it Step 2}. For the next parts of the proof we need an auxiliary step. Suppose that $(P,M,Q)$
satisfies (ii).

Note that thanks to (ii), $\overline{\rho}^{-1}(\fL_j)$ for $j\neq k$ is a symplectic leaf in $\widetilde{Z}$ mapping isomorphically to $\fL_j$. The symplectic leaves in $\widetilde{Z}$ are in bijection with those in $\widetilde{X}_M$ in the following way. If $\fL^M$ is a symplectic leaf in $\widetilde{X}_M$, then $G\times^Q(\fL_M\times \mathfrak{q}^\perp)$ is a symplectic leaf in $\widetilde{Z}$. And if $\fL'$ is a symplectic leaf in $\widetilde{Z}$, then its intersection with $\widetilde{X}_M\subset \widetilde{Z}$ is a symplectic leaf in $\widetilde{X}_M$. These operations are well-defined and give mutually inverse bijections between symplectic leaves because $U^-\times_{G/Q}\widetilde{Z}\cong T^*U^-\times \widetilde{X}_M$.

Let $\fL^M_j$ denote the symplectic leaf in $\widetilde{Y}_M$ corresponding to $\fL_j\subset \widetilde{Z}$, where $j\neq k$, so that $\fL_j\cong G\times^Q(\fL_j^M\times \mathfrak{q}^\perp)$. Thanks to the latter isomorphism, the homomorphism $\pi_1(\fL^M_j)\rightarrow \pi_1(\fL_j)$ is an isomorphism. Let $\mathfrak{S}_j$ denote the minimal resolution of $\Sigma_j$. It follows that the monodromy actions on $H^2(\mathfrak{S}_j,\CC)$ for $\widetilde{Y}_M$ and $\widetilde{Y}$ coincide hence under (\ref{eq:pullback_maps}), 
the subspace $\fP_j\subset H^2(\widetilde{Y}^{\mathrm{reg}},\CC)$ maps isomorphically onto
the subspace $\fP^{\widetilde{X}_M}_j\subset \fP^{\widetilde{X}_M}$ corresponding 
to the leaf $\fL^M_j\subset \widetilde{X}_M$.

{\it Step 3}. Now let $(P,M_k,Q)$ be an $\fL_k$-adapted resolution datum. Set $M:=M_k$. In particular, $\fX(\fm)\xrightarrow{\sim} \fP_k$.
Note that $\fX(\fm)$ is the kernel of (\ref{eq:pullback_maps}). 
Combining Steps 1,2 and the assumption that $H^2(\widetilde{\OO},\C)=\{0\}$, we see that
the image of (\ref{eq:pullback_maps}) is $\bigoplus_{j\neq k} \fP^{\widetilde{X}_M}_j$. Since (\ref{eq:pullback_maps}) is surjective by the construction, we conclude that $H^2(\widetilde{\OO}_M,\CC)=\{0\}$, which is (iii).  

{\it Step 4}. We claim that $\fL^M_j,j\neq k,$ exhaust all codimension $2$ leaves in $\widetilde{X}_M$. Indeed, let $\fL^M$ be some other codimension $2$ leaf in $\widetilde{X}_M$, and let $\fL^Z$ be the corresponding (automatically codimension $2$) leaf in $\widetilde{Z}$. It image in $\widetilde{X}$ under $\overline{\rho}$ is the closure of a leaf because $\overline{\rho}$ is a proper Poisson morphism. If we prove that this image is of codimension $2$, then we arrive at a contradiction with (i) and (ii). 

Let $\OO'_M,\OO'$ be the open orbits in the images of $\fL^M,\fL^Z$ in $\overline{\OO}_M,\overline{\OO}$, respectively. 
Note that $\OO'_M$ has codimension $2$ in $\overline{\OO}_M$ because $\widetilde{X}_M\rightarrow \overline{\OO}_M$ is finite. But $\OO'$ is induced from $\OO'_M$, hence also had codimension $2$ in $\overline{\OO}$. Hence $\overline{\rho}(\fL^Z)$ must have codimension $2$ in $\widetilde{X}$, a contradiction.

{\it Step 5}. The isomorphism of Lemma \ref{lem:compatibility1} yields the decomposition $\fX(\fm)\oplus H^2(\widetilde{Y}_M^{\mathrm{reg}},\CC)\cong H^2(\widetilde{Y}^{reg},\CC)$. (\ref{eq:pullback_maps}) is the projection to the second summand and thanks to Steps 2,4 and condition (iii), it restricts to an isomorphism 
between $\bigoplus_{j\neq k} \fP_j$ and 
$H^2(\widetilde{Y}_M,\CC)$. It remains to show that $\eta(\fX(\fm))\subset H^2(\widetilde{Y}^{reg},\CC)$ lies in $\fP_k$.
Then the dimension count shows that $\eta(\fX(\fm))=\fP_k$, i.e., $(P,M,Q)$ is adapted to $\fL_k$.

 Consider the map
$$\eta_j': \fX(\fm) \hookrightarrow \fX(\fl) \xrightarrow{\sim} \fP \twoheadrightarrow \fP_j,$$
where the middle isomorphism is $\eta$. Note that $\eta_j'$ coincides with the composition 
$$\fX(\fm) \simeq H^2(\widetilde{Z},\CC) \to H^2(\widetilde{Y}, \CC)\to H^2(\mathfrak{S}_j, \CC).$$ 
Note that $\fX(\fm)\simeq H^2(\widetilde{Z},\CC)$ because $\widetilde{Z}$ is a fiber bundle over $G/Q$ with contractible fiber.

If $j\neq k$, then $\overline{\rho}^{-1}(\Sigma_j)\simeq \Sigma_j$ and therefore $\eta_j'=0$. Hence $\eta(\fX(\fm)) \subset \fP_k$, as asserted. This finishes the proof.
 \end{proof}

Suppose $(P,M_k,Q)$ is an $\fL_k$-adapted resolution datum. By Proposition \ref{prop:part_resol_slice}, there is a closed embedding $\mathfrak{S}_k \subset \widetilde{Z}_k$. This gives rise to a restriction map $\Pic(\widetilde{Z}_k) \to \Pic(\mathfrak{S}_k)$ whose image lies in $ \Pic(\mathfrak{S}_k)^{\pi_1(\fL_k)}$. The following lemma is useful for computing $\eta_k$.

 \begin{lemma}\label{lem:etakviapic}
The isomorphism $\eta_k: \fX(\fm_k) \xrightarrow{\sim} \fP_{k}$ restricts to a group homomorphism $\eta_k: \fX(M_k) \to \Lambda_k^{\pi_1(\fL_k)}$. This homomorphism corresponds to the restriction map $\Pic(\widetilde{Z}_k) \to \Pic(\mathfrak{S}_k)^{\pi_1(\fL_k)}$ under the natural identifications $\fX(M_k) \simeq \Pic(\widetilde{Z}_k)$ and $\Lambda_k^{\pi_1(\fL_k)} \simeq \Pic(\mathfrak{S}_k)^{\pi_1(\fL_k)}$.
\end{lemma}

\begin{proof}
Recall, Lemma \ref{lem:c1complexification}, that $\operatorname{Pic}(\widetilde{Y}^{\mathrm{reg}})\otimes \CC\xrightarrow{\sim} \fP$.
Note that $\eta: \fX(\fl) \xrightarrow{\sim} \fP$ is obtained by base change to $\CC$ from the abelian group homomorphism
$$\fX(L) \simeq \Pic(\widetilde{Y}) \to \Pic(\widetilde{Y}^{\mathrm{reg}}).$$
Hence, $\eta_k$ restricts to the map
\begin{equation}\label{eq:XMtopic}\fX(M_k) \subset \fX(L) \to \Pic(\widetilde{Y}) \to \Pic(\mathfrak{S}_k)^{\pi_1(\fL_k)} \simeq \Lambda_k^{\pi_1(\fL_k)}\end{equation}
By Remark \ref{rmk:XMXL}, the inclusion $\fX(M_k) \subset \fX(L)$ corresponds to the pullback map $\Pic(\widetilde{Z}_k) \to \Pic(\widetilde{Y})$ under the natural identifications $\Pic(\widetilde{Z}_k) \simeq \fX(M_k)$ and $\Pic(\widetilde{Y}) \simeq \fX(L)$. Hence, the homomorphism (\ref{eq:XMtopic}) corresponds to the restriction map $\Pic(\widetilde{Z}_k) \to \Pic(\mathfrak{S}_k)^{\pi_1(\fL_k)}$. This completes the proof.
\end{proof}

The conditions appearing in Proposition \ref{prop:part_resol_slice} are difficult to check in practice. In the remainder of this subsection, we will develop some additional tools for computing $M_k$. 

\begin{lemma}\label{lem:characterizationMk}
Suppose $(P,M,Q)$ is a resolution datum for $L$. Then $M=M_k$ if and only if the following conditions are satisfied:
\begin{itemize}
    \item[(i)] The semisimple corank of $M$ equals the dimension of $\fP_k$.
    \item[(ii)] For every $j \neq k$, the map $\overline{\rho}: \overline{\rho}^{-1}(\Sigma_j) \to \Sigma_j$ is an isomorphism.
\end{itemize}
\end{lemma}

\begin{proof}
For $M=M_k$, condition (i) is satisfied by Remark \ref{rmk:Mkunique}, while (ii) is satisfied by Proposition \ref{prop:part_resol_slice}. Next we prove that (i) and (ii) imply $M=M_k$.
Define $\eta_j': \fX(\fm) \to \fP_j$ as in Step 5 of the proof of Proposition \ref{prop:part_resol_slice}. Arguing as in that step, condition (ii) of Lemma \ref{lem:characterizationMk} implies that $\eta_j'=0$ for $j\neq k$. By \cref{prop:partialdecomp} and Lemma \ref{lem:computeH2} there is a decomposition $\fP \simeq \bigoplus \fP_j$. {Since $\eta'_j=0$ for $j\neq k$},  the image of the embedding $\fX(\fm) \subset \fX(\fl) \simeq \fP$ intersects trivially with $\bigoplus_{j\neq k} \fP_j$. {Hence} $\eta_k'$ is injective. By condition (i), this implies that $\eta'_k$ is an isomorphism, and thus that $\eta(\fX(\fm)) = \fP_k$, as desired.
\end{proof}

Our final characterization of $M_k$ is applicable only in certain special cases. Suppose first of all that $\widetilde{\OO}$ is a nilpotent orbit $\OO$ (rather than a cover). We say that $\overline{\OO}$ is \emph{normal in codimension 2} if $\overline{\OO}$ is normal along every codimension 2 orbit $\OO' \subset \overline{\OO}$. This is equivalent to the condition that every dimension 2 slice in $\overline{\OO}$ is normal (and hence a Kleinian singularity). Consider the moment map $\mu: \Spec(\CC[\OO]) \to \overline{\OO}$.

\begin{lemma}\label{lem:surjectionleaves}
If $\fL \subset \Spec(\CC[\OO])$ is a codimension 2 leaf, then $\mu(\overline{\fL})$ is the closure of a codimension 2 orbit $\OO' \subset \overline{\OO}$. This defines a {map}
\begin{equation}\label{eq:leaftoorbit}\{\text{codimension 2 leaves in } \Spec(\CC[\OO])\} \twoheadrightarrow \{\text{codimension 2 orbits in }  \overline{\OO}\}. \end{equation}
{Let $\OO_k$ denote the image of $\fL_k$.}
If $\overline{\OO}$ is normal in codimension 2, then this map is a bijection, and the singularity of $\fL_k$ is equivalent to that of $\OO_k$.
\end{lemma}

\begin{proof}
{Every leaf in $\Spec(\CC[\OO])$ contains a unique dense $G$-orbit. Since $\mu$ is finite,} (\ref{eq:leaftoorbit}) is well-defined. If $\overline{\OO}$ is normal in codimension 2, then {$\mu: \Spec(\CC[\OO]) \to \overline{\OO}$} is an isomorphism over all codimension $2$ orbits. In particular, the map (\ref{eq:leaftoorbit}) is {bijective} and preserves singularities.
\end{proof}

\begin{lemma}\label{lem:findMk}
Suppose $M \subset G$ is a Levi subgroup containing $L$. Then $M=M_k$ if the following conditions are satisfied:
\begin{itemize}
    \item[(i)] The semisimple corank of $M$ equals the dimension of $\fP_k$.
    \item[(ii)] Both $\overline{\OO}$ and $\overline{\OO}_M$ are normal in codimension 2.
    \item[(iii)] For every $j \neq k$, there is a nilpotent $M$-orbit $\OO_{M,j} \subset \overline{\OO}_M$ of codimension 2 such that
    $$\OO_j = \mathrm{Ind}^G_M \OO_{M,j}$$
    \item[(iv)] For every $j \neq k$, the singularity of $\OO_j \subset \overline{\OO}$ is equivalent to that of $\OO_{M,j} \subset \overline{\OO}_M$. 
\end{itemize}
\end{lemma}
 
\begin{proof}
Assume (i). We will show that conditions (ii)-(iv) of Lemma \ref{lem:findMk} imply condition (ii) of Lemma \ref{lem:characterizationMk}. Then Lemma \ref{lem:findMk} will follow from Lemma \ref{lem:characterizationMk}. Fix a resolution datum $(P,M,Q)$ and consider the partial resolution $\overline{\rho}:Z=G\times^Q(X_M\times \fq^\perp) \to {X}$. The restriction $\overline{\rho}: \overline{\rho}^{-1}(\Sigma_j) \to \Sigma_j$ is a partial resolution of the Kleinian singularity $\Sigma_j$. We will show that it is an isomorphism for $j\neq k$. Let $x$ be a generic point in $G\times^Q({\OO}_{M,j}\times \fq^\perp)\subset {Z}$, and let $\fL_j\subset X$ be a symplectic leaf corresponding to $\OO_j$ under (\ref{eq:leaftoorbit}). Since $\OO_j = \mathrm{Ind}^G_M \OO_{M,j}$, we have $\overline{\rho}(x)\in \fL_j$. Note that $Z^{\mathrm{reg}}\simeq G\times^Q(X_M^{\mathrm{reg}}\times \fq^\perp)$, and ${\OO}_{M,j}\subset X_M^{\mathrm{sing}}$. Note that $x$ belongs to a codimension $2$ leaf $\fL \subset {Z}$. The singularity of $\fL \subset {Z}$ is equivalent to that of $\fL_{M,j} \subset {X}_M$, and hence to $\Sigma_j$. Hence, the singularity of $\overline{\rho}^{-1}(\Sigma_j)$ is equivalent to $\Sigma_j$. This implies that  $\overline{\rho}^{-1}(\Sigma_j)\to \Sigma_j$ is an isomorphism. This is exactly condition (ii) of Lemma \ref{lem:characterizationMk}.
 \end{proof}
 
\begin{rmk}\label{rmk:findMk}
There is a version of Lemma \ref{lem:findMk} which holds for arbitrary nilpotent covers. The statement is more technical, but the proof is analogous (and we omit it). Let $\widetilde{\OO}$ be a $G$-equivariant nilpotent cover and suppose $M \subset G$ is a Levi subgroup containing $L$. Then $M=M_k$ if the following conditions are satisfied:
\begin{itemize}
    \item[(i)] The semisimple corank of $M$ equals the dimension of $\fP_k$.
    \item[(ii)] For every codimension 2 orbit $\OO' \subset \overline{\OO}$ (resp. $\OO' \subset \overline{\OO}_M$), there is at most one codimension 2 leaf $\fL \subset \Spec(\CC[\widetilde{\OO}])$ (resp. $\fL_M \subset \Spec(\CC[\widetilde{\OO}_M])$ which maps to it under $ \Spec(\CC[\widetilde{\OO}]) \to \overline{\OO}$ (resp. $ \Spec(\CC[\widetilde{\OO}_M]) \to \overline{\OO}_M$). 
    \item[(iii)] If $\OO_j \subset \overline{\OO}$ is a codimension 2 $G$-orbit which is the image of a codimension 2 leaf $\fL_j \subset \widetilde{X}$ with $j\neq k$ there is a codimension 2 $M$-orbit $\OO_{M,j} \subset \overline{\OO}_M$ which corresponds to a codimension 2 leaf $\fL_{M,j} \subset \Spec(\CC[\widetilde{\OO}_M])$ such that
    $$\OO_j = \mathrm{Ind}^G_M \OO_{M,j}$$
    \item[(iv)] For every such $\OO_j$, the singularity of $\fL_j \subset \widetilde{X}$ is equivalent to that of $\fL_{M,j} \subset \widetilde{X}_M$. 
\end{itemize}
\end{rmk}

\section{Birationally rigid and 2-leafless covers}\label{subsec:codim2leaves}

In this section, we will collect some facts about orbits which admit birationally rigid covers. As explained in the introduction, such orbits play an important role in the computation of unipotent infinitesimal characters. The following result provides a useful criterion for birational rigidity.

\begin{cor}\label{cor:criterionbirigid}
The following conditions are equivalent:
\begin{itemize}
    \item[(i)] $\widetilde{\mathbb{O}}$ is birationally rigid.
    \item[(ii)] $\CC[\widetilde{\mathbb{O}}]$ admits a unique filtered quantization.
    \item[(iii)] $\fP^{\widetilde{X}}=0$.
    \item[(iv)]  $H^2(\widetilde{\mathbb{O}},\CC)=0$ and $\Spec(\CC[\widetilde{\OO}])$ has no codimension $2$ leaves.
\end{itemize}
\end{cor}

\begin{proof}
(i) and (iii) are equivalent by Proposition \ref{prop:namikawacovers}. (iii) and (iv) are equivalent by Proposition \ref{prop:partialdecomp} and the isomorphism $\fP^{\widetilde{X}}_0 \simeq H^2(\widetilde{\OO},\CC)$ established in Lemma \ref{lem:computeH2}. Finally, (ii) and (iii) are equivalent by Proposition \ref{thm:quantsofsymplectic}.
\end{proof}

Sometimes it will be convenient to consider a larger class of covers.

\begin{definition}\label{defi:2_leafless}
Let $\widetilde{\mathbb{O}}$ be a nilpotent cover. We say that $\widetilde{\mathbb{O}}$ is \emph{2-leafless}
if $\Spec(\CC[\widetilde{\mathbb{O}}])$ has no codimension 2 leaves.\index{cover!2-leafless}
\end{definition}

\subsection{Linear classical groups}

Using Corollary \ref{cor:criterionbirigid} and results of Namikawa, we can give a complete classification of birationally rigid orbits in classical types. 

\begin{prop}\label{prop:birigidorbitclassical}
Suppose $\fg$ is classical and let $\mathbb{O} \subset \fg^*$ be a nilpotent orbit corresponding to a partition $p$. Then $\OO$ is birationally rigid if and only if one of the following is true:
\begin{itemize}
    \item[(i)] $\fg = \mathfrak{sl}(n)$ and $\mathbb{O} = \{0\}$.
    \item[(ii)] $\fg = \mathfrak{so}(2n+1)$ or $\mathfrak{sp}(2n)$ and $p$ satisfies
    $$p_i \leq p_{i+1}+1, \qquad \forall i.$$
    \item[(iii)] $\fg = \mathfrak{so}(2n)$, $p$ satisfies
    $$p_i \leq p_{i+1} +1, \qquad \forall i,$$
    and $p$ is not of the form $(2^m,1^2)$ for some $m$.
\end{itemize}
\end{prop}
\begin{proof}
If $\fg=\fs\fl(n)$ all induction is birational. So (i) follows from (i) of \cref{cor:rigid}. Suppose $\fg$ is of type B, C or D.  Lemma \ref{lem:surjectionleaves} {combined with the description of the codimension $2$ singularities in $\overline{\OO}$ obtained in \cite{Kraft-Procesi}} shows that $\OO$ is 2-leafless if and only if $\overline{\OO}$ has no codimension 2 orbits. So by \cref{cor:criterionbirigid} $\OO$ is birationally rigid if and only if:
    \begin{itemize}
        \item[(1)] $\overline{\OO}$ has no codimension 2 orbits, and
        \item[(2)] $H^2(\OO,\CC)=0$.
    \end{itemize}
Let $p$ be the partition corresponding to $\OO$. By the results \cite{Kraft-Procesi} (see also \cite[Prop 1.3.2]{Namikawa2009}) $\OO$ satisfies (1) if and only if
\begin{equation}\label{eq:step1}
p_{i}-p_{i+1} \leq 1, \qquad \forall i.\end{equation}
On the other hand, by \cite[Thms 5.5,5.6]{BISWAS}, $H^2(\OO,\CC) \neq 0$ if and only if $\fg=\fs\fo(2n)$ and $p$ has the following property: $p$ contains an odd part with multiplicity 2, and no other odd parts. Under the condition (\ref{eq:step1}) this can only happen if $p$ is of the form $(2^m,1^2)$ for some $m$.
\end{proof}

We now turn our attention to nilpotent covers. For $G=\mathrm{SL}(n)$, we can give a complete classification of birationally rigid covers.

\begin{prop}\label{prop:birigidcoverA}
Suppose $G=\mathrm{SL}(n)$ and let $\widetilde{\mathbb{O}}$ be a $G$-equivariant nilpotent cover. Then $\widetilde{\mathbb{O}}$ is birationally rigid if and only if it is the universal cover of a nilpotent orbit $\mathbb{O}$ corresponding to a partition $p=(d^m)$ for $m,d \in \mathbb{Z}_{>0}$ satisfying $md=n$.
\end{prop}

\begin{proof}
Let $\widetilde{\OO}$ be the universal cover of the orbit corresponding to the partition $(d^m)$. By \cite[Prop 1.9]{NamikawaQ}, $\widetilde{\mathbb{O}}$ is 2-leafless. Recall that $\mathfrak{r}$ denotes the Lie algebra of the maximal reductive subgroup in the stabilizer of an element of $\widetilde{\OO}$. By \cite[Thm 6.1.3]{CM}, $\mathfrak{r} \simeq \mathfrak{sl}(m)$. So $H^2(\widetilde{\OO},\CC)=0$ by Lemma \ref{lem:computeH2}. Thus, $\widetilde{\OO}$ is birationally rigid by Corollary \ref{cor:criterionbirigid}. Conversely, Namikawa shows in the proof of \cite[Claim 1.10.1]{NamikawaQ} that every $G$-equivariant nilpotent cover not of this form is birationally induced.
\end{proof}

\begin{rmk}\label{rmk:SLlevisbirigid}
Let $L$ be a standard Levi subgroup of $\mathrm{SL}(n)$, i.e.
$$L = \mathrm{S}(\mathrm{GL}(a_1) \times ... \times \mathrm{GL}(a_t)), \qquad \sum a_i = n$$
Adapting the proof of Proposition \ref{prop:birigidcoverA}, one can classify all birationally rigid $L$-equivariant nilpotent covers. The result is as follows: an $L$-equivariant nilpotent cover $\widetilde{\mathbb{O}}$ is birationally rigid if and only if it is the universal $L$-equivariant cover of a nilpotent orbit
$$\mathbb{O} = \mathbb{O}^1 \times  ... \times \mathbb{O}^t \subset \cN_{\mathrm{GL}(a_1)} \times ... \times \cN_{\mathrm{GL}(a_t)} = \cN_L$$
such that each $\mathbb{O}^i$ corresponds to a partition $(d^{m_i})$ of $a_i$ for a fixed integer $d$.
\end{rmk}


Next we consider the groups $\operatorname{SO}(n)$ and $\operatorname{Sp}(2n)$. By Corollary \ref{cor:criterionbirigid}, every birationally rigid cover is 2-leafless. The following is an immediate consequence of Lemma \ref{lem:leaftoleaf}.

\begin{lemma}\label{lem:cover2leafless}
Let $\widetilde{\OO}$ and $\widecheck{\OO}$ be $G$-equivariant nilpotent covers such that $\widecheck{\OO}$ covers $\widetilde{\OO}$. If $\widetilde{\OO}$ is 2-leafless, then $\widecheck{\OO}$ is 2-leafless.
\end{lemma}

The next proposition provides a classification of nilpotent orbits admitting 2-leafless covers for $G=\mathrm{SO}(n)$ or $\mathrm{Sp}(2n)$. A result like Proposition \ref{prop:birigidcoverA} for these groups is possible, but inconvenient to state and not necessary for our purposes.

\begin{prop}\label{prop:nocodim2leaves}
Suppose $G=\mathrm{SO}(n)$ or $\mathrm{Sp}(2n)$, and let $\mathbb{O}$ be a nilpotent $G$-orbit. Write $p$ for the partition corresponding to $\mathbb{O}$. Then $\mathbb{O}$ admits a 2-leafless $G$-equivariant cover if and only if one of the following is true:
\begin{itemize}
\item[(i)] $G = \mathrm{SO}(n)$, $p$ satisfies 
$$p_i \leq p_{i+1} + 2, \qquad \forall i,$$
and the inequality is strict whenever $p_i$ is even.
\item[(ii)] $G = \mathrm{Sp}(2n)$, $p$ satisfies 
$$p_i \leq p_{i+1} + 2, \qquad \forall i,$$
and the inequality is strict whenever $p_i$ is odd. 
\end{itemize}
\end{prop}

\begin{proof}
First, suppose that ${\OO}$ is an orbit satisfying (i) or (ii). Let $\widehat{\mathbb{O}}$ be the universal $G$-equivariant cover of $\OO$. By \cite[Prop 2.3 and Prop 3.4]{NamikawaQ}, $\widehat{\mathbb{O}}$ is 2-leafless. 

Conversely, suppose $\mathbb{O}$ admits a 2-leafless $G$-equivariant cover. Then the universal $G$-equivariant cover $\widehat{\mathbb{O}}$ of $\mathbb{O}$ is 2-leafless. Suppose there is an index $i \in \mathbb{Z}_{>0}$ such that $p_i \geq p_{i+1}+2$. Then by \cite[Section 3]{Kraft-Procesi}, there is a codimension 2 orbit $\mathbb{O}' \subset \overline{\mathbb{O}}$ such that the corresponding partition $q$ satisfies
$$p_i > q_i \quad \text{and} \quad p_{i+1} < q_{i+1}.$$
Write $\widehat{\mathbb{O}}'\subset \Spec(\CC[\widehat{\mathbb{O}}])$ for the preimage of $\mathbb{O}'$ under the $G$-equivariant map $\Spec(\CC[\widehat{\mathbb{O}}])\to \overline{\mathbb{O}}$. Since $\widehat{\mathbb{O}}$ is 2-leafless, $\widehat{\mathbb{O}}'\subset \Spec(\CC[\widehat{\mathbb{O}}])^{\mathrm{reg}}$, and hence by \cite[Theorem 2.6]{Mitya2020}, $p_i=q_i+1=q_{i+1}+1=p_{i+1}+2$, and  $p_j=q_j$ for $j\neq i,i+1$. In particular, the multiplicity of $p_i$ in $p$ is of different parity than its multiplicity in $q$, and thus $p_i$ is odd if $G=\operatorname{SO}(n)$ and even in $G=\operatorname{Sp}(n)$. So $p$ is of the form described in (i) or (ii).

\end{proof}

Now suppose $G$ is linear classical, and choose a nilpotent orbit $\OO$ which admits a birationally rigid cover. In the proposition below, we will find a standard Levi subgroup $L \subset G$ and a birationally rigid $L$-orbit $\OO_L$ such that $\OO=\mathrm{Bind}^G_L \OO_L$. We will also give a parameterization of the codimension 2 leaves $\fL_k \subset X$. For each $\fL_k$, we will describe the Kleinian singularity $\Sigma_k$, as well as the codimension 2 orbit $\OO_k \subset \overline{\OO}$, see Lemma \ref{lem:surjectionleaves}. 
These data will be used in \cref{sec:centralchars} to compute the unipotent infinitesimal character corresponding to a birationally rigid cover $\widetilde{\OO}$. Let us briefly explain the idea of the computation. Assume that we know the unipotent infinitesimal character corresponding to $\OO$. To compute unipotent infinitesimal character corresponding to $\widetilde{\OO}$, one needs to find the preimage of the weighted barycenter parameter (\cref{ex:barycentersymplectic}) under the isomorphism $\eta$ of (\ref{eq:defofeta}), or equivalently the preimage of the barycenter parameter (\cref{ex:barycenterKleinian}) under the isomorphism $\eta_k$ of (\ref{eq:etakdef}). This isomorphism is given in terms of the adapted Levi subgroup $M_k$. In the proposition below we describe the codimension $2$ orbits $\OO_k\subset \fL_k$, in \cref{prop:nocodim2leavesLM} we use it to compute $M_k$. For these purposes we introduce a bit of additional notation. 

\begin{definition}\label{def:S2p}
If $p$ is a partition, define
$$S_2(p) = \{i: p_i = p_{i+1}+2\}.$$
Regard $S_2(p)$ as a partition (with no repeated parts). 
\begin{itemize}
    \item If $i \in S_2(p)$, let $c_i(p)$ be the partition formed by deleting one box from $p_i$ and adding one box to $p_{i+1}$ (`c' stands for `collapse').
    \item  If $x=(x_1,...,x_r)$ is a subpartition of $S_2(p)$, let $p \# x$ be the partition formed by deleting the columns in $p$ numbered $p_{x_1}, p_{x_1}-1$, $p_{x_2}$, $p_{x_2}-1$, ..., $p_{x_r}$, $p_{x_r}-1$.
\end{itemize}
\end{definition}
For example, if $p=(6^2,4,3,2)$, then $S_2(p) = (5,2)$, $c_5(p) = (6^2,4,3,1,1)$, and $p \# S_2(p) = (2^3,1)$.

\begin{prop}\label{prop:Lforbirigidcover}
Suppose $G$ is linear classical, and let $\mathbb{O}$ be a nilpotent orbit which admits a birationally rigid $G$-equivariant cover. Write $p$ for the partition corresponding to $\mathbb{O}$. Then $p$ is of the form described in Propositions \ref{prop:birigidcoverA} (if $G=\mathrm{SL}(n)$) or Proposition \ref{prop:nocodim2leaves} (otherwise).

\begin{itemize}
    \item[(i)] Suppose $G=\mathrm{SL}(n)$. Then
    $$L=\mathrm{S}(\mathrm{GL}(m)^d), \qquad \mathbb{O}_L=\{0\}.$$
    There is a single codimension 2 leaf $\fL_1 \subset X$, and
    $$\Sigma_1 \simeq A_{d-1}, \qquad \OO_1 = \OO_{(d^{m-1},d-1,1)}.$$
    
    \item[(ii)] Suppose $G=\mathrm{SO}(2n+1), \mathrm{Sp}(2n)$, or $\mathrm{SO}(2n)$. Assume $\mathbb{O}$ is \emph{not} of the form $\mathbb{O}_{(4^{2m},3,1)}$ for $\mathrm{SO}(8m+4)$. Then 
    $$L = \prod_{k \in S_2(p)} \mathrm{GL}(k) \times G(n-|S_2(p)|), \qquad \mathbb{O}_L =  \{0\} \times ... \times \{0\} \times \mathbb{O}_{p \# S_2(p)}.$$
    The codimension 2 leaves in $X$ are parameterized by the parts of the partition $S_2(p)$. For $k \in S_2(p)$
    $$\Sigma_k \simeq A_1, \qquad \OO_k = \OO_{c_k(p)}.$$
    \item[(iii)] Suppose $G=\mathrm{SO}(8m+4)$ and $\mathbb{O} = \mathbb{O}_{(4^{2m},3,1)}$. Then 
    $$L = \mathrm{GL}(2m+1) \times \mathrm{GL}(2m+1), \qquad \mathbb{O}_L = \{0\} \times \{0\}.$$
    There are two codimension 2 leaves $\fL_1,\fL_2 \subset X$, and
    $$\Sigma_1 \simeq \Sigma_2 \simeq A_1, \qquad \OO_1 = \OO_{(4^{2m},2^2)}^{\mathrm{I}}, \qquad \OO_2 = \OO_{(4^{2m},2^2)}^{\mathrm{II}}.$$
\end{itemize}
\end{prop}

\begin{proof}
For the computation of $L$ and $\OO_L$, we refer the reader to  \cite[Cor 4.13]{Mitya2020}. The description of singularities and codimension 2 orbits is immediate from \cite[Table I]{Kraft-Procesi}.
\end{proof}

If $\OO$ is an orbit which admits a birationally rigid cover, then $H^2(\OO,\CC)=0$. Indeed, $H^2(\OO,\CC)$ embeds into the second cohomology of any cover of $\OO$ and by Corollary \ref{cor:criterionbirigid}, the second cohomology of a birationally rigid cover is $0$. Thus we can define, for each codimension 2 leaf $\fL_k \subset X$, a Levi subgroup $M_k \subset G$ containing $L$ adapted to $\fL_k$ and a nilpotent $M_k$-orbit $\OO_{M_k}$ with $\OO=\mathrm{Bind}^G_{M_k} \OO_{M_k}$, see (\ref{eq:defofmk}). Recall, Remark \ref{rmk:Mkunique}, that once $(L,\OO_L)$ is fixed, $(M_k,\OO_{M_k})$ is \emph{uniquely} determined by the codimension 2 leaf $\fL_k$. In the next proposition, we will compute $(M_k,\OO_{M_k})$ for each $\fL_k$. For these computations, we will always take $L$ to be the (standard) Levi subgroup described in Proposition \ref{prop:Lforbirigidcover}. For this choice of $L$, typically $M_k$ is \emph{not} standard. We will compute the standard Levi subgroup to which it is conjugate under $G$ (we will use the symbol `$\simeq^G$' to indicate $G$-conjugacy), as well as the partition corresponding to the nilpotent orbit $\OO_{M_k}$. We will also compute a complete set of generators $\tau_i(k)$ for the free abelian group $\fX(M_k)$, written in standard coordinates on $\fh^*$. Of course, the elements $\tau_i(k)$ determine $M_k$ uniquely (i.e. not just up to conjugacy). Having these coordinates allows us to describe the isomorphism $\eta_k$ of (\ref{eq:etakdef}) explicitly, which is done in Section \ref{subsec:identification}. This is an important step in computation the unipotent infinitesimal character corresponding to a birationally rigid cover $\widetilde{\OO}$ in \cref{sec:centralchars}. 

\begin{prop}\label{prop:nocodim2leavesLM}
Suppose $G$ is linear classical, and let $\mathbb{O}$ be a nilpotent orbit which admits a birationally rigid $G$-equivariant cover. Write $p$ for the partition corresponding to $\mathbb{O}$. Then $p$ is as described in Propositions \ref{prop:birigidcoverA} (if $G=\mathrm{SL}(n)$) or Proposition \ref{prop:nocodim2leaves} (otherwise), and $L$ is as described in Proposition \ref{prop:Lforbirigidcover}.

\begin{itemize}
    \item[(a)] Suppose $G=\mathrm{SL}(n)$. Then
    $$M_1 = L = \mathrm{S}(\mathrm{GL}(m)^d), \qquad \mathbb{O}_{M_1} = \mathbb{O}_L = \{0\},$$
    and
    $$\tau_i(1) = \frac{1}{d}(\underbrace{d-i,d-i,...,d-i}_{mi}, \underbrace{-i,-i,...,-i}_{m(d-i)}) \in \fh^*, \qquad 1 \leq i \leq d-1.$$
    \item[(b)] Suppose $G=\mathrm{SO}(2n+1), \mathrm{Sp}(2n)$, or $\mathrm{SO}(2n)$. Assume $\mathbb{O}$ is \emph{not} of the form $\mathbb{O}_{(4^{2m},3,1)}$ for $\mathrm{SO}(8m+4)$. For $k \in S_2(p)$,  
    $$M_k \simeq^G  \mathrm{GL}(k) \times G(n-k), \qquad \mathbb{O}_{M_k} =  \{0\} \times \mathbb{O}_{p \# (k)},$$
    and
    $$\tau_1(k) = (\underbrace{0,...,0}_{l(k)},\underbrace{1,...,1}_k,\underbrace{0,...,0}_{n-l(k)-k}) \in \fh^*, \qquad l(k) := \sum_{\substack{j \in S_2(p)\\ j <k}} j.$$
    \item[(c)] Suppose $G=\mathrm{SO}(8m+4)$ and $\mathbb{O} = \mathbb{O}_{(4^{2m},3,1)}$. Then
    $$M_1 \simeq^G \mathrm{GL}(4m+2)^{\mathrm{I}}, \qquad \mathbb{O}_{M_1} = \mathbb{O}_{(2^{2m+1})}, \qquad \tau_1(1) = \frac{1}{2}(\underbrace{1,...,1}_{2m+1},\underbrace{-1,...,-1}_{2m+1}) \in \fh^*,$$
    and
    $$M_2 \simeq^G \mathrm{GL}(4m+2)^{\mathrm{II}}, \qquad \mathbb{O}_{M_2} = \mathbb{O}_{(2^{2m+1})}, \qquad \tau_1(2) = \frac{1}{2}(1,...,1) \in \fh^*.$$
\end{itemize}
\end{prop}

\begin{proof}
As explained in Remark \ref{rmk:Mkunique}, the weights $\tau_i(k)$ determine a Levi subgroup $M \subset G$ ($M$ is the Levi subgroup corresponding to the coroots which vanish on each $\tau_i(k)$). In each case, it is trivial to check that $M$ contains $L$ and is conjugate to the indicated (standard) Levi subgroup.
We will show that $M=M_k$ using Lemma \ref{lem:findMk}. Condition (i) of that lemma is trivial in all cases. Condition (ii) can be easily verified using \cite[Table I]{Kraft-Procesi}. This leaves conditions (iii) and (iv). In case (a), $\Spec(\CC[\OO])$ contains a unique codimension 2 leaf. So conditions (iii) and (iv) are vacuous. For case (b), we argue as follows. By Proposition \ref{prop:Lforbirigidcover}, $\OO_j = \OO_{c_j(p)}$. Let $\OO_{M_k,j} = \{0\} \times \OO_{c_j(p\#(k))}$. The expression $c_j(p\#(k))$ makes sense since $j \in  S_2(p) - \{k\} = S_2(p\#(k))$.  By \cite[Table I]{Kraft-Procesi}, $\OO_{M_k,j}$ is a codimension 2 orbit in $\overline{\OO}_{M_k}$, and the singularity of $\OO_{M_k,j} \subset \overline{\OO}_{M_k}$ is of type $A_1$, and hence equivalent to that of $\OO_j \subset \overline{\OO}$. This proves condition (iv) of Lemma \ref{lem:findMk}. Condition (iii) becomes
$$\OO_{c_j(p)} = \mathrm{Ind}^G_{\mathrm{GL}(k) \times G(n-k)} \{0\} \times \OO_{c_j(p \#(k))}$$
By Proposition \ref{prop:inductionclassical}, the right hand side corresponds to the partition
$$c_j(p \# (k)) + 2(\underbrace{1,...,1}_{k}) = c_j[p \# (k) + 2(\underbrace{1,...,1}_{k})] =  c_j(p),$$
which proves condition (iii). 
We proceed to case (c). Suppose $k=1$, $j=2$ (the other case, namely $k=2$, $j=1$, is analogous). Note that the Levi subgroup $M_1$ determined by $\tau_1(1)$ is $G$-conjugate to the standard Levi $\mathrm{GL}(4m+2)^{\mathrm{II}} \subset \mathrm{SO}(8m+4)$. By Proposition \ref{prop:Lforbirigidcover}, $\OO_2 = \OO_{(4^{2m},2^2)}^{\mathrm{II}}$. Set $\OO_{M_1,2} = \OO_{(2^{2m},1^2)}$. By \cite[Table I]{Kraft-Procesi}, $\OO_{M_1,2}$ is a codimension 2 orbit in $\overline{\OO}_{M_1}$, and the singularity of $\OO_{M_1,2} \subset \overline{\OO}_{M_1}$ is of type $A_1$, and hence equivalent to that of $\OO_2 \subset \overline{\OO}$. It remains to show that
\begin{equation}\label{eq:check}\OO_{(4^{2m},2^2)}^{\mathrm{II}} = \mathrm{Ind}^G_{GL(4m+2)^{\mathrm{II}}} \OO_{(2^{2m},1^2)}.\end{equation}
This follows at once from (iii) of \cite[Thm 7.3.3]{CM}.
\end{proof}

\subsection{Spin groups}

Next, we describe all nilpotent orbits for $\mathrm{SO}(n)$ which admit $\mathrm{Spin}(n)$-equivariant 2-leafless covers. Of course, all orbits described in Proposition \ref{prop:nocodim2leaves}(i) have this property. But there are others.

\begin{prop}\label{prop:nocodim2leavesspin}
Let $G=\mathrm{SO}(n)$ and let $\mathbb{O}$ be a nilpotent $G$-orbit corresponding to a partition $p$ of $n$. Then $\mathbb{O}$ admits a 2-leafless cover, which is \emph{not} $G$-equivariant, if and only if the following conditions are satisfied:
\begin{itemize}
    \item[(i)] $p$ is rather odd (i.e. every odd part occurs with multiplicity 1).
    \item[(ii)] $p_i \leq p_{i+1}+1$ if $p_i$ is even, and $p_i \leq p_{i+1}+4$ if $p_i$ is odd.
   \item[(iii)] $p_i \neq p_{i+1}+3$ for all $i$.
\end{itemize}
\end{prop}

\begin{proof}
First, suppose $p$ satisfies conditions (i)-(iii). Condition (i) implies that $\pi_1(\OO) \neq \pi_1^G(\OO)$, see Section \ref{subsec:nilpcovers}. Hence, the universal cover $\widehat{\OO}$ of $\OO$ cannot be $G$-equivariant. And by \cite[Prop 3.6]{NamikawaQ}, $\widehat{\OO}$ is 2-leafless. 

Conversely, suppose $\OO$ admits a 2-leafless cover which is not $G$-equivariant. Then $p$ is rather odd by Section \ref{subsec:nilpcovers}.  Suppose, for contradiction, that either (ii) or (iii) is false. Then it follows from the discussion after the proof of \cite[Lem 3.9]{NamikawaQ} that the universal cover $\widehat{\OO}$ is birationally induced from a proper Levi subgroup. Note that Lemma \ref{lem:leaftoleaf} implies that $\widehat{\OO}$ is 2-leafless. Thus by Corollary \ref{cor:criterionbirigid}, $H^2(\widehat{\OO},\CC)$ is nonzero. By Lemma \ref{lem:computeH2}, $H^2(\widehat{\OO},\CC)$ is identified with $\fX(\mathfrak{r})$, where $\mathfrak{r}$ is the reductive part of the centralizer of $e \in \OO$. For orbits in classical types, $\mathfrak{r}$ is described in \cite[Thm 6.1.3]{CM}. For $\fg=\mathfrak{so}(n)$, $\mathfrak{r}$ is semisimple unless $p$ contains an odd part of multiplicity 2. This is impossible in our case, since $p$ is rather odd. So $\mathfrak{r}$ is semisimple and $H^2(\widehat{\OO},\CC)=0$, a contradiction.
\end{proof}

If $p$ is a partition, define
$$S_4(p) = \{i: p_i = p_{i+1} + 4\}.$$
We will regard $S_4(p)$ as a partition (with no repeated parts). If $x \subset S_4(p)$ is a subpartition, define $p\# x$ as in Definition \ref{def:S2p}. For example, if $p=(6,2,1^2)$, then $S_4(p) = (1)$, and $p\#S_4(p)$ is obtained from $p$ by deleting columns numbered $5$ and $6$, i.e. $p\#S_4(p) = (4,2,1^2)$.

\begin{prop}\label{prop:nocodim2leavesspinLM}
Let $G = \mathrm{SO}(n)$ and suppose $\mathbb{O}$ is a nilpotent $G$-orbit which admits a birationally rigid $\mathrm{Spin}(n)$-equivariant cover which is \emph{not} $G$-equivariant. Let $p$ be the partition corresponding to $\OO$. Then $p$ is of the form described in Proposition \ref{prop:nocodim2leavesspin}. Let $\widetilde{\OO}$ denote the universal $G$-equivariant cover of $\OO$, and let
$$L = \prod_{k \in S_4(p)} \mathrm{GL}(k) \times \mathrm{SO}(n-2|S_4(p)|), \qquad \widetilde{\mathbb{O}}_L =  \{0\} \times ... \{0\} \times \widehat{\mathbb{O}}_{p \# S_4(p)},$$
where $\widehat{\mathbb{O}}_{p \# S_4(p)}$ denotes the universal $\mathrm{SO}(n-2|S_4(p)|)$-equivariant cover of $\OO_{p \# S_4(p)}$. Then $\widetilde{\OO}_L$ is birationally rigid and $\widetilde{\OO}=\mathrm{Bind}^G_L \widetilde{\OO}_L$.

The codimension 2 leaves in $\widetilde{X}$ are parameterized by $S_4(p)$, and all corresponding singularities are of type $A_1$. If $k \in S_4(p)$, then
$$M_k \simeq^G  \mathrm{GL}(k) \times \mathrm{SO}(n-2k), \qquad \widetilde{\mathbb{O}}_{M_k} =  \{0\} \times \widehat{\mathbb{O}}_{p \# (k)},$$
and
$$\tau_1(k) = (\underbrace{0,...,0}_{l(k)},\underbrace{1,...,1}_k, \underbrace{0, ...,0}_{\lfloor \frac{n}{2} \rfloor -l(k)-k}) \in \fh^*, \qquad l(k) := \sum_{\substack{j \in S_4(p)\\ j <k}} j.$$
\end{prop}

\begin{proof}
The claim that the singularities of $\widetilde{X}$ are all of type $A_1$ and parameterized by $S_4(p)$ is a special case of \cite[Theorem 2.6]{Mitya2020}. The description of $L$ and $\widetilde{\OO}_L$ is a special case of \cite[Theorem 4.17]{Mitya2020}. For the description of $M_k$ we apply Lemma \ref{rmk:findMk} and Remark \ref{rmk:findMk}. Condition (i) of Remark \ref{rmk:findMk} is trivial. Note that $p \#(k)$ satisfies
\begin{itemize}
    \item $q_i=p_i$ if $i>k$.
    \item $q_i=p_i-2$ if $i\le k$.
\end{itemize}
Conditions (ii) and (iv) now follow from \cite[Thm 2.6]{Mitya2020}. Condition (iii) can be verified by inspecting the partitions corresponding to the boundary orbits in $\overline{\OO}$ and $\overline{\OO}_M$, see the proof of Proposition \ref{prop:nocodim2leavesLM}.
\end{proof}

\subsubsection{Exceptional groups}

Let $\OO$ be a nilpotent orbit and let $\widehat{\OO} \to \OO$ denote the universal $G$-equivariant cover. As usual, let $X=\Spec(\CC[\OO])$ and $\widehat{X} = \Spec(\CC[\widehat{\OO}])$. Below, we will develop a method for determining whether $\widehat{\OO}$ is birationally rigid. Although in the cases where we apply it, $G$ will be a simple exceptional group, the argument below is general and works for arbitrary $G$.

For the next lemma, choose a codimension 2 leaf $\fL_k \subset X$. Let $\OO_k \subset \overline{\OO}$ denote the corresponding codimension 2 orbit, see Lemma \ref{lem:surjectionleaves}.

\begin{lemma}\label{lem:inequality1}
Suppose that the covering map $\widehat{X} \to X$ is \'etale over the locus $\OO\cup \fL_k$. Then
\begin{equation}\label{eq:inequalityPX}\dim \fP^{\widehat{X}} \geq |\pi_1^G(\OO)| |\pi_1^G(\OO_k)|^{-1}.\end{equation}
\end{lemma}

\begin{proof}
Let $\widehat{\OO}_k \subset \widehat{X}$ denote the preimage of $\OO_k$ under the finite map $\widehat{X} \to \overline{\OO}$. Since $G$ acts on $\fL_k$ with finitely many orbits, there is an open $G$-orbit $\fL_k' \subset \fL_k$, which is a finite connected $G$-equivariant cover of $\OO_k$. Choose $x \in \fL_k'$. Since the map $\widehat{X} \to X$ is \'{e}tale over $x$, the preimage of $x$ contains $|\pi_1^G(\OO)|$ elements. Each connected component of $\widehat{\OO}_k$ is a $G$-equivariant cover of $\OO_k$ and therefore contains at most $|\pi_1^G(\OO_k)|$ elements in the preimage of $x$. Thus, $\widehat{\OO}_k$ contains at least $|\pi_1^G(\OO)| |\pi_1^G(\OO_k)|^{-1}$ connected components. Each such component is a codimension 2 leaf in $\widehat{X}$. Now (\ref{eq:inequalityPX}) follows from Proposition \ref{prop:partialdecomp}. 
\end{proof}

For the next lemma, define
$$\mathcal{P}_{\mathrm{rig}}(\mathbb{O}) := \{(L,\mathbb{O}_L) = (\text{Levi subgroup}, \text{rigid orbit}) \mid \mathbb{O} = \Ind^G_L \mathbb{O}_L\}/G,$$
and let $m(\OO)$ be the integer
$$m(\OO) := \mathrm{max}\{\dim(\fX(\fl)) \mid (L,\OO_L) \in \mathcal{P}_{\mathrm{rig}}(\OO)\}.$$

\begin{lemma}\label{lem:inequality2}
Let $\widetilde{\OO} \to \OO$ be a $G$-equivariant cover and let $\widetilde{X}=\Spec(\CC[\widetilde{\OO}])$. Then
\begin{equation}\label{eq:inequalityPX2}\dim \fP^{\widetilde{X}} \leq m(\OO).\end{equation}
\end{lemma}

\begin{proof}
We can assume in the proof that $G$ is semisimple. Choose a Levi subgroup $M \subset G$ and a birationally rigid $M$-equivariant nilpotent cover $\widetilde{\OO}_M$ such that $\widetilde{\OO}=\mathrm{Bind}^G_M \widetilde{\OO}_M$. Since $\OO=\Ind^G_M \OO_M$, there is a pair $(L,\OO_L) \in \mathcal{P}_{\mathrm{rig}}(\OO)$ such that $L \subset M$. In particular, $\dim (\fX(\fm)) \leq m(\OO)$. By \cref{prop:namikawacovers}, $\dim \fP^{\widetilde{X}} = \dim(\fX(\fm))$. The lemma follows.
\end{proof}

In many cases, Lemmas \ref{lem:inequality1} and \ref{lem:inequality2} can be applied in conjunction to prove that $\widehat{\OO}$ is birationally rigid. 

\begin{prop}\label{prop:birrigidcriterion}
Suppose
\begin{itemize}
    \item[(i)] The reductive part $\mathfrak{r}$ of the centralizer of $e \in \OO$ is semisimple.
    \item[(ii)] The finite subgroups of $\mathrm{Sp}(2)$ corresponding to the Kleinian singularities in $X$ are cyclic of prime order.
    \item[(iii)] For each codimension 2 orbit $\OO' \subset \overline{\OO}$, there is a strict inequality
    $$|\pi_1^G(\OO)||\pi_1^G(\OO')|^{-1} > m(\OO).$$
\end{itemize}
Then $\widehat{\OO}$ is birationally rigid. 
\end{prop}

\begin{proof}
Suppose for contradiction that $\widehat{\OO}$ is birationally induced.
Property (i), combined with Lemma \ref{lem:computeH2}(ii), implies that $H^2(\widehat{\OO},\CC)=0$.  Thus, by Corollary \ref{cor:criterionbirigid}, there is a codimension 2 leaf $\widehat{\fL} \subset \widehat{X}$. Let $\fL \subset X$ denote the image of $\widehat{\fL}$ under the covering map $\widehat{X} \to X$ and consider the codimension 2 orbit $\OO' \subset \overline{\OO}$ corresponding to $\fL$. Write $\Gamma', \widehat{\Gamma}' \subset \mathrm{Sp}(2)$ for the (nontrivial) finite subgroups corresponding to $\fL$ and $\widehat{\fL}$, respectively.  By property (ii),  $\widehat{\Gamma}'=\Gamma'$. It follows that the covering map $\widehat{X} \to X$ is \'{e}tale over $\OO \cup \fL$, see Proposition \ref{prop:almostetale}. From \cref{lem:inequality1} we deduce that $\dim \fP^{\widehat{X}}\ge |\pi_1^G(\OO)||\pi_1^G(\OO')|^{-1}$. On the other hand, Lemma \ref{lem:inequality2}, combined with (iii), implies that $\dim \fP^{\widehat{X}} < |\pi_1^G(\OO)||\pi_1^G(\OO')|^{-1}$. This is a contradiction.
\end{proof}

In exceptional types, the Lie algebras $\mathfrak{r}$ are computed in \cite[Sec 13.1]{Carter1993}. The boundary orbits $\OO' \subset \overline{\OO}$ and the corresponding singularities can be found in \cite[Tables]{fuetal2015}. For $G$ simply connected, the fundamental groups $\pi_1^G(\OO)\simeq \pi_1(\OO)$ are listed in \cite[Sec 8.4]{CM}. Finally, $\mathcal{P}_{\mathrm{rig}}(\OO)$, and hence the integer $m(\OO)$, is deducible from \cite[Sec 4]{deGraafElashvili}.

The next result describes the Kleinian singularities which can appear in $\Spec(\CC[\mathbb{O}])$ for nilpotent orbits $\mathbb{O}$ admitting 2-leafless covers. 

\begin{prop}\label{prop:almost all A1}
Suppose $G$ is a simple exceptional group. Assume $\mathbb{O}$ admits a 2-leafless $G$-equivariant cover. Then all codimension 2 singularities in $X$ are of type $A_1$, except in the following cases:
\begin{itemize}
    \item[(i)] $G = E_6$ (simply connected form) and $\mathbb{O} = 2A_2$. There is a unique codimension 2 leaf, and the corresponding singularity is of type $A_2$.
    \item[(ii)] $G=E_6$ (simply connected form) and $\mathbb{O}=A_5$. There is a unique codimension 2 leaf, and the corresponding singularity is of type $A_2$.
    \item[(iii)]$G=E_6$ (simply connected form) and $\mathbb{O}=E_6(a_3)$. There are two codimension 2 leaves, and the corresponding singularities are of types $A_1$ and $A_2$.
    \item[(iv)] $G=E_8$ and $\mathbb{O} = E_8(b_6)$. There are two codimension 2 leaves, and the corresponding singularities are of types $A_1$ and $A_2$.
\end{itemize}
In each case listed above, the universal $G$-equivariant cover $\widehat{\mathbb{O}} \to \mathbb{O}$ is birationally rigid.
\end{prop}

\begin{proof}
First, we will show that the cases listed above are the only possible exceptions. Fix a codimension 2 leaf $\fL_k \subset X=\Spec(\CC[\OO])$. We claim first of all that the natural map $\phi_k: \Gamma_k\to \pi_1(\OO)$ defined in (\ref{eq:defofphik}) is an injection. 

Recall that $\Sigma_k$ embeds into $X$ (see the discussion at the beginning of Section \ref{subsec:descriptionpartial}). 
Choose a small neighborhood $N$ of $0$ in $\Sigma_k \subset X$ and a point $x\in N- \{0\}$. Let $\widehat{\mathbb{O}} \to \mathbb{O}$ be the universal cover and let $p: \widehat{X} =\Spec(\CC[\widehat{\OO}]) \to  X$ be the induced map of affine varieties. Fix a preimage $\widehat{x} \in p^{-1}(x)$, and let $\widehat{N}$ denote the connected component of $p^{-1}(N)$ containing $\hat{x}$. Since $\widehat{X}$ is 2-leafless, we can assume that $\widehat{N}$ is a disc and the map $p:\widehat{N}\rightarrow N$ is the quotient map for the $\Gamma_k$-action on $\widehat{N}$.
For each $\widehat{x}_i \in p^{-1}(x)\cap \widehat{N}$, choose a continuous path $\gamma_i$ inside $\widehat{N}$ connecting $\widehat{x}$ to $\widehat{x}_i$. Note that there are $|\Gamma_k|$ distinct points $\widehat{x}_i$. By construction, the loops $p(\gamma_{i})$ exhaust the subgroup $\phi_k(\Gamma_k)$. Since there are $|\Gamma_k|$ distinct elements, the homomorphism $\phi_k$ is injective, as asserted. 

For simple exceptional $G$, the fundamental groups $\pi_1(\mathbb{O})$ are listed in \cite[Sec 8.4]{CM}. The dimension 2 singularities $\Gamma_k$ are described in \cite[Sec 13]{fuetal2015}. The four orbits listed in the statement of the proposition are the only orbits with the following two properties:
\begin{itemize}
    \item $X$ contains a dimension 2 singularity $\Sigma_k$ not of type $A_1$, and
    \item $\Gamma_k$ admits an embedding into $\pi_1(\OO)$.
\end{itemize}
To complete the proof, we will show that for each of these four orbits, the universal $G$-equivariant cover $\widehat{\mathbb{O}}$ is birationally rigid. For orbits (i), (ii), and (iv) we will do so using a straightforward application of Proposition \ref{prop:birrigidcriterion}. Orbit (iii) will require a slightly more elaborate argument, given below.
 
\vspace{3mm}
$(i)\ \underline{G=E_6, \mathbb{O}=2A_2.}$ We have $\pi_1(\OO) \simeq \ZZ_3$ and $\mathfrak{r} = G_2$. By \cite[Sec 4]{deGraafElashvili}, $\mathcal{P}_{\mathrm{rig}}(\OO) = \{(D_4,\{0\})\}$, and hence $m(\OO)=2$. There is one codimension 2 orbit in $\overline{\OO}$, indicated below.
\vspace{3mm}
\begin{center}
\begin{tabular}{|c|c|c|c|} \hline
$k$ & $\OO_k$ & $\Sigma_k$ & $\pi_1(\OO_k)$ \\ \hline
$1$ & $A_2+A_1$ & $A_2$& $1$\\ \hline
\end{tabular}
\end{center}
\vspace{3mm}
Note that 
$$|\pi_1(\OO)||\pi_1(\OO_1)|^{-1}=3>2=m(\OO).$$
Thus, $\widehat{\OO}$ is birationally rigid by \cref{prop:birrigidcriterion}.
 
\vspace{3mm}
$(ii)\ \underline{G=E_6, \mathbb{O}=A_5.}$ We have $\pi_1(\OO) \simeq \ZZ_3$ and $\mathfrak{r} = A_1$. By \cite[Sec 4]{deGraafElashvili}, $\mathcal{P}_{\mathrm{rig}}(\OO) = \{(D_4,  (3,2^2,1)\}$, and hence $m(\OO)=2$. There is one codimension 2 orbit in $\overline{\OO}$, indicated below.
\vspace{3mm}
\begin{center}
\begin{tabular}{|c|c|c|c|} \hline
$k$ & $\OO_k$ & $\Sigma_k$ & $\pi_1(\OO_k)$ \\ \hline
$1$ & $A_4+A_1$ & $A_2$& $1$\\ \hline
\end{tabular}
\end{center}
\vspace{3mm}
Note that 
$$|\pi_1(\OO)||\pi_1(\OO_1)|^{-1}=3>2=m(\OO).$$
Thus, $\widehat{\OO}$ is birationally rigid by \cref{prop:birrigidcriterion}.

 \vspace{3mm}
$(iii)\ \underline{G=E_6, \mathbb{O}=E_6(a_3).}$ We have $\pi_1(\OO) \simeq S_2 \times \ZZ_3$ and $\mathfrak{r} = \{0\}$. By \cite[Sec 4]{deGraafElashvili}, $\mathcal{P}_{\mathrm{rig}}(\OO) = \{(3A_1,\{0\}),(A_2,\{0\})\}$, and hence $m(\OO)=4$. There are two codimension 2 orbits in $\overline{\OO}$, indicated below.
\vspace{3mm}
\begin{center}
\begin{tabular}{|c|c|c|c|} \hline
$k$ & $\OO_k$ & $\Sigma_k$ & $\pi_1(\OO_k)$ \\ \hline
$1$ & $D_5(a_1)$ & $A_2$ & $1$\\ \hline
$2$ & $A_5$ & $A_1$ & $\ZZ_3$\\ \hline
\end{tabular}
\end{center}
\vspace{3mm}
By Lemma \ref{lem:computeH2}(ii), $H^2(\widehat{\OO},\CC)=0$. Thus, by Corollary \ref{cor:criterionbirigid}, it suffices to show that $\widehat{\OO}$ is 2-leafless, i.e. that $\Sigma_1$ and $\Sigma_2$ are smoothened under $\widehat{X} \to X$. Note that 
$$|\pi_1(\OO)||\pi_1(\OO_1)|^{-1}=6>4=m(\OO).$$
Thus, $\Sigma_1$ is smoothened under $\widehat{X} \to X$ by the argument given in the proof of \cref{prop:birrigidcriterion}. For $\Sigma_2$, the analogous inequality fails so a separate argument is needed. Suppose for contradiction that $\Sigma_2$ is not smoothened under $\widehat{X}\to X$. Consider the 2-fold cover $\widetilde{\OO} \to \OO$, and let $\widetilde{X} = \Spec(\CC[\widetilde{\OO}])$. Since $\Sigma_2$ is not smoothened under $\widehat{X}\rightarrow X$, it is not smoothened under the intermediate cover $\widetilde{X}\rightarrow X$, see \cref{lem:cover2leafless}. Since $\pi_1(\OO_2)=\ZZ_3$, the preimage $\widetilde{\OO}_2\subset \widetilde{X}$ of $\OO_2$ has 2 connected components. Thus, $\fP^{\widetilde{X}}$ contains two 1-dimensional partial Namikawa spaces. 

Since $\Sigma_1 \simeq \CC^2/\ZZ_3$ does not admit a non-trivial $2$-fold cover, the map $\widetilde{X}\to X$ is \'{e}tale over the locus $\OO\cup \fL_1$. Since $\pi_1(\OO_1)=1$, the preimage $\widetilde{\OO}_1\subset \widetilde{X}$ of $\OO_1$ has 2 irreducible components $\widetilde{\fL}_1^1$ and $\widetilde{\fL}_1^2$, each corresponding to a singularity of type $A_2$. Since for $i=1,2$ the real codimension of $\widetilde{\fL}_1^i-\OO_1$ in $\widetilde{\fL}_1^i$ is $\geq 4$, we have $\pi_1(\widetilde{\fL}_1^i)=\pi_1(\OO_1)=1$. Thus, $\fP^{\widetilde{X}}$ contains two 2-dimensional partial Namikawa spaces. This, combined with the previous paragraph, gives an inequality $\dim \fP^{\widetilde{X}} \geq 6$. On the other hand $\dim \fP^{\widetilde{X}} \leq m(\OO) = 4$ by \cref{lem:inequality2}. This is a contradiction, proving $\widehat{\OO}$ is birationally rigid.

 \vspace{3mm}
$(iv)\ \underline{G=E_8, \mathbb{O}=E_8(b_6)}.$ We have $\pi_1(\OO) \simeq S_3$ and $\mathfrak{r} = \{0\}$. By \cite[Sec 4]{deGraafElashvili}
$$\mathcal{P}_{\mathrm{rig}}(\OO) = \{(A_1+E_6,\{0\} \times 2A_2+A_1), (A_1+A_2+A_3,\{0\}), (A_2+D_4, \{0\} \times (2^2,1^4))\},$$
and hence $m(\OO)=2$. There are two codimension 2 orbits in $\overline{\OO}$, indicated below.
\vspace{3mm}
\begin{center}
\begin{tabular}{|c|c|c|c|} \hline
$k$ & $\OO_k$ & $\Sigma_k$ & $\pi_1(\OO_k)$ \\ \hline
$1$ & $E_6(a_1)+A_1$ & $A_2$& $S_2$\\ \hline 
$2$ & $A_7$ & $A_1$ & $1$ \\ \hline
\end{tabular}
\end{center}
\vspace{3mm}
Note that 
$$|\pi_1(\OO)||\pi_1(\OO_1)|^{-1}=3>2=m(\OO),$$
and
$$|\pi_1(\OO)||\pi_1(\OO_2)|^{-1}=6>2=m(\OO).$$
Thus, $\widehat{\OO}$ is birationally rigid by \cref{prop:birrigidcriterion}.

\end{proof}

\section{Computation of $\eta_k$}\label{subsec:identification}

Suppose $G$ is semisimple and simply connected. Let $\widetilde{\mathbb{O}}$ be a $G$-equivariant nilpotent cover. Fix the notation of Section \ref{subsec:terminalizationcover}, i.e. $\widetilde{X}$, $L$, $P$, $\widetilde{\OO}_L$, $\widetilde{Y}$, $\rho: \widetilde{Y} \to \widetilde{X}$, and so on. For the remainder of this subsection, we will also assume
\begin{itemize}
    \item[(a1)] $\widetilde{\OO}$ admits a birationally rigid cover. 
\end{itemize}
As explained in the introduction {of this chapter}, such covers play an important role in the computation of unipotent infinitesimal characters. 

Choose a codimension 2 leaf $\fL_k \subset X$ and let $\Sigma_k = \CC^2/\Gamma_k$ be the corresponding Kleinian singularity. Assumption (a1) guarantees that $H^2(\widetilde{\OO},\CC)=0$, see Corollary \ref{cor:criterionbirigid}. Thus, we can define a Levi subgroup $M_k \subset G$ which is adapted to $\fL_k$, see (\ref{eq:defofmk}). Recall that $L \subset M_k$ and $\eta: \fX(\fl) \xrightarrow{\sim} \fP$ restricts to an isomorphism
\begin{equation}\label{eq:etak}\eta_k: \fX(\fm_k) \xrightarrow{\sim} \fP_{k}.\end{equation}
For simplicity, we will impose the following additional assumption:
\begin{itemize}
    \item[(a2)] $\pi_1(\fL_k)$ acts trivially on $H^2(\mathfrak{S}_k,\CC)$.
\end{itemize}

This is not an unreasonable condition. {First, we claim that (a1) always implies (a2) when $G$ is linear classical.}  
If $\fg$ is of type $A$, by \cref{prop:birigidcoverA} orbit $\OO$ corresponds to the partition $(d^m)$. The unique codimension $2$ orbit in $\Spec(\CC[\OO])$ corresponds to the partition $(d^{m-1}, d-1, 1)$ and by the discussion in Section \ref{subsec:nilpcovers}, $\pi_1(\fL_k)$ is trivial. Suppose now that $G=\mathrm{Sp}(2n)$ or $G=\mathrm{Spin}(n)$, and $\OO$ admits a birationally rigid $\mathrm{SO}(n)$-cover. Let $p$ be the partition corresponding to the orbit $\OO$. \cref{prop:nocodim2leaves} shows that $p_i\le p_{i+1}+2$. \cite[Section 3]{Kraft-Procesi} implies that all dimension $2$ singularities of $\Spec(\CC[\OO])$ are of type $A_1$. The action of $\pi_1(\fL_k)$ on $H^2(\fS_k, \CC)$ comes from the action on the root system $\Delta_k$ of type $A_1$ by diagram automorphisms, and therefore is trivial.

Similarly, in exceptional types, by \cref{prop:almost all A1}, in almost all cases (a1) implies that all dimension $2$ singularities of $\Spec(\CC[\OO])$ are of type $A_1$. Since $X$ is a finite covering of $\Spec(\CC[\OO])$, all dimension $2$ singularities in $X$ are of type $A_1$ as well. 
So, (a1) implies (a2) for the same reason as for $G=\operatorname{SO}(n)$ and $\operatorname{Sp}(2n)$.

If $G=\mathrm{Spin}(n)$, the condition (a1) does not necessarily imply (a2). Let $\OO$ be an orbit that admits a birationally rigid $\mathrm{Spin}(n)$-cover but has no birationally rigid $\mathrm{SO}(n)$-covers, and let $p$ be the partition corresponding to $\OO$. Then \cref{prop:nocodim2leavesspin} implies that $p_i\le p_{i+1}+4$ for all $i$. It follows that all dimension $2$ singularities of $\Spec(\CC[\OO])$ are of type $A_1$ or $A_3$, and the latter are parameterized by the set $S_4(p)$. Let $\fL_k$ be a codimension $2$ leaf corresponding to an $A_3$ singularity, parameterized by $k\in S_4(p)$. By \cite[Theorem 2.3.1]{Mitya2020}, (a2) is satisfied if and only if $p_k$ and $p_{k+1}$ are the only odd members of $p$. Together with the conditions of \cref{prop:nocodim2leavesspin}, that implies $p=(6^{2m}, 5, 1)$.

{It remains to analyze the $4$ orbits listed in \cref{prop:almost all A1}. For each of the orbits $2A_2$, $A_5$, and $E_6(a_3)\subset E_6$, its closure has unique codimension $2$ orbit with $A_2$ singularity, see \cref{prop:almost all A1}. This orbit ($A_2+A_1$, $A_4+A_1$, and $D_5(a_1)$, respectively) has trivial fundamental group, and therefore $\pi_1(\fL_k)$ is trivial. 
For simple $G$ the orbit $E_8(b_6)\subset E_8$ is the only such that (a1) holds, but (a2) does not. For the codimension $2$ orbit $E_6(a_1)+A_1$ in the closure of the orbit $E_8(b_6)$, the corresponding Kleinian singularity is of type $A_2$, and the monodromy action is non-trivial, see \cite[Section 13]{fuetal2015}.}

Assuming condition (a2), $\fP_{k}$ can be identified with the vector space $\fh_{k}^*$, i.e. the dual Cartan subalgebra corresponding to the Kleinian singularity $\Sigma_k$. In particular, $\fP_{k}$ admits a natural basis consisting of fundamental weights, denoted $\{\omega_i(k) \mid 1 \leq i \leq n(k)\}$. On the other hand, $\fX(\fm_k)$ admits a natural basis consisting of dominant generators for the free abelian group $\fX(M_k)$, denoted $\{\tau_i(k) \mid 1 \leq i \leq n(k)\}$. Our goal in this section is to describe $\eta_k$ in terms of these bases. We will achieve this goal for a large class of covers, including:
\begin{itemize}
    \item All nilpotent orbits for linear classical groups which admit birationally rigid covers. See Corollary \ref{cor:lambdaA} for type $A$ and Corollary \ref{cor:lambdaBCD} for types $B$, $C$, and $D$.
    \item All universal $\mathrm{SO}(n)$-equivariant covers which admit birationally rigid $\mathrm{Spin}(n)$-equivariant covers. See Corollary \ref{cor:lambdaspin}.
    \item All nilpotent orbits for simple exceptional groups which admit birationally rigid covers, except for the 4 orbits listed in Proposition \ref{prop:nocodim2leaves}. See Proposition \ref{prop:identification2}, as well as Example \ref{ex:exceptional1}.
\end{itemize}
In all cases above, condition (a2) is satisfied. 

{
Let us explain how our computation is structured. A typical situation is when the slice $\Sigma_k$ is of type $A_1$.
Then our map is given by a scalar. Proposition \ref{subsec:A1} expresses this scalar in terms of orders of certain Picard groups. Then we carry the computations in three different cases. In Section \ref{SSS_case2_1}, we handle the case when $G=\operatorname{SO}(n)$ or $\operatorname{Sp}(2n)$. In Section \ref{SSS_case2_2}, we consider the case when $\fg=\mathfrak{so}(n)$, and $\widetilde{\OO}$ is equivariant with respect to $\operatorname{Spin}(n)$ but not $\operatorname{SO}(n)$. Finally, in Section \ref{SSS_case2_3} we handle the computations for exceptional groups (except the four orbits mentioned above). 
}

{The typical case mentioned above is referred to as Case 2, because Case 1 is what happens in type $A$. This case is handled in Proposition \ref{prop:identification1} in an ``abstract situation'', then the specific case of $\operatorname{SL}(n)$ is Corollary \ref{cor:lambdaA}.}

\subsubsection{Case 1}\label{subsec:rectangularA}

\begin{prop}\label{prop:identification1}
Assume conditions (a1)-(a2) as well as
\begin{itemize} 
    \item[(b1)] $\fL_k=\fL_1$ is the unique codimension 2 leaf in $\widetilde{X}$. 
    \item[(b2)] $\Pic(\widetilde{\mathbb{O}}) \simeq \fX(\Gamma_1)$.
    \item[(b3)] $\Pic(\widetilde{\mathbb{O}}_{M_1}) = 0$.
    \item[(b4)] Up to the action of $W^{\widetilde{X}}$, there is a unique parabolic subgroup $P \subset G$ with Levi factor $L$.
\end{itemize}
Then $M_1=L$, $\eta_1=\eta$, and $\eta_1\{\tau_i(1) \mid 1 \leq i \leq n(1)\}  = \{\omega_i(1) \mid 1 \leq i \leq n(1)\}$
\end{prop}

By condition (b1) of Proposition \ref{prop:identification1} combined with $\fP_0 \simeq H^2(\widetilde{\OO},\CC)=0$, we have that $\fP = \fP_1$. Hence, $\eta_1=\eta$. By Remark \ref{rmk:Mkunique}, $M_1$ is the (unique) Levi subgroup corresponding to the subspace $\eta_1^{-1}(\fP_1) = \fX(\fl)$. Thus, $M_1=L$ and $\widetilde{\OO}_{M_1}= \widetilde{\OO}_L$. It remains to show that $\eta\{\tau_i(1)\} = \{\omega_i(1)\}$. We will do so using Lemma \ref{lem:etakviapic}. The proof {of the proposition} will come after a lemma.

Let $E_{\widetilde{X}}$ and $E_{\Sigma}$ denote the exceptional divisors of the (partial) resolutions $\rho:\widetilde{Y} \to \widetilde{X}$ and $\mathfrak{S}_1 \to \Sigma_1$. Write $U_1,...,U_d$ for the irreducible components of $E_{\widetilde{X}}$ and $V_1,...,V_{d'}$ for the irreducible components of $E_{\Sigma}$. Since $\fL_1 \subset \widetilde{X}$ is the unique codimension 2 leaf, $E_{\widetilde{X}}$ is the closure of a fiber bundle over $\fL_1$, with fiber equal to $E_{\Sigma}$. Since the monodromy action of $\pi_1(\fL_1)$ on $\{V_i\}$ is trivial, we have that $d=d'$ and $V_i = U_i \cap E_{\Sigma}$ (up to reordering).

\begin{lemma}\label{lem:restrictionisisomorphism}
Assume the conditions of Proposition \ref{prop:identification1}. Then the restriction map $\Pic(\widetilde{Y})\to \Pic(\mathfrak{S}_1)$ is an isomorphism.
\end{lemma}

\begin{proof} In the proof we can assume that $G$ is semisimple and simply connected.

Consider the map $\ZZ^d \to \Cl(\widetilde{Y})$ which takes the $i$th generator of $\ZZ^d$ to $[U_i] \in \Cl(\widetilde{Y})$. This map induces a short exact sequence
\begin{equation}\label{eqn:ses1}
\ZZ^d \to \Cl(\widetilde{Y}) \to \Cl(\widetilde{Y} - \bigcup_{i=1}^d U_i) \to 0.
\end{equation}
Note that $\widetilde{\OO}$ is an open subset of $\widetilde{Y} - \bigcup_{i=1}^d U_i$ with complement of codimension $\geq 2$. Thus, $\Cl(\widetilde{Y} - \bigcup_{i=1}^d U_i)  \xrightarrow{\sim} \Cl(\widetilde{\OO})$. Since $\widetilde{\mathbb{O}}$ is smooth, $\Cl(\widetilde{\mathbb{O}}) \simeq \Pic(\widetilde{\mathbb{O}})$. And by condition (b2) of Proposition \ref{prop:identification1}, $\Pic(\widetilde{\mathbb{O}}) \simeq \fX(\Gamma_1)$. Hence, (\ref{eqn:ses1}) becomes
\begin{equation}\label{eqn:ses2}
\ZZ^d \to \Cl(\widetilde{Y}) \to \fX(\Gamma_1) \to 0.
\end{equation}
Now consider the map $\ZZ^d \to \Cl(\mathfrak{S}_1)$ which takes the $i$th generator of $\ZZ^d$ to $[V_i]$. This map induces a short exact sequence
\begin{equation}\label{eqn:ses3}
\ZZ^d \to \Cl(\mathfrak{S}_1) \to \Cl(\Sigma_1 ^{\times}) \to 0,
\end{equation}
where $\Sigma_1^{\times} := \Sigma_1 - \{0\}$. Since $\Sigma_1^{\times}$ is smooth and isomorphic to  $(\CC^2-0)/\Gamma_1$, there are group isomorphisms $\Cl(\Sigma_1^{\times}) \simeq \Pic(\Sigma_1^{\times})  \simeq \fX(\Gamma_1)$. Hence, (\ref{eqn:ses3}) becomes
\begin{equation}\label{eqn:ses4}
   \ZZ^d \to \Cl(\mathfrak{S}_1) \to \fX(\Gamma_1) \to 0.
\end{equation}
Since $U_i$ intersects $\mathfrak{S}_1$ transversely, $[U_i]$ maps to $[V_i]$ under the restriction map $\Cl(\widetilde{Y}) \to \Cl(\mathfrak{S}_1)$. Thus, sequences (\ref{eqn:ses2}) and (\ref{eqn:ses4}) form a commutative diagram
        
\begin{center}
\begin{tikzcd}[
  ar symbol/.style = {draw=none,"#1" description,sloped},
  isomorphic/.style = {ar symbol={\simeq}},
  equals/.style = {ar symbol={=}},
  ]
  \ZZ^d \ar[d, equal] \ar[r]& \Cl(\widetilde{Y})\ar[d, "|_{\mathfrak{S}}"] \ar[r] & \fX(\Gamma_1) \ar[r] \ar[d, "\nu"]& 0\\
   \ZZ^d \ar[r]& \Cl(\mathfrak{S}_1)\ar[r] & \fX(\Gamma_1)\ar[r] & 0
\end{tikzcd}
\end{center}
where $\nu$ denotes the induced map on quotients. By Proposition \ref{prop:descriptionofpic}, there is a short exact sequence
$$0 \to \Pic(\widetilde{Y}) \overset{\mathrm{div}}{\to} \Cl(\widetilde{Y}) \to \Pic(\widetilde{\OO}_L) \to 0$$
and $\Pic(\widetilde{Y})\simeq \fX(L)$. By condition (b3) of Proposition \ref{prop:identification1}, $\Pic(\widetilde{\OO}_L)=0$. So $\fX(L) \simeq \Pic(\widetilde{Y}) \xrightarrow{\sim} \Cl(\widetilde{Y})$. Hence, the complexification of the restriction map $\Cl(\widetilde{Y}) \to \Cl(\mathfrak{S}_1)$ coincides with the isomorphism $\eta: \fX(\fl) \xrightarrow{\sim}\fP$. Any homomorphism of free abelian groups which is an isomorphism after complexification is injective. Hence, $\nu$ is injective by the commutativity of the diagram and so an isomorphism, since $\fX(\Gamma_1)$ is finite. By the commutativity of the diagram, this implies that the restriction map $\Cl(\widetilde{Y}) \to \Cl(\mathfrak{S}_1)$ is an isomorphism.
\end{proof}

\begin{proof}[Proof of \cref{prop:identification1}]
By Lemma \ref{lem:etakviapic}, $\eta$ restricts to a group homomorphism $\eta: \fX(L) \to \Lambda_1$, and this homomorphism corresponds to the restriction map $\Pic(\widetilde{Y}) \to \Pic(\mathfrak{S}_1)$ under the natural identifications $\fX(L) \simeq \Pic(\widetilde{Y})$ and $\Lambda_1 \simeq \Pic(\mathfrak{S}_1)$. By Lemma \ref{lem:restrictionisisomorphism}, this restriction map is an isomorphism. Hence, $\eta_1$ restricts to an isomorphism of free abelian groups $\eta: \fX(L) \xrightarrow{\sim} \Lambda_1$. 

By \cref{thm:Q-term cones}(iv), $\fP_\RR^{\ge 0}$ is the union of the ample cones of all $\QQ$-terminalizations of $\widetilde{X}$. By \cref{lem:prelimfaces}(iii) and condition (b4), $\widetilde{Y}$ is the unique such $\QQ$-terminalization. Thus, $\eta: \fX(\fl)_{\RR}^{\ge 0} \xrightarrow{\sim} \fP_\RR^{\ge 0}$, and hence $\eta$ restricts to an isomorphism of free commutative monoids $\eta: \fX(L)^{\geq 0} \xrightarrow{\sim} \Lambda_1^{\geq 0}$. The source is generated by $\{\tau_i(1)\}$, and the target by $\{\omega_i(1)\}$. So $\eta$ maps $\{\tau_i(1)\}$ to $\{\omega_i(1)\}$, as asserted.
\end{proof}

\subsubsection{Case 1 for $G=\operatorname{SL}(n)$}

Let $G=\mathrm{SL}(n)$ and let $\mathbb{O}$ be a nilpotent $G$-orbit admitting a birationally rigid cover. By Proposition \ref{prop:nocodim2leaves}, $\mathbb{O}$ corresponds to a partition of the form $(d^m)$ of $n$, and by Proposition \ref{prop:Lforbirigidcover}
$$L=\mathrm{S}(\mathrm{GL}(m)^d), \qquad \mathbb{O}_L=\{0\}.$$
Recall the barycenter parameter $\epsilon \in \fP^X$ constructed in Example \ref{ex:barycentersymplectic}, and let $\delta := \eta^{-1}(\epsilon) \in \fX(\fl)$. The following corollary describes $\delta$ in standard coordinates on $\fh^*$. We can describe $\eta$ completely (up to a diagram automorphism), but Corollary \ref{cor:lambdaA} is sufficient (and a bit easier to state).

\begin{cor}\label{cor:lambdaA}
In standard coordinates on $\fh^*$
$$\delta = (\underbrace{\frac{d-1}{2d},\ldots, \frac{d-1}{2d}}_{m}, \underbrace{\frac{d-3}{2d},\ldots, \frac{d-3}{2d}}_{m},\ldots, \underbrace{\frac{1-d}{2d},\ldots, \frac{1-d}{2d}}_{m}) \in \fX(\fl).$$
\end{cor}

\begin{proof}
We first verify that conditions of Proposition \ref{prop:identification1} are satisfied. (a1) is one of our assumptions on $\OO$. By Proposition \ref{prop:Lforbirigidcover}, $X$ contains a single codimension 2 leaf $\fL_1 \subset X$ and $\Sigma_1 \simeq \CC^2/\ZZ_d$ (this proves (b1)). Since $\Pic(\OO)$ is finite, it is identified with $\Hom(\pi_1(\OO),\CC^{\times})$, which is $\ZZ_d$ by the results of Section \ref{subsec:nilpcovers} (this proves (b2)). (b3) is trivial since $\OO_{M_1}=\OO_L=\{0\}$. Note that
$$\dim(\fX(\fl))=d-1= \dim(\fh_{\mathfrak{sl}(d)}^*) = \dim H^2(\mathfrak{S}_1,\CC),$$
Thus, the monodromy action of $\pi_1(\fL_1)$ on $\fh_{\mathfrak{sl}(d)}^*$ is trivial (this proves (a2)). Finally, note that $W^{\widetilde{X}} = W_1^{\widetilde{X}} \simeq S_d \simeq N_G(L)/L$. Note that the parabolic subgroups with Levi factor $L$ are in bijection with elements of $S_d$. Thus, $P$ is the unique parabolic with Levi factor $L$ up to the action of $W^{\widetilde{X}}$  (this proves (b4)). 

Now Proposition \ref{prop:identification1} implies
$$\eta\{\tau_i(1) \mid 1 \leq i \leq d-1\} = \{\omega_i(1) \mid 1 \leq i \leq d-1\}.$$
By Proposition \ref{prop:nocodim2leavesLM}, 
$$\tau_i(1) = \frac{1}{d}(\underbrace{d-i,d-i,...,d-i}_{mi}, \underbrace{-i,-i,...,-i}_{m(d-i)}) \in \fh^*, \qquad 1 \leq i \leq d-1.$$
Hence
$$\delta = \eta^{-1}(\epsilon) = \frac{1}{d}\sum_{i=1}^{d-1} \tau_i(1) = (\underbrace{\frac{d-1}{2d},\ldots, \frac{d-1}{2d}}_{m}, \underbrace{\frac{d-3}{2d},\ldots, \frac{d-3}{2d}}_{m},\ldots, \underbrace{\frac{1-d}{2d},\ldots, \frac{1-d}{2d}}_{m})$$
as asserted.
\end{proof}

\subsubsection{Case 2}\label{subsec:A1}

By (ii) of Lemma \ref{lem:computeH2}, there is an isomorphism $\Pic(\widetilde{\OO}) \otimes_{\ZZ} \CC \xrightarrow{\sim} H^2(\widetilde{\OO},\CC)$. Recall that condition (a1) implies that $H^2(\widetilde{\OO},\CC)=0$ and hence that $\Pic(\widetilde{\OO})$ is finite.

\begin{prop}\label{prop:identification2}
Assume conditions (a1)-(a2) as well as
\begin{itemize}
    \item[(c1)] $\Sigma_k$ is of type $A_1$. 
\end{itemize}
Then $\Pic(\widetilde{\mathbb{O}}_{M_k})$ is finite, 
\begin{equation}\label{eq:defofck}c_k := 2|\Pic(\widetilde{\mathbb{O}}_{M_k})||\Pic(\widetilde{\mathbb{O}})|^{-1}\end{equation}
is an integer, and $\eta_k(\tau_1(k)) = c_k\omega_1(k)$.
\end{prop}

As in the proof of Proposition \ref{prop:identification1}, we will compute $\eta_k$ using Lemma \ref{lem:etakviapic}. To apply this result, we must first fix a resolution datum $(P',M_k,Q)$ for $L$. It may not be possible to arrange so that $P'=P$, but we can always choose $P'$ so that the resulting identification $\fX(\fl) \xrightarrow{\sim} \fP$ coincides with $\eta$, see Remark \ref{rmk:resdatum}. Thus, we can assume for convenience that $P=P'$. As usual, write $\rho_k: \widetilde{Z}_k = G \times^Q (\widetilde{X}_{M_k} \times \fq^{\perp}) \to \widetilde{X}$ for the partial resolution and $\rho_k': \widetilde{Y} \to \widetilde{Z}_k$ for the $\QQ$-terminalization. 

Since $\Sigma_k$ is of type $A_1$, the monodromy action of $\pi_1(\fL_k)$ on $\fh_k^*$ is trivial. Hence, by Lemma \ref{lem:etakviapic}(ii), $\eta_k$ restricts to a group homomorphism $\fX(M_k) \to \Lambda_k$. By Proposition \ref{prop:part_resol_slice}, the closed embedding $\Sigma_k\hookrightarrow \widetilde{X}$ lifts to a closed embedding $\mathfrak{S}_k \hookrightarrow \widetilde{Z}_k$ and by Lemma \ref{lem:etakviapic} the homomorphism  $\fX(M_k) \to \Lambda_k$ corresponds to the restriction map  $\Pic(\widetilde{Z}_k) \to \Pic(\mathfrak{S}_k)$ under the natural identifications $\fX(M_k) \simeq \Pic(\widetilde{Z}_k)$ and $\sigma:\Lambda_k \xrightarrow{\sim} \Pic(\mathfrak{S}_k)$. Since $\fP_k$ is one-dimensional, the Levi subgroup $M_k \subset G$ has semisimple corank 1, see Lemma \ref{lem:findMk}. Hence, by Proposition \ref{prop:descriptionofpic}, $\ZZ \simeq \Pic(\widetilde{Z}_k)$ via $1 \mapsto \pi_k^*\mathcal{L}_{G/Q}(\tau_1)$. Also, $\ZZ \simeq \Pic(\mathfrak{S}_k)$ via $1 \mapsto \sigma(\omega_1)$. Let $\varphi: \ZZ \to \ZZ$ be the group homomorphism induced by the restriction map $\Pic(\widetilde{Z}_k) \to \Pic(\mathfrak{S}_k)$. To prove the proposition, it suffices to show that $\Pic(\widetilde{\OO}_{M_k})$ is finite and $\varphi(1) = c_k$.

Let $E_{\widetilde{X}}$ be the exceptional divisor of the partial resolution $\rho_k: \widetilde{Z}_k \to \widetilde{X}$. 

\begin{lemma}\label{lem:2ses}
Assume the conditions of Proposition \ref{prop:identification2}. The following are true:
\begin{itemize}
    \item[(i)] $\mathrm{rank}(\Cl(\widetilde{Z}_k))=1$.
    \item[(ii)] $\Pic(\widetilde{\mathbb{O}}_{M_k})$ is finite.
    \item[(iii)] The map $a: \ZZ \to \Cl(\widetilde{Z}_k)$ defined by $a(1) = \mathrm{div}(\pi_k^*\mathcal{L}_{G/Q}(\tau_1))$ gives rise to a short exact sequence
    $$0 \to \ZZ \overset{a}{\to} \Cl(\widetilde{Z}_k) \to \Pic(\widetilde{\mathbb{O}}_{M_k}) \to 0.$$
    \item[(iv)] The map $b: \ZZ \to \Cl(\widetilde{Z}_k)$ defined by $b(1) = [E_X]$ gives rise to a short exact sequence
    $$0 \to \ZZ \overset{b}{\to} \Cl(\widetilde{Z}_k) \to \Pic(\widetilde{\mathbb{O}}) \to 0.$$
\end{itemize}
\end{lemma}

\begin{proof}
(iii) is a special case of Proposition \ref{prop:descriptionofpic}. We proceed to proving (i). Consider the restriction map $\Cl(\widetilde{Z}_k) \to \Cl(\widetilde{\mathbb{O}})$. Since $\widetilde{\mathbb{O}}$ is an open subset of $\widetilde{Z}_k$, this map is surjective. Let $\widetilde{Z}_k^2 \subset \widetilde{Z}_k$ be the preimage under $\rho_k$ of $\widetilde{X}^{\mathrm{reg}} \cup \bigcup \fL_j \subset \widetilde{X}$. By Lemma \ref{lem:codim2}, $\codim(\widetilde{Z}_k - \widetilde{Z}_k^2,\widetilde{Z}_k) \geq 2$. Hence, there is a natural isomorphism $\Cl(\widetilde{Z}_k) \simeq \Cl(\widetilde{Z}_k^2)$. 
By Proposition \ref{prop:part_resol_slice}, $\widetilde{Z}_k^2 \cap E_{\widetilde{X}}$ is a fiber bundle over $\fL_k$  with fiber equal to $E_{\Sigma}$, and hence an irreducible variety. Thus, the kernel of the restriction map $\Cl(\widetilde{Z}_k) \to \Cl(\widetilde{\mathbb{O}})$ is generated by the class $[E_{\widetilde{X}}]$ and the following sequence is exact
\begin{equation}\label{eq:clses1}
\ZZ \overset{b}{\to} \Cl(\widetilde{Z}_k) \to \Cl(\widetilde{\mathbb{O}}) \to 0.\end{equation}
Since $\widetilde{\mathbb{O}}$ is smooth, $\Cl(\widetilde{\mathbb{O}}) \simeq \Pic(\widetilde{\mathbb{O}})$. Hence, (\ref{eq:clses1}) becomes
\begin{equation}\label{eq:clses2}
\ZZ \overset{b}{\to} \Cl(\widetilde{Z}_k) \to \Pic(\widetilde{\mathbb{O}}) \to 0.\end{equation}
Since $\Pic(\widetilde{\OO})$ is finite, the exactness of (\ref{eq:clses2}) implies $\mathrm{rank}(\Cl(\widetilde{Z}_k)) \leq 1$. The exactness of (iii) implies the opposite inequality. Hence, $\mathrm{rank}(\Cl(\widetilde{Z}_k))=1$, proving (i). Now, (i) and (iii) imply $\mathrm{rank}(\Pic(\widetilde{\mathbb{O}}_{M_k})) = \mathrm{rank}(\Cl(\widetilde{Z}_k)) - \mathrm{rank}(\ZZ) =0$. Hence, $\Pic(\widetilde{\mathbb{O}}_{M_k})$ is finite, proving (ii). Similarly, (i) and (\ref{eq:clses2}) imply $\mathrm{rank}(\ker(b)) = \mathrm{rank}(\ZZ)-\mathrm{rank}(\Cl(\widetilde{Z}_k)) + \mathrm{rank}(\Pic(\widetilde{\mathbb{O}})) =0$. Hence, $\ker(b) = 0$, proving (iv). 
\end{proof}

\begin{proof}[Proof of \cref{prop:identification2}]
By Lemma \ref{lem:2ses}, $\Pic(\widetilde{\mathbb{O}}_{M_k})$ is finite. To show that $\eta_k(\tau_1(k))=c_k\omega_1(k)$, we must show that $\varphi(1)=c_k$. We will do so by comparing the injections $a,b: \ZZ \hookrightarrow \Cl(\widetilde{Z}_k)$ constructed in Lemma \ref{lem:2ses}.

Since $\widetilde{Z}_k$ is normal, there is a natural isomorphism $\Cl(\widetilde{Z}_k) \simeq \Pic(\widetilde{Z}_k^{\mathrm{reg}})$. Let $\phi$ denote the composition
\begin{equation}\label{eq:defofphi}
\phi: \Cl(\widetilde{Z}_k) \simeq \Pic(\widetilde{Z}_k^{\mathrm{reg}}) \overset{|_{\mathfrak{S}_k}}{\to} \Pic(\mathfrak{S}_k) \simeq \ZZ.\end{equation}
Since the composition
$$\Pic(\widetilde{Z}_k) \overset{\mathrm{div}}{\to} \Cl(\widetilde{Z}_k) \simeq \Pic(\widetilde{Z}_k^{\mathrm{reg}})$$
coincides with the restriction map $\Pic(\widetilde{Z}_k) \to \Pic(\widetilde{Z}_k^{\mathrm{reg}})$, we have 
$$(\phi \circ a)(1) = \phi(\mathrm{div}(\pi_k^*\mathcal{L}_{G/Q}(\tau_1))) = (\pi_k^*\mathcal{L}_{G/Q}(\tau_1))|_{\mathfrak{S}_k}.$$
Hence, by definition, $\phi \circ a = \varphi$.

On the other hand, $(\phi \circ b)(1)$ corresponds to the line bundle $\cO(E_{\Sigma}) \in \Pic(\mathfrak{S}_k)$ under the canonical identification $\ZZ \simeq \Pic(\mathfrak{S}_k)$. Since $\Sigma_k \simeq \CC^2/\ZZ_2$ and hence $\mathfrak{S}_k\simeq T^*\mathbb{P}^1$, this integer must be 2. 

To complete the proof, we appeal to the following general fact: let $H$ be a finitely-generated abelian group of rank 1 and let $\phi: H \to \ZZ$ be an arbitrary homomorphism. Suppose there are exact sequences
$$0 \to \ZZ \overset{a}{\to} H \to A \to 0, \qquad 0 \to \ZZ \overset{b}{\to} H \to B \to 0.$$
Then 
$$(\phi \circ a)(1) = |A||B|^{-1}(\phi \circ b)(1).$$
The proof of this fact is easy and is left to the reader. Setting $H=\Cl(\widetilde{Z}_k)$, $A=\Pic(\widetilde{\mathbb{O}}_{M_k})$, and $B=\Pic(\widetilde{\mathbb{O}})$, we deduce
$$\varphi(1) = (\phi \circ a)(1) = |\Pic(\mathbb{O}_{M_k})||\Pic(\mathbb{O})|^{-1}(\phi \circ b)(1) = 2|\Pic(\mathbb{O}_{M_k})||\Pic(\mathbb{O})|^{-1} = c_k.$$
This completes the proof.
\end{proof}


\subsubsection{Case 2 for linear classical groups}\label{SSS_case2_1}

Let $G=\mathrm{Sp}(2n)$, $\mathrm{SO}(2n)$ or $\mathrm{SO}(2n+1)$ and let $\mathbb{O}$ be a nilpotent $G$-orbit admitting a birationally rigid $G$-equivariant cover. Such orbits are {among those} described in Proposition \ref{prop:nocodim2leaves}. Recall the barycenter parameter $\epsilon \in \fP^X$ constructed in Example \ref{ex:barycentersymplectic}, and let $\delta := \eta^{-1}(\epsilon) \in \fX(\fl)$. The following corollary describes $\delta$ in standard coordinates on $\fh^*$. 

\begin{cor}\label{cor:lambdaBCD}
The following are true
\begin{itemize}
    \item[(i)] Suppose $\mathbb{O}$ is \emph{not} of the form $\mathbb{O}_{(4^{2m},3,1)}$ for $G=\mathrm{SO}(8m+4)$. Then 
    $$L = \prod_{k \in S_2(p)} \mathrm{GL}(k) \times  G(n-|S_2(p)|), 
 \qquad \mathbb{O}_L =  \prod_{k \in S_2(p)} \{0\} \times \mathbb{O}_{p\#S_2(p)}.$$
    In standard coordinates on $\fh^*$
    $$\delta = \frac{1}{2}(\underbrace{1,...,1}_{|S_2(p)|},\underbrace{0,...,0}_{n-|S_2(p)|}).$$
    \item[(ii)] Suppose $G=\mathrm{SO}(8m+4)$ and $\mathbb{O}=\mathbb{O}_{(4^{2m},3,1)}$. Then 
    $$L= \mathrm{GL}(2m+1)\times \mathrm{GL}(2m+1), \qquad \mathbb{O}_L=\{0\} \times \{0\}.$$
    In standard coordinates on $\fh^*$
    $$\delta= \frac{1}{2}(\underbrace{1,...,1}_{2m+1},\underbrace{0,...,0}_{2m+1}).$$
\end{itemize}
\end{cor}

\begin{proof}
First assume $G=\mathrm{Sp}(2n)$. By Proposition \ref{prop:nocodim2leaves}, $\mathbb{O}$ corresponds to a partition $p$ of the following form
\begin{equation}\label{eq:descriptionofp}
\text{For all }i, \ p_i \leq p_{i+1}+2 \text{ and } p_i < p_{i+1}+2 \text{ if } p_i \text{ is odd},\end{equation}
By Proposition \ref{prop:Lforbirigidcover}, the codimension 2 leaves $\fL_k \subset X$ 
are parameterized by the set $S_2(p)$, and the corresponding singularities are of type $A_1$. Choose $k \in S_2(p)$. Then by Proposition \ref{prop:nocodim2leavesLM},
$$M_k \simeq^G  \mathrm{GL}(k) \times \mathrm{Sp}(2n-2k), \qquad \mathbb{O}_{M_k} =  \{0\} \times \mathbb{O}_{p \# (k)},$$
and
$$\tau_1(k) = (\underbrace{0,...,0}_{l(k)},\underbrace{1,...,1}_k,\underbrace{0,...,0}_{n-l(k)-k}) \in \fh^*, \qquad l(k) := \sum_{\substack{j \in S_2(p)\\ j <k}} j.$$
By Proposition \ref{prop:identification2}, $\eta_k(\tau_1(k))= c_k\omega_1(k)$, where $c_k = 2|\Pic(\OO_{M_k})||\Pic(\OO)|^{-1}$. Since $\Pic(\mathbb{O})$ is finite, $\Pic(\mathbb{O}) \simeq \pi_1(\mathbb{O})_{\mathrm{ab}}$. Similarly, $\Pic(\mathbb{O}_{M_k}) \simeq \pi_1(\mathbb{O}_{p\#(k)})_{\mathrm{ab}}$. Both $\pi_1(\mathbb{O})$ and $\pi_1(\mathbb{O}_{p\#(k)})$ are elementary abelian 2-groups of rank equal to the number of distinct even parts in the corresponding partition, see Section \ref{subsec:nilpcovers}.  By the description of $p$ (see (\ref{eq:descriptionofp})) and the definition of $p\#(k)$, it is clear that $p$ contains one more distinct even part than $p \#(k)$. Hence, $c_k=1$ and $\eta_k(\tau_1(k)) = \omega_1(k)$. Now the computation of $\delta$ is immediate:
$$\delta=\eta^{-1}(\epsilon) = \sum_{k \in S_2(p)} \eta_k^{-1}(\epsilon_k) = \sum_{k \in S_2(p)} \frac{1}{2}\tau_1(k) = \frac{1}{2}(\underbrace{1,...,1}_{|S_2(p)|},\underbrace{0,...,0}_{n-|S_2(p)|}) \in \fh^*.$$
If $G=\mathrm{SO}(2n+1)$ or $\mathrm{SO}(2n)$ and $\mathbb{O}$ is not of the form $\mathbb{O}_{(4^{2m},3,1)}$ for $G=\mathrm{SO}(8m+4)$, the computation is analogous. We leave the details to the reader.

Finally, assume $G=\mathrm{SO}(8m+4)$ and $\mathbb{O}=\mathbb{O}_{(4^{2m},3,1)}$. By Proposition \ref{prop:Lforbirigidcover}, there are two codimension 2 leaves $\fL_1,\fL_2 \subset X$, and the corresponding singularities are both of type $A_1$. By Proposition \ref{prop:nocodim2leavesLM},
$$M_1 \simeq^G \mathrm{GL}(4m+2)^{\mathrm{II}}, \qquad \mathbb{O}_{M_1} = \mathbb{O}_{(2^{2m+1})}, \qquad \tau_1(1) = \frac{1}{2}(\underbrace{1,...,1}_{2m+1},\underbrace{-1,...,-1}_{2m+1}) \in \fh^*,$$
    and
$$M_2 \simeq^G \mathrm{GL}(4m+2)^{\mathrm{I}}, \qquad \mathbb{O}_{M_2} = \mathbb{O}_{(2^{2m+1})}, \qquad \tau_1(2) = \frac{1}{2}(1,...,1) \in \fh^*.$$
For $k=1,2$, we have $|\Pic(\mathbb{O}_{M_k})| = |\pi_1(\mathbb{O}_{M_k})_{\mathrm{ab}}| = 2$ and $|\Pic(\mathbb{O})| = |\pi_1(\mathbb{O})_{\mathrm{ab}}|= 4$ (see Section \ref{subsec:nilpcovers} for the description of the fundamental groups). Hence $c_k=1$ and by Proposition \ref{prop:identification2} $\eta_k(\tau_1(k))=\omega_1(k)$. Thus
$$\delta = \eta^{-1}(\epsilon)  =  \frac{1}{2}\eta_1^{-1}(\omega_1(1)) + \frac{1}{2}\eta_2^{-1}(\omega_1(2)) = \frac{1}{2}\tau_1(1)+\frac{1}{2}\tau_1(2) = \frac{1}{2}(\underbrace{1,...,1}_{2m+1},\underbrace{0,...,0}_{2m+1})$$
as asserted.
\end{proof}

\subsubsection{Case 2 for Spin groups}\label{SSS_case2_2}

For the next corollary, let $G=\mathrm{SO}(n)$. Suppose $\mathbb{O}$ is a nilpotent $G$-orbit which admits a $\mathrm{Spin}(n)$-equivariant birationally rigid cover which is \emph{not} $\mathrm{SO}(n)$-equivariant. Hence, $\OO$ is one of the orbits described in Proposition \ref{prop:nocodim2leavesspin}. Let $\widetilde{\mathbb{O}} \to \mathbb{O}$ be the universal $G$-equivariant cover. Then the universal $\mathrm{Spin}(n)$-equivariant cover $\widehat{\OO}$ is birationally rigid (and a two-fold cover of $\widetilde{\OO}$). 

Let $\epsilon \in \fP^{\widetilde{X}}$ be the barycenter parameter and let $\delta:=\eta^{-1}(\epsilon)$. By Proposition \ref{prop:nocodim2leavesspinLM}, 
$$L =  \prod_{k \in S_4(p)} \mathrm{GL}(k) \times \mathrm{SO}(n-2|S_4(p)|), \qquad \widetilde{\mathbb{O}}_L = \prod_{k \in S_4(p)} \{0\} \times \widehat{\mathbb{O}}_{p\#S_4(p)},$$
where $\widehat{\mathbb{O}}_{p\#S_4(p)}$ denotes the universal $\mathrm{SO}(n-2|S_4(p)|)$-equivariant cover. 

\begin{cor}\label{cor:lambdaspin}
In standard coordinates on $\fh^*$
$$\delta = \frac{1}{4}(\underbrace{1,...,1}_{|S_4(p)|},\underbrace{0,...,0}_{\lfloor \frac{n}{2} \rfloor-|S_4(p)|}) \in \fh^*.$$
\end{cor}

\begin{proof}
By Proposition \ref{prop:nocodim2leavesspinLM}, the codimension 2 leaves $\fL_k \subset \widetilde{X}$ are parameterized by $S_4(p)$ and the corresponding singularities are of type $A_1$. If $k \in S_4(p)$, then by the same proposition
$$M_k \simeq^G  \mathrm{GL}(k) \times \mathrm{SO}(n-2k), \qquad \widetilde{\mathbb{O}}_{M_k} =  \{0\} \times \widehat{\mathbb{O}}_{p \# (k)},$$
where $\widehat{\mathbb{O}}_{p \# (k)}$ denotes the universal $\mathrm{SO}(n-2k)$-equivariant cover. Note that $p$ is rather odd by Proposition \ref{prop:nocodim2leavesspin}. Hence $p \#(k)$ is rather odd by the construction of $p\#(k)$. So $\pi_1(\widetilde{\mathbb{O}}) \simeq \pi_1(\widetilde{\mathbb{O}}_{M_k}) \simeq \ZZ_2$, see Section \ref{subsec:nilpcovers}, and 
$$c_k = 2|\Pic(\widetilde{\mathbb{O}}_{M_k})||\Pic(\widetilde{\mathbb{O}})|^{-1} = 2.$$
Hence by Proposition \ref{prop:identification2}, $\eta_k(\tau_1(k)) = 2\omega_1(k)$. By Proposition \ref{prop:nocodim2leavesLM}
$$\tau_1(k) = (\underbrace{0,...,0}_{l(k)},\underbrace{1,...,1}_k, \underbrace{0, ...,0}_{n-l(k)-k}) \in \fh^*, \qquad l(k) := \sum_{\substack{j \in S_4(p)\\ j <k}} j.$$
Hence,
$$\eta_k^{-1}(\epsilon_k) = \eta_k^{-1}(\frac{1}{2}\omega_1(k)) = \frac{1}{4}\tau_1(k) =\frac{1}{4} (\underbrace{0,...,0}_{l(k)},\underbrace{1,...,1}_k, \underbrace{0, ...,0}_{n-l(k)-k}) \in \fh^*.$$
Summing over all $k\in S_4(p)$, we obtain the desired formula for $\delta$,
$$\delta = \eta^{-1}(\epsilon) = \sum_{k \in S_4(p)} \eta_k^{-1}(\epsilon_k) =  \frac{1}{4}(\underbrace{1,...,1}_{|S_4(p)|},\underbrace{0,...,0}_{n-|S_4(p)|}) \in \fh^*.$$
\end{proof}

\subsubsection{Case 2 for exceptional groups}\label{SSS_case2_3}

\begin{example}\label{ex:exceptional0}
Let $G$ be the (unique) simple group of type $G_2$, and let $\mathbb{O} = G_2(a_1)$. By \cite[Sec 13]{fuetal2015}, $\overline{\OO}$ is normal in codimension 2, and there is a unique orbit $\OO_1 \subset \overline{\OO}$ of codimension 2, namely $\OO_1=\widetilde{A}_1$. The corresponding singularity is of type $A_1$. In this example, we will show that the universal cover $\widehat{\OO}$ of $\OO$ is birationally rigid. We will also compute the Levi subgroup $M_1 \subset G$ adapted to the codimension 2 leaf $\fL_1 \subset X=\Spec(\CC[\OO])$ as well as the identification $\eta_1: \fX(\fm_1) \xrightarrow{\sim} \fP_1^X=\fP^X$. 

Write $\alpha_1,\alpha_2$ for the simple roots for $\fg$, with $\alpha_2$ the short root, and write $\varpi_1,\varpi_2$ for the corresponding fundamental weights. Since $\OO$ is distinguished, it is birationally induced from the $\{0\}$-orbit of its Jacobson-Morozov Levi $L_{\OO}$, see Proposition \ref{prop:even} and Remark \ref{rmk:distinguishedbirational}. In this case, $L_{\OO}$ is the Levi subgroup (of type $A_1$) corresponding to the short simple root $\alpha_2$.

\emph{$\widehat{\OO}$ is birationally rigid}. Since $\OO$ is distinguished, the reductive part $\mathfrak{r}$ of the centralizer of $e \in \OO$ is $0$. Thus $H^2(\widehat{\OO},\CC) \simeq \fX(\mathfrak{r})=0$ by Lemma \ref{lem:computeH2}. So by Corollary \ref{cor:criterionbirigid}, it suffices to show that $\widehat{\OO}$ is 2-leafless. Note that $\pi_1(\OO) \simeq S_3$ and $\pi_1(\OO_1) =1$. By \cite[Sec 4]{deGraafElashvili}
$$\mathcal{P}_{\mathrm{rig}}(\OO) = \{(A_1, \{0\}), (\widetilde{A}_1,\{0\})\}$$
In particular $m(\OO)=1$. We have
$$|\pi_1(\OO)||\pi_1(\OO_1)|^{-1} = 6 > 1 = m(\OO),$$
Thus, $\widehat{\OO}$ is birationally rigid by Proposition \ref{prop:birrigidcriterion}.

\emph{Computation of $M_1$}. Since $\fL_1 \subset X$ is the unique codimension 2 leaf, it follows immediately from Lemma \ref{lem:characterizationMk} that $L=M_1$.

\emph{Computation of $\eta_1$}. Since $\OO$ admits a birationally rigid cover and $\Sigma_1 \simeq A_1$, Proposition \ref{prop:identification2} is applicable. Note that $\OO_{M_1} = \OO_L = \{0\}$. Hence
$$c_1 = 2|\Pic(\OO_{M_1})||\Pic(\OO)|^{-1} = 2|\pi_1(\OO)_{\mathrm{ab}}|^{-1} = 1,$$
and therefore by Proposition \ref{prop:identification2}
\begin{equation}\label{eq:G21}\eta_1(\tau_1(1)) = \omega_1(1).\end{equation}
By definition $\tau_1(1)$ is the dominant generator of the free abelian group $\fX(M_1)$, i.e. $\tau_1(1) = \varpi_1$. So (\ref{eq:G21}) becomes
\begin{equation}\label{eq:G22}
    \eta_1(\varpi_1) = \omega_1(1).
\end{equation}
\end{example}

\begin{example}\label{ex:exceptional1}
Let $G$ be the simply connected group of type $E_7$ and let $\mathbb{O}=E_7(a_4)$. By \cite[Sec 13]{fuetal2015}, $\overline{\OO}$ is normal in codimension 2. The codimension 2 orbits $\OO_k \subset \overline{\OO}$ as well as the singularities $\Sigma_k$ and fundamental groups $\pi_1(\OO_k)$ are indicated below
\vspace{3mm}
\begin{center}
\begin{tabular}{|c|c|c|c|} \hline
$k$ & $\OO_k$ & $\Sigma_k$ & $\pi_1(\OO_k)$ \\ \hline
$1$ & $A_6$ & $A_1$ & $1$\\ \hline
$2$ & $D_5+A_1$ & $A_1$ & $\ZZ_2$\\ \hline
$3$ & $D_6(a_1)$ & $A_1$ & $\ZZ_2$ \\ \hline
\end{tabular}
\end{center}
\vspace{3mm}
In this example, we will show that the universal cover $\widehat{\OO}$ of $\OO$ is birationally rigid. We will also compute the Levi subgroup $M_1 \subset G$ adapted to the codimension 2 leaf $\fL_1 \subset X=\Spec(\CC[\OO])$ as well as the identification $\eta_1: \fX(\fm_1) \xrightarrow{\sim} \fP_1^X$. 

Label the simple roots as follows
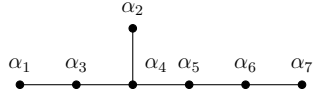
\begin{center}
$$\dynkin[labels*={\alpha_1,\alpha_2,\alpha_3,\alpha_4,\alpha_5,\alpha_6,\alpha_7},edge
length=.75cm] E7$$
\end{center}
and write $\varpi_1,...,\varpi_7$ for the fundamental weights. Arguing as in Example \ref{ex:exceptional0}, we see that $\OO = \mathrm{Bind}^G_L\{0\}$, where $L$ is the standard Levi subgroup (of Lie type $2A_1+A_2$) corresponding to the simple roots $\{\alpha_2,\alpha_3,\alpha_5,\alpha_6\}$. 

\emph{$\widehat{\OO}$ is birationally rigid}. Since $\OO$ is distinguished, the reductive part $\mathfrak{r}$ of the centralizer of $e \in \OO$ is $0$. Thus $H^2(\widehat{\OO},\CC) \simeq \fX(\mathfrak{r})=0$ by Lemma \ref{lem:computeH2}. So by Corollary \ref{cor:criterionbirigid}, it suffices to show that $\widehat{\OO}$ is 2-leafless.

Note that $\pi_1(\OO) \simeq S_2 \times \ZZ_2$ and by \cite[Sec 4]{deGraafElashvili}
\begin{equation}\label{eq:Prig}\mathcal{P}_{\mathrm{rig}}(\OO) = \{(A_1+D_4, \{0\} \times (3,2^2,1)), (2A_1+A_2,\{0\})\}\end{equation}
In particular $m(\OO)=3$. Note that
$$|\pi_1(\OO)||\pi_1(\OO_1)|^{-1}=4>3=m(\OO)$$
Thus, by the argument given in the proof of \cref{prop:birrigidcriterion}., $\Sigma_1$ is smoothened under $\widehat{X} \to X$.
 
For $\Sigma_2$ and $\Sigma_3$, a separate argument is needed. We supply one for $\Sigma_2$ (the argument for $\Sigma_3$ is analogous). Suppose that $\Sigma_2$ is \emph{not} smoothened under $\widehat{X} \to X$. Since $\Sigma_1\simeq \CC^2/\ZZ_2$, and $\widehat{X}\to X$ is generically 4 to 1, the preimage of $\Sigma_1$ under this map consists of two copies of $\CC^2$. Write $\widehat{\OO}_1$ for the preimage of $\OO_1$ under $\widehat{X} \to X$. Since $\pi_1(\OO_1) \simeq 1$, $\widehat{\OO}_1$ has two connected components. The action of $\pi_1(\OO)$ on $\widehat{X}$ permutes the preimages of a point $e\in \OO_1$. Consider the resulting map $\pi_1(\OO)\to S_2$, and let $K\simeq \ZZ_2$ be its kernel. Let $\widetilde{X}\to X$ be the cover corresponding to the subgroup $K\subset \pi_1(\OO)$. By construction, the preimage of $\OO_1$ under $\widetilde{X} \to X$ has $2$ connected components permuted by the action of $K$, and the corresponding singularities are of type $A_1$. \cref{lem:cover2leafless} implies that the singularity $\Sigma_2$ is not smoothened under the map $\widetilde{X}\to X$. Therefore, $\dim \fP^{\widetilde{X}}\ge 3$. So by (\ref{eq:Prig}), we see that both ${\OO}$ and $\widetilde{\OO}$  are birationally induced from $(2A_1+A_2,\{0\})$. This is a contradiction. Thus $\Sigma_2$ is smoothened under $\widehat{X}\to X$.

\emph{Computation of $M_1$}. Let $M$ denote the standard Levi  of Lie type $D_6$ corresponding to the simple roots $\{\alpha_2,\alpha_3,\alpha_4,\alpha_5,\alpha_6,\alpha_7\}$. Note that $L \subset M$ and by Proposition \ref{prop:inductionclassical}
$$\OO_M = \mathrm{Ind}^M_L \{0\} = \OO_{(5,3^2,1)}.$$
Using Lemma \ref{lem:findMk}, we will show that $M=M_1$. By \cite[Table 1]{Kraft-Procesi}, $\overline{\OO}_M$ is normal in codimension 2. There are two boundary orbits $\OO_{M,2}, \OO_{M,3} \subset \overline{\OO}_{M}$ of codimension 2, namely
$$\OO_{M,2} = \OO_{(4^2,3,1)}, \qquad \OO_{M,3} = \OO_{(5,3,2^2)}, $$
and the corresponding singularities are of type $A_1$. Note that $\OO_{M,2} =\Ind^{D_6}_{2A_2} \{0\}$ and $\OO_{M,3} = \Ind^{D_6}_{A_3} \{0\}$. Now by \cite[Sec 4]{deGraafElashvili}, we have $\Ind^G_{M_2} \OO_{M,2} = D_5+A_1 = \OO_2$ and $\Ind^G_{M_2} \OO_{M,3} = D_6(a_1) = \OO_3$. So $M=M_1$ by Lemma \ref{lem:findMk}.

\emph{Computation of $\eta_1$}. Since $\OO$ admits a birationally rigid cover and $\Sigma_1 \simeq A_1$, the hypotheses of Proposition \ref{prop:identification2} are satisfied. Thus
$$\eta_1 (\tau_1(1)) = c_1\omega_1(1), \qquad c_1 = 2|\Pic(\OO_{M_1})||\Pic(\OO)|^{-1}.$$
Note that $\Pic(\OO) \simeq \pi_1(\OO)_{\mathrm{ab}} \simeq S_2 \times \ZZ_2$ and $\Pic(\OO_{M_1}) \simeq \pi_1(\OO_{M_1})_{\mathrm{ab}} \simeq \ZZ_2 \times \ZZ_2$. Thus, $c_1=2(4)/4=2$. By definition, $\tau_1(1)$ is the dominant generator of the free abelian group $\fX(M_1)$, i.e. $\tau_1(1) = \varpi_1$. It follows that
$$\eta_1(\varpi_1) = 2 \omega_1(1)$$
\end{example}

\begin{example}\label{ex:exceptional2}
Let $G$ be the (unique) simple group of type $E_8$, and let $\mathbb{O}=E_8(b_6)$. The universal cover $\widehat{\OO}$ of $\OO$ is birationally rigid by Proposition \ref{prop:almost all A1}. By \cite[Sec 13]{fuetal2015}, $\overline{\OO}$ is normal in codimension 2. The codimension 2 orbits $\OO_k \subset \overline{\OO}$ and the singularities $\Sigma_k$ are indicated below
\vspace{3mm}
\begin{center}
    \begin{tabular}{|c|c|c|}\hline
        $k$ & $\Sigma_k$ & $\mathbb{O}_k$\\ \hline
        $1$ & $A_1$ & $A_7$ \\ \hline
        $2$ & $A_2$ & $E_6(a_1)+A_1$ \\ \hline
    \end{tabular}
\end{center}
\vspace{3mm}
In this example, we will compute the Levi subgroups $M_1,M_2 \subset G$ adapted to the codimension 2 leaves $\fL_1,\fL_2 \subset \Spec(\CC[\OO])=X$ as well as the identifications $\eta_1,\eta_2$.

Label the simple roots as follows

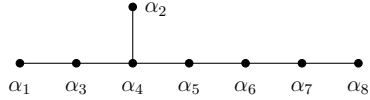
\begin{center}
$$\dynkin[labels={\alpha_1,\alpha_2,\alpha_3,\alpha_4,\alpha_5,\alpha_6,\alpha_7,\alpha_8},edge
length=.75cm] E8$$
\end{center}
and write $\varpi_1,...,\varpi_8$ for the fundamental weights. Arguing as in Example \ref{ex:exceptional0}, we see that $\OO=\mathrm{Bind}^G_L \{0\}$, where $L$ is the standard Levi subgroup (of Lie type $A_3+A_2+A_1$) corresponding to the simple roots $\{\alpha_1,\alpha_2,\alpha_3,\alpha_5,\alpha_6,\alpha_7\}$. 

\emph{Computation of $M_1$}. Let $M$ denote the standard Levi subgroup (of Lie type $E_7$) corresponding to the simple roots $\alpha_1,\alpha_2,\alpha_3,\alpha_4,\alpha_5,\alpha_6,\alpha_7$. Note that $L \subset M_1$ and by \cite[Sec 4]{deGraafElashvili}
$$\OO_M = \Ind^M_L \{0\} = A_4+A_2$$
Using Lemma \ref{lem:findMk}, we will show that $M=M_1$. By \cite[Sec 13]{fuetal2015}, $\overline{\OO}_M$ is normal in codimension 2. There is one boundary orbit $\OO_{M,2} \subset \overline{\OO}_M$, namely $\OO_{M,2}=A_4+A_1$, and the corresponding singularity is of type $A_2$. By \cite[Sec 13]{deGraaf2013}, $\OO_{M,2} = \Ind^M_{A_4+A_1} \{0\}$. Hence, $\Ind^G_M A_4+A_1 = \Ind^G_{A_4+A_1} \{0\} = E_6(a_1)+A_1$. So $M=M_1$ by Lemma \ref{lem:findMk}.

\emph{Computation of $\eta_1$}. Since $\OO$ admits a birationally rigid cover and $\Sigma_1 \simeq A_1$, the hypotheses of Proposition \ref{prop:identification2} are satisfied. Thus
$$\eta_1(\tau_1(1)) = c_1\omega_1(1), \qquad c_1=2|\Pic(\OO_{M_1})||\Pic(\OO)|^{-1}.$$
Note that $\Pic(\OO) \simeq \pi_1(\OO)_{\mathrm{ab}} \simeq (S_3)_{\mathrm{ab}} \simeq \ZZ_2$ and $\Pic(\OO_{M_1}) \simeq \pi_1(\OO_{M_1})_{\mathrm{ab}} \simeq 1$. Thus, $c_1= 2(1)/2=1$. By definition, $\tau_1(1)$ is the dominant generator of the free abelian group $\fX(M)$, i.e. $\tau_1(1)=\varpi_8$. Hence
\begin{equation}\label{eq:eta1final}\eta_1(\varpi_8) = \omega_1(1).\end{equation}

\emph{Computation of $M_2$}. For our given choice of $L$, $M_2$ is non-standard. So the computation of $M_2$ is more involved than it was for $M_1$. Let $L'$ denote the standard Levi subgroup (of Lie type $A_3+A_2+A_1$) corresponding to the simple roots $\alpha_1,\alpha_2,\alpha_4,\alpha_5,\alpha_7,\alpha_8$. In $E_8$, a Levi subgroup is determined up to conjugacy by the isomorphism type of its Lie algebra. So $L$ and $L'$ are $G$-conjugate. Furthermore, {\tt atlas} provides us with a Weyl group element $w \in W$ conjugating $L'$ onto $L$. 

Let $K'$ denote the standard Levi subgroup (of Lie type $D_5+A_2$) corresponding to the simple roots $\alpha_1,\alpha_2,\alpha_3,\alpha_4,\alpha_5,\alpha_7,\alpha_8$ and let $K=wK'$. Note that $L' \subset K'$ and therefore $L \subset K$. Using Lemma \ref{lem:findMk}, we will show that $K=M_2$. Note that
$$\OO_{K} = \Ind^{D_5+A_2}_{D_3+A_1+A_2} \{0\} = \OO_{(3^2,1^4)} \times \{0\}$$
(here $D_3+A_1+A_2$ embeds in $D_5+A_2$, with $D_3+A_1 \subset D_5$). By \cite[Table 1]{Kraft-Procesi}, $\overline{\OO}_{K}$ is normal in codimension 2. There is one boundary orbit $\OO_{K,1} \subset \overline{\OO}_{K}$, namely $\OO_{K,1} = \OO_{(3,2^2,1)}$, and the corresponding singularity is of type $A_1$. By \cite[Sec 4]{deGraafElashvili}, $\Ind^G_{K} \OO_{(3,2^2,1^3)} = A_7$. So $K=M_2$ by Lemma \ref{lem:findMk}. 

\emph{Computation of $\eta_2$}. Since $\dim(\fX(\fm_2))=1$ and $\dim(\fh_2^*)=2$, $\pi_1(\fL_2)$ acts on $\fh_2^*$ by the nontrivial diagram automorphism. In particular, $\fP_2^X = (\fh_2^*)^{\pi_1(\fL_2)}$ is spanned by the element $\omega_1(2)+\omega_2(2)$. Since the monodromy action is nontrivial, neither Proposition \ref{prop:identification1} nor Proposition \ref{prop:identification2} can be straightforwardly applied. We will adapt the argument of Proposition \ref{prop:identification1} to show that
$$\eta_2(\tau_1(2)) = \omega_1(2)+\omega_2(2).$$
Note that $X_{M_2}$ has a single codimension 2 leaf $\fL \subset X_{M_2}$, and the corresponding singularity $\Sigma \subset X$ is of type $A_1$. Let $Z_{M_2}$ denote the cotangent bundle of the partial flag variety for $M_2$ corresponding to $L$---it is a $\QQ$-terminalization (in fact, a symplectic resolution) of $X_{M_2}$. Furthermore, there is a natural identification $\fP_2^X \simeq H^2(\mathfrak{S}_k,\CC)^{\pi_1(\fL)}$ which takes $\omega_1(2) + \omega_2(2) \in \fP_2^X$ to $\omega_1+\omega_2 \in H^2(\mathfrak{S}_k,\CC)^{\pi_1(\fL)}$. Thus, it suffices to show that the identification $\eta: \fX(\mathfrak{l}) \simeq H^2(\mathfrak{S}_k,\CC)^{\pi_1(\fL)}$ induced by the restriction map $\Pic(Z_{M_2}) \to \Pic(\mathfrak{S}_k)$ takes $\tau_1(2)$ to $\omega_1(2)+\omega_2(2)$. By Proposition \ref{prop:descriptionofpic}, we have $\Pic(Z_{M_2})\simeq \ZZ$ (with $\overline{\pi}^*\mathcal{L}_{G/Q}(\tau_1(2)) \leftrightarrow 1$). Let $E$ denote the exceptional divisor of $\overline{\rho}: Z_{M_2} \to X_{M_2}$. Since $X_{M_2}$ has only one singularity (of type $A_2$) and the monodromy action is nontrivial, $E$ is irreducible, and the mapping $1 \mapsto [E]$ induces a short exact sequence as in Lemma \ref{lem:2ses}
$$\ZZ\to \Cl(Z_{M_2}) \to \Pic(\mathbb{O}_{M_2})\to 0.$$
Since $\Pic(\mathbb{O}_{M_2})\simeq \ZZ/2\ZZ$, and $[E]$ is anti-dominant, we have  $[E]= \mathrm{div}(\mathcal{L}_{Z_{M_2}}(-2\tau_1(2)))$.
   
Note that $E_{\Sigma} := E \cap \mathfrak{S}_k$ is the union of the two irreducible components of the exceptional divisor of $\mathfrak{S}_k \to \Sigma$. For simple roots $\alpha_1, \alpha_2 \in \mathfrak{h}_{\mathfrak{sl}(3)}^* \simeq H^2(\mathfrak{S}_k,\CC)$, we have $\langle [E_\Sigma], \alpha_i \rangle=-2+1=-1$, and therefore $[E_\Sigma]=-\omega_1(2)-\omega_2(2)$. Arguing as in the proof of \cref{prop:identification1}, we deduce that  $\eta_2(\tau_1(2)) = \omega_1(2)+\omega_2(2)$, as asserted.

By definition, $\tau_1(2)$ is the dominant generator of the free abelian group $\fX(M_2)$. Thus,in the notation of the previous step, $\tau_1(2) = w\varpi_6$. Using {\tt atlas} we compute $\tau_1(2) = \varpi_4-3\varpi_8$. Hence
\begin{equation}\label{eq:eta2final}\eta_2(\varpi_4-3\varpi_8) = \omega_1(2)+\omega_2(2).\end{equation}
\end{example}

\section{A proof of the Dixmier conjecture}\label{subsec:Dixmier}

In this section, we will prove the following result, resolving a conjecture of Vogan (see \cite[Conj 2.3]{Vogan1990}).

\begin{theorem}\label{thm:Dixmier}
Suppose $G$ is semisimple.
The following are true:
\begin{enumerate}
\item 
There is an injective correspondence
$$\mathrm{Dix}: \{G\text{-eqvt covers of coadjoint orbits}\} \hookrightarrow \{\text{Dixmier algebras for }G\}$$
such that for any cover $\widetilde{\OO}^1$, there is an isomorphism of $G$-representations
$$\mathrm{Dix}(\widetilde{\mathbb{O}}^1) \simeq_G \CC[\widetilde{\OO}^1].$$
\item The image of $\mathrm{Dix}$ is the set of Dixmier algebras which arise as filtered quantizations of $\CC[\widetilde{\OO}]$, where $\widetilde{\OO}$ runs over the set of $G$-equivariant nilpotent covers. 
\item 
If $\widetilde{\OO}$ is a \emph{nilpotent} cover, then $\mathrm{Dix}(\widetilde{\OO}) = \cA_0^{\widetilde{X}}$ (where the latter is regarded as a Dixmier algebra via the uniquely defined co-moment map $\Phi_0^{\widetilde{X}}: U(\fg) \to \cA_0^{\widetilde{X}}$).
\end{enumerate}
\end{theorem}
    
Note that Theorem \ref{thm:Dixmier} generalizes  \cite[Thm 5.3(1)]{Losev4}, and is proved by similar methods. 

\begin{proof}
We will establish bijections between the following four sets. 
\begin{itemize}
    \item[(a)] The set of $G$-equivariant covers of coadjoint $G$-orbits (up to isomorphisms of $G$-equivariant covers).
    \item[(b)] The set of $G$-equivariant filtered Poisson deformations of algebras $\CC[\widetilde{\OO}]$ (up to $G$-equivariant isomorphisms of filtered Poisson algebras).
    \item[(c)] The set of $G$-equivariant filtered quantizations of algebras
    $\CC[\widetilde{\OO}]$ (up to $G$-equivariant isomorphisms of filtered algebras).
    \item[(d)] The set of $G$-equivariant filtered quantizations of algebras $\CC[\widetilde{\OO}]$ (up to Dixmier algebra isomorphisms). 
\end{itemize}

A bijection between (b) and (c) is constructed in \cite[Rmk 3.24]{Losev4}: for $\widetilde{X}:=\operatorname{Spec}(\CC[\widetilde{\OO}])$ both sets are classified by $\fP^{\widetilde{X}}/N_G(L,\widetilde{\OO}_L)$ and the bijection preserves the parameter. In particular, $\cA^{\widetilde{X}}_0$ corresponds to the trivial deformation $\CC[\widetilde{\OO}]$. 

A bijection between (c) and (d) is as follows. 
Every $G$-equivariant filtered quantization $\cA$ of $\CC[\widetilde{\OO}]$ admits a unique quantum co-moment map $\Phi: U(\fg) \to \cA$, see Lemma \ref{lem:co-momentunique}. Note that $\Phi$ turns $\cA$ into a Dixmier algebra for $G$. This defines a (tautologically) surjective map from (c) to (d). To prove that this map is injective we will show that the filtration on $\cA$ is uniquely recovered up to a $G$-equivariant algebra isomorphism. Consider the functor $\bullet_{\dagger}: \HC^G(U(\fg)) \to \HC^R(\cW)$ associated to the orbit $\OO$ covered by $\widetilde{\OO}$, see Section 
\ref{subsec:W}. By Lemma \ref{lem:Adagger}, $\cA_{\dagger}$ is concentrated in degree $0$ and there is a $G$-equivariant filtered algebra isomorphism $\cA\xrightarrow{\sim} (\cA_\dagger)^\dagger$. So the filtration is recovered uniquely. This shows that the map from (c) to (d) is injective.

It remains to establish a bijection between (a) and (b). To go from (b) to (a), we take a filtered Poisson deformation $\cA^0$ of $\CC[\widetilde{\OO}]$ and form its spectrum $\widetilde{X}^1$. The group $G$ acts on $\widetilde{X}^1$ in a Hamiltonian fashion. The moment map is unique because $G$ is semisimple. The $G$-variety $\widetilde{X}^1$ contains a (unique) open $G$-orbit since $\widetilde{X}$ does, and the moment map realizes this open orbit as a cover of a coadjoint orbit. We map $\cA^0$ to this open orbit.  Taking the disjoint union over all possible $\widetilde{\OO}$, we get a map from (b) to (a).

To prove that this map is a bijection 
we will need a parameterization of $G$-equivariant covers of coadjoint orbits. Generalizing \cite[Definition 1.2]{Losev4}, we define the notion of a birationally minimal induction datum for a cover of a coadjoint orbit. This is a triple $(L,\widetilde{\OO}_L,\xi)$ consisting of a Levi subgroup $L \subset G$, a {birationally rigid} $L$-equivariant nilpotent cover $\widetilde{\OO}_L$, and an element $\xi \in \fX(\fl)$. To each such triple, we can associate a $G$-equivariant cover of a coadjoint orbit, namely the open $G$-orbit in the variety $G\times^P(\{\xi\}\times \widetilde{X}_L\times \mathfrak{p}^\perp)$, where $\widetilde{X}_L=\Spec(\CC[\widetilde{\OO}_L])$. This defines a bijection between the set of all $G$-equivariant covers of coadjoint orbits and the set of all birationally minimal induction data up to conjugation by $G$. Indeed, it suffices to show that two birationally minimal induction data for $\OO^1$ are $G$-conjugate. This easily reduces to the case when $\widetilde{\OO}^1$ is a nilpotent cover. In this case, the conjugacy follows from (iii) of Proposition \ref{prop:birationalinduction}. 

Now let $\widetilde{\OO}^1$ be a $G$-equivariant cover with birationally minimal induction datum
$(L,\widetilde{\OO}_L,\xi)$. Let $\widetilde{\OO} = \mathrm{Bind}^G_L \widetilde{\OO}_L$, and let $\mathsf{Cov}(\widetilde{\OO})$ denote the set of covers $\widetilde{\OO}^1$
which give rise to $\widetilde{\OO}$ via this construction. We will show that 
\begin{itemize}
\item[(*)]
the map from (b) to (a) defines a bijection between the set of $G$-equivariant Poisson deformations of $\CC[\widetilde{\Orb}]$ and
$\mathsf{Cov}(\widetilde{\OO})$.
\end{itemize}

Recall, Theorem \ref{thm:defsofsymplectic}, that the set $\mathrm{PDef}(\CC[\widetilde{\OO}])$ of Poisson deformations of $\CC[\widetilde{\OO}]$ up to isomorphisms of Poisson deformations is in natural bijection with $\fP^{\widetilde{X}}/W^{\widetilde{X}}$. By Remark \ref{rmk:equiv_automor}, every such deformation has a natural $G$-action. By the same remark, the group of graded $G$-equivariant automorphisms of $\CC[\widetilde{\OO}]$ acts on $\fP^{\widetilde{X}}/W^{\widetilde{X}}$.
And the filtered Poisson deformations up to $G$-equivariant isomorphism of filtered Poisson algebras are classified by orbits of this group action. 

By Proposition \ref{prop:namikawacovers}, there is a natural isomorphism $\fP^{\widetilde{X}} \simeq \fX(\fl)$. The Galois group $\Pi$ of the cover $\widetilde{\OO}\rightarrow \OO$ is  naturally identified with the 
group of $G$-equivariant graded Poisson automorphisms of $\CC[\widetilde{\OO}]$, and, by (ii) of Proposition \ref{prop:namikawacovers}, the Weyl group $W^{\widetilde{X}}$ is the kernel of the natural epimorphism $N_G(L,\widetilde{\OO}_L)/L\twoheadrightarrow \Pi$. So the filtered Poisson deformations up to filtered $G$-equivariant Poisson isomorphism are parameterized by $\fX(\fl)/N_G(L,\widetilde{\OO}_L)$. The deformation corresponding to the orbit of $\lambda \in \fX(\fl)$ is recovered as $\CC[G\times^P(\{\xi\}\times \widetilde{X}_L\times \fp^\perp)]$, this follows from (iii) of Proposition \ref{prop:namikawacovers}. The open $G$-orbit therein coincides with the cover birationally induced from $(\fl,\widetilde{\OO}_L,\xi)$ because the natural morphism
$G\times^P(\{\xi\}\times \widetilde{X}_L\times \fp^\perp)\rightarrow \Spec(\CC[G\times^P(\{\xi\}\times \widetilde{X}_L\times \fp^\perp)])$
is an isomorphism over the open orbit.  To finish the proof of (*) note that the $G$-conjugacy classes of birationally minimal induction data with $\fl,\widetilde{\OO}_L$ fixed are in bijection with $\fX(\fl)/N_G(L,\widetilde{\OO}_L)$.
Under the identifications of both sets with $\fX(\fl)/N_G(L,\widetilde{\OO}_L)$ the map from (b) to (a) is the identity. 

This proves claim  (2) of the theorem. Claim (3) follows because the bijection between (b) and (c) sends $\CC[\widetilde{\OO}]$ to $\cA^{\widetilde{X}}_0$. For (1), it remains to establish an isomorphism $\operatorname{Dix}(\widetilde{\OO}^1)\simeq_G \CC[\widetilde{\OO}^1]$. 
Since $\mathrm{Dix}(\widetilde{\OO}^1)$ is a filtered quantization of $\CC[\widetilde{\OO}]$, there is an isomorphism of $G$-representations
\begin{equation}\label{eq:G_mod_iso_Dixmier}
\operatorname{Dix}(\widetilde{\OO}^1)\simeq_G \CC[\widetilde{\OO}].
\end{equation}
The variety 
$\Spec(\CC[G\times^P(\{\xi\}\times \widetilde{X}_L\times \fp^\perp)])$ is normal and, since it admits a finite $G$-equivariant map to a coadjoint orbit closure, contains finitely many orbits, all of even dimension. Since $\widetilde{\OO}^1$ is an open orbit, we see that $\CC[G\times^P(\{\xi\}\times \widetilde{X}_L\times \fp^\perp)]=\CC[\widetilde{\OO}^1]$. The left hand side is a filtered Poisson deformation of $\CC[\widetilde{\OO}]$. 
So $\CC[\widetilde{\OO}^1]\simeq_G \CC[\widetilde{\OO}]$. Combined with (\ref{eq:G_mod_iso_Dixmier}) this completes the proof of (1).
\end{proof}

\chapter{Infinitesimal characters of unipotent ideals}\label{sec:centralchars}

Let $G$ be a connected reductive algebraic group and let $\widetilde{\OO}$ be a $G$-equivariant nilpotent cover. Choose $\lambda \in \overline{\fP}^{\widetilde{X}}$ and consider the corresponding ideal $I_{\lambda}(\widetilde{\OO})$. By Proposition \ref{prop:propsofIbeta}(i), $I_{\lambda}(\widetilde{\OO})$ is primitive. Define
$$\gamma_{\lambda}(\widetilde{\OO}) := \text{infinitesimal character of } I_{\lambda}(\widetilde{\OO}) \in \fh^*/W.$$
In this chapter, we will outline a general strategy for computing the \emph{unipotent infinitesimal character} $\gamma_0(\widetilde{\OO})$\index{infinitesimal character!unipotent}. It is based on the computation of the map $\eta$ in 
Chapter \ref{sec:nilpquant}. This strategy is implemented in Section \ref{subsec:centralcharclassical} in the case of linear classical groups. In Sections \ref{subsec:centralcharspin} and \ref{subsec:centralcharexceptional}, we discuss the application to spin and exceptional groups.

\section{Unipotent infinitesimal characters}\label{subsec:unipotentcentralchars}

Let $\widetilde{\OO}$ be a $G$-equivariant nilpotent cover. Choose a Levi subgroup $L \subset G$ and a birationally rigid $L$-equivariant nilpotent cover $\widetilde{\mathbb{O}}_L$ such that $\widetilde{\mathbb{O}} = \mathrm{Bind}^G_L \widetilde{\mathbb{O}}_L$. Let $M \subset G$ be a Levi subgroup containing $L$ and let $\widetilde{\mathbb{O}}_M = \mathrm{Bind}^M_L \widetilde{\mathbb{O}}_L$. Suppose $\gamma$ is an infinitesimal character for $\fm$, regarded as an element of $\fh^*/W_M$. Since $W_M$ is a subgroup of $W$, $\gamma$ defines an element in $\fh^*/W$, and hence an infinitesimal character for $\fg$. 

\begin{prop}\label{prop:Ibetainduced}
Let $\beta \in \fX(\fl)$. There is an equality in $\fh^*/W$
$$\gamma_{\beta}(\widetilde{\mathbb{O}}) = \gamma_{\beta}(\widetilde{\mathbb{O}}_M).$$
\end{prop}

\begin{proof}
As in the proof of Lemma \ref{lem:induct_transitivity} we can assume that $G$ is semisimple. Choose a parabolic subgroup $Q\subset G$ with Levi subgroup $M$ and a Borel subgroup $B \subset G$ contained in $Q$. The Borel subgroup $B$ determines a system of positive roots in $M$ and in $G$. Let $\rho(\fm)$ denote the half-sum of the positive roots in $\fm$ and choose an anti-dominant representative for the infinitesimal character $\gamma_{\beta}(\widetilde{\OO}_M)$. Consider the sheaf $\mathfrak{D}_{Q/B}^{\gamma_{\beta}(\widetilde{\OO}_M)+\rho(\fm)}$ of twisted differential operators on $Q/B$ with twist $\gamma_{\beta}(\widetilde{\OO}_M)+\rho(\fm)$. Taking global sections, we get an algebra isomorphism
$$ \Gamma(Q/B, \mathfrak{D}_{Q/B}^{\gamma_{\beta}(\widetilde{\OO}_M)+\rho(\fm)}) \simeq U(\fm)/(I_{\beta}(\widetilde{\OO}_M) \cap \fz(\fm)) =: \mathcal{U}_M,$$
see \cite[Theorem 11.2.2]{HTT}. Furthermore, $\Phi_{\beta}^{\widetilde{X}_M}: U(\fm) \to \cA_{\beta}^{\widetilde{X}_M}$ factors through an algebra homomorphism $\Phi_M: \mathcal{U}_M \to \cA_{\beta}^{\widetilde{X}_M}$. We can repeat the construction of \cref{subsec:inductionquantizations} to produce an algebra $\Ind_M^G \mathcal{U}_M$. 
By Proposition \ref{prop:quantizationparaminduction} applied to $\widetilde{X}_M:=\mathcal{N}_M$, we have
$$\Ind_M^G \mathcal{U}_M \simeq \Gamma(G/B, \mathfrak{D}_{G/B}^{\gamma_\beta(\widetilde{\OO}_M)+\rho(\fg)}),$$
Let $\Psi: U(\fg) \to \Ind^G_M \mathcal{U}_M$ denote the quantum co-moment map. Note that $\ker{\Psi}$ has infinitesimal character $\gamma_{\beta}(\widetilde{\OO}_M)$, see again \cite[Theorem 11.2.2]{HTT}. Furthermore, $\Phi_M$ gives rise to an algebra homomorphism $\Phi:\Ind^G_M\mathcal{U}_M \to \Ind^G_M\cA_{\beta}^{\widetilde{X}_M} \simeq \cA_{\beta}^{\widetilde{X}}$, where the isomorphism $\Ind^G_M\cA_{\beta}^{\widetilde{X}_M} \simeq \cA_{\beta}^{\widetilde{X}}$ follows from \cref{prop:quantizationparaminduction}. By Lemma \ref{lem:co-momentunique}, $\Phi_{\beta}^{\widetilde{X}} = \Phi \circ \Psi$. In particular, $\ker{\Phi_{\beta}^{\widetilde{X}}} = I_{\beta}(\widetilde{\OO})$ has infinitesimal character $\gamma_{\beta}(\widetilde{\OO}_M)$.
\end{proof}

\begin{rmk}
There is a classical notion of parabolic induction for 2-sided ideals in enveloping algebras  (see e.g. \cite[Chapter 5]{Dixmier}). It is not difficult to show that $I_{\beta}(\widetilde{\mathbb{O}})$ coincides with induced ideal $\Ind^G_M I_{\beta}(\widetilde{\mathbb{O}}_M)$. Since parabolic induction preserves infinitesimal character, this is strictly stronger than Proposition \ref{prop:Ibetainduced}. Since we will not use this fact, we omit the proof. 
\end{rmk}
Now choose a Levi subgroup $K \subset L$ and a birationally rigid $K$-orbit such that $\mathbb{O}_L = \Ind^L_K \mathbb{O}_K$ (the induction need not be birational). Let $\widecheck{\mathbb{O}}_L = \mathrm{Bind}^L_K \mathbb{O}_K$ and $\widecheck{X}_L = \Spec(\CC[\widecheck{\mathbb{O}}_L])$. Note that the universal $L$-equivariant cover $\widehat{\OO}_L$ of $\OO_L$ is a Galois cover of $\widecheck{\OO}_L$. By Corollary \ref{cor:criterionbirigid}, $\widetilde{\OO}_L$ is 2-leafless, and thus $\widehat{\OO}_L$ is 2-leafless by Lemma \ref{lem:cover2leafless}. Consider the barycenter parameter $\epsilon \in \fP^{\widecheck{X}_L}$ (for the Galois cover $\widehat{X}_L\twoheadrightarrow \widecheck{X}_L$) defined in Example \ref{ex:barycentersymplectic} and let $\delta = \eta^{-1}(\epsilon) \in \fX([\fl,\fl]\cap\mathfrak{k})$. We can view $\delta$ as an element of $\fX(\mathfrak{k})$ via the direct sum decomposition $\fX(\mathfrak{k})=\fX(\mathfrak{l})\oplus \fX([\fl,\fl]\cap\mathfrak{k})$.

Suppose $\gamma$ is an infinitesimal character for $\mathfrak{k}$, regarded as an element of $\fh^*/W_K$. Then the infinitesimal character $\gamma+\delta \in \fh^*/W$ is defined in the usual way: choose a representative for $\gamma$ in $\fh^*$ and form the sum $\gamma+\delta \in \fh^*$. Since $W_K \subset W$ and $\delta$ is a fixed point for the $W_K$-action, $\gamma+\delta$ is well-defined modulo $W$, and independent of the choice of representative. The unipotent infinitesimal character $\gamma_0(\widetilde{\mathbb{O}})$ can now be computed as follows.

\begin{prop}\label{prop:unipotentcentralchars}
There is an equality in $\fh^*/W$
\begin{equation}\label{eq:desired_equality}\gamma_0(\widetilde{\mathbb{O}}) = \gamma_0(\mathbb{O}_K) + \delta.\end{equation}
\end{prop}

\begin{proof}
By Proposition \ref{prop:Ibetainduced}, there is an equality
\begin{equation}\label{eq:centralchar1}
    \gamma_0(\widetilde{\mathbb{O}}) = \gamma_0(\widetilde{\mathbb{O}}_L).
\end{equation}
Note that $\widetilde{\OO}_L$ is equivalent to $\widehat{\OO}_L$, since both covers are 2-leafless. By Proposition \ref{prop:equivalencerelation}, 
$I_0(\widehat{\mathbb{O}}_L) = I_0(\widetilde{\mathbb{O}}_L)$, and therefore $\gamma_0(\widehat{\OO}_L) = \gamma_0(\widetilde{\OO}_L)$. Let $\Pi = \pi_1^L(\widecheck{\mathbb{O}}_L)$. By Proposition \ref{prop:parameterofinvariantssymplectic}, there is an isomorphism of Hamiltonian quantizations $\cA_{\epsilon}^{\widecheck{X}_L} \simeq (\cA_0^{\widehat{X}_L})^{\Pi}$. Hence,
$$I_{\delta}(\widecheck{\mathbb{O}}_L) = I_0(\widehat{\mathbb{O}}_L) = I_0(\widetilde{\mathbb{O}}_L),$$
and therefore
\begin{equation}\label{eq:centralchar2}
\gamma_{\delta}(\widecheck{\mathbb{O}}_L) = \gamma_0(\widehat{\mathbb{O}}_L) = \gamma_0(\widetilde{\mathbb{O}}_L).\end{equation}
By Proposition \ref{prop:Ibetainduced}
\begin{equation}\label{eq:centralchar3}
\gamma_{\delta}(\widecheck{\mathbb{O}}_L) = \gamma_{\delta}(\mathbb{O}_K) = \gamma_0(\mathbb{O}_K) + \delta.\end{equation}
Combining (\ref{eq:centralchar1}), (\ref{eq:centralchar2}) and (\ref{eq:centralchar3}) gives 
(\ref{eq:desired_equality}).
\end{proof}
Using Proposition \ref{prop:unipotentcentralchars}, we can reduce the calculation of $\gamma_0(\widetilde{\mathbb{O}})$ to the case of \emph{rigid} nilpotent orbits. For such orbits, we will appeal to the following proposition.

\begin{prop}\label{prop:mult1ideals}
Suppose $\fg$ is simple. Let $\mathbb{O}$ be a birationally rigid nilpotent orbit, which is not one of the following:
\begin{center}
    \begin{tabular}{|c|c|c|c|c|c|c|} \hline
        $\fg$ & $G_2$ & $F_4$ & $E_7$ & $E_8$ & $E_8$ & $E_8$ \\ \hline
        $\mathbb{O}$ & $\widetilde{A}_1$ & $\widetilde{A}_2+A_1$ & $(A_1+A_3)'$ & $A_3+A_1$ & $D_5(a_1)+A_2$ & $A_5+A_1$ \\ \hline
    \end{tabular}
\end{center}
Then there is a unique primitive ideal $I \subset U(\fg)$ such that
\begin{itemize}
    \item[(i)] $V(I) = \overline{\mathbb{O}}$.
    \item[(ii)] $m_{\overline{\mathbb{O}}}(U(\fg)/I) = 1$.
\end{itemize}
Furthermore, $I=I_0(\mathbb{O})$.
\end{prop}

\begin{proof}
By Proposition \ref{prop:propsofIbeta}, $I_0(\mathbb{O})$ satisfies (i) and (ii) above. The uniqueness claim is \cite[Lemma 4.6]{Losev5}.
%
\end{proof}

\begin{rmk}\label{Rem:6_orbits}
There is a geometric characterization of the six orbits appearing in Proposition \ref{prop:mult1ideals}: they are exactly the rigid orbits $\OO$ such that $\overline{\OO}$ is not normal in codimension $2$. This can be deduced from the incidence tables in \cite{fuetal2015}. 
\end{rmk}

\section{Unipotent infinitesimal characters: linear classical groups}\label{subsec:centralcharclassical}

Let $G$ be linear classical. In this section, we will compute the unipotent infinitesimal characters attached to $G$-equivariant nilpotent covers. We will follow the recipe outlined in Section \ref{subsec:unipotentcentralchars}.

\subsection{Birationally rigid orbits}\label{subsec:centralcharrigid}

First, assume $\mathbb{O}$ is a birationally rigid orbit. We will compute $\gamma_0(\mathbb{O})$ using Proposition \ref{prop:mult1ideals} and the results of McGovern (\cite{McGovern1994}). It will be convenient to define some combinatorial operations.

\begin{definition}\label{def:rhoplus}
Suppose $q$ is a partition of $n$. Define $\rho(q) \in (\frac{1}{2}\ZZ)^{n}$ by appending
$$(\frac{q_i-1}{2}, \frac{q_i-3}{2}, ..., \frac{3-q_i}{2}, \frac{1-q_i}{2})$$
for each $i \geq 1$.

Define $\rho^+(q) \in (\frac{1}{2}\ZZ)^{\floor{\frac{n}{2}}}$ by appending the \emph{positive} elements of the sequence
$$(\frac{q_i-1}{2}, \frac{q_i-3}{2}, ..., \frac{3-q_i}{2}, \frac{1-q_i}{2})$$
for each $i \geq 1$, and then $0$'s as needed so that $|\rho^+(q)| = \floor{\frac{n}{2}}$.
\end{definition}

For example, if $q = (4,3^2,1)$, then
$$\rho(q) = (\frac{3}{2},\frac{1}{2},-\frac{1}{2},-\frac{3}{2},1,0,-1,1,0,-1,0) \quad \text{and} \quad \rho^+(q) = (\frac{3}{2},\frac{1}{2},1,1,0).$$
\begin{definition}\label{def:fBfC}
Suppose $q$ is a partition. Define $f_B(q)$ as follows: for every odd $i$ with $q_i \geq q_{i+1} +2$, move one box down from $q_i$ to $q_{i+1}$. Define $f_C(q)$ as follows: for every even $i$ with $q_i \geq q_{i+1} + 2$, move one box from $q_i$ to $q_{i+1}$ and then add a single box to $q_1$. If $q = (0)$ is the trivial partition, define $f_C(q)= (1)$.
\end{definition}
For example, if $q = (7^2,4^3,2,1^2)$, then
$$f_B(q) = (7^2,4^2,3^2,1^2) \quad \text{and} \quad f_C(q) = (8,6,5,4^2,2,1^2).$$

\begin{prop}\label{prop:centralcharacterbirigidorbit}
Suppose $\mathbb{O} \subset \fg^*$ is a birationally rigid nilpotent orbit corresponding to a partition $p$. Then $I_0(\mathbb{O})$ is a maximal ideal, and $\gamma_0(\mathbb{O})$ is as follows

\begin{itemize}
    \item[(i)] If $\fg = \mathfrak{sl}(n)$, then $\mathbb{O} = \{0\}$ and $\gamma_0(\mathbb{O}) = \rho$.
    \item[(ii)] If $\fg = \mathfrak{so}(n)$, then
    $$\gamma_0(\mathbb{O}) = \rho^+(f_B(p^t)).$$
    \item[(iii)] If $\fg = \mathfrak{sp}(2n)$, then 
    $$\gamma_0(\mathbb{O}) = \rho^+(f_C(p^t)).$$
\end{itemize}
\end{prop}

\begin{proof}
We will show that for each infinitesimal character $\gamma$ listed above, the maximal ideal $I_{\mathrm{max}}(\gamma)$ satisfies the properties
\begin{equation}\label{eqn:avandmult}
V(I_{\mathrm{max}}(\gamma)) = \overline{\mathbb{O}}, \qquad m_{\overline{\mathbb{O}}}(U(\fg)/I_{\mathrm{max}}(\gamma)) = 1.
\end{equation}
Then by Proposition \ref{prop:mult1ideals}, we have $I_{\mathrm{max}}(\gamma) = I_0(\mathbb{O})$.

If $\gamma=\rho$, then $I_{\mathrm{max}}(\gamma)$ is the augmentation ideal, and (\ref{eqn:avandmult}) is trivial. Thus, we can assume $\fg = \mathfrak{so}(n)$ or $\mathfrak{sp}(2n)$. 

In \cite[Sec 4]{McGovern1994}, McGovern attaches a finite collection of infinitesimal characters $Q(\mathbb{O})$ to every nilpotent orbit $\mathbb{O} \subset \fg^*$. For each $\gamma \in Q(\mathbb{O})$, he proves that 
$$V(I_{\mathrm{max}}(\gamma)) = \overline{\mathbb{O}}$$
(see \cite[Thm 5.1]{McGovern1994}). For the infinitesimal characters above, it is easy to check that $\gamma \in Q(\mathbb{O})$. For a proof of this containment (and a description of $Q(\mathbb{O})$), we refer the reader to Proposition \ref{lem:maximalitybirigidtypeC}. McGovern also gives a formula for $m_{\overline{\mathbb{O}}}(U(\fg)/I_{\mathrm{max}}(\gamma))$, see \cite[Thm 5.14]{McGovern1994}. For the infinitesimal characters above, it is easy to compute that $m_{\overline{\mathbb{O}}}(U(\fg)/I_{\mathrm{max}}(\gamma)) = 1$. We leave the straightforward details to the reader.
\end{proof}

\begin{rmk}\label{rmk:Barbasch}
From the formulas above, it is clear that all unipotent infinitesimal characters attached to rigid orbits in classical types are `extremal' in the sense of \cite[Def 10.3]{Barbasch1989}. 
\end{rmk}

\begin{example}\label{ex:centralcharbirigdorbit}
Let $G=\mathrm{SO}(9)$. Using Proposition \ref{prop:birigidorbitclassical}, we see that there are 4 birationally rigid orbits, corresponding to the partitions
$$(3,2^2,1^2), \ (2^4,1), \ (2^2,1^5), \ (1^9).$$
The corresponding infinitesimal characters are given in the table below
\vspace{5mm}
\begin{center}
\begin{tabular}{|c|c|}\hline
$\mathbb{O}$ & $\gamma_0(\mathbb{O})$\\ \hline
$(3,2^2,1^2)$ & $(\frac{3}{2},\frac{3}{2},\frac{1}{2},\frac{1}{2})$ \\ \hline
$(2^4,1)$ & $(2,\frac{3}{2},1,\frac{1}{2})$\\ \hline
$(2^2,1^5)$ & $(\frac{5}{2},\frac{3}{2},1,\frac{1}{2})$ \\ \hline
$(1^9)$ & $(\frac{7}{2},\frac{5}{2},\frac{3}{2},\frac{1}{2})$ \\ \hline
\end{tabular}
\end{center}
\vspace{5mm}
\end{example}
If $\mathbb{O}$ is both birationally rigid and special, then $\gamma_0(\mathbb{O})$ has an even simpler description. 

\begin{lemma}\label{lem:specialbirigid}
Assume $\fg$ is classical. Suppose $\mathbb{O}$ is birationally rigid and special. Then there is a unique nilpotent orbit $\mathbb{O}^{\vee} \subset (\fg^{\vee})^*$ satisfying $\mathsf{D}(\mathbb{O}^{\vee}) = \mathbb{O}$ (defined by $\mathbb{O}^{\vee} = \mathsf{D}(\mathbb{O})$, for example). Furthermore, $I_0(\mathbb{O}) = I_{\mathrm{max}}(\frac{1}{2}h_{\mathbb{O}^{\vee}})$.
\end{lemma}

\begin{proof}
If $\mathbb{O}^{\vee} \subset \cN^{\vee}$ satisfies $\mathsf{D}(\mathbb{O}^{\vee}) = \mathbb{O}$, then by \cite[Cor. 5.19]{McGovern1994}
$$m_{\overline{\mathbb{O}}}(U(\fg)/I_{\mathrm{max}}(\frac{1}{2}h_{\mathbb{O}^{\vee}})) = 1.$$
Thus, $I_0(\mathbb{O}) = I_{\mathrm{max}}(\frac{1}{2}h_{\mathbb{O}^{\vee}})$ by \cref{prop:mult1ideals}. By the Dynkin classification of nilpotent orbits, $\mathbb{O}^{\vee}$ is uniquely determined by $\frac{1}{2}h_{\mathbb{O}^{\vee}}$. Hence, $\mathbb{O}^{\vee}$ is the unique nilpotent $G^{\vee}$-orbit satisfying $\mathsf{D}(\mathbb{O}^{\vee}) = \mathbb{O}$.
\end{proof}

\subsection{Birationally rigid covers}

Next, assume $\widetilde{\mathbb{O}}$ is a birationally rigid cover of $\OO$. As usual, let $X=\Spec(\CC[\OO])$. Choose a Levi subgroup $K \subset G$ and a birationally rigid nilpotent $K$-orbit $\mathbb{O}_K$ such that $\mathbb{O} = \mathrm{Bind}^G_K \mathbb{O}_K$. The possible orbits $\mathbb{O}$ and pairs $(K,\mathbb{O}_K)$ are described in Proposition \ref{prop:nocodim2leavesLM}. Consider the barycenter parameter $\epsilon \in \fP^X$ and let $\delta = \eta^{-1}(\epsilon) \in \fX(\mathfrak{k})$. By Proposition \ref{prop:unipotentcentralchars}, there is an equality in $\fh^*/W$
\begin{equation}\label{eq:gamma0formulaL}\gamma_0(\widetilde{\mathbb{O}}) = \gamma_0(\mathbb{O}_K) + \delta.\end{equation}
Since $\mathbb{O}_K$ is birationally rigid, $\gamma_0(\mathbb{O}_K)$ can be computed using Proposition \ref{prop:centralcharacterbirigidorbit}. The element $\delta$ was computed in Corollaries \ref{cor:lambdaA} and \ref{cor:lambdaBCD}. Below, we will give a formula for $\gamma_0(\widetilde{\mathbb{O}})$ involving only the partition for $\mathbb{O}$ (this is possible since all birationally rigid covers $\widetilde{\mathbb{O}} \to \mathbb{O}$ are equivalent, see Proposition \ref{prop:equivalencerelation}).

\begin{definition}\label{def:xy}
Suppose $q$ is a partition. Let $x(q) \subset q$ be the subpartition consisting of multiplicity-1 parts and let $y(q) \subset q$ be the subpartition consisting of multiplicity-2 parts. 

Suppose $y$ is a partition such that every part has multiplicity 2. Define a partition $g(y)$ (of the same size as $y$) by replacing every pair $(y_i,y_i)$ with $(y_i+1, y_i-1)$. 
\end{definition}

For example, if $q = (6^2,5,3^3,1)$, then $x(q) = (5,1)$ and $y(q) = (6^2)$. If $y = (5^2,4^2,1^2)$, then $g(y) = (6,5^2,3,2)$.

\begin{prop}\label{prop:centralcharacterbirigidcover}
Let $p$ be the partition corresponding to $\mathbb{O}$. Form the partitions $x=x(p^t)$ and $y=y(p^t)$. Then $\gamma_0(\widetilde{\mathbb{O}})$ is given by the following formulas

\begin{itemize}
    \item[(i)] If $G = \mathrm{SL}(n)$, then $p=(d^m)$ and
$$\gamma_0(\widetilde{\mathbb{O}}) = (\frac{n-1}{2d}, \frac{n-3}{2d}, ..., \frac{1-n}{2d}) = \frac{\rho}{d}.$$
    \item[(ii)] If $G = \mathrm{SO}(2n+1)$ or $\mathrm{SO}(2n)$, then
    $$\gamma_0(\widetilde{\mathbb{O}}) = \rho^+(g(y) \cup f_B(x)),$$
    where $\rho^+$ and $f_B$ are as in Definitions \ref{def:rhoplus} and \ref{def:fBfC}, respectively. 
     \item[(iii)] If $G= \mathrm{Sp}(2n)$, then
    $$\gamma_0(\widetilde{\mathbb{O}}) = \rho^+(g(y) \cup f_C(x)).$$
\end{itemize}
\end{prop}

\begin{proof}
First, suppose $G=\mathrm{SL}(n)$. By Proposition \ref{prop:unipotentcentralchars}
\begin{align*}
\gamma_0(\widetilde{\mathbb{O}}) &= \gamma_0(\{0\}) + \delta.
\end{align*}
By Corollary \ref{cor:lambdaA} we have 
\begin{align*}
\delta=(\underbrace{\frac{d-1}{2d},\ldots, \frac{d-1}{2d}}_{m}, \underbrace{\frac{d-3}{2d},\ldots, \frac{d-3}{2d}}_{m},\ldots, \underbrace{\frac{1-d}{2d},\ldots, \frac{1-d}{2d}}_{m}).
\end{align*}
Combining, we get 
\begin{align*}
    \gamma_0(\widetilde{\mathbb{O}}) &= (\underbrace{\frac{m-1}{2},\ldots, \frac{1-m}{2}}_{m},\ldots, \underbrace{\frac{m-1}{2},\ldots, \frac{1-m}{2}}_{m})\\
                                      &+(\underbrace{\frac{d-1}{2d},\ldots, \frac{d-1}{2d}}_{m}, \underbrace{\frac{d-3}{2d},\ldots, \frac{d-3}{2d}}_{m},\ldots, \underbrace{\frac{1-d}{2d},\ldots, \frac{1-d}{2d}}_{m})
\end{align*}
which coincides up to permutation with $\rho/d$.

{\ Next, suppose $G=G_{\epsilon}(n)$ is the simple classical group of type $\epsilon \in \{B,C,D\}$ and rank $n$. We will assume for the moment that $(G,\OO) \neq (SO(8m+4),\OO_{(4^{2m},31)})$ (this case requires separate treatment and will be handled below). By Proposition \ref{prop:nocodim2leavesLM}, we have
$$K =  \prod_{k \in S_2(p)}\mathrm{GL}(k) \times G_{\epsilon}(n-|S_2(p)|), \qquad \mathbb{O}_K = \{0\} \times ... \times \{0\} \times \mathbb{O}_{p\#S_2(p)}.$$
By the conditions on $p$ in Proposition \ref{prop:nocodim2leaves}, we have $p^t = x \cup y$. Note that $i \in y$ if and only if $i \in S_2(p)$. Hence, $y = S_2(p) \cup S_2(p)$. On the other hand, $i \in x$ if and only if $i$ is a column in $p$ of multiplicity 1. Hence, $x = (p \# S_2(p))^t$. Now by Proposition \ref{prop:centralcharacterbirigidorbit}
$$\gamma_0(\mathbb{O}_K) = (\rho(S_2(p)),\rho^+(f_{\epsilon}(x))),$$
where $\epsilon\in \{B,C,D\}$ and $f_D$ is the same as $f_B$. 
By Corollary \ref{cor:lambdaBCD}
$$\delta = \frac{1}{2}(\underbrace{1,1,...,1}_{|S_2(p)|},\underbrace{0,...,0}_{n-|S_2(p)|}).$$
Hence
\begin{equation}\label{eq:gamma01}
    \gamma_0(\widetilde{\mathbb{O}}) = \gamma_0(\mathbb{O}_K) + \delta =(\rho(S_2(p)),\rho^+(f_{\epsilon}(x)))+ (\underbrace{\frac{1}{2}, ..., \frac{1}{2}}_{|S_2(p)|},\underbrace{0,...,0}_{n-|S_2(p)|}).
\end{equation}
Let $W'$ denote the finite subgroup of $GL_n(\RR)$ generated by arbitrary permutations and sign changes. Note that $W'=W$ for $\epsilon \in \{B,C\}$ and $[W:W']=2$ for $\epsilon=D$. Since, $\rho(S_2(p))+(\frac{1}{2},...,\frac{1}{2})$ agrees with $\rho^+(g(y))$ up to permutations and sign changes, the right hand side of Equation \ref{eq:gamma01} is $W'$-conjugate to $(\rho^+(g(y)),\rho^+(f_{\epsilon}(x))) = \rho^+(g(y) \cup f_{\epsilon}(x))$. If $\epsilon=D$, we make the following additional observation. By construction, all columns in $p\#S_2(p)$ are distinct. Hence, all rows in $x=(p\#S_2(p))^t$ are distinct. If $x_1$ is odd, then $f_B(x)_1 = x_1$. Otherwise, $f_B(x)_1 = x_1-1$. In either case, $f_B(x)$ contains at least one odd part, and therefore $\rho^+(g(y) \cup f_B(x))$ contains at least one entry equal to $0$. Hence $W(\rho^+(g(y) \cup f_B(x))) = W'(\rho^+(g(y) \cup f_B(x)))$.

Finally, suppose $G=\mathrm{SO}(8m+4)$ and $\mathbb{O}=\mathbb{O}_{(4^{2m}31)}$. By Proposition \ref{prop:nocodim2leavesLM}, we have
$$K = \mathrm{GL}(2m+1) \times \mathrm{GL}(2m+1), \qquad \mathbb{O}_K = \{0\} \times \{0\}.$$
By Corollary \ref{cor:lambdaBCD}, we have
$$\delta = \frac{1}{2}(\underbrace{1,...,1}_{2m+1},\underbrace{0,...,0}_{2m+1}).$$
Hence
$$\gamma_0(\widetilde{\mathbb{O}}) = \gamma_0(\mathbb{O}_K) + \delta = (\rho(2m+1),\rho(2m+1)) + \frac{1}{2}(\underbrace{1,...,1}_{2m+1},\underbrace{0,...,0}_{2m+1}).$$
This weight contains (exactly) one entry equal to $0$, and is therefore $W$-conjugate to $\rho^+(2m+2, 2m+1,2m+1,2m)$. On the other hand, $p^t = (2m+2,2m+1,2m+1,2m)$ and therefore $f_B(x) \cup g(y) = (2m+2,2m+1,2m+1,2m)$. Thus, $\gamma_0(\widetilde{\mathbb{O}}) = \rho^+(f_B(x) \cup g(y))$, as asserted. This completes the proof.

}

\end{proof}

\begin{example}
Let $G = \mathrm{SO}(9)$. The orbits admitting birationally rigid covers are given by the partitions
$$(5,3,1), (3,2^2,1^2), (3,1^6), (2^4,1), (2^2,1^5), (1^9).$$
Only $(5,3,1)$ and $(3,1^6)$ are not birationally rigid (see Example \ref{ex:centralcharbirigdorbit}). The orbit $(5,3,1)$ admits two (equivalent) birationally rigid covers, of degrees 2 and 4. The orbit $(3,1^6)$ admits a single birationally rigid cover, of degree 2. The corresponding infinitesimal characters are $(1,\frac{1}{2},\frac{1}{2},0)$ and $(\frac{5}{2},\frac{3}{2},\frac{1}{2},\frac{1}{2})$, respectively.
\end{example}

\subsection{Arbitrary covers}\label{subsec:arbitrarycovers}

Finally, let $\widetilde{\mathbb{O}}$ be an arbitrary $G$-equivariant nilpotent cover. Choose a Levi subgroup $L \subset G$ and a birationally rigid $L$-equivariant nilpotent cover $\widetilde{\mathbb{O}}_L$ such that $\widetilde{\mathbb{O}} = \mathrm{Bind}^G_L \widetilde{\mathbb{O}}_L$. 

If $G = \mathrm{SL}(n)$, then $L$ is of the form 
$$L = \mathrm{S}(\mathrm{GL}(a_1) \times \mathrm{GL}(a_2) \times ... \times \mathrm{GL}(a_t)),$$
and
$$\widetilde{\mathbb{O}}_L = \widetilde{\mathbb{O}}_{\mathrm{SL}(a_1)} \times \widetilde{\mathbb{O}}_{\mathrm{SL}(a_2)} \times ... \times \widetilde{\mathbb{O}}_{\mathrm{SL}(a_t)},$$
where each $\widetilde{\mathbb{O}}_{\mathrm{SL}(a_i)}$ is a birationally rigid $\mathrm{SL}(a_i)$-equivariant nilpotent cover. Hence by Proposition \ref{prop:Ibetainduced}
$$\gamma_0(\widetilde{\mathbb{O}}) = (\gamma_0(\widetilde{\mathbb{O}}_{\mathrm{SL}(a_1)}), \gamma_0(\widetilde{\mathbb{O}}_{\mathrm{SL}(a_2)}), ..., \gamma_0(\widetilde{\mathbb{O}}_{\mathrm{SL}(a_t)}))).$$
Each $\gamma_0(\widetilde{\mathbb{O}}_{\mathrm{SL}(a_i)})$ can be computed using Proposition \ref{prop:centralcharacterbirigidcover}.

If $G=\mathrm{Sp}(2n)$ or $\mathrm{SO}(n)$, then $L$ is of the form
$$L = \mathrm{GL}(a_1) \times ... \times \mathrm{GL}(a_t) \times  G(m),$$
and
$$\widetilde{\mathbb{O}}_L =  \{0\} \times ... \{0\} \times \widetilde{\mathbb{O}}_{G(m)},$$
where $\widetilde{\mathbb{O}}_{G(m)}$ is a birationally rigid $G(m)$-equivariant nilpotent cover. If $m \neq 0$, then
$$\gamma_0(\widetilde{\mathbb{O}}) = ( \rho(a_1),\rho(a_2),...,\rho(a_t),\gamma_0(\widetilde{\mathbb{O}}_{G(m)})).$$
Otherwise, $\mathbb{O}$ is very even. If $G=\operatorname{Sp}(2n)$ or $\operatorname{SO}(2n+1)$, then 
$$\gamma_0(\widetilde{\mathbb{O}}) = (\rho(a_1),...,\rho(a_t))$$
If $G=\operatorname{SO}(2n)$, then 
$$\gamma_0(\widetilde{\mathbb{O}}) = (\rho(a_1),...,\rho(a_t)) \quad  \text{or} \quad \gamma_0(\widetilde{\mathbb{O}}) = (\rho(a_1),...,\rho(a_t)) \text{ with last entry negated},$$
depending on the numeral attached to $\mathbb{O}$.

\section{Unipotent infinitesimal characters: spin groups}\label{subsec:centralcharspin}

In this section, we modify the approach outlined in Section \ref{subsec:centralcharclassical} to compute the unipotent infinitesimal characters for spin-equivariant covers. Let $G=\mathrm{Spin}(n)$, and let $\widetilde{\OO}$ be a ($G$-equivariant) nilpotent cover, maximal in its equivalence class. By Proposition \ref{prop:Ibetainduced}, the computation of $\gamma_0(\widetilde{\mathbb{O}})$ reduces to the case of birationally rigid covers. Thus, we can (and will) assume that $\widetilde{\mathbb{O}}$ is birationally rigid. In particular, since $\widetilde{\OO}$ is assumed to be maximal, $\widetilde{\OO}$ is the universal ($G$-equivariant) cover of $\OO$, see Lemma \ref{lem:cover2leafless}.

If $\widetilde{\mathbb{O}}$ is $\mathrm{SO}(n)$-equivariant, then $\gamma_0(\widetilde{\mathbb{O}})$ can be computed using the results of Section \ref{subsec:centralcharclassical}. So we can also assume that $\widetilde{\mathbb{O}}$ is \emph{not} $\mathrm{SO}(n)$-equivariant. Thus, $\OO$ is is as described in Proposition \ref{prop:nocodim2leavesspin} and $\widetilde{\OO}$ is a 2-fold cover of the universal $\mathrm{SO}(n)$-equivariant cover $\widehat{\OO}$ of $\OO$. Let
$$L =  \prod_{i \in S_4(p)} \mathrm{GL}(i) \times \mathrm{SO}(n-2|S_4(p)|), \qquad \widehat{\mathbb{O}}_L =  \{0\} \times ... \times \{0\} \times \widehat{\mathbb{O}}_{p\#S_4(p)}.$$
where $\widehat{\OO}_{p\# S_4(p)}$ denotes the universal $\mathrm{SO}(n-2|S_4(p)|)$-equivariant cover of $\OO_{p\#S_4(p)}$. By Proposition \ref{prop:nocodim2leavesLM}, $\widehat{\mathbb{O}}_L$ is birationally rigid and $\mathrm{Bind}^{\mathrm{SO}(n)}_L \widehat{\mathbb{O}}_L = \widehat{\OO}$. Arguing as in the proof of Proposition \ref{prop:unipotentcentralchars},
\begin{align*}
\gamma_0(\widetilde{\mathbb{O}}) &= \gamma_0(\widehat{\OO}_L) + \delta,
\end{align*}
where $\delta \in \fX(\fl)$ is the element corresponding to the barycenter parameter $\epsilon \in \fX^{\widehat{X}}$. The element $\delta$ was computed in Corollary \ref{cor:lambdaspin}:
\begin{align*}
\delta= (\underbrace{\frac{1}{4},...,\frac{1}{4}}_{|S_4(p)|},0,...,0).
\end{align*}
Combining, we get 
\begin{align*}
    \gamma_0(\widetilde{\mathbb{O}}) &= (\rho(S_4(p)),\gamma_0(\widehat{\mathbb{O}}_{p\#S_4(p)})) + (\underbrace{\frac{1}{4},...,\frac{1}{4}}_{|S_4(p)|},0,...,0).
\end{align*}
The infinitesimal character $\gamma_0(\widehat{\mathbb{O}}_{p\#S_4(p)})$, and hence $\gamma_0(\widetilde{\mathbb{O}})$, can be computed using the results of Section \ref{subsec:centralcharclassical}.

\begin{example}
Let $G=\mathrm{Spin}(6)$, $\mathbb{O}=\mathbb{O}_{(5,1)}$, and let $\widetilde{\mathbb{O}}$ be the universal cover of $\OO$. By Proposition \ref{prop:nocodim2leavesspin}, $\widetilde{\mathbb{O}}$ is birationally rigid. We have
$$S_4(p) = (1), \qquad p\#S_4(p) = (3,1).$$
Hence,
$$L = \mathrm{GL}(1) \times  \mathrm{SO}(4), \qquad \widehat{\mathbb{O}}_L =  \{0\} \times \widehat{\mathbb{O}}_{(3,1)}.$$
By Proposition \ref{prop:centralcharacterbirigidcover}, we have $\gamma_0(\widehat{\mathbb{O}}_{(3,1)}) = (\frac{1}{2},0)$. So
$$\gamma_0(\widetilde{\mathbb{O}}) = \gamma_0(\widehat{\OO}_L) + \delta = (0,\frac{1}{2},0) + (\frac{1}{4},0,0) = (\frac{1}{4},\frac{1}{2},0),$$
which is $W$-conjugate to $\rho/4$. Note that there is an isomorphism $\mathrm{Spin}(6) \simeq \mathrm{SL}(4)$. And indeed, the unipotent infinitesimal character attached to the universal cover of the principal orbit in $\mathfrak{sl}(4)$ is $\rho/4$, see Proposition \ref{prop:centralcharacterbirigidcover}(i).
\end{example}

\section{Unipotent infinitesimal characters: exceptional groups}\label{subsec:centralcharexceptional}

Let $G$ be a simple exceptional group, and let $\widetilde{\OO}$ be a $G$-equivariant nilpotent cover. Following the recipe outlined in Section \ref{subsec:unipotentcentralchars} one can, in principle, compute the infinitesimal character of $I_0(\widetilde{\OO})$. Here, we will sketch the procedure and give some examples. In \cite{MBM}, we carry out this procedure in all cases to produce a complete list of unipotent infinitesimal characters in exceptional types.

\begin{enumerate}
    \item Find a Levi subgroup $L \subset G$ and a birationally rigid cover $\widetilde{\mathbb{O}}_L$ such that $\widetilde{\mathbb{O}} = \mathrm{Bind}^G_L \widetilde{\mathbb{O}}_L$.
    \item Find a Levi subgroup $K \subset L$ and a rigid orbit $\mathbb{O}_K$ such that $\mathbb{O}_L = \Ind^L_K \mathbb{O}_K$ (for this, one can use the results of \cite{deGraafElashvili}).
    \item Compute $\gamma_0(\mathbb{O}_K)$. For the six `bad' orbits in Proposition \ref{prop:mult1ideals}, this computation is carried out in Appendix \ref{SS_m_sing}. In all other cases, there is a unique multiplicity 1 primitive ideal with $V(I) = \overline{\mathbb{O}}_K$, computed by Premet in \cite[Secs 3-5]{Premet2013}, which must coincide with $I_0(\OO_K)$, see Proposition \ref{prop:mult1ideals}. For the reader's convenience, a full list of unipotent infinitesimal characters attached to rigid orbits in exceptional types is provided in Table \ref{table:exceptional} after Appendix \ref{SS_m_sing}. 

    \item Compute the barycenter parameter $\delta \in \fX(\mathfrak{k})$. Except in 4 cases, all dimension 2 singularities in $\Spec(\CC[\mathbb{O}_L])$ are of type $A_1$ (see Proposition \ref{prop:almost all A1}), and hence $\delta$ can be computed as in Example \ref{ex:exceptional1}.
    \item Use Proposition \ref{prop:unipotentcentralchars} and steps (3) and (4) to compute $\gamma_0(\widetilde{\mathbb{O}}) = \gamma_0(\mathbb{O}_K) + \delta$.
\end{enumerate}

\begin{example}\label{ex: G2covers}
Let $G=G_2$ and let $\mathbb{O} = G_2(a_1)$. In this example, we will compute the unipotent infinitesimal characters associated to every cover of $\OO$.

Since $\pi_1(\mathbb{O}) \simeq S_3$, $\mathbb{O}$ admits 3 non-isomorphic covers (apart from $\mathbb{O}$ itself): 
\begin{itemize}
    \item $\widetilde{\mathbb{O}}_2$: Galois 2-fold cover, $\pi_1(\widetilde{\mathbb{O}}_2) \simeq \ZZ_3$.
    \item $\widetilde{\mathbb{O}}_3$: non-Galois 3-fold cover, $\pi_1(\widetilde{\mathbb{O}}_3) \simeq \ZZ_2$.
    \item $\widehat{\mathbb{O}}$: Galois 6-fold (universal) cover, $\pi_1(\widehat{\mathbb{O}}) \simeq 1$.
\end{itemize}
Recall from Example \ref{ex:exceptional0} that there is a unique codimension 2 leaf $\fL_1 \subset \Spec(\CC[\OO])$ and the corresponding singularity is of type $A_1$. This $A_1$ singularity cannot be smoothened under the 3-fold cover $\Spec(\CC[\widetilde{\OO}_3]) \to \Spec(\CC[\OO])$. Hence, $\widetilde{\OO}_3$ is birationally induced. Apart from $T$ and $G$, there are two non-conjugate Levi subgroups of $G$, corresponding to $\alpha_1$ and $\alpha_2$. Denote them by $L_1$ and $L_2$, respectively. As explained in Example \ref{ex:exceptional0}, $\OO=\mathrm{Bind}^G_{L_2} \{0\}$. So $\widetilde{\OO}_3 = \mathrm{Bind}^G_{L_1} \{0\}$. The remaining two covers of $\mathbb{O}$, namely $\widetilde{\mathbb{O}}_2$ and $\widehat{\mathbb{O}}$, are birationally rigid, and hence equivalent. 

We first compute the infinitesimal characters $\gamma_0(\OO)$ and $\gamma_0(\widetilde{\OO}_3)$. By Proposition \ref{prop:Ibetainduced}, we have
\begin{align*}
\gamma_0(\OO) &= \rho(\mathfrak{l}_2) = \frac{1}{2}\alpha_2\\
\gamma_0(\widetilde{\OO}_3) &= \rho(\mathfrak{l}_1) = \frac{1}{2}\alpha_1
\end{align*}
So that our answers are consistent with {\tt atlas} conventions, we will write these expressions in the basis $\{\varpi_1,\varpi_2\}$ of fundamental weights (see, e.g., \cite[Plates I-IX]{Bourbaki46} for the change-of-basis matrix)
\begin{align*}
\gamma_0(\OO) &= -\frac{3}{2}\varpi_1+ \varpi_2\\
\gamma_0(\widetilde{\OO}_3) &= \varpi_1 -\frac{1}{2} \varpi_2\end{align*}
In these coordinates, it is clear that neither weight is dominant. We use {\tt atlas} to compute their dominant representatives
$$\gamma_0(\OO) = \frac{1}{2}\varpi_2, \qquad \gamma_0(\widetilde{\OO}_3) = \frac{1}{2}\varpi_1$$
For $\widetilde{\OO}=\widetilde{\OO}_2$ or $\widetilde{\OO}=\widehat{\OO}$, we will compute $\gamma_0(\widetilde{\OO})$ using Proposition \ref{prop:unipotentcentralchars}. By equation (\ref{eq:G22}) of Example \ref{ex:exceptional0},
$$\delta = \eta_1^{-1}(\epsilon_1) = \eta_1^{-1}(\frac{1}{2}\omega_1(1)) = \frac{1}{2}\varpi_1. $$
Now by Proposition \ref{prop:unipotentcentralchars}
$$\gamma_0(\widetilde{\mathbb{O}}) = \gamma_0(\OO_{L_2}) + \delta = \frac{1}{2}\alpha_2 + \frac{1}{2}\varpi_1 = \varpi_1.$$
We record all of these computations in the table below.
\vspace{3mm}
\begin{center}
\begin{tabular}{|c|c|c|c|c|c|} \hline
 $\widetilde{\mathbb{O}}$ & $(L,\widetilde{\mathbb{O}}_L)$ & $(K,\mathbb{O}_K)$ & $\gamma_0(\mathbb{O}_K)$ & $\delta$ & $\gamma_0(\widetilde{\mathbb{O}})$\\ \hline
 
  $\mathbb{O}$  & $(L_2,\{0\})$ & $(L_2,\{0\})$ & $\frac{1}{2}\alpha_2$ & $0$ & $\frac{1}{2}\varpi_2$\\ \hline
  
 $\widetilde{\mathbb{O}}_2$ & $(G_2, \widetilde{\mathbb{O}}_2)$ & $(L_2,\{0\})$ & $\frac{1}{2}\alpha_2$ & $\frac{1}{2}\varpi_1$ & $\varpi_1$\\ \hline
 
 $\widetilde{\mathbb{O}}_3$ & $(L_1,\{0\})$ & $(L_1,\{0\})$ & $\frac{1}{2}\alpha_1$ & $0$ & $\frac{1}{2}\varpi_1$ \\ \hline
 
 $\widehat{\mathbb{O}}$ & $(G_2, \widetilde{\mathbb{O}}_2)$ & $(L_2,\{0\})$ & $\frac{1}{2}\alpha_2$ & $\frac{1}{2}\varpi_1$ & $\varpi_1$\\ \hline
\end{tabular}
\end{center}
\vspace{3mm}
We note that all three infinitesimal characters $\gamma_0(\widetilde{\mathbb{O}})$ are special unipotent (cf. Definition \ref{def:spec_unipotent}). This is not an accident (see Chapter \ref{sec:duality} for a general statement and explanation).
\end{example}

If $\mathbb{O}_L$ is one of the 4 exceptions listed in Proposition \ref{prop:almost all A1}, then Proposition \ref{prop:identification2} is not, on its own, sufficient for computing the barycenter parameter $\delta \in \fX(\mathfrak{k})$. In these cases, one can apply (a modification of) Proposition \ref{prop:identification1}, as illustrated below.

\begin{example}
Let $G=E_8$, $\mathbb{O}=E_8(b_6)$, and let $\widehat{\mathbb{O}} \to \mathbb{O}$ be the universal (6-fold) $G$-equivariant cover. We will compute the unipotent infinitesimal character $\gamma_0(\widehat{\OO})$. 

Label the simple roots as follows
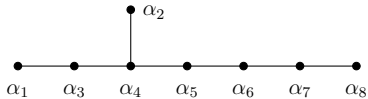
\begin{center}
$$\dynkin[labels={\alpha_1,\alpha_2,\alpha_3,\alpha_4,\alpha_5,\alpha_6,\alpha_7,\alpha_8},edge
length=.75cm] E8$$
\end{center}
and write $\varpi_1,...,\varpi_8$ for the corresponding fundamental weights. By Proposition \ref{prop:almost all A1}, $\widehat{\mathbb{O}}$ is birationally rigid. Hence, $(L,\widetilde{\mathbb{O}}_L) = (G,\widehat{\mathbb{O}})$. On the other hand, $\mathbb{O}$ is birationally induced from $(K,\mathbb{O}_K) = (A_3+A_2+A_1,\{0\})$. Choose $K$ corresponding to the simple roots $\{\alpha_1,\alpha_3,\alpha_2,\alpha_5,\alpha_6,\alpha_7\}$. Then
\begin{equation}\label{eq:E81}
\gamma_0(\mathbb{O}_K) = \rho(\mathfrak{k}) = \alpha_1+ \frac{1}{2}\alpha_2 + \alpha_3 + \frac{3}{2}\alpha_5 + 2\alpha_6+\frac{3}{2}\alpha_7.
\end{equation}
Converting once again into fundamental weight coordinates, we get
$$\gamma_0(\OO_K) = \varpi_1 + \varpi_2 +\varpi_3 - 3\varpi_4 + \varpi_5 +\varpi_6 + \varpi_7 - \frac{3}{2}\varpi_8.$$
By equations (\ref{eq:eta1final}) and (\ref{eq:eta2final}) of Example \ref{ex:exceptional2}
$$\eta_1^{-1}(\epsilon_1) = \frac{1}{2}\eta_1^{-1}(\omega_1(1)) = \frac{1}{2}\varpi_8, \qquad \eta_2^{-1}(\epsilon_2) = \frac{1}{3}\eta_2^{-1}(\omega_1(2)+\omega_2(2)) = \frac{1}{3}\varpi_4 - \varpi_8.$$
Now by Proposition \ref{prop:unipotentcentralchars}
$$\gamma_0(\widehat{\OO}) = \gamma_0(\OO_K) + \delta = \gamma_0(\OO_K) + \eta_1^{-1}(\epsilon_1) + \eta_2^{-1}(\epsilon_2)= \varpi_1 + \varpi_2 +\varpi_3 - \frac{8}{3}\varpi_4 + \varpi_5 +\varpi_6 + \varpi_7 - 2\varpi_8.$$
Conjugating by $W$ to obtain a dominant representative
$$\gamma_0(\widehat{\OO}) = \frac{1}{3}\varpi_1 + \frac{1}{3}\varpi_2 + \frac{1}{3}\varpi_6 + \frac{1}{3}\varpi_8$$
As a sanity check, we conclude by showing that the unipotent ideal $I_0(\widehat{\OO})$ is maximal. We will do so using Proposition \ref{prop:maximalitycriterion}. Fix the notation of Section \ref{subsec:maximality}. Using {\tt atlas} we compute reductive subgroups $L^{\vee}_{\gamma}$ and $L^{\vee}_{\gamma,0}$ of $G^{\vee}$. We see that $L^{\vee}_{\gamma}$ is of type $E_6+A_2$ and $L^{\vee}_{\gamma,0}$ is of type $A_3+A_1$, embedded in the $E_6$ factor of $L^{\vee}_{\gamma}$. Thus, $\OO_{\gamma}^{\vee} = \Ind^{E_6+A_2}_{A_1+A_3} = D_4(a_1) \times (3)$, see \cite[Sec 3]{deGraafElashvili} and $\OO_{\gamma} = \mathsf{D}(\OO_{\gamma}^{\vee}) = D_4(a_1) \times \{0\}$. Now it is easy to see that
$$\codim(\OO,\cN) = 20 = \codim(\OO_{\gamma},\cN_{L_{\gamma}})$$
Thus $I_0(\widehat{\OO})$ is maximal by Proposition \ref{prop:maximalitycriterion}.
\end{example}

\section{Maximality of unipotent ideals}\label{subsec:maximality}

In Appendix \ref{sec:maximality}, we prove the following result.

\begin{theorem}\label{thm:maximality1}
Suppose that $G$ is linear classical, and let $\widetilde{\OO}$ be a $G$-equivariant nilpotent cover. Then $I_0(\widetilde{\OO})$ is a maximal ideal.
\end{theorem}
The idea of the proof is as follows. Let $\gamma:=\gamma_0(\widetilde{\OO})$ and define subgroups
$$L_{\gamma,0}^{\vee} \subset L_{\gamma}^{\vee} \subset G^{\vee}$$
as in Section \ref{subsec:assvarmax}. Let $$\OO_{\gamma} := \mathsf{D}(\Ind^{L^{\vee}_{\gamma}}_{L^{\vee}_{\gamma,0}} \{0\}),$$ a nilpotent orbit for $L_{\gamma}$, the dual group of $L_{\gamma}^{\vee}$. By Proposition \ref{prop:maximalitycriterion}, $I_0(\widetilde{\OO})$ is maximal if and only if
$$\codim(\OO,\cN) = \codim(\OO_{\gamma},\cN_{L_{\gamma}})$$
Checking this condition amounts to some standard (but tedious) casework, which is carried out in Sections \ref{subsec:maximalityA}-\ref{subsec:maximalityD}. 

In \cite{MBM}, the second and third-named authors extend Theorem \ref{thm:maximality1} to arbitrary $G$.

\chapter{A refinement of BVLS duality}\label{sec:duality}

Let $G$ be a connected reductive algebraic group. Recall from Section \ref{subsec:BVduality} there is an order-reversing map
$$\mathsf{D}: \{\mathbb{O}^{\vee} \subset \cN^{\vee}\}\to \{\mathbb{O} \subset \cN\}, \qquad \mathsf{D}(\mathbb{O}^{\vee}) := \text{open orbit in } V(I_{\mathrm{max}}(\frac{1}{2}h_{\mathbb{O}^{\vee}})).$$
called BVLS duality. In this section, for each nilpotent $G^{\vee}$-orbit $\mathbb{O}^{\vee}$, we will construct an equivalence class of covers $\widetilde{\mathsf{D}}(\mathbb{O}^{\vee})$ of $\mathsf{D}(\mathbb{O}^{\vee})$ such that
$$I_0(\widetilde{\mathsf{D}}(\mathbb{O}^{\vee}))= I_{\mathrm{max}}(\frac{1}{2}h_{\mathbb{O}^{\vee}}).$$
In particular, this shows that every special unipotent ideal is unipotent.

The assignment $\mathbb{O}^{\vee} \mapsto \widetilde{\mathsf{D}}(\mathbb{O}^{\vee})$ defines an injective map
$$\widetilde{\mathsf{D}}: \{\text{nilpotent orbits for } G^{\vee}\} \hookrightarrow \{\text{equivalence classes of nilpotent covers for } G\}$$
which we call \emph{refined BVLS duality}\index{duality!refined BVLS}. This map is constructed in Sections \ref{subsec:dualdistinguished}-\ref{subsec:refinedBVLS}. Some motivation for our construction is provided in Section \ref{subsec:motivationsymplectic}.

\section{Distinguished orbits and their duals}\label{subsec:dualdistinguished}

\begin{prop}\label{prop:dualtodistinguished}
Let $\mathbb{O}^{\vee}$ be a distinguished nilpotent $G^{\vee}$-orbit. Let $\mathbb{O} = \mathsf{D}(\mathbb{O}^{\vee})$ and $\widehat{\mathbb{O}} \to \mathbb{O}$ the universal $G$-equivariant cover. Then 
\begin{itemize}
    \item[(i)] $\widehat{\mathbb{O}}$ is birationally rigid.
    \item[(ii)] $\frac{1}{2}h_{\mathbb{O}^{\vee}} = \gamma_0(\widehat{\mathbb{O}})$.
\end{itemize}
\end{prop}
Our proof of this proposition is a case-by-case computation, and will occupy the remainder of this section.

\subsection*{Classical types} If $\fg = \mathfrak{sl}(n)$, then Proposition \ref{prop:dualtodistinguished} is trivial---only the principal orbit $\OO^\vee$ is distinguished, its dual orbit is $\OO:=\{0\}$, and
$$\frac{1}{2}h_{\mathbb{O}^{\vee}} = \rho = \gamma_0(\mathbb{O}).$$

Suppose $\fg = \mathfrak{so}(2n+1), \mathfrak{sp}(2n)$ or $\mathfrak{so}(2n)$, and let $q$ be the partition corresponding to $\mathbb{O}^{\vee}$. By Corollary \ref{corollary: distinguished}, $q$ is of the following form

\begin{itemize}
   \item If $\fg = \mathfrak{so}(2n+1)$, then $\fg^{\vee} = \mathfrak{sp}(2n)$, and
   $$q = (2k_1,2k_2,...,2k_m), \qquad k_1 > k_2 >...>k_m >0.$$
    \item If $\fg = \mathfrak{sp}(2n)$, then $\fg^{\vee} = \mathfrak{so}(2n+1)$, and
    $$q = (2k_1-1, 2k_2-1,...,2k_m-1), \qquad k_1 > k_2 >...>k_m >0, \ m  \text{ is odd}.$$

   \item If $\fg = \mathfrak{so}(2n)$, then $\fg^{\vee} = \mathfrak{so}(2n)$, and
   $$q = (2k_1-1, 2k_2 - 1,...,2k_m-1), \qquad k_1 > k_2 >...>k_m >0, \ m  \text{ is even}.$$
\end{itemize}

Let $p$ be the partition corresponding to $\mathbb{O}$. By Proposition \ref{prop:BVduality}, 
\begin{equation}\label{eq:pai}p = ((m+1)^{a_{m+1}}, m^{a_m}, ..., 1^{a_1}),\end{equation}
where the multiplicities $a_1,a_2,...,a_m$ are as follows.

\begin{itemize}
    \item If $\fg = \mathfrak{so}(2n+1)$ and $m$ is odd, then
    $$a_i = \begin{cases} 
             0 & i=m+1\\
             2k_m+1 & i = m \\
             2(k_i-k_{i+1}) +2 & i \text{ odd, }i\neq m \\
             2(k_i-k_{i+1}) -2 & i \text{ even, }i\neq m+1. \\
             
   \end{cases}$$
   
   \item If $\fg = \mathfrak{so}(2n+1)$ and $m$ is even, then

    $$a_i = \begin{cases}
             1 & i = m+1 \\
             2(k_i-k_{i+1}) +2 & i \text{ odd, }i\neq m+1 \\
             2(k_i-k_{i+1}) -2 & i \text{ even}. \\
             
   \end{cases}$$
   
   \item If $\fg = \mathfrak{sp}(2n)$, then
    $$a_i = \begin{cases} 
             0 & i=m+1\\
             2(k_i-k_{i+1})-2 & i \text{ odd} \\
             2(k_i-k_{i+1})+2 & i \text{ even, }i\neq m+1. \\
   \end{cases}$$
  
   \item If $\fg = \mathfrak{so}(2n)$, then
    $$a_i = \begin{cases} 
             0 & i=m+1\\
             2(k_i-k_{i+1}) +2 & i \text{ odd, }i\neq m+1  \\
             2(k_i-k_{i+1}) -2 & i \text{ even}. \\
   \end{cases}$$
\end{itemize}
We leave the (straightforward) verification to the reader. By Corollary \ref{cor:criterionbirigid} {to prove (i)} it is enough to show that 
\begin{itemize}
    \item[(a)] $\widehat{\OO}$ is $2$-leafless;
    \item[(b)] $H^2(\widehat{\OO},\CC)=0$. 
\end{itemize}
For (a) by \cref{lem:cover2leafless} we may assume that $G$ is adjoint group, and let $G^{cl}$ be the corresponding classical group (i.e. $SO(n)$ or $Sp(n)$). One can check using \cite[Cor 6.1.6]{CM} that in each case $\pi_1^G(\OO)\simeq \pi_1^{G^{cl}}(\OO)$. Proposition \ref{prop:nocodim2leaves} implies (a). To prove (b) note that for orbits $\OO$ in classical types, the reductive part $\mathfrak{r}$ of the centralizer of $e \in \OO$ is given in \cite[Thm 6.1.3]{CM}. If $\mathfrak{g} = \mathfrak{so}(2n+1)$ or $\mathfrak{so}(2n)$ (resp. $\mathfrak{sp}(2n)$), then $\mathfrak{r}$ is semisimple if and only if the corresponding partition contains no odd (resp. even) parts of multiplicity 2. For the partitions above, it is easy to check that this condition is satisfied. Thus, $H^2(\widehat{\OO},\CC)=0$ by Lemma \ref{lem:computeH2}. This completes the proof of (i) in classical types. 

We now proceed to proving (ii). First, assume $\fg = \mathfrak{sp}(2n)$. Let $x=x(p^t)$ and $y=y(p^t)$, cf. Definition \ref{def:xy}. Then by Proposition \ref{prop:centralcharacterbirigidorbit}
\begin{equation}\label{eq:gamma0p}\gamma_0(\widehat{\mathbb{O}}) = \rho^+(g(y) \cup f_C(x)),\end{equation}
where $f_C(x)$ and $g(y)$ are the partitions constructed in Definitions \ref{def:fBfC} and \ref{def:xy} and $\rho^+$ is the infinitesimal character of Definition \ref{def:rhoplus}. The weight $\frac{1}{2}h_{\mathbb{O}^{\vee}}$ can be determined from the explicit triples in \cite[Chp 5.2]{CM}. For $\fg=\mathfrak{sp}(2n)$, the conclusion is as follows
\begin{equation}\label{eq:halfhcheck}
\frac{1}{2}h_{\mathbb{O}^{\vee}} = \rho^+(q).\end{equation}
Thus, comparing (\ref{eq:gamma0p}) and (\ref{eq:halfhcheck}), it suffices to show that $g(y) \cup f_C(x) = q$. By (\ref{eq:pai}), we have
$$p^t = (\sum_{i \geq 1} a_i, \sum_{i \geq 2}a_i, ..., \sum_{i\geq m}a_i).$$
Using the formulas above for $a_1,...,a_m$, we can rewrite this in terms of $k_1,...,k_m$
$$p^t = (2k_1-2,2k_2,2k_3-2,2k_4,...,2k_m-2).$$
Let $S = \{i \text{ odd} \mid k_i = k_{i+1}+1\}$ and $T = \{i \text{ odd} \mid k_i \geq k_{i+1}+2 \}$. Then
$$y= y(p^t) = \bigcup_{i \in S} (2k_i-2, 2k_{i+1}), \qquad x = x(p^t) = \bigcup_{i \in T} (2k_i-2,2k_{i+1}).$$
Hence by definition 
$$g(y) = \bigcup_{i \in S} (2k_i-1,2k_{i+1}-1), \qquad f_C(x) = \bigcup_{i \in T} (2k_i-1,2k_{i+1}-1),$$
and therefore
$$g(y) \cup f_C(x)  = \bigcup_{i \in S \cup T} (2k_i-1,2k_{i+1}-1) = (2k_1-1,2k_2-1,...,2k_m-1) = q$$
as desired. For $\fg=\mathfrak{so}(2n+1)$ or $\mathfrak{so}(2n)$, the proof is analogous. 

\subsection*{Exceptional types}

There are a total of 26 distinguished nilpotent orbits in the simple exceptional Lie algebras. For each $\OO^{\vee}$, the dual orbit $\OO=\mathsf{D}(\OO^{\vee})$ is given in \cite[Sec 13.4]{Carter1993}. In the majority of cases, $\OO$ is rigid (see \cite[Sec 4]{deGraafElashvili} for a complete list of rigid orbits in simple exceptional types). In the table below, we list all such $\OO^{\vee}$:

\vspace{3mm}
\begin{center}
    \begin{tabular}{|c|c|c|} \hline
        $\fg$ & $\OO^{\vee}$ & $\OO$\\ \hline
        $G_2$ & $G_2$ & $\{0\}$ \\ \hline
        $F_4$ & $F_4$ & $\{0\}$ \\ \hline
        $F_4$ & $F_4(a_1)$ & $\widetilde{A}_1$ \\ \hline
        $F_4$ & $F_4(a_2)$ & $\widetilde{A}_1+A_1$ \\ \hline
        $E_6$ & $E_6$ & $\{0\}$ \\ \hline
        $E_6$ & $E_6(a_1)$ & $A_1$ \\ \hline
        $E_7$ & $E_7$ & $\{0\}$ \\ \hline
        $E_7$ & $E_7(a_1)$ & $A_1$ \\ \hline
        $E_7$ & $E_7(a_2)$ & $2A_1$ \\ \hline
        $E_7$ & $E_7(a_4)$ & $A_2+2A_1$ \\ \hline
        $E_8$ & $E_8$ & $\{0\}$ \\ \hline
        $E_8$ & $E_8(a_1)$ & $A_1$ \\ \hline
        $E_8$ & $E_8(a_2)$ & $2A_1$ \\ \hline
        $E_8$ & $E_8(a_4)$ & $A_2+A_1$ \\ \hline
        $E_8$ & $E_8(b_4)$ & $A_2+2A_1$ \\ \hline
        $E_8$ & $E_8(a_6)$ & $D_4(a_1)+A_1$ \\ \hline
    \end{tabular}
\end{center}
\vspace{3mm}
In these cases, part (i) of Proposition \ref{prop:dualtodistinguished} is automatic: any $G$-equivariant cover of a rigid nilpotent orbit is birationally rigid. For (ii), we argue as follows. Since $\OO$ is rigid, we have $\gamma_0(\widehat{\OO})=\gamma_0(\OO)$. Then $\gamma_0(\OO)=\frac{1}{2}h_{\OO^{\vee}}$ by \cite[Thm B]{Premet2013}. 

Next, we handle the orbits $\OO^{\vee}$ for which $\OO=\mathsf{D}(\OO^{\vee})$ is induced. For each such $\OO^{\vee}$, we verify (i) and (ii) of Proposition \ref{prop:dualtodistinguished} using essentially the same arguments. We will summarize these arguments now before carrying out the casework. 

For each $\OO^{\vee}$ with $\OO=\mathsf{D}(\OO^{\vee})$ induced we record the following information in tables:

\begin{itemize}
    \item The dominant weight $\frac{1}{2}h_{\mathbb{O}^{\vee}} \in \fh^*$ corresponding to $\mathbb{O}^{\vee}$. This weight is denoted in fundamental weight coordinates, and is computed using the weighted Dynkin diagram of $\mathbb{O}^{\vee}$.
    \item The dual orbit $\mathbb{O} = \mathsf{D}(\mathbb{O}^{\vee})$. 
    \item The reductive part $\mathfrak{r}$ of the centralizer of $e \in \mathbb{O}$, see \cite[Sec 13.1]{Carter1993}. 
    \item The fundamental group $\pi_1(\mathbb{O})$, see \cite[Sec 8.4]{CM}.
    \item The Levi subgroup $L :=L_{\OO}$ corresponding the weighted diagram of $\OO$, see (\ref{eq:JMlevi}).
    \item The integer $m(\OO)$, see the remarks preceding Lemma \ref{lem:inequality2}. This is deducible from \cite[Sec 4]{deGraafElashvili}.
    \item The codimension 2 orbits $\OO' \subset \overline{\OO}$. For each $\OO$, we find that there is unique such orbit by inspection of the tables in \cite[Sec 13]{fuetal2015}.
    \item The fundamental group $\pi_1(\OO')$, see \cite[Sec 8.4]{CM}.
\end{itemize}
In each case, we find that $\OO$ satisfies the following four properties:
\begin{enumerate}
    \item $\OO$ is even.
    \item $X=\Spec(\CC[\OO])$ contains a unique codimension 2 leaf $\fL_1 \subset X$.
    \item The Kleinian singularity $\Sigma_1$ is of type $A_1$.
    \item $\mathfrak{r}$ is semisimple.
\end{enumerate}
Property (1) is immediate from the weighted Dynkin diagram. Properties (2) and (3) are verified using \cite[Sec 13]{fuetal2015}. (4) is of course immediate from the determination of $\mathfrak{r}$. Note that (1) implies that $\OO=\Ind^G_L\{0\}$, see Proposition \ref{prop:even}. Properties (2)-(4) imply that $\dim(\fX(\fl))=\dim(\fP^X)=1$ and hence that $M_1=L$. 

In each case, we will show that $\widehat{\OO}$ is birationally rigid using Proposition \ref{prop:birrigidcriterion}. By that proposition, it is sufficient to exhibit a strict inequality
$$|\pi_1(\OO)||\pi_1(\OO_1)|^{-1} > m(\OO).$$
Once we have determined that $\widehat{\OO}$ is birationally rigid, we proceed to computing the infinitesimal character $\gamma_0(\widehat{\OO})$. By Proposition \ref{prop:unipotentcentralchars}, we have
\begin{equation}\label{eq:gamma0eqprelim}\gamma_0(\widehat{\mathbb{O}}) = \rho(\fl) + \delta,\end{equation}
where $\delta \in \fX(\fl)$ is the weight corresponding to the barycenter parameter $\epsilon \in \fP^X_1=\fP^X$ under the natural identification $\eta: \fX(\fl) \simeq \fP^X$. To compute this weight, we use Proposition \ref{prop:identification2}. In each case, we find that $|\pi_1(\mathbb{O})_{\mathrm{ab}}|=2$. Thus, 
$$c_1 = 2|\Pic(\OO_L)||\Pic(\mathbb{O})|^{-1} = 2|\Pic(\{0\})||\pi_1(\mathbb{O})_{\mathrm{ab}}|^{-1} = 1.$$
Let $\alpha_k$ be the simple root \emph{not} contained in $L$ and let $\varpi_k$ be the corresponding fundamental weight. Then by Proposition \ref{prop:identification2}, equation (\ref{eq:gamma0eqprelim}) becomes
\begin{equation}\label{eq:gamm0eq}\gamma_0(\widehat{\mathbb{O}}) = \rho(\Pi - \{\alpha_k\}) + \frac{1}{2}\varpi_k,\end{equation}
where $\rho(\Pi - \{\alpha_k\})$ is the half-sum of the positive roots in the Levi defined by $\Pi-\{\alpha_k\}$. We write this sum in fundamental weight coordinates and compute its dominant $W$-conjugate using {\tt atlas}. In each case, this dominant weight coincides with $\frac{1}{2}h_{\mathbb{O}^{\vee}}$, proving (ii) of Proposition \ref{prop:dualtodistinguished}.
\vspace{3mm}

\paragraph{Type $G_2$}

\begin{center}
$$\dynkin[labels={\alpha_1,\alpha_2},edge
length=.75cm] G2$$

\begin{tabular}{|c|c|c|c|c|c|c|c|c|} \hline
$\mathbb{O}^{\vee}$ & $\frac{1}{2}h_{\mathbb{O}^{\vee}}$ & $\mathbb{O}$ & $\mathfrak{r}$ & $\pi_1(\mathbb{O})$ & $L$ & $m(\OO)$ & $\OO'$ & $\pi_1(\OO')$\\ \hline
$G_2(a_1)$ & $(1,0)$ & $G_2(a_1)$ & $\{0\}$ & $S_3$ & $\Pi - \{\alpha_2\}$ & $1$ & $\widetilde{A}_1$ & $1$\\ \hline
\end{tabular}
\end{center}

\vspace{3mm}

\begin{itemize}
    \item \underline{$\mathbb{O}=G_2(a_1)$}. We have
$$|\pi_1(\OO)||\pi_1(\OO')|^{-1}=6>1=m(\OO).$$
Thus, $\widehat{\OO}$ is birationally rigid by \cref{prop:birrigidcriterion}. By (\ref{eq:gamm0eq}),
    $$\gamma_0(\widehat{\mathbb{O}}) = \frac{1}{2}\alpha_2 + \frac{1}{2}\varpi_1 = (1,0).$$
\end{itemize}

\vspace{3mm}

\paragraph{Type $F_4$}

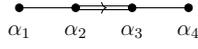
\begin{center}

$$\dynkin[labels={\alpha_1,\alpha_2,\alpha_3,\alpha_4},edge
length=.75cm] F4$$

\begin{tabular}{|c|c|c|c|c|c|c|c|c|} \hline
$\mathbb{O}^{\vee}$ & $\frac{1}{2}h_{\mathbb{O}^{\vee}}$ & $\mathbb{O}$ & $\mathfrak{r}$ & $\pi_1(\mathbb{O})$ & $L$ & $m(\OO)$ & $\OO'$ & $\pi_1(\OO')$\\ \hline
$F_4(a_3)$ & $(0,1,0,0)$ & $F_4(a_3)$ & $\{0\}$ & $S_4$ & $\Pi - \{\alpha_2\}$ & $2$ & $C_3(a_1)$ & $S_2$\\ \hline
\end{tabular}
\end{center}

\vspace{3mm}

\begin{itemize}
    \item \underline{$\mathbb{O}=F_4(a_3)$}. We have
$$|\pi_1(\OO)||\pi_1(\OO')|^{-1}=12>2=m(\OO).$$
Thus, $\widehat{\OO}$ is birationally rigid by \cref{prop:birrigidcriterion}. By (\ref{eq:gamm0eq}),
    $$\gamma_0(\widehat{\mathbb{O}}) = \rho(\Pi - \{\alpha_2\}) + \frac{1}{2}\varpi_2 = (1,-2,1,1).$$
    This is $W$-conjugate to $(0,1,0,0)$. 
\end{itemize}

\vspace{3mm}

\paragraph{Type $E_6$}

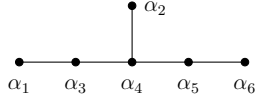
\begin{center}

$$\dynkin[labels={\alpha_1,\alpha_2,\alpha_3,\alpha_4,\alpha_5,\alpha_6},edge
length=.75cm] E6$$

\begin{tabular}{|c|c|c|c|c|c|c|c|c|} \hline
$\mathbb{O}^{\vee}$ & $\frac{1}{2}h_{\mathbb{O}^{\vee}}$ & $\mathbb{O}$ & $\mathfrak{r}$ & $\pi_1(\mathbb{O})$ & $L$ & $m(\OO)$ & $\OO'$ & $\pi_1(\OO')$\\ \hline

$E_6(a_3)$ & $(1,0,0,1,0,1)$ & $A_2$ & $2A_2$ & $2$ & $\Pi - \{\alpha_2\}$ & $1$ & $3A_1$ & $1$\\ \hline
\end{tabular}
\end{center}

\vspace{3mm}

\begin{itemize}
    \item \underline{$\mathbb{O}=A_2$}. We have
    $$|\pi_1(\OO)||\pi_1(\OO')|^{-1}=2>1=m(\OO).$$
    Thus, $\widehat{\OO}$ is birationally rigid by \cref{prop:birrigidcriterion}. By (\ref{eq:gamm0eq}),
    $$\gamma_0(\widehat{\mathbb{O}}) = \rho(\Pi - \{\alpha_2\}) + \frac{1}{2}\varpi_2 = (1,-4,1,1,1,1).$$
    This is $W$-conjugate to $(1,0,0,1,0,1)$.
\end{itemize}

\vspace{3mm}

\paragraph{Type $E_7$}

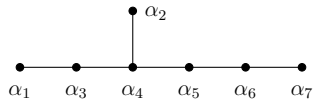
\begin{center}
$$\dynkin[labels={\alpha_1,\alpha_2,\alpha_3,\alpha_4,\alpha_5,\alpha_6,\alpha_7},edge
length=.75cm] E7$$

\begin{tabular}{|c|c|c|c|c|c|c|c|c|} \hline
$\mathbb{O}^{\vee}$ & $\frac{1}{2}h_{\mathbb{O}^{\vee}}$ & $\mathbb{O}$ & $\mathfrak{r}$ & $\pi_1(\mathbb{O})$ & $L$ & $m(\OO)$ & $\OO'$ & $\pi_1(\OO')$\\ \hline

$E_7(a_3)$ & $(1,0,0,1,0,1,1)$ & $A_2$ & $A_5$ & $S_2$ & $\Pi - \{\alpha_1\}$ & $1$ & $(3A_1)'$ & $1$\\ \hline

$E_7(a_5)$ & $(0,0,0,1,0,0,1)$ & $D_4(a_1)$ & $3A_1$ & $S_3$ & $\Pi - \{\alpha_3\}$ & $1$ & $(A_3+A_1)'$ & $1$\\ \hline
\end{tabular}
\end{center}

\vspace{3mm}

\begin{itemize}
    \item \underline{$\mathbb{O}=A_2$}. We have
    $$|\pi_1(\OO)||\pi_1(\OO')|^{-1}=2>1=m(\OO).$$
    Thus, $\widehat{\OO}$ is birationally rigid by \cref{prop:birrigidcriterion}. By (\ref{eq:gamm0eq}),
    $$\gamma_0(\widehat{\mathbb{O}}) = \rho(\Pi - \{\alpha_1\}) + \frac{1}{2}\varpi_1 = (-7,1,1,1,2,-1,2).$$
    This is  $W$-conjugate to $(1,0,0,1,0,1,1)$.
    
    \item \underline{$\mathbb{O}=D_4(a_1)$}. We have
    $$|\pi_1(\OO)||\pi_1(\OO')|^{-1}=6>1=m(\OO).$$
    Thus, $\widehat{\OO}$ is birationally rigid by \cref{prop:birrigidcriterion}. By (\ref{eq:gamm0eq})
    $$\gamma_0(\widehat{\mathbb{O}}) = \rho(\Pi- \{\alpha_3\}) + \frac{1}{2}\varpi_3 = (1,1,-4,1,2,-1,2)$$
    This is $W$-conjugate to $(0,0,0,1,0,0,1)$.
\end{itemize}

\vspace{3mm}

\paragraph{Type $E_8$}

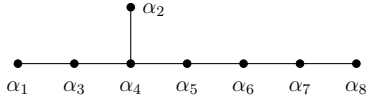
\begin{center}
$$\dynkin[labels={\alpha_1,\alpha_2,\alpha_3,\alpha_4,\alpha_5,\alpha_6,\alpha_7,\alpha_8},edge
length=.75cm] E8$$

\begin{tabular}{|c|c|c|c|c|c|c|c|c|} \hline
$\mathbb{O}^{\vee}$ & $\frac{1}{2}h_{\mathbb{O}^{\vee}}$ & $\mathbb{O}$ & $\mathfrak{r}$ & $\pi_1(\mathbb{O})$ & $L$ & $m(\OO)$ & $\OO'$ & $\pi_1(\OO')$\\ \hline

$E_8(a_3)$ & $(1,0,0,1,0,1,1,1)$ & $A_2$ & $E_6$ & $S_2$ & $\Pi - \{\alpha_8\}$ & $1$ & $3A_1$ & $1$\\ \hline

$E_8(a_5)$ & $(1,0,0,1,0,0,1,0)$ & $2A_2$ & $2G_2$ & $S_2$ & $\Pi - \{\alpha_1\}$ & $1$ & $A_2+3A_1$ & $1$\\ \hline

$E_8(b_5)$ & $(0,0,0,1,0,0,1,1)$ & $D_4(a_1)$ & $D_4$ & $S_3$ & $\Pi - \{\alpha_7\}$ & $1$ & $A_3+A_1$ & $1$\\ \hline

$E_8(b_6)$ & $(0,0,0,1,0,0,0,1)$ & $D_4(a_1)+A_2$ & $A_2$ & $S_2$ & $\Pi - \{\alpha_2\}$ & $1$ & $A_3+A_2+A_1$ & $1$\\ \hline

$E_8(a_6)$ & $(0,0,0,0,1,0,0,0)$ & $E_8(a_7)$ & $0$ & $S_5$ & $\Pi - \{\alpha_5\}$ & $2$ & $E_7(a_5)$ & $S_3$\\ \hline
\end{tabular}
\end{center}

\vspace{3mm}

\begin{itemize}
    \item \underline{$\mathbb{O} = A_2$}. We have
    $$|\pi_1(\OO)||\pi_1(\OO')|^{-1}=2>1=m(\OO).$$
    Thus, $\widehat{\OO}$ is birationally rigid by \cref{prop:birrigidcriterion}. By (\ref{eq:gamm0eq}),
    $$\gamma_0(\widehat{\mathbb{O}}) = \rho(\Pi - \{\alpha_8\}) + \frac{1}{2}\varpi_8 = (1,1,1,1,1,1,1,-13).$$
    This is $W$-conjugate to $(1,0,0,1,0,1,1,1)$.
    
    \item \underline{$\mathbb{O} = 2A_2$}. We have 
    $$|\pi_1(\OO)||\pi_1(\OO')|^{-1}=2>1=m(\OO).$$
    Thus, $\widehat{\OO}$ is birationally rigid by Proposition \ref{prop:birrigidcriterion}. By (\ref{eq:gamm0eq})
    $$\gamma_0(\widehat{\mathbb{O}}) = \rho(\Pi - \{\alpha_1\}) + \frac{1}{2}\varpi_1 = (-10,1,1,1,1,1,1,1),$$
    This is $W$-conjugate to $(1,0,0,1,0,0,1,0)$.

    \item \underline{$\mathbb{O}=D_4(a_1)$}. We have
    $$|\pi_1(\OO)||\pi_1(\OO')|^{-1}=6>1=m(\OO).$$
    Thus, $\widehat{\OO}$ is birationally rigid by \cref{prop:birrigidcriterion}. By (\ref{eq:gamm0eq}),
    $$\gamma_0(\widehat{\mathbb{O}}) = \rho(\Pi - \{\alpha_7\}) + \frac{1}{2}\varpi_7 = (1,1,1,1,1,1,-8,1).$$
    This is $W$-conjugate to $(0,0,0,1,0,0,1,1)$.
    
    \item \underline{$\mathbb{O}=D_4(a_1)+A_2$}. We have
    $$|\pi_1(\OO)||\pi_1(\OO')|^{-1}=2>1=m(\OO).$$
    Thus, $\widehat{\OO}$ is birationally rigid by \cref{prop:birrigidcriterion}. By (\ref{eq:gamm0eq})
    $$\gamma_0(\widehat{\mathbb{O}}) = \rho(\Pi - \{\alpha_2\}) + \frac{1}{2}\varpi_2 =  (1,-7,1,1,1,1,1,1)$$
    This is $W$-conjugate to $(0,0,0,1,0,0,0,1)$.

    \item \underline{$\mathbb{O}=E_8(a_7)$}. We have 
    $$|\pi_1(\OO)||\pi_1(\OO_1)|^{-1}=20>2=m(\OO).$$
    Thus, $\widehat{\OO}$ is birationally rigid by \cref{prop:birrigidcriterion}. By (\ref{eq:gamm0eq}),
    $$\gamma_0(\widehat{\mathbb{O}}) = \rho(\Pi - \{\alpha_5\}) + \frac{1}{2}\varpi_5=  (1,1,1,1,-4,1,1,1).$$
    This is $W$-conjugate to $(0,0,0,0,1,0,0,0)$.
\end{itemize}

\section{Refined BVLS duality}\label{subsec:refinedBVLS}

Let $\mathbb{O}^{\vee} \subset (\fg^{\vee})^*$ be a nilpotent $G^{\vee}$-orbit. By \cref{prop:birationalinduction}, there is a pair $(L^{\vee}, \mathbb{O}^{\vee}_0)$, unique up to conjugation by $G^{\vee}$, consisting of a Levi subgroup $L^{\vee} \subset G^{\vee}$ and a distinguished $L^{\vee}$-orbit $\mathbb{O}^{\vee}_0$, such that
$$\mathbb{O}^{\vee} = \mathrm{Sat}^{G^{\vee}}_{L^{\vee}} \mathbb{O}^{\vee}_0.$$
Let $\mathbb{O}_0 := \mathsf{D}(\mathbb{O}_0^{\vee})$. By Proposition \ref{prop:dualtodistinguished}, the universal $L$-equivariant cover $\widehat{\mathbb{O}}_0 \to \mathbb{O}_0$ is birationally rigid. Consider the equivalence class
\textbf{}\begin{equation}\label{eq:defoftilded}
\widetilde{\mathsf{D}}(\mathbb{O}^{\vee}) :=  [\mathrm{Bind}^G_L \widehat{\mathbb{O}}_0].\end{equation}
This defines a map
$$\widetilde{\mathsf{D}}: \{\text{nilpotent orbits for } G^{\vee}\} \to \{\text{equivalence classes of nilpotent covers for } G\}.$$

\begin{prop}
The map $\widetilde{\mathsf{D}}$ has the following properties:
\begin{itemize}
    \item[(i)] For every $\OO^{\vee}$, $\widetilde{\mathsf{D}}(\mathbb{O}^{\vee})$ covers $\mathsf{D}(\mathbb{O}^{\vee})$.
    \item[(ii)] If $\OO^{\vee}$ is distinguished, then $\widetilde{\mathsf{D}}(\mathbb{O}^{\vee})$ consists of birationally rigid covers. 
    
    \item[(iii)] $\widetilde{\mathsf{D}}$ is injective.

    \item[(iv)] For every $\OO^{\vee}$, there is an equality of ideals
    $$I_{\mathrm{max}}(\frac{1}{2}h_{\mathbb{O}^{\vee}}) = I_0(\widetilde{\mathsf{D}}(\mathbb{O}^{\vee})).$$
\end{itemize}
\end{prop}

\begin{proof}
Let $\OO^{\vee}=\mathrm{Sat}^{G^{\vee}}_{L^{\vee}}\OO_0^{\vee}$ for $\OO^{\vee}_0$ distinguished. By  (\ref{eq:defoftilded}),   $\widetilde{\mathsf{D}}(\OO^{\vee})$ covers $\Ind^G_L \mathsf{D}(\OO_0^{\vee})$, and $\Ind^G_L\mathsf{D}(\OO_0^{\vee}) = \mathsf{D}(\OO^{\vee})$ by Proposition \ref{prop:inclusioninduction}. This proves (i). (ii) is an easy consequence of Proposition \ref{prop:dualtodistinguished}(i). For (iii), first recall that all distinguished orbits are special (distinguished orbits are Richardson by Remark \ref{rmk:distinguishedbirational}, and Richardson orbits are special by Proposition \ref{prop:inclusioninduction}). Thus, $\mathsf{D}$ is injective when restricted to distinguished nilpotent orbits. Since $\widetilde{\mathsf{D}}(\OO^{\vee})$ covers $\mathsf{D}(\OO^{\vee})$, the restriction of  $\widetilde{\mathsf{D}}$ to distinguished orbits is injective as well. Now (iii) follows immediately from Proposition \ref{prop:birationalinduction}(iii). For (iv), note that
$$\frac{1}{2}h_{\mathbb{O}^{\vee}} = \frac{1}{2}h_{\mathbb{O}^{\vee}_0}
= \gamma_0[\widehat{\mathbb{O}}_0] 
= \gamma_0[\mathrm{Bind}^G_L \widehat{\mathbb{O}}_0] 
= \gamma_0(\widetilde{\mathsf{D}}(\mathbb{O}^{\vee})).$$
The second equality follows from Proposition \ref{prop:dualtodistinguished}, and the third from Proposition \ref{prop:Ibetainduced}. Thus by Proposition \ref{prop:propertiesofprim}(ii), there is an inclusion of primitive ideals $I_0(\widetilde{\mathsf{D}}(\mathbb{O}^{\vee})) \subseteq I_{\mathrm{max}}(\frac{1}{2}h_{\mathbb{O}^{\vee}})$. By Proposition \ref{prop:propsofIbeta}(i)
$$V(I_0(\widetilde{\mathsf{D}}(\mathbb{O}^{\vee}))) = \overline{\mathsf{D}(\mathbb{O}^{\vee})} = V(I_{\mathrm{max}}(\frac{1}{2}h_{\mathbb{O}^{\vee}})).$$
Hence, the inclusion $I_0(\widetilde{\mathsf{D}}(\mathbb{O}^{\vee})) \subseteq I_{\mathrm{max}}(\frac{1}{2}h_{\mathbb{O}^{\vee}})$ is an equality by Proposition \ref{prop:propertiesofprim}(iii).
\end{proof}

As an immediate consequence, we obtain the following result (Theorem \ref{thm:spec_unip} from the introduction).

\begin{cor}\label{cor:specialimpliesunipotent}
Every special unipotent bimodule (cf. Definition \ref{def:spec_unipotent}) is unipotent (cf. Definition \ref{def:unipotentbimods}).
\end{cor}

\begin{rmk}\label{rmk:classificationspecial}
In light of Corollary \ref{cor:specialimpliesunipotent}, our classification of unipotent bimodules (Theorem \ref{thm:classificationbimods}) gives, in particular, a classification of special unipotent bimodules (even in cases when $\mathbb{O}^{\vee}$ is not special and the results of \cite{BarbaschVogan1985} do not apply). The result is as follows. Let $\widetilde{\mathbb{O}}$ be the maximal cover in the equivalence class $\widetilde{\mathsf{D}}(\mathbb{O}^{\vee})$ and let $\Pi = \Aut_{\mathbb{O}}(\widetilde{\mathbb{O}})$. Then there is a natural bijection
$$\{\text{irreducible representations of } \Pi\} \xrightarrow{\sim} \unip_{\mathbb{O}^{\vee}}^s(G).$$
\end{rmk}

It is well-known that BVLS duality is order-reversing with respect to the closure ordering on orbits, see Proposition \cite[Prop A2]{BarbaschVogan1985}. We conjecture that refined BVLS duality satisfies an analogous property.

Define a partial order on the set of equivalence classes of equivariant nilpotent covers as follows. Write $[\widetilde{\OO}_1] \gg [\widetilde{\OO}_2]$ if 
\begin{itemize}
    \item $\OO_1 > \OO_2$ with respect to the closure ordering on orbits, or
    \item $\OO_1=\OO_2$ and the maximal cover in  $[\widetilde{\OO}_2]$ covers the maximal cover in $[\widetilde{\OO}_1]$.
\end{itemize} 

\begin{conj}\label{conj:orderreversing}
    The map 
    $$\widetilde{\mathsf{D}}: \{\text{nilpotent orbits for } G^{\vee}\} \to \{\text{equivalence classes of nilpotent covers for } G\}.$$
    is order-reversing.
\end{conj}

We can prove this conjecture for classical groups. 
The proof amounts to a tedious calculation with partitions, and we do not include it here. In the example below, we verify the conjecture for $G=G_2$. We expect similar arguments to work for other exceptional groups.

\begin{example}
For $G=G_2$, $\mathsf{D}$ is as follows (see \cite[Sec 13.4]{Carter1993}).
\begin{center}
\begin{tabular}{|c|c|} \hline
 $\mathbb{O}^{\vee}$ &  $\mathsf{D}(\OO^\vee)$ \\ \hline
$0$ & $G_2$ \\ \hline
$A_1$ & $G_2(a_1)$  \\ \hline
$\widetilde{A}_1$ & $G_2(a_1)$  \\ \hline
$G_2(a_1)$ & $G_2(a_1)$  \\ \hline
$G_2$ & $\{0\}$ \\ \hline
\end{tabular}
\end{center}
Let $\OO_1^\vee=G_2(a_1)$, $\OO_2^\vee=\widetilde{A}_1$, and $\OO_3^\vee=A_1$. Note that $\OO_1^\vee>\OO_2^\vee>\OO_3^\vee$. Thus, it suffices to show that $\widetilde{\mathsf{D}}(\OO_1^{\vee}) \ll \widetilde{\mathsf{D}}(\OO_2^{\vee}) \ll \widetilde{\mathsf{D}}(\OO_3^{\vee})$. 

Define the Levi subgroups $L_1,L_2 \subset G$ as in Example \ref{ex: G2covers}. Since $\OO_1^{\vee}$ is distinguished, $\widetilde{\mathsf{D}}(\OO_1^{\vee})$ contains the universal cover of $\mathsf{D}(\OO_1^{\vee}) = G_2(a_1)$. On the other hand, $\OO_2^{\vee}$ and $\OO_3^{\vee}$ are saturated from the principal orbits for the short and long root Levis. Thus, 
$$\widetilde{\mathsf{D}}(\OO_2^{\vee}) = [\mathrm{Bind}^G_{A_1} \{0\}], \qquad \widetilde{\mathsf{D}}(\OO_3^{\vee}) = [\mathrm{Bind}^G_{\widetilde{A}_1} \{0\}],$$
Using the computations in \cref{ex: G2covers} we deduce
$$\widetilde{\mathsf{D}}(\OO_2^\vee) = \{\widetilde{\OO}_3\}, \qquad \widetilde{\mathsf{D}}(\OO_3^\vee)=\{\OO\}.$$
In particular, $\widetilde{\mathsf{D}}(\OO_1^\vee)\ll\widetilde{\mathsf{D}}(\OO_2^\vee)\ll \widetilde{\mathsf{D}}(\OO_3^\vee)$, as asserted.
\end{example}

\section{Motivation: symplectic duality}\label{subsec:motivationsymplectic}

Our construction of refined BVLS duality should be viewed as a special case of a more general (but still largely conjectural) duality known as {\it symplectic duality}\index{duality!symplectic} (a closely related duality is known as the {\it 3D mirror symmetry}), see \cite[Section 10]{BPWII}. Since there is no rigorous definition of symplectic duality in the cases that concern us, our exposition in this section will be largely speculative. 

Let $X$ be a conical symplectic singularity, cf. Definition \ref{def:conicalsymplecticsing}. To $X$ one can attach numerous invariants which come in pairs. The most basic example is the pair $(\mathfrak{P}^X, \mathfrak{t}^X)$. As usual, $\mathfrak{P}^X$ is the (complex) Namikawa space, parameterizing the Poisson deformations of a $\QQ$-factorial terminalization $Y \to X$. The second space, $\mathfrak{t}^X$, is constructed as follows. Consider the group $H^X$ of graded Hamiltonian automorphisms of $\CC[X]$. This is an algebraic group. Choose a maximal torus $T^X \subset H^X$, and let $\mathfrak{t}^X := \mathrm{Lie}(T^X)$.

\begin{example}\label{Ex:sd_cover}
Let $X$ be the affinization of a $G$-equivariant nilpotent cover. Then $\mathfrak{P}^X \simeq \mathfrak{X}(\mathfrak{l})$, where $\fl$ is a Levi subalgebra of $\fg$ such that $X$ is birationally induced from a birationally rigid $L$-equivariant nilpotent cover, see Proposition \ref{prop:namikawacovers}. In most cases, $H^X=G$ (assuming $G$ acts faithfully on $X$), and therefore $\mathfrak{t}^X$ is a Cartan subalgebra of $\fg$. However, in some cases, the maximal reductive subgroup of $H^X$ can be larger, see \cite[Thm 4]{BrylinskiKostant1994}.
\end{example}

\begin{example}\label{Ex:sd_slice}
Let $e^\vee \in \mathcal{N}^{\vee}$ be a nilpotent element and choose an $\mathfrak{sl}(2)$-triple $(e^{\vee},f^{\vee},h^{\vee})$. For $X^\vee$ take the intersection of $\mathcal{N}^\vee$ with the Slodowy slice $S^\vee$ to $e^\vee$.  
In this case, $Y^\vee$ is the preimage of $S^\vee$ under the Springer resolution.
In most cases, we have $\mathfrak{P}^{X^\vee} \simeq \mathfrak{h}$, see \cite[Theorem 1.3]{LehnNamikawaSorger}, and the group $H^{X^\vee}$ coincides with $Z_{G^\vee}(e^\vee,h^\vee,f^\vee)$. So $\mathfrak{t}^{X^\vee}$ is identified with connected component of the center of a minimal Levi subalgebra containing $e^\vee$.       
\end{example}

To a first approximation, symplectic duality is a conjectural duality between conical symplectic singularities. To get a precise correspondence, the conical symplectic singularities, on both sides, should be equipped with certain `decorations', whose precise nature is unclear in the cases that concern us. We will ignore this complication. Write $X^{\vee}$ for the symplectic dual of $X$. Some expected properties of symplectic duality are as follows:

\begin{itemize}
    \item $(X^{\vee})^{\vee} \simeq X$.
    \item There are natural isomorphisms
    $$\mathfrak{t}^X \simeq \fP^{X^{\vee}}, \qquad \fP^X \simeq \mathfrak{t}^{X^{\vee}}.$$
    \item The following conditions are equivalent:
    \begin{itemize}
        \item[(i)] $T^X$ acts on $X$ with unique fixed point, $0$. 
        \item[(ii)] $X^\vee$ admits a symplectic resolution (equivalently, all $\mathbb{Q}$-factorial terminalizations are smooth). 
    \end{itemize}
\end{itemize}

Let $Y$ be a $\QQ$-factorial terminalization of $X$. It is expected that the slices to the symplectic leaves in $Y$ are the formal neighborhoods of $0$ in conical symplectic singularities. Let $X_1,...,X_k$ be the conical symplectic singularities corresponding to the \emph{minimal} symplectic leaves. It is also expected that the connected components of the normalization of the $T^{X^{\vee}}$-fixed point locus in $X^{\vee}$ are conical symplectic singularities (compare to \cite[Section 5.1]{LosevCatO}), denoted $X_1',...,X_{\ell}'$. Then a general expectation is that $k=\ell$ and (up to permutation of indices), $X_i^{\vee} \simeq X_i'$ for all $i$. Note that this expectation generalizes the third bullet above.

Now let $X^\vee$ be as in Example \ref{Ex:sd_slice}. It is natural to expect that $X$ should have something to do with nilpotent orbits in $\fg$, see, e.g. \cite[Section 10.4]{BPWII}. Suppose $e^\vee$ is distinguished. Then $T^{X^\vee}=\{1\}$. Since $\mathfrak{t}^{X^\vee} \simeq \fP^{X}$, $X$ should have no Poisson deformations. If we assume that $X^\vee$ is the affinization of a nilpotent cover, this cover should be birationally rigid, see Corollary \ref{cor:criterionbirigid}. This motivates the first part of Proposition \ref{prop:dualtodistinguished}. 

Now consider an arbitrary element 
$e^{\vee}\in \cN^{\vee}$ and the corresponding variety $X^\vee$. Let $\mathbb{O}^\vee$ denote the $G^{\vee}$-orbit of $e^{\vee}$, and $\mathfrak{l}^\vee$ be a minimal Levi subalgebra containing $e^{\vee}$. Let $\underline{S}^\vee$ denote the Slodowy slice to $e^\vee$ in $\mathfrak{l}^\vee$ and let $\underline{X}^\vee$ be the intersection of $\underline{S}^\vee$
with the nilpotent cone in $\fl^\vee$. Thus $T^{X^\vee}=Z(L^\vee)^{\circ}$. The fixed point locus of $T^{X^\vee}$ in $S^\vee$ is $\underline{S}^\vee$. It follows that the fixed point locus of $T^{X^\vee}$ in $X^\vee$ is 
$\underline{X}^\vee$. 

Now let $X$ be the affinization of $\widetilde{\mathsf{D}}(\Orb^\vee)$. We have $\fP^X \simeq \mathfrak{t}^{X^\vee}$ and, in most cases, $\mathfrak{t}^X\simeq\fP^{X^\vee}\simeq\mathfrak{h}$. The slice to the minimal leaf in $Y$ is identified with the symplectic dual of $\underline{X}^\vee$. We know that $Y$ is of the form $G\times^P(\fp^{\perp} \times X_L)$. So the slice in question is $X_L$. This discussion motivates the construction of refined BVLS duality in Section \ref{subsec:refinedBVLS}. 

Another justification for viewing refined BVLS duality as a special case of symplectic duality can be found in the ongoing work of Finkelberg, Hanany, and Nakajima. If $G$ is either $\operatorname{SO}(n)$ or $\operatorname{Sp}(2n)$, for many $\Orb^\vee$ they construct an orthosymplectic quiver gauge theory with the Higgs branch isomorphic to $X^\vee$, and the Coulomb branch isomorphic to an affinization of a certain $G$-equivariant cover $\widetilde{\OO}$ of the BVLS dual orbit $\OO=\mathsf{D}(\Orb^\vee)$. We expect that $\widetilde{\OO}$ is equivalent to $\widetilde{\mathsf{D}}(\Orb^\vee)$.

Finally, we speculate on how the special unipotent infinitesimal character $\frac{1}{2}h^\vee$ relates to symplectic duality. Let $X$ be a conical symplectic singularity. Write $\mathfrak{P}$ for $\mathfrak{P}^X$ to simplify the notation. 
Suppose that $0$ is the unique $T^X$-fixed point in $X$. Let $\mathcal{A}_{\fP,\hbar}$ denote the Rees algebra of $\Gamma(\calD^{Y,\mathrm{univ}})$, where
$\calD^{Y,\mathrm{univ}}$ is the universal quantization of $Y$, see Proposition \ref{prop:universal for Q-term}. $\mathcal{A}_{\fP,\hbar}$ is a graded algebra over $\C[\fP,\hbar]$, with a Hamiltonian $T^X$-action. Choose a generic one-parameter subgroup $\nu:\C^\times\rightarrow T^X$, where `generic' means that $\nu(\C^\times)$ has a unique fixed point in $X$. The choice of $\nu$ defines a grading on $\mathcal{A}_{\fP,\hbar}$. Write $\mathcal{A}_{\fP,\hbar}^i$ for the $i$th graded component. Now consider the so called \emph{Cartan subquotient}\index{Cartan subquotient} of $\mathcal{A}_{\fP,\hbar}$ (called the \emph{B-algebra} in \cite{BPWII}) 
$$\mathsf{C}_\nu(\mathcal{A}_{\fP,\hbar}):=
\mathcal{A}_{\fP,\hbar}^0/\sum_{i>0}
\mathcal{A}_{\fP,\hbar}^{-i}\mathcal{A}_{\fP,\hbar}^i.$$
Note that $\mathsf{C}_\nu(\mathcal{A}_{\fP,\hbar}) $ has the structure of a graded $\C[\fP,\hbar]$-algebra. The condition that $0$ is the unique $T^X$-fixed point in $X$ is equivalent to the condition that
$\mathsf{C}_\nu(\mathcal{A}_{\fP,\hbar})$ is a finitely generated $\CC[\fP,\hbar]$-module. On the other hand, the quantum co-moment map $S(\mathfrak{t}^X)\rightarrow \mathcal{A}_{\fP,\hbar}$ gives rise to a homomorphism $S(\mathfrak{t}^X)\rightarrow 
\mathsf{C}_\nu(\mathcal{A}_{\fP,\hbar})$. In fact, there is a canonical choice of co-moment map, cf. \cite[Section 5.4]{Losev_isofquant}. So $\mathsf{C}_\nu(\A_{\fP,\hbar})$ has a distinguished $\C[(\mathfrak{t}^X)^*\oplus \fP,\hbar]$-algebra structure. 

One can write down a similar algebra on the $X^\vee$-side motivated by   a conjecture of Hikita \cite{Hikita} and its extension due to Nakajima, see \cite[Conjecture 8.9]{KamnitzerTingleyWebsterWeeksYacobi}. We call this the {\it deformed Hikita conjecture}. Consider the equivariant cohomology algebra $H^*_{T^{X^\vee}\times \C^\times}(Y^\vee)$, where $\C^\times$ is a contracting torus. This is an algebra over $\C[\operatorname{Lie}(T^{X^\vee}\times \C^\times)]$. For $\hbar$ we take $d\in \operatorname{Lie}(\C^\times)$, where $-d$ is the weight of the Poisson bracket on $X^\vee$. 

In a related but different setup, Nakajima has conjectured
a graded $\C[\fP,\hbar]$-algebra isomorphism
\begin{equation}\label{eq:deformed_Hikita}
    \mathsf{C}_\nu(\A_{\fP,\hbar})\xrightarrow{\sim}H^*_{T^{X^\vee}\times \C^\times}(Y^\vee).
\end{equation}

In \cite{KMH} Hoang, Krylov, and the third named author show that this conjecture is not true as stated, see Section 6.4 of \emph{loc.cit.} for explicit counterexamples. A modified version of the conjecture is proposed in Sections 8 and 9. To avoid the technicalities of the modified conjecture, we assume that (\ref{eq:deformed_Hikita}) is an isomorphism.


%
In general, we do not know how to define a $\CC[\mathfrak{t}^{X*}]$-algebra structure  on the right hand side of (\ref{eq:deformed_Hikita}) or how to choose the contracting action on $Y^\vee$. However, either of these is not an issue when $X^\vee$ is the nilpotent part of a Slodowy slice. Namely, there is a natural contracting $\C^\times$-action (recalled in 
Section \ref{subsec:W} and known as the Kazhdan action) on $Y^\vee$. In the case when $e^\vee$ is distinguished, this is the unique action such that the weight of the Poisson bracket is $-2$. For $Y^\vee$ we take the preimage of $S^\vee$ under the Springer resolution $T^*\mathcal{B}^\vee \to\cN^{\vee}$. The $\C[\mathfrak{t}^{X*}\oplus \fP,\hbar]$-algebra structure on $H^*_{T^{X^\vee}\times \C^\times}(Y^\vee)$
comes by pullback from $Y^\vee\rightarrow \mathcal{B}^\vee$. Namely, there is a pullback homomorphism in equivariant cohomology, $$H^*_{T^{X^\vee}\times \C^\times}(\mathcal{B}^\vee)\rightarrow H^*_{T^{X^\vee}\times \CC^\times}(Y^\vee),$$ and the source algebra is a quotient of $\CC[(\mathfrak{t}^X)^*\oplus \fP,\hbar]$. It is natural to conjecture that (\ref{eq:deformed_Hikita}) is $\CC[(\mathfrak{t}^X)^*\oplus \fP,\hbar]$-linear. We note that this is consistent with the modifications to (\ref{eq:deformed_Hikita}) proposed in \cite{KMH}. The modified conjecture replaces the right hand side of (\ref{eq:deformed_Hikita}) with the image of the pullback homomorphism in equivariant cohomology. The modification of the left hand side is a bit more subtle and will not be important for what follows.

Let us return to the issue of the special unipotent infinitesimal character. By construction, the homomorphism $\CC[(\mathfrak{t}^X)^*\oplus \fP,\hbar]\rightarrow H^*_{T^{X^\vee}\times \CC^\times}(Y^\vee)$ factors through $H^*_{T^{X^\vee}\times \CC^\times}(\mathcal{B}^\vee)$, where $\CC^\times$ acts on $\mathcal{B}^\vee$ via the Kazhdan action. The formula for the Kazhdan action is $t\mapsto t^{-2}t^{h^\vee}$, where $t^{h^\vee}$ denotes the 1-parameter subgroup of $G^\vee$ corresponding to the element $h^\vee$. 

Since we are interested in the canonical quantization, we specialize $\fP$ to $0$ in (\ref{eq:deformed_Hikita}). Thus, we obtain an isomorphism
\begin{equation}\label{eq:deformed_Hikita1}
    \mathsf{C}_\nu(\A_{0,\hbar})\xrightarrow{\sim}H^*_{\CC^\times}(X^\vee).
\end{equation}
We see that (\ref{eq:deformed_Hikita1}) factors through the algebra $H^*_{\mathbb{C}^\times}(\mathcal{B}^\vee)=\mathbb{C}[\mathfrak{t}^\vee]\otimes_{\mathbb{C}[\mathfrak{t^\vee}/W]}\mathbb{C}[\hbar]$, where the homomorphism $\mathbb{C}[\mathfrak{t}^\vee/W]\rightarrow \mathbb{C}[\hbar]$ comes from the morphism $\mathbb{C}\rightarrow \mathfrak{t}/W, z\mapsto W(zh^\vee/2)$. It follows that the infinitesimal character of the kernel of $U(\mathfrak{g})\rightarrow \mathcal{A}_0$ is $\frac{1}{2}h^\vee$.


Our definition of refined BVLS duality raises the following question: what is the symplectic dual of the affinization of a more general nilpotent cover? This question is further explored in the paper \cite{MBMatYu}.

\chapter{Construction and unitarity of unipotent bimodules}\label{sec:unipbimod}

The main goal of this chapter is to construct all unipotent bimodules for linear classical groups and to deduce that they are unitary. Our approach is as follows. Suppose $G$ is linear classical, and let $\widetilde{\OO}$ be a $G$-equivariant nilpotent cover. We will construct a Levi subgroup $L \subset G$ and a rigid orbit $\OO_L$ such that all bimodules in $\unip_{\widetilde{\OO}}(G)$ are built out of bimodules in $\unip_{\OO_L}(L)$ via the following operations:
\begin{itemize}
    \item {Tensoring with unitary characters (i.e., in the notation of (\ref{eqn:1dimbimod}), bimodules of the form $\CC(-\frac{\chi}{2},\frac{\chi}{2})$ for $\chi \in \mathfrak{X}(L)$).}
    \item Unitary induction.
    \item Complementary series.
    \item Extraction of direct summands.
\end{itemize}
Our proofs rely on the classification of unipotent bimodules from Section \ref{subsec:classificationbimods} and a notion of parabolic induction for Harish-Chandra bimodules over Hamiltonian quantizations, which is developed in Section \ref{subsec:bimodinductioncovers}. 

Since all four operations in the list above preserve unitarity, our construction reduces the question of unitarity to the case of rigid nilpotent orbits. For such orbits, the unitarity of $\unip_{\OO}(G)$
follows from a classical result of Barbasch (\cite[Prop 10.6]{Barbasch1989}). In Section \ref{subsec:spinexceptional}, we will discuss the applicability of this strategy to spin and exceptional groups.

We note that Dougal Davis and and the second-named author have recently proved the unitarity of \emph{all} unipotent bimodules, including for spin and exceptional groups (\cite[Corollary 5.23]{DavisMasonBrown}), using completely different methods. Their proof is uniform across types, but it does not give a construction of unipotent bimodules by parabolic induction, which may be of independent interest.

\section{Unitarity of Harish-Chandra bimodules}\label{subsec:unitarydef}

In this section, we will recall what it means for a Harish-Chandra bimodule $\cB \in \HC^G(U(\fg))$ to be \emph{unitary}.

Fix a maximal compact subgroup $K \subset G$. There is an anti-holomorphic involution $\sigma$ of $G$ such that $K = G^{\sigma}$. The differential of $\sigma$ is a conjugate-linear Lie algebra involution of $\mathfrak{g}$ (still denoted by $\sigma$). Let $\cB \in \HC^G(U(\fg))$. A Hermitian form on $\cB$ is a sesquilinear pairing
$$\langle \ , \ \rangle: \cB \times \cB \to \CC$$
such that
$$\langle v, w \rangle = \overline{\langle w, v\rangle}, \qquad v,w \in \cB.$$
This form is said to be $\sigma$-\emph{invariant} if the following condition is satisfied
$$
\langle Xv, w\rangle = \langle v,  w\sigma(X) \rangle, \qquad 
X \in \fg, \quad v,w \in \cB.
$$
We say that $\cB$ is \emph{Hermitian}\index{Harish-Chandra bimodule!Hermitian} if it admits a non-degenerate $\sigma$-invariant Hermitian form. If $\cB$ is irreducible and Hermitian, then by a version of Schur's lemma, this form is unique up to multiplication by $\mathbb{R}^{\times}$. On the level of Langlands parameters (cf. Theorem \ref{thm:Langlands}) there is a simple criterion for deciding whether an irreducible bimodule is Hermitian. Recall from Theorem \ref{thm:Langlands} that if $(\lambda_{\ell},\lambda_r) \in \fh^* \times \fh^*$ is a Langlands parameter, the induced bimodule $I^G_H(\lambda_{\ell},\lambda_r) := \Ind^G_H \CC(\lambda_{\ell},\lambda_r)$ has a unique irreducible subquotient $\overline{I}(\lambda_{\ell},\lambda_r)$ containing the irreducible $G$-representation of extremal weight $\lambda_{\ell}-\lambda_r$.

\begin{prop}[\cite{Duflo1979}, Sec 3, see also \cite{Knapp}, Chp XVI]\label{prop:Hermitian}
Let $(\lambda_{\ell},\lambda_r)$ be a Langlands parameter for $G$. Then the following conditions are equivalent
\begin{itemize}
    \item[(i)] $I^G_H(\lambda_{\ell},\lambda_r)$ is Hermitian.
    \item[(ii)] $\overline{I}^G_H(\lambda_{\ell},\lambda_r)$ is Hermitian.
    \item[(iii)] There is an element $w \in W$ such that
$$w(\lambda_{\ell} - \lambda_r) = \lambda_{\ell} - \lambda_r, \qquad w(\lambda_{\ell} + \lambda_r) = - \overline{(\lambda_{\ell}+\lambda_r)}.$$
\end{itemize}
\end{prop}

A $\sigma$-invariant Hermitian form is \emph{positive-definite} if it satisfies the following additional condition
$$\langle v, v\rangle > 0, \qquad 0 \neq v \in \cB.$$
We say that $\cB$ is \emph{unitary}\index{Harish-Chandra bimodule!unitary} if it admits a positive-definite $\sigma$-invariant Hermitian form. There is no simple condition like Proposition \ref{prop:Hermitian} for deciding whether an irreducible bimodule is unitary. For classical groups, a classification of unitary Harish-Chandra bimodules was obtained in \cite{Barbasch1989}. For general groups, there is no known classification.

Finally, we explain the behavior of unitarity under parabolic induction. Let $M \subset G$ be a Levi subgroup of $G$. Recall the functor $\Ind^G_M: \HC^M(U(\fm)) \to \HC^G(U(\fg))$ defined in Section \ref{subsec:HCbimodsclassical}. 

\begin{prop}[\cite{GelfandNaimark}, see also \cite{Knapp}, Chp XVI]\label{prop:inductionunitary}
Choose continuous functions 
$$\lambda_{\ell}, \lambda_r: [0,1] \to \fX(\fm) \text{ such that } \lambda_{\ell}(t) - \lambda_r(t) \in \fX(M), \quad \forall t \in [0,1],$$
and let $\cB' \in \HC^M(U(\fm))$. Consider the one-parameter families of Harish-Chandra bimodules
$$\cB'(t) := \cB' \otimes \CC(\lambda_{\ell}(t),\lambda_r(t)) \in \HC^M(U(\fm)), \qquad \cB(t) := \Ind^G_M \cB'(t) \in \HC^G(U(\fg)), \qquad t \in [0,1].$$
Suppose $\cB'(0)$ is unitary. Then
\begin{enumerate}
    \item $\B(0)$ is unitary.
    \item If $\B(t)$ is irreducible and Hermitian for every $t \in [0,1]$, then $\B(1)$ is unitary.
\end{enumerate}
\end{prop}
In the setting of the proposition above, it is traditional to say that $\cB(0)$ is \emph{unitarily induced}\index{induction!unitary} from $\cB'(0)$. Indeed, under the correspondence between Harish-Chandra bimodules and continuous representations of $G$, the functor $\Ind^G_M$ corresponds to (normalized) parabolic induction and $\cB'(0)$ corresponds to a unitary representation of $M$. If $\CC(\lambda_{\ell}(1),\lambda_r(1))$ is a non-unitary character (as it will be in applications), then $\cB(1)$ is \emph{not} unitarily induced from $\cB'(1)$, but is unitary nonetheless. In this instance, it is traditional to say that $\cB(1)$ is obtained from $\cB'(0)$ through a \emph{complementary series construction}, see \cite[Chp XVI]{Knapp}.\index{complementary series}

\section{Parabolic induction of bimodules for Hamiltonian quantizations}\label{subsec:bimodinductioncovers}

In this section, we will define the notion of parabolic induction for Harish-Chandra bimodules for Hamiltonian quantizations of nilpotent covers. The goal is to construct 
unipotent bimodules attached to covers of induced orbits via parabolic induction. We note that our construction is closely related to the classical notion of parabolic induction of Harish-Chandra bimodules, see Section
\ref{subsec:HCbimodsclassical}. In Appendix \ref{sec:coincidence} we will show that, under suitable assumptions, these two constructions coincide. Our version has the advantage of making certain geometric properties of parabolic induction much more transparent, see e.g. Proposition \ref{prop:twinductiondagger}.

Choose a parabolic subgroup $P = LN \subset G$ and a birationally rigid $L$-equivariant nilpotent cover $\widetilde{\OO}_L$ such that $\widetilde{\OO} = \mathrm{Bind}^G_L\widetilde{\OO}_L$. Choose also a parabolic subgroup $Q=MU \subset G$ such that $P \subset Q$ and $L \subset M$. Let $\widetilde{\OO}_M = \mathrm{Bind}^M_L \widetilde{\OO}_L$, $\widetilde{X}_M = \Spec(\CC[\widetilde{\OO}_M])$, and form the partial resolutions 
$$\rho: \widetilde{Y}:=G\times^P (\widetilde{X}_L \times \fp^{\perp}) \to \widetilde{X}, \qquad \overline{\rho}: \widetilde{Z}:=G \times^Q (\widetilde{X}_M \times \fq^{\perp}) \to \widetilde{X},$$ 
and the projections
$$\pi: \widetilde{Y} \to G/P, \qquad \overline{\pi}: \widetilde{Z} \to G/Q.$$
Choose parameters $\beta_1,\beta_2 \in \fX(\fl)$ and let $\cA_{\beta_1}^{\widetilde{X}_M}$, $\cA_{\beta_2}^{\widetilde{X}_M}$ be the corresponding Hamiltonian quantizations of $\widetilde{X}_M$. A $\cA_{\beta_1}^{\widetilde{X}_M}$-$\cA_{\beta_2}^{\widetilde{X}_M}$-bimodule $\cB$ is \emph{Harish-Chandra} if it admits a good filtration, see the discussion preceding Definition \ref{def:HCbimodfiltered}. {This is equivalent to the definition used in Section \ref{SS_HC_different}}. Let $\HC(\cA_{\beta_1}^{\widetilde{X}_M},\cA_{\beta_2}^{\widetilde{X}_M})$ denote the category of Harish-Chandra $\cA_{\beta_1}^{\widetilde{X}_M}$-$\cA_{\beta_2}^{\widetilde{X}_M}$-bimodules, and define the full subcategory $\HC^M(\cA_{\beta_1}^{\widetilde{X}_M},\cA_{\beta_2}^{\widetilde{X}_M})$ of $M$-equivariant bimodules analogously to Definition \ref{def:eqvtbimods}. There is a forgetful functor
$$\HC^M(\cA_{\beta_1}^{\widetilde{X}_M},\cA_{\beta_2}^{\widetilde{X}_M}) \to \HC^M(U(\fm))$$
defined using the co-moment maps $\Phi^{\widetilde{X}_M}_{\beta_1}$ and $\Phi^{\widetilde{X}_M}_{\beta_2}$. 

Given a Harish-Chandra bimodule $\cB \in \HC^M(\cA_{\beta_1}^{\widetilde{X}_M},\cA_{\beta_2}^{\widetilde{X}_M})$, we will produce an induced bimodule
$$\Ind^G_M \cB \in \HC^G(\cA_{\beta_1}^{\widetilde{X}},\cA_{\beta_2}^{\widetilde{X}}),$$
thus defining a functor $\HC^M(\cA_{\beta_1}^{\widetilde{X}_M},\cA_{\beta_2}^{\widetilde{X}_M}) \to \HC^G(\cA_{\beta_1}^{\widetilde{X}},\cA_{\beta_2}^{\widetilde{X}})$ called \emph{parabolic induction}. 

The construction requires a bit of preparation. Let $\mathcal{D}^{\widetilde{X}_M}_{\beta_i}$ denote the microlocalization of $\cA_{\beta_i}^{\widetilde{X}_M}$ over $\widetilde{X}_M$ (for $i=1,2$), and let $\mathcal{E}$ denote the microlocalization of $\cB$ over $\widetilde{X}_M$ (the microlocalization of a filtered quantization was defined in Section \ref{subsec:quant}; the definition is analogous for Harish-Chandra bimodules). Note that $\mathcal{E}$ is an $M$-equivariant Harish-Chandra $\mathcal{D}^{\widetilde{X}_M}_{\beta_1}$-$\mathcal{D}^{\widetilde{X}_M}_{\beta_2}$-bimodule (the sheaf-theoretic version of a Harish-Chandra bimodule was discussed in Remark 
\ref{rmk:HC_sheaf}). As in Section  \ref{subsec:inductionquantizations}, form the completed tensor product 
$\mathcal{D}_{\beta_i}':=\mathfrak{D}_{G/U} \  \widehat{\otimes} \ \mathcal{D}^{\widetilde{X}_M}_{\beta_i}$.
Then let
$$\mathcal{E}' := \mathfrak{D}_{G/U} \  \widehat{\otimes} \ \mathcal{E}.$$
Note that $\mathcal{E}'$ is an $M$-equivariant Harish-Chandra $\mathcal{D}_{\beta_1}'$-$\mathcal{D}_{\beta_2}'$-bimodule. A good filtration can be constructed as follows. Pick a good filtration $\cB_{\le i}$ on $\cB$, and consider the corresponding filtration $\cE_{\le i}$ on $\cE$. Consider the doubled order filtration $\mathfrak{D}_{G/U,\le i}$ (so that vector fields lie in degree $\mathfrak{D}_{G/U,\le 2}$). Consider the tensor product filtration $\cE'_{\leqslant i}$ on $\cE'$. It is $G\times M$-stable and good. Then there is a natural identification
\begin{equation}\label{eq:product_assoc_graded}
\gr \cE'\simeq \mathcal{O}_{T^*(G/U)}\boxtimes \gr \mathcal{E}.
\end{equation}

Let $\mu: T^*(G/U)\times \widetilde{X}_M\to \fm^*$ be the moment map for the $M$-action on $T^*(G/U) \times \widetilde{X}_M$, and $\Phi_{\beta_1}'$, $\Phi_{\beta_2}'$ be the quantum comoment maps for $\calD_{\beta_1}'$ and $\calD_{\beta_1}'$, respectively. Here we use the same shift as in Section 
\ref{subsec:inductionquantizations}. The quotient $\mathcal{E}'/\mathcal{E}'\Phi_{\beta_2}(\fm)$ is a (weakly) $M$-equivariant $\mathcal{D}'_{\beta_1}$-module, set-theoretically supported on $\mu^{-1}(0) \subset T^*(G/U) \times \widetilde{X}_M$. Consider the quotient morphism $q: \mu^{-1}(0) \to G \times^Q (\widetilde{X}_M \times \fq^{\perp}) = \widetilde{Z}$ and set
$$\Ind^G_M \mathcal{E} := \left(q_*[\mathcal{E}'/\mathcal{E}'\Phi_{\beta_2}(\fm)]\right)^M.$$
By construction, $\Ind^G_M\mathcal{E}$ has the structure of a right module for $\Ind^G_M \mathcal{D}_{\beta_2}^{\widetilde{X}_M}$. It is also a left module for $\Ind^G_M \mathcal{D}_{\beta_1}^{\widetilde{X}_M}$ (this follows from the $M$-equivariance of $\mathcal{E}'$). These two actions commute, i.e. $\Ind^G_M\mathcal{E}$ has the structure of an $\Ind^G_M \mathcal{D}_{\beta_1}^{\widetilde{X}_M}$-$\Ind^G_M \mathcal{D}_{\beta_2}^{\widetilde{X}_M}$ bimodule. This bimodule has a $G$-action coming from the $G$-action on $\mathfrak{D}_{G/U}$. It also inherits a filtration from $\cE'$, which is automatically complete and separated. Our next task is to describe its associated graded.  

For any $M$-equivariant coherent sheaf $\cF$ on $\widetilde{X}_M$, we will construct a $G$-equivariant coherent sheaf $\Ind^G_M \cF$ on $\widetilde{Z}$. View $\cF$ as a $Q$-equivariant sheaf via the natural map $Q \to M$ and let $p_1$ denote the projection $\widetilde{X}_M\times \fq^\perp\rightarrow \widetilde{X}_M$. Note that the pullback $p_1^*\mathcal{F}$ is $Q$-equivariant. Restriction to the fiber over $1Q\in G/Q$ induces a category equivalence $\Coh^G(\widetilde{Z})\xrightarrow{\sim} \Coh^Q(\widetilde{X}_M\times \fq^\perp)$. 
Let $\Ind^G_M\cF$ be the object in
$\Coh^G(\widetilde{Z})$ corresponding to $p_1^*\mathcal{F}$ under this equivalence.

\begin{rmk}\label{rmk:indofsheaves}
Note that the above is a special case of a more general construction. Let $Z_M$ be a variety with $M$-action, and let $Z=G\times^Q(Z_M\times \fq^\perp)$. Then for $\cF \in \Coh^M(Z_M)$, the construction above gives rise to a sheaf $\Ind_M^G \cF \in \Coh^G(\widetilde{Z})$. Later in this section, we will apply this construction to $Z_M = \widetilde{\OO}_M$ and $Z = G\times^Q(\widetilde{\OO}_M\times \fq^\perp)$.
\end{rmk}

\begin{lemma}\label{lem:grind}
There is a $G$ and $\C^\times$-equivariant isomorphism of coherent sheaves on $\widetilde{Z}$
$$\gr \Ind^G_M \mathcal{E} \xrightarrow{\sim} \Ind^G_M(\gr \cE).$$
\end{lemma}
    
\begin{proof}
The proof has three steps.

{\it Step 1}.  Recall that $\mu$ stands for the moment map $T^*(G/U)\times \widetilde{X}_M\rightarrow \fm^*$. We claim that the natural epimorphism
\begin{equation}\label{eq:product_assoc_graded1}
[\mathcal{O}_{T^*(G/U)}\boxtimes \gr \mathcal{E}]/[\mathcal{O}_{T^*(G/U)}\boxtimes \gr \mathcal{E}]\mu^*(\mathfrak{m}) \twoheadrightarrow \gr 
(\cE'/\cE'\Phi_{\beta_2}(\fm)) 
\end{equation}
induced by (\ref{eq:product_assoc_graded}) is an isomorphism. 

Choose a basis $x_1,\ldots,x_k\in \mathfrak{m}$.
The action of $M$ on $T^*(G/U)$ is free, and hence the elements  $\mu_{T^*(G/U)}^*(x_i)$ form a regular sequence in the sheaf of algebras $\mathcal{O}_{T^*(G/U)}$. It follows that the elements  $\mu^*(x_i)$ form a regular sequence for the sheaf of modules $\mathcal{O}_{T^*(G/U)}\boxtimes \gr \mathcal{E}$. Hence the Koszul complex associated to this sheaf of modules and the specified elements is exact (in positive degrees). The Chevalley-Eilenberg complex for the Lie algebra $\fm$ acting on $\cE'$ via $\Phi_{\beta_2}$ is a filtered deformation of the Koszul complex above. Since the filtration on $\mathcal{E}'$ is complete and separated, the exactness of the Koszul complex implies the exactness of the Chevalley-Eilenberg complex and the claim that the associated graded of the 0th homology of the latter is the 0th homology of the former (the argument is by induction on the filtration degree, compare to the proof of Lemma \ref{lem:quant_pushforward}). So 
(\ref{eq:product_assoc_graded1}) is an isomorphism. 

{\it Step 2}. Since $M$ is reductive, taking $\gr$ commutes with taking $M$-invariants. So
$$
\gr\Ind^G_M \mathcal{E} = \gr\left(q_*[\mathcal{E}'/\mathcal{E}'\Phi_{\beta_2}(\fm)]\right)^M=\left(\gr q_*[\mathcal{E}'/\mathcal{E}'\Phi_{\beta_2}(\fm)]\right)^M$$
Since $q$ is $\CC^\times$-equivariant, taking $\gr$ commutes with $q_*$. Thanks to
(\ref{eq:product_assoc_graded1}), we see that
\begin{equation}\label{eq:product_assoc_graded2}
\gr\Ind^G_M \mathcal{E}\xrightarrow{\sim}
(q_*\left([\mathcal{O}_{T^*(G/U)}\boxtimes \gr \mathcal{E}]/[\mathcal{O}_{T^*(G/U)}\boxtimes \gr \mathcal{E}]\mu^*(\mathfrak{m})\right))^M.
\end{equation}

{\it Step 3}. It remains to identify the target of (\ref{eq:product_assoc_graded2}) with $\Ind^G_M(\gr\cE)$. First of all, 
$[\mathcal{O}_{T^*(G/U)}\boxtimes \gr \mathcal{E}]/[\mathcal{O}_{T^*(G/U)}\boxtimes \gr \mathcal{E}]\mu^*(\mathfrak{m})$ is nothing but the restriction of $\mathcal{O}_{T^*(G/U)}\boxtimes \gr \mathcal{E}$ to $\mu^{-1}(0)$. This restriction is identified with the pullback 
of $\gr\cE$ to $\mu^{-1}(0)=G\times^U (\widetilde{X}_M\times \fq^\perp)$ under the natural projection $G\times^U (\widetilde{X}_M\times\fq^\perp)\twoheadrightarrow \widetilde{X}_M$. This is a $G\times M$-equivariant sheaf. It is easy to see that it coincides with  $q^*\Ind^G_M\gr\cE$.
So the target of (\ref{eq:product_assoc_graded2}) is 
$\Ind^G_M\gr\cE$. This completes the proof. 
\end{proof}

To complete our construction of parabolic induction, define
$$\Ind^G_M\cB := \Gamma(\widetilde{Z}, \Ind^G_M\mathcal{E}).$$
By Proposition \ref{prop:quantizationparaminduction}, there are isomorphisms $\Ind^G_M \cA_{\beta_i}^{\widetilde{X}_M}\simeq \cA_{\beta_i}^{\widetilde{X}}$ for $i=1,2$. Hence, $\Ind^G_M\cB$ is an $\cA_{\beta_1}^{\widetilde{X}}$-$\cA_{\beta_2}^{\widetilde{X}}$ bimodule. It inherits a filtration from $\Ind^G_M \mathcal{E}$. This filtration is compatible with the filtrations on $\cA^{\widetilde{X}}_{\beta_1},\cA^{\widetilde{X}}_{\beta_2}$. There is a $G$-action on $\Ind^G_M \cB$ which preserves this filtration, making  $\Ind^G_M \cB$  a $G$-equivariant $\cA_{\beta_1}^{\widetilde{X}}$-$\cA_{\beta_2}^{\widetilde{X}}$ bimodule.

By Lemma \ref{lem:grind}, there is an injective homomorphism of $G$-equivariant $\CC[\widetilde{X}]$-modules
\begin{equation}\label{eq:grind to indgr}
  \gr(\Ind^G_M \cB) \hookrightarrow \Gamma(\widetilde{Z}, \gr(\Ind^G_M \mathcal{E})) \simeq \Gamma(\widetilde{Z},\Ind^G_M\gr(\mathcal{E})).  
\end{equation}
The latter module is finitely-generated over $\CC[\widetilde{X}]$. Thus, $\Ind^G_M\cB \in \HC^G(\cA_{\beta_1}^{\widetilde{X}}, \cA_{\beta_2}^{\widetilde{X}})$.

\begin{prop}\label{prop:twistedinductiondefined}
 $\Ind^G_M$ defines a left exact functor
\begin{equation}\label{eq:ind_functor}
    \Ind^G_M: \HC^M(\cA_{\beta_1}^{\widetilde{X}_M},\cA_{\beta_2}^{\widetilde{X}_M}) \to \HC^G(\cA_{\beta_1}^{\widetilde{X}},\cA_{\beta_2}^{\widetilde{X}}).
\end{equation}    
    This functor sends $\HC^M_{\partial}(\cA_{\beta_1}^{\widetilde{X}_M},\cA_{\beta_2}^{\widetilde{X}_M})$ to $\HC^G_{\partial}(\cA_{\beta_1}^{\widetilde{X}},\cA_{\beta_2}^{\widetilde{X}})$ and descends to a functor between the quotient categories
    $$\Ind^G_M: \overline{\HC}^M(\cA_{\beta_1}^{\widetilde{X}_M},\cA_{\beta_2}^{\widetilde{X}_M}) \to \overline{\HC}^G(\cA_{\beta_1}^{\widetilde{X}},\cA_{\beta_2}^{\widetilde{X}}).$$
\end{prop}

\begin{proof}
That $\Ind_M^G$ defines a functor is immediate from the construction. Next we show that this functor is left exact. 
Let 
$$0\rightarrow \cB_1\rightarrow \cB_2
\rightarrow \cB_3\rightarrow 0$$
be an exact sequence in $\HC^M(\cA^{\widetilde{X}_M}_{\beta_1}, 
\cA^{\widetilde{X}_M}_{\beta_2})$. Equip the terms with good filtrations so that 
there is an exact sequence of graded modules
$$0\rightarrow \gr\cB_1\rightarrow \gr\cB_2\rightarrow \gr\cB_3\rightarrow 0.$$
The induction functor for coherent sheaves is manifestly exact. Thus, there is an exact sequence of Harish-Chandra $\Ind^G_M \calD^{\widetilde{X}}_{\beta_1}$-
$\Ind^G_M \calD^{\widetilde{X}}_{\beta_2}$-
bimodules
$$0\rightarrow \Ind^G_M\cE_1\rightarrow 
\Ind^G_M\cE_2\rightarrow \Ind^G_M \cE_3\rightarrow 0.$$
Since $\Gamma$ is left exact so is
$$\Ind^G_M: \HC^M(\cA_{\beta_1}^{\widetilde{X}_M},\cA_{\beta_2}^{\widetilde{X}_M}) \to \HC^G(\cA_{\beta_1}^{\widetilde{X}},\cA_{\beta_2}^{\widetilde{X}}).$$
Next we show that  
\begin{equation}\label{eq:boundary_cat_incl}
\Ind^G_M\left(
\HC^M_{\partial}(\cA_{\beta_1}^{\widetilde{X}_M},\cA_{\beta_2}^{\widetilde{X}_M})\right) \subset\HC^G_{\partial}(\cA_{\beta_1}^{\widetilde{X}},\cA_{\beta_2}^{\widetilde{X}}) 
\end{equation}

Since $\overline{\rho}$ is proper, $\mathcal{V}(\Ind^G_M \cB)\subseteq
\overline{\rho}\mathcal{V}(\Ind^G_M \cE)$ (cf. \cite[Lemma 2.18]{BezLosev}). By (\ref{eq:grind to indgr}), there is an inclusion
$$\mathcal{V}(\Ind^G_M\cE) \subseteq \overline{\rho} [\mathrm{Supp}(\Ind^G_M \gr(\mathcal{E}))] \subseteq \overline{\rho}(G \times^Q (\mathcal{V}(\mathcal{E}) \times \fq^{\perp})) = \overline{\rho}(G \times^Q (\mathcal{V}(\cB) \times \fq^{\perp})).$$
Since $\overline{\rho}: \widetilde{Z} \to \widetilde{X}$ maps $\widetilde{Z} - \widetilde{\OO}$ to $\widetilde{X} - \widetilde{\OO}$, the condition $\mathcal{V}(\cB) \subseteq \widetilde{X}_M - \widetilde{\mathbb{O}}_M$ implies $\mathcal{V}(\Ind^G_M\cB) \subseteq \widetilde{X} - \widetilde{\mathbb{O}}$.
This proves (\ref{eq:boundary_cat_incl}).

Finally, we show that $\Ind^G_M$ descends to the quotient categories. Let $\theta_M$ denote the quotient functor 
$\HC^M(\cA_{\beta_1}^{\widetilde{X}_M},\cA_{\beta_2}^{\widetilde{X}_M})\twoheadrightarrow
\overline{\HC}^M(\cA_{\beta_1}^{\widetilde{X}_M},\cA_{\beta_2}^{\widetilde{X}_M})$ (and define $\theta_G$ similarly). Note that $\theta_M$ has a right adjoint (and left inverse) $\theta_M^*$: first microlocalize to $\widetilde{\OO}_M$ and then take global sections. Define the functor 
\begin{equation}\label{eq:induction_quot_cat}
\Ind^G_M: \overline{\HC}^M(\cA_{\beta_1}^{\widetilde{X}_M},\cA_{\beta_2}^{\widetilde{X}_M}) \to \overline{\HC}^G(\cA_{\beta_1}^{\widetilde{X}},\cA_{\beta_2}^{\widetilde{X}})
\end{equation} 
as the composition $\theta_G \circ \Ind^G_M \circ \theta_M^*$. 
Since the cokernel of the adjunction unit $\operatorname{id}\rightarrow \theta_M^*\circ \theta_M$ lies in $\HC_{\partial}^M(\cA_{\beta_1}^{\widetilde{X}_M},\cA_{\beta_2}^{\widetilde{X}_M})$, we see that $\theta_G\circ \Ind^G_M\simeq \Ind^G_M\circ \theta_M$. Hence, (\ref{eq:ind_functor})
descends to (\ref{eq:induction_quot_cat}).
\end{proof}

Let $\chi \in \fX(M)$ and consider the one-dimensional bimodule $\CC(\chi,0)$, see (\ref{eqn:1dimbimod}). If $\cB \in \HC^M(\cA_{\beta}^{\widetilde{X}_M})$, then $\cB \otimes \CC(\chi,0) \in \HC^M(\cA_{\beta + \chi}^{\widetilde{X}_M}, \cA_{\beta}^{\widetilde{X}_M})$. This defines an equivalence
$$\otimes \ \CC(\chi,0): \HC^M(\cA_{\beta}^{\widetilde{X}_M}) \to \HC^M(\cA_{\beta + \chi}^{\widetilde{X}_M}, \cA_{\beta}^{\widetilde{X}_M})$$
with inverse $\otimes \ \CC(-\chi,0)$.

For the remainder of this section, we will impose the following condition on $\chi$:
\begin{equation}\label{eq:conjugacycondition}\beta \text{ and } \beta+\chi \text{ are conjugate under } W^{\widetilde{X}}\end{equation}
Under this condition, $\cA_{\beta}^{\widetilde{X}} \simeq \cA_{\beta+\chi}^{\widetilde{X}}$ as Hamiltonian quantizations and thus $\Ind^G_M$ defines a functor
$$\Ind^G_M: \HC^M(\cA_{\beta + \chi}^{\widetilde{X}_M}, \cA_{\beta}^{\widetilde{X}_M}) \to \HC^G(\cA_{\beta}^{\widetilde{X}})$$
Consider the composition
$$\Ind^G_M[\chi]: \HC^M(\cA_{\beta}^{\widetilde{X}_M}) \overset{\otimes \CC(\chi,0)}{\xrightarrow{\sim}} \HC^M(\cA_{\beta + \chi}^{\widetilde{X}_M}, \cA_{\beta}^{\widetilde{X}_M}) \overset{\Ind^G_M}{\to} \HC^G(\cA_{\beta}^{\widetilde{X}}).$$
By Proposition \ref{prop:twistedinductiondefined}, this functor descends to a functor (still denoted by $\Ind^G_M[\chi]$)
$$\Ind^G_M[\chi]: \overline{\HC}^M(\cA_{\beta}^{\widetilde{X}_M}) \to \overline{\HC}^G(\cA_{\beta}^{\widetilde{X}})$$
These functors should be viewed as `twisted' versions of parabolic induction. We will now give $\Ind^G_M[\chi]$ an alternative description involving restriction functors. 

Choose $x \in \widetilde{\OO}$ and let $\Omega = \pi_1^G(\widetilde{\OO}) \simeq G_x/G_x^{\circ}$. By Proposition \ref{prop:classificationeqvtbimods}, there is a full monoidal embedding
$$\bullet_{\dagger}: \overline{\HC}^G(\cA_{\beta}^{\widetilde{X}}) \hookrightarrow \Omega\modd.$$
Consider the group homomorphisms $\mathcal{L}_{\widetilde{\OO}}: \fX(G_x) \to \Pic(\widetilde{\OO})$ and $\mathcal{L}_{G/Q}: \fX(M) \to \Pic(G/Q)$, as well as the embedding $\fX(\Omega) \hookrightarrow \fX(G_x)$ induced from the quotient homomorphism $G_x \twoheadrightarrow \Omega$. Let $\widetilde{\psi}$ denote the composition $\widetilde{\OO}\hookrightarrow\widetilde{Z}=G\times^Q(\widetilde{X}_M\times \fq^\perp) \twoheadrightarrow G/Q$.

\begin{lemma}\label{lem:mutilde}
Assume $\chi$ satisfies (\ref{eq:conjugacycondition}). There is a uniquely defined character $\varphi_M^G(\chi) \in \fX(\Omega)$ such that
$$\widetilde{\psi}^*\mathcal{L}_{G/Q}(\chi) \simeq \mathcal{L}_{\widetilde{\mathbb{O}}}(\varphi_M^G(\chi))$$
as $G$-equivariant line bundles on $\widetilde{\OO}$.
\end{lemma}

\begin{proof}
    The proof has two steps. First, we show that the line bundle $\widetilde{\psi}^*\mathcal{L}_{G/Q}(\chi)$ has trivial first Chern class. Then we use this to deduce the statement of the lemma.

        \emph{Step 1.} We can assume that $G$ is semisimple. The projection $\overline{\pi}:\widetilde{Z} \to G/Q$ induces a homomorphism
$$\fX(\fm) \simeq H^2(G/Q,\CC) \overset{\overline{\pi}^*}{\to} H^2(\widetilde{Z},\CC).$$
Since taking Chern classes commutes with pullbacks, the Chern class of $\overline{\pi}^*\mathcal{L}_{G/Q}(\chi)$ corresponds to $\chi$ under the isomorphism above. Thus, the Chern class of $\widetilde{\psi}^*\mathcal{L}_{G/Q}(\chi)=[\overline{\pi}^*\mathcal{L}_{G/Q}(\chi)]|_{\widetilde{\mathbb{O}}}$ corresponds to the image of $\chi$ under the composition $\fX(\fm) \to H^2(\widetilde{Z},\CC) \to H^2(\widetilde{\mathbb{O}},\CC)$. We claim that this map coincides with the composition
\begin{equation}\label{eq:c1trivial}\fX(\fm) \hookrightarrow \fX(\fl) \xrightarrow{\sim} \fP^{\widetilde{X}} \twoheadrightarrow H^2(\widetilde{\mathbb{O}},\CC),\end{equation}
where middle map is the isomorphism $\eta$, see (\ref{eq:defofeta}), and the final map is projection onto $\fP_0^{\widetilde{X}} \simeq H^2(\widetilde{\OO},\CC)$ (the latter identification is by Lemma \ref{lem:computeH2}). Recall from Proposition \ref{prop:independence} that $\eta$ is well-defined defined up to the $W^{\widetilde{X}}$-action on $\fP^{\widetilde{X}}$. Since $W^{\widetilde{X}}$ acts trivially on  $\fP_0^{\widetilde{X}}$, we see that $\fX(\fl)\twoheadrightarrow H^2(\widetilde{\OO},\CC)$ is independent of $P$.

To show that the composition $\fX(\fm) \to H^2(\widetilde{Z},\CC) \to H^2(\widetilde{\mathbb{O}},\CC)$ coincides with (\ref{eq:c1trivial}), pick a $\QQ$-terminalization $\widetilde{Y}_M$ of $\widetilde{X}_M$. It has the form $M\times^{P\cap M} (\widetilde{X}_L\times (\fm^*\cap \fp^\perp))$ for some choice of a parabolic subgroup $P$ contained in $Q$.
Consider the $\QQ$-terminalization $\widetilde{Y}:=G\times^Q(\widetilde{Y}_M\times \fq^\perp)\simeq G\times^P(\widetilde{X}_L\times \fp^\perp)$ of $\widetilde{X}$ (the isomorphism was established in Lemma \ref{lem:compatibility1}). Let $\pi': \widetilde{Y}\twoheadrightarrow G/Q$
denote the projection. The composition $\fX(\fm) \hookrightarrow \fX(\fl) \xrightarrow{\sim} \fP^{\widetilde{X}}$ coincides with the map $\pi'^*: H^2(G/Q,\CC)\to H^2(\widetilde{Y}, \CC)$. Since the natural map $\widetilde{Y}\to \widetilde{Z}$ is an isomorphism over $\widetilde{\OO}$, the restrictions of $\overline{\pi}^*\mathcal{L}_{G/Q}(\chi)$ and $\pi'^*\mathcal{L}_{G/Q}(\chi)$ to $\widetilde{\OO}$ are isomorphic. Thus the composition $\fX(\fm) \to H^2(\widetilde{Z},\CC) \to H^2(\widetilde{\mathbb{O}},\CC)$ indeed coincides with (\ref{eq:c1trivial}).

Since $\beta+\chi$ and $\beta$ are conjugate under $W^{\widetilde{X}}$, and $W^{\widetilde{X}}$ acts trivially on $H^2(\widetilde{\mathbb{O}},\CC) \hookrightarrow \fP^{\widetilde{X}} \simeq \fX(\fl)$, the image of $\chi$ under (\ref{eq:c1trivial}) is $0$. Thus, the Chern class of the line bundle $\widetilde{\psi}^*\mathcal{L}_{G/Q}(\chi)$ is trivial, as asserted. 

\emph{Step 2.} By Step 1,  
\begin{equation}\label{eq:trivial_c1}
c_1(\widetilde{\psi}^*\mathcal{L}_{G/Q}(\chi))=0.\end{equation}
We will use this to deduce the statement of the lemma. For any algebraic group $H$, the torsion subgroup $\fX(H)_{\mathrm{tor}} \subset \fX(H)$ coincides with $\fX(H/H^{\circ})$. In particular, $\fX(G_x)_{\mathrm{tor}} = \fX(G_x/G_x^{\circ}) \simeq \fX(\Omega)$. Thus, $\mathcal{L}_{\widetilde{\mathbb{O}}}: \fX(G_x) \simeq \Pic^G(\widetilde{\mathbb{O}})$ restricts to an isomorphism
\begin{equation}\label{eq:tor}
\mathcal{L}_{\widetilde{\mathbb{O}}}: \fX(\Omega) \simeq \Pic^G(\widetilde{\mathbb{O}})_{\mathrm{tor}}.\end{equation}
On the other hand, the Chern class map $c_1: \Pic(\widetilde{\mathbb{O}}) \to H^2(\widetilde{\mathbb{O}},\CC)$ induces a linear isomorphism $\Pic(\widetilde{\mathbb{O}}) \otimes_{\ZZ} \CC \simeq H^2(\widetilde{\mathbb{O}},\CC)$, see Lemma \ref{lem:computeH2}. Hence, $\ker{c_1} = \Pic(\widetilde{\mathbb{O}})_{\mathrm{tor}}$. By (\ref{eq:trivial_c1}), the line bundle $\overline{\pi}^*\mathcal{L}_{G/Q}(\chi)|_{\widetilde{\mathbb{O}}}$ belongs to $\ker{c_1}$, and is $G$-equivariant by construction. Thus by (\ref{eq:tor}), there is a unique character $\varphi_M^G(\chi) \in \fX(\Omega)$ such that $\widetilde{\psi}^*\mathcal{L}_{G/Q}(\chi) \simeq \mathcal{L}_{\widetilde{\mathbb{O}}}(\varphi_M^G(\chi))$.
\end{proof}

For the next lemma, let $P_M:=P\cap M$, a parabolic in $M$. Consider the partial resolution $\widetilde{Y}_M:=M\times^{P_M} (\widetilde{X}_L\times \fp_M^{\perp}) \to \widetilde{X}_M$ and the projection map ${\pi}_M: \widetilde{Y}_M \to M/P_M$. There is a $G$-equivariant identification $\widetilde{Y}=G\times^Q(\widetilde{Y}_M\times \fq^\perp)$, see Lemma \ref{lem:compatibility1}.

Set $\Omega' = \pi_1^M(\widetilde{\OO}_M)$. Recall the surjective homomorphism $f: \Omega \to \Omega'$ defined in Lemma \ref{lem:mappi1}, and let $\Omega_0=\ker{f}$. Write $\widehat{\mathbb{O}}_M \to \widetilde{\OO}_M$ for the universal $M$-equivariant cover---its Galois group is $\Omega'$. Set 
$$\widecheck{\mathbb{O}}:= \mathrm{Bind}^G_M \widehat{\mathbb{O}}_M, \qquad  \widecheck{Z}=G\times^Q(\widecheck{\OO}\times \fq^\perp).$$
Observe that $\widecheck{\OO}$ is a Galois cover of $\widetilde{\OO}$ with Galois group $\Omega'$.
So by the construction of $f:\Omega\twoheadrightarrow \Omega'$ in Section \ref{subsec:inductionpi1}, we have $\widecheck{\OO}\simeq 
\widehat{\OO}/\Omega_0$, where $\widehat{\OO}$ is the universal $G$-equivariant cover of $\OO$.
Finally, let $\widecheck{\psi}$ denote the composition $\widecheck{\OO}\hookrightarrow \widecheck{Z}\twoheadrightarrow G/Q$.

\begin{lemma}\label{lem:coherentind}
The following are true:
\begin{itemize}
    \item[(i)] Let $\chi \in \fX(L)$. Then there is an isomorphism of $G$-equivariant lines bundles on $\widetilde{Y}$
        $$\pi^*\calL_{G/P}(\chi)\simeq \Ind_M^G[{\pi}_M^*\calL_{M/P_M}(\chi)].$$
        \item[(ii)] Let $\chi' \in \fX(\Omega')$. Then there is an isomorphism of $G$-equivariant line bundles on $\widetilde{\OO}$
        \begin{equation}\label{eq:line_bundle_iso1}
        \calL_{\widetilde{\OO}}(f^*\chi')\simeq \Ind_M^G[\calL_{\widetilde{\OO}_M}(\chi')]|_{\widetilde{\OO}},\end{equation}
        and also an isomorphism of $G\times \Omega'$-equivariant line bundles on $\widecheck{\OO}$
        \begin{equation}\label{eq:line_bundle_iso2}
        \chi'\otimes \mathcal{O}_{\widecheck{\OO}}\simeq \Ind_M^G[\calL_{\widehat{\OO}_M}(\chi')]|_{\widecheck{\OO}}.\end{equation}
\end{itemize}
\end{lemma}

\begin{proof}
We first prove (i). We can assume that $G$ is semisimple and simply connected.
Note that both sides are line bundles on $\widetilde{Y}$ so it suffices to show that their classes in $\operatorname{Pic}(\widetilde{Y})$ coincide. Recall, Proposition \ref{prop:descriptionofpic}, that $\pi^*: \Pic(G/P)\xrightarrow{\sim} 
\Pic(\widetilde{Y})$ is an isomorphism. Its inverse is the restriction to the zero section $G/P\subset \widetilde{Y}$. So it is enough to show that 
$$\pi^*\calL_{G/P}(\chi)|_{G/P}\simeq \Ind_M^G[{\pi}_M^*\calL_{M/P_M}(\chi)]|_{G/P}.$$
The left hand side is $\mathcal{L}_{G/P}(\chi)$.
It follows easily from the construction of $\operatorname{Ind}^G_M$ that the right hand side is the same.

We proceed to proving (ii).  By construction, the left (resp. right) side of (\ref{eq:line_bundle_iso1})
is obtained from the left (resp. right) side of (\ref{eq:line_bundle_iso2}) by equivariant descent for the action of $\Omega'$. So (\ref{eq:line_bundle_iso1}) and (\ref{eq:line_bundle_iso2}) are equivalent. To prove (\ref{eq:line_bundle_iso2}) we note 
    $$\Ind_M^G[\calL_{\widehat{\OO}_M}(\chi')]|_{\widecheck{\OO}}\simeq \chi'\otimes \Ind_M^G[\mathcal{O}_{\widehat{\OO}_M}]|_{\widecheck{\OO}}= \chi'\otimes \cO_{\widecheck{\OO}}.$$
This completes the proof of (ii).
\end{proof}

The next proposition is the main result of this section. In it, we describe parabolic induction of Harish-Chandra bimodules in terms of the equivalence of \cref{prop:classificationeqvtbimods}. 

\begin{prop}\label{prop:twinductiondagger}
Recall that $\widetilde{\OO}$ is birationally induced from $\widetilde{\OO}_M$, and we write $\Omega',\Omega$ for $\pi_1^M(\widetilde{\OO}_M),\pi^G_1(\widetilde{\OO})$, respectively. 
Assume $\chi$ satisfies (\ref{eq:conjugacycondition}), and define $\varphi_M^G(\chi) \in \fX(\Omega)$ as in Lemma \ref{lem:mutilde}. Then the following diagram of functors commutes:

\begin{center}
    \begin{tikzcd}
      \overline{\HC}^M(\cA_{\beta}^{\widetilde{X}_M}) \ar[d,hookrightarrow,"\bullet_{\dagger}"] \ar[rr, "{\Ind^G_M[\chi]}"]  & & \overline{\HC}^G(\cA_{\beta}^{\widetilde{X}}) \ar[d,hookrightarrow,"\bullet_{\dagger}"] \\
      \Omega'\modd \ar[rr,"\varphi_M^G(\chi) \otimes  f^*(\bullet)"] & & \Omega\modd
    \end{tikzcd}
\end{center}

\end{prop}

\begin{proof}
Since all categories in question are semisimple, it is enough to show that  
\begin{equation}\label{eq:ind_iso_objects}
\varphi_M^G(\chi) \otimes  f^*(\cB_\dagger)\simeq \left(\Ind^G_M[\chi](\cB)\right)_\dagger
\end{equation}
 for all $\cB \in \HC(\cA_{\beta}^{\widetilde{X}_M})$. Set $V := \cB_{\dagger} \in \Omega'\modd$. By \cref{lem:determineddagger} , it suffices to show that
\begin{equation}\label{eqn:inductiondagger2}
(\gr \Ind^G_M [\cB\otimes \CC(\chi, 0)])|_{\widetilde{\mathbb{O}}} \simeq ( \varphi^G_M(\chi)\otimes f^*V \otimes \cO_{\widehat{\mathbb{O}}})^{\Omega}
\end{equation}
as $G$-equivariant coherent sheaves, where $\widehat{\mathbb{O}} \to \widetilde{\mathbb{O}}$ denotes the universal $G$-equivariant cover. 

Using \cref{lem:grind}, we have 
$$(\gr \Ind^G_M [\cB\otimes \CC(\chi, 0)])|_{\widetilde{\mathbb{O}}} \simeq (\Ind^G_M [\gr \cB\otimes \CC(\chi, 0)])|_{\widetilde{\mathbb{O}}}.$$ 

Let $i:\widetilde{\OO}_M\to \widetilde{X}_M$ be the inclusion. Since $\widetilde{\OO}\subset G\times^Q(\widetilde{\OO}_M\times \fq^\perp)$, there is an isomorphism of $G$-equivariant coherent sheaves on $\widetilde{\OO}$

$$(\Ind^G_M [\gr \cB\otimes \CC(\chi, 0)])|_{\widetilde{\mathbb{O}}}\simeq (\Ind^G_M i^*[\gr \cB\otimes \CC(\chi, 0)])|_{\widetilde{\mathbb{O}}}.$$ 
By (\ref{eq:local_system_iso2}), $i^*\gr \cB\simeq (V\otimes \cO_{\widehat{\OO}_M})^{\Omega'}$. Thus, we get  
\begin{equation}\label{eq:commrephrasing}
(\gr \Ind^G_M [\cB\otimes \CC(\chi, 0)])|_{\widetilde{\mathbb{O}}} \simeq \Ind^G_M (V\otimes \cO_{\widehat{\OO}_M} \otimes \CC(\chi, 0)^{\Omega'})|_{\widetilde{\mathbb{O}}}
=\Ind^G_M([\chi\otimes V \otimes \cO_{\widehat{\mathbb{O}}_M}]^{\Omega'})|_{\widetilde{\mathbb{O}}}.
\end{equation}
We proceed to computing $\Ind^G_M(\chi\otimes V \otimes \cO_{\widehat{\mathbb{O}}_M})$. 
%
%
We claim that for any representation $V$ of $\Omega'$ we have an isomorphism
\begin{equation}\label{eq:another_induced_iso}
V\otimes \widecheck{\psi}^*\calL_{G/Q}(\chi)\simeq 
\Ind^G_M(\chi\otimes V \otimes \cO_{\widehat{\mathbb{O}}_M})|_{\widecheck{\OO}}
\end{equation}
of $G\times \Omega'$-equivariant coherent sheaves on $\widecheck{\OO}$ (recall that $\widecheck{\psi}$ is the morphism $\widecheck{\OO} \hookrightarrow \widecheck{Z} \to G/Q$). Note that $\Ind^G_M(\chi\otimes V \otimes \cO_{\widehat{\mathbb{O}}_M})\simeq V\otimes \Ind^G_M(\chi \otimes \cO_{\widehat{\mathbb{O}}_M})$, so we can assume that $V$ is trivial. Let $\overline{\pi}$
denote the projection $\widetilde{Z}\twoheadrightarrow G/Q$.  Now (\ref{eq:another_induced_iso}) 
follows from $\overline{\pi}^*\calL_{G/Q}(\chi)\simeq \Ind^G_M(\chi\otimes \cO_{\widetilde{\OO}_M})$ (a tautological special case of (i) of Lemma \ref{lem:coherentind}): to get 
(\ref{eq:another_induced_iso}) from the latter isomorphism we pull back along $\widecheck{\OO}\rightarrow \widetilde{Z}$.

Let $\calL_{\widecheck{\OO}}(\varphi^G_M(\chi))$ denote the $G\times \Omega'$-equivariant line bundle on $\widecheck{\OO}$ corresponding to $\varphi^G_M(\chi)$, this is the pullback of $\calL_{\widetilde{\OO}}(\varphi^G_M(\chi))$. By Lemma \ref{lem:mutilde}, we have an isomorphism $\calL_{\widecheck{\OO}}(\varphi^G_M(\chi))\simeq \widecheck{\psi}^*\calL_{G/Q}(\chi)$ of $G\times \Omega'$-equivariant line bundles. Combining this with (\ref{eq:another_induced_iso}) we see that 
$$\Ind^G_M(\chi\otimes V \otimes \cO_{\widehat{\mathbb{O}}_M})|_{\widecheck{\OO}}\simeq V\otimes \calL_{\widecheck{\OO}}(\varphi^G_M(\chi)),$$
equivalently, by $\Omega'$-equivariant descent,
$$\Ind^G_M([\chi\otimes V \otimes \cO_{\widehat{\mathbb{O}}_M}]^{\Omega'})|_{\widetilde{\OO}}\simeq [V\otimes \calL_{\widecheck{\OO}}(\varphi^G_M(\chi))]^{\Omega'}\simeq  [f^*(V)\otimes \varphi^G_M(\chi)\otimes \cO_{\widehat{\OO}}]^\Omega.$$
Combining the composed isomorphism with 
(\ref{eq:commrephrasing}) we arrive at 
(\ref{eqn:inductiondagger2}).
\end{proof}

\begin{cor} \label{cor:propsofbimodind}
Assume $\chi$ satisfies (\ref{eq:conjugacycondition}). Then the functor
$$\Ind^G_M[\chi]: \overline{\HC}^M(\cA_{\beta}^{\widetilde{X}_M}) \to \overline{\HC}^G(\cA_{\beta}^{\widetilde{X}})$$ 
enjoys the following properties:
\begin{itemize}
    \item[(i)] Suppose $L\subset M$ is a Levi subgroup and choose a character $\chi' \in \fX(L)$ such that 
    \begin{itemize}
        \item[$\bullet$] $\beta$ and $\beta+\chi'$ are conjugate under $W^{\widetilde{X}_M}$,
        \item[$\bullet$] $\beta$ and $\beta+\chi+\chi'$ are conjugate under $W^{\widetilde{X}}$.
    \end{itemize}
    Then there is a natural isomorphism
    $$\Ind^G_L[\chi+\chi'] \simeq \Ind^G_M[\chi] \circ \Ind^M_L[\chi'].$$
    \item[(ii)] $\Ind^G_M[\chi]$ is conservative, i.e. takes non-isomorphic objects to non-isomorphic objects.
    \item[(iii)] $\Ind^G_M[\chi]$ takes irreducibles to irreducibles.
\end{itemize}
\end{cor}

\begin{proof}
By Proposition \ref{prop:twinductiondagger}, it suffices to establish the analogous properties for the functor $\varphi^G_M(\chi) \otimes f^*(\bullet): \Omega'\modd \to \Omega\modd$
(and the similar functors in (i)).  Properties (ii) and (iii) follow from the surjectivity of $f$. 
 
 For (i), we need some notation. Let $\Omega'' =\pi_1^L(\widetilde{\OO}_L)$ and let
 $g:\Omega'\twoheadrightarrow \Omega''$ denote the analog of $f:\Omega\twoheadrightarrow \Omega'$ for $L \subset M$. Then $g\circ f$ is the analog of $f$ for $L \subset G$. For (i), it suffices to show that
\begin{equation}\label{eq:varphi1}\varphi_L^G(\chi+\chi') \otimes (g\circ f)^*(\cB_{\dagger}) = \varphi_M^G(\chi) \otimes f^*[\varphi_L^M(\chi') \otimes g^*(\cB_{\dagger})],\end{equation}
We will need the following facts about the maps $\varphi_\bullet^\bullet$:
\begin{enumerate}
    \item $\varphi_L^G(\chi) = \varphi_M^G(\chi)$.
    \item $\varphi_L^G(\chi+\chi') = \varphi_L^G(\chi) \otimes \varphi_L^G(\chi')$.
    \item $\varphi_L^G(\chi') = f^*\varphi_L^M(\chi')$.
\end{enumerate}
(1) follows from the isomorphism $\pi^*\mathcal{L}_{G/P}(\chi) \simeq \overline{\pi}^*\mathcal{L}_{G/Q}(\chi)$. (2) follows from the fact that $\pi^*\mathcal{L}_{G/P}^*(\bullet)|_{\overline{\OO}}$ is a homomorphism $\fX(L) \to \Pic^G(\widetilde{\OO})$. For (3) we argue as follows. By Lemma \ref{lem:coherentind}(i) (using the notation therein), there is an isomorphism $$\Ind^G_M(\pi_M^*\mathcal{L}_{M/P_M}(\chi')) \simeq \pi^*\mathcal{L}_{G/P}(\chi').$$ Restricting to $\widetilde{\OO}$ we get a $G$-equivariant isomorphism
$$\mathcal{L}_{\widetilde{\OO}}(\varphi_L^G(\chi')) \simeq \Ind^G_M[\mathcal{L}_{\widetilde{\OO}_M}(\varphi_L^M(\chi')]|_{\widetilde{\OO}}.$$
Now (3) follows from Lemma \ref{lem:coherentind}(ii). Equation (\ref{eq:varphi1}) is an immediate consequence of facts (1)-(3).
\end{proof}

\section{Unitarity of bimodules for Hamiltonian quantizations}

Let $\widetilde{\OO}$ be a $G$-equivariant nilpotent cover and choose a Hamiltonian quantization $\cA$ of $\CC[\widetilde{\OO}]$. As explained in Section \ref{subsec:classificationbimods}, there is a faithful embedding
$$\Phi^*: \HC^G(\cA) \to \HC^G(U(\fg)).$$
A bimodule $\B \in \HC^G(\cA)$ is said to be unitary if $\Phi^*(\B)$ is unitary. Proposition \ref{prop:inductionunitary} concerns the relation between unitarity and parabolic induction for Harish-Chandra $U(\fg)$-bimodules. In this section, we will prove analogous results for Harish-Chandra bimodules for Hamiltonian quantizations of nilpotent covers. 

Choose Levi subgroups $L \subset M \subset G$ and a birationally rigid $L$-equivariant nilpotent cover $\widetilde{\OO}_L$ such that $\widetilde{\OO} = \mathrm{Bind}^G_L \widetilde{\OO}_L$. Let $\widetilde{\OO}_M = \mathrm{Bind}^M_L \widetilde{\OO}_L$. Let $\beta \in \fX(\fl)$ and form the Hamiltonian quantizations $\cA_{\beta}^{\widetilde{X}_M}$ and $\cA_{\beta}^{\widetilde{X}}$. Choose a character $\chi \in \fX(M)$ such that $\beta+\frac{\chi}{2}$ and $\beta-\frac{\chi}{2}$ are conjugate under $W^{\widetilde{X}}$. If $\cB' \in \HC^M(\cA_{\beta}^{\widetilde{X}_M})$, then $\Ind^G_M[\cB' \otimes \CC(-\frac{\chi}{2},\frac{\chi}{2})] \in \HC^G(\cA_{\beta+\frac{\chi}{2}}^{\widetilde{X}})$.

\begin{prop}\label{prop:inductionunitarygeom}
Suppose $\B' \in \HC^M(\cA_{\beta}^{\widetilde{X}_M})$ is unitary, and $I_{\beta+\frac{\chi}{2}}(\widetilde{\mathbb{O}})= I_{\beta-\frac{\chi}{2}}(\widetilde{\mathbb{O}})$ as well as $I_{\beta}(\widetilde{\OO}_M)$ are  maximal ideals. Then $\Ind^G_M[\B' \otimes \CC(-\frac{\chi}{2},\frac{\chi}{2})]$ is unitary.
\end{prop}

\begin{proof}
Since $I_{\beta+\frac{\chi}{2}}(\widetilde{\mathbb{O}})= I_{\beta-\frac{\chi}{2}}(\widetilde{\mathbb{O}})$ and $I_{\beta}(\widetilde{\OO}_M)$ are  maximal ideals,  Proposition \ref{prop:HC_induction} is applicable.  So there is an isomorphism in $\HC^G(U(\fg))$
$$\Phi^*\Ind^G_M[\B' \otimes \CC(-\frac{\chi}{2},\frac{\chi}{2})] \simeq \Ind^G_M[\Phi^*(\B') \otimes \CC(-\frac{\chi}{2},\frac{\chi}{2})].$$
By assumption,  $\Phi^*(\cB')$ is unitary. Hence $\Phi^*(\B') \otimes \CC(-\frac{\chi}{2},\frac{\chi}{2})$ is unitary, since $\CC(-\frac{\chi}{2},\frac{\chi}{2})$ corresponds to an algebraic (thus, unitary) character of $L$. Now apply Proposition \ref{prop:inductionunitary}(1).
\end{proof}

For the next proposition, we will assume $M$ decomposes as a direct product of reductive subgroups $M_1,M_2 \subset M$. Let $\OO_1$ be a nilpotent $M_1$-orbit, and $\OO_M = \mathbb{O}_1 \times \{0\}$. Let $W_1$ and $W_2$ denote the Weyl groups of $M_1,M_2$, respectively. Choose $\nu \in \fX(\fm_2)$ and $\chi \in \fX(M_2)$ such that $\beta+\nu+\chi$ and $\beta+\nu$ are conjugate under $W^{\widetilde{X}}$. If $\cB' \in \HC^M(\cA_{\beta}^{X_M})$, then, $\Ind^G_M[\cB' \otimes \CC(\nu+\chi,\nu)] \in \HC^G(\cA_{\beta+\nu}^X)$.

To simplify notation, let
\begin{equation}\label{eq:lambdallambdar}\lambda_{\ell}(t) := \nu +\chi + \frac{1}{2}(t-1)(\chi+2\nu),  \qquad \lambda_r(t) :=  \nu + \frac{1}{2}(t-1)(\chi+2\nu), \qquad t \in [0,1].\end{equation}
\begin{prop}\label{prop:deformationgeom}
Suppose $\OO$ is birationally induced from $\OO_M$.
Let $\B' \in \HC^M(\cA_{\beta}^{X_M})$ be an irreducible unitary Harish-Chandra bimodule. Suppose, further, that
\begin{itemize}
    \item[(i)] The following ideals are maximal
    $$I_\beta(\OO_M), I_{\beta+\lambda_{\ell}(t)}(\mathbb{O}), I_{\beta+\lambda_r(t)}(\mathbb{O}), \qquad \forall \ t \in [0,1].$$
    \item[(ii)] There is an element $w \in W_2$ such that
    $$w \chi = \chi, \qquad w(\chi+2\nu) = -\overline{(\chi+2\nu)}.$$
\end{itemize}
Then $\Ind^G_M[\B' \otimes \CC(\nu+\chi, \nu)]$ is unitary.
\end{prop}

\begin{proof}
Consider the one-parameter family of Harish-Chandra bimodules
$$\B(t) := \Phi^*\Ind^G_M[\B' \otimes \CC(\lambda_{\ell}(t),\lambda_r(t))], \qquad t \in [0,1].$$
Note that $\B(0) = \Phi^*\Ind^G_M[\B' \otimes \CC(\frac{\chi}{2}, -\frac{\chi}{2})]$ and $\cB(1) = \Phi^*\Ind^G_M[\cB' \otimes \CC(\nu+\chi,\nu)]$. The first module is unitary by Proposition \ref{prop:inductionunitarygeom}. Thus by Proposition \ref{prop:inductionunitary} it suffices to show that $\B(t)$ is irreducible and Hermitian for every $t \in [0,1]$. 

\vspace{3mm}
$\underline{\B(t) \text{ is irreducible}}$: let $\cW$ be the $W$-algebra associated to $\mathbb{O}$ and let $R$ be the reductive part of the centralizer of $e \in \OO$. Recall the functors constructed in Section \ref{subsec:W}
$$\bullet_{\dagger}: \HC^G_{\overline{\OO}}(U(\fg)) \to \HC^R_{\mathrm{fin}}(\cW), \qquad \bullet^{\dagger}: \HC^R_{\mathrm{fin}}(\cW) \to \HC^G_{\overline{\OO}}(U(\fg)).$$
For every $t$, there is the adjunction unit homomorphism $\cB(t) \to (\cB(t)_{\dagger})^{\dagger}$. Since $I_{\beta+\lambda_{\ell}(t)}(\mathbb{O})$ is maximal, this map is an isomorphism, see \cref{lem:adjunctionmorphismiso}. From this, we will deduce the following implication
\begin{equation}\label{eq:daggerimplies}\cB(t)_{\dagger} \text{ irreducible} \implies \cB(t) \text{ irreducible}.\end{equation}
Indeed, suppose $\cB(t)_{\dagger}$ is irreducible and choose a nonzero subbimodule $C \subseteq \cB(t)$. There are inclusions $I_{\beta+\lambda_{\ell}(t)}(\mathbb{O}) \subseteq \mathrm{LAnn}_{U(\fg)}(\cB(t)) \subseteq \mathrm{LAnn}_{U(\fg)}(C)$. Since $I_{\beta+\lambda_{\ell}(t)}(\mathbb{O})$ is maximal, these inclusion are equalities. In particular, $\mathcal{V}(C) = V(I_{\beta+\lambda_{\ell}(t)}(\mathbb{O})) = \overline{\OO}$ and therefore $C_{\dagger}\neq 0$ by Proposition \ref{prop:propsofdagger}(ii). Since $B(t)_{\dagger}$ is irreducible, we must have $C_{\dagger} = B(t)_{\dagger}$. Applying Lemma \ref{lem:adjunctionmorphismiso} to $C$, we see that $C$ coincides with $(C_{\dagger})^{\dagger} = (\cB(t)_{\dagger})^{\dagger}$ under the natural isomorphism $\cB(t) \simeq (\cB(t)_{\dagger})^{\dagger}$. Hence, $C = \cB(t)$ and $\cB(t)$ is irreducible. 

In view of (\ref{eq:daggerimplies}), it is enough to show that $\cB(t)_{\dagger}$ is irreducible, for every $t \in [0,1]$. In fact, we will show that $F\cB(t)_{\dagger}$ is irreducible, for every $t \in [0,1]$, where $F: \HC^R_{\mathrm{fin}}(\cW) \to R\modd$ is the forgetful functor.

First, we will show that $F\cB(t)_{\dagger}$ is independent of $t$. If we choose a good filtration on $\cB'$, we get a good filtration on $(\B' \otimes \CC(\lambda_{\ell}(t),\lambda_r(t)))_\dagger$, and the associated graded sheaf $\gr(\B' \otimes \CC(\lambda_{\ell}(t),\lambda_r(t)))_\dagger$ is independent of $t$. Let $\mathcal{F} := \gr(\B' \otimes \CC(\lambda_{\ell}(t),\lambda_r(t)))_\dagger$. By the construction of $\Ind_M^G$ in Section \ref{subsec:bimodinductioncovers}, the bimodule $\cB(t)$ comes with a good filtration
such that $\gr \cB(t)|_{\OO}\simeq \Ind^G_M\mathcal{F}|_{\OO}$. By (ii) of Proposition 
\ref{prop:propsofdagger}, $\gr ([\Phi^*\cB(t)]_\dagger)$  is $R$-equivariantly isomorphic to the pullback {of $\Ind^G_M\mathcal{F}$ to the Slodowy slice $S$ attached to $\OO$}. Since $\gr ([\Phi^*\cB(t)]_\dagger)$ and $[\Phi^*\cB(t)]_\dagger$ are isomorphic as $R$-modules, $F\cB(t)_\dagger$ is independent of $t$, as asserted above. Thus, it is enough to prove the irreducibility of $F\cB(1)_{\dagger}$. 

Since $\OO$ is birationally induced from $\OO_M$, $\A_{\beta+\nu}$ is a quantization of $\CC[\OO]$. 
By Proposition \ref{prop:propsofdagger}, $\underline{J}:=I_{\beta+\nu}(\OO)_\dagger$ is an $R$-stable codimension $1$ ideal in $\cW$. 
 Recall the functor $\mathsf{B}_\bullet: \Omega\modd\to \HC^R_{\mathrm{fin}}(\cW/\underline{J}) $ from \cref{lem:twodaggers}. In our case, $\Omega=\pi_1^G(\OO)=R/R^\circ$. Hence $\mathsf{B}_\bullet$ is a category equivalence. 
 By (iii) of Corollary \ref{cor:propsofbimodind}, 
 $\mathcal{B}(1)$ is an irreducible object in 
 $\overline{\HC}(\cA_{\beta+\nu}^X)$. Hence by Proposition \ref{prop:classificationeqvtbimods}, 
 $\mathcal{B}(1)_\dagger$ is an irreducible object in $\Omega\operatorname{-mod}$. The $R$-module $F\cB(1)_\dagger$ is obtained from $\cB(1)_\dagger$ by pullback along the surjective homomorphism
  $R\twoheadrightarrow R/R^\circ=\Omega$. We conclude that $F\cB(1)_\dagger$ is irreducible. This finishes the proof of the claim that $\cB(t)$ is irreducible for all $t\in [0,1]$.

\vspace{3mm}
$\underline{\B(t) \text{ is Hermitian}}$: thanks to condition (i) and Proposition \ref{prop:HC_induction}, there is an isomorphism in $\HC^G(U(\fg))$
$$\cB(t) = \Phi^*\Ind^G_M[\cB' \otimes \CC(\lambda_{\ell}(t),\lambda_r(t))] \simeq \Ind^G_M[\Phi^*\cB' \otimes \CC(\lambda_{\ell}(t),\lambda_r(t))]$$
Let $(\lambda_{\ell}',\lambda_r') \in \fh^* \times \fh^*$ denote the Langlands parameters for $\Phi^*\B' \in \HC^M(U(\fm))$ and define
$$\\underline{\lambda}_{\ell}(t) = \lambda'_{\ell}+\lambda_{\ell}(t) \in \fh^*, \qquad \underline{\lambda}_r(t) =  \lambda_r' + \lambda_r(t) \in \fh^*, \qquad \forall t \in [0,1].$$
Then $\Phi^*\cB' \simeq \overline{I}^M_H(\lambda'_{\ell},\lambda_r')$ and
$$\Phi^*\cB' \otimes \CC(\lambda_{\ell}(t),\lambda_r(t)) \simeq \overline{I}^M_H(\underline{\lambda}_{\ell}(t),\underline{\lambda}_r(t)), \qquad \forall t \in [0,1].$$
By the previous step, $\cB(t) \simeq \Ind^G_M \overline{I}^M_H(\underline{\lambda}_{\ell}(t),\underline{\lambda}_r(t))$ is irreducible for all $t \in [0,1]$. Thus by Lemma \ref{lem:Langlands}
$$\cB(t) \simeq \overline{I}^G_H(\underline{\lambda}_{\ell}(t),\underline{\lambda}_r(t)), \qquad \forall t \in [0,1].$$
To show that $\B(t)$ is Hermitian, we use Proposition \ref{prop:Hermitian}. Since $\Phi^*\B' \simeq \overline{I}^M_H(\lambda'_{\ell},\lambda_r')$ is Hermitian, there is an element $w' \in W_1$ such that 
$$w'(\lambda'_{\ell}-\lambda_r') = \lambda'_{\ell} - \lambda_r', \qquad w'(\lambda'_{\ell}+\lambda_r') = -\overline{(\lambda_{\ell}'+\lambda_r')}.$$
Take $w \in W_2$ as described in the statement of the proposition. Then
$$ww'(\underline{\lambda}_{\ell}(t)-\underline{\lambda}_r(t)) = ww'(\lambda_{\ell}'-\lambda_r') + ww'(\chi) = \lambda_{\ell}'-\lambda_r' + \chi = \underline{\lambda}_{\ell}(t)-\underline{\lambda}_r(t),$$
and
$$ww'(\underline{\lambda}_{\ell}(t)+\underline{\lambda}_r(t)) = ww'(\lambda_{\ell}'+\lambda_r') + tww'(\chi+2\nu)  = -\overline{(\lambda_{\ell}'+\lambda_r')} - t\overline{(\chi+2\nu)} = - \overline{(\underline{\lambda}_{\ell}(t) + \underline{\lambda}_r(t))}.$$
Hence, $\B(t)$ is Hermitian by Proposition \ref{prop:Hermitian}.
\end{proof}

\section{Construction of unipotent bimodules: some general results}\label{subsec:bimodinductionstep}

In Section \ref{subsec:constructionclassical}, we will construct all unipotent bimodules for linear classical groups via parabolic induction. More precisely, for each nilpotent cover $\widetilde{\OO}$, we will construct a Levi subgroup $L \subset G$ and a rigid orbit $\OO_L$ such that all bimodules in $\unip(\widetilde{\OO})$ are constructed from bimodules in $\unip_{\OO_L}(L)$ through several unitarity-preserving operations: unitary induction, complementary series, and extraction of direct summands. The proof of this result is essentially an induction on the number of codimension 2 leaves in $X=\Spec(\CC[\OO])$. In this section, we will work out the induction step. 

Suppose for convenience that $G$ is semisimple, and let $\mathbb{O}$ be a nilpotent $G$-orbit. Choose a Levi subgroup $L \subset G$ and a birationally rigid $L$-orbit $\mathbb{O}_L$ such that $\mathbb{O} = \mathrm{Bind}^G_L \mathbb{O}_L$. For the remainder of this section, we will assume
\begin{itemize}
    \item[(a1)] $\OO$ admits a birationally rigid cover.
\end{itemize}
Choose a codimension 2 leaf $\fL_k \subset X$ and let $\Sigma_k = \CC^2/\Gamma_k$ be the corresponding singularity. Assumption (a1) guarantees that $H^2(\OO,\CC)=0$, see Corollary \ref{cor:criterionbirigid}. Let $M_k \subset G$ be the (unique) Levi subgroup adapted to $\fL_k$, see (\ref{eq:defofmk}). Then $L \subset M_k$ and $\eta$ restricts to an isomorphism $\eta_k: \fX(\fm_k) \xrightarrow{\sim} \fP_k^X$. As usual, let $\OO_{M_k} = \mathrm{Bind}^{M_k}_L \OO_L$. {As in Section \ref{subsec:identification},} we will impose {the following} additional condition on $\OO$ and $\fL_k$:
\begin{itemize}
    \item[(a2)] $\pi_1(\fL_k)$ acts trivially on $H^2(\mathfrak{S}_k,\CC)$.
\end{itemize}
Here $\mathfrak{S}_k$ denotes, as usual, the minimal resolution of $\Sigma_k$. 

Consider the barycenter parameters $\epsilon \in \fP^X$ and $\epsilon' \in \fP^{X_{M_k}}$, cf. Example \ref{ex:barycentersymplectic}. Define the elements
$$\delta := \eta^{-1}(\epsilon) \in \fX(\fl), \qquad \delta' := \eta_{M_k}^{-1}(\epsilon') \in \fX(\fl\cap [\fm_k,\fm_k]), \qquad 
\delta_k: = \eta_k^{-1}(\epsilon_k) \in \fX(\fm_k).$$
\begin{lemma}\label{lem:delta_induction} 
$\delta = \delta' + \delta_k$.
\end{lemma}
\begin{proof}
The space $\fX(\fl)$ inherits a nondegenerate symmetric bilinear form $B$ from the dual Cartan subalgebra $\fh^*$. Note that $B$ induces a form $B^{\eta}$ on $\fP^X$ via the identification $\eta: \fX(\fl) \xrightarrow{\sim} \fP^X$. We claim that the decomposition $\fP^X = \bigoplus_k \fP_k^X$ is orthogonal with respect to $B^{\eta}$. Note that $B$ is invariant under the natural action of $N_G(L)$. By (ii) of Proposition \ref{prop:namikawacovers} $W^X=N_G(L,\OO_L)/L$. Hence, $B^{\eta}$ is invariant under $W^X$. The decomposition $\fP^X = \bigoplus_k \fP_k^X$ is orthogonal with respect to any nondegenerate $W^X$-invariant symmetric form, and, in particular, with respect to $B^\eta$, as asserted. Under $\eta$, $\fP^X_k$ corresponds to $\fX(\fm_k)$, see
(\ref{eq:etakdef}). Since 
$\fX(\fl\cap [\fm_k,\fm_k]) \subset \fX(\fl)$ is the orthogonal complement of $\fX(\fm_k) \subset \fX(\fl)$, it follows that $\eta \left(\fX(\fl\cap [\fm_k,\fm_k])\right) \subset \fP^X$ is the orthogonal complement of $\fP_k^X$. By the discussion above, $\bigoplus_{i \neq k} \fP_i^X$ is the orthogonal complement of $\fP_k^X$. So $\eta( \fX(\fl\cap [\fm_k,\fm_k])) = \bigoplus_{i \neq k} \fP_i^X$.

Now it suffices to show that the orthogonal projections of $\delta$ to $\fX(\fm_k)$ and $\fX(\fl\cap [\fm_k,\fm_k])$ coincide with $\delta_k$ and $\delta'$, respectively. Note that the projection of $\epsilon$
to $\fP^X_k$ is determined by the isomorphism type of $\Sigma_k$. By Step 4 of the proof of Proposition \ref{prop:part_resol_slice}, the codimension $2$ leaves in $Z_k:=G\times^{Q_k}(X_{M_k}\times \fq_k^\perp)$ are in bijection with the leaves $\mathfrak{L}_i, i\neq k$, via the partial resolution $Z_k\twoheadrightarrow X$. 

We can describe these symplectic leaves explicitly using the geometry of $X_{M_k}$, Namely, let $\mathfrak{L}_i'\subset X_{M_k}$ be a codimension $2$ leaf. Then $\fL_i=G\times^{Q_k}(\mathfrak{L}'_i\times \fq_k^\perp)$ is a symplectic leaf of $Z_k$. Pick a point $x\in \fL_i'\subset \fL_i$, and let $\Sigma_i\subset X_{M_k}$ be the corresponding Kleinian singularity. We have the following commutative diagram:
\begin{center}
        \begin{tikzcd}
          \mathfrak{S}_i\ar[d] \ar[r] & {Y}_{M_k}^{\mathrm{reg}} \ar[d] \ar[r] & {Y}^{\mathrm{reg}} \ar[d] \\
          {\Sigma}_i \ar[r] \ar[rrd]& X_{M_k} \ar[r] & Z_k \ar[d]\\
          && X
        \end{tikzcd}
    \end{center}
    It implies that for every $i\neq k$ the pullback map $H^2({Y}^{\mathrm{reg}}, \CC)\to H^2(\fS_i, \CC)$ factors through the pullback map $H^2(Y^{\mathrm{reg}},\CC)\twoheadrightarrow H^2(Y_{M_k}^{\mathrm{reg}},\CC)$. Therefore, the projection $\fP^X\twoheadrightarrow \bigoplus_{i\neq k}\fP^X_i$ coincides with the pullback map $H^2(Y^{\mathrm{reg}},\CC)\twoheadrightarrow H^2(Y_{M_k}^{\mathrm{reg}},\CC)$. By Proposition \ref{prop:compatibility2}, the latter coincides with the projection $\fX(\fl)\twoheadrightarrow \fX(\fl\cap [\fm_k,\fm_k])$. We conclude that the projection of $\delta$ to $\fX(\fl\cap [\fm_k,\fm_k])$ coincides with $\delta'$. 

It remains to show that the orthogonal projection of $\delta$ to $\fX(\fm_k)$ coincides with $\delta_k$.  This is equivalent to the claim that the unique $W^X$-equivariant projection $\fP^X\twoheadrightarrow \fP^X_{k}$ sends $\epsilon$ to $\epsilon_k$. This follows directly from the definition of $\epsilon$.
\end{proof}

Below, we will construct the irreducible objects in  $\overline{\HC}^G(\cA_{\delta}^X)$ from the irreducible objects in $\overline{\HC}^M(\cA_{\delta'}^{X_M})$ using parabolic induction.

We write $M$ for $M_1$ and $\underline{\delta}$ for $\delta_1$. We still use the notations $\fh_1,W_1$ to avoid a confusion with the objects associated to the Lie algebra $\fg$.  

Let $\{\omega_i\}$ denote the fundamental weights in $\fP_1^X \simeq \fh_1^*$.

Our main result is the following:

\begin{prop}\label{prop:bimods_construction}
In addition to conditions (a1) and (a2) recalled earlier in Section \ref{subsec:bimodinductionstep}, assume that the following conditions hold:
\begin{itemize}
\item[(d1)] $\OO_{M}$ admits a 2-leafless $M$-equivariant cover.
\item[(d2)] $\pi_1^G(\OO)$ and $\pi_1^{M}(\OO_{M})$ are abelian.
\item[(d3)] $\Sigma(=\Sigma_1)$ is of type $A_{d-1}$.
\item[(d4)] $|\pi_1^G(\OO)| = d |\pi_1^M(\OO_M)|$.
\item[(d5)] There are elements $\{\tau_i\}$ of $\mathfrak{X}(M)$ forming a basis in $\mathfrak{X}(\fm)$ such that $\eta_1\{\tau_i\}=\{\omega_i\}$.
\end{itemize}
Then the following claims are true:
\begin{enumerate}
    \item If $\cB \in \overline{\HC}^{M}(\cA_{\delta'}^{X_{M}})$, then the induced bimodules
    $$\Ind^G_{M} [\cB \otimes \CC(\underline{\delta}, \underline{\delta})] \quad \text{and} \quad \Ind^G_{M} [\cB \otimes \CC(\underline{\delta} - \tau_i, \underline{\delta})], \quad 1 \leq i \leq d-1,$$ are  in $\overline{\HC}^G(\cA_{\delta}^X)$.
    \item If $\cB$ is irreducible, then so are the objects in (1).
    \item Every irreducible object in $\overline{\HC}^G(\cA_{\delta}^X)$ is isomorphic to exactly one bimodule of the form described in (1) (for a unique irreducible $\cB$).
\end{enumerate}
\end{prop}
\begin{proof}
For (1), we must show that $\delta$ and $\delta-\tau_i$ are conjugate under $W^X$, for $1 \leq i \leq d-1$. This will, of course, follow if we show that these elements are conjugate under the subgroup $W_1\subset W^X$. Note that $\eta(\delta)=\epsilon$ and $\eta(\delta-\tau_i) = \epsilon-\omega_i$.  By conditions (d3) and (a2), $W_1 \simeq S_d$ with its standard action on $\fh_1^*$. Thus, it suffices to show that $\epsilon$ and $\epsilon-\omega_i$, when written in standard coordinates on $\fh_1^*$, differ by permutation of entries. In standard coordinates
$$\epsilon = \frac{1}{d}(\omega_1+...+\omega_{d-1}) = \frac{1}{d}\sum_{i=1}^{d-1}(\underbrace{d-i,...,d-i}_{i},\underbrace{-i,...,-i}_{d-i}) = \frac{1}{2d}(d-1,d-3,...,3-d,1-d)$$
and
\begin{align*}
\epsilon - \omega_i   &= \frac{1}{2d}(d-1,d-3,...,3-d,1-d) - \frac{1}{d}(\underbrace{d-i,...,d-i}_{i},\underbrace{-i,...,-i}_{d-i})\\
                      &=\frac{1}{2d}(\underbrace{2i-1-d,2i-3-d,...,1-d}_{i},\underbrace{d-1,d-3,...,2i+1-d}_{d-i})
\end{align*}
Let $\sigma_i \in S_d$ denote the permutation given by $\sigma_i(j)=d -i+j$ for $j\leqslant i$ and $\sigma_i(j)=j-i$ for $j>i$. Then clearly $\epsilon-\omega_i = \sigma_i(\epsilon)$. This completes the proof of (1). 

Since $\cB \in \overline{\HC}^{M}(\cA_{\delta'}^{X_{M}})$ is irreducible, the induced bimodules $\Ind^G_{M} [\cB \otimes \CC(\underline{\delta}, \underline{\delta})],\Ind^G_M [\cB \otimes \CC(\underline{\delta}-\tau_i,\underline{\delta})] \in \overline{\HC}^G(\cA_{\delta}^X)$ are irreducible by Corollary \ref{cor:propsofbimodind}(iii). 
This completes the proof of (2). 

Now we proceed to (3).
By condition (a1), the universal $G$-equivariant cover $\widehat{\OO}$ of $\OO$ is 2-leafless. By Proposition \ref{prop:parameterofinvariantssymplectic}, $\cA_{\delta}^X$ is identified with the $\pi_1^G(\OO)$-invariants  in $\cA^{\widehat{X}}_0$, where, as usual $\widehat{X}:=\Spec(\CC[\widehat{\OO}])$. 
Similarly, by Proposition \ref{prop:parameterofinvariantssymplectic} and (d1),  $\cA_{\delta'}^{X_M}$ is  the  $\pi_1^M(\OO_M)$-invariants in $\cA_0^{\widehat{X}_M}$. Hence by Corollary \ref{cor:isotypiceqvt}, there are equivalences
$$\overline{\HC}^G(\cA_{\delta}^X) \simeq \pi_1^G(\OO)\modd, \qquad \overline{\HC}^{M}(\cA_0^{X_M}) \simeq \pi_1^M(\OO_M)\modd.$$
Thus, by condition (d2),  
$$|\mathrm{Irr}[\overline{\HC}^G(\cA_{\delta}^X)]| = |\pi_1^G(\OO)|,   \ |\mathrm{Irr}[\overline{\HC}^{M}(\cA_{\delta'}^{X_M})]|=|\pi_1^M(\OO_M)|.$$
And, by (d4),
$$|\mathrm{Irr}[\overline{\HC}^G(\cA_{\delta}^X)]| = d \ |\mathrm{Irr}[\overline{\HC}^{M}(\cA_{\delta'}^{X_M})]|.$$
To prove (3) it suffices to show that the induced bimodules
$$\{\Ind^G_M[\cB \otimes \CC(\underline{\delta},\underline{\delta})] \mid \cB \text{ irreducible}\} \cup \{\Ind^G_M[\cB \otimes \CC(\underline{\delta},\underline{\delta}-\tau_i)] \mid \cB \text{ irreducible}, \ 1 \leq i \leq d-1\}$$
are pairwise non-isomorphic. Recall that $\tau_i$ determines a character $\varphi^G_M(\tau_i)$ of $\pi_1^G(\OO)$, see Lemma \ref{lem:mutilde}. By \cref{prop:twinductiondagger} we have
$$\left(\Ind^G_M[\cB \otimes \CC(\underline{\delta},\underline{\delta})]\right)_{\dagger} \simeq f^*(\cB_{\dagger}), \quad \left(\Ind^G_M[\cB \otimes \CC(\underline{\delta},\underline{\delta}-\tau_i)]\right)_{\dagger} \simeq  (\varphi^G_M\tau_i)^{-1} \otimes f^*(\cB_{\dagger})$$
for $1 \leq i \leq d-1$. Thus, it is enough to show that 
\begin{itemize}
\item[(*)]
the restrictions of the characters $\varphi^G_M(\tau_i)$ to $\ker[\pi^G(\OO)\twoheadrightarrow \pi^M(\OO_M)]$ are nontrivial and pairwise distinct. 
\end{itemize}

Let $\Sigma^{\times} = \Sigma - \{0\}$.
By Lemma \ref{lem:mutilde}, there is an isomorphism $\mathcal{L}_{\mathbb{O}}(\varphi^G_M(\tau_i))|_{\Sigma^{\times}} \simeq \bar{\pi}^*\mathcal{L}_{G/Q}(\tau_i)|_{\Sigma^{\times}}$ of line bundles on $\Sigma^\times$. By \cref{prop:identification2}, $\bar{\pi}^*\mathcal{L}_{G/Q}(\tau_i)|_{\mathfrak{S}}\simeq \sigma(\omega_i)$, where $\sigma(\omega_i)$ is the line bundle on $\mathfrak{S}$ corresponding to $\omega_i$ via the isomorphism $\sigma: \Lambda \simeq \Pic(\mathfrak{S})$, where we write $\Lambda$ for the weight lattice of $\operatorname{SL}(d)$, see Proposition 
\ref{prop:PicSigma}.

We claim that the restrictions $\sigma(\omega_i)|_{\Sigma^\times}$ are distinct (and nontrivial) for $1\leq i \leq d-1$. 
The restriction map $\Pic(\mathfrak{S})\to \Pic(\Sigma^\times)$ corresponds to the restriction map $\Lambda \twoheadrightarrow \fX(Z(\mathrm{SL}(d)))$, under the identifications $\Lambda \simeq \Pic(\mathfrak{S})$ and $\fX(Z(\mathrm{SL}(d))) \simeq \Pic(\Sigma^{\times})$, see Section \ref{subsec:Mckay}. It is trivial to check that the fundamental weights $\{\omega_i\}$ restrict to distinct nontrivial characters of $Z(\mathrm{SL}(d))$.

Now we identify $\ker[\pi_1^G(\OO)\twoheadrightarrow \pi_1^M(\OO_M)]$ with $\pi_1(\Sigma^\times)$. First, we observe that $\pi_1(\Sigma^\times)$ maps naturally to the kernel. Indeed,  $\pi_1^M(\OO_M)\cong \pi_1^G(G\times^{Q}(\OO_M\times \fq^\perp))$ and under this identification, the homomorphism $\pi_1^G(\OO)\twoheadrightarrow \pi_1^M(\OO_M)$ is induced by the inclusion $\OO\hookrightarrow 
G\times^Q(\OO_M\times \fq^\perp)$, see Section \ref{subsec:birationalinduction}. The inclusion $\Sigma^\times\hookrightarrow \OO$ gives rise to a homomorphism $\pi_1(\Sigma^\times)\rightarrow \pi_1(\OO)$. On the other hand, the composition $\Sigma^\times\hookrightarrow \OO\hookrightarrow G\times^Q(\OO_M\times \fq^\perp)$ factors through $\mathfrak{S}$, which is simply connected. To see the latter, notice that $\mathfrak{S}$ comes with an action of $\mathbb{G}_m^2$ with finitely many fixed points. From here one sees that there is a Zariski open subset of $\mathfrak{S}$ isomorphic to $\mathbb{A}^2$. 

So the image of $\pi_1(\Sigma^\times)\rightarrow \pi_1(\OO)$ lies in the kernel of $\pi_1(\OO)\twoheadrightarrow \pi_1(G\times^Q(\OO_M\times \fq^\perp))$. The $G$-equivariant fundamental groups are quotients of the usual ones, so we get a homomorphism 
$\pi_1(\Sigma^\times)\rightarrow \ker[\pi_1^G(\OO)\rightarrow \pi_1^M(\OO_M)]$. 
The cardinalities of these groups are the same (by (d3) and (d4)), so to prove their isomorphism, it suffices to show that the homomorphism is injective. Since the groups in question are abelian, this is equivalent to saying that the homomorphism $\mathfrak{X}(\pi_1^G(\OO))\rightarrow \mathfrak{X}(\pi_1(\Sigma^\times))$ is surjective. But these character groups are nothing else but $\operatorname{Pic}^G(\OO)$ and $\operatorname{Pic}(\Sigma^\times)$, and the homomorphism is the pullback under $\Sigma^\times\hookrightarrow \OO$. We have seen that it is surjective in the previous paragraph.
\end{proof}

\section{Unipotent bimodules attached to birationally rigid covers: classical case}\label{subsec:constructionclassical}

{In Sections \ref{subsec:bimodsA} and \ref{subsec:bimodsBCD} below, we will give a construction of the unipotent bimodules attached to birationally rigid nilpotent covers for classical groups. The proofs in all cases follow a (more or less) uniform pattern, which we will now briefly outline for the reader's convenience. Let $\widetilde{\OO}$ be a birationally rigid $G$-equivariant nilpotent cover and suppose that $\OO = \mathrm{Bind}^G_L \OO_L$ for a Levi subgroup $L$ and a rigid nilpotent $L$-orbit $\OO_L$. Then $\pi_1(\OO)$ acts on the canonical quantization of $\CC[\widetilde{\OO}]$, and the $\pi_1(\OO)$-invariants are a (non-canonical) quantization $\cA_{\delta}^X$ of $\CC[\OO]$. First, we check that the ideal $I_0(\widetilde{\OO}) = I_{\delta}(\OO)$ is maximal (using the results of Appendix \ref{sec:maximality}). Then we deduce (using Lemma \ref{lem:simplenoboundary} and Corollary \ref{Cor:check_A_description}) that $\overline{\HC}^G(\cA_{\delta}^X) \xrightarrow{\sim} \HC^G(U(\fg)/I_{\delta}(\OO))$. After checking the relevant hypotheses, we apply Proposition \ref{prop:bimods_construction} to construct the irreducibles in $\overline{\HC}^G(\cA_{\delta}^X)$ via parabolic induction from irreducibles in $\overline{\HC}^L(\cA_0^{X_L})$. Finally, to verify that these inductions preserve unitarity, we appeal to Propositions \ref{prop:inductionunitarygeom} and \ref{prop:deformationgeom}.}

\subsection{Type A}\label{subsec:bimodsA}

Let $G=\mathrm{SL}(n)$ and let $\widetilde{\OO}$ be a birationally rigid $G$-equivariant nilpotent cover. By Proposition \ref{prop:birigidcoverA}, $\widetilde{\OO}$ is the universal cover of an orbit $\OO$ corresponding to a partition of the form $(d^m)$ of $n$. Let 
$$L=\mathrm{S}(\mathrm{GL}(m)^d), \qquad \mathbb{O}_L = \{0\},$$
so that $\mathbb{O} = \mathrm{Bind}^G_L \mathbb{O}_L$, and fix $\tau_1(1),\tau_2(1),...,\tau_{d-1}(1), \delta \in \fh^*$ as in Corollary \ref{cor:lambdaA} (and its proof). Thus, $\cA_{\delta}^X$ is isomorphic to the $\pi_1(\OO)$-invariants in $\cA_0^{\widetilde{X}}$ and $\pi_1(\OO) \simeq \ZZ_d$.

\begin{prop}\label{prop:bimodsA}
The set $\unip_{\widetilde{\OO}}(G)$ consists of the following $d$ bimodules
$$\{\Ind^G_L \CC(\delta,\delta)\} \cup \{\Ind^G_L\CC(\delta-\tau_i(1),\delta) \mid 1\leq i \leq d-1\}.$$
All are unitary.
\end{prop}

\begin{proof}
Let $I =I_0(\widetilde{\OO})= I_{\delta}(\OO)$. By Proposition \ref{prop:maximalityA}, $I$ is maximal, so the algebra $\cA_{\delta}^X$ is simple, see Proposition \ref{prop:simplemaximal}. Thus by Lemma \ref{lem:simplenoboundary}, 
$$\HC^G(\cA_{\delta}^X) \xrightarrow{\sim} \overline{\HC}^G(\cA_{\delta}^{X}),$$
By Corollary \ref{Cor:check_A_description}, $\Phi: U(\fg)/I \xrightarrow{\sim} (\cA_{\delta}^X)^{\Pi}$. {Combining this with Lemma \ref{Lem:quant_bimodule_easy}, we see that}
$$\Phi^*: \HC^G(\cA_{\delta}^{X}) \xrightarrow{\sim} \HC^G(U(\fg)/I).$$
By Proposition \ref{prop:Lforbirigidcover}, there is a unique codimension 2 leaf $\fL_1 \subset X$ and, by Proposition \ref{prop:nocodim2leavesLM},
$$M = L, \qquad \OO_{M} = \{0\}.$$
Conditions (a1)-(a2) as well as (d5) were verified in the proof of Corollary \ref{cor:lambdaA}. To apply Proposition \ref{prop:bimods_construction} it remains to check (d1)-(d4). Since $\OO_{M}=\{0\}$, condition (d1) is trivial. Condition (d2) and (d4) follow from $\pi_1^G(\OO) \simeq \ZZ_d$. (d3) was established in Proposition \ref{prop:Lforbirigidcover}(i). So, we can apply Proposition \ref{prop:bimods_construction}.

If $\Phi^*\cB \in \unip_{\widetilde{\OO}}(G)$, then, by Proposition \ref{prop:bimods_construction},
$$\cB \simeq \Ind^G_L \CC(\delta,\delta) \quad \text{or} \quad \cB \simeq \Ind^G_L \CC(\delta-\tau_i(1),\delta), \quad 1 \leq i \leq d-1.$$
First, assume $\B \simeq \Ind^G_L \CC(\delta,\delta)$. We will show that $\B$ is unitary using Proposition \ref{prop:deformationgeom}. Set
$$M_1 = 1, \qquad M_2 = L, \qquad \mathbb{O}_M = \{0\}, \qquad \beta =0, \qquad \nu = \delta, \qquad \chi=0.$$
Condition (i) of Proposition \ref{prop:deformationgeom} follows from Lemma \ref{lem:auxidealsA1}. It remains to exhibit an element $w \in S_n$ such that $w(\delta) = - \delta$. For this, we choose $w=w_0$ (the longest element of $S_n$).

Next, assume $\B \simeq \Ind^G_L\CC(\delta -\tau_i(1),\delta)$. We will show that $\B$ is unitary in the same way. Set
$$M_1 = L, \qquad M_2 = \{1\}, \qquad \widetilde{\mathbb{O}}_M = \{0\}, \qquad \beta = 0, \qquad \nu = \delta, \qquad \chi=-\tau_i(1).$$
Condition (i) of Proposition \ref{prop:deformationgeom} follows from Lemma \ref{lem:auxidealsA2}. It remains to exhibit an element $w \in S_n$ such that
\begin{equation}\label{eq:wcond}
w\tau_i(1) = \tau_i(1), \qquad w(2\delta -\tau_i(1)) = \tau_i(1)- 2\delta.\end{equation}
Note that
\begin{align*}
\tau_i(1) &= \frac{1}{d}(\underbrace{d-i,d-i,...,d-i}_{mi}, \underbrace{-i,-i,...,-i}_{m(d-i)}),\\
2\delta-\tau_i(1) &= \frac{1}{d}(\underbrace{i-1,...}_m, \underbrace{i-3,...}_m,...,\underbrace{1-i,...}_m, \underbrace{d-i-1}_m,...,\underbrace{1+i-d}_m).\end{align*}
Write $w'_0 \in S_{mi}$ and $w_0'' \in S_{m(d-i)}$ for the longest elements, and take $w$ to be the image of $(w'_0,w''_0) \in S_{mi} \times S_{m(d-i)}$ under the natural embedding $S_{mi} \times S_{m(d-i)} \subset S_n$. Then an easy calculation shows that $w$ satisfies (\ref{eq:wcond}). 
\end{proof}

\subsection{Types B,C, and D}\label{subsec:bimodsBCD}

Let $G=\mathrm{Sp}(2n)$, $\mathrm{SO}(2n)$, or $\mathrm{SO}(2n+1)$, and let 
$\widetilde{\OO}$ be a birationally rigid $G$-equivariant nilpotent cover. We will first assume that $\widetilde{\OO}$ covers an orbit $\OO$ which is \emph{not} of the form $\mathbb{O}_{(4^{2m},3,1)}$ for $G=\mathrm{SO}(2n)$. The orbits $\OO$ which can appear are described in Proposition \ref{prop:nocodim2leaves}.

Let
$$L= \prod_{k \in S_2(p)} \mathrm{GL}(k) \times G(n-|S_2(p)|), \qquad \mathbb{O}_L = \prod_{k \in S_2(p)} \{0\} \times \mathbb{O}_{p\#S_2(p)} $$
so that $\mathbb{O} = \mathrm{Bind}^G_L \mathbb{O}_L$, and fix $\{\tau_1(k) \mid k \in S_2(p)\} \subset \fh^*$ and $\delta \in \fh^*$ as in Corollary \ref{cor:lambdaBCD}(i) (and its proof). Thus, $\cA_{\delta}^X$ is identified with the $\pi_1(\OO)$-invariants in $\cA_0^{\widetilde{X}}$.

\begin{prop}\label{prop:bimodsBCD}
For every bimodule $\cB \in \unip_{\widetilde{\OO}}(G)$ there is a unique subset $S \subseteq S_2(p)$ and a unique bimodule $\cB' \in \unip_{\OO_L}(L)$ such that
$$\cB\cong\Ind^G_L [\cB' \otimes \CC(\delta-\sum_{k \in S}\tau_1(k),\delta)].$$
If $\cB'$ is unitary, then $\cB$ is unitary.
\end{prop}

\begin{proof}
Let $I =I_0(\widetilde{\OO})= I_{\delta}(\OO)$. By Proposition \ref{prop:maximalitytypeC}, $I$ is maximal, so the algebra $\cA_{\delta}^{X}$ is simple, see Proposition \ref{prop:simplemaximal}. Thus by Lemma \ref{lem:simplenoboundary}, 
$$\HC^G(\cA_{\delta}^{X}) \xrightarrow{\sim} \overline{\HC}^G(\cA_{\delta}^{X}),$$
By Corollary \ref{Cor:check_A_description}, $\Phi: U(\fg)/I \xrightarrow{\sim} \cA_{\delta}^X$. {Combining this with Lemma \ref{Lem:quant_bimodule_easy}, we see that}
$$\Phi^*: \HC^G(\cA_{\delta}^{X}) \xrightarrow{\sim} \HC^G(U(\fg)/I).$$
Hence, $\cB$ is the image under $\Phi^*$ of a (unique) irreducible object in $\HC^G(\cA_{\delta}^X)$ (which we will continue to denote by $\cB$).

We will prove the proposition by induction on $|S_2(p)|$. If $|S_2(p)|=0$, then $\mathbb{O}$ is birationally rigid by Proposition \ref{prop:birigidorbitclassical} and the assertion is vacuous.

Now suppose $|S_2(p)| \geq 1$, and choose $k \in S_2(p)$. Recall, see Proposition \ref{prop:Lforbirigidcover}, that the set $S_2(p)$ parameterizes codimension 2 leaves $\fL \subset X$. 
Let $M_k\subset G$ be the Levi subgroup adapted to the leaf $\fL_k$, see (\ref{eq:defofmk}), and let $\OO_{M_k} = \mathrm{Bind}^{M_k}_L \OO_L$. We claim that conditions (a1),(a2) and (d1)-(d5) hold for $\OO$ and $\fL_k$. Of course, (a1) holds by assumption. By Propositions \ref{prop:Lforbirigidcover} and \ref{prop:nocodim2leavesLM}, $\Sigma_k \simeq \CC^2/\ZZ_2$ and
$$M_k \simeq GL(k) \times G(n-k), \qquad \OO_{M_k} = \{0\} \times \OO_{p \# (k)}.$$
Conditions (a2) and (d3) follow immediately. Condition (d1) follows from Proposition \ref{prop:nocodim2leaves}. Condition (d2) is automatic since $G$ is linear classical. Conditions (d4) was verified in the proof of Corollary \ref{cor:lambdaBCD}. In the proof of Corollary \ref{cor:lambdaBCD} we have also seen that $c_k$ defined by (\ref{eq:defofck}) is equal to $1$. Then (d5) follows from Proposition \ref{prop:identification2}. Now Proposition \ref{prop:bimods_construction} implies
$$\text{either }\cB \simeq \Ind^G_{M_k}[\cB^1 \otimes \CC(\frac{1}{2}\tau_1(k),\frac{1}{2}\tau_1(k))] \quad \text{or} \quad \cB \simeq \Ind^G_{M_k}[\cB^1 \otimes \CC(-\frac{1}{2}\tau_1(k),\frac{1}{2}\tau_1(k))]$$
for a uniquely determined irreducible object $\cB^1 \in \overline{\HC}^{M_k}(\cA_{\delta-\frac{1}{2}\tau_1(k)}^{X_{M_k}})$. By the induction hypothesis,
$$\cB^1 \simeq \Ind^{M_k}_{L}[\cB' \otimes \CC(\delta-\frac{1}{2}\tau_1(k)-\sum_{j \in S'}\tau_1(j),\delta-\frac{1}{2}\tau_1(k))]$$
for a uniquely determined irreducible object $\cB' \in \overline{\HC}^L(\cA_0^{X_L})$ and a uniquely determined subset $S' \subset S_2(p)$ not containing $k$. By Corollary \ref{cor:propsofbimodind}(ii),
$$\cB \simeq \Ind^G_L [\cB' \otimes \CC(\delta - \sum_{k \in S} \tau_1(k),\delta)].$$
It remains to show that $\cB$ is unitary if $\cB'$ is unitary. If $\cB'$ is unitary, then $\cB^1$ is unitary by the induction hypothesis. To show that that $\cB$ is unitary, we apply Propositions \ref{prop:inductionunitarygeom} and \ref{prop:deformationgeom}. 

First, suppose $\B \simeq \Ind^G_{M_k}[\B^1 \otimes \CC(\frac{1}{2}\tau_1(k),\frac{1}{2}\tau_1(k))]$. Condition (i) of Proposition \ref{prop:deformationgeom} follows from Lemma \ref{lem:auxilaryidealstypeC}. It remains to exhibit an element $w \in W$ such that $w\delta = -\delta$. For this, we choose $w = w_0$ (the longest element of $W$).

Next, suppose $\B \simeq \Ind^G_{M_k}[\B^1 \otimes \CC(-\frac{1}{2}\tau_1(k),\frac{1}{2}\tau_1(k))]$. By Lemmas \ref{lem:maximalitybirigidtypeC}, \ref{lem:maximalitybirigidtypeB}, and \ref{lem:maximalitybirigidtypeD}, the ideal $I_{\delta}(\mathbb{O}) \subset U(\fg)$ is maximal. Thus, $\B$ is unitary by Proposition \ref{prop:inductionunitarygeom}.
\end{proof}

Finally, suppose $G=\mathrm{SO}(8m+4)$, $\mathbb{O} = \mathbb{O}_{(4^{2m},3,1)}$, and $\widetilde{\mathbb{O}} \to \mathbb{O}$ is the universal $G$-equivariant cover. Let
$$L=\mathrm{GL}(2m+1) \times \mathrm{GL}(2m+1), \qquad \mathbb{O}_L = \{0\},$$
so that $\mathbb{O} = \mathrm{Bind}^G_L \mathbb{O}_L$. Fix $\tau_1(1), \tau_1(2), \delta \in \fh^*$ as in Corollary \ref{cor:lambdaBCD}(ii). Thus, $\cA_{\delta}^X$ is identified with the $\pi_1(\OO)$-invariants in $\cA_0^{\widetilde{X}}$ and $\pi_1(\OO) \simeq \ZZ_2$.

\begin{prop}\label{prop:bimodsBCD2}
The set $\unip_{\widetilde{\OO}}(G)$ consists of the following 2 bimodules
$$\{\Ind^G_L\CC(\delta,\delta), \  \Ind^G_L\CC(\delta-\tau_1(1),\delta)\}.$$
Both are unitary.
\end{prop}

\begin{proof}
Let $I =I_0(\widetilde{\OO})= I_{\delta}(\OO)$. By Proposition \ref{prop:maximalitytypeC}, $I$ is maximal, so the algebra $\cA_{\delta}^X$ is simple, see Proposition \ref{prop:simplemaximal}. Thus by Lemma \ref{lem:simplenoboundary},
$$\HC^G(\cA_{\delta}^X) \xrightarrow{\sim} \overline{\HC}^G(\cA_{\delta}^X).$$
By Corollary \ref{Cor:check_A_description}, $\Phi: U(\fg)/I \xrightarrow{\sim} \cA_{\delta}^X$. In particular
$$\Phi^*: \HC^G(\cA_{\delta}^X) \xrightarrow{\sim} \HC^G(U(\fg)/I).$$
We note that (d1) of Proposition \ref{prop:bimods_construction} is not satisfied, this follows from Corollary \ref{cor:lambdaBCD}. 
So we cannot apply that proposition directly, however techniques similar to its proof work.

Since $\pi_1^G(\mathbb{O}) \simeq \ZZ_2$, the category $\overline{\HC}^G(\cA_{\delta}^X) \simeq \HC(U(\fg)/I)$ contains two irreducible objects, see Theorem \ref{thm:classificationbimods}. So for the first claim it suffices to show that the induced bimodules $\cB:=\Ind^G_L\CC(\delta,\delta)$, $\cB':=\Ind^G_L \CC(\delta-\tau_1(1),\delta)$ are irreducible and distinct. Irreducibility is immediate from Corollary \ref{cor:propsofbimodind}(iii). By Proposition \ref{prop:twinductiondagger},
$$\cB_{\dagger} \simeq \operatorname{triv}, \qquad \cB'_{\dagger} \simeq \varphi^G_L( \tau_1(1)).$$
Thus in order to prove that $\cB \not\simeq \cB'$, it is enough to show that $\varphi^G_L (\tau_1(1))$ is a nontrivial character of $\pi_1^G(\OO)$. Let $\Sigma_1^{\times} = \Sigma_1 - \{0\}$, and let $\mathfrak{S}_1 \to \Sigma_1$ be the minimal resolution. By Lemma \ref{lem:mutilde}, there is an isomorphism $\mathcal{L}_{\mathbb{O}}\left(\varphi^G_L(\tau_1(1))\right)|_{\Sigma_1^{\times}} \simeq \pi^*\mathcal{L}_{G/P}(\tau_1(1))|_{\Sigma^{\times}_1}$ of line bundles on $\Sigma_1^\times$. And by \cref{prop:identification2}, $\pi^*\mathcal{L}_{G/P}(\tau_1(1))|_{\mathfrak{S}_1}\simeq \sigma_1(\omega_1(1))$, where $\sigma_1(\omega_1(1))$ is the line bundle on $\mathfrak{S}_1$ corresponding to $\omega_1(1)$ via the natural identification $\sigma_1: \Lambda_1 \simeq \Pic(\mathfrak{S}_1)$, see Proposition 
\ref{prop:PicSigma}. Thus, it suffices to show that $\sigma_1(\omega_1(1))$ restricts to a nontrivial line bundle on $\Sigma^{\times}$. This follows as in the penultimate paragraph of the proof of Proposition \cref{prop:bimods_construction} (setting $d=2$). 

 It remains to show that $\cB$ and $\cB'$ are unitary. To prove $\cB$ is unitary, we use Proposition \ref{prop:deformationgeom}. Condition (i) of that proposition follows from Lemma \ref{lem:auxilaryidealstypeC}. It remains to exhibit an element $w \in W$ such that $w\delta = -\delta$. For this, we choose $w = w_0$ (the longest element of $W$). The unitarity of $\cB'$ is proved by a similar argument.
\end{proof}

\section{Unitarity of unipotent bimodules: classical case}\label{subsec:bimodsclassical}

Let $G$ be a linear classical group, and let $\widetilde{\OO}$ be a $G$-equivariant nilpotent cover, maximal in its equivalence class. In this section, we will prove the following result.

\begin{theorem}\label{thm:unipotentunitary}
The set $\unip_{\widetilde{\OO}}(G)$ consists of unitary bimodules.
\end{theorem}

The proof will require several preliminary results. The first is an easy consequence of Barbasch's classification.

\begin{lemma}\label{lem:Barbasch}
Suppose $\OO$ is rigid. Then the set $\unip_{\OO}(G)$ consists of unitary bimodules.
\end{lemma}

\begin{proof}
Let $\cB \in \unip_{\OO}(G)$. By Proposition \ref{prop:propsofIbeta}, $V(\cB) = \overline{\mathbb{O}}$, and the ideal $\mathrm{LAnn}(\B) = \mathrm{RAnn}(\B) = I_0(\mathbb{O}) \subset U(\mathfrak{g})$ is maximal by Lemmas \ref{lem:maximalitybirigidtypeC}, \ref{lem:maximalitybirigidtypeB}, and \ref{lem:maximalitybirigidtypeD}. Thus, $\cB$ is unitary by \cite[Prop 10.6]{Barbasch1989} and Remark \ref{rmk:Barbasch}.
\end{proof}

The next two lemmas hold in complete generality. For the first lemma, choose a Levi subgroup $M \subset G$ and an $M$-equivariant nilpotent cover $\widetilde{\OO}_M$ such that $\widetilde{\OO}=\mathrm{Bind}^G_M \widetilde{\OO}_M$.

\begin{lemma}\label{lem:uniptounip}
Suppose the ideals $I_0(\widetilde{\OO}_M)$, $I_0(\widetilde{\OO})$ are maximal and $\unip_{\widetilde{\OO}_M}(M)$ consists of unitary bimodules.
Then $\unip_{\widetilde{\OO}}(G)$ consists of unitary bimodules.
\end{lemma}

\begin{proof}
Let $\widehat{\OO}_M$ denote the maximal cover in the equivalence class $[\widetilde{\OO}_M]$. By Theorem \ref{thm:classificationbimods}, there are category equivalences
$$\Aut_{\OO}(\widetilde{\OO})\modd \xrightarrow{\sim} \HC^G(U(\fg)/I_0(\widetilde{\OO})), \qquad \Aut_{\OO_M}(\widehat{\OO}_M)\modd \xrightarrow{\sim} \HC^M(U(\fm)/I_0(\widetilde{\OO}_M)).$$
In particular, the categories $\HC(U(\fg)/I_0(\widetilde{\OO}))$ and $\HC(U(\fm)/I_0(\widetilde{\OO}_M))$ are semisimple. This means that $\cA_0^{\widetilde{X}_M}$, regarded as an object in $\HC(U(\fm)/I_0(\widetilde{\OO}_M))$, is a direct sum of bimodules in $\unip_{\widetilde{\OO}_M}(M)$ and hence unitary.  By Proposition \ref{prop:quantizationparaminduction}, $\Ind^G_M \cA_0^{\widetilde{X}_M} \simeq \cA_0^{\widetilde{X}}$. So, $\cA_0^{\widetilde{X}}$ is unitary by Proposition \ref{prop:inductionunitarygeom} (setting $\chi=0$). By Theorem \ref{thm:classificationbimods}, the bimodules in $\unip_{\widetilde{\OO}}(G)$ are isotypic components in $\cA_0^{\widetilde{X}}$, and thus unitary as well.

\end{proof}

We will also need the following general fact about birationally rigid orbits.

\begin{lemma}\label{lem:rigidtobirigid}
Suppose $\OO$ is a birationally rigid orbit. Then there is a Levi subgroup $L \subset G$ and a rigid orbit $\mathbb{O}_L$ such that $\widehat{\mathbb{O}} = \mathrm{Bind}^G_L \widehat{\mathbb{O}}_L$, where $\widehat{\OO} \to \OO$, and $\widehat{\mathbb{O}}_L \to \OO_L$ are the universal equivariant covers.
\end{lemma}

\begin{proof}
{\it Step 1}. Choose a Levi subgroup $L\subset G$ and a birationally rigid $L$-equivariant nilpotent cover $\widetilde{\mathbb{O}}_L$ such that $\widehat{\mathbb{O}} = \mathrm{Bind}^G_L \widetilde{\mathbb{O}}_L$. By Lemma \ref{lem:mappi1}, there is a surjective homomorphism $\{1\}= \pi_1^G(\widehat{\mathbb{O}}) \twoheadrightarrow \pi_1^L(\widetilde{\mathbb{O}}_L)$. Thus, $\widetilde{\mathbb{O}}_L$ is the universal $L$-equivariant cover $\widehat{\OO}_L$ of $\OO_L$.

{\it Step 2}. Next, we show that $\OO_L$ is 2-leafless, cf. Definition \ref{defi:2_leafless}. Let $\widetilde{\OO} = \mathrm{Bind}^G_L \OO_L$.
Set $\widetilde{X}:=\operatorname{Spec}(\CC[\widetilde{\OO}]), 
\widetilde{X}_L:=\operatorname{Spec}(\CC[\widetilde{\OO}_L])$.
By Corollary \ref{cor:criterionbirigid}, $\OO$ is 2-leafless. Hence, $\widetilde{\OO}$ is 2-leafless, see Lemma \ref{lem:cover2leafless}. Suppose $\fL' \subset \Spec(\CC[\OO_L])$ is a codimension 2 leaf, and let $\OO_L' \subset \overline{\OO}_L$ be the corresponding codimension 2 $L$-orbit, see Lemma \ref{lem:surjectionleaves}. Note that $\fL_Y:=G\times ^P(\fL' \times \fp^\perp)\subset 
\widetilde{Y}:=G\times ^P(\widetilde{X}_L \times \fp^\perp)$ is a codimension $2$ leaf. Since $\rho: \widetilde{Y}\to \widetilde{X}$ is a partial Poisson reslution, $\rho(\fL_Y)$ is a proper closed Poisson subvariety in $\widetilde{X}$, hence it has to be contained in $X-\widetilde{X}^{\mathrm{reg}}$. On the other hand,  consider the codimesion 2 $G$-orbit $\OO'=\Ind_L^G(\OO'_L) \subset \overline{\OO}$, and let $\widetilde{\OO}' \subset \widetilde{X}$ be the preimage of $\OO'$ under the cover $\widetilde{X} \to \overline{\OO}$. Since $\widetilde{\OO}$ is $2$-leafless, $\widetilde{\OO}'$ is contained in $\widetilde{X}^{\mathrm{reg}}$. But since $\OO'$ is induced from $\OO'_L$, $\widetilde{\OO}'$ intersects nontrivially with $\rho(\fL_Y)$. This is a contradiction.


{\it Step 3} Suppose $\OO_L = \Ind^L_{L'} \OO_{L'}$ for a proper Levi subgroup $L' \subset L$, and let $\widecheck{\OO}_L = \mathrm{Bind}^L_{L'} \OO_{L'}$. Since $\widecheck{\OO}_L$ covers $\OO_L$ and $\OO_L$ is 2-leafless, $\widecheck{\OO}_L$ is 2-leafless, see Lemma \ref{lem:cover2leafless}. Thus, Corollary \ref{cor:criterionbirigid} implies that $H^2(\widecheck{\OO}_L,\CC)\neq 0$, since $\widecheck{\OO}_L$ is, by definition, birationally induced. Hence, $\fP^{\widehat{X}_L} \neq 0$, since the pullback map $H^2(\widecheck{\OO}_L,\CC) \to H^2(\widehat{\OO}_L,\CC)$ is an embedding. It follows that $\widehat{\OO}_L$ is birationally induced. This contradicts Step 1.
\end{proof}

\begin{proof}[Proof of Theorem \ref{thm:unipotentunitary}]
Choose a Levi subgroup $K \subset G$ and a birationally rigid $K$-equivariant nilpotent cover $\widetilde{\OO}_K$ such that $\widetilde{\OO} = \mathrm{Bind}^G_K \widetilde{\OO}_K$. Choose a Levi subgroup $M \subset K$ and a birationally rigid orbit $\OO_M$ such that $\OO_K = \mathrm{Bind}^K_M \OO_M$. Finally, choose (using Lemma \ref{lem:rigidtobirigid}) a Levi subgroup $L \subset M$ and a rigid orbit $\OO_L$ such that $\widehat{\OO}_M = \mathrm{Bind}^M_L \widehat{\OO}_L$, where $\widehat{\OO}_M$ (resp. $\widehat{\OO}_L$) is the universal equivariant cover of $\OO_M$ (resp. $\OO_L$).

By Lemma \ref{lem:Barbasch}, the set $\unip_{\OO_L}(L)$ consists of unitary bimodules. Hence by Lemma \ref{lem:uniptounip}, applied to $\widehat{\OO}_M = \mathrm{Bind}^M_L \widehat{\OO}_L$, the same is true of $\unip_{\widehat{\OO}_M}(M)$. Note that $\unip_{\widehat{\OO}_M}(M) = \unip_{\OO_M}(M)$, since $[\widehat{\OO}_M]=[\OO_M]$. Hence by Propositions \ref{prop:bimodsA}, \ref{prop:bimodsBCD}, and \ref{prop:bimodsBCD2} the set $\unip_{\widetilde{\OO}_K}(K)$ consists of unitary bimodules. Finally, Lemma \ref{lem:uniptounip}, applied to $\widetilde{\OO}=\mathrm{Bind}^G_K \widetilde{\OO}_K$, shows that $\unip_{\widetilde{\OO}}(G)$ consists of unitary bimodules. This completes the proof.
\end{proof}

\section{Spin and exceptional groups}\label{subsec:spinexceptional}

The strategy outlined in the beginning of Section \ref{sec:unipbimod} for proving unitarity has two basic components:

\begin{enumerate}
    \item Show that every unipotent bimodule is obtained from one attached to a rigid nilpotent orbit via unitarity-preserving operations (i.e. parabolic induction, complementary series, and extraction of direct summands).
    \item Show that every unipotent bimodule attached to a rigid orbit is unitary. 
\end{enumerate}

For linear classical groups, this strategy was successfully implemented in Section \ref{subsec:bimodsclassical}. For (1), we use an exhaustion argument (involving a detailed understanding of fundamental groups and birational induction). For (2), we appeal to a classical result of Barbasch. 

For spin and exceptional groups, several problems arise. For the exhaustion arguments in Section \ref{subsec:bimodinductionstep} we assume that the fundamental group $\pi_1^G(\OO)$ is abelian. For spin and exceptional groups, this is often not the case. We emphasize, however, that this condition on $\pi_1^G(\OO)$ is an artifact of the proof. It should be possible to weaken this condition (or dispense with it altogether). A second, more fundamental issue is the result of Barbasch applies only to classical groups. 

As noted previously, Dougal Davis and the second-named author have given a uniform conceptual proof of the unitarity of \emph{all} unipotent bimodules (\cite[Corollary 5.23]{DavisMasonBrown}). This result supersedes the unitarity results in this chapter. However, we note that the argument in \cite{DavisMasonBrown} does not give a construction of unipotent bimodules via parabolic induction, which may be of independent interest.

\addtocontents{toc}{\string\vspace{-5pt}}

\appendix

\chapter{Coincidence of Inductions for Harish-Chandra Bimodules}\label{sec:coincidence}

In Section \ref{subsec:HCbimodsclassical}, we recalled the classical construction of parabolic induction for Harish-Chandra $U(\fg)$-bimodules. In Section \ref{subsec:bimodinductioncovers}, we defined parabolic induction for Harish-Chandra bimodules for Hamiltonian quantizations of nilpotent covers. In this appendix, we will show that under certain conditions these two constructions coincide (more precisely, are intertwined by forgetful functors). A similar result was obtained by Vogan in \cite{Vogan1990}, see, for example, \cite[Corollary 6.16]{Vogan1990}. Checking the conditions of Vogan's result in our setting is difficult, so we will pursue a different approach. 

We will need some notation:

\begin{itemize}
    \item Let $Q \subset G$ be a parabolic subgroup with Levi decomposition $Q=MU$. Let $Q^-=MU^-$ denote the opposite parabolic.
    \item Let $\widetilde{\OO}_M$ be an $M$-equivariant nilpotent cover and $\widetilde{\OO} = \mathrm{Bind}^G_M \widetilde{\OO}_M$.
    \item Let $\cA_M^{\ell}$ and $\cA_M^r$ be Hamiltonian quantizations of $\widetilde{X}_M=\Spec(\CC[\widetilde{\OO}_M])$.
    \item Let $\cA^{\ell} := \Ind^G_M \cA_M^{\ell}$ and $\cA^r := \Ind^G_M \cA_M^r$, Hamiltonian quantizations of $\widetilde{X} = \Spec(\CC[\widetilde{\OO}])$.

    \item Let
    $$\Phi^{\ell}: U(\fg) \to \cA^{\ell}, \quad \Phi^r: U(\fg) \to \cA^r, \quad \Phi^{\ell}_M: U(\fm) \to \cA_M^{\ell}, \quad \Phi^r_M: U(\fm) \to \cA_M^r$$
    denote the quantum co-moment maps.
    
    \item Set
    $$I^{\ell}:=\ker{\Phi^{\ell}}, \quad I^r:=\ker{\Phi^r}, \quad  I^{\ell}_M:=\ker{\Phi^{\ell}_M}, \quad  I^r_M:=\ker{\Phi^r_M}.$$
    \item Let $\gamma^{\ell},\gamma^r,\gamma^{\ell}_M, \gamma^R_M$ denote the respective infinitesimal characters. {We assume that $\gamma^\ell_M,\gamma^r_M$ differ by a character of $M$.}
    \item Let 
    $$\Phi^*: \HC^G(\cA^{\ell},\cA^r) \to \HC^G(U(\fg)), \qquad \Phi_M^*: \HC^M(\cA_M^{\ell},\cA_M^r)\to \HC^M(U(\fm))$$
    denote the forgetful functors.
    \item Let $\cB_M$ be an irreducible object in $\HC^M(\cA^{\ell}_M,\cA^r_M)$ and $\cB:=\Ind^G_M\cB_M \in \HC^G(\cA^{\ell},\cA^r)$. 
\end{itemize}

\begin{prop}\label{prop:HC_induction}
Suppose that 
$I^\ell, I^r,I^r_M$
are maximal ideals.
Then there is an isomorphism in $\HC^G(U(\fg))$
$$\Phi^*\cB \simeq \Ind^G_M \Phi_M^* \cB_M$$
\end{prop}

The proof scheme is classical, compare to \cite{Vogan1990}. We embed the left hand side into the right hand side and show that the $G$-types coincide. 

The proof requires some preparation. Recall from Section \ref{subsec:bimodinductioncovers} that $\cA^{\ell}$,$\cA^r$, $\cB$ are obtained as the global sections of sheaves on $G/Q$, denoted here by $\mathcal{D}^{\ell}$, $\mathcal{D}^r$, and $\mathcal{B}^{\mathrm{loc}}$, respectively. Write
$\bar{\cA}^\ell$, $\bar{\cA}^r$, and $\bar{\mathcal{B}}$ for the sections of these sheaves over the open Bruhat cell
$U^-\hookrightarrow G/Q$. Note that $\bar{\mathcal{B}}$ is an $\bar{\cA}^\ell$-$\bar{\cA}^r$-bimodule.
There are $Q^-$-equivariant decompositions
\begin{equation}\label{eq:Bruhat_decomp}
\bar{\cA}^\ell=D(U^-)\otimes \cA_M^\ell, \qquad 
\bar{\cA}^r=D(U^-)\otimes \cA_M^r, \qquad \bar{\mathcal{B}}=D(U^-)\otimes \mathcal{B}_M,
\end{equation}
where $M$ acts 
diagonally and $U^-$ acts on the first factor only. The isomorphisms in (\ref{eq:Bruhat_decomp}) are compatible with algebra and bimodule structures.

We will need to define several auxilary categories of bimodules. Let $\mathcal{U}_{\gamma^\ell}, \mathcal{U}_{\gamma^r}$ denote the
central reductions of $U(\fg)$ corresponding to $\gamma^{\ell}, \gamma^r$. Choose a one-parameter subgroup $\nu:\C^\times \rightarrow Z(M)$ such that $\nu(t)$ acts on $\fu:=\operatorname{Lie}(U)$ by positive powers of $t$, and let $h$ denote the element of $\mathfrak{z}(\fm)$ corresponding to $\nu$. Consider the category
$\mathrm{OHC}(\mathcal{U}_{\gamma^\ell},\cA^r_M)$ consisting of all  $\mathcal{U}_{\gamma^\ell}$-$\cA^r_M$-bimodules
$V$ such that
\begin{enumerate}
\item $h$ acts diagonalizably on $V$,
and the set of eigenvalues is bounded from above, i.e. there are $c_1,\ldots,c_k\in \CC$ depending on $V$ such that every eigenvalue is of the form $c_i-n$ for $n\in \mathbb{Z}_{\geqslant 0}$ for some $i$.
\item Each $h$-eigenspace is a finite length Harish-Chandra $U(\fm)$-bimodule, where the action on the right is via the quantum co-moment map $\Phi_M^r: U(\mathfrak{m})\rightarrow \mathcal{A}_M^r$.
\end{enumerate}
Define the categories $\mathrm{OHC}(\mathcal{U}_{\gamma^r},\cA^r_M), \mathrm{OHC}(\bar{\cA}^\ell,\cA^r_M),
\mathrm{OHC}(\bar{\cA}^r, \cA^r_M)$, etc., in a similar fashion. These categories interpolate between a version of category $\mathcal{O}$
and the category of Harish-Chandra bimodules (justifying our notation).

Consider the grading on $\bar{\cA}^r$ defined by $\nu$. We can decompose
$\bar{\cA}^r$ as the direct sum of vector spaces $\bar{\cA}^r_{<0}\oplus \bar{\cA}^r_0\oplus \bar{\cA}^r_{>0}$,
where the subscript $<0$ (resp. $>0$) indicates the direct sum of homogeneous components of
negative (resp. positive) degree. It is easy to
see that $\bar{\cA}^r/\bar{\cA}^r\bar{\cA}^r_{>0}$ is naturally identified with $U(\mathfrak{u}^-)\otimes \cA^r_M$ as an $M$-equivariant $\left(U(\mathfrak{u}^-)\otimes \cA^r_M\right)-\cA^r_M$-bimodule. Also $\bar{\cA}^r/\bar{\cA}^r_{<0}\bar{\cA}^r\simeq \cA^r_M\otimes \CC[U^-]$, an $M$-equivariant isomorphism of $\cA^r_M$-$(\cA^r_M\otimes \C[U^-])$-bimodules.  

Similarly, we can consider the decompositions
$$\bar{\cA}^\ell=\bar{\cA}^\ell_{<0}\oplus \bar{\cA}^\ell_0\oplus \bar{\cA}^\ell_{>0}, \qquad 
U(\fg)=U(\fg)_{<0}\oplus U(\fg)_0\oplus U(\fg)_{>0}$$
as well as  the induced decompositions for $\mathcal{U}_{\gamma_\ell},\mathcal{U}_{\gamma_r}$.
Note that $U(\fg)/U(\fg)U(\fg)_{>0}\simeq U(\mathfrak{q}_-)$ and $U(\fg)/U(\fg)_{<0}U(\fg)\simeq U(\mathfrak{q})$.

We will need the following objects
\begin{align*}
&\Delta^{\A,r}(\A^r_M):=(\bar{\A}^r/\bar{\A}^r\bar{\A}^r_{>0})\otimes_{\A^r_M}\A^r_M \in \mathrm{OHC}(\bar{\A}^r, \A^r_M),\\
&\nabla^{\A,r}(\A^r_M):=\Hom^{\mathrm{fin}}_{\A^r_L}(\bar{\A}^r/\bar{\A}^r_{<0}\bar{\A}^r,\A^r_M)
\in \mathrm{OHC}(\bar{\A}^r, \A^r_M),\\
&\Delta^{\A,\ell}(\mathcal{B}_M):=(\bar{\A}^\ell/\bar{\A}^\ell\bar{\A}^\ell_{>0})\otimes_{\A^\ell_M}\mathcal{B}_M
\in \mathrm{OHC}(\bar{\A}^\ell, \A^r_M),\\
&\nabla^{\A,\ell}(\mathcal{B}_M):=\Hom^{\mathrm{fin}}_{\A^\ell_M}(\bar{\A}^\ell/\bar{\A}^\ell_{<0}\bar{\A}^\ell,\mathcal{B}_M)
\in \mathrm{OHC}(\bar{\A}^\ell, \A^r_M),\\
&\Delta^{\mathcal{U},r}(\A^r_M):=(U(\fg)/U(\fg)U(\fg)_{>0})\otimes_{U(\fm)}\A^r_M\in \mathrm{OHC}(\mathcal{U}_{\gamma^r}, \A^r_M),\\
&\nabla^{\mathcal{U},r}(\A^r_M):=\Hom^{\mathrm{fin}}_{U(\fm)}(U(\fg)/U(\fg)_{<0}U(\fg),\A^r_M)
\in \mathrm{OHC}(\mathcal{U}_{\gamma^r}, \A^r_M),\\
& \Delta^{\mathcal{U},\ell}(\mathcal{B}_M):=(U(\fg)/U(\fg)U(\fg)_{>0})\otimes_{U(\fm)}\mathcal{B}_M
\in \mathrm{OHC}(\mathcal{U}_{\gamma^\ell}, \A^r_M),\\
& \nabla^{\mathcal{U},\ell}(\mathcal{B}_M):=\Hom^{\mathrm{fin}}_{U(\fm)}(U(\fg)/U(\fg)_{<0}U(\fg),\mathcal{B}_M)
\in \mathrm{OHC}(\mathcal{U}_{\gamma^\ell}, \A^r_M).
\end{align*}
%
Here $\Hom^{\mathrm{fin}}$ denotes the direct sum of Hom's from graded components as in the definition
of (parabolic) dual Verma modules. It is easy to check that the objects above lie in the specified categories.

\begin{proof}[Proof of Proposition \ref{prop:HC_induction}]
The proof has several steps.

{\it Step 1}. The composition $U(\g)\rightarrow \A^r\rightarrow \bar{\A}^r$
gives rise to a forgetful functor $\mathrm{OHC}(\bar{\A}^r,\A^r_M)\rightarrow
\mathrm{OHC}(\mathcal{U}_{\gamma_r},\A^r_M)$. Since the homomorphism $U(\fg)\rightarrow \bar{\A}^r$
is $M$- and hence $\nu$-equivariant, there are natural
homomorphisms
\begin{equation}\label{eq:OHC_homomorphisms}
\Delta^{\mathcal{U},r}(\A^r_M)\rightarrow \Delta^{\A,r}(\A^r_M), \qquad 
\nabla^{\A,r}(\A^r_M)\rightarrow \nabla^{\mathcal{U},r}(\A^r_M)
\end{equation}
in
$\mathrm{OHC}(\mathcal{U}_{\gamma_r},\A^r_M)$. Note that
both $\Delta^{\mathcal{U},r}(\A^r_M),\Delta^{\A,r}(\A^r_M)$ are identified with
$U(\mathfrak{u}^-)\otimes \A^r_M$ and, under this identification, (\ref{eq:OHC_homomorphisms}) is the identity. So
$\Delta^{\mathcal{U},r}(\A^r_M)\xrightarrow{\sim}\Delta^{\A,r}(\A^r_M)$. For similar
reasons, $\nabla^{\A,r}(\A^r_M)\xrightarrow{\sim} \nabla^{\mathcal{U},r}(\A^r_M)$, and
$$\Delta^{\mathcal{U},\ell}(\mathcal{B}_M)\xrightarrow{\sim}\Delta^{\A,\ell}(\B_M),
\nabla^{\A,\ell}(\B_M)\xrightarrow{\sim} \nabla^{\U,\ell}(\B_M).$$

{\it Step 2}. We can define $\Delta^{\U,r}(\mathcal{B}'), \nabla^{\U,r}(\mathcal{B}')$
for every  Harish-Chandra $U(\fm)$-$\A_M^r$-bimodule $\mathcal{B}'$. These are objects in
$\OHC(U(\g),\A_M^r)$. They enjoy usual properties of Verma and dual Verma modules including:
\begin{enumerate}
\item If $\mathcal{B}'$ is simple, then $\Delta^{\U,r}(\mathcal{B}')$ has a unique simple quotient, 
$\nabla^{\U,r}(\mathcal{B}')$ has a unique simple sub, and these
simple objects coincide.
\item $\Delta^{\U,r}(\mathcal{B}')$ 
has finite length. 
\item $\Delta^{\U,r}(\mathcal{B}'), \nabla^{\U,r}(\mathcal{B}')$
are isomorphic as $U(\fm)$-$\A^r_M$-bimodules (to $U(\mathfrak{u}^-)\otimes \mathcal{B}'$) and hence have the same
composition series in $\OHC(U(\g),\A_M^r)$. 
\end{enumerate}
The proofs are standard. 

{\it Step 3}. We claim that $\Delta^{\U,r}(\A_M^r)$ is a simple object
in $\OHC(\U_{\gamma_r},\A^r_M)$. Assume the contrary. We first note that 
$I^r\otimes U(\fm)+ U(\g)\otimes I^r_M$ is a maximal ideal in 
$U(\g)\otimes U(\fm)$. 
Indeed, it is primitive because it coincides with the kernel of 
$U(\fg)\otimes U(\mathfrak{m})\rightarrow \cA^r\otimes \cA^r_M$ and the target is a quantization of a nilpotent cover. Every such kernel is primitive, see Proposition \ref{prop:propsofIbeta}. The associated variety of $I^r\otimes U(\fm)+ U(\g)\otimes I^r_M$ has the same dimension as the associated variety of the maximal ideal containing it. This can be deduced from 
Proposition \ref{prop:maximalitycriterion}. So 
$I^r\otimes U(\fm)+ U(\g)\otimes I^r_M$  is indeed maximal.

Thanks to Step 1, the action of $\U_{\gamma_r}\otimes 
\A^r_M$ on $\Delta^{\U,r}(\A_M^r)$ factors through $U(\g)/I^r\otimes \A_M^r$. 

It makes sense to talk about good filtrations on objects in $\OHC(U(\g),\A_M^r)$.
The object $\Delta^{\U,r}(\A_M^r)$ comes with a good filtration induced from the filtration on $\A_M^r$.
The associated graded for this filtration is the $\mathbb{C}[\mathcal{N}]\otimes \mathbb{C}[\widetilde{X}_M]$-algebra $\mathbb{C}[\mathfrak{u}^-]\otimes 
\mathbb{C}[\widetilde{X}_M]$, where the module structure is as follows: the homomorphism $\mathbb{C}[\mathcal{N}]\otimes \mathbb{C}[\widetilde{X}_M]\rightarrow \mathbb{C}[\mathfrak{u}^-]$ comes from the inclusion $\mathfrak{u}^-\hookrightarrow \mathcal{N}$, while the homomorphism
$\mathbb{C}[\mathcal{N}]\otimes \mathbb{C}[\widetilde{X}_M]\rightarrow \mathbb{C}[\widetilde{X}_M]$ comes from the diagonal embedding $\widetilde{X}_M\hookrightarrow \mathcal{N}\times\widetilde{X}_M$. 
It follows that the associated variety of $\Delta^{\U,r}(\A_M^r)$ in $\mathcal{N}\times
\widetilde{X}_M$ is the subvariety  $\mathfrak{u}^-\times \widetilde{X}_M\subset \mathcal{N}\times\widetilde{X}_M$. The multiplicity of $\Delta^{\U,r}(\A_M^r)$ is equal to $1$. 
Recall, see (2) from Step 2, that $\Delta^{\U,r}(\A_M^r)$ has finite length. 

From our assumption that  $\Delta^{\U,r}(\A_M^r)$ is not simple, it follows that 
there is a simple constituent $S$ of $\Delta^{\U,r}(\A_M^r)$ whose associated 
variety is properly contained in $\mathfrak{u}^-\times \widetilde{X}_M$. {It follows that the associated variety of $\Delta^{\U,r}(\A_M^r)$ viewed as a $U(\fg)$-$U(\fm)$-bimodule is properly contained in $\mathfrak{u}^-\times \overline{\OO}_M$, where, as usual, $\OO_M$ denotes the orbit covered by $\widetilde{\OO}_M$. In particular, the GK dimension of $S$ is less than $\dim \mathfrak{u}^-+\dim \OO_M$. On the other hand, the annihilator of $S$ in $U(\fg)\otimes U(\fm)$ must coincide with $I^r\otimes U(\fm)+ U(\g)\otimes I^r_M$ because the latter is a maximal ideal annihilating $S$. The GK dimension of the quotient of $U(\fg)\otimes U(\fm)$ by this annihilator is $\dim \OO+\dim \OO_M=2(\dim \mathfrak{u}^-+\dim \OO_M)$. We get a contradiction with Gabber's theorem, see, e.g., \cite[Theorem 9.11]{Krause_Lenagan}: for a finite length module over a reductive Lie algebra its GK dimension is at least half of the GK dimension of the quotient of the universal enveloping algebra by the annihilator. This contradiction finishes the proof of the claim that $\Delta^{\U,r}(\A_M^r)$ is a simple object
in $\OHC(\U_{\gamma_r},\A^r_M)$.}



{\it Step 4}. Combining the conclusion of Step 3 with (1) and (3) of Step 2, we see that 
there is an isomorphism $\Delta^{\U,r}(\A_M^r)\xrightarrow{\sim} \nabla^{\U,r}(\A_M^r)$.

{\it Step 5}. 
We now produce an injective $\bar{\A}^\ell$-$\bar{\A}^r$-bilinear map 
\begin{equation}\label{eq:bimodule_map}
\bar{\B}\rightarrow \Hom_{\A_M^r}(\nabla^{\A,r}(\A_M^r),
\nabla^{\A,\ell}(\B_M)).\end{equation} 
%
Recall that we can identify
$\bar{\B}$ with $D(U^-)\otimes \B_M$, $\bar{\A}^\ell$ with $D(U^-)\otimes \A_M^\ell$,
and $\bar{\A}^r$ with $D(U^-)\otimes \A_M^r$. Under these identifications,
we have  $\nabla^{\A,r}(\A_M^r)\simeq \Hom^{\mathrm{fin}}(\C[U^-], \A_M^r), 
\nabla^{\A,\ell}(\B_M)\simeq \Hom^{\mathrm{fin}}(\C[U^-],\B_M)$. 
The homomorphism (\ref{eq:bimodule_map}) sends $d\otimes b\in D(U^-)\otimes \B_M$ to the map 
$\varphi\mapsto b[\varphi\circ d]$.
This description
shows that (\ref{eq:bimodule_map}) is injective. 

{\it Step 6}. Note that $\B$ embeds into the $\operatorname{ad}(\g)$-finite part 
$\bar{\B}^{\g-\mathrm{fin}}$ of $\bar{\B}$. So $\B$ also embeds into 
the $\operatorname{ad}(\g)$-finite part of the target of (\ref{eq:bimodule_map}). 
Thanks to Step 4, we get an embedding
\begin{equation}\label{eq:bimod_map2}
\B\hookrightarrow \Hom_{\A_M^r}^{\g-\mathrm{fin}}(\Delta^{\U,r}(\A_M^r),\nabla^{\U,\ell}(\B_M))
\end{equation}
By Lemma \ref{lem:bimodinduction}, the target, viewed as a Harish-Chandra $U(\g)$-bimodule, coincides with $\operatorname{Ind}^G_M(\mathcal{B}_M)$. 

{\it Step 7}. We will show that (\ref{eq:bimod_map2}) is an 
isomorphism. For this it is sufficient to show that the $G$-multiplicities of its source and target coincide. The target is isomorphic as a $G$-representation to 
$\operatorname{AlgInd}^G_M(\B_M)$, see Proposition \ref{prop:propsofind}. Recall that $\B=\Gamma(\mathcal{B}^{\mathrm{loc}})$. The higher cohomology 
of this bimodule is a Harish-Chandra $\A^\ell$-$\A^r$-bimodule, supported 
away from the locus where the partial resolution $G\times^Q(\widetilde{X}_M\times \mathfrak{u})
\rightarrow \widetilde{X}$ is an isomorphism. In particular, the associated variety of the higher cohomology is a proper subvariety in $X\times \widetilde{X}_M$. Since 
$I^\ell,I^r$ are assumed to be maximal, this implies that the higher cohomology is $0$. So the $G$-multiplicities of $\B$ are equal to the $G$-multiplicities of the Euler characteristic of $R\Gamma(\mathcal{B}^{\mathrm{loc}})$ and thus of $R\Gamma(\operatorname{gr}\mathcal{B}^{\mathrm{loc}})$. So it suffices to show that $G$-multiplicities of the complex $R\Gamma(\operatorname{gr}\mathcal{B}^{\mathrm{loc}})$ are equal to the $G$-multiplicities of $\operatorname{AlgInd}^G_M(\mathcal{B}_M)$. This amounts to showing that for any finite dimensional $M$-representation $V$ there is an isomorphism of $G$-representations
\begin{equation}\label{eq:G_type_equality}
R\Gamma(G\times^Q \mathcal{V})\simeq \operatorname{AlgInd}^G_M V,
\end{equation} 
where $\mathcal{V}$ is, by definition, the $M$-equivariant sheaf on $\widetilde{X}_M \times \fq^{\perp}$ obtained by pull back from the sky-scraper sheaf with fiber $V$ at $0 \in \widetilde{X}_M$. For $\lambda\in \mathfrak{X}(\mathfrak{m})$, let
$\mathcal{V}_\lambda$ denote the sheaf on
$G\times^Q(\widetilde{X}_M\times\{\lambda\}\times \mathfrak{u})$ defined similarly to $\mathcal{V}$.
The $G$-multiplicities of 
$R\Gamma(G\times^Q \mathcal{V})$ are constant under flat deformation. In particular, the $G$-multiplicities of 
$R\Gamma(G\times^Q \mathcal{V}_\lambda)$ are independent of $\lambda$. For a generic $\lambda$, the sheaf $G\times^Q \mathcal{V}_\lambda$ coincides with the homogeneous bundle on $G/M$ with fiber $V$. There is a $G$-module isomorphism $\Gamma(G\times^Q \mathcal{V}_\lambda)\simeq \operatorname{AlgInd}^G_M(\mathcal{B}_M)$. This completes the proof.
\end{proof}

\chapter{Maximality calculations}\label{sec:maximality}

In this appendix, we show that for linear classical groups all unipotent ideals are maximal. We will also prove the maximality of certain auxilary ideals which appear in the deformation arguments in Chapter \ref{sec:unipbimod}. All of these proofs proceed in more or less the same fashion. Let $\widetilde{\mathbb{O}}$ be a $G$-equivariant nilpotent cover. Let $I:=I_{\beta}(\widetilde{\mathbb{O}}) \subset U(\fg)$, and let $\gamma := \gamma_{\beta}(\widetilde{\mathbb{O}}) \in \fh^*/W$ (cf. Definition \ref{def:unipotentideals}). By Proposition \ref{prop:propsofIbeta}, $I$ is a primitive ideal with infinitesimal character $\gamma$ and associated variety $\overline{\mathbb{O}}$. As explained in Section \ref{subsec:assvarmax}, $\gamma$ can be used to define reductive subgroups
$$L^{\vee}_{\gamma,0} \subset L^{\vee}_{\gamma} \subset G^{\vee}.$$
Consider the Richardson orbit $\mathbb{O}_{\gamma}^{\vee} := \Ind^{L^{\vee}_{\gamma}}_{L^{\vee}_{\gamma,0}} \{0\} \subset (\fl_{\gamma}^{\vee})^*$, and its dual orbit $\mathbb{O}_{\gamma} := \mathsf{D}(\mathbb{O}^{\vee}_{\gamma}) \subset \fl_{\gamma}^*$. By Proposition \ref{prop:maximalitycriterion}, $I$ is maximal if and only if the following condition is satisfied
\begin{equation}\label{eq:codimcodim}\codim(\mathbb{O},\cN) = \codim(\mathbb{O}_{\gamma},\cN_{L_{\gamma}}).\end{equation}
In special cases, $L_{\gamma}$ is a Levi subgroup of $G$ (for example, this is always the case in type $A$). In these cases, (\ref{eq:codimcodim}) is equivalent to the condition (see Remark \ref{rmk:maximalitycriterion})
$$\mathbb{O} = \Ind^G_{L_{\gamma}}\mathbb{O}_{\gamma}.$$
In practice, the latter condition is often easier to check.

Although the proofs in all types follow the same basic outline, there are a number of small differences which would make case-free proofs difficult. Examples of these differences include the parameterization of nilpotent orbits and covers, the combinatorial formulas for unipotent infinitesimal characters, and the combinatorics of BVLS duality. Instead, we give detailed proofs in types A and C and then indicate the necessary changes for types B and D.

\section{Calculations in type A}\label{subsec:maximalityA}

Suppose $G = \mathrm{SL}(n)$. The following result is used in the proof of Proposition \ref{prop:bimodsA}.
\begin{prop}\label{prop:maximalityA}
Suppose $\widetilde{\mathbb{O}}$ is a $G$-equivariant nilpotent cover. Then $I_0(\widetilde{\mathbb{O}}) \subset U(\mathfrak{g})$ is a maximal ideal.
\end{prop}

\begin{proof}
Choose a Levi subgroup $M = \mathrm{S}(\mathrm{GL}(a_1) \times ... \times \mathrm{GL}(a_t)) \subset G$ and a birationally rigid cover
$$\widetilde{\mathbb{O}}_M = \widetilde{\mathbb{O}}_M^1 \times ... \times \widetilde{\mathbb{O}}_M^t$$
such that $\widetilde{\mathbb{O}} = \mathrm{Bind}^G_M \widetilde{\mathbb{O}}_M$. Write $p^i$ for the partitions of $a^i$ corresponding to $\mathbb{O}_M^i$ and write $p$ for the partition of $n$ corresponding to $\mathbb{O}$. By Proposition \ref{prop:nocodim2leaves} and Remark \ref{rmk:SLlevisbirigid}, there is a positive integer $d \mid n$ and a partition $m$ of $\frac{n}{d}$ such that $p^i = (d^{m_i})$ for $1 \leq i \leq t$. Since $\mathbb{O} = \mathrm{Ind}^G_M \mathbb{O}_M$, this implies $p = dm$.

By Propositions \ref{prop:Ibetainduced} and \ref{prop:centralcharacterbirigidcover}
$$\gamma:=\gamma_0(\widetilde{\mathbb{O}})  = \gamma_0(\widetilde{\mathbb{O}}_M)= (\frac{\rho(a_1)}{d}, \frac{\rho(a_2)}{d}, ..., \frac{\rho(a_t)}{d}).$$
We will consider separately the cases of odd and even $d$.

\begin{description}
    \item[d odd] Let $e$ (resp. $o$) denote the transpose of the subpartition of $a$ consisting of even (resp. odd) parts. For example, if $a=(15,12^2,9)$, then $e=(2^{12})$, and $o=(2^9,1^6)$. Let $S=\{0,1,...,2d-1\}$ and let $S_e$ (resp. $S_o$) denote the set of even (resp. odd) elements of $S$. By an easy computation
    $$L^{\vee}_{\gamma} = \mathrm{S}\left(\prod_{x \in S_o} \mathrm{GL}\left(\frac{|e|}{d}\right) \times \prod_{x \in S_e} \mathrm{GL}\left(\frac{|o|}{d}\right)\right)$$
    and
    $$L^{\vee}_{\gamma,0} = \mathrm{S}\left(\prod_{x \in S_o} \left(\prod_{y \equiv x \text{ mod }d} \mathrm{GL}(e_y)\right) \times \prod_{x \in S_e} \left(\prod_{y \equiv x \text{ mod }d} \mathrm{GL}(o_y)\right) \right).$$
    Note that $L_{\gamma} \simeq L_{\gamma}^{\vee}$ is a Levi subgroup of $G$. Thus it suffices to show that $\mathbb{O} = \Ind^G_{L_{\gamma}} \mathbb{O}_{\gamma}$. 
    
    For each $x \in S_o$, let $r_x(e)$ denote the subpartition of $e$ consisting of parts $e_y$ with $y \equiv x \text{ mod }d$. For example, if $e=(3^6,2^6,1^3)$ and $d=3$, then $r_1(e) = (3^2,2^2,1)$. For each $x \in S_e$, define $r_x(o)$ similarly. Then by Proposition \ref{prop:inductionclassical}
    $$\mathbb{O}_{\gamma}^{\vee} = \mathrm{Ind}^{L^{\vee}_{\gamma}}_{L^{\vee}_{\gamma,0}} \{0\} = \prod_{x \in S_o} r_x(e)^t \times \prod_{x \in S_e} r_x(o)^t \subset \cN_{L^{\vee}_{\gamma}},$$
    where we have identified orbits (here and for the remainder of this section) with their corresponding partitions. By Proposition \ref{prop:BVduality}
    $$\mathbb{O}_{\gamma} = \mathsf{D}(\mathbb{O}_{\gamma}^{\vee}) = \prod_{x \in S_o} r_x(e) \times \prod_{x \in S_e} r_x(o) \subset \cN_{L_{\gamma}},$$
    and therefore
    $$\mathrm{Ind}^G_{L_{\gamma}} \mathbb{O}_{\gamma} = \mathrm{Ind}^G_{\mathrm{S}(\mathrm{GL}(|e|) \times \mathrm{GL}(|o|))} e \times o = e + o = p$$
    as desired.
    
    \item[d even] Let $S = \{0,1,...d-1\}$. Then
    $$L^{\vee}_{\gamma} = \mathrm{S}\left(\prod_{x \in S} \mathrm{GL}\left(|m|\right) \right), \qquad L^{\vee}_{\gamma,0} = \mathrm{S}\left(\prod_{x \in S} \left(\prod_{y \equiv x \text{ mod }d} \mathrm{GL}((a^t)_y) \right) \right).$$
    For each $x \in S$, let $r_x(a^t)$ denote the subpartition of $a^t$ consisting of parts $a^t_y$ with $y \equiv x \text{ mod }d$. Hence
    $$\mathbb{O}_{\gamma}^{\vee} = \mathrm{Ind}^{L^{\vee}_{\gamma}}_{L^{\vee}_{\gamma,0}} \{0\} = \prod_{x \in S} r_x(a^t)^t \subset \cN_{L^{\vee}_{\gamma}},$$
    and therefore
    $$\mathbb{O}_{\gamma} = \mathsf{D}\left(\mathbb{O}_{\gamma}^{\vee}\right) =  \prod_{x \in S} r_x(a^t) \subset \cN_{L_{\gamma}}.$$
    Now
    $$\mathrm{Ind}^G_{L_{\gamma}} \mathbb{O}_{\gamma} = \mathsf{D}\left(\frac{a}{d}\right)^t
    = p$$
as desired.
\end{description}
\end{proof}

\subsection{Auxilary ideals}

We now proceed to proving the maximality of certain auxilary ideals attached to birationally rigid covers. These results are used in the proof of Proposition \ref{prop:bimodsA}.

Let $\mathbb{O}$ be an orbit corresponding to a partition $p = (d^m)$ of $n$. Recall that the universal $G$-equivariant cover $\widetilde{\mathbb{O}} \to \mathbb{O}$ is birationally rigid and $I_0(\widetilde{\mathbb{O}}) = I_{\delta}(\mathbb{O})$, where
   $$\delta = (\underbrace{\frac{d-1}{2d},\ldots, \frac{d-1}{2d}}_{m}, \underbrace{\frac{d-3}{2d},\ldots, \frac{d-3}{2d}}_{m},\ldots, \underbrace{\frac{1-d}{2d},\ldots, \frac{1-d}{2d}}_{m}).$$
To prove the unitarity of the spherical bimodule $U(\fg)/I_{\delta}(\mathbb{O}) \in \unip_{\widetilde{\mathbb{O}}}(G)$, we need to prove the maximality not only of $I_{\delta}(\mathbb{O})$, but of all ideals in the one-parameter family $\{I_{t\delta}(\mathbb{O})\}_{t \in [0,1]}$.

\begin{lemma}\label{lem:auxidealsA1}
For every $t \in [0,1]$, $I_{t\delta}(\mathbb{O}) \subset U(\fg)$ is a maximal ideal.
\end{lemma}

\begin{proof}
If $t \in \{0,1\}$, $I_{t\delta}(\mathbb{O})$ is maximal by Proposition \ref{prop:maximalityA}. Assume $t \in (0,1)$. By Proposition \ref{prop:unipotentcentralchars} and a direct computation
\begin{align*}
\gamma:=\gamma_{t\delta}(\mathbb{O}) &= \rho(\fl) + t\delta\\
                   &=(\underbrace{\frac{m-1}{2},\ldots, \frac{1-m}{2}}_{m}, \underbrace{\frac{m-1}{2},\ldots, \frac{1-m}{2}}_{m},\ldots, \underbrace{\frac{m-1}{2},\ldots, \frac{1-m}{2}}_{m})\\
                   &+ t(\underbrace{\frac{d-1}{2d},\ldots, \frac{d-1}{2d}}_{m}, \underbrace{\frac{d-3}{2d},\ldots, \frac{d-3}{2d}}_{m},\ldots, \underbrace{\frac{1-d}{2d},\ldots, \frac{1-d}{2d}}_{m}).
\end{align*}
Hence, the difference between two entries of $\gamma$ is of the form $a + t\frac{b}{d}$ for nonnegative integers $a \leq m-1$ and $b \leq d-1$. The difference is integral if and only if $b=0$ and zero if and only if $a=b=0$. Hence,
$$L^{\vee}_{\gamma} = \mathrm{S}(\mathrm{GL}(m)^d), \qquad L^{\vee}_{\gamma,0} = H^{\vee},$$
where $H^{\vee}$ is the (dual) maximal torus. Thus,
$$\mathbb{O}_{\gamma}^{\vee} = \Ind^{L^{\vee}_{\gamma}}_{L^{\vee}_{\gamma,0}} \{0\} = (m) \times ... \times (m),$$
and therefore
$$\mathbb{O}_{\gamma} = \mathsf{D}(\mathbb{O}_{\gamma}^{\vee}) = \{0\} \times ... \times \{0\}.$$
Now
$$\Ind^G_{L_{\gamma}} \mathbb{O}_{\gamma} = (d^m) = p$$
as desired.
\end{proof}

There are $d-1$ bimodules in $\unip_{\widetilde{\mathbb{O}}}(G)$ which are distinct from the spherical bimodule $U(\fg)/I_{\delta}(\mathbb{O})$. To prove that these bimodules are unitary, we will need to consider a different set of auxilary ideals. Fix $j \in \{1,...,d-1\}$ and let
$$\tau_j = (\underbrace{\frac{d-j}{d},...,\frac{d-j}{d}}_{mj},\underbrace{-j,...,-j}_{m(d-j)}).$$
Consider the weights
$$\delta_{\ell}(t) := \delta -\tau_j + \frac{1}{2}(t-1)(2\delta-\tau_j), \qquad \delta_r(t) := \delta +\frac{1}{2}(t-1)(2\delta-\tau_j).$$

\begin{lemma}\label{lem:auxidealsA2}
For every $t \in [0,1]$ and $j\in \{1,...,d-1\}$, both  $I_{\delta_{\ell}(t)}(\mathbb{O}) \subset U(\fg)$ and $I_{\delta_r(t)}(\mathbb{O}) \subset U(\fg)$ are maximal ideals.
\end{lemma}

\begin{proof}
We will prove the assertion only for $\delta_r(t)$. The proof for $\delta_{\ell}(t)$ is analogous. Note that $\delta_r(1) = \delta$. Hence, $I_{\delta_r(1))}(\mathbb{O}) = I_0(\widetilde{\mathbb{O}})$ is maximal by Proposition \ref{prop:maximalityA}. On the other extreme, $\delta_r(0) = \frac{\tau_j}{2}$. By Proposition \ref{prop:unipotentcentralchars} and a direct computation
\begin{align*}
\gamma:=\gamma_{\delta_r(0)} &= \rho(\fl) + \frac{\tau_j}{2}\\
                   &= (\underbrace{\frac{n-j}{2d},\frac{n-2d-j}{2d},...,\frac{2d-n-j}{2d},\ldots,\frac{n-j}{2d},\frac{n-2d-j}{2d},\ldots,\frac{2d-n-j}{2d}}_{mj},\\
                   &\quad \underbrace{\frac{n-d-j}{2d}, \frac{n-3d-j}{2d},\ldots,\frac{d-n-j}{2d},\ldots, \frac{n-d-j}{2d}, \frac{n-3d-j}{2d},\ldots,\frac{d-n-j}{2d}}_{m(d-j)}).
\end{align*}
Hence, 
$$L^{\vee}_{\gamma} = \mathrm{S}(\mathrm{GL}(mj) \times \mathrm{GL}(m(d-j))), \qquad 
L^{\vee}_{\gamma,0} = \mathrm{S}(\mathrm{GL}(j)^m \times \mathrm{GL}(d-j)^m).$$
So,
$$\mathbb{O}^{\vee}_{\gamma} = \Ind^{L_{\gamma}^{\vee}}_{L^{\vee}_{\gamma,0}} \{0\} = (m^j) \times (m^{d-j}),$$
and therefore
$$\mathbb{O}_{\gamma} = \mathsf{D}(\mathbb{O}^{\vee}_{\gamma}) = (j^m) \times ((d-j)^m).$$
Now,
$$\Ind^G_{L_{\gamma}} \mathbb{O}_{\gamma} = (j^m) + ((d-j)^m) = (d^m) =p$$
as desired. 

Next, suppose $t \in (0,1)$. By Proposition \ref{prop:unipotentcentralchars} and a direct computation
\begin{align*}
\gamma_{\delta_r(t)}:=\gamma_{\delta_r(t)}(\mathbb{O}) &= \rho(\fl) + (1-t)\frac{\tau_j}{2}  + t\delta                             \\                                &=(\underbrace{\frac{m-1}{2},\ldots, \frac{1-m}{2}}_{m}, \underbrace{\frac{m-1}{2},\ldots, \frac{1-m}{2}}_{m},\ldots, \underbrace{\frac{m-1}{2},\ldots, \frac{1-m}{2}}_{m})\\
                             &+ (\underbrace{\frac{(1-t)(d-j)}{2d},\ldots, \frac{(1-t)(d-j)}{2d}}_{mj}, \underbrace{\frac{-(1-t)j}{2d},...,\frac{-(1-t)j}{2d}}_{m(d-j)})\\    &+(\underbrace{\frac{t(d-1)}{2d},\ldots,\frac{t(d-1)}{2d}}_m, \ldots,\underbrace{\frac{t(1-d)}{2d},\ldots,\frac{t(1-d)}{2d}}_m).
\end{align*}
The difference between two of the first $mj$ entries of $\gamma_{\delta_r(t)}$ is of the form
$$a+\frac{tb}{d}, \qquad a,b \in \ZZ \quad |a| \leq m-1 \quad |b|\leq j-1.$$
The difference between two of the last $m(d-j)$ entries is of the form
$$a+\frac{tb}{d}, \qquad a,b\in \ZZ \quad |a| \leq m-1 \quad |b| \leq d-j-1.$$
The difference between one of the first $mj$ entries and one of the last $m(d-j)$ entries is of the form
$$a+\frac{tb}{d} + \frac{(1-t)}{2}, \qquad a,b\in \ZZ \quad |a| \leq m-1 \quad |b| \leq d-1.$$
From these formulas we deduce
$$L^{\vee}_{\gamma_{\delta_r(t)}} = \mathrm{S}(\mathrm{GL}(m)^d) \qquad L^{\vee}_{\gamma_{\delta_r(t)},0} = H^{\vee},$$
where $H^{\vee}$ denotes the (dual) maximal torus. Now, proceed as in the proof of Lemma \ref{lem:auxidealsA1}.
\end{proof}

\section{Calculations in type C}\label{subsec:maximalityC}

Suppose $G = \mathrm{Sp}(2n)$. Our main result in this section is the following. It is used in the proof of Proposition \ref{prop:bimodsBCD}.

\begin{prop}\label{prop:maximalitytypeC}
Suppose $\widetilde{\mathbb{O}}$ is a $G$-equivariant cover. Then $I_0(\widetilde{\mathbb{O}}) \subset U(\fg)$ is a maximal ideal.
\end{prop}

We will prove this proposition by induction (on the number of $\mathrm{GL}$-factors in the Levi subgroup of $G$ from which $\widetilde{\mathbb{O}}$ is birationally induced). The base case is the following.

\begin{lemma}\label{lem:maximalitybirigidtypeC}
Suppose $\widetilde{\mathbb{O}}$ is a birationally rigid $G$-equivariant nilpotent cover. Then $I_0(\widetilde{\mathbb{O}}) \subset U(\fg)$ is a maximal ideal.
\end{lemma}

\begin{proof}
Let $p$ be the partition of $2n$ corresponding to $\mathbb{O}$ and let $q:=p^t$. By Proposition \ref{prop:nocodim2leaves}, $p$ satisfies
\begin{itemize}
    \item $p_i \leq p_{i+1}+2$ for all $i$.
    \item If $p_i$ is odd, then $p_i \leq p_{i+1}+1$.
\end{itemize}
Thus, $q$ satisfies
\begin{itemize}
    \item all parts in $q$ occur with multiplicity $\leq 2$.
    \item If $i$ is even, then $q_i \neq q_{i+1}$.
\end{itemize}
Let $x \subset q$ (resp. $y \subset q$) be the subpartition of multiplicity 1 (resp. 2) parts. By Proposition \ref{prop:centralcharacterbirigidcover}
$$\gamma_0(\widetilde{\mathbb{O}}) = \rho^+(g(y) \cup f_C(x))$$
(the notations $f_C(x)$, $g(y)$, and $\rho^+$ are explained in Definitions \ref{def:xy}, \ref{def:fBfC}, and \ref{def:rhoplus}). For any partition $r$, let $e(r)$ (resp. $l(r)$) denote the partitions obtained from $r$ by adding (resp. deleting) a single box from the first (resp. last) part in $r$. In \cite{McGovern1994}, Mcgovern considers the finite set of infinitesimal characters
$$Q(\mathbb{O}) = \{\rho^+(q) \mid C(l(q)) = p\},$$
where $C$ denotes the $C$-collapse. By \cite[Thm 5.1]{McGovern1994}), $V(I_{\mathrm{max}}(\gamma)) = \overline{\mathbb{O}}$ for every $\gamma \in Q(\mathbb{O})$. Thus, it suffices to show
\begin{equation}\label{eq:Clp}
C(l(g(y) \cup f_C(x))^t) = p.\end{equation}
We begin by providing an alternative description of the partition $g(y) \cup f_C(x)$. First, decompose $q$ into consecutive subpartitions (called \emph{blocks})
$$q = (\underbrace{q_1,...q_{l_1}}_{b^1}, \underbrace{q_{l_1+1}, ..., q_{l_1+l_2+1}}_{b^2},...,\underbrace{q_{2n-l_t+1},...,q_{2n}}_{b^s})$$
such that
\begin{itemize}
    \item For every $i$ ($1 \leq i \leq s$) and even $j$ ($1 \leq j \leq l_i-1$), $b^i_j = b^i_{j+1}+1$.
    \item For every $i$ ($1 \leq i \leq s-1$), $b^i_{l_i} \geq b^{i+1}_1+2$.
\end{itemize}
For example, if $p=(6,5^4,4^3,2,1^2)$, then $q=(11,9,8^2,5,1)$, $b^1=(11,9,8^2)$, and $b^2=(5,1)$. Note that if $\#q$ is even, then $l_i$ is even for $1 \leq i \leq s$. Otherwise, $l_i$ is even for $1 \leq i \leq s-1$, $l_s$ is odd, and $b^s$ contains (an odd number of) parts of multiplicity 1, including the smallest part of $b^s$, namely $b^s_{l_s}$. For $1 \leq i \leq s$, define 
$$\widetilde{b}^i := e(l(b^i)).$$
Since $b^i_{l_i} \geq b_1^{i+1} + 2$, we have $\widetilde{b}_{l_i}^i \geq \widetilde{b}_1^{i+1}$ for $1 \leq i \leq s-1$. So we can define the partition
$$\widetilde{b}:= (\widetilde{b}^1,\widetilde{b}^2,...,\widetilde{b}^s).$$
We will define one final partition $\vardbtilde{b}$, as follows: if $\#q$ is even, define $\vardbtilde{b}=(\widetilde{b},1)$. Otherwise, define
$$\vardbtilde{b}(q) = (\widetilde{b}^1,...,\widetilde{b}^{s-1},e(b^s))$$
(i.e. add a single box to the smallest part in $\widetilde{b}$). For example, if we choose $q=(11,9,8^2,5,1)$ as above, then
$$b^1=(11,9,8^2), \ b^2=(5,1), \qquad \widetilde{b}^1 = (12,9,8,7), \ \widetilde{b}^2 = (6), \qquad \vardbtilde{b} = (12,9,8,7,6,1).$$
If, on the other hand, $q=(6,4,2)$, then 
$$b^1 = (6,4), \ b^2 = (2), \qquad \widetilde{b}^1 = (7,3), \ \widetilde{b}^2 = (2), \qquad \vardbtilde{b} = (7,3^2).$$
We will sometimes write $\vardbtilde{b}^i(q)$ and $\vardbtilde{b}(q)$ to indicate dependence on $q$. We claim, first of all, that $g(y) \cup f_C(x) = \vardbtilde{b}(q)$. We proceed by induction on $s$, the number of blocks in $q$. For notational convenience, let $h(q) := f_C(x) \cup g(y)$. The base case (namely, $s=1$) is an easy exercise, which we leave to the reader. Now suppose $s \geq 2$ and let $q' := (b^2,...,b^s)$. If $b^1$ contains no parts of multiplicity 1, then by definition $h(b^1)$ contains 1 with positive multiplicity. Let $h'(b^1)$ be the partition obtained from $h(b^1)$ by deleting one such part. If, on the other hand, $b^1$ contains at least one part with multiplicity 1, let $b^1_r$ be the maximal such and define $h'(b^1)$ be removing one box from $h(b^1)_r$. Then clearly
$$h(q) = h(b^1,q') = (h'(b^1),h(q')).$$
By the $s=1$ case, discussed above, $h(b^1) = \vardbtilde{b}(b^1)$, and hence $h'(b^1) = \widetilde{b}(b^1) = \widetilde{b}^1$. By induction, $h(q') = \vardbtilde{b}(q')$. Thus,
\begin{equation}\label{eq:hvardbtilde}
h(q) = (\widetilde{b}^1, \vardbtilde{b}(q')) = (\widetilde{b}^1,..., \widetilde{b}^{s-1},e(b^t)) = \vardbtilde{b}(q)\end{equation}
as asserted. Now, comparing (\ref{eq:Clp}) and (\ref{eq:hvardbtilde}), it suffices to show that $C(l((\vardbtilde{b}(q))^t)) = p$. We will do so, again, by induction on $s$. First, suppose $s=1$. There are two cases to consider. If $\#q$ is even, then $\vardbtilde{b}(q) = (\widetilde{b}(q),1) = (e(l(q)), 1)$. Thus, $l(\vardbtilde{b}(q)^t)$ is obtained from $p$ by adding a box to the first row and deleting a box from the first column. By the algorithm for the $C$-collapse (see \cite[Lem 6.3.3]{McGovern1994}) it is clear that $C(l(\vardbtilde{b}(q)^t)) = p$. If, on the other hand, $\#q$ is odd, then $\vardbtilde{b}(q) = p$. So again, $C(l(\vardbtilde{b}(q)^t)) = p$. 

Now suppose $s \geq 2$. As above, let $q' := (b^2, ..., b^s)$ and recall that
\begin{equation}\label{eq:reminder}
\vardbtilde{b}(q) = (\widetilde{b}^1(q),\vardbtilde{b}(q')).\end{equation}
Write $\mathrm{col}(\vardbtilde{b}(q'))$ for the number of columns in $\vardbtilde{b}(q')$ (note that $\mathrm{col}(\vardbtilde{b}(q')) = \mathrm{col}(q')+1$). Then from (\ref{eq:reminder}), we obtain
$$\vardbtilde{b}(q)^t_i = \vardbtilde{b}(q')^t_i + l_1, \qquad i \leq \mathrm{col}(\vardbtilde{b}(q')).$$
Since $l_1$ is even, this implies
$$C(l(\vardbtilde{b}(q)^t))_i = C(l(\vardbtilde{b}(q')^t))_i + l_1, \qquad i \leq \mathrm{col}(\vardbtilde{b}(q')) -1.$$
Hence by the induction hypothesis (applied to $q'$) we obtain
$$C(l(\vardbtilde{b}(q)^t))_i = (q')^t_i + l_1 = p_i, \qquad i \leq \mathrm{col}(\vardbtilde{b}(q')) -1.$$
For $i = \mathrm{col}(\vardbtilde{b}(q'))$, the row $\vardbtilde{b}(q)^t_i$ is odd, appearing with odd multiplicity. Hence,
$$C(l(\vardbtilde{b}(q)^t))_i = \vardbtilde{b}(q)^t -1 = p_i, \qquad i = \mathrm{col}(\vardbtilde{b}(q')).$$
It remains to show that 
\begin{equation}\label{eqn:maximalityeqn1}
C(l(\vardbtilde{b}(q)^t))_i =  p_i, \qquad i > \mathrm{col}(\vardbtilde{b}(q')).\end{equation}
But indeed
$$\vardbtilde{b}(q)^t_i = (\vardbtilde{b}^1)^t_i, \qquad i>\mathrm{col}(\vardbtilde{b}(q')).$$
So (\ref{eqn:maximalityeqn1}) follows from the base case of the induction (applied to $q = b^1$).
\end{proof}

The induction step in our proof of Proposition \ref{prop:maximalitytypeC} will depend on the following definition.

\begin{definition}\label{def:triangular}
Let $\mu=(\mu_1,...,\mu_r)$ and $\nu = (\nu_1,...,\nu_s)$ be partitions and assume $\mu_1$ is even. Define the $(|\mu|+|\nu|-\mu_1)$-tuple
\begin{align*}v(\mu,\nu) = (&\underbrace{r-1,...,r-1}_{\mu_r},\underbrace{r-2,...,r-2}_{\mu_{r-1}}, ...,\underbrace{1,...,1}_{\mu_2},\underbrace{0,...,0}_{\frac{\mu_1}{2}},\\
&\underbrace{\frac{2s-1}{2},...,\frac{2s-1}{2}}_{\nu_s},\underbrace{\frac{2s-3}{2},...,\frac{2s-3}{2}}_{\nu_{s-1}},...,\underbrace{\frac{1}{2},...,\frac{1}{2}}_{\nu_1}).\end{align*}
We say that a tuple $v$ is \emph{triangular} if there are partitions $\mu$ and $\nu$ as above such that $v=v(\mu,\nu)$.
\end{definition}

Now, choose a Levi subgroup of the form
$$M = \mathrm{GL}(a) \times G(m) \subset G, \qquad a+m=n.$$
Let $\gamma' = v(\mu',\nu')$ be a triangular $m$-tuple, and let
$$\gamma = (\frac{a-1}{2}, \frac{a-3}{2}, ..., \frac{1-a}{2},\gamma').$$
Regard $\gamma$ (resp. $\gamma'$) as a weight (in standard coordinates) for $G$ (resp. $G(m)$). Form the reductive groups $L_{\gamma}$ and $L_{\gamma'}$, and the nilpotent orbits $\mathbb{O}_{\gamma} \subset \fl_{\gamma}^*$ and $\mathbb{O}_{\gamma'} \subset \fl_{\gamma'}^*$. 

\begin{lemma}\label{lem:inductionsteptypeC}
The $n$-tuple $\gamma$ is ($W$-conjugate to) a triangular one, and there is an equality
$$\codim(\mathbb{O}_{\gamma}, \cN_{L_{\gamma}}) = \codim(\mathbb{O}_{\gamma'},\cN_{L_{\gamma'}}) + \dim(\cN_{GL(a)}).$$
\end{lemma}

\begin{proof}
First, suppose $a$ is odd. Define partitions $\mu$ and $\nu$ by

\begin{equation}\label{eqn:codim3}
\mu_i = \left\{
        \begin{array}{ll}
            \mu'_i+2 & \quad i=1\\
            \mu'_i+1 & \quad 2 \leq i \leq a\\
            \mu'_i & \quad \text{else}
        \end{array}\right. \qquad \nu = \nu'.
\end{equation}
It is easy to see that $\gamma$ is $W$-conjugate to $v(\mu,\nu)$. This proves the first part of the proposition.

Define the integers
$$R' = |\mu'| - \frac{1}{2}\mu'_1, \qquad R = |\mu| - \frac{1}{2}\mu_1, \qquad S' = |\nu'| = |\nu|.$$
Then
$$
L^{\vee}_{\gamma'} = \mathrm{SO}(2R'+1) \times \mathrm{SO}(2S'), \qquad L^{\vee}_{\gamma',0} = \mathrm{SO}(\mu'_1+1) \times \prod_{i=2}^r \mathrm{GL}(\mu'_i) \times \prod_{j=1}^s \mathrm{GL}(\nu'_j),
$$
and
$$
L^{\vee}_{\gamma} = \mathrm{SO}(2R+1) \times \mathrm{SO}(2S'), \qquad L^{\vee}_{\gamma,0} = \mathrm{SO}(\mu_1+1) \times \prod_{i=2}^r \mathrm{GL}(\mu_i) \times \prod_{j=1}^s \mathrm{GL}(\nu'_i).
$$
Hence, by (\ref{eqn:codim3}) and Proposition \ref{prop:inclusionclassical}
$$\mathbb{O}^{\vee}_{\gamma} = \Ind^{L^{\vee}_{\gamma}}_{L^{\vee}_{\gamma,0}}\{0\} =  \mathrm{Sat}^{L^{\vee}_{\gamma}}_{\mathrm{GL}(a) \times L^{\vee}_{\gamma'} } \left(\mathbb{O}_{\mathrm{GL}(a)}^{\mathrm{prin}} \times \mathbb{O}_{\gamma'}^{\vee}  \right).$$
Using Proposition \ref{prop:inclusioninduction} we obtain
$$\mathbb{O}_{\gamma} = \mathsf{D}(\mathbb{O}^{\vee}_{\gamma}) = \mathrm{Ind}^{L_{\gamma}}_{\mathrm{GL}(a) \times L_{\gamma'} } \left(\{0\} \times \mathbb{O}_{\gamma} \right).$$
Since induction preserves codimension (see Proposition \ref{prop:propertiesofInd}(iii)), this implies
$$\codim(\mathbb{O}_{\gamma}, \cN_{L_{\gamma}}) = \codim(\mathbb{O}_{\gamma'},\cN_{L_{\gamma'}}) + \dim(\cN_{GL(a)})$$
as asserted. Next, suppose $a$ is even. Define
\begin{equation}\label{eqn:codim4}
\mu = \mu', \qquad 
\nu_j = \left\{
        \begin{array}{ll}
            \nu'_j+1 & \quad 1 \leq j \leq a\\
            \nu'_j & \quad \text{else}.
        \end{array}\right. 
\end{equation}
Then $\gamma$ is $W$-conjugate to $v(\mu,\nu)$. Set
$$R' = |\mu'|-\frac{1}{2}\mu_1' = |\mu|-\frac{1}{2}\mu_1, \qquad S' = |\nu'|, \qquad S = |\nu|.$$
Then
$$
L^{\vee}_{\gamma} = \mathrm{SO}(2R'+1) \times \mathrm{SO}(2S), \qquad L^{\vee}_{\gamma,0} = \mathrm{SO}(\mu'_1+1) \times \prod_{i=2}^r \mathrm{GL}(\mu'_i) \times \prod_{j=1}^s \mathrm{GL}(\nu_j).
$$
So again
$$\mathbb{O}^{\vee}_{\gamma} = \mathrm{Sat}^{L^{\vee}_{\gamma}}_{\mathrm{GL}(a) \times L^{\vee}_{\gamma'} } \left(\mathbb{O}_{\mathrm{GL}(a)}^{\mathrm{prin}} \times \mathbb{O}^{\vee}_{\gamma'}  \right).$$
Now the equality follows exactly as above.
\end{proof}

We are now prepared to prove Proposition \ref{prop:maximalitytypeC}.

\begin{proof}[Proof of Proposition \ref{prop:maximalitytypeC}]
Choose a Levi subgroup $M = \mathrm{GL}(a_1) \times ... \times \mathrm{GL}(a_t) \times G(m)$ and a birationally rigid cover
$$\widetilde{\mathbb{O}}_M =  \{0\} \times ... \{0\} \times \widetilde{\mathbb{O}}_{G(m)}$$
such that $\widetilde{\mathbb{O}} = \mathrm{Bind}^G_M \widetilde{\mathbb{O}}_M$. We will prove the following statement by induction on $t$:
$$I_0(\widetilde{\mathbb{O}}) \subset U(\fg) \text{ is maximal, and } \gamma_0(\widetilde{\mathbb{O}}) \text{ is the } W\text{-orbit of a triangular weight.}$$
If $t=0$, then $\widetilde{\mathbb{O}}$ is birationally rigid. Hence $I_0(\widetilde{\mathbb{O}}) \subset U(\fg)$ is maximal by Lemma \ref{lem:maximalitybirigidtypeC}. Any weight of the form $\rho^+(q)$ (for some partition $q$) is triangular. Thus, $\gamma_0(\widetilde{\mathbb{O}})$ is triangular by Proposition \ref{prop:centralcharacterbirigidcover}. 

Now suppose $t \geq 1$. Let $G':=G(n-a_1)$. Consider the Levi subgroup 
$$M' = GL(a_2) \times ... \times GL(a_t) \times G(m) \subset G'$$
and the $M'$-equivariant nilpotent cover 
$$\widetilde{\mathbb{O}}_{M'} :=  \{0\} \times ... \times \{0\} \times \widetilde{\mathbb{O}}_{G(m)}.$$
Consider the $G'$-equivariant nilpotent cover $\widetilde{\mathbb{O}}' = \mathrm{Bind}^{G'}_{M'} \widetilde{\mathbb{O}}_{M'}$. By the transitivity of induction 
$$\widetilde{\mathbb{O}} = \mathrm{Bind}^G_{\mathrm{GL}(a_1) \times G'} \left(\{0\} \times \widetilde{\mathbb{O}}' \right)$$
and hence (writing $\mathbb{O}$ and $\mathbb{O}'$ for the underlying orbits)
$$\mathbb{O} = \Ind^G_{\mathrm{GL}(a_1) \times G'} \left(\{0\} \times \mathbb{O}'\right).$$
For notational convenience, let $\gamma := \gamma_0(\widetilde{\mathbb{O}})$ and $\gamma' := \gamma_0(\widetilde{\mathbb{O}}')$. By Proposition \ref{prop:Ibetainduced}
$$\gamma = (\gamma_0(\{0\}),\gamma') = ( \frac{a_1-1}{2},\frac{a_1-3}{2},...,\frac{1-a_1}{2},\gamma').$$
By the induction hypothesis, $\gamma'$ is triangular and $I_0(\widetilde{\mathbb{O}}') \subset U(\fg')$ is maximal. Thus, by Lemma \ref{lem:inductionsteptypeC}, $\gamma$ is triangular, and
\begin{align*}
\codim(\mathbb{O},\cN) &= \codim(\mathbb{O}',\cN_{G'}) + \dim(\cN_{\mathrm{GL}(a_t)}) & \text{(Proposition \ref{prop:propertiesofInd}(iii))}\\
&= \codim(\mathbb{O}_{\gamma'}, \cN_{L_{\gamma'}}) + \dim(\cN_{\mathrm{GL}(a_t)}) & \text{(Proposition \ref{prop:maximalitycriterion})}\\
&= \codim(\mathbb{O}_{\gamma}, \cN_{L_{\gamma}}) & \text{(Lemma \ref{lem:inductionsteptypeC})}.
\end{align*}
Hence, $I_0(\widetilde{\mathbb{O}})$ is maximal by a second application of Proposition \ref{prop:maximalitycriterion}. 
\end{proof}

\subsection{Auxilary ideals}

We now proceed to proving the maximality of certain auxilary ideals attached to birationally rigid covers. These results are used in the proof of Proposition \ref{prop:bimodsBCD}. 

Suppose $\widetilde{O}$ is a birationally rigid cover. Choose a Levi subgroup $L= \mathrm{GL}(a_1) \times ... \times \mathrm{GL}(a_s) \times G(m) \subset G$ and a birationally rigid $L$-orbit 
$$\mathbb{O}_L =  \{0\} \times ... \times \{0\} \times \mathbb{O}_{G(m)} \subset \cN_L$$
such that $\mathbb{O} = \mathrm{Bind}^G_L \mathbb{O}_L$. For each $k \in \{1,...,s\}$ let
$$\tau_1(k) = \frac{1}{2}(\underbrace{0,...,0}_{a_1+...+a_{k-1}},\underbrace{1,...,1}_k,\underbrace{0,...,0}_{a_{k+1}+...+a_s+m}) \in \fX(\fl)$$
and let
$$\delta(t) := \frac{t}{2}\tau_1(s) + \sum_{k=1}^{s-1} \frac{1}{2}\tau_1(k) \in \fX(\fl).$$
\begin{lemma}\label{lem:auxilaryidealstypeC}
For every $t \in (0,1)$, $I_{\delta(t)}(\mathbb{O}) \subset U(\fg)$ is a maximal ideal. 
\end{lemma}

\begin{proof}
Let $G' = G(n-a_1)$, let 
$$L' = \mathrm{GL}(a_2) \times ... \mathrm{GL}(a_{s}) \times G(m) \subset G', \qquad \mathbb{O}_{L'} = \underbrace{\{0\} \times ... \{0\}}_{s-1} \times \mathbb{O}_{G(m)},$$
and let $\mathbb{O}_{G'} = \mathrm{Bind}^{G'}_{L'} \mathbb{O}_{L'}$. Note that by Proposition \ref{prop:nocodim2leaves}, $\mathbb{O}_{G'}$ admits a 2-leafless $G'$-equivariant cover $\widetilde{\mathbb{O}}_{G'} \to \mathbb{O}_{G'}$. Let $\gamma' = \gamma_0(\widetilde{\mathbb{O}}_{G'})$ and $\gamma(t) = \gamma_{\delta(t)}(\mathbb{O}_G)$. Then
$$L^{\vee}_{\gamma(t)} = \mathrm{GL}(a_1) \times L^{\vee}_{\gamma'},   \qquad L^{\vee}_{\gamma(t),0} =  H_{a_1} \times M^{\vee}_{\gamma'},$$
where $H_{a_1} \subset \mathrm{GL}(a_1)$ is the maximal torus. Hence
$$\mathbb{O}_{\gamma(t)} = \mathsf{D}(\mathbb{O}^{\vee}_{\gamma(t)}) =  \{0\} \times \mathsf{D}(\mathbb{O}^{\vee}_{\gamma'})  =  \{0\} \times \mathbb{O}_{\gamma'}.$$
Using Proposition \ref{prop:maximalitytypeC}, we deduce
\begin{align*}
\codim(\mathbb{O}_{\gamma(t)}, \cN_{L_{\gamma(t)}}) &= \dim(\cN_{\mathrm{GL}(a_1)}) + \codim(\mathbb{O}_{\gamma'}, \cN_{L_{\gamma'}}) \\
&= \dim(\cN_{\mathrm{GL}(a_1)}) + \codim(\mathbb{O}_{G'}, \cN_{G'}) \\
&= \codim(\mathbb{O}, \cN).
\end{align*}
Hence, $I_{\delta(t)}(\mathbb{O})$ is maximal by Proposition \ref{prop:maximalitycriterion}.
\end{proof}

\section{Modifications for type $B$}\label{subsec:maximalityB}

Suppose $G = \mathrm{SO}(2n+1)$. The statements of Lemmas \ref{lem:maximalitybirigidtypeC}, \ref{lem:inductionsteptypeC}, \ref{lem:auxilaryidealstypeC}, and Proposition \ref{prop:maximalitytypeC} remain true. The proofs of Proposition \ref{prop:maximalitytypeC} and Lemma \ref{lem:auxilaryidealstypeC} hold word for word. We indicate below how the proofs of Lemmas \ref{lem:maximalitybirigidtypeC} and \ref{lem:inductionsteptypeC} should be modified.

\begin{lemma}\label{lem:maximalitybirigidtypeB}
Suppose $\widetilde{\mathbb{O}}$ is a birationally rigid $G$-equivariant nilpotent cover. Then $I_0(\widetilde{\mathbb{O}}) \subset U(\fg)$ is a maximal ideal.
\end{lemma}

\begin{proof}
Define partitions $p$,$q$,$x$, and $y$ as in the proof of Lemma \ref{lem:maximalitybirigidtypeC}. By Proposition \ref{prop:centralcharacterbirigidcover} and \cite[Thm 5.1]{McGovern1994}, it suffices to show
$$B(((g(y) \cup f_B(x))^t ) = p.$$
Our approach will be similar to Lemma \ref{lem:maximalitybirigidtypeC}, with a few modifications. 

We decompose $q$ into blocks
$$q = (\underbrace{q_1,...q_{l_1}}_{b^1}, \underbrace{q_{l_1+1}, ..., q_{l_1+l_2+1}}_{b^2},...,\underbrace{q_{2n-l_t+1},...,q_{2n}}_{b^s})$$
such that
\begin{itemize}
    \item For every $i$ ($1 \leq i \leq s$) and $j$ ($1 \leq j \leq l_i-1$) such that $l_1+...+l_{i-1}+j$ is odd, we have $b^i_j = b^i_{j+1}+1$.
    \item For every $i$ ($1 \leq i \leq s-1$), $b^i_{l_i} \geq b^{i+1}_1+2$.
\end{itemize}
(note, the definition of `block' in this case is different than in the proof of Lemma \ref{lem:maximalitybirigidtypeC}). For example, if $p=(6^4,5,3^3,2^2,1)$, then $q=(11,10,8,5^2,4)$, $b^1=(11,10,8)$, and $b^2=(5,5,4)$. For $2 \leq i \leq s-1$, put
$$\vardbtilde{b}^i = e(l(b^i)).$$
If $s=1$, put $\vardbtilde{b}^1=b^1$, else put
$$\vardbtilde{b}^1=l(b^1), \qquad \vardbtilde{b}^s=e(b^s).$$
Put
$$\vardbtilde{b}(q) = (\vardbtilde{b}^1,...,\vardbtilde{b}^s).$$
For example, if $q=(11,10,8,5^2,4)$ as above, then
$$b^1=(11,10,8), \ b^2 = (5,5,4), \qquad \vardbtilde{b}^1 = (11,10,7), \ \vardbtilde{b}^2 = (6,5,4), \qquad \vardbtilde{b} = (11,10,7,6,5,4).$$
Arguing by induction on $s$ (similarly to the proof of Lemma \ref{lem:maximalitybirigidtypeC}), we see that $g(y) \cup f_B(x) = \vardbtilde{b}(q)$ and $B(\vardbtilde{b}(q)^t)=p$. 
\end{proof}

Define $\gamma'$ and $\gamma$ as in Section \ref{subsec:maximalityC}.

\begin{lemma}\label{lem:inductionsteptypeB}
The $n$-tuple $\gamma$ is ($W$-conjugate to) a triangular one, and there is an equality
$$\codim(\mathbb{O}_{\gamma}, \cN_{L_{\gamma}}) = \codim(\mathbb{O}_{\gamma'},\cN_{L_{\gamma'}}) + \dim(\cN_{GL(a)}).$$
\end{lemma}

\begin{proof}
Repeat the proof of Lemma \ref{lem:inductionsteptypeC}, replacing all odd $\mathrm{SO}$ factors with $\mathrm{Sp}$ factors in $L^{\vee}_{\gamma}$, $L^{\vee}_{\gamma,0}$, $L^{\vee}_{\gamma'}$, and $L^{\vee}_{\gamma',0}$.
\end{proof}

\section{Modifications for type $D$}\label{subsec:maximalityD}

Suppose $G=\mathrm{SO}(2n)$. The statements of Lemma \ref{lem:maximalitybirigidtypeC}, Proposition \ref{prop:maximalitytypeC}, and Lemma \ref{lem:auxilaryidealstypeC} remain true. The proof of Lemma \ref{lem:auxilaryidealstypeC} holds word for word. The proofs of Lemma \ref{lem:maximalitybirigidtypeC} and Proposition \ref{prop:maximalitytypeC} must be altered slightly. Lemma \ref{lem:inductionsteptypeC} should be replaced by a slightly more elaborate argument. We indicate below how the statements and proofs should be modified.

\begin{lemma}\label{lem:maximalitybirigidtypeD}
Suppose $\widetilde{\mathbb{O}}$ is a birationally rigid $G$-equivariant cover. Then $I_0(\widetilde{\mathbb{O}}) \subset U(\fg)$ is a maximal ideal.
\end{lemma}

\begin{proof}
By Proposition \ref{prop:centralcharacterbirigidcover} and \cite[Thm. 5.1]{McGovern1994}, it suffices to show
$$B((g(y) \cup f_B(x))^t) = p.$$
Proceed exactly as in Lemma \ref{lem:maximalitybirigidtypeB}.
\end{proof}

Lemma \ref{lem:inductionsteptypeC} should be replaced by a pair of lemmas. Choose a Levi subgroup of the form
$$M = \mathrm{GL}(a) \times G(m) \subset G, \qquad a+m=n, \ m >0.$$
Suppose $\mu \neq (0)$, and let $\gamma' = v(\mu',\nu')$ be a triangular $m$-tuple. Let
$$\gamma = (\frac{a-1}{2}, \frac{a-3}{2}, ..., \frac{1-a}{2},\gamma').$$

\begin{lemma}\label{lem:inductionsteptypeD}
The $n$-tuple $\gamma$ is ($W$-conjugate to) a triangular one, and there is an equality
$$\codim(\mathbb{O}_{\gamma}, \cN_{L_{\gamma}}) = \codim(\mathbb{O}_{\gamma'},\cN_{L_{\gamma'}}) + \dim(\cN_{GL(a)}).$$
\end{lemma}

\begin{proof}
Proceed exactly as in the proof of Lemma \ref{lem:inductionsteptypeC}. 
\end{proof}

Next, choose a Levi subgroup of the form
$$M \simeq \mathrm{GL}(a_1) \times ... \times \mathrm{GL}(a_t) \subset G, \qquad a_1+...+a_t=n,$$
and let $\mathbb{O} = \mathrm{Bind}^G_M \{0\}$. Note that $\mathbb{O}$ corresponds to a very even partition and has a decoration determined by the conjugacy class of $M$. 

\begin{lemma}\label{lem:inductionsteptypeDextra}
$I_0(\mathbb{O}) \subset U(\fg)$ is a maximal ideal.
\end{lemma}

\begin{proof}
Let $\gamma:=\gamma_0(\mathbb{O})$. By Proposition \ref{prop:Ibetainduced}
$$\gamma_0(\mathbb{O}) = (\gamma_0(\{0\}), ..., \gamma_0(\{0\})) = (\rho(a_1),...,\rho(a_t)),$$
or possibly (if all $a_i$ are even)
$$\gamma_0(\mathbb{O}) = (\rho(a_1),...,\rho(a_t))',$$
where the prime indicates that the last entry is negated. We will show that
$$\codim(\mathbb{O},\cN) = \codim(\mathbb{O}_{\gamma},\cN_{L_{\gamma}}).$$
Write $e \subset a$ and $o  \subset a$ for the subpartitions consisting of even and odd parts. Define partitions $\mu = (\mu_1,...,\mu_r)$ and $\nu = (\nu_1,...,\nu_s)$ by
$$\mu_i = 2(o^t)_{2i-1}, \qquad \nu_j = 2(e^t)_{2j},$$
and put $R = |\mu| - \frac{1}{2}\mu_1$ and $S = |\nu|$.

Note that
$$\gamma = v(\mu,\nu) \text{ or } v(\mu,\nu)'$$
(see Definition \ref{def:triangular}). First suppose $\gamma = v(\mu,\nu)$. Then
$$L^{\vee}_{\gamma} = \mathrm{SO}(2R) \times \mathrm{SO}(2S), \qquad L^{\vee}_{\gamma,0} = \mathrm{SO}(\mu_1) \times \prod_{i=2}^r \mathrm{GL}(\mu_i) \times \prod_{j=1}^s \mathrm{GL}(\nu_j).$$
Hence $\mathbb{O}^{\vee}_{\gamma} = \mathbb{O}^{\vee,o}_{\gamma} \times \mathbb{O}^{\vee,e}_{\gamma}$ where $ \mathbb{O}^{\vee,o}_{\gamma}$ is the nilpotent $\mathrm{SO}(2R)$-orbit corresponding to the partition $(2o^t)^t$ and $\mathbb{O}^{\vee,e}_{\gamma}$ is the nilpotent $\mathrm{SO}(2S)$-orbit corresponding to the partition $(2e^t)^t$. Note that $(2o^t)^t$ is just $o$, but with multiplicities doubled. Similarly for $(2e^t)^t$. Hence, by Proposition \ref{prop:inclusionclassical} we have
$$\mathbb{O}^{\vee,o}_{\gamma} = \mathrm{Sat}^{\mathrm{SO}(2R)}_{\prod \mathrm{GL}(o_i)} \prod \mathbb{O}^{\mathrm{prin}}_{\mathrm{GL}(o_i)}, \qquad \mathbb{O}^{\vee,e}_{\gamma} = \mathrm{Sat}^{\mathrm{SO}(2S)}_{\prod \mathrm{GL}(e_j)} \prod \mathbb{O}^{\mathrm{prin}}_{\mathrm{GL}(e_j)}.$$
Using Proposition \ref{prop:inclusioninduction} we obtain
$$\mathbb{O}_{\gamma} = \mathsf{D}(\mathbb{O}^{\vee}_{\gamma}) = \mathsf{D}(\mathbb{O}^{\vee,o}_{\gamma}) \times \mathsf{D}(\mathbb{O}^{\vee,e}_{\gamma}) = \mathrm{Ind}^{\mathrm{SO}(2R)}_{\prod \mathrm{GL}(o_i)} \{0\} \times  \mathrm{Ind}^{\mathrm{SO}(2S)}_{\prod \mathrm{GL}(e_j)} \{0\}.$$
Since induction preserves codimension
\begin{align*}\codim(\mathbb{O}_{\gamma}, \cN_{L_{\gamma}}) &= \codim(\mathrm{Ind}^{\mathrm{SO}(2R)}_{\prod \mathrm{GL}(o_i)}  \{0\}, \cN_{\mathrm{SO}(2R)}) + \codim(\mathrm{Ind}^{\mathrm{SO}(2S)}_{\prod \mathrm{GL}(e_j)} \{0\}, \cN_{\mathrm{SO}(2S)}) \\
&= \sum_i \dim(\cN_{\mathrm{GL}(a_i)})\\
&= \codim(\mathbb{O}, \cN)
\end{align*}
as desired.

Next, suppose $\gamma = v'(\mu,\nu)$. In this case, we have $o = (0)$ and $a=(e)$. Compute
$$L^{\vee}_{\gamma} = \mathrm{SO}(2n), \qquad L^{\vee}_{\gamma,0} = \prod_{j=1}^s \mathrm{GL}(\nu_j)'.$$
Arguing as above (with $\mu = (0)$) we deduce
$$\codim(\mathbb{O}_{\gamma}, \cN_{L_{\gamma}}) = \codim(\mathbb{O}, \cN).$$
\end{proof}

\begin{cor}\label{cor:maximalitytypeD}
Suppose $\widetilde{\mathbb{O}}$ is a $G$-equivariant cover. Then $I_0(\widetilde{\mathbb{O}}) \subset U(\fg)$ is a maximal ideal.
\end{cor}

\begin{proof}
Choose a Levi subgroup $M \simeq  \mathrm{GL}(a_1) \times ... \times \mathrm{GL}(a_t) \times G(m)$ of $G$ and a birationally rigid $M$-equivariant cover
$$\widetilde{\mathbb{O}}_M =  \{0\} \times ... \{0\} \times \widetilde{\mathbb{O}}_{G(m)}$$
such that $\widetilde{\mathbb{O}} = \mathrm{Bind}^G_M \widetilde{\mathbb{O}}_M$. If $m=0$, use Lemma \ref{lem:inductionsteptypeDextra}. Otherwise, proceed by induction on $t$. The base case is Lemma \ref{lem:maximalitybirigidtypeD}. The induction step is Lemma \ref{lem:inductionsteptypeD}.
\end{proof}

\chapter{Singularities of type $m$}\label{SS_m_sing}

Let $\fg$ be a simple exceptional Lie algebra and let $\OO \subset \fg^*$ be a rigid nilpotent orbit. Premet shows in \cite{Premet2013} that in almost all cases, there is a unique multiplicity 1 primitive ideal $I \subset U(\fg)$ with $V(I) = \overline{\OO}$. In six cases, however, there are several such ideals. The orbits in question are: 
\begin{equation}\label{eq:sixorbits}\widetilde{A}_1 \subset G_2, \quad \widetilde{A}_2+A_1 \subset F_4, \quad (A_3+A_1)' \subset E_7, \quad A_3+A_1, A_5+A_1, D_5(a_1)+A_2 \subset E_8.\end{equation}
For each of these orbits, Premet constructs two multiplicity 1 primitive ideals with associated variety $\overline{\OO}$. From Premet's computations, it is not at all clear which of these ideals is unipotent (in fact, the situation is even murkier---in two cases, Premet does not claim that his list of ideals is exhaustive. In these cases, it is not even clear that $I_0(\OO)$ is among the two ideals he constructs). In this appendix, we will determine the unipotent ideals attached to the six orbits in (\ref{eq:sixorbits}), thus completing the determination of the unipotent spectrum for rigid nilpotent orbits.

There is a geometric feature which these orbits have in common: in each case, $\overline{\OO}$ contains a (unique) codimension 2 orbit $\OO' \subset \overline{\OO}$ and the transverse slice for $(\OO', \overline{\OO})$ is a singularity of type $m$, see the incidence tables in \cite{fuetal2015}. Recall, a type $m$ singularity is a a non-normal conical singularity with an $\mathrm{SL}(2)$-action admitting an open orbit isomorphic
to $\CC^2 - \{0\}$. See \cite[Sec 1.8.4]{fuetal2015} for details. 

Let $R'$ denote the reductive part of the centralizer
of an element $e' \in \OO'$. In the following table, we indicate the orbit $\OO'$ and the isomorphism type of $\mathfrak{r}'$. As usual, we use \cite[Sec 13.1]{Carter1993} for $\mathfrak{r}'$ and \cite[Tables]{fuetal2015} for $\OO'$.

\smallskip

\begin{center}
\begin{tabular}{|c|c|c|c|}\hline
$\fg$&$\OO$&$\OO'$&$R'$\\\hline
$G_2$&$\widetilde{A}_1$&$A_1$&$A_1$\\\hline
$F_4$&$\widetilde{A}_2+A_1$&$A_2+\widetilde{A}_1$&$A_1$\\\hline
$E_7$&$(A_3+A_1)'$&$2A_2+A_1$&$2A_1$\\\hline
$E_8$&$A_3+A_1$&$2A_2+A_1$&$G_2+A_1$\\\hline
$E_8$&$A_5+A_1$&$A_4+A_3$&$A_1$\\\hline
$E_8$&$D_5(a_1)+A_2$&$A_4+A_3$&$A_1$\\\hline
\end{tabular}
\end{center}
\smallskip
By inspection, we arrive at the following.

\begin{lemma}\label{lem:O_prime_properties}
The orbit $\OO'$ is rigid and of principal Levi type. The Lie algebra
$\mathfrak{r}'$ is semisimple.
\end{lemma}

Choose an $\mathfrak{sl}(2)$-triple $(e',h',f')$ and let $R'=Z_G(e',h',f')$. Let $S'$ be the Slodowy
slice $e'+\fg_{f'}$. The intersection $S'\cap \overline{\OO}$ is a
singularity of type $m$, by our choice of $(\OO,\OO')$. The following claim was
established in \cite[Sec 1.8.4]{fuetal2015}.

\begin{lemma}\label{lem:slice_action}
The action of $R'$ on $S'\cap \overline{\OO}'$ factors through an epimorphism
$R'\twoheadrightarrow \mathrm{SL}(2)$. The induced action of $\mathrm{SL}(2)$ on $S'\cap \overline{\OO}'$ coincides with the usual $\mathrm{SL}(2)$-action on a type $m$ singularity.
\end{lemma}

It follows that $R'=R_1'\times R_0'$, where $R_1'\simeq \mathrm{SL}(2)$ and $R_0'$
is an algebraic group with semisimple identity component.

\section{Algebras $\cA$ and $\cA_\dagger$}

Now let $I \subset U(\fg)$ be a multiplicity 1 primitive ideal with $V(I)=\overline{\OO}$. Set $\cA:=U(\fg)/I$. Combining Lemmas \ref{lem:adjunctionmorphismiso} and \ref{lem:upper_dag_1dim}, we equip $\cA$ with a good algebra filtration such that $\gr\cA\hookrightarrow \CC[\OO]$. Since $\OO$ is rigid, $\CC[\OO]$ admits a unique filtered quantization, see Corollary \ref{cor:criterionbirigid}. Thus, $\gr\cA\xrightarrow{\sim} \CC[\OO]$ if and only if $I=I_0(\OO)$. 

Let $\cW'$ denote the
W-algebra for the orbit $\OO'$ and let $\bullet_{\dagger}$ denote the
corresponding restriction functor, see Section \ref{subsec:W}. Consider the ideal $I_\dagger\subset \cW'$
and the quotient $\cA_\dagger:=\cW'/I_{\dagger}$.
By (ii) of Proposition \ref{prop:propsofdagger}, $\gr\cA_\dagger$ coincides with the pullback of $\gr\cA$ to $S'$.
Therefore $\gr\cA_\dagger\hookrightarrow
\CC[S'\cap \OO]=\CC[\CC^2]$. There is a Hamiltonian action of $R'$ on $\cW'$ and on $\cA_\dagger$.
Note that the action of $R_0'$  on $\cA_\dagger$ is trivial. So the quantum comoment map for the $R'$-action on $\cA_\dagger$ can be thought of as a homomorphism $U(\mathfrak{sl}_2)\rightarrow \cA_\dagger$.

The following lemma is the first step towards computing $I_0(\OO)$. 

\begin{lemma}\label{lem:sl_2_central_character}
The following claims are true:
\begin{enumerate}
\item
The kernel of $U(\mathfrak{sl}(2))\rightarrow \cA_\dagger$ has infinitesimal character in $\ZZ+\frac{1}{2}$.
\item This infinitesimal character is $\pm \frac{1}{2}$ (after $\rho$-shift) if and only if
$I=I_0(\OO)$.
\end{enumerate}
\end{lemma}
\begin{proof}
The algebra $\cA_{\dagger}$ is a Harish-Chandra bimodule for $U(\mathfrak{sl}(2))$ and its associated variety
is the image of $S'\cap \overline{\OO}$ in $\mathfrak{sl}(2)^*$. This image is the nilpotent cone. Moreover,
the kernel of $U(\mathfrak{sl}(2))\rightarrow \cA_\dagger$ is a completely prime ideal, hence primitive. Note that $\gr\cA$ embeds into $\CC[\CC^2]$ as an $\mathfrak{sl}(2)$-module with finite dimensional cokernel.
This is possible if and only if the infinitesimal character is in $\ZZ+\frac{1}{2}$.

If $I=I_0(\OO)$, then $\cA/I$ is the (unique) quantization of $\CC[\OO]$, and hence $\cA_\dagger$
quantizes $\CC[S'\cap \OO]=\CC[\CC^2]$. It follows that $\cA_{\dagger}$ is the Weyl algebra in two variables. The Hamiltonian $\mathrm{SL}(2)$-action on $\cA_\dagger$ lifts the standard action on
$\CC[\CC^2]$ and so is determined uniquely. We get the standard action of $\mathrm{SL}(2)$ on the Weyl algebra. The claim about its infinitesimal character is classical (see, e.g., (\ref{eq:barycenter}) for the case of $n=1$).

Now suppose $I$ is not unipotent. This implies that the inclusion $\gr\cA_\dagger\hookrightarrow \CC[\CC^2]$ is proper. Indeed, assume (for contradiction) that $\gr\cA_\dagger\xrightarrow{\sim} \CC[\CC^2]$. By (ii) of Proposition \ref{prop:propsofdagger}, this implies that the cokernel of $\gr \cA\hookrightarrow \CC[\OO]$ is supported away from $\Orb'$, i.e. has support of codimension at least $4$. The microlocalization of $U(\fg)/I$ to the complement of this support is a filtered quantization thereof. Similarly to \cite[Proposition 3.1]{Losev4}, the global sections of this microlocalization is a filtered quantization of $\CC[\OO]$. So $I$ is the kernel of the quantum comoment map to a filtered quantization of $\CC[\OO]$. Since $\OO$ is rigid (and hence birationally rigid) this quantization is unique, see Corollary \ref{cor:criterionbirigid}. So $I$ is the unipotent ideal, a contradiction. 

Since $\gr\cA_\dagger \hookrightarrow \CC[\CC^2]$ is proper, the infinitesimal character is different from $\pm \frac{1}{2}$.
\end{proof}


\section{Category $\mathsf{O}$ for $\cW'$}

Our method for computing the infinitesimal characters of unipotent ideals involves some additional machinery from \cite{BGK} and \cite{LosevStructureO}, which we will now briefly review.

Choose a minimal Levi subgroup $\underline{G} \subset G$ such that $e' \in \underline{\fg}$ (by Lemma \ref{lem:O_prime_properties}, $e'$ is then principal in $\underline{\fg}$). We can assume without loss that $\underline{G}$ is a standard Levi sugroup. Choose a dominant one-parameter subgroup $\nu:\CC^\times \rightarrow R'$ such that the centralizer in $G$ of $\nu(\CC^{\times})$ coincides with $\underline{G}$. Write $\fh$ for the Cartan subalgebra and $\fh^0$ for the center of $\underline{\fg}$. We can regard $\nu$ as an element of $\fh^0$ (by evaluating $d\nu(1)$). Recall, see Section \ref{subsec:W}, that $\mathfrak{r}'$ embeds into $\cW'$, so we can also view $\nu$ as an element of $\cW'$.

Following \cite[Sec 4.4]{BGK},\cite{LosevStructureO}, we will consider the category $\mathsf{O}_\nu(\cW')$ consisting
of finitely generated $\cW'$-modules $M$ such that
\begin{itemize}
\item $\nu\in \cW'$ acts locally finitely on $M$ with finite dimensional generalized
eigenspaces.
\item The eigenvalues are bounded from above, i.e. there are complex
numbers $z_1,\ldots,z_k$ (depending on $M$) such that any eigenvalue $z$ of $M$
satisfies $z_i-z\in \ZZ_{\geqslant 0}$ for some $i$.
\end{itemize}

We will now recall the classification of irreducible objects in $\mathsf{O}_\nu(\cW')$,
see \cite[Thm 4.5]{BGK}. The one-parameter subgroup $\nu$ gives rise to a $\ZZ$-grading on $\cW'$, $\cW'=\bigoplus_{i\in \ZZ}\cW'_i$. We will consider the subspace $\cW'_{>0} \subset \cW'$ and the
{\it Cartan subquotient}
$$\mathsf{C}_\nu(\cW'):=\cW'_0/\sum_{i>0}\cW'_{-i}\cW'_i.$$
For $M\in \mathsf{O}_\nu(\cW')$ consider the subspace of singular vectors,
$\operatorname{Ann}_{\cW'_{>0}}(M)$. This is a $\mathsf{C}_\nu(\cW')$-module.
The assignment $L\mapsto \operatorname{Ann}_{\cW'_{>0}}(L)$
defines a bijection between the set of simples in $\mathsf{O}_\nu(\cW')$
and the set of simple finite dimensional $\mathsf{C}_\nu(\cW')$-modules.

A description of $\mathsf{C}_\nu(\cW')$ was given in \cite[Section 4.1]{BGK}.
Let $\underline{\cW}'$ denote the W-algebra for $(\underline{\fg},e')$.
Since $e'$ is principal, the algebra $\underline{\cW}'$ is isomorphic
to the center of $U(\underline{\fg})$, i.e., to $\CC[\fh^*]^{\underline{W}}$,
where $\underline{W}$ is the Weyl group of $\underline{\fg}$. By
\cite[Theorem 4.3]{BGK}, see also \cite[Rmk 5.5]{LosevStructureO}, there is an isomorphism $\iota:\underline{\cW}'\xrightarrow{\sim}\mathsf{C}_\nu(\cW')$. We note that there are natural
embeddings $\fh^0\hookrightarrow \underline{\cW}',\mathsf{C}_\nu(\cW')$, the latter
induced from the embedding $\fh^0\hookrightarrow \cW'_0$. The next lemma describes the restriction of $\iota$ to the subspace $\fh^0 \subset \underline{\cW}'$. 

For each $n \in \ZZ$, define $\Delta^+(n) := \{\alpha \in \Delta^+ \mid \alpha(h') = n\}$ and consider the element
\begin{equation}\label{eq:defofrhoe}\rho_{e'} := \frac{1}{2}(\sum_{\alpha \in \Delta^+(0)}\alpha + \sum_{\alpha\in \Delta^+(1)}\alpha) \in \fh^*\end{equation}
%

The next lemma follows from \cite[Theorem 4.3]{BGK}, see also \cite[Sec 2.6]{Premet2013}.

\begin{lemma}\label{lem:shift_formula}
$\iota(x)=x+\langle \rho-\rho_{e'},x\rangle$ for every $x \in \fh^0$.
\end{lemma}

The isomorphism $\iota: \underline{\cW}'\xrightarrow{\sim} \mathsf{C}_\nu(\cW')$ shows
that the space of singular vectors in any simple module $L$ is one-dimensional
and the action is given by a character of $\underline{\cW}'$, i.e. an element in
$\fh^*/\underline{W}$. We will represent this element by an
{\it integrally anti-dominant} (for $\underline{\fg}$) weight $\lambda\in \fh^*$. This means
$$\langle \lambda, \alpha^{\vee} \rangle \notin \ZZ_{>0}, \qquad \forall \alpha \in \Delta^+(\underline{\fg},\fh).$$

Write $L_\lambda$ for the irreducible module in $\mathsf{O}_\nu(\cW')$ corresponding
to $\lambda$. For what follows we will need a formula for the action of $\fh^0\hookrightarrow \mathsf{C}_\nu(\cW')$
on $L_\lambda$. The next Corollary is a direct consequence of Lemma \ref{lem:shift_formula}.

\begin{cor}\label{Cor:highest_weight}
The subspace $\fh^0 \hookrightarrow \mathsf{C}_\nu(\cW')$ acts on $L_{\lambda}$ by the character $(\lambda-\rho_{e'})|_{\fh^0}$. This character will be called the \emph{highest weight} of $L_{\lambda}$.
\end{cor}

Next, consider the category $\mathsf{Wh}$ of \emph{generalized Whittaker modules} for $\fg$, see \cite[Sec 4]{LosevStructureO}. These are $U(\fg)$-modules $M$ such that
\begin{itemize}
\item $\nu$ acts locally finitely on $M$ {and the eigenvalues are bounded from above}.
\item Each generalized eigenspace for the $\nu$-action is a Whittaker module for
$\underline{\fg}$ in the sense of Kostant, \cite{KostantWhittaker}.
\end{itemize}
We note that the irreducible objects in $\mathsf{Wh}$ are parameterized by their
highest weight spaces for $\nu$. Each highest weight space is an irreducible
Whittaker module for $\underline{\fg}$, and they are in bijection with $\fh^*/\underline{W}$
(by taking infinitesimal character). If $\lambda \in \fh^*$ is integrally anti-dominant, we write
$\widetilde{L}_\lambda$ for the irreducible object in $\mathsf{Wh}$ corresponding
to $\lambda$. The following is a direct corollary of \cite[Theorem 4.1]{LosevStructureO}.

\begin{prop}\label{Prop:category_equivalence}
There is a category equivalence $\Psi: \mathsf{O}_\nu(\cW')\xrightarrow{\sim}\mathsf{Wh}$ such that $\Psi(L_\lambda)\simeq \widetilde{L}_\lambda$
for all $\lambda\in \fh^*/\underline{W}$.
\end{prop}

We will need a result concerning the behavior of annihilator of $L_\lambda$
in a special case. For $\lambda\in \fh^*$ we write $I(\lambda)$ for the annihilator of the irreducible
$U(\fg)$-module with highest weight $\lambda-\rho$.


There is a map $\bullet^{\ddag}$ from the set of two-sided ideals in $\cW'$ to the set of
two-sided ideals in $U(\fg)$, see \cite[Theorem 1.2.2]{Losev3}, extending the map constructed
in Section \ref{subsec:W}.

\begin{prop}\label{Prop:annihilators}
Let $\lambda$ be 
integrally anti-dominant for $\underline{\fg}$.
Then $\operatorname{Ann}_{\cW'}(L_\lambda)^{\ddag}=I(\lambda)$.
\end{prop}
\begin{proof}
 By \cite[Thm 4.1]{LosevStructureO}, $\operatorname{Ann}_{U(\fg)}(\Psi(L_\lambda))=
\operatorname{Ann}_{\cW'}(L_\lambda)^\ddag$. By Proposition \ref{Prop:category_equivalence},
$\Psi(L_\lambda)=\widetilde{L}_\lambda$. By \cite[Thm 5.1.1]{Losev_parabolic}
(and its proof),  $\operatorname{Ann}_{U(\fg)}(\widetilde{L}_\lambda)=I(\lambda)$.
\end{proof}

\section{Category $\mathsf{O}$ for $\cA_\dagger$}
Consider the full subcategory $\mathsf{O}_\nu(\cA_\dagger)$ of $\mathsf{O}_\nu(\cW')$ consisting
of all modules annihilated by $I_\dagger$. We can form the Cartan subquotient $\mathsf{C}_\nu(\cA_\dagger)$ of $\cA_\dagger$
similarly to $\mathsf{C}_\nu(\cW')$ in the previous section.  Taking the subspace of singular vectors
defines a bijection between the set of isomorphism classes of irreducible modules in $\mathsf{O}_\nu(\cA_\dagger)$ and the set of isomorphism classes of irreducible $\mathsf{C}_\nu(\cA_\dagger)$-modules.

The following is our main result on the structure of irreducible modules in $\mathsf{O}_\nu(\cA_\dagger)$.
Note that the choice of $\nu$ gives rise to a system of simple roots in the Lie algebra $\mathfrak{r}'$.

\begin{prop}\label{Prop:Cat_O_A_dagger}
The following are true:
\begin{enumerate}
\item There is a single irreducible module in $\mathsf{O}_\nu(\cA_\dagger)$.
\item  Its highest weight, see Corollary \ref{Cor:highest_weight}, is of the form
$z\omega_1$, where $z\in \frac{1}{2}+\ZZ$ and $\omega_1$ is the fundamental weight of the factor
$\mathfrak{r}_1'\simeq \mathfrak{sl}(2)$ in $\mathfrak{r}'$ (see the discussion after Lemma \ref{lem:slice_action}).
\item We have $z=-\frac{1}{2}$ if $I$ is a unipotent ideal and $z\in \frac{1}{2}+\ZZ_{\geqslant 0}$
otherwise.
\end{enumerate}
\end{prop}
\begin{proof}
Recall that there is the quantum comoment map $U(\mathfrak{r}')\rightarrow \cA_\dagger$.
Its restriction to the semisimple factor $\mathfrak{r}'_0$ is zero, while the kernel of
$U(\mathfrak{r}_1')\rightarrow \cA_\dagger$ has infinitesimal character in $\frac{1}{2}+\ZZ$
by Lemma \ref{lem:sl_2_central_character}. It follows that a central reduction
of $U(\mathfrak{sl}(2))$ has a filtered embedding into $\cA_\dagger$. The image has
a complimentary nonzero Harish-Chandra bimodule summand. Let $h$ denote the image of
the standard basis element of $\mathfrak{sl}(2)$ in $\cA_\dagger$. It defines a grading such that
the sum of even graded components is the central reduction of $U(\mathfrak{sl}(2))$, while there
are also odd components. This grading is proportional to the grading defined by $\nu$.

We claim that $\mathsf{C}_\nu(\cA_\dagger)$ is generated by the image of $h$. Indeed, the degree $0$ component of $\cA_\dagger$ surjects onto $\mathsf{C}_\nu(\cA_\dagger)$, so it is enough to show that it coincides with $\CC[h]$. The proof of that reduces to the analogous claim for $\gr\cA_\dagger$, which in turn, reduces to $\CC[\CC^2]$, thanks to the embedding  $\gr\cA_\dagger\hookrightarrow\CC[\CC^2]$.
The claim that the degree $0$ part of $\CC[\CC^2]$ is generated by $h$ is clear. 

It follows that
an irreducible module in $\mathsf{O}_\nu(\cA_\dagger)$ is determined by its highest weight up to an isomorphism. The highest weight
must have the form $z\omega_1$ because $\mathfrak{r}'_0$ acts by $0$. By Lemma \ref{lem:sl_2_central_character}, $z\in \frac{1}{2}+\ZZ$. This proves (2).

Now we prove (1). Take an irreducible object $L\in \mathsf{O}_\nu(\cA_\dagger)$ and restrict
it to $U(\mathfrak{sl}(2))$. This restriction lies in the category $\mathsf{O}$ for $\mathfrak{sl}(2)$. Its
infinitesimal character is half-integral, so it is semisimple. We claim that both irreducible
objects in the block occur in the restriction. Indeed, if only one simple occurs,
then the elements of odd degree in $\cA_\dagger$ act by $0$ on $L$. It follows that
$L$ is annihilated by a proper two-sided ideal in $\cA_\dagger$. Since $\gr\cA_\dagger\hookrightarrow
\CC[\CC^2]$ with finite dimensional cokernel, we see that $\operatorname{Spec}(\gr\cA_\dagger)$
has two symplectic leaves, $\{0\}$ and the open leaf. It follows that any proper two-sided
ideal in $\cA_\dagger$ must have finite codimension. This is impossible because its central
character for $\mathfrak{sl}(2)$ is in $\frac{1}{2}+\ZZ$. It follows that the highest weight
of $L$ is the larger of the highest weights of the two irreducible $\mathfrak{sl}(2)$-modules.
This determines the highest weight of $L$ uniquely. Since the highest weight of $L$
determines it uniquely, (1) is proved.

Now we prove (3). By Lemma \ref{lem:sl_2_central_character}, $I$ is unipotent if and only if
the infinitesimal character of the $U(\mathfrak{sl}(2))$-module $L$ is $\pm \frac{1}{2}$. The description of
the highest weight of $L$ in the previous paragraph implies that it is $-\frac{1}{2}$.
The case when $I$ is not unipotent is handled similarly.
\end{proof}

The following technical corollary is our main computational tool for determining $I_0(\OO)$ for $\OO$ in (\ref{eq:sixorbits}). A weight $\lambda\in \fh^*$ is {\it integrally dominant} (for $\fg$) if
$$\langle\lambda,\alpha^\vee\rangle\not\in \ZZ_{\leqslant 0}, \qquad \forall \alpha \in \Delta^+(\fg,\fh).$$ 
It is a classical fact that $\lambda$ is integrally dominant if and only if $I(\lambda)$ is maximal, see \cite[Thm 7.6.24]{Dixmier}. Recall the weight $\rho_{e'} \in \fh^*$ defined in (\ref{eq:defofrhoe}).

\begin{cor}\label{Cor:how_to_compute}
Let $\lambda\in \fh^*$ be integrally dominant for $\fg$ and integrally
anti-dominant for $\underline{\fg}$. Suppose that $I=I(\lambda)$ is a multiplicity 1 primitive ideal such that $V(I) = \overline{\OO}$. 
Then the following are true:
\begin{enumerate}
\item The kernel of the operator $(\lambda-\rho_{e'})|_{\fh^0}$ coincides with the Cartan subalgebra
of $\mathfrak{r}'_0$.
\item Choose $h_1 \in \fh^0$ in the nonnegative span of the positive coroots and orthogonal to the subspace $\ker{(\lambda-\rho_{e'})|_{\fh^0}}$. If
$$\langle \lambda - \rho_{e'}, h_1 \rangle <0,$$
then $I(\lambda) = I_0(\OO)$.
\end{enumerate}
\end{cor}
\begin{proof}
We claim that the unique irreducible in  $\mathsf{O}_\nu(\cA_{\dagger})$
is $L_{\lambda}$. Indeed, $\operatorname{Ann}_{\cW'}(L_{\lambda})^\ddag=I$
by Proposition \ref{Prop:annihilators}. By \cite[Theorem 1.2.2]{Losev3}, for any  ideal
$J\subset \cW'$, we have $(J^\ddag)_\dagger\subset J$. Applying this to $J:=\operatorname{Ann}_{\cW'}(L_{\lambda})$ we see that $L_{\lambda}$ is annihilated by
$I_\dagger$, which implies our claim.

By Corollary \ref{Cor:highest_weight}, the highest weight of $L_{\lambda}$ is
$(\lambda - \rho_{e'})|_{\fh^0}$. By (2) of Proposition \ref{Prop:Cat_O_A_dagger},
this highest weight vanishes on $\fh^0\cap \mathfrak{r}'_0$, which is a codimension $1$
subspace in $\fh^0$, and is nonzero. Hence $\fh^0\cap \mathfrak{r}'_0$ is the kernel. This proves (1).

To prove (2) we observe that $h_1$ must be a positive multiple of the standard basis element
$h\in \mathfrak{r}_1'\simeq \mathfrak{sl}(2)$. By Proposition \ref{Prop:Cat_O_A_dagger}, the unipotent
ideal is uniquely characterized by the property that the corresponding highest
weight $(\lambda-\rho_{e'})|_{\fh^0}$ is negative on $h$. (2) follows.
\end{proof}

\section{Computations}

Let $\OO$ be one the six orbits in (\ref{eq:sixorbits}), and let $I(\lambda)$ be one of the (two) multiplicity 1 primitive ideals constructed by Premet (note that Premet's highest weights do not include the $\rho$-shift, unlike ours). Assume that $\lambda$ is integrally dominant for $\fg$. We will use the following algorithm to decide whether $I(\lambda)=I(\OO)$.

\begin{enumerate}
    \item Fix a standard Levi subgroup $\underline{G} \subset G$ in which $\OO'$ is principal. Then $\fh^0$ is spanned by the fundamental coweights corresponding to the simple roots which are not contained in $\underline{G}$.
    \item Replace $\lambda$ with a $W$-conjugate which is both integrally dominant for $\fg$ and integrally anti-dominant for $\underline{\fg}$.
    \item Copy the character $\rho_{e'} \in \fh^*$ from Premet, and compute the difference $\lambda - \rho_{e'}$. 
    \item If $\dim(\mathfrak{h}^0)=1$, then there is a unique simple root $\alpha_i$ not contained in $\underline{G}$. If 
    $$\langle \lambda - \rho_{e'}, \varpi_i^{\vee}\rangle < 0,$$
    then $I(\lambda) = I(\OO)$.
    \item If $\dim(\mathfrak{h}^0) > 1$ (this happens in exacly two cases), compute the kernel of the operator  $(\lambda-\rho_{e'})|_{\fh^0}$. In both cases, this kernel is the span of a collection of fundamental coweights $\varpi_j^{\vee} \in \fh^0$. 
    \item In the setting of (5), let $\varpi_i^{\vee}$ denote the unique fundamental coweight in $\fh^0$ not contained in $\ker{(\lambda-\rho_{e'})|_{\fh^0}}$. Note that $\varpi_i^{\vee}$ differs from $h_1$ in Corollary \ref{Cor:how_to_compute} by an element of $\ker{(\lambda-\rho_{e'})|_{\fh^0}}$. If
    $$\langle \lambda-\rho_{e'}, \varpi_i^{\vee}\rangle <0,$$
    then $I(\lambda)=I_0(\OO)$. 
\end{enumerate}

For the computations below, we will use the standard (Bourbaki) numbering of simple roots in $\Delta^+$, see the Dynkin diagrams in Section \ref{subsec:dualdistinguished}. This is also numbering used by Premet in \cite{Premet2013}. As usual, we will write elements of $\fh^*$ (resp. $\fh$) in the basis of fundamental weights (resp. simple coroots). One advantage of this convention is that  the pairing $\langle\lambda,x \rangle$ is the dot product. Very often, we will need to convert from the bases of simple roots (resp. fundamental co-weights). For this, we use \cite[Planches I-IX]{Bourbaki46} as a reference (usually without comment).

\subsubsection{$\widetilde{A}_1 \subset G_2$}

Premet's ideals correspond to the weights
$$\lambda_1 = \frac{1}{2}(1,1), \qquad \lambda_2 = \frac{1}{2}(5,-1).$$
We claim that $I(\lambda_1) = I_0(\OO)$. 

Note that $\OO'=A_1$ and $\mathfrak{r}'=A_1$. $\underline{G}$ is the standard Levi subgroup corresponding to the short simple root $\alpha_2$, and therefore $\fh^0$ is spanned by the fundamental co-weight $\varpi_1^{\vee}$. Note that $\lambda_1$ is (integrally) dominant for $G$ and integrally anti-dominant for $\underline{\fg}$, so Corollary \ref{Cor:how_to_compute} is applicable. In our (fundamental weight) coordinates, $\rho_{e'} = \frac{1}{2}(3,0)$, see \cite[Sec 5.6]{Premet2013}. Thus $\lambda_1 - \rho_{e'} = \frac{1}{2}(-2,1)$. Note that
$$\langle \lambda_1-\rho_{e'},\varpi_1^{\vee}\rangle = \langle \frac{1}{2}(-2,1), (2,3)\rangle = -\frac{1}{2}$$
Hence, $I(\lambda_1)=I_0(\OO)$ by Corollary \ref{Cor:how_to_compute}. Since $\lambda_1$ is integrally dominant, $I_0(\OO)$ is thus the maximal ideal of infinitesimal character $\frac{1}{2}(1,1)$.

\subsubsection{$\widetilde{A}_2+A_1 \subset F_4$}

Premet's ideals correspond to the weights
$$\lambda_1 = \frac{1}{3}(1,1,1,1), \qquad \lambda_2 = \frac{1}{3}(-1,-1,5,2).$$
We claim that $I(\lambda_1) = I_0(\OO)$. 

Note that $\OO'=A_2+\widetilde{A}_1$ and $\mathfrak{r}'=A_1$. $\underline{G}$ is the standard Levi subgroup corresponding to the simple roots $\{\alpha_1,\alpha_2,\alpha_4\}$, and therefore $\fh^0$ is spanned by the fundamental co-weight $\varpi_3^{\vee}$. Note that $\lambda_1$ is (integrally) dominant for $G$ and integrally anti-dominant for $\underline{\fg}$, so Corollary \ref{Cor:how_to_compute} is applicable. In our (fundamental weight) coordinates, $\rho_{e'} = \frac{1}{2}(-1,2,0,1)$, see \cite[Sec 5.4]{Premet2013}. Thus $\lambda_1 - \rho_{e'} = \frac{1}{6}(5,-4,2,-1)$. Note that
$$\langle \lambda_1-\rho_{e'},\varpi_3^{\vee}\rangle = \langle \frac{1}{6}(5,-4,2,-1), (4,8,6,3)\rangle = -\frac{1}{2}$$
Hence, $I(\lambda_1)=I_0(\OO)$ by Corollary \ref{Cor:how_to_compute}. Since $\lambda_1$ is integrally dominant, $I_0(\OO)$ is thus the maximal ideal of infinitesimal character $\frac{1}{3}(1,1,1,1)$.

\subsubsection{$(A_3+A_1)' \subset E_7$}

Premet's ideals correspond to the weights
$$\lambda_1 = \frac{1}{2}(1,1,0,1,0,1,1), \qquad \lambda_2 = \frac{1}{2}(3,4,-4,-1,3,1,1).$$
We claim that $I(\lambda_1) = I_0(\OO)$. 

Note that $\OO'=2A_2+A_1$ and $\mathfrak{r}'=2A_1$. $\underline{G}$ is the standard Levi subgroup corresponding to the simple roots $\{\alpha_1,\alpha_2,\alpha_3,\alpha_5,\alpha_6\}$, and therefore $\fh^0$ is spanned by the fundamental co-weights $\varpi_4^{\vee},\varpi_7^{\vee}$. Note that $\lambda_1$ is (integrally) dominant for $G$ and integrally anti-dominant for $\underline{\fg}$, so Corollary \ref{Cor:how_to_compute} is applicable. In our (fundamental weight) coordinates, $\rho_{e'} = \frac{1}{2}(-1,2,3,-1,0,2,0)$, see \cite[Sec 4.10]{Premet2013}. Thus $\lambda_1 - \rho_{e'} = \frac{1}{2}(2,-1,-3,2,0,-1,1)$. Note that
\begin{align*}
\langle \lambda_1-\rho_{e'},\varpi_4^{\vee}\rangle &= \langle\frac{1}{2}(2,-1,-3,2,0,-1,1), (4,6,8,12,9,6,3)\rangle = -\frac{1}{2}\\
\langle \lambda_1-\rho_{e'},\varpi_7^{\vee}\rangle &= \langle \frac{1}{2}(2,-1,-3,2,0,-1,1), \frac{1}{2}(2,3,4,6,5,4,3)\rangle = 0
\end{align*}
Hence, $I(\lambda_1)=I_0(\OO)$ by Corollary \ref{Cor:how_to_compute}. Since $\lambda_1$ is integrally dominant, $I_0(\OO)$ is thus the maximal ideal of infinitesimal character $\frac{1}{2}(1,1,0,1,0,1,1)$.

\subsubsection{$A_3+A_1 \subset E_8$}

Premet's ideals correspond to the weights
$$\lambda_1 = \frac{1}{2}(1,1,0,1,0,1,1,2), \qquad \lambda_2 = \frac{1}{2}(-3,3,-2,1,2,1,1,2).$$
We claim that $I(\lambda_1) = I_0(\OO)$. 

Note that $\OO'=2A_2+A_1$ and $\mathfrak{r}'=A_1+G_2$. So $\underline{G}$ is the standard Levi subgroup corresponding to the simple roots $\{\alpha_1,\alpha_2,\alpha_3,\alpha_5,\alpha_6\}$, and therefore $\fh^0$ is spanned by the fundamental co-weights $\varpi_4^{\vee},\varpi_7^{\vee},\varpi_8^{\vee}$. Note that $\lambda_1$ is (integrally) dominant for $G$ and integrally anti-dominant for $\underline{\fg}$, so Corollary \ref{Cor:how_to_compute} is applicable. In our (fundamental weight) coordinates, $\rho_{e'} =  \frac{1}{2}(0,2,4,-2,-1,4,-1,2)$, see \cite[Sec 3.8]{Premet2013}. Thus $\lambda_1 - \rho_{e'} = \frac{1}{2}(1,-1,-4,3,1,-3,2,0)$. Note that
\begin{align*}
\langle \lambda_1-\rho_{e'},\varpi_4^{\vee}\rangle &= \langle \frac{1}{2}(1,-1,-4,3,1,-3,2,0), (10,15,20,30,24,18,12,6) \rangle = -\frac{1}{2}\\
\langle \lambda_1-\rho_{e'},\varpi_7^{\vee}\rangle &= \langle \frac{1}{2}(1,-1,-4,3,1,-3,2,0), (4,6,8,12,10,8,6,3)\rangle = 0\\
\langle \lambda_1-\rho_{e'},\varpi_8^{\vee}\rangle &= \langle \frac{1}{2}(1,-1,-4,3,1,-3,2,0), (2,3,4,6,5,4,3,2)\rangle = 0
\end{align*}
Hence, $I(\lambda_1)=I_0(\OO)$ by Corollary \ref{Cor:how_to_compute}. Since $\lambda_1$ is integrally dominant, $I_0(\OO)$ is thus the maximal ideal of infinitesimal character $\frac{1}{2}(1,1,0,1,0,1,1,2)$.

\subsubsection{$A_5+A_1 \subset E_8$}

Premet's ideals correspond to the weights
$$\lambda_1 = \frac{1}{6}(2,2,1,1,1,1,1,1), \qquad \lambda_2 = \frac{1}{6}(2,2,1,7,-11,7,1,1).$$
We claim that $I(\lambda_1-\rho) = I_0(\OO)$. 

Note that $\OO'=A_4+A_3$ and $\mathfrak{r}'=A_1$. So $\underline{G}$ is the standard Levi subgroup corresponding to the simple roots $\{\alpha_1,\alpha_2,\alpha_3,\alpha_4,\alpha_6,\alpha_7,\alpha_8\}$, and therefore $\fh^0$ is spanned by the fundamental co-weight $\varpi_5^{\vee}$. Although $\lambda_1$ is clearly (integrally) dominant for $G$, it is not integrally anti-dominant for $\underline{\fg}$. A bit of tinkering shows that $\lambda_1$ is $W$-conjugate to the weight $\lambda_1' = \frac{1}{6}(2,2,1,2,-1,2,1,1)$ (the simple reflection $s_{\alpha_5}$ takes $\lambda_1$ to $\lambda_1'$). It is straightforward to check that $\lambda_1'$ is integrally dominant for $\fg$ and integrally anti-dominant for $\underline{\fg}$. 

In our (fundamental weight) coordinates, $\rho_{e'} = \frac{1}{2}(0,2,0,0,0,1,1,-1)$, see \cite[Sec 3.15]{Premet2013}. Thus $\lambda_1' - \rho_{e'} = \frac{1}{6}(2,-4,1,2,-1,-1,-2,4)$. Note that
$$\langle \lambda_1'-\rho_{e'}, \varpi_5^{\vee}\rangle =  \langle \frac{1}{6}(2,-4,1,2,-1,-1,-2,4), (8,12,16,24,20,15,10,5)\rangle = -\frac{1}{2}.$$
Hence $I(\lambda_1-\rho)=I(\lambda_1'-\rho)=I_0(\OO)$ by Corollary \ref{Cor:how_to_compute}. Since $\lambda_1$ is integrally dominant, $I_0(\OO)$ is thus the maximal ideal of infinitesimal character $\frac{1}{6}(2,2,1,1,1,1,1,1)$.

\subsubsection{$D_5(a_1)+A_2 \subset E_8$}

Premet's ideals correspond to the weights
$$\lambda_1 = \frac{1}{4}(-1,-1,-1,4,-1,4,-1,-1), \qquad \lambda_2 = \frac{1}{4}(-1,-1,-1,8,-9,8,-1,-1).$$
We claim that $I(\lambda_1) = I_0(\OO)$. 

Note that $\OO'=A_4+A_3$ and $\mathfrak{r}'=A_1$. So $\underline{G}$ is the standard Levi subgroup corresponding to the simple roots $\{\alpha_1,\alpha_2,\alpha_3,\alpha_4,\alpha_6,\alpha_7,\alpha_8\}$, and therefore $\fh^0$ is spanned by the fundamental co-weight $\varpi_5^{\vee}$. Although $\lambda_1$ is integrally dominant for $G$, it is not integrally anti-dominant for $\underline{\fg}$. Using {\tt atlas}, we see that $\lambda_1$ is $W$-conjugate to the dominant weight $\lambda_1' = \frac{1}{4}(1,1,1,0,1,1,1,1)$ (the Weyl group element $s_{\alpha_1}s_{\alpha_2}s_{\alpha_3}s_{\alpha_1}s_{\alpha_5}s_{\alpha_7}s_{\alpha_8}s_{\alpha_7}$ takes $\lambda_1$ to $\lambda_1'$). It is straightforward to check that $\lambda_1'$ is integrally anti-dominant for $\underline{\fg}$. 

Again, $\rho_{e'} = \frac{1}{2}(0,2,0,0,0,1,1,-1)$. Thus $\lambda_1' - \rho_{e'} = \frac{1}{4}(1,-3,1,0,1,-1,-1,3)$. Note that
$$\langle \lambda_1'-\rho_{e'}, \varpi_5^{\vee}\rangle =  \langle \frac{1}{4}(1,-3,1,0,1,-1,-1,3),  (8,12,16,24,20,15,10,5)\rangle = -\frac{1}{2}$$
Hence $I(\lambda_1)=I(\lambda_1'-\rho)=I_0(\OO)$ by Corollary \ref{Cor:how_to_compute}. Thus, $I_0(\OO)$ is the maximal ideal of infinitesimal character $\frac{1}{4}(1,1,1,0,1,1,1,1)$.

\chapter{Tables}

\begin{table}[hbt!]
  \small
  \caption{Unipotent infinitesimal characters attached to rigid orbits in exceptional types. Characters are denoted in fundamental weight coordinates. Special unipotent infinitesimal characters are highlighted in blue. See Section \ref{subsec:centralcharexceptional} for further explanation.}    \label{table:exceptional}
    \begin{tabular}{|c|c|c|c|} \hline
       $\fg$ & $\OO$ & $\gamma_0(\OO)$ \\ \hline 
       
     $G_2$ & $\{0\}$ & $\cellcolor{blue!20} (1,1)$ \\ \hline
     $G_2$ & $A_1$  & $\frac{1}{3}(3,1)$  \\ \hline
     $G_2$ & $\widetilde{A}_1$ & $\frac{1}{2}(1,1)$ \\ \hline

    $F_4$ & $\{0\}$ & \cellcolor{blue!20}$(1,1,1,1)$\\ \hline
    $F_4$ & $A_1$  & $\frac{1}{2}(1,1,2,2)$  \\ \hline
    $F_4$ & $\widetilde{A}_1$ & \cellcolor{blue!20}$(1,0,1,1)$  \\ \hline
    $F_4$ & $A_1+\widetilde{A}_1$ & \cellcolor{blue!20}$(1,0,1,0)$  \\ \hline
    $F_4$ & $A_2+\widetilde{A}_1$ & $\frac{1}{4}(1,1,2,2)$  \\ \hline
    $F_4$ & $\widetilde{A}_2+A_1$ & $\frac{1}{3}(1,1,1,1)$  \\ \hline

    $E_6$ & $\{0\}$ & \cellcolor{blue!20}$(1,1,1,1,1,1)$ \\ \hline
    $E_6$ & $A_1$ & \cellcolor{blue!20}$(1,1,1,0,1,1)$  \\ \hline
    $E_6$ & $3A_1$ & $\frac{1}{2}(1,1,1,1,1,1)$  \\ \hline
    $E_6$ & $2A_2+A_1$ & $\frac{1}{3}(1,1,1,1,1,1)$ \\ \hline

    $E_7$ & $\{0\}$ & \cellcolor{blue!20}$(1,1,1,1,1,1,1)$  \\ \hline
    $E_7$ & $A_1$ &  \cellcolor{blue!20}$(1,1,1,0,1,1,1)$\\ \hline
    $E_7$ & $2A_1$ & \cellcolor{blue!20}$(1,1,1,0,1,0,1)$  \\ \hline
    $E_7$ & $(3A_1)'$ & $\frac{1}{2}(1,1,1,1,1,1,2)$ \\ \hline
    $E_7$ & $4A_1$ & $\frac{1}{2}(1,1,1,1,1,1,1)$ \\ \hline
    $E_7$ & $A_2+2A_1$ & \cellcolor{blue!20}$(1,0,0,1,0,0,1)$  \\ \hline
    $E_7$ & $2A_2+A_1$ & $\frac{1}{3}(1,1,1,1,1,1,1)$  \\ \hline
    $E_7$ & $(A_3+A_1)'$ & $\frac{1}{2}(1,1,0,1,0,1,1)$  \\ \hline

    $E_8$ & $\{0\}$ & \cellcolor{blue!20}$(1,1,1,1,1,1,1,1)$ \\ \hline
    $E_8$ & $A_1$ & \cellcolor{blue!20}$(1,1,1,0,1,1,1,1)$  \\ \hline
    $E_8$ & $2A_1$ & \cellcolor{blue!20}$(1,1,1,0,1,0,1,1)$ \\ \hline
    $E_8$ & $3A_1$ & $\frac{1}{2}(1,1,1,1,1,1,2,2)$ \\ \hline
    $E_8$ & $4A_1$ & $\frac{1}{2}(1,1,1,1,1,1,1,1)$  \\ \hline
    $E_8$ & $A_2+A_1$ & \cellcolor{blue!20}$(1,0,0,1,0,1,0,1)$  \\ \hline
    $E_8$ & $A_2+2A_1$ & \cellcolor{blue!20}$(1,0,0,1,0,0,1,1)$ \\ \hline
    $E_8$ & $A_2+3A_1$ & $\frac{1}{2}(1,1,1,0,1,1,1,1)$  \\ \hline
    $E_8$ & $2A_2+A_1$ & $\frac{1}{3}(1,1,1,1,1,1,1,3)$ \\ \hline
    $E_8$ & $A_3+A_1$ & $\frac{1}{2}(1,1,0,1,0,1,1,2)$  \\ \hline
    $E_8$ & $2A_2+2A_1$ & $\frac{1}{3}(1,1,1,1,1,1,1,1)$  \\ \hline
    $E_8$ & $A_3+2A_1$ & $\frac{1}{2}(1,1,1,0,1,0,1,1)$ \\ \hline
    $E_8$ & $D_4(a_1)+A_1$ & \cellcolor{blue!20}$(0,0,0,1,0,0,1,0)$\\ \hline
    $E_8$ & $A_3+A_2+A_1$ & $\frac{1}{2}(1,0,0,1,0,1,1,1)$  \\ \hline
    $E_8$ & $2A_3$ & $\frac{1}{4}(1,1,1,1,1,1,1,1)$  \\ \hline
    $E_8$ & $A_4+A_3$ & $\frac{1}{5}(1,1,1,1,1,1,1,1)$ \\ \hline
    $E_8$ & $A_5+A_1$ & $\frac{1}{6}(2,2,1,1,1,1,1,1)$  \\ \hline
    $E_8$ & $D_5(a_1)+A_2$ & $\frac{1}{4}(1,1,1,0,1,1,1,1)$  \\ \hline
    \end{tabular}

\end{table}

\begin{table}[hbt!]
\tiny
\caption{Unipotent infinitesimal characters attached to nilpotent covers for $\mathrm{Sp}(8)$. Characters are denoted in standard coordinates. Subscripts indicate the degree of the cover. Special unipotent infinitesimal characters are highlighted in blue. See Section \ref{subsec:unipotentcentralchars} for further explanation.}\label{table:Sp8}
\begin{tabular}{|c|c|c|c|c|c|c|c|}
\hline
$\widetilde{\mathbb{O}}$ & $L$ & $\widetilde{\mathbb{O}}_L$ & $K$ & $\mathbb{O}_K$ & $\gamma_0(\mathbb{O}_K)$ & $\delta$ & $\gamma_0(\widetilde{\mathbb{O}})$ \\ \hline
$(8)$ & $\mathrm{GL}(1)^4$ & $\{0\}$ & $\mathrm{GL}(1)^4$ & $\{0\}$ & $(0,0,0,0)$ & $(0,0,0,0)$ & $\cellcolor{blue!20} (0,0,0,0)$ \\ \hline
$(8)_2$ & $\mathrm{GL}(1)^3 \times \mathrm{Sp}(2)$ & $(2)_2$ & $\mathrm{GL}(1)^4$ & $\{0\}$ & $(0,0,0,0)$ & $(0,0,0,\frac{1}{2})$ & $(\frac{1}{2},0,0,0)$ \\ \hline
$(6,2)$ & $\mathrm{GL}(2) \times \mathrm{GL}(1)^2$ & $\{0\}$ & $\mathrm{GL}(2) \times \mathrm{GL}(1)^2$ & $\{0\}$ & $(\frac{1}{2}, -\frac{1}{2},0,0)$ & $(0,0,0,0)$ & \cellcolor{blue!20} $(\frac{1}{2},\frac{1}{2},0,0)$ \\ \hline
$(6,2)_2$ & $\mathrm{GL}(2) \times \mathrm{GL}(1) \times \mathrm{Sp}(2)$ & $(2)_2$ & $\mathrm{GL}(2) \times \mathrm{GL}(1)^2$ & $\{0\}$ & $(\frac{1}{2},-\frac{1}{2},0,0)$ & $(0,0,0,\frac{1}{2})$ & $(\frac{1}{2},\frac{1}{2}, \frac{1}{2},0)$ \\ \hline
$(6,2)_2$ & $\mathrm{GL}(1)^3 \times \mathrm{Sp}(2)$ & $\{0\}$ & $\mathrm{GL}(1)^3 \times \mathrm{Sp}(2)$ & $\{0\}$ & $(0,0,0,1)$ & $(0,0,0,0)$ & \cellcolor{blue!20} $(1,0,0,0)$ \\ \hline
$(6,2)_2$ & $\mathrm{GL}(1) \times \mathrm{Sp}(6)$ & $(4,2)_2$ & $\mathrm{GL}(1)^2 \times \mathrm{GL}(2)$ & $\{0\}$ & $(0,0,\frac{1}{2},-\frac{1}{2})$ & $(0,\frac{1}{2},\frac{1}{2},\frac{1}{2})$ &  $(1,\frac{1}{2},0,0)$ \\ \hline
$(6,2)_4$ & $\mathrm{GL}(1) \times \mathrm{Sp}(6)$ & $(4,2)_4$ & $\mathrm{GL}(1)^2 \times \mathrm{GL}(2)$ & $\{0\}$ & $(0,0,\frac{1}{2},-\frac{1}{2})$ & $(0, \frac{1}{2},\frac{1}{2},\frac{1}{2})$ &   $(1,\frac{1}{2}, 0, 0)$ \\ \hline
$(6,1^2)$ & $\mathrm{GL}(1)^2 \times \mathrm{Sp}(4)$ & $(2,1^2)$ & $\mathrm{GL}(1)^2 \times \mathrm{Sp}(4)$ & $(2,1^2)$ & $(0,0,\frac{3}{2},\frac{1}{2})$ & $(0,0,0,0)$ & $(\frac{3}{2},\frac{1}{2},0,0)$ \\ \hline
$(6,1^2)_2$ & $\mathrm{GL}(1)^2 \times \mathrm{Sp}(4)$ & $(2,1^2)_2$ & $\mathrm{GL}(1)^2 \times \mathrm{Sp}(4)$ & $(2,1^2)$ & $(0,0,\frac{3}{2},\frac{1}{2})$ & $(0,0,0,0)$ & $(\frac{3}{2},\frac{1}{2},0,0)$ \\
\hline

$(4^2)$ & $\mathrm{GL}(2)^2$ & $\{0\}$ & $\mathrm{GL}(2)^2$ & $\{0\}$ & $(\frac{1}{2},-\frac{1}{2},\frac{1}{2},-\frac{1}{2})$ & $(0,0,0,0)$ & \cellcolor{blue!20} $(\frac{1}{2},\frac{1}{2},\frac{1}{2},\frac{1}{2})$ \\ \hline
$(4^2)_2$ & $\mathrm{GL}(2) \times \mathrm{GL}(1) \times \mathrm{Sp}(2)$ & $\{0\}$ & $\mathrm{GL}(2) \times \mathrm{GL}(1) \times \mathrm{Sp}(2)$ & $\{0\}$ & $(\frac{1}{2},-\frac{1}{2}, 0, 1)$ & $(0,0,0,0)$ & \cellcolor{blue!20} $(1,\frac{1}{2},\frac{1}{2},0)$ \\ \hline
$(4,2^2)$ & $\mathrm{GL}(3) \times \mathrm{GL}(1)$ & $\{0\}$ & $\mathrm{GL}(3) \times \mathrm{GL}(1)$ & $\{0\}$ & $(1,0,-1,0)$ & $(0,0,0,0)$ & \cellcolor{blue!20} $(1,1,0,0)$\\ \hline

$(4,2^2)_2$ & $\mathrm{GL}(3) \times \mathrm{Sp}(2)$ & $(2)_2$ & $\mathrm{GL}(3) \times \mathrm{GL}(1)$ & $\{0\}$ & $(1,0,-1,0)$ & $(0,0,0,\frac{1}{2})$ & $(1,1,\frac{1}{2},0)$\\ \hline

$(4,2^2)_2$ & $\mathrm{GL}(2) \times \mathrm{Sp}(4)$ & $(2,1^2)$ & $\mathrm{GL}(2) \times \mathrm{Sp}(4)$ & $(2,1^2)$ & $(\frac{1}{2},-\frac{1}{2},\frac{3}{2},\frac{1}{2})$ & $(0,0,0,0)$ & $(\frac{3}{2},\frac{1}{2},\frac{1}{2},\frac{1}{2})$\\ \hline

$(4,2^2)_2$ & $\mathrm{GL}(1) \times \mathrm{Sp}(6)$ & $(2^3)_2$ & $\mathrm{GL}(1) \times \mathrm{GL}(3)$ & $\{0\}$ & $(0,1,0,-1)$ & $(0,\frac{1}{2},\frac{1}{2},\frac{1}{2})$ & $(\frac{3}{2},\frac{1}{2},\frac{1}{2},0)$ \\ \hline

$(4,2^2)_4$ & $\mathrm{GL}(2) \times \mathrm{Sp}(4)$ & $(2,1^2)_2$ & $\mathrm{GL}(2) \times \mathrm{Sp}(4)$ & $(2,1^2)$ & $(\frac{1}{2},-\frac{1}{2},\frac{3}{2},\frac{1}{2})$ & $(0,0,0,0)$ &  $(\frac{3}{2},\frac{1}{2},\frac{1}{2},\frac{1}{2})$\\ \hline

$(4,2,1^2)$ & $\mathrm{GL}(1) \times \mathrm{Sp}(6)$ & $(2^2,1^2)$ & $\mathrm{GL}(1) \times \mathrm{Sp}(6)$ & $(2^2,1^2)$ & $(0,2,1,0)$ & $(0,0,0,0)$ & \cellcolor{blue!20} $(2,1,0,0)$\\ \hline

$(4,2,1^2)_2$ & $\mathrm{GL}(1) \times \mathrm{Sp}(6)$ & $(2^2,1^2)_2$ & $\mathrm{GL}(1) \times \mathrm{Sp}(6)$ & $(2^2,1^2)$ & $(0,2,1,0)$ & $(0,0,0,0)$ & \cellcolor{blue!20} $(2,1,0,0)$\\ \hline

$(4,2,1^2)_2$ & $\mathrm{GL}(1)^2 \times \mathrm{Sp}(4)$ & $\{0\}$ & $\mathrm{GL}(1)^2 \times \mathrm{Sp}(4)$ & $\{0\}$ & $(0,0,2,1)$ & $(0,0,0,0)$ & \cellcolor{blue!20} $(2,1,0,0)$\\ \hline

$(4,2,1^2)_2$ & $\mathrm{Sp}(8)$ & $(4,2,1^2)_2$ & $\mathrm{GL}(1) \times \mathrm{Sp}(6)$ & $(2^2,1^2)$ & $(0,2,1,0)$ & $(\frac{1}{2},0,0,0)$ & $(2,1,\frac{1}{2},0)$\\ \hline

$(4,2,1^2)_4$ & $\mathrm{Sp}(8)$ & $(4,2,1^2)_4$ & $\mathrm{GL}(1) \times \mathrm{Sp}(6)$ & $(2^2,1^2)$ & $(0,2,1,0)$ & $(\frac{1}{2},0,0,0)$ & $(2,1,\frac{1}{2},0)$\\ \hline

$(4,1^4)$ & $\mathrm{GL}(1) \times \mathrm{Sp}(6)$ & $(2,1^4)$ & $\mathrm{GL}(1) \times \mathrm{Sp}(6)$ & $(2,1^4)$ & $(0,\frac{5}{2},\frac{3}{2},\frac{1}{2})$ & $(0,0,0,0)$ & $(\frac{5}{2},\frac{3}{2},\frac{1}{2},0)$\\ \hline

$(4,1^4)_2$ & $\mathrm{GL}(1) \times \mathrm{Sp}(6)$ & $(2,1^4)_2$ & $\mathrm{GL}(1) \times \mathrm{Sp}(6)$ & $(2,1^4)$ & $(0,\frac{5}{2},\frac{3}{2},\frac{1}{2})$ & $(0,0,0,0)$ & $(\frac{5}{2},\frac{3}{2},\frac{1}{2},0)$\\ \hline

$(3^2,2)$ & $\mathrm{GL}(3) \times \mathrm{Sp}(2)$ & $\{0\}$ & $\mathrm{GL}(3) \times \mathrm{Sp}(2)$ & $\{0\}$ & $(1,0,-1,1)$ & $(0,0,0,0)$ & \cellcolor{blue!20} $(1,1,1,0)$\\ \hline

$(3^2,2)_2$ & $\mathrm{Sp}(8)$ & $(3^2,2)_2$ & $\mathrm{GL}(3) \times \mathrm{Sp}(2)$ & $\{0\}$ & $(1,0,-1,1)$ & $(\frac{1}{2},\frac{1}{2},\frac{1}{2},0)$ & $(\frac{3}{2},1,\frac{1}{2},\frac{1}{2})$\\ \hline

$(3^2,1^2)$ & $\mathrm{GL}(2) \times \mathrm{Sp}(4)$ & $\{0\}$ & $\mathrm{GL}(2) \times \mathrm{Sp}(4)$ & $\{0\}$ & $(\frac{1}{2},-\frac{1}{2},2,1)$ & $(0,0,0,0)$ & \cellcolor{blue!20} $(2,1,\frac{1}{2},\frac{1}{2})$\\ \hline

$(2^4)$ & $\mathrm{GL}(4)$ & $\{0\}$ & $\mathrm{GL}(4)$ & $\{0\}$ & $(\frac{3}{2},\frac{1}{2},-\frac{1}{2},-\frac{3}{2})$ & $(0,0,0,0)$ & \cellcolor{blue!20} $(\frac{3}{2},\frac{3}{2},\frac{1}{2},\frac{1}{2})$\\ \hline

$(2^4)_2$ & $\mathrm{Sp}(8)$ & $(2^4)_2$ & $\mathrm{GL}(4)$ & $\{0\}$ & $(\frac{3}{2},\frac{1}{2},-\frac{1}{2},-\frac{3}{2})$ & $(\frac{1}{2},\frac{1}{2},\frac{1}{2},\frac{1}{2})$ & \cellcolor{blue!20} $(2,1,1,0)$\\ \hline

$(2^3,1^2)$ & $\mathrm{Sp}(8)$ & $(2^3,1^2)$ & $\mathrm{Sp}(8)$ & $(2^3,1^2)$ & $(\frac{5}{2},\frac{3}{2},\frac{1}{2},\frac{1}{2})$ & $(0,0,0,0)$ & $(\frac{5}{2},\frac{3}{2},\frac{1}{2},\frac{1}{2})$\\ \hline

$(2^3,1^2)_2$ & $\mathrm{Sp}(8)$ & $(2^3,1^2)_2$ & $\mathrm{Sp}(8)$ & $(2^3,1^2)$ & $(\frac{5}{2},\frac{3}{2},\frac{1}{2},\frac{1}{2})$ & $(0,0,0,0)$ & $(\frac{5}{2},\frac{3}{2},\frac{1}{2},\frac{1}{2})$ \\ \hline

$(2^2,1^4)$ & $\mathrm{Sp}(8)$ & $(2^2,1^2)$ & $\mathrm{Sp}(8)$ & $(2^2,1^2)$ & $(3,2,1,0)$ & $(0,0,0,0)$ & \cellcolor{blue!20} $(3,2,1,0)$\\ \hline

$(2^2,1^4)_2$ & $\mathrm{Sp}(8)$ & $(2^2,1^2)_2$ & $\mathrm{Sp}(8)$ & $(2^2,1^2)$ & $(3,2,1,0)$ & $(0,0,0,0)$ & \cellcolor{blue!20} $(3,2,1,0)$\\ \hline

$(2,1^6)$ & $\mathrm{Sp}(8)$ & $(2,1^6)$ & $\mathrm{Sp}(8)$ & $(2,1^6)$ & $(\frac{7}{2},\frac{5}{2},\frac{3}{2},\frac{1}{2})$ & $(0,0,0,0)$ & $(\frac{7}{2},\frac{5}{2},\frac{3}{2},\frac{1}{2})$\\ \hline

$(2,1^6)_2$ & $\mathrm{Sp}(8)$ & $(2,1^6)_2$ & $\mathrm{Sp}(8)$ & $(2,1^6)$ & $(\frac{7}{2},\frac{5}{2},\frac{3}{2},\frac{1}{2})$ & $(0,0,0,0)$ & $(\frac{7}{2},\frac{5}{2},\frac{3}{2},\frac{1}{2})$\\ \hline

$(1^8)$ & $\mathrm{Sp}(8)$ & $\{0\}$ & $\mathrm{Sp}(8)$ & $\{0\}$ & $(4,3,2,1)$ & $(0,0,0,0)$ & \cellcolor{blue!20} $(4,3,2,1)$\\ \hline
\end{tabular}

\end{table}

\chapter{Index of Notation}\label{sec:notation}

\begin{longtable}{l l l}
$G^{\circ}$ & Identity component of algebraic group $G$ & \S \ref{subsec:nilpcovers} \\
$G_x$ & Stabilizer of $x$ under group action of $G$ & \S \ref{subsec:nilpcovers} \\
$\Aut_{\widetilde{\OO}}(\widehat{\OO})$ & Galois group of covering map $\widehat{\OO} \to \widetilde{\OO}$ & \S \ref{subsec:nilpcovers} \\
$\geq$ & Partial order on nilpotent covers defined by covering relation & \S \ref{subsec:nilpcovers} \\
$\pi_1^G(\widetilde{\OO})$ & $G$-equivariant fundamental group of nilpotent cover $\widetilde{\OO}$  & \S \ref{subsec:nilpcovers} \\
$\Ind^G_M$ & Lusztig-Spaltenstein induction of nilpotent orbits  &  \S \ref{subsec:LSinduction} \\
$p+q$ & Row-wise sum of partitions $p$ and $q$ & \S \ref{subsec:LSinduction} \\
$\mathrm{Bind}^G_M$ & Birational induction of nilpotent covers & \S \ref{subsec:birationalinduction} \\
$\mathrm{Sat}^G_M$ & Bala-Carter inclusion of nilpotent orbits & \S \ref{subsec:BCinclusion} \\
$P_{\OO}, L_{\OO}$ & Jacobson-Morozov parabolic, Levi attached to a nilpotent orbit $\OO$ & \S \ref{subsec:BCinclusion} \\
$p \cup q$ & Concatenation of partitions $p$ and $q$ & \S \ref{subsec:BCinclusion} \\
$\fZ(\fg)$ & Center of universal enveloping algebra $U(\fg)$ & \S \ref{subsec:assvar} \\
$V(I)$ & Associated variety of primitive ideal $I$ & \S \ref{subsec:assvar} \\
$\Prim_{\gamma}(U(\fg))$ & Set of primitive ideals in $U(\fg)$ with infinitesimal character $\gamma \in \fh^*/W$ & \S \ref{subsec:assvar} \\ 
$I_{\mathrm{max}}(\gamma)$ & Maximal ideal in $U(\fg)$ with infinitesimal character $\gamma \in \fh^*/W$ & \S \ref{subsec:BVduality}\\
$G^{\vee}, \fg^{\vee}$ & Langlands dual of $G$, $\fg$ & \S \ref{subsec:BVduality} \\
$\mathsf{D}$ & BVLS duality & \S \ref{subsec:BVduality} \\
$p^t$ & Transpose of partition $p$ & \S \ref{subsec:BVduality} \\
$C(p), B(p)$ & C-collapse, B-collapse of partition $p$ & \S \ref{subsec:BVduality} \\
$l(p)$ & Partition obtained by removing a single box from the last row of $p$ & \S \ref{subsec:BVduality} \\
$e(p)$ & Partition obtained by appending $1$ to $p$ & \S \ref{subsec:BVduality} \\
$\fl^{\vee}_{\gamma}$ & Reductive subalgebra of $\fg^{\vee}$ corresponding to integral roots for $\gamma \in \fh^*$ & \S \ref{subsec:assvarmax}\\
$\fl^{\vee}_{\gamma,0}$ & Levi subalgebra of $\fg^{\vee}$ corresponding to singular roots for $\gamma \in \fh^*$ & \S \ref{subsec:assvarmax}\\
$J^G_M$ & Truncated induction of nilpotent orbits from endoscopic group $M$ & \S \ref{subsec:assvarmax}\\
$\HC^G(U(\fg))$ & Category of Harish-Chandra bimodules for reductive group $G$ & \S \ref{subsec:HCbimodsclassical} \\
$\mathcal{V}(\cB)$ & Associated variety of Harish-Chandra bimodule $\cB$ & \S \ref{subsec:HCbimodsclassical} \\
$\mathrm{LAnn}(\cB), \mathrm{RAnn}(\cB)$ & Left and right annihilators of Harish-Chandra bimodule $\cB$ & \S \ref{subsec:HCbimodsclassical} \\
$m_Z(\cB)$ & Multiplicity of Harish-Chandra bimodule $\cB$ along &  \S \ref{subsec:HCbimodsclassical} \\
& irreducible component $Z \subset \mathcal{V}(\cB)$ & \\
$\fX(\bullet)$ & One-dimensional representations of group or Lie algebra & \S \ref{subsec:HCbimodsclassical} \\ 
$\CC(\lambda,\nu)$ & 1-dimensional bimodule corresponding to weights $\lambda,\nu \in \fX(\fg)$ & \S \ref{subsec:HCbimodsclassical} \\
$\Ind^G_M$ & (Normalized) parabolic induction of Harish-Chandra bimodules & \S \ref{subsec:HCbimodsclassical} \\
$I(\lambda_{\ell},\lambda_r)$ & Induced bimodule corresponding to Langlands parameter $(\lambda_{\ell},\lambda_r)$ & \S \ref{subsec:HCbimodsclassical} \\
$\overline{I}(\lambda_{\ell},\lambda_r)$ & Langlands subquotient of $I(\lambda_{\ell},\lambda_r)$ & \S \ref{subsec:HCbimodsclassical} \\
$R$ & Reductive part of centralizer of $e \in \OO$ & \S \ref{subsec:W}\\
$\HC^R_{\mathrm{fin}}(\cW)$ & Category of $R$-equivariant Harish-Chandra bimodules for & \S \ref{subsec:W} \\
&  $W$-algebra $\cW$ & \\
$\mathrm{Prim}_{\mathrm{fin}}(\cW)$ & Set of primitive ideals of finite codimension in $\cW$ & \S \ref{subsec:W} \\
$\mathrm{Id}_{\mathrm{fin}}(\cW)$ & Set of ideals of finite-codimension in $\cW$ & \S \ref{subsec:W} \\
$\bullet_{\dagger}$ & Restriction functor $\HC^G_{\overline{\OO}}(U(\fg)) \to \HC^R_{\mathrm{fin}}(\cW)$ & \S \ref{subsec:W} \\
$\bullet^{\dagger}$ & Extension functor $\HC^R_{\mathrm{fin}}(\cW) \to \HC^G_{\overline{\OO}}(U(\fg))$ & \S \ref{subsec:W} \\
$I_{\dagger} \in \mathrm{Id}_{\mathrm{fin}}(\cW)$ & Restriction of primitive ideal $I \in \Prim_{\overline{\OO}}(U(\fg))$ & \S \ref{subsec:W} \\
$J^{\ddag} \in \Prim_{\overline{\OO}}(U(\fg))$ & Extension of primitive ideal $J \in \Prim_{\mathrm{fin}}(\cW)$ & \S \ref{subsec:W} \\
$\cW-\dim(I)$ & $\cW$-dimension of primitive ideal $I \in \Prim_{\overline{\OO}}(U(\fg))$ & \S \ref{subsec:W} \\
$\mathrm{Quant}(A)$ & Set of isomorphism classes of filtered quantizations  & \S \ref{subsec:quant} \\
& of graded Poisson algebra $A$ & \\
$\mathrm{PDef}(A)$ &  Set of isomorphism classes of Poisson deformations & \S \ref{subsec:quant} \\
& of graded Poisson algebra variety $A$ & \\
$\Pic(X)$ & Picard group of algebraic variety $X$ & \S \ref{subsec:Mckay} \\
$X^{\mathrm{reg}},X^{\mathrm{sing}}$ & Regular, singular loci of algebraic variety $X$ & \S \ref{subsec:symplectic} \\
$\fP_{\RR}^X, \fP^X$ & Real, complex Namikawa space associated to & \S \ref{subsec:Qfactorial} \\
& conical symplectic singularity $X$ & \\
$c_1(\mathcal{L}) \in H^2(X,\RR)$ & First Chern class of line bundle $\mathcal{L}$ on smooth manifold $X$ & \S \ref{subsec:Qfactorial}\\
$\fL_k \subset X$ & Codimension 2 leaf of conical symplectic singularity $X$ & \S \ref{subsec:structurenamikawa} \\
$\Sigma_k,\fg_k,\fh_k,\Lambda_k,\Delta_k$ & Kleinian singularity, complex simple Lie algebra, Cartan subalgebra, & \S \ref{subsec:structurenamikawa}\\
& weight lattice, and root system associated to codimension 2 leaf &  \\
& $\fL_k \subset X$ & \\
$\fP_{\RR,k}^X,\fP_k^X$ & Real, complex partial Namikawa space associated to & \S \ref{subsec:structurenamikawa} \\
& codimension 2 leaf $\fL_k \subset X$ in conical symplectic singularity $X$ & \\
$W_k^X$ & Partial Namikawa Weyl group associated to codimension 2 leaf & \S \ref{subsec:structurenamikawa} \\
& $\fL_k \subset X$ in conical symplectic singularity $X$ & \\
$W^X$ & Namikawa Weyl group associated to conical symplectic singularity $X$ & \S \ref{subsec:structurenamikawa} \\
$X_{\mathrm{univ}}$ & Universal Poisson deformation of conical symplectic singularity $X$ & \S \ref{subsec:quantsymplectic} \\
$\mathcal{D}^{Y,\mathrm{univ}}$ & Canonical quantization of universal deformation $Y_{\mathrm{univ}}$ & \S \ref{subsec:quantsymplectic} \\
$\cA^0_{\lambda}, \cA_{\lambda}$ & Poisson deformation/filtered quantization of conical symplectic & \S \ref{subsec:quantsymplectic} \\
& singularity $X$ with parameter $\lambda \in \fP^X$ & \\
$\cN \subset \fg^*$ & Nilpotent cone of complex reductive Lie algebra $\fg$ & \S \ref{subsec:quantsymplectic} \\
$eH_ce$ & Filtered quantization of Kleinian singularity $\CC^2/\Gamma$ with & \S \ref{subsec:quantKleinian} \\
&  CBH parameter $c \in \CC[\Gamma]^{\Gamma}$ & \\
$\lambda^c \in \fh^*$ & Parameter in $\fP^{\Sigma}$ associated to CBH parameter $c$ & \S \ref{subsec:quantKleinian} \\
$\fP_{\RR}^{\geq 0}$ & Fundamental chamber of Namikawa space & \S \ref{amplecones} \\
$\Pic^a(Y)$ & Semigroup of relatively ample line bundles on $\QQ$-factorial & \S \ref{amplecones}\\
&  terminalization $Y \to X$ & \\
$\mathrm{Amp}(Y) \subset \fP_{\RR}^{\geq 0}$ & Ample cone of $\QQ$-factorial terminalization $Y \to X$& \S \ref{amplecones}\\
$X_{\lambda}$ & Fiber of universal deformation of conical symplectic singularity $X$ & \S \ref{amplecones}\\
& over $\lambda \in \fP^X/W$ & \\
$\fP^{\mathrm{sing}}$ & Subset of $\fP^X$ consisting of $\lambda$ for which $Y_{\lambda} \to X_{\lambda}$ not an isomorphism & \S \ref{amplecones}\\
$\mathrm{Der}(A)$ & Lie algebra of derivations of Lie algebra $A$ & \S \ref{subsec:equivariant} \\
$\mathrm{Quant}^G(A)$ & Set of isomorphism classes of Hamiltonian quantizations of & \S \ref{subsec:equivariant} \\
& $G$-equivariant graded Poisson algebra $A$ & \\
$\overline{\fP}^X$ & Extended Namikawa space of conical symplectic singularity $X$ with & \S \ref{subsec:equivariant} \\
& Hamiltonian $G$-action & \\
$\HC(\cA)$ & Category of Harish-Chandra bimodules for filtered quantization $\cA$  & \S \ref{subsec:HCbimods} \\
& of conical symplectic singularity & \\
$\HC_{\partial}(\cA)$ & Subcategory of $\HC(\cA)$ consisting of $\cA$-bimodules of proper support & \S \ref{subsec:HCbimods} \\
$\overline{\HC}(\cA)$ & Quotient category $\HC(\cA)/\HC_{\partial}(\cA)$ & \S \ref{subsec:HCbimods} \\
$\bullet_{\dagger}$ & Functor $\HC(\cA) \to \Gamma\modd$ for quantization $\cA$ of conical symplectic & \S \ref{subsec:daggers} \\
& singularity $X$ with $\Gamma = \pi_1(X^{\mathrm{reg}})$ & \\
$\bullet^{\dagger}$ & Functor $\Gamma\modd \to \HC(\cA)$ & \S \ref{subsec:daggers} \\
$\Gamma(\lambda) \subset \Gamma$ & Normal subgroup corresponding to quantization parameter $\lambda \in \fP^X$ & \S \ref{subsec:daggers} \\
$\epsilon(\Gamma') \in \fh^*$ & Dominant weight for simple Lie algebra associated to Kleinian  & \S \ref{subsec:invariantsKleinian} \\
& singularity $\CC^2/\Gamma$ corresponding to normal subgroup $\Gamma' \subset \Gamma$ & \\
$\epsilon \in \fP^X$ & Element of Namikawa space corresponding to Galois cover $\widetilde{X} \to X$  & \S \ref{sec:invariantssymplectic} \\
& of conical symplectic singularities & \\
$I_{\lambda}(\widetilde{\OO})$ & Kernel of quantum comoment map $\Phi_{\lambda}: U(\fg) \to \cA_{\lambda}^{\widetilde{X}}$ corresponding to & \S \ref{sec:unipotent} \\
&  Hamiltonian quantization of $\widetilde{X}=\Spec(\CC[\widetilde{\OO}])$ with parameter $\lambda \in \overline{\fP}^{\widetilde{X}}$ & \\
$\succeq$ & Partial order on nilpotent covers defined by almost \'{e}tale relation & \S \ref{subsec:classificationideals} \\
$\sim$ & Equivalence relation on nilpotent covers generated by $\succeq$ & \S \ref{subsec:classificationideals} \\
$[\widetilde{\OO}]$ & Equivalence class of nilpotent cover $\widetilde{\OO}$ & \S \ref{subsec:classificationideals} \\
$\Cl(X)$ & Divisor class group of algebraic variety $X$ & \S \ref{subsec:picard} \\
$\mathrm{div}$ & Natural map $\Pic(X) \to \Cl(X)$ & \S \ref{subsec:picard} \\
$\eta$ & Isomorphism $\eta: \fX(\fl) \xrightarrow{\sim} \overline{\fP}^{\widetilde{X}}$ & \S \ref{subsec:terminalizationcover} \\
$\mathfrak{D}_X$ & Sheaf of differential operators on algebraic variety $X$ & \S \ref{subsec:inductionquantizations} \\
$\Ind^G_M$ & (Normalized) parabolic induction of Hamiltonian quantizations  & \S \ref{subsec:inductionquantizations} \\
& of nilpotent covers & \\
$M_k$ & Levi subgroup adapted to codimension 2 leaf $\fL_k \subset \Spec(\CC[\widetilde{\OO}])$ & \S \ref{subsec:descriptionpartial} \\
$\eta_k$ & Isomorphism $\fX(\fm_k) \xrightarrow{\sim} \fP_k$ associated to codimension 2 leaf  & \S \ref{subsec:descriptionpartial} \\
& $\fL_k \subset \Spec(\CC[\widetilde{\OO}])$ & \\
$\OO_k$ & Codimension 2 orbit in $\overline{\OO}$ corresponding to codimension 2 leaf & \S \ref{subsec:descriptionpartial} \\
&  $\fL_k \subset \Spec(\CC[\widetilde{\OO}])$ & \\
$\simeq^G$ & Conjugacy relation on Levi subgroups of $G$ & \S \ref{subsec:codim2leaves} \\
$S_2(p)$ & Set of jumps of size 2 in partition $p$ & \S \ref{subsec:codim2leaves} \\
$c_i(p)$ & Collapse of partition $p$ at $p_i$ & \S \ref{subsec:codim2leaves} \\
$S_4(p)$ & Set of jumps of size 4 in partition $p$ & \S \ref{subsec:codim2leaves} \\
$p \# x$ & Partition obtained from $p$ by deleting columns numbered & \S \ref{subsec:codim2leaves} \\
& $p_{x_1}, p_{x_1}-1, p_{x_2}, p_{x_2}-1, ...$ & \\
$\tau_1(k),...,\tau_{n(k)}(k)$ & Generators for free abelian group $\fX(M_k)$  & \S \ref{subsec:codim2leaves} \\
$G_{\mathrm{ab}}$ & Abelianization of group $G$ & \S \ref{subsec:codim2leaves} \\
$\mathcal{P}_{\mathrm{rig}}(\OO)$ & Conjugacy classes of pairs $(L,\OO_L)$ with $\OO_L$ rigid and $\OO=\mathrm{Ind}^G_L \OO_L$ & \S \ref{subsec:codim2leaves} \\
$m(\OO)$ & Maximum value of $\dim \fX(\fl)$ for $(L,\OO_L) \in \mathcal{P}_{\mathrm{rig}}(\OO)$ & \S \ref{subsec:codim2leaves} \\
$\omega_1(k),...,\omega_{n(k)}(k)$ & Fundamental weights for $\fg_k$, simple Lie algebra corresponding to $\Sigma_k$ & \S \ref{subsec:identification} \\
$\delta \in \fX(\fl)$ & Element of $\fX(\fl)$ corresponding to barycenter parameter $\epsilon \in \fP^X$  & \S \ref{subsec:identification} \\
& for $X=\Spec(\CC[\OO])$ & \\
$\gamma_{\lambda}(\widetilde{\OO})$ & infinitesimal character of primitive ideal $I_{\lambda}(\widetilde{\OO}) \subset U(\fg)$ & \S \ref{sec:centralchars} \\
$\rho(q),\rho^+(q)$ & Tuples attached to partition $q$ & \S \ref{subsec:centralcharclassical} \\
$Q(\OO)$ & $Q$-unipotent infinitesimal characters attached to $\OO$ & \S \ref{subsec:centralcharclassical} \\
$x(q),y(q)$ & Subpartitions of $q$ consisting of multiplicity 1 (resp. 2) parts & \S \ref{subsec:centralcharclassical} \\
$|p|$ & Size of partition $p$ & \S \ref{subsec:centralcharclassical} \\
$\widetilde{\mathsf{D}}$ & Refined BVLS duality & \S \ref{sec:duality} \\
$\Ind^G_M[\chi]$ & Parabolic induction of Harish-Chandra bimodules for Hamiltonian & \S \ref{subsec:bimodinductioncovers} \\
& quantizations of nilpotent covers & \\
$G_{\mathrm{tor}}$ & Torsion subgroup of an abelian group $G$ & \S \ref{subsec:bimodinductioncovers} \\
\end{longtable}

\newpage
\thispagestyle{empty}

\backmatter

\addtocontents{toc}{{\string\vskip1\baselineskip}}%
\addtocontents{toc}{\string\enlargethispage{-24pt}}%

%
%

\bibliography{bibsamp}

\printindex

\end{document}